\documentclass[12pt,oldfontcommands, openany]{memoir}
\usepackage{mystyle}

\makeindex
\shipout\null
\frontmatter
\title{The Ultrapower Axiom}
\date{}
\author{ \\ A dissertation presented\\
	by\\
	Gabriel Goldberg\\
to\\
The Harvard Department of Mathematics\\\\
in partial fulfillment of the requirements\ for the degree of\\
Doctor of Philosophy\\
in the subject of\\
Mathematics\\\\
Harvard University\\
Cambridge, Massachusetts\\
April 2019}
\begin{document}
\pagestyle{plain}
\maketitle
\thispagestyle{empty}
\newpage
\noindent \copyright\ 2019 Gabriel Goldberg\\
All rights reserved.
\thispagestyle{empty}
\newpage
Dissertation Advisor: Professor Hugh Woodin \hfill Gabriel Goldberg\\\\\\
\begin{center}{\large The Ultrapower Axiom}\end{center}
\begin{abstract}
The inner model problem for supercompact cardinals, one of the central open
problems in modern set theory, asks whether there is a canonical model of set
theory with a supercompact cardinal. The problem is closely related to the more
precise question of the equiconsistency of strongly compact cardinals and
supercompact cardinals. This dissertation approaches these two problems
abstractly by introducing a principle called the Ultrapower Axiom which is
expected to hold in all known canonical models of set theory. By investigating
the consequences of the Ultrapower Axiom under the hypothesis that there is a
supercompact cardinal, we provide evidence that the inner model problem can be
solved. Moreover, we establish that under the Ultrapower Axiom, strong
compactness and supercompactness are essentially equivalent.
\end{abstract}
\newpage
\tableofcontents

\section*{Acknowledgements}
I want to thank Peter Koellner, who taught my first course in set theory. Since
then, Peter has profoundly shaped my view of the subject. I would not have made
it through graduate school without his support and friendship.

Any reader will see that this work owes a great debt to the mathematics of Hugh
Woodin. I should also acknowledge my personal debt. I could not have asked for a
better advisor. Hugh was incredibly generous with his time, and was a constant
source of optimism and encouragement over the course of this research.

I want to thank Doug Blue for helping me secure housing a week before this
school year began, Akihiro Kanamori for carefully reading this dissertation and
pointing out numerous errors, and John and Colleen Steel for hosting me for
three weeks at Berkeley.
\mainmatter
\chapter{Introduction}
The goal of inner model theory is to construct and analyze canonical models of set theory. The simplest example of such a model is G\"odel's constructible universe \(L\), the smallest model of set theory that contains every ordinal number. One sense in which \(L\) is canonical is that seemingly every question about the internal structure of \(L\) can be answered. For example, G\"odel proved that \(L\) satisfies the Generalized Continuum Hypothesis. This stands in stark contrast with the universe of sets \(V\), many of whose most basic properties (for example, whether the Continuum Hypothesis holds) cannot be determined in any commonly accepted axiomatic system. 

To what extent does \(L\) provide a good approximation to the universe of sets? On the one hand, the principle that every set belongs to \(L\) (or in other words, \(V= L\)) cannot be refuted using the ZFC axioms, since \(L\) itself is a model of the theory \(\text{ZFC}+ V = L\). If \(V = L\), then \(L\) approximates \(V\) {\it very }well. 
On the other hand, the model \(L\) fails to satisfy relatively weak large cardinal axioms. If one takes the stance that these large cardinal axioms are true in the universe of sets, one must conclude that \(V \neq L\). Moreover, it follows from large cardinal axioms that \(L\) constitutes only a tiny fragment of the universe of sets. For example, assuming large cardinal axioms, the set of real numbers that lie in \(L\) is countable.

Are there canonical models generalizing \(L\) that yield better approximations to \(V\)?  It turns out that there is a hierarchy of canonical models beyond \(L\), satisfying stronger and stronger large cardinal axioms. The program of building such models has met striking success, reaching large cardinal axioms as strong as a Woodin limit of Woodin cardinals. Based on the pattern that has emerged so far, it is plausible every large cardinal axiom has a canonical model. 

At present, however, a vast expanse of large cardinal axioms are not yet known to admit canonical models. A key target problem for inner model theory is the construction of a canonical model with a supercompact cardinal. Work of Woodin suggests that the solution to this problem alone will yield an \emph{ultimate canonical model} that inherits essentially all large cardinals present in the universe. There is therefore hope that the goal of constructing inner models for all large cardinal axioms could might be achieved in a single stroke. If this is possible, the resulting model would be of enormous set theoretic interest, since it would closely approximate the universe of sets and yet admit an analysis that is as detailed as that of G\"odel's \(L\).

This dissertation investigates whether there can be a canonical model with a supercompact cardinal. To do this, we develop an abstract approach to inner model theory. This is accomplished by introducing a combinatorial principle called the Ultrapower Axiom, which is expected to hold in all canonical models. If one could show that the Ultrapower Axiom is inconsistent with a supercompact cardinal, one would arguably have to conclude that there can be no canonical model with a supercompact cardinal.

Supplemented with large cardinal axioms, the Ultrapower Axiom turns out to have surprisingly strong and coherent consequences for the structure of the upper reaches of the universe of sets, particularly above the first supercompact cardinal. These consequences are entirely consistent with what one would expect to hold in a canonical model, yet are proved by methods that are completely different from the usual techniques of inner model theory. The coherence of this theory provides compelling evidence that the Ultrapower Axiom is consistent with a supercompact cardinal. If this is the case, it seems that the only possible explanation is that the canonical model for a supercompact cardinal does indeed exist. Optimistically, studying the consequences of the Ultrapower Axiom will shed light on how this model should be constructed. 
\subsection{Outline}
We now describe the main results of this dissertation. \\

\noindent {\bf \sc \cref{MOChapter}.} In this introductory chapter, we introduce UA in the context of the problem of the linearity of the Mitchell order on normal ultrafilters. We  show first that UA holds in all canonical inner models, a result that is philosophically central to this dissertation. More precisely, we prove that UA is a consequence of Woodin's Weak Comparison principle:

\begin{repthm}{WCtoUA}
	Assume that \(V = \textnormal{HOD}\) and there is a \(\Sigma_2\)-correct worldly cardinal. If Weak Comparison holds, then the Ultrapower Axiom holds.
\end{repthm}

We then show that UA implies the linearity of the Mitchell order:
\begin{repthm}{UAMO}[UA]
	The Mitchell order is linear.
\end{repthm}

Two applications of this result to longstanding problems of Solovay-Reinhardt-Kanamori \cite{Kanamori} are explained in the introduction to \cref{MOChapter}.\\

\noindent  {\bf \sc \cref{KetonenChapter}.} This chapter introduces the Ketonen order, a generalization of the Mitchell order to all countably complete ultrafilters on ordinals. The restriction of this order to weakly normal ultrafilters was originally introduced by Ketonen. The first proof of the wellfoundedness of the generalization of this order to countably complete ultrafilters is due to the author:
\begin{repthm}{KOWellfounded}
	The Ketonen order is wellfounded.
\end{repthm}

The main theorem of this chapter explains the fundamental role of the Ketonen order in applications of the Ultrapower Axiom:
\begin{repthm}{LinKetThm}
	The Ultrapower Axiom is equivalent to the linearity of the Ketonen order.
\end{repthm}
In addition, we analyze the relationship between the Ketonen order and various well-known orders like the Rudin-Keisler order and the Mitchell order. \\

\noindent {\bf \sc \cref{GMOChapter}.} The topic of this chapter is the generalized Mitchell order, which is defined in exactly the same way as the usual Mitchell order on normal ultrafilters but removing the requirement that the ultrafilters involved be normal. This order is not linear (assuming there is a measurable cardinal), and in fact it is quite pathological when considered on ultrafilters in general. The two main results of this chapter generalize the linearity of the Mitchell order to nice classes of ultrafilters:
\begin{repthm}{DoddMO}[UA]
	The generalized Mitchell order is linear on Dodd sound ultrafilters.
\end{repthm}
Dodd soundness is a generalization of normality that was first isolated in the context of inner model theory by Steel \cite{Schimmerling}. A uniform ultrafilter \(U\) on a cardinal \(\lambda\) is Dodd sound if the map \(h : P(\lambda)\to M_U\) defined by \(h(X) = j_U(X)\cap [\text{id}]_U\) belongs to \(M_U\). The concept is discussed at great length in \cref{DoddSection}.

A better-known generalization of normality is the concept of a normal fine ultrafilter (\cref{NormalFineDef}), introduced by Solovay, and underpinning the theory of supercompact cardinals. The second result of this chapter generalizes the linearity of the Mitchell order to this class of ultrafilters:
\begin{repthm}{GCHLinear}[UA]
	Suppose \(\lambda\) is a cardinal such that \(2^{<\lambda} = \lambda\). Then the generalized Mitchell order is linear on normal fine ultrafilters on \(P_{\textnormal{bd}}(\lambda)\).
\end{repthm}
Here \(P_{\textnormal{bd}}(\lambda)\) denotes the set of bounded subsets of \(\lambda\).\\

\noindent {\bf \sc \cref{RFChapter}.} We turn to another fundamental order on ultrafilters, the Rudin-Frol\'ik order. The structure of the Rudin-Frol\'ik order on countably complete ultrafilters is intimately related to the Ultrapower Axiom. For example, we point out the following simple connection:
\begin{repcor}{RFDirected}\index{Rudin-Frol\'ik order!directedness}
	The Ultrapower Axiom holds if and only if the Rudin-Frol\'ik order is directed on countably complete ultrafilters.
\end{repcor}
On the other hand, it is well-known that the Rudin-Frol\'ik order is not directed on ultrafilters on \(\omega\).

The chapter is devoted to deriving deeper structural features of the Rudin-Frol\'ik order from UA. The most interesting one is that it is locally finite:
\begin{repthm}{RFFinite}[UA]
	A countably complete ultrafilter has at most finitely many predecessors in the Rudin-Frol\'ik order up to isomorphism.
\end{repthm}
Given the finiteness of the Rudin-Frol\'ik order, it turns out to be possible to represent every ultrafilter as a finite iterated ultrapower consisting of \emph{ irreducible ultrafilters}, ultrafilters whose ultrapowers cannot be factored as an iterated ultrapower (\cref{FactorizationThm}). We apply this to analyze ultrafilters on the least measurable cardinal under UA:
\begin{repthm}{KunenUA}[UA]
	Every countably complete ultrafilter on the least measurable cardinal \(\kappa\) is isomorphic to \(U^n\) where \(U\) is the unique normal ultrafilter on \(\kappa\) and \(n\) is a natural number.
\end{repthm}
This generalizes a classic theorem of Kunen \cite{KunenLU}.\\

\noindent {\bf \sc \cref{GCHChapter}.} This chapter exposits two inner model principles that follow abstractly from UA in the presence of a supercompact cardinal:
\begin{repthm}{HODGen}[UA]
	Assume there is a supercompact cardinal. Then \(V\) is a generic extension of \(\textnormal{HOD}\).
\end{repthm}
Thus UA almost implies \(V = \text{HOD}\). This is best possible in the sense that it is consistent that UA holds and \(V\) is a \emph{ nontrivial} generic extension of HOD.

The main result of the chapter is that UA implies the Generalized Continuum Hypothesis:
\begin{repthm}{MainTheorem}[UA]
	Suppose \(\kappa\) is supercompact. Then for all \(\lambda\geq \kappa\), \(2^\lambda = \lambda^+\).
\end{repthm}
Thus UA almost implies the GCH. This is best possible in the sense that it is consistent that UA holds but CH fails.\\

\noindent {\bf \sc \cref{SCChapter1}.} This chapter initiates an analysis of strongly compact and supercompact cardinals under UA. In this chapter, we investigate the structure of the least strongly compact cardinal, introducing the ``least ultrafilters" \(\mathscr K_\lambda\), and proving that they witness its supercompactness:
\begin{repthm}{StrongSuper}[UA]
	The least strongly compact cardinal is supercompact.\\
\end{repthm}

\noindent {\bf \sc \cref{SCChapter2}.} The main result of this chapter is the \emph{ Irreducibility Theorem} (\cref{SuccessorIrredThm}, \cref{LimitIrredThm}) that relates supercompactness and irreducibility (that is, Rudin-Frol\'ik minimality) under UA. The original impetus for proving this theorem was to analyze all larger strongly compact cardinals:
\begin{repthm}{MenasUA}[UA]
	A cardinal \(\kappa\) is strongly compact if and only if it is supercompact or a measurable limit of supercompact cardinals.
\end{repthm}
We also analyze various other large cardinals using UA. For example, we consider huge cardinals (\cref{HugeThm}) and rank-into-rank cardinals (\cref{I3Thm}).
\chapter{The Linearity of the Mitchell Order}\label{MOChapter}
\section{Introduction}\label{MOIntroduction}
\subsection{Normal ultrafilters and the Mitchell Order}
Normal ultrafilters are among the simplest objects that arise from modern large
cardinal axioms, yet despite their apparent simplicity, and despite the past six
decades of remarkable progress in the theory of large cardinals, the class of
normal ultrafilters remains mysterious, its underlying structure inextricably
bound up with some of the deepest and most difficult problems in set theory. The
following questions, posed by
Solovay-Reinhardt-Kanamori\index{Solovay-Reinhardt-Kanamori questions}
\cite{Kanamori} in the 1970s, are among the most prominent open questions in
this subject:

\begin{qst}\label{SRK1}
	Assume \(\kappa\) is \(2^\kappa\)-supercompact. Must there be more than one
	normal ultrafilter on \(\kappa\) concentrating on nonmeasurable cardinals?
\end{qst}

\begin{qst}\label{SRK2}
	Assume \(\kappa\) is strongly compact. Must \(\kappa\) carry more than one
	normal ultrafilter?
\end{qst}

These questions turn out to be merely the most concrete instances of a sequence
of more and more general structural questions in the theory of large cardinals.
Let us start down this path by stating a conjecture that would answer both
questions at once:

\begin{conj}\label{UniqueConj}
	It is consistent with an extendible cardinal that every measurable cardinal
	carries a unique normal ultrafilter concentrating on nonmeasurable
	cardinals.
\end{conj}

This would obviously answer \cref{SRK1} negatively, but what bearing does it
have on \cref{SRK2}? Assume \(\kappa\) is extendible and every measurable
cardinal carries a unique normal ultrafilter concentrating on nonmeasurable
cardinals. Consider the least strongly compact cardinal \(\kappa\) that is a
limit of strongly compact cardinals. By a theorem of Menas \cite{Menas} (proved
here as \cref{Menas}), the set of measurable cardinals below \(\kappa\) is
nonstationary. It follows that every normal ultrafilter on \(\kappa\)
concentrates on nonmeasurable cardinals. Since we assumed there is only one such
ultrafilter, \(\kappa\) is a strongly compact cardinal that carries a unique
ultrafilter. \cref{UniqueConj} thus supplies a negative answer to \cref{SRK2} as
well.

Why would someone make \cref{UniqueConj}? To answer this question, we must
consider the broader question of the structure of the Mitchell order under large
cardinal hypotheses. Recall that if \(U\) and \(W\) are normal ultrafilters, the
{\it Mitchell order} is defined by setting \(U\mo W\) if \(U\) belongs to the
ultrapower of the universe by \(W\). It is not hard to see that a normal
ultrafilter \(U\) on a cardinal \(\kappa\) concentrates on nonmeasurable
cardinals if and only if \(U\) is a minimal element in the Mitchell order on
normal ultrafilters on \(\kappa\). The following conjecture therefore
generalizes \cref{UniqueConj}:

\begin{conj}\label{MOLinearQ}
It is consistent with an extendible cardinal that the Mitchell order is
linear.\index{Mitchell order!linearity}
\end{conj}
How could one possibly prove \cref{MOLinearQ}? The most general technique for
proving consistency results, namely forcing, seems to be powerless in this
instance. To force the linearity of the Mitchell order, one would in particular
have to force that the least measurable cardinal carries a unique normal
ultrafilter, but even this much more basic problem remains open.

Kunen \cite{KunenLU} famously did prove that it is consistent for the least
measurable cardinal to carry a unique normal ultrafilter, not by forcing but
instead by building an inner model. In fact, he showed that if \(U\) is a normal
ultrafilter on a cardinal \(\kappa\), then the inner model \(L[U]\) satisfies
that \(\kappa\) is the unique measurable cardinal and
\(U\cap L[U]\) is the unique normal ultrafilter on \(\kappa\).
In an attempt to generalize Kunen's results, Mitchell
\cite{Mitchell} isolated the Mitchell order and used it to guide the 
construction of generalizations of \(L[U]\) that can have many measurable cardinals. 
Mitchell's proof of the linearity of the Mitchell order in these models
proceeds as follows:
\begin{itemize}
	\item Consider the model \(M = L[\langle U_\alpha: \alpha < \gamma\rangle]\) 
	built from a coherent sequence \(\langle U_\alpha: \alpha < \gamma\rangle\) 
	of normal ultrafilters.\footnote{
		Coherence is a key technical definition that includes the assumption 
		that \(\langle U_\alpha: \alpha < \gamma\rangle\) is increasing in the Mitchell order.}
		\begin{itemize}
			\item As a consequence of coherence, the sequence
			\(\langle U_\alpha\cap M: \alpha < \gamma\rangle\)
			is linearly ordered by the Mitchell order in \(M\).
		\end{itemize}
	\item Show that every normal ultrafilter of \(M\)
	appears on the sequence \(\langle U_\alpha\cap M: \alpha < \gamma\rangle\).
\end{itemize}

In the decades since Mitchell's result, inner model theory has ascended much
farther into the large cardinal hierarchy.  Combining the results of many
researchers (especially Neeman \cite{Neeman} and Schlutzenberg
\cite{Schlutzenberg}), the following is the best partial result towards
\cref{MOLinearQ} to date:
\begin{thm*}
	If it is consistent that there is a Woodin limit of Woodin cardinals, then
	the linearity of the Mitchell order is consistent with a Woodin limit of
	Woodin cardinals.
\end{thm*}
The linearity proof, due to Schlutzenberg, is much harder, but the argument
still roughly follows Mitchell's:
\begin{itemize}
	\item Consider the model \(M = L[\langle E_\alpha: \alpha < \gamma\rangle]\) built from a coherent extender sequence \(\langle E_\alpha: \alpha < \gamma\rangle\).\begin{itemize}
	\item By the definition of a coherent extender sequence, \(\langle E_\alpha: \alpha < \gamma\rangle\) is linearly ordered by the Mitchell order in \(M\).\end{itemize}
	\item Show that every normal ultrafilter of \(M\) appears on the sequence
	\(\langle E_\alpha: \alpha < \gamma\rangle\).
\end{itemize}

By now, it may appear that \cref{MOLinearQ} itself is subsumed by the far more
important (but far less precise) Inner Model Problem:\index{Inner Model Problem}
\begin{conj}\label{IMP}
	There is a canonical inner model with an extendible cardinal.
\end{conj}

The relationship between \cref{MOLinearQ} and \cref{IMP} is actually not as
straightforward as one might expect, because if inner model theory can be
extended to the level of extendible cardinals, the models must be significantly
different from the current models. For example, take the Woodin and Neeman-Steel
models with long extenders, which are canonical inner models designed to accommodate
large cardinals at the finite levels of supercompactness. It is not known
whether the constructions actually succeed, but the following conjecture is
plausible:
\begin{conj}\label{WoodinConj}
If for all \(n < \omega\), there is a cardinal \(\kappa\) that is
\(\kappa^{+n}\)-supercompact, then for all \(n < \omega\), there is an iterable
Woodin model with a cardinal \(\kappa\) that is \(\kappa^{+n}\)-supercompact.
\end{conj}
Given the pattern described above, one might expect to generalize Mitchell and
Schlutzenberg's results to the Woodin models, and therefore obtain for any \(n <
\omega\), the consistency of the linearity of the Mitchell order with a cardinal
\(\kappa\) that is \(\kappa^{+n}\)-supercompact. But there is a catch: the
proofs of these theorems cannot generalize verbatim to this level.
\begin{prp}
If \(L[\mathbb E]\) is an iterable Woodin model satisfying that \(\kappa\) is
\(\kappa^{++}\)-supercompact, then in \(L[\mathbb E]\), there is a normal
ultrafilter on \(\kappa\) that does not lie on the coherent sequence \(\mathbb
E\).
\end{prp}
(The same result applies to the Neeman-Steel models at this level.)
Therefore Mitchell's proof of the consistency of the linearity of the Mitchell
order cannot extend to the level of a cardinal \(\kappa\) that is
\(\kappa^{++}\)-supercompact. This might have been taken as a reason for
skepticism about \cref{MOLinearQ}.
\subsection{The Ultrapower Axiom}
The problem of generalizing the linearity of the Mitchell order to canonical
inner models at the finite levels of supercompactness was the original
inspiration for all the work in this dissertation. Our initial discovery was a
new argument that shows that {\it any} canonical inner model built by the
methodology of modern inner model theory must satisfy that the Mitchell order is
linear. The argument is extremely simple and relies on a single fundamental
property of the known canonical inner models: the {\it Comparison
Lemma}.\index{Comparison Lemma} The Comparison Lemma roughly states that any two
canonical inner models at the same large cardinal level can be embedded into a
common model. The inner model constructions are perhaps best viewed as an
attempt to build models satisfying the Comparison Lemma and accommodating large
cardinals.

Upon further reflection, it became clear that this argument relied solely on a
a single simple consequence of the Comparison Lemma, which could be distilled
into an abstract combinatorial principle.
This principle is called the {\it Ultrapower Axiom} (UA). 
(The definition appears in \cref{UAMOSection}.) The
Comparison Lemma implies that UA holds in all known canonical inner models.
Since the Comparison Lemma is fundamental to the methodology of inner model
theory, UA is expected to hold in any canonical inner model that will ever be
built.

Our theorem on the linearity of the Mitchell order now reads:
\begin{repthm}{UAMO}
Assume the Ultrapower Axiom. Then the Mitchell order is linear.
\end{repthm}
Granting our contention that UA holds in every canonical inner model, we have
reduced \cref{MOLinearQ} to the Inner Model Problem (for example, \cref{IMP}).
Moreover, we can state a perfectly precise test question that seems to capture
the essence of the Inner Model Problem:
\begin{conj}\label{UAQ}
The Ultrapower Axiom is consistent with an extendible cardinal.
\end{conj}
It is our expectation that neither this conjecture nor even \cref{MOLinearQ}
will be proved without first solving the Inner Model Problem.
What sets \cref{UAQ} apart from similar test questions like \cref{MOLinearQ} 
is that UA turns has a host of structural consequences in the theory of large cardinals. By
studying UA, one can hope to glean insight into the inner models that have not
yet been built, or perhaps to refute their existence by refuting UA from a large
cardinal hypothesis. The latter has not happened. Instead a remarkable theory of
large cardinals under UA has emerged which in our opinion provides evidence for
\cref{UAQ} and hence for the existence of inner models for very large cardinals.

\subsection{Outline of \cref{MOChapter}}\label{MOOutline}
We now briefly outline the contents of the rest of this chapter.\\

\noindent {\sc \cref{PrelimsSection}.} This section contains preliminary definitions most
of which are standard or self-explanatory. The topics we cover are ultrapowers,
close embeddings, uniform ultrafilters, and normal ultrafilters.\\

\noindent {\sc \cref{LinMOSection}.} This section contains proofs of the linearity of the
Mitchell order and motivation for the Ultrapower Axiom. We begin in
\cref{WCSection} by introducing and motivating Woodin's Weak Comparison axiom.
Then in \cref{WCMOSection}, we give our original argument for the linearity of
the Mitchell order under Weak Comparison (\cref{WCtoMO}). In \cref{WCUASection},
we abstract from this argument the Ultrapower Axiom, the central principle in
this dissertation and prove UA from Weak Comparison (\cref{WCtoUA}). This proof
is incomplete in the sense that several technical lemmas are deferred until the
end of the chapter. In \cref{UAMOSection}, we give the proof of the linearity of
the Mitchell order from UA (\cref{UAMO}), which is actually a simplification of
the proof in \cref{WCMOSection}. We also prove a sort of converse: UA restricted
to normal ultrafilters is equivalent to the linearity of the Mitchell order.
Finally, in \cref{LemmaProofs}, we prove the technical lemmas we had set aside
in \cref{WCUASection}.
\section{Preliminary definitions}\label{PrelimsSection}
\subsection{Ultrapowers}
We briefly put down our conventions on ultrapowers.\index{Ultrapower} If \(U\)
is an ultrafilter, we denote by \[j_U : V\to M_U\] the ultrapower of the
universe by \(U\). If \(M_U\) is wellfounded, or equivalently if \(U\) is
countably complete, our convention is that \(M_U\) denotes the unique transitive
class isomorphic to the ultrapower of the universe by \(U\). The ultrafilters we
consider will almost always be countably complete.

Many arguments in this dissertation proceed by applying an ultrafilter to a
model to which it does not belong. This involves a taking {\it relativized
ultrapower}.\index{Ultrapower!relativized ultrapower} If \(N\) is a transitive
model of ZFC and \(X\in N\), an {\it \(N\)-ultrafilter on \(X\)} is a set
\(U\subseteq P(X)\cap N\) such that \((N,U)\vDash U\text{ is an ultrafilter}\).
Equivalently, \(U\) is an ultrafilter on the Boolean algebra \(P(X)\cap N\). One
can form the ultrapower of \(N\) by \(U\), denoted \[j_U^N : N\to M_U^N\] using
a modified ultrapower construction that uses only functions that belong to
\(N\). For any function \(f\in N\) that is defined \(U\)-almost everywhere, we
denote by \([f]_U^N\) the point in \(M_U^N\) represented by \(f\). Since the
point \([\text{id}]_U^N\) comes up so often, we introduce special notation for
it:
\begin{defn}\label{SeedNotation}
If \(U\) is an \(N\)-ultrafilter, \(\id_U^N\) denotes the point
\([\text{id}]_U^N\).\index{\(\id_U\) (seed of an ultrapower embedding)}
\end{defn}
We will drop the superscript \(N\) when it is convenient and unambiguous.

Derived ultrafilters\index{Ultrafilter!derived}\index{Derived ultrafilter} allow
us to extract combinatorial content from elementary embeddings:
\begin{defn}
Suppose \(N\) and \(M\) are transitive models of ZFC and \(j : N\to M\) is an
elementary embedding. Suppose \(X\in N\) and \(a\in j(X)\). The {\it
\(N\)-ultrafilter on \(X\) derived from \(j\) using \(a\)} is the
\(N\)-ultrafilter \(\{A\in P(X)\cap N : a\in j(A)\}\).
\end{defn}

What is the relationship between an elementary embedding and the ultrapowers by
its derived ultrafilters? The answer is contained in the following lemma:
\begin{lma}\label{FactorEmbedding}
	Suppose \(N\) and \(M\) are transitive models of ZFC and \(j : N\to M\) is
	an elementary embedding. Suppose \(X\in N\) and \(a\in j(X)\). Let \(U\) be
	the \(N\)-ultrafilter on \(X\) derived from \(j\) using \(a\). Then there is
	a unique embedding \(k : M_U^N\to M\) such that \(k\circ j_U^N = j\) and
	\(k(\id_U) = a\).\qed
\end{lma}
We refer to the embedding \(k\) as the {\it factor embedding}\index{Factor
embedding} associated to the derived ultrafilter \(U\).

Often we will wish to discuss an ultrapower embedding without the need to choose
any particular ultrafilter giving rise to it, so we introduce the following
terminology:

\begin{defn} If \(N\) and \(M\) are transitive models of ZFC, an elementary
embedding \(j : N\to M\) is an {\it ultrapower
embedding}\index{Ultrapower!ultrapower embedding} if there is an
\(N\)-ultrafilter \(U\) such that \(M = M_U^N\) and \(j = j_U^N\). 
\end{defn}

\begin{defn}
If \(N\) is a transitive model of ZFC, a {\it countably complete ultrafilter of
\(N\)}\index{Ultrafilter!of an inner model} is a point \(U\in N\) such that
\(N\) satisfies that \(U\) is a countably complete ultrafilter.
\end{defn}
An \(N\)-ultrafilter \(U\) is a countably complete ultrafilter of \(N\) if and
only if its corresponding ultrapower \(j: N\to M\) is wellfounded and definable
over \(N\).
\begin{defn}
An ultrapower embedding \(j : N \to M\) is an {\it internal ultrapower
embedding}\index{Ultrapower!internal ultrapower embedding} of \(N\) if there is
a countably complete ultrafilter \(U\) of \(N\) such that \(j =
j_U^N\).\end{defn}

An important point is that for our purposes, when we speak of ultrapower
embeddings, we only mean ultrapower embeddings between wellfounded models. For
example, if \(U\) is an ultrafilter on \(\omega\), then the embedding \(j_U :
V\to M_U\) does not count as an ultrapower embedding. 

There is a characterization of ultrapower embeddings that does not refer to
ultrafilters at all.

\begin{defn}
Suppose \(N\) and \(M\) are transitive set models of ZFC. An elementary
embedding \(j : N\to M\) is {\it cofinal}\index{Cofinal embedding} if for all
\(a\in M\), there is some \(X\in N\) such that \(a\in j(X)\).
\end{defn}
Equivalently, \(j\) is cofinal if \(j[\text{Ord}\cap N]\) is cofinal in
\(\text{Ord}\cap M\).
\begin{defn}
Suppose \(N\) and \(M\) are transitive set models of ZFC. An elementary
embedding \(j : N\to M\) is a {\it weak ultrapower
embedding}\index{Ultrapower!ultrapower embedding!weak ultrapower embedding} if
there is some \(a\in M\) such that every element of \(M\) is definable in \(M\)
from parameters in \(j[N]\cup \{a\}\).
\end{defn}
For metamathematical reasons (namely, the undefinability of definability), we
cannot define the concept of a weak ultrapower embeddings when \(M\) is a proper
class.

\begin{lma}
Suppose \(N\) and \(M\) are transitive set models of \textnormal{ZFC}. An
elementary embedding \(j: N\to M\) is an ultrapower embedding if and only if
\(j\) is a cofinal weak ultrapower embedding.\qed
\end{lma}

The following notation will be extremely important in our analysis of elementary
embeddings:
\begin{defn}
Suppose \(N\) and \(M\) are transitive models of \textnormal{ZFC}. Suppose \(j :
N\to M\) is a cofinal elementary embedding and \(S\) is a subclass of \(M\).
Then the {\it hull of \(S\) in \(M\) over \(j[N]\)}\index{\(H^M\) (hull inside
\(M\))} is the class \(H^M(j[N]\cup S) = \{j(f)(x_1,\dots,x_n) :
x_1,\dots,x_n\in S\}\).
\end{defn}

The fundamental theorem about these hulls, which we use repeatedly and
implicitly, is the following:

\begin{lma}
	Suppose \(N\) and \(M\) are transitive models of \textnormal{ZFC}. Suppose
	\(j : N\to M\) is a cofinal elementary embedding and \(S\) is a subclass of
	\(M\). Then the hull of \(S\) in \(M\) over \(j[N]\) is the minimum
	elementary substructure of \(M\) containing \(j[N]\cup S\).\qed
\end{lma}

For more on hulls, see \cite{Larson} Chapter 1 Lemma 1.1.18. (Larson's lemma
should use a stronger theory than ZFC \(-\) Replacement;
\(\Sigma_2\)-Replacement suffices.) Using hulls, we can give a
metamathematically unproblematic model theoretic characterization of ultrapower
embeddings between transitive models that are not assumed to be sets:
\begin{lma}\label{UltraChar}
	Suppose \(N\) and \(M\) are transitive models of \textnormal{ZFC}. A cofinal
	elementary embedding \(j : N\to M\) is an ultrapower embedding if \(M =
	H^M(j[N]\cup \{a\})\) for some \(a\in M\).\qed
\end{lma}
\subsection{Close embeddings}
The property of being an internal ultrapower embedding is a very stringent
requirement. {\it Closeness} is a natural weakening that originated in inner
model theory:
\begin{defn}
Suppose \(N\) and \(M\) are transitive models of ZFC and \(j: N\to M\) is an
elementary embedding. Then \(j\) is {\it close}\index{Close!Embedding} to \(N\)
if \(j\) is cofinal and for all \(X\in N\) and \(a\in j(X)\), the
\(N\)-ultrafilter on \(X\) derived from \(j\) using \(a\) belongs to \(N\).
\end{defn}

Close embeddings have a very natural model theoretic definition that makes no
reference to ultrafilters:
\begin{lma}\label{CloseLemma0}
Suppose \(N\) and \(P\) are transitive models of ZFC and \(j: N\to P\) is an
elementary embedding. Then the following are equivalent:
\begin{enumerate}[(1)]
	\item \(j\) is close to \(N\).
	\item For any \(a\in P\), \(j\) factors as \(N\stackrel{i}\longrightarrow M
	\stackrel{k}\longrightarrow P\) where \(i : N\to M\) is an internal
	ultrapower embedding, \(k : M\to P\) is an elementary embedding, and \(a\in
	k[M]\).
	\item For any set \(A\in P\), the inverse image \(j^{-1}[A]\) belongs to
	\(N\).
\end{enumerate}
\begin{proof}
	{\it (1) implies (2):} Immediate from the factor embedding construction
	\cref{FactorEmbedding}.
	
	{\it (2) implies (3):} Fix \(A\in P\), and we will show \(j^{-1}[A]\in N\).
	Factor \(j\) as \(N\stackrel{i}\longrightarrow M \stackrel{k}\longrightarrow
	P\) where \(i : N\to M\) is an internal ultrapower embedding, \(k : M\to N\)
	is an elementary embedding, and \(A\in k[M]\). Fix \(B\in M\) such that
	\(k(B) = A\). Now \(i^{-1}[B] \in N\) since \(i\) is an internal ultrapower
	embedding of \(N\). We finish by showing \(i^{-1}[B] = j^{-1}[A]\). First,
	by the elementarity of \(k\), \(B = k^{-1}[A]\). Therefore \(i^{-1}[B] =
	i^{-1}[k^{-1}[A]] = (k\circ i)^{-1}[A] = j^{-1}[A]\).
	
	{\it (3) implies (1):} 	We first show that \(j\) is cofinal. Assume not,
	towards a contradiction. Then there is an ordinal \(\alpha\in P\) that lies
	above all ordinals in the range of \(j\). Therefore \(j^{-1}[\alpha] =
	\text{Ord}\cap N\notin N\), which is a contradiction.
	
	Finally, fix \(X\in N\) and \(a\in P\) with \(a\in j(X)\). We must show that
	the \(N\)-ultrafilter on \(X\) derived from \(j\) using \(a\) belongs to
	\(N\). Let \(\pr a {j(X)}\) denote the principal \(N\)-ultrafilter on
	\(j(X)\) concentrated at \(a\). Then the \(N\)-ultrafilter on \(X\) derived
	from \(j\) using \(a\) is precisely \(j^{-1}[\pr a {j(X)}]\), which belongs
	to \(N\) by assumption.
\end{proof}
\end{lma}

Most texts on inner model theory define {\it close
extenders}\index{Close!Extender} rather than close embeddings, so we briefly
describe the relationship between these two concepts. If \(N\) is a transitive
model of ZFC and \(E\) is an \(N\)-extender of length \(\lambda\), then \(E\) is
{\it close to \(N\)} if \(E_a\in M\) for all \(a\in [\lambda]^{<\omega}\). 
\begin{lma}An \(N\)-extender \(E\) is close to \(N\) if and only if the elementary embedding \(j_E^N\) is close to \(N\).\qed\end{lma} 
The fact that the comparison process gives rise to close embeddings is less
well-known than the fact that all extenders applied in a normal iteration tree
are close, which for example is proved in \cite{MitchellSteel}. Given that each
of the individual extenders that are applied are close, the following fact shows
that all the embeddings between models of ZFC in a normal iteration tree are
close:
\begin{lma}\label{CloseCompositions}
\begin{enumerate}[(1)]
\item If \(N\stackrel{i}\longrightarrow M \stackrel{k}\longrightarrow P\) are
close embeddings, then the composition \(k\circ i\) is close to \(N\).
\item Suppose \(\mathcal D = \{M_p,j_{pq} : p\leq q\in I\}\) is a directed
system of transitive models of \textnormal{ZFC} and elementary embeddings.
Suppose \(p\in I\) is an index such that for all \(q\geq p\) in \(I\), \(j_{pq}:
M_p\to M_q\) is close to \(M_p\). Let \(N\) be the direct limit of \(\mathcal
D\), and assume \(N\) is transitive. Then the direct limit embedding
\(j_{p\infty} : M_p\to N\)  is close to \(M_p\).
\end{enumerate}
\begin{proof}
	(1) is immediate from \cref{CloseLemma0} (3). (2) is clear from
	\cref{CloseLemma0} (2).
\end{proof}
\end{lma}

An often useful trivial fact about close embeddings is that their right-factors
are close:
\begin{lma}\label{CloseLemma}
	If \(j : N \to P\) is a close embedding and \(j = k\circ i\) where
	\(N\stackrel i\longrightarrow M\stackrel k\longrightarrow P\) are elementary
	embeddings. Then \(i\) is close to \(N\).\qed
\end{lma}

Another fact which is almost tautological is that an ultrapower embedding is
internal if and only if it is close:
\begin{lma}\label{CloseUltrapower}
	If \(j  :N\to M\) is an ultrapower embedding, then \(j\) is an internal
	ultrapower embedding of \(N\) if and only if \(j\) is close to \(N\).\qed
\end{lma}
\subsection{Uniform ultrafilters}
One of the most basic notions from ultrafilter theory is that of a uniform
ultrafilter:
\begin{defn}
An ultrafilter \(U\) is {\it uniform}\index{Ultrafilter!uniform} if every set in
\(U\) has the same cardinality. If \(U\) is an ultrafilter, the {\it
size}\index{\(\lambda_U\) (size of an ultrafilter)}\index{Size of an ultrafilter
(\(\lambda_U\))} of \(U\), denoted \(\lambda_U\), is the least cardinality of a
set in \(U\).
\end{defn}
The cardinals \(\lambda_U\) for \(U\) a countably complete ultrafilter will
become very important in \cref{SCChapter1}.

Equivalently, \(U\) is a uniform ultrafilter on \(X\) if it extends the
Fr\'echet filter on \(X\), the collection of \(A\subseteq X\) such that
\(|X\setminus A| < |X|\).\index{Fr\'echet filter} It will be important to
distinguish between the notion of a uniform ultrafilter and the similar but
distinct notion of a {\it tail uniform} ultrafilter on an ordinal, defined in
\cref{TailUniformSection}. These notions are often confused in the literature.

\begin{defn}
Suppose \(U\) and \(W\) are ultrafilters. Then \(U\) and \(W\) are {\it
isomorphic},\index{Isomorphism of ultrafilters} denoted \(U\cong W\), if there
exist \(X\in U\), \(Y\in W\), and a bijection \(f : X\to Y\) such that for all
\(A\subseteq X\), \(A\in U\) if and only if \(f[A]\in W\).
\end{defn}

Ultrafilter isomorphism is equivalent to the following model theoretic property:

\begin{defn}\label{ModelIsoDef}
Suppose \(j_0 : N\to M_0\) and \(j_1 : N\to M_1\) are elementary embeddings. We
write \((M_0,j_0)\cong (M_1,j_1)\) to denote that there is an isomorphism \(k :
M_0\to M_1\) such that \(k\circ j_0 = j_1\).
\end{defn}

The following lemma (due to Rudin-Keisler) is explained in \cref{RKSection}:

\begin{lma}\label{IsoChar}
If \(U\) and \(W\) are ultrafilters, then \(U\) and \(W\) are isomorphic if and
only if \((M_U,j_U)\cong (M_W,j_W)\).\qed
\end{lma}
For countably complete ultrafilters, there is a much simpler model theoretic
characterization of ultrafilter isomorphism (so we will not really need the
notation from \cref{ModelIsoDef}):
\begin{cor}
If \(U\) and \(W\) are countably complete ultrafilters, then \(U\) and \(W\) are
isomorphic if and only if \(j_U = j_W\).
\begin{proof}
Since \(M_U\) and \(M_W\) are transitive, the only possible isomorphism between
\(M_U\) and \(M_W\) is the identity. Hence \((M_U,j_U)\cong (M_W,j_W)\) if and
only if \(j_U= j_W\). 
\end{proof}
\end{cor}

Notice that if \(U\cong W\) then \(\lambda_U = \lambda_W\). Since we are mostly
interested in ultrapower embeddings and not ultrafilters themselves, the
following lemma lets us focus our attention on uniform ultrafilters that lie on
cardinals:

\begin{lma}\label{UniformIso}
Any ultrafilter \(U\) is isomorphic to a uniform ultrafilter \(W\) on
\(\lambda_U\).
\begin{proof}
Fix \(X\in U\) such that \(|X| = \lambda_U\). Let \(f : X\to \lambda_U\) be a
bijection. Let \(W = \{A\subseteq \lambda_U : f^{-1}[A]\in U\}\). Then \(U\cong
W\). Moreover \(W\) is uniform since \(\lambda_W = \lambda_U\).
\end{proof}
\end{lma}

Let us also mention a basic generalization of uniformity to the relativized
case:\index{\(\lambda_U\) (size of an ultrafilter)!relativized}
\begin{defn}\label{RelativizedLambda}
	Suppose \(M\) is a transitive model of ZFC and \(U\) is an
	\(M\)-ultrafilter. Then the {\it size } of \(U\) is the \(M\)-cardinal
	\(\lambda_U = \min\{|X|^M  : X\in U\}\).
\end{defn}
\subsection{Normal ultrafilters and the Mitchell order}
\begin{defn}
Suppose \(\langle X_\alpha : \alpha < \delta\rangle\) is a sequence of subsets
of \(\delta\). The {\it diagonal intersection}\index{Diagonal
intersection}\index{\(\triangle_{\alpha < \delta}\) (diagonal intersection)} of
\(\langle X_\alpha : \alpha < \delta\rangle\) is the set \[\triangle_{\alpha <
\delta} X_\alpha = \{\alpha < \delta : \alpha\in \textstyle\bigcap_{\beta <
\alpha} X_\beta\}\]
\end{defn}

\begin{defn}
A uniform ultrafilter on an infinite cardinal \(\kappa\) is {\it
normal}\index{Normal ultrafilter}\index{Normal ultrafilter|seealso{Normal fine
ultrafilter}}\index{Ultrafilter!normal} if it is closed under diagonal
intersections.
\end{defn}
\begin{lma}\label{NormalChar}
Suppose \(U\) is an ultrafilter on an ordinal \(\kappa\). The following are
equivalent:
\begin{enumerate}[(1)]
\item \(U\) is normal.
\item \(U\) is \(\kappa\)-complete and \(a_U = \kappa\).\qed
\end{enumerate}
\end{lma}

The Mitchell order was introduced by Mitchell in \cite{Mitchell}.
\begin{defn}
Suppose \(U\) and \(W\) are normal ultrafilters. The {\it Mitchell
order}\index{Mitchell order|seealso{Generalized Mitchell order}} is defined by
setting \(U\mo W\) if \(U\in M_W\).
\end{defn}
This definition makes sense by our convention that the ultrapower of the
universe by a countably complete ultrafilter is taken to be transitive.
\begin{lma}
The Mitchell order is a wellfounded partial order.\qed
\end{lma}
Actually, the interested reader will find several generalizations of this fact
scattered throughout this dissertation. For example, \cref{KOWellfounded} and
\cref{MOWF} come to mind.
\begin{defn}
If \(U\) is a normal ultrafilter on a cardinal \(\kappa\), then
\(o(U)\)\index{Mitchell order!rank (\(o(U)\))}\index{\(o(U)\), \(o(\kappa)\)
(rank in the Mitchell order)} denotes the rank of \(U\) in the restriction of
the Mitchell order to normal ultrafilters on \(\kappa\).  For any ordinal
\(\kappa\), \(o(\kappa)\) denotes the rank of the restriction of the Mitchell
order to normal ultrafilters on \(\kappa\).
\end{defn}

\section{The linearity of the Mitchell order}\label{LinMOSection}
Our original proof of the linearity of the Mitchell order did not use the
Ultrapower Axiom as a hypothesis. Instead, it used a principle called {\it Weak
Comparison} that was introduced by Woodin \cite{Woodin} in his work on the Inner
Model Problem for supercompact cardinals. 

Weak Comparison is directly motivated by the Comparison Lemma of inner model
theory, and it is immediately clear that Weak Comparison holds in all known
canonical inner models. On the other hand, although the Ultrapower Axiom is a
more elegant principle than Weak Comparison, the fact that the Ultrapower Axiom
holds in all known canonical inner models is not nearly as obvious. But our
proof of the linearity of the Mitchell order from Weak Comparison actually shows
that the Ultrapower Axiom follows from Weak Comparison, and this is how the
Ultrapower Axiom was originally isolated.

In this section, we will introduce Weak Comparison and then prove that Weak
Comparison implies the linearity of the Mitchell order. We then isolate the
Ultrapower Axiom by remarking that this proof breaks naturally into two
implications: first, that Weak Comparison implies the Ultrapower Axiom, and
second, that the Ultrapower Axiom implies the linearity of the Mitchell order.
We hope that this ``genetic approach" will help motivate the formulation of the
Ultrapower Axiom. The reader who does not want to learn about Weak Comparison
can skip ahead to \cref{UAMOSection}. We emphasize, however, that the fact that
Weak Comparison implies UA is central to the motivation for this work.
\subsection{Weak Comparison}\label{WCSection}
Stating Weak Comparison requires a number of definitions. The following
notational convention will make many of our arguments easier to read:
\begin{defn}
Suppose \(N_0,N_1,P\) are transitive models of ZFC. We write \[(j_0,j_1) :
(N_0,N_1)\to P\] to mean that \(j_0: N_0\to P\) and \(j_1 : N_1\to P\) are
elementary embeddings.
\end{defn}

Weak Comparison is a comparison principle for a class of structures. The next
two definitions specify this class.
\begin{defn}
Suppose \(M\) is a model of ZFC. Then \(M\) is {\it finitely
generated}\index{Finitely generated model} if there is some \(a\in M\) such that
every point in \(M\) is definable in \(M\) using \(a\) as a parameter.
\end{defn}

\begin{defn}
Suppose \(M\) is a transitive set that satisfies ZFC. Then \(M\) is a {\it
\(\Sigma_2\)-hull} if there is a \(\Sigma_2\)-elementary embedding \(\pi : M \to
V\).
\end{defn}

We can now state Weak Comparison:
\begin{wc}
If \(M_0\) and \(M_1\) are finitely generated \(\Sigma_2\)-hulls such that
\(P(\omega)\cap M_0 = P(\omega)\cap M_1\), then there are close embeddings
\((k_0,k_1) : (M_0,M_1)\to N\).\index{Weak Comparison}
\end{wc}
We conclude this section by sketching Woodin's argument that Weak Comparison
holds in all known canonical inner models. Assume that \(V\) itself is a
canonical inner model, so that there is some sort of Comparison
Lemma\index{Comparison Lemma} for countable sufficiently elementary
substructures of \(V\). Assume \(M_0\) and \(M_1\) are finitely generated
\(\Sigma_2\)-hulls. We will show that there are close embeddings \((k_0,k_1) :
(M_0,M_1)\to N\).

The fact that \(M_0\) and \(M_1\) are countable \(\Sigma_2\)-hulls implies that
the Comparison Lemma applies to them. The comparison process therefore produces
transitive structures \(N_0\) and \(N_1\) such that one of the following holds:
\begin{case}\label{GoodComp}\(N_0 = N_1\) and there are close embeddings \((k_0,k_1) : (M_0,M_1)\to N_0\).\end{case}
\begin{case}\label{BadComp1}\(N_0\in N_1\), \(P(\omega)\cap N_1\subseteq M_1\), and there is a close embedding \(k_0 :M_0\to N_0\). \end{case}
\begin{case}\label{BadComp2}  \(N_1\in N_0\), \(P(\omega)\cap N_0\subseteq M_0\), and there is a close embedding \(k_1 :M_1\to N_1\). \end{case}
\cref{GoodComp} is the result of ``coiteration," while in \cref{BadComp1} and
\cref{BadComp2}, one of the models has ``outiterated" the other. To obtain weak
comparison for the pair \(M_0\) and \(M_1\), it suffices to show that
\cref{GoodComp} holds. To do this we show that \cref{BadComp1} and
\cref{BadComp2} cannot hold. 

Assume towards a contradiction that \cref{BadComp1} holds. Since \(M_0\) is
finitely generated, there is some \(a\in M_0\) such that every point in \(M_0\)
is definable in \(M_0\) from \(a\). Therefore \(k_0[M_0]\) is equal to the set
of points in \(N_0\) that are definable in \(N_0\) from \(k_0(a)\). Since
\(N_0\in N_1\), it follows that \(k_0[M_0]\in N_1\) and \(k_0[M_0]\) is
countable in \(N_1\). Therefore its transitive collapse, namely \(M_0\), is in
\(N_1\) and is countable in \(N_1\). Let \(x\in P(\omega)\cap N_1\) code \(M_0\)
in the sense that any transitive model \(H\) of ZFC with \(x\in H\) has \(M_0\in
H\). Then \(x\in P(\omega)\cap N_1\subseteq P(\omega)\cap M_1 = P(\omega)\cap
M_0\). It follows that \(x\in M_0\). Since \(x\) codes \(M_0\), this implies
\(M_0\in M_0\), which is a contradiction.

A similar argument shows that \cref{BadComp2} does not hold. Therefore
\cref{GoodComp} must hold.

This argument actually shows that a slight strengthening of Weak Comparison is
true in all known canonical inner models:
\begin{wc2}[Strong Version] If \(M_0\) and \(M_1\) are finitely generated
\(\Sigma_2\)-hulls, either \(M_0\in M_1\), \(M_1\in M_0\), or there are close
embeddings \((k_0,k_1) : (M_0,M_1)\to N\).
\end{wc2}
The strong version of Weak Comparison implies the Continuum
Hypothesis.\footnote{Here one must assume in addition to the strong version of
Weak Comparison that \(V = \text{HOD}\) and there is a \(\Sigma_2\)-correct
worldly cardinal. In fact, these hypotheses are necessary for all our
consequences of Weak Comparison.} It is not clear if it has any other
advantages.
\subsection{Weak Comparison and the Mitchell order}\label{WCMOSection}
In the interest of full disclosure, we admit that we cannot actually prove the
linearity of the Mitchell order from Weak Comparison. Rather we will need some
auxiliary hypotheses:
\begin{thm}\label{WCtoMO}
Assume that \(V = \textnormal{HOD}\) and there is a \(\Sigma_2\)-correct worldly
cardinal. Assume Weak Comparison holds. Then the Mitchell order is linear.
\end{thm}
The need for these auxiliary hypotheses is one of the quirks of Weak Comparison,
and it is part of the reason we think the Ultrapower Axiom is a more elegant
principle.

Here a cardinal \(\kappa\) is {\it worldly}\index{Worldly cardinal} if
\(V_\kappa\) satisfies ZFC, and {\it \(\Sigma_2\)-correct} if
\(V_\kappa\prec_{\Sigma_2} V\). This is a very weak large cardinal hypothesis.
For example, if \(\kappa\) is inaccessible, then in \(V_\kappa\) there is a
\(\Sigma_2\)-correct worldly cardinal; indeed, Morse-Kelley set theory implies
the existence of a \(\Sigma_2\)-correct worldly cardinal. If \(\kappa\) is a
strong cardinal, then \(\kappa\) itself is a \(\Sigma_2\)-correct worldly
cardinal. The hypothesis is motivated by the following lemma, which we defer to
a later section:
\begin{replma}{Worldliness}
The existence of a \(\Sigma_2\)-hull is equivalent to the existence of a
\(\Sigma_2\)-correct worldly cardinal.
\end{replma}
If one wants to apply Weak Comparison at all, at the very least, one needs the
existence of a \(\Sigma_2\)-hull, and therefore one needs a \(\Sigma_2\)-correct
worldly cardinal. One also needs {\it finitely generated} models, and this is
where we use the principle \(V = \text{HOD}\):
\begin{defn}
Suppose \(M\) is a model of ZFC. Then \(M\) is {\it pointwise
definable}\index{Pointwise definable model} if every point in \(M\) is definable
in \(M\) without parameters.
\end{defn}
\begin{lma}\label{Pointwise} Assume \(V = \textnormal{HOD}\). If there is a
\(\Sigma_2\)-hull, then there is a pointwise definable \(\Sigma_2\)-hull.\qed
\end{lma}
The principle \(V = \textnormal{HOD}\) arguably does not hold in all canonical
inner models. (The standard counterexample is \(L[M_1^\#]\), though one might
instead argue that this is not a canonical inner model.) The proof that the
Mitchell order is linear, however, really does work in any inner model. The fact
that we must assume \(V=\text{HOD}\) is again a quirk of the formulation of Weak
Comparison.

The key technical lemma of \cref{WCtoMO} is the following closure property:
\begin{replma}{BigLemma}
	The set of finitely generated \(\Sigma_2\)-hulls is closed under internal
	ultrapowers.
\end{replma}

We defer the proof to \cref{LemmaProofs}. We now proceed to the proof of
\cref{WCtoMO} granting the lemmas.
\begin{proof}[Proof of \cref{WCtoMO}] Since there is a \(\Sigma_2\)-correct
worldly cardinal and \(V = \textnormal{HOD}\), we can fix a pointwise definable
\(\Sigma_2\)-hull \(H\) (by \cref{Pointwise}). It suffices to show that the
Mitchell order is linear in \(H\), since this is a \(\Pi_2\)-statement and
\(H\equiv_{\Pi_2} V\). 

Suppose \(U_0\) and \(U_1\) are normal ultrafilters of \(H\). We must show that
in \(H\), either \(U_0 = U_1\), \(U_0\mo U_1\), or \(U_0\gmo U_1\). We might as
well assume that \(U_0\) and \(U_1\) are normal ultrafilters on the same
cardinal \(\kappa\), since otherwise it is immediate that \(U_0\mo U_1\) or
\(U_0\gmo U_1\).

Let \(j_0 : H\to M_0\) be the ultrapower of \(H\) by \(U_0\) and let \(j_1 :
H\to M_1\) be the ultrapower of \(H\) by \(U_1\). By the closure of finitely
generated \(\Sigma_2\)-hulls under internal ultrapowers (\cref{BigLemma}),
\(M_0\) and \(M_1\) are finitely generated \(\Sigma_2\)-hulls. Since \(M_0\) and
\(M_1\) are internal ultrapowers of \(H\), \(P(\omega)\cap M_0 = P(\omega)\cap
M_1\). Therefore, by Weak Comparison there are close embeddings \[(k_0,k_1) :
(M_0,M_1)\to N\] Since \(H\) is pointwise definable,
\begin{equation}\label{Commutativity}k_0\circ j_0 = k_1\circ j_1\end{equation}
This is because \(k_0\circ j_0,k_1\circ j_1\) are both elementary embeddings
from \(H\) to \(N\), and therefore must shift all parameter-free definable
points in the same way.

The proof now splits into three cases.
\begin{case}\label{EqualCase}\(k_0(\kappa) = k_1(\kappa)\).\end{case}
\begin{case}\label{LessCase}\(k_0(\kappa) < k_1(\kappa)\).\end{case}
\begin{case}\label{GreatCase}\(k_0(\kappa) > k_1(\kappa)\).\end{case}
In \cref{EqualCase}, we will show \(U_0 = U_1\), in \cref{LessCase}, we will
show \(U_0\mo U_1\), and in \cref{GreatCase}, we will show \(U_0\gmo U_1\). This
will complete the proof.
\begin{proof}[Proof in \cref{EqualCase}] Suppose \(A\subseteq \kappa\) and
\(A\in H\). We have 
\begin{align}
A\in U_0&\iff\kappa\in j_0(A)\nonumber\\
&\iff k_0(\kappa)\in k_0(j_0(A))\nonumber\\
&\iff k_0(\kappa)\in k_1(j_1(A))\label{ComUse0}\\
&\iff k_1(\kappa)\in k_1(j_1(A))\label{CaseHyp0}\\
&\iff \kappa\in j_1(A)\nonumber\\
&\iff A\in U_1\nonumber
\end{align}
To obtain \cref{ComUse0}, we use \cref{Commutativity} above. To obtain
\cref{CaseHyp0}, we use the case hypothesis that \(k_0(\kappa) = k_1(\kappa)\).
It follows that \(U_0 = U_1\).
\end{proof}
\begin{proof}[Proof in \cref{LessCase}] Suppose \(A\subseteq \kappa\) and \(A\in
H\). We have 
\begin{align}
A\in U_0&\iff \kappa\in j_0(A)\nonumber\\
&\iff k_0(\kappa)\in k_0(j_0(A))\nonumber\\
&\iff k_0(\kappa)\in k_1(j_1(A))\label{ComUse1}\\
&\iff k_0(\kappa)\in k_1(j_1(A))\cap k_1(\kappa)\label{CaseHyp1}\\
&\iff k_0(\kappa)\in k_1(j_1(A) \cap \kappa)\nonumber\\
&\iff k_0(\kappa)\in k_1(A)\label{Complete}
\end{align}
To obtain \cref{ComUse1}, we use \cref{Commutativity} above. To obtain
\cref{CaseHyp1}, we use the case hypothesis that \(k_0(\kappa) < k_1(\kappa)\).
To prove \cref{Complete}, we use that \(U_1\) is \(\kappa\)-complete, so
\(\textsc{crt}(j_1) = \kappa\) and hence \(j_1(A)\cap \kappa = A\) for any
\(A\subseteq \kappa\). 

It follows from this calculation that \(U_0\) is the \(M_1\)-ultrafilter on
\(\kappa\) derived from \(k_1\) using \(k_0(\kappa)\). (Here we use that
\(P(\kappa)\cap M_1 = P(\kappa)\cap H\).) Since \(k_1\) is close to \(M_1\), it
follows that \(U_0\in M_1\). Since \(M_1 = M_{U_1}^H\), this means that \(U_0\mo
U_1\) in \(H\).
\end{proof}
\begin{proof}[Proof in \cref{GreatCase}] The proof in this case is just like the
proof in \cref{LessCase} but with \(U_0\) and \(U_1\) swapped.
\end{proof}
This completes the proof of \cref{WCtoMO}.
\end{proof}
\subsection{Weak Comparison and the Ultrapower Axiom}\label{WCUASection}
We now define the Ultrapower Axiom, which arises naturally from the proof of
\cref{WCtoMO}. Notice that the first half of this proof, which justifies our
application of Weak Comparison to the ultrapowers \(M_0\) and \(M_1\), does not
actually require that \(U_0\) and \(U_1\) are normal ultrafilters. Instead, it
simply requires that they are countably complete. 

In order to state UA succinctly, we make the following definitions.

\begin{defn}\label{ComparisonDef}
Suppose \(N,M_0,M_1,P\) are transitive models of ZFC and \(j_0: N\to M_0\),
\(j_1 : N\to M_1\), and \((k_0,k_1) : (M_0,M_1)\to P\) are elementary
embeddings. 
\begin{itemize}
\item \((k_0,k_1)\) is a {\it comparison} of \((j_0,j_1)\) if \(k_0\circ j_0 = k_1\circ j_1\).\index{Comparison}
\item \((k_0,k_1)\) is an {\it internal ultrapower comparison} if \(k_0\) is an internal ultrapower embedding of \(M_0\) and \(k_1\) is an internal ultrapower embedding of \(M_1\).\index{Internal ultrapower comparison}\index{Comparison!internal ultrapower comparison}
\item \((k_0,k_1)\) is a {\it close comparison} if \(k_0\) is close to \(M_0\) and \(k_1\) is close to \(M_1\).\index{Close comparison}
\end{itemize}
\end{defn}

\begin{ua}
Every pair of ultrapower embeddings of the universe of sets has an internal
ultrapower comparison.\index{Ultrapower Axiom}
\end{ua}

On the face of it, the statement that every pair of ultrapowers has a comparison
by internal ultrapowers looks {\it much} stronger than the conclusion of Weak
Comparison, which only supplies close comparisons. But this is an illusion.

\begin{lma}\label{ClosetoUltra}
Suppose \(N,M_0,M_1\) are transitive set models of \textnormal{ZFC} and \(j_0:
N\to M_0\) and \(j_1 : N\to M_1\) are weak ultrapower embeddings. If
\((j_0,j_1)\) has a comparison by close embeddings, then \((j_0,j_1)\) has a
comparison by internal ultrapowers. 
\begin{proof}
Suppose \((k_0,k_1) : (M_0,M_1)\to P\) is a comparison by close embeddings. Let
\(H\prec P\) be defined by \[H = H^P(k_0[M_0]\cup k_1[M_1])\] Let \(Q\) be the
transitive collapse of \(H\) and let \(h : Q\to P\) be the inverse of the
transitive collapse embedding. Let \(i_0 = h^{-1}\circ k_0\) and \(i_1 = h^{-1}
\circ k_1\).

Obviously \((i_0,i_1) : (M_0,M_1)\to Q\) is a comparison of \((j_0,j_1)\) and
\[Q= H^Q(i_0[M_0]\cup i_1[M_1])\] We need to show it is a comparison by internal
ultrapowers, or in other words that \(i_0\) is an internal ultrapower embedding
of \(M_0\) and \(i_1\) is an internal ultrapower embedding of \(M_1\).

We first show that \(i_0\) is an ultrapower embedding of \(M_0\). Since
\(j_1:N\to M_1\) is a weak ultrapower embedding, there is some \(a\in M_1\) such
that every element of \(M_1\) is definable in \(M_1\) from parameters in
\(j_1[N]\cup \{a\}\). It follows easily that \(Q = H^Q(i_0[M_0]\cup
\{i_1(a)\})\). Therefore \(i_0\) is an ultrapower embedding by \cref{UltraChar}.

Next, we show that \(i_0\) is an {\it internal} ultrapower embedding. Since
\(h\circ i_0 = k_0\) and \(k_0\) is close, in fact, \(i_0\) is close to \(M_0\)
(\cref{CloseLemma}). Since \(i_0\) is a close ultrapower embedding of \(M_0\),
in fact, \(i_0\) is an internal ultrapower embedding of \(M_0\)
(\cref{CloseUltrapower}).

A symmetric argument shows that \(i_1\) is an internal ultrapower embedding of
\(M_1\), completing the proof.
\end{proof}
\end{lma}

This yields a strengthening of Weak Comparison:

\begin{thm}\label{WCStrengthening}Assume Weak Comparison and \(V =
\textnormal{HOD}\). Suppose \(M_0\) and \(M_1\) are finitely generated
\(\Sigma_2\)-hulls such that \(P(\omega)\cap M_0 = P(\omega)\cap M_1\). Then
there are internal ultrapower embeddings \((i_0,i_1) : (M_0,M_1)\to Q\).
\begin{proof}
Applying Weak Comparison, fix close embeddings \((k_0,k_1) : (M_0,M_1)\to P\).

Since \(M_0\) and \(M_1\) are \(\Sigma_2\)-hulls, they satisfy any \(\Pi_3\)
sentence true in \(V\). Therefore they both satisfy \(V = \textnormal{HOD}\).
Let \(H_0\prec M_0\) be the set of points that are definable without parameters
in \(M_0\). Let \(H_1\prec M_1\) be the set of points that are definable without
parameters in \(M_1\). Then \(k_0[H_0] = k_1[H_1]\) is the set of points that
are definable without parameters in \(P\). It follows that \(H_0 \cong H_1\).
Let \(N\) be the common transitive collapse of \(H_0\) and \(H_1\), and let
\(j_0 :N\to M_0\) and \(j_1 : N\to M_1\) be the inverses of the transitive
collapse maps. Note that \(j_0\) and \(j_1\) are weak ultrapower embeddings, and
since \(N\) is pointwise definable, \(k_0\circ j_0 = k_1\circ j_1\).

The weak ultrapower embeddings \((j_0,j_1)\) therefore have a comparison by
close embeddings, namely \((k_0,k_1)\). It follows from \cref{ClosetoUltra} that
they have a comparison by internal ultrapower embeddings.
\end{proof}
\end{thm}

\cref{ClosetoUltra} also yields a proof of the Ultrapower Axiom from the same
hypotheses as \cref{WCtoMO}:
\begin{thm}\label{WCtoUA}
Assume that \(V = \textnormal{HOD}\) and there is a \(\Sigma_2\)-correct worldly
cardinal. If Weak Comparison holds, then the Ultrapower Axiom holds.
\begin{proof}
Since there is a \(\Sigma_2\)-correct worldly cardinal and \(V = \text{HOD}\),
we can fix a pointwise definable \(\Sigma_2\)-hull \(H\) (by \cref{Pointwise}).
Since UA is a \(\Pi_2\)-statement and \(H\equiv_{\Pi_2} V\), it suffices to show
that \(H\) satisfies UA.

Suppose \(j_0 : H\to M_0\) and \(j_1 : H\to M_1\) are internal ultrapower
embeddings of \(H\). We must show that \(H\) satisfies that \((j_0,j_1)\) has an
internal ultrapower comparison. 

By the closure of finitely generated \(\Sigma_2\)-hulls under internal
ultrapowers (\cref{BigLemma}), \(M_0\) and \(M_1\) are finitely generated
\(\Sigma_2\)-hulls. Moreover, since \(M_0\) and \(M_1\) are internal ultrapowers
of \(H\), \(P(\omega)\cap M_0 = P(\omega)\cap H = P(\omega)\cap M_1\). Therefore
by \cref{WCStrengthening}, there are internal ultrapower embeddings \((i_0,i_1)
: (M_0,M_1)\to Q\). Moreover since \(H\) is finitely generated, \(i_0\circ j_0 =
i_1\circ j_1\). It follows that \((i_0,i_1)\) is an internal ultrapower
comparison of \((j_0,j_1)\). This is absolute to \(H\), and therefore \(H\)
satisfies that \((j_0,j_1)\) has an internal ultrapower comparison, as desired.
\end{proof}
\end{thm}
\subsection{The Ultrapower Axiom and the Mitchell order}\label{UAMOSection}
In this subsection, we prove the linearity of the Mitchell order from the
Ultrapower Axiom. We include this proof largely for the benefit of the reader
who would prefer to skip over our discussions of Weak Comparison, since the
proof is very similar to that of \cref{WCtoMO}. The reader who has followed
\cref{WCtoMO} will no doubt notice that both the statement and proof of
\cref{UAMO} below are much simpler and more elegant than those of \cref{WCtoMO}.
It is a general pattern that UA is easier to use than Weak Comparison. In fact,
almost every known consequence of Weak Comparison is a consequence of UA.

\begin{thm}[UA]\label{UAMO}
The Mitchell order is linear.\index{Mitchell order!linearity}
\begin{proof}
Suppose \(U_0\) and \(U_1\) are normal ultrafilters. We must show that either
\(U_0 = U_1\), \(U_0\mo U_1\), or \(U_0\gmo U_1\). We may assume without loss of
generality that \(U_0\) and \(U_1\) are normal ultrafilters on the same cardinal
\(\kappa\), since otherwise it is obvious that either  \(U_0\mo U_1\) or
\(U_0\gmo U_1\).

Let \(j_0 : V\to M_0\) be the ultrapower of the universe by \(U_0\). Let \(j_1 :
V\to M_1\) be the ultrapower of the universe by \(U_1\). Applying UA, there is
an internal ultrapower comparison \((i_0,i_1) : (M_0,M_1)\to P\) of
\((j_0,j_1)\).

The proof now breaks into three cases.
\begin{case}\label{Equal}
\(i_0(\kappa) = i_1(\kappa)\).
\end{case}
\begin{case}\label{Less}
\(i_0(\kappa) < i_1(\kappa)\).
\end{case}
\begin{case}\label{Greater}
\(i_0(\kappa) > i_1(\kappa)\).
\end{case}
In \cref{Equal}, we will prove \(U_0 = U_1\). In \cref{Less}, we will prove
\(U_0\mo U_1\). In \cref{Greater}, we will prove \(U_0\gmo U_1\).
\begin{proof}[Proof in \cref{Equal}] Suppose \(A\subseteq \kappa\). Then
\begin{align}
A\in U_0 &\iff \kappa\in j_0(A)\nonumber\\
&\iff i_0(\kappa)\in i_0(j_0(A))\nonumber\\
&\iff i_0(\kappa)\in i_1(j_1(A))\label{Comp}\\
&\iff i_1(\kappa)\in i_1(j_1(A))\label{CaseHyp2}\\
&\iff \kappa\in j_1(A)\nonumber\\
&\iff A\in U_1\nonumber
\end{align}
To obtain \cref{Comp}, we use the fact that \((i_0,i_1)\) is a comparison, and
in particular that \(i_1\circ j_1 = i_0\circ j_0\). To obtain \cref{CaseHyp2},
we use the case hypothesis that \(i_0(\kappa) = i_1(\kappa)\). It follows that
\(U_0 = U_1\).
\end{proof}
\begin{proof}[Proof in \cref{Less}] Suppose \(A\subseteq \kappa\). Then
\begin{align}
A\in U_0 &\iff \kappa\in j_0(A)\nonumber\\
&\iff i_0(\kappa)\in i_0(j_0(A))\nonumber\\
&\iff i_0(\kappa)\in i_1(j_1(A))\nonumber\\
&\iff i_0(\kappa)\in i_1(j_1(A))\cap i_1(\kappa)\label{CaseHyp3}\\
&\iff i_0(\kappa)\in i_1(j_1(A)\cap \kappa)\nonumber\\
&\iff i_0(\kappa)\in i_1(A)\label{Complete2}
\end{align}
To obtain \cref{CaseHyp3}, we use the case hypothesis that \(i_0(\kappa) <
i_1(\kappa)\). To obtain \cref{Complete2}, we use that \(U_1\) is
\(\kappa\)-complete; therefore \(\textsc{crt}(j_1) = \kappa\) so \(j_1(A)\cap
\kappa = A\) for any \(A\subseteq \kappa\).

It follows that \(U_0\) is the \(M_1\)-ultrafilter derived from \(i_1\) using
\(i_0(\kappa)\). (Here we use that \(P(\kappa)\subseteq M_1\).) Since \(i_1\) is
an internal ultrapower embedding of \(M_1\), \(i_1\) is definable over \(M_1\),
and therefore \(U_0\) is definable over \(M_1\) from \(i_1\) and
\(i_0(\kappa)\). It follows that \(U_0\in M_1\). Since \(M_1 = M_{U_1}\), this
means \(U_0\mo U_1\), as desired.
\end{proof}
\begin{proof}[Proof in \cref{Greater}] The proof in this case is identical to
the proof in \cref{Less} but with \(U_0\) and \(U_1\) swapped.
\end{proof}
Thus no matter which of the cases hold, either \(U_0 = U_1\), \(U_0\mo U_1\), or
\(U_0\gmo U_1\). This completes the proof.
\end{proof}
\end{thm}
There is a partial converse to \cref{UAMO} that helps explain the motivation for
the proof of \cref{UAMO}. To state this converse, we first defines a restricted
version of the Ultrapower Axiom for ultrapower embeddings coming from normal
ultrafilters:
\begin{defn}
The {\it Normal Ultrapower Axiom} is the statement that any pair of ultrapower
embeddings of the universe of sets associated to normal ultrafilters have a
comparison by internal ultrapowers.
\end{defn}

\begin{prp}\label{UAMOConverse}
The Normal Ultrapower Axiom is equivalent to the linearity of the Mitchell
order.
\begin{proof}
The proof that the Normal Ultrapower Axiom implies the linearity of the Mitchell
order is immediate from the proof of \cref{UAMO}. 

Conversely, assume the Mitchell order is linear. Suppose \(U_0\) and \(U_1\) are
normal ultrafilters, and let \(j_0 : V\to M_0\) and \(j_1 : V\to M_1\) be their
ultrapowers. We will show \((j_0,j_1)\) has a comparison by internal
ultrapowers. Assume without loss of generality that \(U_0\mo U_1\). Let \(i_0 :
M_0\to P_0\) be the ultrapower of \(M_0\) by \(j_0(U_1)\). Let \(i_1: M_1\to
P_1\) be the ultrapower of \(M_1\) by \(U_0\). Then \(i_0\) and \(i_1\) are
internal ultrapowers of \(M_0\) and \(M_1\) respectively.
Moreover\footnote{Suppose \(M,N,\) and \(P\) are transitive models of ZFC.
Suppose \(j: M\to N\) and \(i : M\to P\) are elementary embeddings. Assume
\(j\restriction x\in M\) for all \(x\in M\). Assume \(i\) is a cofinal
embedding. Then \(i(j) = \bigcup_{X\in M} i(j\restriction X)\). } \(i_0 =
j_0(j_1)\) and \(i_1 = j_0\restriction M_1\), so \[i_0\circ j_0 = j_0(j_1)\circ
j_0 = j_0\circ j_1 = i_1\circ j_1\] It follows that \((i_0,i_1)\) is a
comparison of \((j_0,j_1)\) by internal ultrapowers.
\end{proof}
\end{prp}
The proof of \cref{UAMOConverse} is local in the sense that it shows that the
comparability of two normal ultrafilters in the Mitchell order is equivalent to
their comparability by internal ultrapowers. This is a special feature of the
Mitchell order on normal ultrafilters. For the generalized Mitchell order
(defined in \cref{GMOChapter}), neither implication is provable. Motivated by
this issue, we develop in \cref{InternalSection} a variant of the generalized
Mitchell order called the {\it internal relation}.
\subsection{Technical lemmas related to Weak Comparison}\label{LemmaProofs}
In this section, we prove several lemmas promised in \cref{WCMOSection}.

\begin{lma}\label{GenClosure}
Suppose \(N\) is a finitely generated model of \textnormal{ZFC} and \(U\) is an
\(N\)-ultrafilter. Then \(M_U^N\) is finitely generated.
\begin{proof}
Fix \(b\in N\) such that every element of \(N\) is definable in \(N\) using
\(b\) as a parameter. Obviously every element of \(j_U[N]\) is definable in
\(M_U^N\) using \(j_U(b)\) as a parameter. But \(M_U^N = \{j_U(f)(a_U) : f\in
N\} = \{g(a_U) : g\in j_U[N]\}\). Therefore every element of \(M_U^N\) is
definable using \(j_U(b)\) and \(a_U\) as parameters. 
\end{proof}
\end{lma}

The next lemma, standard in the case of fully elementary embeddings, is the key
to our analysis of \(\Sigma_2\)-hulls:

\begin{lma}\label{Sigma2Factor}
Suppose \(j: N\to M\) is a \(\Sigma_2\)-elementary embedding between transitive
models of \textnormal{ZFC}. Suppose \(X\in N\), and \(a\in j(X)\). Let \(U\) be
the \(N\)-ultrafilter on \(X\) derived from \(j\) using \(a\). Then there is a
unique \(\Sigma_2\)-elementary embedding \(k : M_U^N\to M\) such that \(k\circ
j_U^N = j\) and \(k(a_U) = a\).
\begin{proof}
We begin with a simple remark. Suppose \(\varphi(v_1,\dots,v_n)\) is a
\(\Sigma_2\)-formula and \(f_1,\dots,f_n\) are functions in \(N\) that are
defined \(U\)-almost everywhere. The statement \(S = \{x\in X : N\vDash
\varphi(f_1(x),\dots,f_n(x))\}\) can be written as a Boolean combination of
\(\Sigma_2\) formulas in the variables \(S\) and \(f_1,\dots,f_n\). It follows
that \[j(\{x\in X : N\vDash \varphi(f_1(x),\dots,f_n(x))\}) = \{x\in j(X) :
M\vDash \varphi(j(f_1)(x),\dots,j(f_n)(x))\}\]

For any function \(f\in N\) defined \(U\)-almost everywhere, set \[k([f]_U) =
j(f)(a)\] Fix a \(\Sigma_2\)-formula \(\varphi(x_1,\dots,x_n)\). The following
calculation shows that \(k\) is a well-defined \(\Sigma_2\)-elementary embedding
from \(M^N_U\) to \(M\): 
\begin{align*}M_U^N\vDash \varphi([f_1]_U,\dots,[f_n]_U)&\iff \{x\in X : N\vDash \varphi(f_1(x),\dots,f_n(x))\}\in U\\
&\iff M\vDash a\in j(\{x\in X : N\vDash \varphi(f_1(x),\dots,f_n(x))\})\\
&\iff M\vDash a\in \{x\in j(X) : M\vDash \varphi(j(f_1)(x),\dots,j(f_n)(x))\}\\
&\iff M\vDash \varphi(j(f_1)(a),\dots,j(f_n)(a))\\
&\iff M\vDash \varphi(k([f_1]_U),\dots,k([f_n]_U))\qedhere
\end{align*}
\end{proof}
\end{lma}

\cref{Sigma2Factor} yields a \(\Sigma_2\)-elementary generalization of the
standard Realizability Lemma:

\begin{lma}\label{HullClosure}
	Suppose \(N\) is a countable \(\Sigma_2\)-hull and \(N\vDash U\) is a
	countably complete ultrafilter. Then \(M_U^N\) is a \(\Sigma_2\)-hull.
	\begin{proof}
		Let \(\pi : N\to V\) be a \(\Sigma_2\)-elementary embedding. Let \(U' =
		\pi(U)\), so \(U'\) is a countably complete ultrafilter. Since
		\(\pi[U]\subseteq U'\) is countable, there is some \(a\in \bigcap
		\pi[U]\). Note that \(U = \pi^-\{a\}\). Therefore by
		\cref{Sigma2Factor}, there is a \(\Sigma_2\)-elementary embedding \(k :
		M_U^N\to V\), so \(M_U^N\) is a \(\Sigma_2\)-hull.
	\end{proof}
\end{lma}

\begin{lma}\label{BigLemma}
	The set of finitely generated \(\Sigma_2\)-hulls is closed under internal
	ultrapowers.
	\begin{proof}
		Immediate from the conjunction of \cref{GenClosure} and
		\cref{HullClosure}.
	\end{proof}
\end{lma}

%Now suppose \(j : N\to M\) is a \(\Sigma_2\)-elementary embedding. Let
%\(H\subseteq M\) be the collection of \(a\in M\) such that there is some \(X\in
%N\) with \(a\in j(X)\). For each \(a\in H\), letting \(U = j^-\{a\}\) and \(k :
%M_U^N\to M\) be the factor embedding, define \[H^M(j[N]\cup \{a\}) =
%\text{ran}(k)\] so \(H^M(j[N]\cup \{a\})\) is a \(\Sigma_2\)-elementary
%substructure of \(M\), \(j : N\to H^M(j[N]\cup \{a\})\) is fully elementary,
%and \(H^M(j[N]\cup \{a\}) = \{j(f)(a) : f\in N_X\}\) for any \(X\) with \(a\in
%j(X)\). If \(a_1,\dots,a_n\in H\), let \[H^M(j[N]\cup \{a_1,\dots,a_n\}) =
%H^M(j[N]\cup \{\{a_1,\dots,a_n\}\})\] For any \(S\subseteq H\), let
%\[H^M(j[N]\cup S) = \bigcup_{a_1,\dots,a_n\in S} H^M(j[N]\cup
%\{a_1,\dots,a_n\})\] The structures \(H^M(j[N]\cup\{a_1,\dots,a_n\})\), for
%\(a_1,\dots,a_n\in S\), form a directed system of \(\Sigma_2\)-elementary
%substructures of \(M\), so their union \(H^M(j[N]\cup S)\) is a
%\(\Sigma_2\)-elementary substructure of \(M\) as well. Similarly, \(j : N\to
%H^M(j[N]\cup S)\) is fully elementary. For example, this construction yields:
\cref{Sigma2Factor} can also be used to prove the following fact:
\begin{lma}\label{Sigma2Factor2}
Suppose \(N\) is a set model of \textnormal{ZFC} and \(j : N\to M\) is a
\(\Sigma_2\)-elementary embedding. Then \(j\) factors as a cofinal elementary
embedding followed by a \(\Sigma_2\)-elementary end extension.\qed
%\begin{proof} Let \(H\subseteq N\) be the collection of \(a\in N\) such that
%there is some \(X\in M\) with \(a\in j(X)\). It is easy to see that \(M\) is an
%end-extension of \(N\). Moreover \(H = H^M(j[N]\cup H)\), so \(j : N\to H\) is
%fully elementary and \(H\) is a \(\Sigma_2\)-elementary substructure of \(M\).
%This proves the lemma. \end{proof}
\end{lma}

\begin{prp}\label{Worldliness}
There is a \(\Sigma_2\)-hull if and only if there is a \(\Sigma_2\)-correct
worldly cardinal. 
\end{prp}
\begin{proof}
Suppose \(N\) is a \(\Sigma_2\)-hull. Let \(\pi : N\to V\) be a
\(\Sigma_2\)-elementary embedding. By \cref{Sigma2Factor}, \(\pi\) factors as a
cofinal elementary embedding \(\pi : N\to H\) followed by a
\(\Sigma_2\)-elementary end extension \(H\prec_{\Sigma_2} V\). Since
\(H\prec_{\Sigma_2} V\), \(H = V_\kappa\) for some cardinal \(\kappa\). Since
\(\pi : N\to V_\kappa\) is fully elementary, \(V_\kappa\) satisfies ZFC. Thus
\(\kappa\) is a \(\Sigma_2\)-correct worldly cardinal.
\end{proof}

\chapter{The Ketonen Order}\label{KetonenChapter}
\section{Introduction}
\subsection{Ketonen's order}
Central to \cref{MOChapter} was an argument that the Mitchell order is linear in
all known canonical inner models. In \cref{WCUASection}, we delved deeper into
the first half of this proof, extracting from it a general inner model principle
called the Ultrapower Axiom. It turns out that a closer look at the second half
of the proof also yields more information: it shows that the Ultrapower Axiom
implies not only the linearity of the Mitchell order, but also the linearity of
a much more general order on countably complete ultrafilters. 

This order dates back to the early 1970s. A remarkable theorem of Ketonen
\cite{Ketonen} from this period states that if every regular cardinal \(\lambda
\geq \kappa\) carries a \(\kappa\)-complete uniform ultrafilter, then \(\kappa\)
is strongly compact. Ketonen\index{Ketonen} gave two proofs of this theorem. The
first is an induction. The second is not as well-known, but is of much greater
interest here. Ketonen introduced a wellfounded order on countably complete
weakly normal ultrafilters, and showed that certain minimal elements in this
order witness the strong compactness of \(\kappa\). (We give this proof in
\cref{KetonenThm} since generalizations of the proof form a key component of our
analysis of strong compactness and supercompactness under UA.)

Independently of Ketonen's work, and a rather long time after, we extracted from
the proof of the linearity of the Mitchell order under UA (\cref{UAMO}) a more
general version of Ketonen's order, which we now call the {\it Ketonen order}.
The Ketonen order is a wellfounded partial order on countably complete
ultrafilters concentrating on ordinals. The key realization, which distinguishes
our work from Ketonen's, is that the Ketonen order can be {\it linear.} In fact,
the totality of the Ketonen order is an immediate consequence of UA
(\cref{Totality}). In fact, the linearity of the Ketonen order is equivalent to
the Ultrapower Axiom. This equivalence is \cref{LinKetThm}, which is probably
the hardest theorem of this chapter. The Ketonen order will be our main tools in
the investigation of the structure of countably complete ultrafilters under UA. 

\subsection{Outline of \cref{KetonenChapter}}
Let us outline the rest of \cref{KetonenChapter}.\\

\noindent {\sc\cref{PrelimsSection2}.} We introduce some more preliminary
definitions that will be used throughout the rest of this dissertation.
Especially important are limits of ultrafilters, which we introduce both in the
traditional ultrafilter theoretic sense and in a generalized setting in terms of
inverse images.\\

\noindent {\sc\cref{KetonenOrderSection}.} We introduce the main object of study of this
chapter, a fundamental tool in the theory of the Ultrapower Axiom: a wellfounded
partial order on countably complete ultrafilters called the Ketonen order. In
\cref{KOCharsSection}, we define the Ketonen order and give various alternate
characterizations. The most important characterization is given by
\cref{KOChars}, which shows that the Ketonen order can be reformulated in terms
of comparisons. This immediately leads to the observation that the Ketonen order
is linear under the Ultrapower Axiom. In \cref{BasicSection}, we establish the
basic order-theoretic properties of the Ketonen order: it is a preorder on the
class of countably complete ultrafilters concentrating on ordinals. Restricted
to tail uniform ultrafilters, it is a partial order. \cref{SpaceLemma} shows
that the Ketonen order is {\it graded} in the sense that if \(\alpha < \beta\),
then the tail uniform ultrafilters on \(\alpha\) all lie below those on
\(\beta\). In particular, the Ketonen order is setlike. We finally prove the
wellfoundedness of the Ketonen order (\cref{KOWellfounded}). The general proof
of the wellfoundedness of the Ketonen order is due to the author.\\

\noindent {\sc\cref{OtherOrderSection}.} We explore the relationship between the Ketonen
order and two well-known orders on ultrafilters. \cref{MOKetSection} shows that
the restriction of the Ketonen order to normal ultrafilters is precisely the
Mitchell order. In this sense the Ketonen order is a generalized Mitchell order.
In \cref{RKSection} we turn to perhaps the best-known order on ultrafilters: the
Rudin-Keisler order. We take this opportunity to set down some basic facts about
this order, sometimes with proofs. The Ketonen order is not isomorphism
invariant, so it cannot extend the Rudin-Keisler order. To explain these orders'
relationship better, we define an auxiliary order called the {\it revised
Rudin-Keisler order} which is contained in the intersection of the Rudin-Keisler
order and the Ketonen order. Moreover we introduce the concept of an {\it
incompressible ultrafilter}, an ultrafilter \(U\) whose generator \(\id_U\) is
as small as possible (see \cref{IncomPower}). An argument due to Solovay shows
that the strict Rudin-Keisler order and the revised Rudin-Keisler order coincide
on incompressible ultrafilters. Thus the Ketonen order extends the strict
Rudin-Keisler order on countably complete incompressible ultrafilters.\\

\noindent {\sc\cref{VariantSection}.} We study several variants of the Ketonen order. In
\cref{ExtendedKetonen}, we investigate the relationship between the Ketonen
order and notions from inner model theory. We introduce a model theoretic
generalization of the Ketonen order whose domain is the class of pointed models
of ZFC, structures \((M,\xi)\) where \(M\) is a transitive model of ZFC and
\(\xi\) is an ordinal of \(M\). This defines a  coarse analog of the Dodd-Jensen
order, the canonical prewellorder on mice. We give a generalized wellfoundedness
proof for this order (\cref{GenWellfounded}) that is closely related to the
proof of the wellfoundedness of the Dodd-Jensen order. We use this to prove an
often useful lemma that is a coarse analog of the Dodd-Jensen lemma:
\cref{MinDefEmb} shows that definable embeddings are pointwise minimal on the
ordinals.

In the next subsection, \cref{SOSection}, we introduce the {\it seed order}. In
early versions of this work, we mainly used the seed order where we now use the
Ketonen order. There is no substantive difference between these approaches since
under UA the two orders coincide. In ZFC, however, one cannot prove that the
seed order is transitive: indeed, we show by a silly argument that the
transitivity of the seed order implies the Ultrapower Axiom. We also introduce
an extension of the seed order to pointed ultrapowers. The next subsection
\cref{MInftySection} is spent relating this order to the structure of the direct
limit of all ultrapowers under UA.

The next two subsections are devoted to combinatorial generalizations of the
Ketonen order. One does not need to read them to understand the rest of this
dissertation. In  \cref{LipschitzSection}, we introduce a generalized version of
the Lipschitz order, and show that this order extends the Ketonen order on
countably complete ultrafilters. Therefore under UA, the two orders coincide,
which gives a strange analog of the linearity of the Lipschitz order in
determinacy theory. \cref{FilterSection} introduces a combinatorial
generalization of the Ketonen order to filters, which demonstrates a
relationship between the Ketonen order and the canonical order on stationary
sets due to Jech \cite{JechStationary}.\\

\noindent {\sc\cref{LinKetSection}.} This section contains \cref{LinKetThm}, the most
substantive result of the chapter:  the linearity of the Ketonen order is
equivalent to the Ultrapower Axiom. The fact that UA implies the linearity of
the Ketonen order is immediate. (The proof appears in \cref{KOCharsSection}.)
The converse, however, is subtle. Since we will mostly work under the assumption
of UA, this equivalence is itself not that important (although it does show that
all of our results can be proved from an a priori weaker premise). More
important is the proof, which identifies a canonical way to compare a pair of
ultrafilters assuming the linearity of the Ketonen order. 

\section{Preliminary definitions}\label{PrelimsSection2}
\subsection{Tail uniform ultrafilters}\label{TailUniformSection}
A common notational issue we will encounter in this dissertation is that two
ultrafilters may differ only in the sense that they have different underlying
sets. The change-of-space relation, defined below, articulates our tendency to
identify such ultrafilters.
\begin{defn}\label{UltrafilterConcentration}
	Suppose \(U\) is an ultrafilter on \(X\) and \(C\) is a class. We say \(U\)
	{\it concentrates on \(C\)} if \(C\cap X\in U\). If \(C\) is a set and \(U\)
	concentrates on \(C\), the {\it projection of \(U\) on
	\(C\)}\index{Projection of an ultrafilter}\index{\(U\mid C\) (projection of
	an ultrafilter on a set)} is the ultrafilter \(U\mid C = \{A\subseteq C :
	A\cap X\in U\}\).\index{Ultrafilter!concentrating on a class}
\end{defn}

\begin{defn}\label{EmbeddingEquivalenceDef}
	The {\it change-of-space relation}\index{\(\KE\) (change-of-space relation)}
	is defined on ultrafilters \(U\) and \(W\) by setting \(U \KE W\) if \(U = W
	\mid X\) where \(X\) is the underlying set of \(U\).
\end{defn}

\begin{lma}\label{EmbeddingEquivalence}
	Suppose \(U\) and \(W\) are ultrafilters. Then the following are equivalent:
	\begin{enumerate}[(1)]
		\item \(U\KE W\).
		\item For some set \(S\in U\cap W\), \(U\cap P(S) = W\cap P(S)\).
		\item For all sets \(A\), \(\id_U \in j_U(A)\) if and only if \(\id_W\in
		j_W(A)\).
		\item There is a comparison \((k,h)\) of \((j_U,j_W)\) such that
		\(k(\id_U) = h(\id_W)\).\qed
	\end{enumerate}
\end{lma}
The change-of-space relation is therefore an equivalence relation on
ultrafilters.

The Ketonen order will be a partial order on the class of ultrafilters on
ordinals. On such general ultrafilters, however, the (nonstrict) Ketonen order
is only a preorder, due to the existence of \(\KE\)-equivalent ultrafilters (see
\cref{StrictvsNonStrict2}). Thus we sometimes restrict further to those
ultrafilters that are uniform in a slightly nonstandard sense:

\begin{defn}A filter \(F\) on an ordinal \(\delta\) is {\it tail uniform}\index{Tail uniform} if it contains \(\delta \setminus\alpha\) for every \(\alpha < \delta\).\end{defn}

For any ordinal \(\delta\), the {\it tail filter}\index{Tail filter} on
\(\delta\) is the filter generated by sets of the form \(\delta\setminus
\alpha\) for \(\alpha < \delta\). A filter is therefore tail uniform if it
extends the tail filter. Equivalently, \(F\) is tail uniform if every element of
\(F^+\) is cofinal in \(\alpha\). 

For example, the principal ultrafilter on \(\alpha+1\) concentrated at
\(\alpha\) is uniform.

\begin{defn}\label{SpaceDef} If \(U\) is an ultrafilter that concentrates on
ordinals, then \(\delta_U\)\index{\(\delta_U\)} denotes the least ordinal
\(\delta\) on which \(U\) concentrates.
\end{defn}

\begin{lma} If \(U\) is an ultrafilter that concentrates on ordinals, then \(U\) is tail uniform if and only if \(\delta_U\) is the underlying set of \(U\).\qed\end{lma}

The key property of tail uniform ultrafilters, which is quite obvious, is that
they yield canonical representatives of \(\KE\) equivalence classes of
ultrafilters concentrating on ordinals.
\begin{lma}
	For any ultrafilter \(U\) that concentrates on ordinals, then \(U\mid
	\delta_U\) is the unique tail uniform ultrafilter \(W\) such that \(U \KE
	W\). In particular, if \(U\) and \(W\) are tail uniform ultrafilters such
	that \(U\KE W\), then \(U = W\).\qed
\end{lma}

There is an obvious but useful characterization of \(\delta_U\) in terms of
elementary embeddings:
\begin{lma}\label{SpaceChar}
If \(U\) is an ultrafilter that concentrates on ordinals, then \(\delta_U\) is
the least ordinal \(\delta\) such that \(\id_U < j_U(\delta)\).\qed
\end{lma} 

\begin{defn}
The class of countably complete tail uniform ultrafilters\index{\(\Un\)
(countably complete tail uniform ultrafilters)} is denoted by \(\Un\).
\end{defn}

Let us just point out that tail uniformity and uniformity are not the same
concept, and moreover neither is a strengthening of the other. The simplest way
to separate these concepts is by considering the Fr\'echet and tail filters
themselves. For any set \(X\), let \(F_X\) denote the Fr\'echet filter on \(X\).
For any ordinal \(\alpha\), let \(T_\alpha\) denote the tail filter on
\(\alpha\).

\begin{lma}
Suppose \(\lambda\) is an ordinal.
\begin{itemize}
\item \(T_\lambda \subseteq F_\lambda\) if and only if \(\lambda\) is a
cardinal.
\item \(F_\lambda \subseteq T_\lambda\) if and only if \(|\lambda| =
\textnormal{cf}(\lambda)\). 
\end{itemize}
Thus \(T_\lambda = F_\lambda\) if and only if \(\lambda\) is a regular cardinal.
If \(\lambda\) is a singular cardinal, \(T_\lambda\) is tail uniform but not
uniform. If \(\lambda\) is not a cardinal, then \(F_\lambda\) is uniform but not
tail uniform.\qed
\end{lma}

One can easily obtain ultrafilters that are counterexamples to the equivalence
of tail uniformity and true uniformity by combining the previous lemma with the
Ultrafilter Lemma.

\subsection{Limits of ultrafilters}\label{LimitSection}
The following definition comes from classical ultrafilter theory:
\begin{defn}\label{LimitDef}
Suppose \(W\) is an ultrafilter, \(I\) is a set in \(W\), and \(\langle U_i :
i\in I\rangle\) is a sequence of ultrafilters on a set \(X\). The {\it
\(W\)-limit of \(\langle U_i : i\in I\rangle\)}
\index{Limit}
\index{\(U\text{-}\lim\) (\(U\)-limit)} is the ultrafilter \[W\text{-}\lim_{i\in
I} U_i = \{A\subseteq X: \{i\in I : A\in U_i\}\in W\}\]
\end{defn}

It is often easier to think about limits in terms of elementary embeddings:
\begin{lma}\label{LimitInverse}
Suppose \(W\) is an ultrafilter, \(I\) is a set in \(W\), and \(\langle U_i :
i\in I\rangle\) is a sequence of ultrafilters on a fixed set \(X\). Then 
\[W\text{-}\lim_{i\in I} U_i = j_W^{-1}[Z]\] where \(Z = [\langle U_i : i\in
I\rangle]_W\).
\begin{proof}
Suppose \(A\subseteq X\). Then
\begin{align*}
A\in W\text{-}\textstyle\lim_{i\in I} U_i &\iff  A\in U_i\text{ for \(W\)-almost all }i\in I\\
&\iff j_W(A)\in [\langle U_i : i\in I\rangle]_W\\
&\iff A\in j_W^{-1}[Z]
\end{align*}
where the middle equivalence follows from \L o\'s's Theorem.
\end{proof}
\end{lma}

Limits generalize the usual derived ultrafilter and pushforward constructions:
\begin{defn}\label{PrincipalDef}\index{Principal ultrafilter}\index{\(\pr {a} {X}\) (principal ultrafilter)}
	Suppose \(X\) is a set and \(a\in X\). The {\it principal ultrafilter on
	\(X\) concentrated at \(a\)} is the ultrafilter \(\pr a X = \{A\subseteq X :
	a\in A\}\).
\end{defn}

\begin{defn}\label{PushDef}\index{Pushforward}
	Suppose \(W\) is an ultrafilter, \(I\) is a set in \(W\), and \(f : I\to X\)
	is a function. Then the {\it pushforward} of \(W\) by \(f\) is the
	ultrafilter \(f_*(W) = \{A\subseteq X: f^{-1}[A]\in W\}\).
\end{defn}

The following lemmas relate the derived ultrafilter construction to inverse
images, limits, and pushforwards.

\begin{lma}\label{DerivedPrincipal}\index{Derived ultrafilter!as an inverse image}
	Suppose \(N\) and \(P\) are transitive models of \textnormal{ZFC}, \(X\) is
	a set in \(N\), \(i : N \to P\) is an elementary embedding, and \(a\in
	i(X)\). Then the \(N\)-ultrafilter on \(X\) derived from \(i\) using \(a\)
	is simply \(i^{-1}[\pr a {i(X)}]\).\qed
\end{lma}

\begin{lma}\label{PushDerived}\index{Derived ultrafilter!as a pushforward}
	Suppose \(W\) is an ultrafilter, \(I\) is a set in \(W\), and \(f : I\to X\)
	is a function. Then 
	\[f_*(W) = W\text{-}\lim_{i\in I} \pr {f(i)} {X} = j_W^{-1}[\pr {[f]_W}
	{j_W(X)}]\] In other words, \(f_*(W)\) is the  ultrafilter on \(X\) derived
	from \(j_W\) using \([f]_W\).\qed
\end{lma}

These lemmas are trivial, but it turns out that many calculations are
significantly simpler when one treats limits and derived ultrafilters uniformly
as inverse images.

To be really pedantic, the reader might point out that for example in
\cref{PushDerived}, \(\pr {[f]_W} {j_W(X)}\) is not an \(M_W\)-ultrafilter but a
\(V\)-ultrafilter. Moreover if \(M_W\) is not wellfounded, then the statement
\([f]_W\in j_W(X)\) so \(\pr {[f]_W} {j_W(X)}\) is not well-defined. Of course,
\(\pr {[f]_W} {j_W(X)}\)  really denotes \((\pr {[f]_W} {j_W(X)})^{M_W}\). For
the reader's own sake, we will try to omit all these superscripts in our
notation for principal ultrafilters when they can be guessed from context. For
example, in \cref{PushDerived}, we would usually write: \[f_*(W) =
W\text{-}\lim_{i\in I} \pr {f(i)} {} = j_W^{-1}[\pr {[f]_W} {}]\]

The key to understanding derived ultrafilters is to consider the natural factor
embeddings associated to them. There is a generalization of the factor embedding
construction to the case of limits. In fact, this works somewhat more generally
for arbitrary inverse images of ultrafilters:

\begin{lma}\label{LimitFactor}\index{Factor embedding!associated to a limit}
Suppose \(N\) and \(P\) are transitive models of \textnormal{ZFC}, \(X\) is a
set in \(N\), \(i : N \to P\) is an elementary embedding, and \(D\) is a
\(P\)-ultrafilter on \(i(X)\). Let \(U = i^{-1}[D]\). There is a unique
elementary embedding \(k : M_U^N\to M^P_D\) such that \(k(\id_U) = \id_D\) and
\(k\circ j_U^N = j_D^P\circ i\).
\begin{proof}
For any function \(f\in N\) defined on a set in \(U\), set \[k([f]^N_U) =
[i(f)]^P_D\] It is immediate from this definiton that \(k(\id_U) = \id_D\) and
\(k\circ j_U^N = j_D^P\circ i\). We must show that \(k\) is well-defined and
elementary. This follows from the usual calculation:
\begin{align*}
M_U^N\vDash \varphi([f_1]^N_U,\dots,[f_n]^N_U)&\iff N\vDash \varphi(f_1(x),\dots,f_n(x))\text{ for \(U\)-almost all \(x\)}\\
&\iff P\vDash \varphi(i(f_1)(x),\dots,i(f_n)(x))\text{ for \(D\)-almost all \(x\)}\\
&\iff M_U^P\vDash \varphi([i(f_1)]_D^P,\dots,[i(f_n)]_D^P)\qedhere
\end{align*}
\end{proof}
\end{lma}
\section{The Ketonen order}\label{KetonenOrderSection}
\subsection{Characterizations of the Ketonen order}\label{KOCharsSection}
Let us begin our investigation of the Ketonen order with a purely combinatorial
definition.
\begin{defn}
	Suppose \(X\) is a set and \(A\) is a class. Then \(\mathscr B(X)\) denotes
	the set of countably complete ultrafilters on \(X\), and \(\mathscr
	B(X,A)\)\index{\(\mathscr B(X),\mathscr B(X,A)\) (countably complete
	ultrafilters on \(X\), concentrating on \(A\))} denotes the set of countably
	complete ultrafilters on \(X\) that concentrate on \(A\).
\end{defn}

\begin{defn}\label{KODef}
	Suppose \(\delta\) is an ordinal. The {\it Ketonen order}\index{Ketonen
	order} is defined on \(\mathscr B(\delta)\) as follows. For \(U,W\in
	\mathscr B(\delta)\):
	\begin{itemize} 
		\item \(U\sE W\) if there is a set \(I\in W\) and a sequence \(\langle
		U_\alpha : \alpha \in I\rangle\in \prod_{\alpha\in I}\mathscr
		B(\delta,\alpha)\) such that \(U = W\text{-}\lim_{i\in I} U_\alpha\).
		\item \(U\E W\) if there is a set \(I\in W\) and a sequence \(\langle
		U_\alpha : \alpha \in I\rangle\in \prod_{\alpha\in I}\mathscr
		B(\delta,\alpha+1)\) such that \(U = W\text{-}\lim_{i\in I} U_\alpha\).
	\end{itemize}
\end{defn}
We refer to \(\sE\) and \(\E\) as the {\it strict} and {\it non-strict} Ketonen
orders.

Of course one could take \(I = \delta\setminus \{0\}\) in the first bullet-point
and \(I = \delta\) in the second, but the combinatorics are typically clearer if
one does not make this demand.

There is perhaps a potential ambiguity in our notation, since the order \(\sE\)
depends on the ordinal \(\delta\), which we suppress in our notation. This
dependence is always immaterial, however, since there are canonical embeddings
between the various Ketonen orders. These embeddings allow us to spin all these
orders together into one (\cref{GlobalKODef}).

Let us first explain the straightforward relationship between the strict and
nonstrict Ketonen orders.
\begin{prp}\label{StrictvsNonStrict}
	Suppose \(\delta\) is an ordinal and \(U,W\in \mathscr B(\delta)\). Then
	\(U\E W\) if and only if \(U \sE W\) or \(U = W\).
	\begin{proof}
		Let \(\langle U_\alpha : \alpha \in I\rangle\in \prod_{\alpha\in
		I}\mathscr B(\delta,\alpha+1)\) witness \(U\E W\). Let \[J = \{i\in I :
		U_\alpha\in \mathscr B(\delta,\alpha)\}\] If \(J\in W\), then \(\langle
		U_\alpha : \alpha \in J\rangle\in \prod_{\alpha\in J}\mathscr
		B(\delta,\alpha)\) witnesses \(U\sE W\). 
		
		Assume therefore that \(J\notin W\). For all \(\alpha\in I\setminus J\),
		\(U_\alpha \in \mathscr B(\delta,\alpha+1)\setminus \mathscr
		B(\delta,\alpha)\). Note that \(\mathscr B(\delta,\alpha+1)\setminus
		\mathscr B(\delta,\alpha)\) contains only the principal ultrafilter
		\(\pr \alpha \delta\), and hence \(U_\alpha = \pr \alpha\delta\) for
		\(\alpha\in I\setminus J\). Thus 
		\[U = W\text{-}\lim_{\alpha\in I}U_\alpha = W\text{-}\lim_{\alpha\in
		I\setminus J}\pr \alpha {} = W\] where the final equality follows easily
		from the definitions (or from \cref{PushDerived}).
	\end{proof}
\end{prp}

We therefore focus our attention on the strict Ketonen order \(\sE\) for now.
Before establishing its basic order-theoretic properties, let us give some
fairly obvious alternate characterizations of it. We think the characterization
\cref{KOChars} (2) is quite elegant in that it demonstrates a basic relationship
between the Ketonen order, the covering properties of ultrapowers, and
extensions of filter bases to countably complete ultrafilters, foreshadowing the
powerful interactions between strong compactness and the Ultrapower Axiom that
we will see in \cref{SCChapter1} and \cref{SCChapter2}. \cref{KOChars} (3) and
(4) are more useful, though, linking the Ketonen order and the Ultrapower Axiom
through the concept of a comparison (\cref{ComparisonDef}).
\begin{lma}\label{KOChars}
	Suppose \(\delta\) is an ordinal and \(U,W\in \mathscr B(\delta)\). The
	following are equivalent:
	\begin{enumerate}[(1)]
		\item \(U\sE W\).
		\item \(j_W[U]\) extends to an \(M_W\)-ultrafilter \(Z \in \mathscr
		B^{M_W}(j_W(\delta),\id_W)\). 
		\item There is a comparison \((k,h) : (M_U,M_W)\to P\) of \((j_U,j_W)\)
		such that \(h\) is an internal ultrapower embedding of \(M_W\) and
		\(k(\id_U) < h(\id_W)\).
		\item There is a comparison \((k,h) : (M_U,M_W)\to P\) of \((j_U,j_W)\)
		such that \(h\) is close to \(M_W\) and \(k(\id_U) < h(\id_W)\).
	\end{enumerate}
\begin{proof}
	{\it (1) implies (2):} Fix \(I\in W\) and \(\langle U_\alpha : \alpha\in
	I\rangle\in \prod_{\alpha\in I}\mathscr B(\delta,\alpha)\) witnessing \(U\sE
	W\). Let \(Z = [\langle U_\alpha: \alpha\in I\rangle]_W\). By \L o\'s's
	Theorem, \(Z\in \mathscr B^{M_W}(j_W(\delta),\id_W)\), and by
	\cref{LimitInverse}, \(j_W^{-1}[Z] = W\text{-}\lim_{i\in I} U_i  = U\). This
	implies \(j_W[U]\subseteq Z\).
	
	{\it (2) implies (1):} Similar. 
	
	{\it (2) implies (3):} Fix \(Z \in \mathscr B^{M_W}(j_W(\delta),\id_W)\)
	such that \(j_W[U]\subseteq Z\). Because of the basic structure of
	ultrafilters, the fact that \(j_W[U]\subseteq Z\) implies that \(j_W^{-1}[Z]
	= U\). Let \(h : M_W\to N\) be the ultrapower of \(M_W\) by \(Z\). Since
	\(Z\) concentrates on \(\id_W\), \(\id_Z < h(\id_W)\). By
	\cref{LimitFactor}, there is a unique elementary embedding \(k : M_U\to N\)
	such that \(k(\id_U) = \id_Z\) and \(k\circ j_U = h\circ j_W\). The former
	equation implies \(k(\id_U) < h(\id_W)\), while the latter equation says
	that \((k,h)\) is a comparison of \((j_U,j_W)\). Therefore (3) holds.
	
	{\it (3) implies (4):} Internal ultrapower embeddings are close.
	
	{\it (4) implies (2):} Let \(Z\) be the \(M_W\)-ultrafilter on
	\(j_W(\delta)\) derived from \(h\) using \(k(\id_U)\). Thus \(Z = h^{-1}[\pr
	{k(\id_U)} {}]\). (Here \(\pr {k(\id_U)} {}\) denotes the principal
	ultrafilter on \(k(j_U(\delta))\) concentrated at \(k(\id_U)\); see
	\cref{PrincipalDef} and the ensuing discussion.) Since \(h\) is close, \(Z\)
	belongs to \(M_W\), and since \(k(\id_U) < h(\id_W)\), \(Z\) concentrates on
	\(\id_W\). Thus \(Z\in \mathscr B^{M_W}(j_W(\delta),\id_W)\). Moreover,
	\[j_W^{-1}[Z] = j_W^{-1}[h^{-1}[\pr {k(\id_U)} {}]] =  j_U^{-1}[k^{-1}[\pr
	{k(\id_U)} {}]] = j_U^{-1}[\pr {\id_U} {}] = U\] In particular,
	\(j_W[U]\subseteq Z\), which shows (2).
\end{proof}
\end{lma}

Of course, there are identical characterizations for the nonstrict Ketonen order
as well:
\begin{lma}\label{NSKOChars}
	Suppose \(\delta\) is an ordinal and \(U,W\in \mathscr B(\delta)\). The
	following are equivalent:
	\begin{enumerate}[(1)]
		\item \(U\E W\).
		\item \(j_W[U]\) extends to an \(M_W\)-ultrafilter \(Z\in \mathscr
		B^{M_W}(j_W(\delta),\id_W+1)\). 
		\item There is a comparison \((k,h) : (M_U,M_W)\to P\) of \((j_U,j_W)\)
		such that \(h\) is an internal ultrapower embedding of \(M_W\) and
		\(k(\id_U) \leq h(\id_W)\).
		\item There is a comparison \((k,h) : (M_U,M_W)\to P\) of \((j_U,j_W)\)
		such that \(h\) is close to \(M_W\) and \(k(\id_U) \leq h(\id_W)\).\qed
	\end{enumerate}
\end{lma}

\cref{KOChars} and \cref{NSKOChars} lead to the central linearity theorem for
the Ketonen order under UA:
\begin{thm}[UA]\label{Totality}\index{Ketonen order!linearity}
	Suppose \(\delta\) is an ordinal and \(U,W\in \mathscr B(\delta)\). Either
	\(U\sE W\) or \(W\E U\).
	\begin{proof}
		Let \((k,h) : (M_U,M_W)\to N\) be an internal ultrapower comparison of
		\((j_U,j_W)\). If \(k(\id_U) < h(\id_W)\), then \cref{KOChars} (3)
		implies \(U\sE W\). Otherwise, \(h(\id_W)\leq k(\id_U)\) and so \(W\E
		U\) by \cref{NSKOChars}.
	\end{proof}
\end{thm}

This linearity theorem is only interesting, of course, if we know that the
Ketonen order is ``well defined": if \(U \sE W\) and \(W\sE U\) held for all
\(U,W\in \mathscr B(\delta)\), it wouldn't be very useful. We now show that in
fact the Ketonen order is a wellfounded partial order.
\subsection{Basic properties of the Ketonen order}\label{BasicSection}
We state the main theorem of this section, which we will prove in pieces:
\begin{thm*}
	For any ordinal \(\delta\), \((\mathscr B(\delta),\sE)\) is a strict
	wellfounded partial order.
\end{thm*}

Thus we must show the following facts:
\begin{prp}\label{KOTransitive}
	For any ordinal \(\delta\), \(\sE\) is a transitive relation on \(\mathscr
	B(\delta)\).
\end{prp}

\begin{thm}\label{KOWellfounded}
	For any ordinal \(\delta\), \(\sE\) is a wellfounded relation on \(\mathscr
	B(\delta)\).
\end{thm}

Let us warm up to this by proving irreflexivity:
\begin{prp}\label{KOIrreflexive}
	For any ordinal \(\delta\), \(\sE\) is an irreflexive relation on \(\mathscr
	B(\delta)\).
	\begin{proof}
		Suppose towards a contradiction that \(U\in \mathscr B(\delta)\)
		satisfies \(U\sE U\). Fix \(I\in U\), and \(\langle U_\alpha : \alpha\in
		I\rangle\in \prod_{\alpha\in I}  \mathscr B(\delta,\alpha)\) such that
		\[U = U\text{-}\lim_{\alpha\in I}U_\alpha\] Define \(A\subseteq \delta\)
		by induction: put \(\alpha\in A\) if and only if \(A\cap\alpha\notin
		U_\alpha\). Then \begin{align*}A\in U&\iff \{\alpha \in I : A\in
		U_\alpha\}\in U\\
			&\iff \{\alpha\in I : A\cap \alpha\in U_\alpha\} \in U\\ &\iff
			\{\alpha\in I : \alpha\notin A\}\in U\\
			&\iff I\setminus A\in U\end{align*} Since \(I\in U\) and \(U\) is an
		ultrafilter, either \(A\) or \(I\setminus A\) must belong to \(U\). Thus
		both belong to \(U\), contradicting that \(U\) is closed under
		intersections.
	\end{proof}
\end{prp}

Notice that the proof does not use the wellfoundedness of \(U\). We now give two
proofs of the transitivity of the Ketonen order.
\begin{proof}[Proof of \cref{KOTransitive}] Suppose \(U\sE W \E Z\). We will
	show that \(U\sE Z\). Fix the following objects:
	\begin{itemize}
		\item A set \(I\in W\) and a sequence \(\langle U_\alpha : \alpha\in
		I\rangle\in \prod_{\alpha\in I}\mathscr B(\delta,\alpha)\) such that \(U
		= W\text{-}\lim_{\alpha\in I} U_\alpha\). 
		\item A set \(J\in Z\) and a sequence \(\langle W_\beta : \beta\in
		J\rangle\in \prod_{\beta\in J}\mathscr B(\delta,\beta)\) such that \(W=
		Z\text{-}\lim_{\beta\in J} Z_\beta\). 
	\end{itemize}
	
	Since \(I\in W = Z\text{-}\lim_{\beta\in J} W_\alpha\), the set \(J' =
	\{\beta\in J : I\in W_\beta\}\) belongs to \(Z\). For \(\beta\in J'\), we
	can define \(U'_\beta = W_\beta\text{-}\lim_{\alpha\in I} U_\alpha\). Thus: 
	\begin{align*}
		U &= W\text{-}\lim_{\alpha\in I} U_\alpha\\
			&= (Z\text{-}\lim_{\beta\in J} W_\alpha)\text{-}\lim_{\alpha\in I} U_\alpha\\
			&= Z\text{-}\lim_{\beta\in J'} (W_\alpha\text{-}\lim_{\alpha\in I} U_\alpha)\\
			&= Z\text{-}\lim_{\beta\in J'}U_\alpha'
	\end{align*}
	
	Finally, if \(\beta\in J'\), then \(\{\alpha\in I : U_\alpha\in \mathscr
	B(\delta,\beta)\}\supseteq I\cap (\beta + 1)\in W_\beta\), so \[\langle
	U'_\beta : \beta\in J'\rangle\in \prod_{\beta\in J'} \mathscr
	B(\delta,\beta)\] Therefore \(\langle U'_\beta : \beta\in J'\rangle\)
	witnesses \(U\sE Z\).
\end{proof}

We are still just warming up, so let us give another proof of the transitivity
of the Ketonen order that is more diagrammatic:
\begin{proof}[Alternate Proof of \cref{KOTransitive}] Using \cref{KOChars}, fix
	the following objects:
	\begin{itemize}
		\item A comparison \((k_0,h_0): (M_U,M_W)\to N_0\) of \((j_U,j_W)\) such
		that \(h_0\) is an internal ultrapower embedding of \(M_W\) and
		\(k_0(\id_U) < h_0(\id_W)\).
		\item A comparison \((k_1,h_1) : (M_W,M_Z)\to N_1\) of \((j_W,j_Z)\)
		such that \(h_1\) is an internal ultrapower embedding of \(M_Z\) and
		\(k_1(\id_W)\leq h_1(\id_Z)\).
	\end{itemize}
	
	The rest of the proof is contained in \cref{KOTransitiveFig}.
	\begin{figure}
	\center
	\includegraphics[scale=.6]{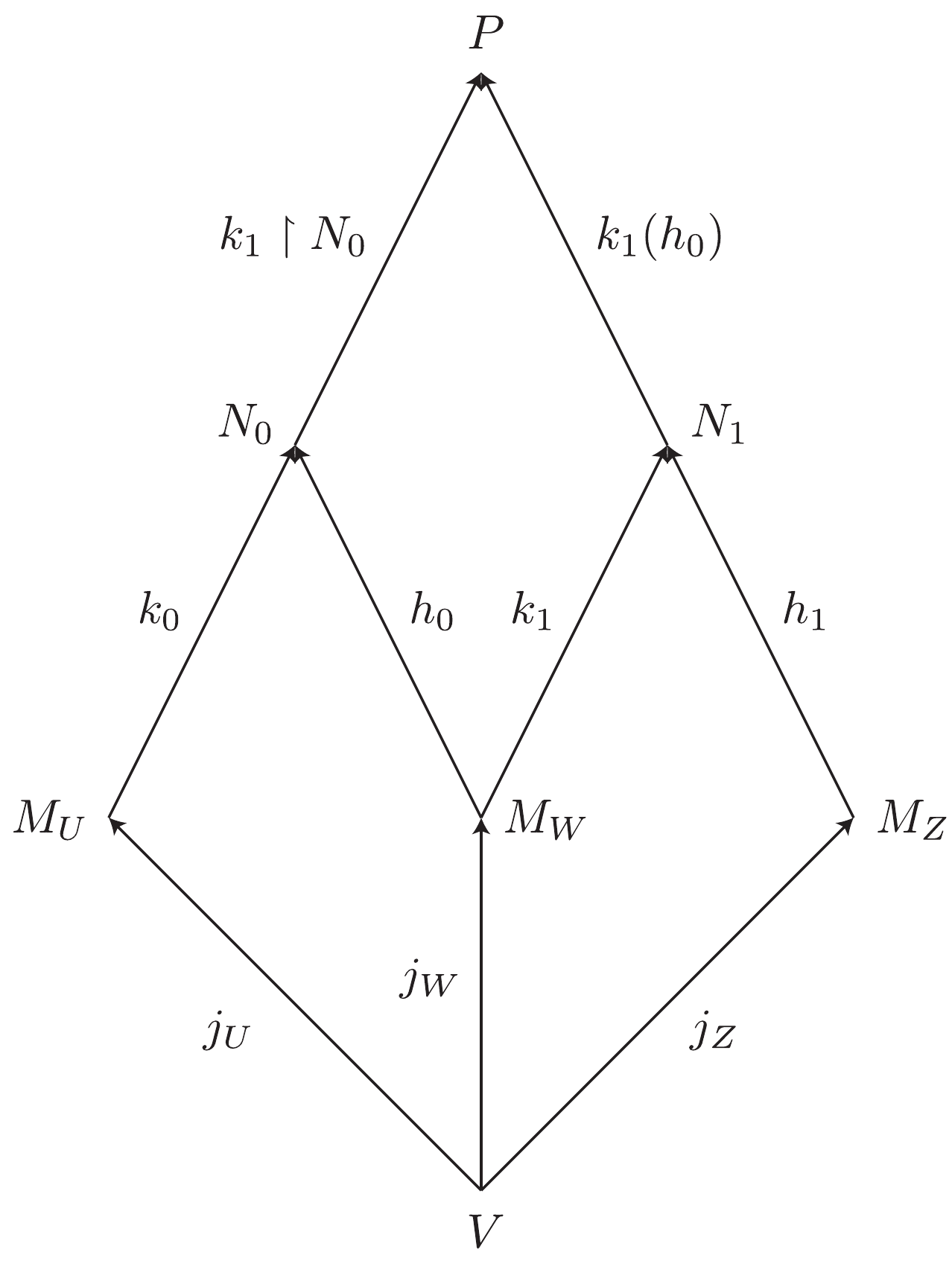}
	\caption{The transitivity of the Ketonen order}\label{KOTransitiveFig}
	\end{figure}
	Consider the embeddings \(h_0 : M_W\to N_0\) and \(k_1 : M_W\to N_1\). There
	is a very general construction that yields a comparison of \((h_0,k_1)\).
	Since \(h_0\) is amenable to \(M_W\), one can define \(k_1(h_0) : N_1\to
	k_1(N_0)\) by shifting the fragments of \(h_0\) using \(k_1\). The
	well-known identity \((k_1\restriction N_0)\circ h_0 = k_1(h_0)\circ k_1\)
	implies that \((k_1\restriction N_0,k_1(h_0)) : (N_0,N_1)\to k_1(N_0)\) is a
	comparison of \((h_0,k_1)\).
	
	It follows easily that \(((k_1\restriction N_0)\circ k_0,k_1(h_0)\circ
	h_1)\) is a comparison of \((j_U,j_Z)\). Easily \(k_1(h_0)\circ h_0\) is an
	internal ultrapower embedding of \(M_Z\). Finally \[(k_1\restriction
	N_0)\circ k_0(\id_U) < (k_1\restriction N_0)\circ h_0(\id_W) = k_1(h_0)
	\circ k_1(\id_W) \leq k_1(h_0)\circ h_1(\id_Z)\] Thus \(U\sE Z\) by
	\cref{KOChars}.
\end{proof}

We finally turn to wellfoundedness. We will give a combinatorial proof, but the
reader can consult \cref{ExtendedKetonen} for a diagrammatic approach in a more
general context. The proof proceeds by iterating the following {\it strong
transitivity lemma} for the Ketonen order, abstracted from the proof of
\cref{KOTransitive}:
\begin{lma}\label{StrongTrans}\index{Ketonen order!strong transitivity}
	Suppose \(\delta\) is an ordinal, \(U,W\in \mathscr B(\delta)\), and \(U\sE
	W\). Suppose \(Z\) is an ultrafilter, \(J\) is a set in \(Z\), and \(\{ W_x
	: x\in J\} \subseteq\mathscr B(\delta)\) is a sequence such that \[W =
	Z\text{-}\lim_{x\in J}W_x\] Then there is a set \(J'\subseteq J\) in \(Z\)
	and \(\{ U_x : x\in J'\}\subseteq \mathscr B(\delta)\) with \(U_x\sE W_x\)
	for all \(x\in J'\) such that 
	\[U = Z\text{-}\lim_{x\in J'} U_x\]
	\begin{proof}[Sketch] Fix \(I\in W\) and \(\langle U_\alpha : \alpha\in
		I\rangle\in \prod_{\alpha\in I} \mathscr B(\delta,\alpha)\) such that
		\(U = W\text{-}\lim_{\alpha\in I} U_\alpha\). Let \(J' = \{x\in J : I\in
		W_x\}\). For \(x\in J'\), let \(U_x = W_x\text{-}\lim_{\alpha\in I}
		U_\alpha\). Then \(\langle U_\alpha : \alpha\in I\rangle\) witnesses
		that \(U_x\sE W_x\). Moreover the calculation in \cref{KOTransitive}
		shows that \(U = Z\text{-}\lim_{x\in J'} U_x\).
	\end{proof}
\end{lma}

This is more elegantly stated using elementary embeddings:
\begin{lma}
	Suppose \(\delta\) is an ordinal, \(U,W\in \mathscr B(\delta)\), and \(U\sE
	W\). Suppose \(j : V\to M\) is an elementary embedding and \(W_*\in
	j(\mathscr B(\delta))\) extends \(j[W]\). Then there is some \(U_*\in
	j(\mathscr B(\delta))\) extending \(j[U]\) such that \(M\vDash U_*\sE
	W_*\).\qed
\end{lma}

Recall now the notation \(U\mid C\) from \cref{UltrafilterConcentration},
denoting the projection of an ultrafilter \(U\) to a set \(C\) on which it
concentrates. We will need the following trivial lemma, which is also implicit
in the proof of \cref{KOTransitive}:
\begin{lma}\label{KetonenProjection}
	Suppose \(\epsilon\) and \(\delta\) are ordinals, \(U\in \mathscr
	B(\delta)\), and \(W\in \mathscr B(\delta,\epsilon)\). If \(U\E W\), then
	\(U\in \mathscr B(\delta,\epsilon)\) and \(U\mid \epsilon \E W\mid
	\epsilon\) in the Ketonen order on \(\mathscr B(\epsilon)\).
	\begin{proof}
		Fix \(I\in W\) and \(\langle U_\alpha : \alpha\in I\rangle\in \prod
		\mathscr B(\delta,\alpha)\) such that \(U = W\text{-}\lim_{\alpha\in I}
		U_\alpha\). Then since \(U = W\text{-}\lim_{\alpha\in I} U_\alpha =
		W\text{-}\lim_{\alpha\in I\cap \epsilon} U_\alpha\) is a limit of
		ultrafilters concentrating on \(\epsilon\), so \(U\) itself concentrates
		on \(\epsilon\). Moreover \(\langle U_\alpha\mid \epsilon : \alpha\in
		I\cap \epsilon\rangle\) witnesses that \(U\mid \epsilon \sE W\mid
		\epsilon\) in the Ketonen order on \(\mathscr B(\epsilon)\).
	\end{proof}
\end{lma}

As we prove \cref{KOWellfounded}, the reader may profit from the observation
that the proof consists of the combinatorial core of the proof of the
wellfoundedness of the Mitchell order on normal ultrafilters, stripped of all
applications of normality and \L o\'s's Theorem.
\begin{proof}[Proof of \cref{KOWellfounded}]\index{Ketonen order!wellfoundedness}
	Assume towards a contradiction that there is an ordinal \(\delta\) such that
	\(\sE\) is illfounded on \(\mathscr B(\delta)\). Fix the least such
	\(\delta\). Choose a sequence \(\{ U_n : n < \omega\}\subseteq  \mathscr
	B(\delta)\) that is \(\sE\)-descending:
	\[U_0\sgE U_1 \sgE U_2 \sgE\cdots\]
	
	For each positive number \(n\), we will define by recursion a set \(J_n\in
	U_0\) and a sequence of ultrafilters \(\langle U^n_\alpha : \alpha\in
	J_n\rangle\in \prod_{\alpha\in J_n} \mathscr B(\delta,\alpha)\) such that
	for all \(n < \omega\), the following hold:
	\begin{itemize}
		\item \(U_n = U\text{-}\lim_{\alpha\in J_n} U^n_\alpha\). 
		\item If \(n > 1\), then \(J_n\subseteq J_{n-1}\) and for all
		\(\alpha\in J_{n}\), \(U^{n}_\alpha \sE U^{n-1}_\alpha\).
	\end{itemize}
	To start, fix \(J_1\in U_0\) and \(\langle U^1_\alpha : \alpha\in
	J_1\rangle\in \prod_{\alpha\in J_1} \mathscr B(\delta,\alpha)\) witnessing
	that \(U_1\sE U_0\); that is,
	\[U_1 = U_0\text{-}\lim_{\alpha\in J_1} U_\alpha^1\] Suppose \(n > 1\) and
	\(J_{n-1}\in U_0\) and \(\langle U^{n-1}_\alpha : \alpha\in
	J_{n-1}\rangle\in \prod_{\alpha\in J_{n-1}} \mathscr B(\delta,\alpha)\) have
	been defined so that \(U_{n-1} = U\text{-}\lim_{\alpha\in J_{n-1}}
	U^{n-1}_\alpha\). \cref{StrongTrans} (with \(U = U_{n}\), \(W = U_{n-1}\),
	and \(Z = U_0\)) yields \(J_{n}\subseteq J_{n-1}\) and \(\{U^{n}_\alpha :
	\alpha \in J_n\}\subseteq \mathscr B(\delta)\) such that the two bullet
	points above are satisfied. We must verify that \(\langle U^{n}_\alpha :
	\alpha\in J_n\rangle\in \prod_{\alpha\in J_n} \mathscr B(\delta,\alpha)\).
	But for any \(\alpha\in J^n\), \(U^n_\alpha \sE U^{n-1}_\alpha\in \mathscr
	B(\delta,\alpha)\), and therefore \(U^n_\alpha\in \mathscr
	B(\delta,\alpha)\) by \cref{KetonenProjection}, as desired. This completes
	the recursive definition.
	
	Now let \(J = \bigcap_{n < \omega} J_n\). For any \(\alpha\in J\), we have
	\[U^1_\alpha\sgE U^2_\alpha \sgE U^3_\alpha \sgE\cdots\] by the second
	bullet point above. Since \(U^n_\alpha\in \mathscr B(\delta,\alpha)\) for
	all \(n < \omega\), \cref{KetonenProjection} implies \[U^1_\alpha\mid \alpha
	\sgE U^2_\alpha\mid \alpha \sgE U^3_\alpha\mid \alpha\sgE \cdots \] Thus the
	restriction of \(\sE\) to \(\mathscr B(\alpha)\) is illfounded. This
	contradicts the minimality of \(\delta\).
\end{proof}
Observe that the proof of \cref{KOWellfounded} goes through in ZF + DC. The
structure of countably complete ultrafilters on ordinals is of great interest in
the context of the Axiom of Determinacy, and so the existence of a combinatorial
analog of the Mitchell order in that context raises a number of very interesting
structural questions that we will not pursue seriously in this dissertation.
\subsection{The global Ketonen order}
\cref{KetonenProjection} above suggests extending the Ketonen order to an order
on ultrafilters that is agnostic about the underlying sets of the ultrafilters
involved:
\begin{defn}\label{GlobalKODef}
	Suppose \(U\) and \(W\) are countably complete ultrafilters on ordinals. The
	(global) {\it Ketonen order}\index{Ketonen order!global} is defined as
	follows:
	\begin{itemize}
		\item \(U\sE W\) if \(U\mid \delta\sE W\mid \delta\).
		\item \(U\E W\) if \(U\mid \delta\E W\mid \delta\).
	\end{itemize}
	where \(\delta\) is any ordinal such that \(U\) and \(W\) both concentrate
	on \(\delta\).
\end{defn}

By \cref{KetonenProjection}, this definition does not conflict with our original
definition of the Ketonen order on \(\mathscr B(\delta)\). In fact, various
characterizations of the Ketonen order from \cref{KOChars} translate smoothly to
this context:

\begin{lma}\label{GlobalKOChars}
	Suppose \(\epsilon\) and \(\delta\) are ordinals, \(U\in \mathscr
	B(\epsilon)\), and \(W\in \mathscr B(\delta)\). Then the following are
	equivalent:
	\begin{enumerate}[(1)]
		\item \(U\sE W\).
		\item There exist \(I\in W\) and \(\langle U_\alpha : \alpha \in
		I\rangle\in \prod _{\alpha \in I}\mathscr B(\epsilon, \alpha)\) such
		that \(U = W\text{-}\lim_{\alpha\in I}U_\alpha\).
		\item \(j_W[U]\subseteq Z\) extends to an \(M_W\)-ultrafilter \(Z\in
		\mathscr B^{M_W}(j_W(\epsilon), \id_W)\).
		\item There is a comparison \((k,h) : (M_U,M_W)\to P\) of \((j_U,j_W)\)
		such that \(h\) is an internal ultrapower embedding of \(M_W\) and
		\(k(\id_U) < h(\id_W)\).
		\item There is a comparison \((k,h) : (M_U,M_W)\to P\) of \((j_U,j_W)\)
		such that \(h\) is close to \(M_W\) and \(k(\id_U) < h(\id_W)\).\qed
	\end{enumerate}
\end{lma}

We have the following simple relationship between the space of an ultrafilter
and its position in the Ketonen order:
\begin{lma}\label{SpaceLemma}
	Suppose \(U\) and \(W\) are countably complete ultrafilters on ordinals. 
	\begin{itemize}
		\item If \(\delta_U < \delta_W\), then \(U\sE W\).
		\item If \(U\E W\), then \(\delta_U \leq \delta_W\).
	\end{itemize}
	\begin{proof}
		To see the first bullet point, note that for any \(\alpha\in
		[\delta_U,\delta_W)\), \(\alpha \geq \delta_U\) and hence \(U\)
		concentrates on \(\alpha\). Thus the constant sequence \(\langle U :
		\alpha \in [\delta_U,\delta_W)\rangle\) belongs to \(\prod _{\alpha \in
		[\delta_U,\delta_W)}\mathscr B(\epsilon, \alpha)\), and clearly \(U =
		W\text{-}\lim_{\alpha\in [\delta_U,\delta_W)}U\). By
		\cref{GlobalKOChars}, \(U\sE W\).
		
		The second bullet point is immediate from \cref{KetonenProjection}.
	\end{proof}
\end{lma}

The one issue with the global Ketonen order, which presents only notational
difficulties, is that in this generalized context, \(\E\) is no longer the
irreflexive part of \(\sE\). Instead we have the following fact, where \(\KE\)
is the change-of-space relation defined in \cref{EmbeddingEquivalenceDef}:

\begin{lma}\label{StrictvsNonStrict2}
	Suppose \(U\) and \(W\) are countably complete ultrafilters on ordinals.
	Then \(U\E W\) if and only if \(U\sE W\) or \(U \KE W\).\qed
\end{lma}

Since the \(\KE\)-relation convenient to restrict the global Ketonen order to
the class of tail uniform ultrafilters \(\Un\):
\begin{lma}\label{StrictvsNonStrictUn}
	Suppose \(U, W\in \Un\). Then \(U\E W\) if and only if \(U\sE W\) or \(U =
	W\). 
\end{lma}

\begin{defn}
	For any ordinal \(\delta\), let \(\Un(\delta)\) denote the set of tail
	uniform ultrafilters \(U\) such that \(\delta_U\leq \delta\).
\end{defn}

\begin{lma}\label{KOUnGood}
	For all ordinals \(\delta\), the map \(\phi : \mathscr B(\delta)\to
	\Un(\delta)\) defined by \(\phi(U) = U\mid \delta_U\) is an isomorphism from
	\((\mathscr B(\delta), \sE, \E)\) to \((\Un(\delta),\sE, \E)\). Thus the
	Ketonen order is a set-like wellfounded partial order on \(\Un\).\qed
\end{lma}

The following easy lemma generalizes our work in this section, showing that not
only do the various Ketonen orders on \(\mathscr B(\delta)\) cohere, but in
fact, order-preserving maps between ordinals induce order-preserving maps on
their associated Ketonen orders:

\begin{lma}\label{Functorial}
	Suppose \(\epsilon\leq \delta\) are ordinals and \(f : \epsilon\to \delta\)
	is an increasing function. For any \(U,W\in \mathscr B(\epsilon)\), \(U\sE
	W\) in the Ketonen order on \(\mathscr B(\epsilon)\) if and only if
	\(f_*(U)\sE f_*(W)\) in the Ketonen order on \(\mathscr B(\delta)\).
	\begin{proof}[Sketch] Fix \(I\in W\) and \(\langle U_\alpha : \alpha \in
		I\rangle\in \prod_{\alpha\in I}\mathscr B(\epsilon,\alpha)\) such that
		\(U = W\text{-}\lim_{\alpha\in I} U_\alpha\). Let \(J = f[I]\), and for
		\(\alpha \in I\), let \(Z_{f(\alpha)}= f_*(W_\alpha)\). Thus \(J\in
		f_*(W)\). Moreover, for all \(\alpha\in I\), \(f(\alpha)\supseteq
		f[\alpha]\in Z_{f(\alpha)}\) since \(f\) is increasing, so
		\(Z_{f(\alpha)}\in \mathscr B(\delta,f(\alpha))\). Thus \(\langle
		Z_\beta:\beta\in J\rangle\in \prod_{\beta\in J}\mathscr
		B(\delta,\beta)\). Finally \[f_*(U) = W\text{-}\lim_{\alpha\in I}
		f_*(W_\alpha) = f_*(W)\text{-}\lim_{\beta\in f[I]} Z_\alpha\] It follows
		that \(f_*(U)\sE f_*(W)\). The other direction is similar.
	\end{proof}
\end{lma}
\section{Orders on ultrafilters}\label{OtherOrderSection}
In this section, we discuss some generalizations of the Ketonen order and
compare the Ketonen order with other well-known orders. 

\subsection{The Mitchell order}\label{MOKetSection}
The Ketonen order can be seen as a combinatorial generalization of the Mitchell
order on normal ultrafilters. We will discuss the relationship between the
Ketonen order and the generalization of the Mitchell order to arbitrary
countably complete ultrafilters at length in \cref{GMOChapter}, but for now, we
satisfy ourselves by proving that the Ketonen and Mitchell orders coincide on
normal ultrafilters.
\begin{thm}\label{KetMO}
Suppose \(U\) and \(W\) are normal ultrafilters. Then \(U\mo W\) if and only if
\(U\sE W\).\index{Ketonen order!vs. the Mitchell order}
\begin{proof}
Suppose first that \(U\) and \(W\) are normal ultrafilters on distinct cardinals
\(\kappa\) and \(\lambda\). Clearly \(U\mo W\) if and only if \(\kappa <
\lambda\). Moreover by \cref{SpaceLemma}, \(U\sE W\) if and only if \(\kappa <
\lambda\). Thus \(U\mo W\) if and only if \(U\sE W\).

Assume instead that \(U\) and \(W\) lie on the same cardinal \(\kappa\). By
\cref{NormalChar}, \(U\) and \(W\) are \(\kappa\)-complete and \(\kappa  = \id_U
= \id_W\). The key fact we use is that since \(\textsc{crt}(j_W) = \kappa\),
\(j_W(A)\cap \kappa = A\) for all \(A\subseteq \kappa\).

Suppose first that \(U\mo W\). Then \(U\in M_W\). Working in \(M_W\), consider
the projection \(Z = U\mid j_W(\kappa)\in \mathscr B^{M_W}(j_W(\kappa),
\kappa)\). For any \(A\subseteq \kappa\), \(j_W(A)\cap \kappa = A\in U\), or in
other words, \(j_W(A)\in Z\). In other words, \(j_W[U]\subseteq Z\), so by
\cref{KOChars}, \(U\mo W\).

Conversely, suppose \(U\sE W\). Fix \(Z\in \mathscr B^{M_W}(j_W(\kappa),
\kappa)\) such that \(U = j_W^{-1}[Z]\). Suppose \(A\subseteq \kappa\). Then
\(A\in U\) if and only if \(j_W(A)\cap \kappa\in Z\) if and only if \(A\in Z\).
Therefore \(U = Z\mid \kappa\), so \(U\in M_W\). This implies \(U\mo W\).
\end{proof}
\end{thm}

Thus the wellfoundedness of the Mitchell order follows from the wellfoundedness
of the Ketonen order. Notice that this theorem gives another proof of the
linearity of the Mitchell order on normal ultrafilters under UA. Finally, the
proof has the following consequence:
\begin{cor}
	Suppose \(\kappa\) is a cardinal, \(U\in \mathscr B(\kappa)\), and \(W\) is
	a normal ultrafilter on \(\kappa\). Then \(U\sE W\) if and only if \(U\mo
	W\).\qed
\end{cor}
Thus the Ketonen predecessors of a normal ultrafilter \(W\) on \(\kappa\) are
precisely \(\mathscr B(\kappa)\cap M_W\). We will see various nontrivial
generalizations of this fact to more general types of ultrafilters than normal
ones.

\subsection{The Rudin-Keisler order}\label{RKSection}
In this section, we briefly recall the theory of the Rudin-Keisler order and
explain its relationship with the Ketonen order. We also introduce the notion of
an {\it incompressible ultrafilter}, which will be a useful technical tool.

The Rudin-Keisler order is defined in terms of pushforward ultrafilters
(\cref{PushDef}).

\begin{defn}
Suppose \(U\) and \(W\) are ultrafilters. The {\it Rudin-Keisler order} is
defined by setting \(U\RK W\) if there is a function \(f : I \to X\) such that
\(f_*(W) = U\) where \(I\in W\) and \(X\) is the underlying set of
\(U\).\index{Rudin-Keisler order}
\end{defn}

We could of course take \(I\) to be the underlying set of \(W\) above. The
Rudin-Keisler order is a (nonstrict) preorder on the class of ultrafilters. For
us, the most important characterization of the Rudin-Keisler order uses
elementary embeddings:
\begin{lma}\label{RKChar}
Suppose \(U\) and \(W\) are ultrafilters. Then \(U\RK W\) if and only if there
is an elementary embedding \(k : M_U\to M_W\) such that \(k\circ j_U = j_W\).
\begin{proof}
Let \(X\) be the underlying set of \(U\).

First assume \(U\RK W\). Fix \(I\in W\) and \(f : I\to X\) such that \(f_*(W) =
U\). Let \(a = [f]_W\), so by \cref{PushDerived}, \(U\) is the ultrafilter on
\(X\) derived from \(j_W\) using \(a\). Let \(k : M_U\to M_W\) be the factor
embedding, so \(k(\id_U) = a\) and \(k\circ j_U = j_W\). Then \(k\) witnesses
the conclusion of the lemma.

Conversely, assume there is an elementary embedding \(k : M_U\to M_W\) such that
 \(k\circ j_U = j_W\). Let \(b = k(\id_U)\). Then \(b\in j_W(X)\). On the one
 hand, \(U\) is equal to the ultrafilter on \(X\) derived from \(j_W\) using
 \(b\). (Explicitly: \(U = j_U^{-1}[\pr {\id_U} {}] = j_U^{-1}[k^{-1}[k(\pr
 {\id_U} {})]] = j_W^{-1}[\pr b {}]\).) Fix \(I\in W\) and \(f : I\to X\) such
 that \([f]_W = b\). Then by \cref{PushDerived}, \(f_*(W)\) is the ultrafilter
 on \(X\) derived from \(j_W\) using \(b\), or in other words \(f_*(W) = U\).
 Thus \(U\RK W\) as desired.
\end{proof}
\end{lma}

A second combinatorial formulation of the Rudin-Keisler order is in terms of
partitions which will become relevant when we study indecomposability
(especially in \cref{Silver}):
\begin{lma}\label{RKPartitionChar}
Suppose \(U\) and \(W\) are ultrafilters. Let \(X\) be the underlying set of
\(U\). Then \(U\RK W\) if and only if there is a sequence of pairwise disjoint
sets \(\langle Y_x : x\in X\rangle\) such that \(U = \textstyle\left\{A\subseteq
X: \bigcup_{x\in A} Y_x\in W\right\}\).\qed
\end{lma}
The following is the fundamental theorem of the Rudin-Keisler order:
\begin{thm}\label{RKThm}
Suppose \(U\) and \(W\) are ultrafilters. Then  \(U\cong W\) if and only if
\(U\RK W\) and \(W\RK U\).
\end{thm}
We sketch the proof even though we do not need it in what follows. This involves
a very interesting rigidity theorem for pushforwards:
\begin{lma}\label{PushRigid}
Suppose \(U\) is an ultrafilter on \(X\) and \(f : X\to X\) is a function. If
\(f_*(U) = U\) then \(f(x) = x\) for \(U\)-almost all \(x\in X\).
\begin{proof}
Assume \(f : X\to X\) is such that \(f(x)\neq x\) for all \(x\in X\). We will
show that \(f_*(U)\neq U\). 

\begin{clm*} There is a partition \(X = A_0\cup A_1\cup A_2\) such that \(f[A_n] \subseteq X\setminus A_n\) for \(n = 0,1,2\). \end{clm*}
\begin{proof}[Sketch]Consider the directed graph \(G\) with vertices \(X\) and a directed edge from \(x\) to \(f(x)\) for each \(x\in X\). Our claim above amounts to the fact that \(G\) is 3-colorable. It suffices to show that each connected subgraph \(H\subseteq G\) is 3-colorable. Therefore suppose \(H\) is a connected subgraph of \(G\). The key point is that \(H\) contains at most one cycle (since \(G\) is a ``functional graph"), so one obtains an acyclic graph \(H'\) by removing an edge of \(H\) if necessary. Since \(H'\) is acyclic, \(H'\) is 2-colorable. By changing the color of at most one vertex in the coloring of \(H'\), one obtains a \(3\)-coloring of \(H\).\end{proof}
Since \(A_0\cup A_1\cup A_2 = X\in U\), either \(A_0,A_1,\) or \(A_2\) belongs
to \(U\). Assume without loss of generality that \(A_0\in U\). Then \(X\setminus
A_0\supseteq f[A_0] \in f_*(U)\), so \(X\setminus A_0\in f_*(U)\) as desired.
\end{proof}
\end{lma}

Let us reformulate this in terms of ultrapowers:
\begin{thm}\label{RKRigid}
	Suppose \(U\) is an ultrafilter and \(k : M_U\to M_U\) is an elementary
	embedding such that \(k\circ j_U = j_U\). Then \(k\) is the identity.
	\begin{proof}
		Let \(X\) be the underlying set of \(U\). Fix \(X\in U\) and a function
		\(f : X\to X\) such that \([f]_U = k(\id_U)\). Then by
		\cref{PushDerived}, \(f_*(U)\) is the ultrafilter on \(X\) derived from
		\(j_U\) using \(k(\id_U)\), which is easily seen to equal \(U\). (Yet
		another inverse image calculation: \(j_U^{-1}[\pr {k(\id_U)} {}] =
		(k\circ j_U)^{-1}[\pr {k(\id_U)} {}] = j_U^{-1}[k^{-1}[\pr {k(\id_U)}
		{}]] = j_U^{-1}[\pr {\id_U} {}] = U\).) Therefore by \cref{PushRigid},
		\(f\restriction I = \text{id}\) for some \(I\in U\). Thus \(k(\id_U) =
		[f]_U = \id_U\). It follows that \(k\restriction j_U[V]\cup \{\id_U\}\)
		is the identity, so \(k\restriction M_U\) is the identity since \(M_U =
		H^{M_U}(j_U[V]\cup \{\id_U\})\).
	\end{proof}
\end{thm}

\cref{PushRigid} immediately implies \cref{RKThm}:

\begin{proof}[Proof of \cref{RKThm}] Let \(X\) be the underlying set of \(U\)
and \(Y\) be the underlying set of \(W\). The trivial direction is to prove that
\(U\cong W\) implies \(U\RK W\) and \(W\RK U\). Fix \(I\in U\), \(J\in W\), and
a bijection \(f : I\to J\) such that for all \(A\subseteq I\), \(A\in U\) if and
only if \(f[A]\in W\). Viewing \(f\) as a function \(p : I\to Y\), we have \(W =
p_*(U)\). Viewing \(f^{-1}\) as a function \(p : J\to X\), we have \(U =
p_*(W)\). This implies implies \(W\RK U\) and \(U\RK W\).

Conversely assume \(U\RK W\) and \(W\RK U\). Fix \(I\in U\) and \(f : I\to Y\)
such that \(f_*(U) = W\). Fix \(J \in W\) and \(g: J\to X\) such that \(g_*(W)=
U\). We claim there is a set \(I'\subseteq I\) such that \(I'\in U\) and
\(g\circ f\restriction I'\) is the identity. To see this, note that \((g\circ
f)_*(U) = g_*(f_*(U)) = g_*(W) = U\). Therefore by \cref{PushRigid}, there is a
set \(I'\subseteq I\) such that \(I'\in U\) and \(g\circ f\) is the identity.
\end{proof}

\cref{RKThm} motivates the following definition:

\begin{defn}\index{Rudin-Keisler order!strict}
The {\it strict Rudin-Keisler order} is defined on ultrafilters \(U\) and \(W\)
by setting \(U\sRK W\) if \(U\RK W\) and \(W\not \cong U\).
\end{defn}

We now discuss the structure of the Rudin-Keisler order on countably complete
ultrafilters and its relationship to the Ketonen order. To facilitate this
discussion, we introduce a revised version of the Rudin-Keisler order. Recall
that a function \(f\) defined on a set of ordinals \(I\) is {\it regressive} if
\(f(\alpha) < \alpha\) for all \(\alpha\in I\).
\begin{defn}\label{rRKDef}\index{Rudin-Keisler order!revised}
Suppose \(U\) and \(W\) are ultrafilters on ordinals. Let \(X\) be the
underlying set of \(U\). The {\it revised Rudin-Keisler order} is defined by
setting \(U\rRK W\) if there is a set \(I\in W\) and a regressive function \(f :
I \to X\) such that \(f_*(W) = U\).
\end{defn}

\begin{lma}\label{rRKEmbedding}
If \(U\) and \(W\) are ultrafilters on ordinals, then \(U\rRK W\) if and only if
there is an elementary embedding \(k : M_U\to M_W\) such that \(k\circ j_U =
j_W\) and \(k(\id_U) < \id_W\).\qed
\end{lma}

\begin{cor}
The Ketonen order and the Rudin-Keisler order extend the revised Rudin-Keisler
order.\qed
\end{cor}

\begin{lma}\label{rRKLinear}
For any ultrafilter \(U\), the collection of tail uniform ultrafilters
isomorphic to \(U\) is linearly ordered by the revised Rudin-Keisler order.
\begin{proof}
Suppose \(W_0\cong U\cong W_1\) are tail uniform ultrafilters. Then \(W_0\rRK
W_1\) if and only if \(M_U\vDash \id_{W_0} < \id_{W_1}\).
\end{proof}
\end{lma}

We now introduce a concept that is very useful in the study of countably
complete ultrafilters. (The same concept was considered by Ketonen
\cite{Ketonen2}, who called them normalized ultrafilters.)

\begin{defn}
A tail uniform ultrafilter \(U\) on an ordinal \(\lambda\) is {\it
incompressible} if for any set \(I\in U\), no regressive function on \(I\) is
one-to-one.\index{Incompressible ultrafilter}
\end{defn}

\begin{lma}\label{IncomRK0}
Suppose \(U\) is tail uniform. The following are equivalent:
\begin{enumerate}[(1)]
\item \(U\) is incompressible.
\item If \(W \rRK U\), then \(W\sRK U\).\qed
\end{enumerate}
\end{lma}

\begin{lma}
A tail uniform ultrafilter \(U\) is incompressible if and only if it is the
\(\rRK\)-minimum element of \(C = \{U'\in \Un : U'\cong U\}\).
\begin{proof}
By \cref{IncomRK0}, \(U\) is an \(\rRK\)-minimal element of \(C\). Since
\(\rRK\) linearly orders \(C\) by \cref{rRKLinear}, \(U\) is the
\(\rRK\)-minimum element of \(C\).
\end{proof}
\end{lma}

\begin{cor}
An ultrafilter is isomorphic to at most one incompressible ultrafilter.\qed
\end{cor}

\begin{lma}\label{IncomPower}
Suppose \(U\) is tail uniform ultrafilter on \(\delta\). Then the following are
equivalent:
\begin{enumerate}[(1)]
\item \(U\) is incompressible.
\item \(\id_U\) is the least ordinal \(a\) of \(M_U\) such that \(M_U =
H^{M_U}(j_U[V]\cup \{a\})\).
\item \(\id_U\) is the largest ordinal \(a\) of \(M_U\) such that \(a\neq
j_U(f)(b)\) for any function \(f :\delta\to \delta\) and \(b < a\).\qed
\end{enumerate}
\end{lma}

If \(U\) is countably complete, then the collection of tail uniform ultrafilters
isomorphic to \(U\) is {\it wellordered} by \(\rRK\), and therefore it has a
minimum element. The following is the key existence theorem for incompressible
ultrafilters:
\begin{lma}
Any countably complete ultrafilter \(U\) is isomorphic to a unique
incompressible ultrafilter \(W\) which can be obtained in any of the following
ways:
\begin{itemize}
\item \(W\) is the \(\rRK\)-minimum element of the isomorphism class of \(U\).
\item \(W = f_*(U)\) where \(f : \delta_U\to\delta_U\) is the least one-to-one
function modulo \(U\).
\item \(W\) is the tail uniform ultrafilter derived from \(j_U\) using
\(\alpha\) where \(\alpha\) is the ordinal defined in either of the following
ways: 
\begin{itemize}
\item \(\alpha\) is least such that \(M_U = H^{M_U}(j_U[V]\cup \{\alpha\})\).
\item \(\alpha\) is largest such that \(\alpha\neq j_U(f)(\beta)\) for any
\(\beta < \alpha\).\qed
\end{itemize}
\end{itemize}
\end{lma} 

What makes incompressible ultrafilters useful is the following dual to
\cref{IncomRK0}:

\begin{prp}\label{IncomRK}
Suppose \(U\) is incompressible and \(W\) is an ultrafilter on an ordinal. If
\(U\sRK W\) then \(U\rRK W\).
\begin{proof}
Assume \(U\sRK W\). Fix \(k : M_U\to M_W\) such that \(k\circ j_U = j_W\). Since
\(U\not\cong W\), \(k\) is not an isomorphism. It follows that \(\id_W\notin
k[M_U]\): otherwise \(j_W[V]\cup\{\id_W\}\subseteq k[M_U]\) and so \(M_W =
H^{M_W}(j_W[V]\cup\{\id_W\})\subseteq k[M_U]\), and therefore \(k\) is
surjective and hence an isomorphism.

To show that \(U\rRK W\), it suffices by \cref{rRKEmbedding} to show that
\(k(\id_U) < \id_W\). Suppose not. Then \(\id_W\leq k(\id_U)\), and since
\(\id_W\notin k[M_U]\), in fact \(\id_W < k(\id_U)\). Since \(M_W =
H^{M_W}(j_W[V]\cup \{\id_W\})\) we can fix a function \(f :\delta_W\to
\delta_W\) such that \(j_W(f)(\id_W) = k(\id_U)\). Since \(\id_W < k(\id_U)\), 
\[M_W\vDash \exists \xi < k(\id_U)\ j_W(f)(\xi) = k(\id_U)\] Since \(j_W(f) =
k(j_U(f))\), the elementarity of \(k: M_U\to M_W\) implies
\[M_U\vDash \exists \xi < \id_U\ j_U(f)(\xi) = \id_U\] This contradicts
\cref{IncomPower} (3), which in particular states that \(\id_U \neq
j_U(f)(\xi)\) for any \(\xi < \id_U\).
\end{proof}
\end{prp}

\begin{cor}\label{RKrRK}
The strict Rudin-Keisler order and the revised Rudin-Keisler order coincide on
incompressible ultrafilters.\qed
\end{cor}

\begin{cor}\label{RKKet}\index{Ketonen order!vs. the Rudin-Keisler order}
The Ketonen order extends the strict Rudin-Keisler order on countably complete
incompressible ultrafilters.\qed
\end{cor}

We remark that given \cref{RKKet}, one might guess that \({\rRK} = {\RK}\cap
{\sE}\), but it is not hard to construct a counterexample under weak large
cardinal assumptions.

\begin{cor}[Solovay]\label{RKWF}
The strict Rudin-Keisler order is wellfounded on countably complete
ultrafilters.\index{Rudin-Keisler order!Wellfoundedness}
\begin{proof}
Suppose towards a contradiction that \[U_0\sgRK U_1\sgRK U_2\sgRK\cdots\] is a
descending sequence of countably complete ultrafilters in the strict
Rudin-Keisler order. For each \(n\), let \(W_n\) be the unique incompressible
ultrafilter isomorphic to \(U_n\). Then \[W_0\sgRK W_1 \sgRK  W_2\sgRK\cdots\]
since the strict Rudin-Keisler order is isomorphism invariant. But by
\cref{RKKet}, the Ketonen order extends the strict Rudin-Keisler order on
countably complete incompressible ultrafilters, and therefore \[W_0\sgE W_1 \sgE
W_2 \sgE\cdots\] This contradicts the wellfoundedness of the Ketonen order
(\cref{KOUnGood}).
\end{proof}
\end{cor}
Note that this yields another proof of \cref{RKThm} in the case that \(U\) and
\(W\) are countably complete.

\section{Variants of the Ketonen order}\label{VariantSection}
\subsection{Minimality of internal embeddings}\label{ExtendedKetonen}
In this subsection, we study an extension of the Ketonen order that provides
some insight into arbitrary extender embeddings (as opposed to just ultrapower
embeddings). This order also clarifies the connection between the Ketonen order
and the mouse order from inner model theory (which is also known as the
Dodd-Jensen order). Using these ideas, we prove a lemma (\cref{MinDefEmb}) that
states that if \(N\) and \(M\) are transitive models and \(j : N\to M\) is an
elementary embedding that is definable over \(N\) from parameters, then
\(j(\alpha) \leq k(\alpha)\) for any other elementary embedding \(k : N \to M\).

\begin{defn}\index{Pointed model}
A {\it pointed model} is a structure \((M,\xi)\) such that \(M\) is a transitive
model of ZFC and \(\xi\in \text{Ord}^{M}\). If \(\mathcal M = (M,\xi)\) is a
pointed model, then \(\xi_{\mathcal M} = \xi\). 
\end{defn}

We allow pointed models \((M,\xi)\) where \(M\) is a proper class. We abuse
notation by confusing a pointed model \(\mathcal M = (M,\xi)\) with its
underlying set \(M\). We therefore sometimes denote \(\xi\) by \(\xi_M\) instead
of \(\xi_\mathcal M\). 

When we discuss elementary embeddings of pointed models, we never impose
elementarity in the language of a pointed model (i.e., with a distinguished
constant for \(\xi\)), only elementarity in the language of set theory. 

\begin{defn}
Suppose \(N\) and \(M\) are transitive models of ZFC. An elementary embedding
\(j : N\to M\) is:
\begin{itemize}
\item an {\it extender embedding} if \(j\) is cofinal and \(M = H^M(j[M]\cup S)\) for some \(S\in M\).\index{Extender embedding}
\item an {\it internal extender embedding} if it is furthermore definable over
\(N\).
\end{itemize}
\end{defn}

\begin{defn}
Suppose \(M\) and \(N\) are pointed models.

The {\it Ketonen order on models}\index{Ketonen order!on models} is defined on
\(M\) and \(N\) by setting \(M\sE N\) if there are embeddings \((k,h) : (M,N)\to
P\) such that \(k(\xi_M) < h(\xi_N)\) and \(h\) is an internal extender
embedding of \(N\).

The {\it nonstrict Ketonen order} is defined on \(M\) and \(N\) by setting \(M\E
N\) if there are embeddings \((k,h) : (M,N)\to P\) such that \(k(\xi_M) \leq
h(\xi_N)\) and \(h\) is an internal extender embedding of \(N\).

{\it Ketonen equivalence} is defined on \(M\) and \(N\) by setting \(M\KE N\) if
\(M\E N\) and \(N\E M\).
\end{defn}

The transitivity for the Ketonen order on pointed models uses a trivial
``comparison lemma" that is provable in ZFC.

\begin{lma}\label{ZFCComparison}
Suppose \(M\), \(N_0\), and \(N_1\) are transitive models of \textnormal{ZFC}.
Suppose \(h: M\to N_0\) is an internal extender embedding and \(k : M\to N_1\)
is a cofinal elementary embedding. Then there is a comparison \((\ell, i) :
(N_0,N_1)\to P\) of \((h,k)\) such that \(i\) is an internal extender embedding
of \(N_1\).
\begin{proof}
Let \(\ell = k\restriction N_0\) and \(i = k(h)\). Then \(i\) is an internal
extender embedding and \(\ell\circ h = k\circ h = k(h)\circ k = i\circ k\), so
\((\ell,i)\) is a comparison of \((h,k)\).
\end{proof}
\end{lma}

\begin{lma}\label{PointedTransitive}
If \(M_0 \sE M_1\E M_2\), then \(M_0\sE M_2\).
\begin{proof}
Suppose \(M_0 \sE M_1\E M_2\). Let \((k_0,h_0) : (M_0,M_1)\to N_0\) witness
\(M_0\sE M_1\). Let \((k_1,h_1) : (M_1,M_2)\to N_1\) witness \(M_1\E M_2\).
Applying \cref{ZFCComparison}, let \[(\ell,i) : (N_0,N_1) \to P\] be a
comparison of \((h_0,k_1)\) such that \(i\) is an internal extender embedding of
\(N_1\). Let \(k = \ell \circ k_0\) and let \(h = i\circ h_1\). Then \((k,h) :
(M_0,M_2)\to P\) and \(h\) is an internal extender embedding of \(M_2\) since it
is the composition of the internal extender embeddings \(i\) and \(h_1\).
Finally, \[k(\xi_{M_0}) = \ell\circ k_0(\xi_{M_0}) < \ell\circ h_0(\xi_{M_1}) =
i\circ k_1(\xi_{M_1}) \leq i\circ h_1(\xi_{M_2}) = h(\xi_{M_2})\qedhere\]
\end{proof}
\end{lma}

We will prove the wellfoundedness of the Ketonen order on pointed models that
satisfy a very weak form of iterability.

\begin{defn}
A transitive model \(M\) of ZFC is {\it \(\omega\)-linearly iterable} if the
following holds. Suppose \[M = M_0\stackrel{h_0}{\longrightarrow}
M_1\stackrel{h_1}{\longrightarrow} M_2 \stackrel{h_2}{\longrightarrow}\cdots\]
is such that for all \(i < \omega\), \(h_i : M_i\to M_{i+1}\) is an internal
extender embedding, then its direct limit is wellfounded.
\end{defn}

The following is a well-known fact, versions of which are due to Gaifman, Kunen,
and Mitchell (see \cite{MitchellHandbook}):

\begin{lma}\label{Iterability}
Suppose \(M\) is a model of \textnormal{ZFC} such that \(\textnormal{Ord}^M\)
has uncountable cofinality. Then \(M\) is \(\omega\)-linearly iterable.
Similarly, any inner model is \(\omega\)-linearly iterable.\qed
\end{lma}

The proof of the wellfoundedness of the Ketonen order on pointed models is based
on the proof of the wellfoundedness of the Dodd-Jensen order.
\begin{thm}\label{GenWellfounded}
The Ketonen order is wellfounded on \(\omega\)-linearly iterable pointed
models.\index{Ketonen order!on models!wellfoundedness}
\begin{proof}
To simplify notation, we isolate the main step of the proof as a lemma:
\begin{lma}\label{WF1Step}
Suppose that \(M_0\sgE M_1\sgE M_2\sgE \cdots\) is a descending sequence of
pointed models. Then there is a descending sequence \(N_0\sgE  N_1\sgE N_2\sgE
\cdots\) of pointed models and an internal extender embedding \(h : M_0\to N_0\)
with \(\xi_{N_0} < h(\xi_{M_0})\).
\begin{proof}
The proof is illustrated by \cref{GenWFFig}. 
\begin{figure}
	\center
	\includegraphics[scale=.75]{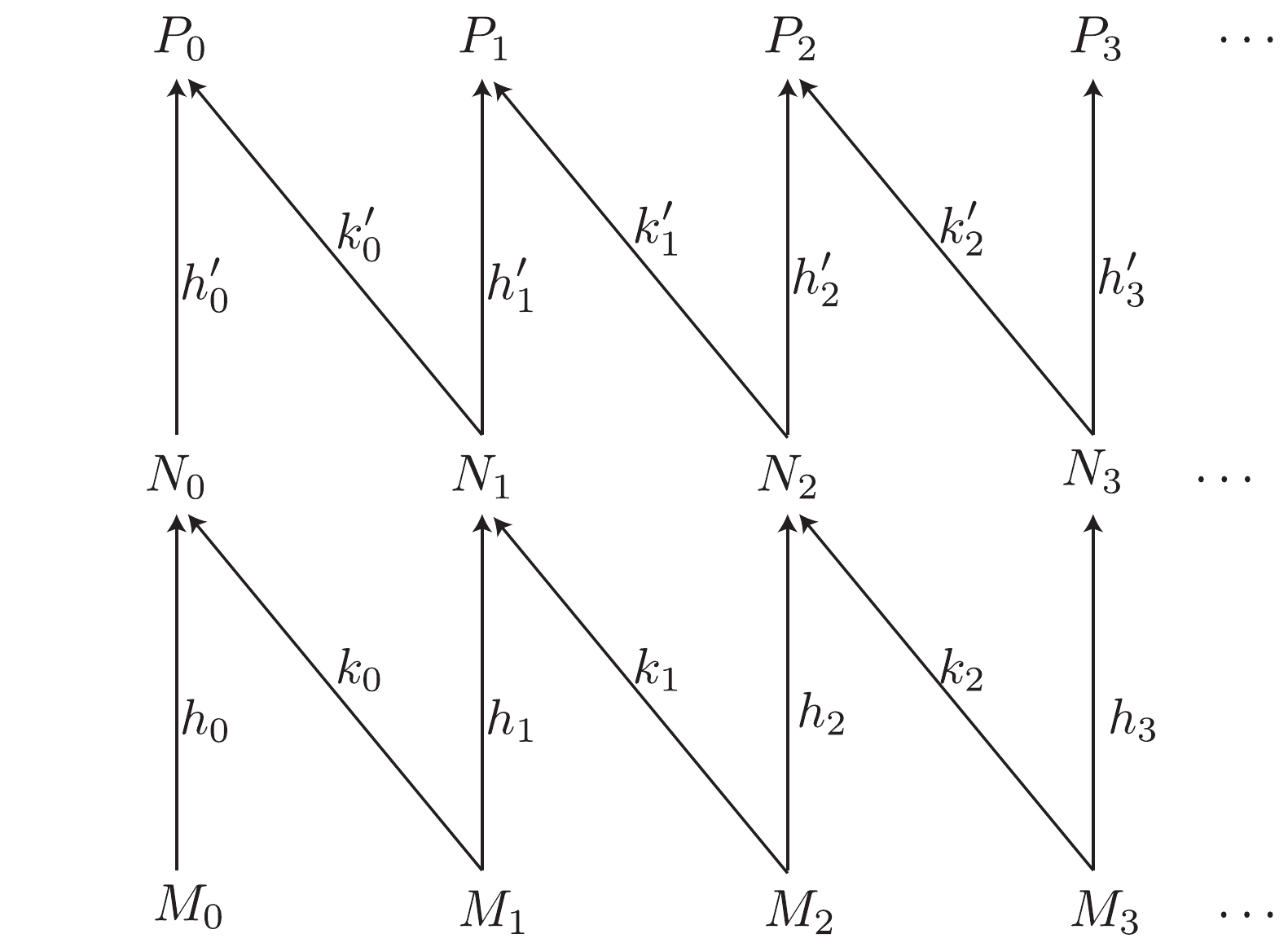}
	\caption{The proof of \cref{WF1Step}}\label{GenWFFig}
\end{figure} Let \((h_i,k_i) : (M_i,M_{i+1})\to N_i\) witness \(M_i\sgE
M_{i+1}\). We endow \(N_i\) with the structure of a pointed model by letting
\(\xi_{N_i} = k_i(\xi_{M_{i+1}})\).

Setting \(h = h_0\), particular, \(h\) is an internal extender embedding and
\(\xi_{N_0} < h(\xi_{M_0})\). It remains to verify that \(N_0\sgE  N_1\sgE
N_2\sgE \cdots\). Fix \(i < \omega\). By \cref{ZFCComparison}, there is a
comparison \((h_i',k_i') : (N_i,N_{i+1})\to P_i\) of \((k_i,h_{i+1})\) such that
\(h_i'\) is an internal extender embedding of \(N_i\). As in
\cref{PointedTransitive}, \[h_i'(\xi_{N_i}) = h_i'(k_i(\xi_{M_{i+1}})) =
k_i'(h_{i+1}(\xi_{M_{i+1}})) > k_i'(k_{i+1}(\xi_{M_{i+1}})) =
k_i'(\xi_{N_{i+1}})\] and hence \(N_i\sgE N_{i+1}\).
\end{proof}
\end{lma}
Now suppose towards a contradiction that  \(M^0_0\sgE M^0_1\sgE M^0_2\sgE
\cdots\) is a descending sequence of \(\omega\)-linearly iterable pointed
models. By recursion, using the lemma, one obtains sequences \(M^i_0\sgE
M^i_1\sgE M^i_2\sgE \cdots\) and internal extender embeddings \(h^i : M^i_0 \to
M^{i+1}_0\) with \(\xi_{M^{i+1}_0} < h^i(\xi_{M^i_0})\) for all \(i < \omega\).
But then the iteration \[M^0_0\stackrel{h^0}{\longrightarrow}
M_0^1\stackrel{h^1}{\longrightarrow}M_0^2\stackrel{h^2}{\longrightarrow}\cdots\]
has an illfounded direct limit, which contradicts that \(M^0_0\) is
\(\omega\)-linearly iterable.
\end{proof}
\end{thm}

The wellfoundedness of the Ketonen order on pointed models has some useful
consequences. Of course, it provides an alternate proof of the wellfoundedness
of the Ketonen order:

\begin{proof}[Alternate Proof of \cref{KOWellfounded}] For \(U\in \Un\), let
\(\Phi(U) = (M_U,\id_U)\). Then for any \(U\in \Un\), \(\Phi(U)\) is an
\(\omega\)-linearly iterable pointed model. Moreover, if \(U\sE W\), then
\(\Phi(U)\sE \Phi(W)\) since internal ultrapower embeddings are internal
extender embeddings. Thus the Ketonen order is wellfounded on \(\Un\) since by
\cref{GenWellfounded}, the Ketonen order is wellfounded on \(\omega\)-linearly
iterable models.
\end{proof}

More interestingly, \cref{GenWellfounded} implies a coarse version of the
Dodd-Jensen Lemma (proved for example, \cite{Steel}):
\begin{thm}\label{MinDefEmb}\index{Minimality of definable embeddings}
Suppose \(M\) is an \(\omega\)-linearly iterable model. Suppose \(h,k : M\to N\)
are elementary embeddings and \(h\) is an internal extender embedding of \(M\).
Then for all \(\alpha\in \textnormal{Ord}^M\), \(h(\alpha) \leq k(\alpha)\).
\begin{proof}
Suppose towards a contradiction that \(k(\alpha) < h(\alpha)\). Then \((k,h) :
(M,M)\to N\) witnesses \((M,\alpha) \sE (M,\alpha)\), contradicting
\cref{GenWellfounded}.
\end{proof}
\end{thm}
The idea of generalizing arguments from inner model theory to prove results like
\cref{MinDefEmb} is due to Woodin \cite{Woodin}, who proved the similar theorem
that if \(M\) and \(N\) are models of ZFC and \(M\) is finitely generated, then
there is at most one close embedding from \(M\) to \(N\). Woodin's theorem
actually follows from the restriction of \cref{MinDefEmb} to ultrapower
embeddings.

By tracing through the proof of this theorem, one can prove the following fact,
which is really a theorem scheme:
\begin{thm}
Suppose \(M\) and \(N\) are inner models, \(h,k : M\to N\) are elementary
embeddings. If \(h\) is definable over \(M\), then \(h(\alpha) \leq k(\alpha)\)
for all ordinals \(\alpha\).\qed
\end{thm}
There are some metamathematical difficulties involving the linear iterability of
an inner model \(M\) by an \(\omega\)-sequence of definable embeddings: it is
not in general clear that this is first-order expressible in the language of set
theory with a predicate for \(M\). The iterability required for the proof of
\cref{MinDefEmb}, however, can be stated and proved. We omit the proof since we
have no applications of this more general theorem.

\subsection{The seed order}\label{SOSection}
We now define the seed order, a variant of the Ketonen order that uses fully
internal ultrapower comparisons.
\begin{defn}
Suppose \(U\) and \(W\) are countably complete ultrafilters on ordinals.
\index{Seed order}

The {\it seed order} is defined by setting \(U\swo W\) if there is an internal
ultrapower comparison \((k,h)\) of \((j_{U},j_{W})\) such that \(k(\id_{U}) <
h(\id_{W})\).

The {\it nonstrict seed order} is defined by setting \(U\wo W\) if there is an
internal ultrapower comparison \((k,h)\) of \((j_{U},j_{W})\) such that
\(k(\id_{U}) \leq h(\id_{W})\).

{\it Seed equivalence} is defined by setting \(U=_S W\) if there is an internal
ultrapower comparison \((k,h)\) of \((j_{U},j_{W})\) such that \(k(\id_{U}) =
h(\id_{W})\).
\end{defn}

\begin{lma}
If \(U\) and \(W\) are countably complete ultrafilters, then \(U =_S W\) if and
only if \(U \KE W\).\qed
\end{lma}

\begin{lma}
Suppose \(U_0\) and \(U_1\) are countably complete ultrafilters concentrating on
ordinals. Then \(U_0\wo U_1\) if and only if \(U_0 \swo U_1\) or \(U_0 =_S
U_1\).\qed
\end{lma}

By the characterization of the Ketonen order in terms of comparisons
(\cref{KOChars}) we have the following fact:
\begin{lma}
The Ketonen order extends the seed order.\qed
\end{lma}

It follows that the seed order is a strict wellfounded set-like relation.
(Transitivity is another story; see \cref{TransUA} below.) 

\begin{cor}\label{SeedAntisymmetry}
Suppose \(U\) and \(W\) are countably complete ultrafilters on ordinals. Then
\(U =_S W\) if and only if \(U\wo W\) and \(W\wo U\).\qed
\end{cor}

\begin{prp}[UA]\label{SeedKet}
The seed order linearly orders \(\Un\).
\begin{proof}
By the definition of the Ultrapower Axiom, the nonstrict seed order is a total
relation on \(\Un\). By \cref{SeedAntisymmetry} and the fact that \(=_S\)
restricts to equality on \(\Un\), the seed order is antisymmetric on \(\Un\).
Thus the seed order linearly orders \(\Un\).
\end{proof}
\end{prp}

The seed order, unlike the Ketonen order, is not provably transitive in ZFC (for
mundane reasons):
\begin{prp}\label{TransUA}
The seed order is transitive if and only if the Ultrapower Axiom
holds.\index{Seed order!transitivity}
\end{prp}

The proof uses the following trivial variant of \cref{SpaceLemma}.

\begin{lma}\label{PrincipalOrder2}
Suppose \(\alpha\) is an ordinal and \(U\) is a countably complete ultrafilter
that concentrates on ordinals. Then \(U\) and the principal ultrafilter \(\pr
\alpha {}\) are comparable in the seed order:
\begin{itemize}
\item \(U \swo \pr \alpha {}\) if and only if \(\delta_U\leq \alpha\).
\item \(U \KE \pr \alpha {}\) if and only if \(\delta_U = \alpha + 1\).
\item \(U\slwo \pr \alpha {}\) if and only if \(\alpha +1 < \delta_U\).\qed
\end{itemize}
\end{lma}

\begin{proof}[Proof of \cref{TransUA}] Suppose \(j_0 : V\to M_0\) and \(j_1 :
V\to M_1\) are ultrapower embeddings. We will show they can be compared. For \(i
= 0,1\), fix ordinals \(\alpha_i\in M_i\) such that \(M_i = H^{M_i}(j_i[V]\cup
\{\alpha_i\})\) with the further property that letting \(U_i\) be the tail
uniform ultrafilter derived from \(j_i\) using \(\alpha_i\),  \(\delta_{U_0} <
\delta_{U_1}\).

By \cref{PrincipalOrder2}, \[U_0 \swo \pr {\delta_{U_0}} {}\wo U_1\] Thus if the
seed order is transitive, \(U_0 \wo U_1\). Since \(j_0 = j_{U_0}\) and \(j_1 =
j_{U_1}\) the fact that \(U_0 \swo U_1\) implies in particular that there is an
internal ultrapower comparison of \((j_0,j_1)\). This verifies the Ultrapower
Axiom for the pair \((j_0,j_1)\).
\end{proof}

We now consider the seed order on pointed models.

\begin{defn}\index{Pointed ultrapower}\index{Pointed ultrapower embedding}
A {\it pointed ultrapower} is a pointed model \(\mathcal M\) whose underlying
class \(M\) is an ultrapower of the universe \(V\). A {\it pointed ultrapower
embedding} is a pair \((j,\xi)\) where \(j\) is an ultrapower embedding and
\(\xi\) is an ordinal.
\end{defn}

There is a natural identification of countably complete ultrafilters with a
certain class of pointed ultrapower embeddings:
\begin{defn}\index{Pointed ultrapower embedding!representing an ultrafilter}
	Suppose \(U\) is a countably complete ultrafilter on an ordinal. Then the
	{\it pointed ultrapower embedding representing \(U\)} is \((j_U,\id_U)\). A
	pointed ultrapower embedding \((j,\xi)\) {\it represents an ultrafilter} if
	it is the pointed ultrapower embedding representing some ultrafilter.
\end{defn}

We apologize for bombarding the reader with definitions. The following
definitions extend the seed order and Ketonen order to pointed ultrapowers and
embeddings.

\begin{defn}
Suppose \(M\) and \(N\) are pointed ultrapowers.

The {\it seed order}\index{Seed order!on pointed ultrapowers} ({\it nonstrict
seed order}) is defined by setting \( M \swo  N\) (\( M \wo  N\)) if there is an
internal ultrapower comparison \((k,h) : ( M, N)\to P\) such that \(k(\xi_ M) <
h(\xi_{ N})\) (\(k(\xi_ M) \leq h(\xi_{ N})\)). 

{\it Seed equivalence}\index{Seed order!seed equivalence}\index{\(=_S\) (seed
equivalence)} is defined by setting \( M \wo  N\) if there is an internal
ultrapower comparison \((k,h) : ( M, N)\to P\) such that \(k(\xi_ M) = h(\xi_{
N})\). 
\end{defn}

\begin{defn}
Suppose \((i,\nu)\) and \((j,\xi)\) are pointed embeddings.

The {\it seed order} ({\it nonstrict seed order}) is defined by setting
\((i,\nu) \swo (j,\xi)\) (\((i,\nu) \wo (j,\xi)\)) if there is an internal
ultrapower comparison \((k,h)\) of \((i,j)\) such that \(k(\nu) < h(\xi)\)
(\(k(\nu) \leq h(\xi)\)). 

{\it Seed equivalence} is defined on by setting \((i,\nu) =_S (j,\xi)\) if there
is an internal ultrapower comparison \((k,h)\) of \((i,j)\) such that \(k(\nu) =
h(\xi)\). 

The {\it Ketonen order} ({\it nonstrict Ketonen order}) is defined on by setting
\((i,\nu) \sE (j,\xi)\) (\((i,\nu)  \E (j,\xi)\)) if there is a comparison
\((k,h)\) of \((i,j)\) such that \(h\) is an internal ultrapower embedding and
\(k(\nu) < h(\xi)\) (\(k(\nu) \leq h(\xi)\)).

{\it Ketonen equivalence}\index{Ketonen order!Ketonen
equivalence}\index{\(=_{\Bbbk}\) (Ketonen equivalence)} is defined by setting
\((i,\nu)=_{\Bbbk} (j,\xi)\) if \((i,\nu)\E(j,\xi)\) and \((j,\xi)\E (i,\nu)\).
\end{defn}

The following is in a sense the strongest consequence of UA for these orders:

\begin{prp}[UA]\label{GenTrichotomy}
Suppose \((i,\nu)\) and \((j,\xi)\)  are pointed ultrapower embeddings. Then
either \((i,\nu) \swo (j,\xi)\), \((i,\nu) =_S (j,\xi)\), or \((i,\nu) \slwo
(j,\xi)\).\qed
\end{prp}

\begin{defn}
If \((j,\xi)\) is a pointed ultrapower embedding, then \(M(j)\) denotes the
target model of \(j\) and \(M(j,\xi)\) denotes the pointed ultrapower
\((M(j),\xi)\).
\end{defn}

\begin{lma}[UA]\label{GenSE}\index{Ketonen order!vs. the seed order}\index{Seed order!vs. the Ketonen order}
Suppose  \((i,\nu)\) and \((j,\xi)\)  are pointed ultrapower embeddings. Then
the following are equivalent:
\begin{enumerate}[(1)]
\item \((i,\nu) \wo (j,\xi)\).
\item \((i,\nu)\E (j,\xi)\).
\item \(M(i,\nu) \wo M(j,\xi)\).
\item \(M(i,\nu) \E M(j,\xi)\). 
\end{enumerate}
\begin{proof}
The implications from (1) to (2) to (3) to (4) are trivial, so to prove the
lemma, it suffices to show that (4) implies (1). Therefore assume \(M(i,\nu)\E
M(j,\xi)\). Assume (4) fails, towards a contradiction, so that by
\cref{GenTrichotomy}, \((j,\xi)\swo (i,\nu)\). Therefore \(M(j,\xi)\swo M
(i,\nu)\) and hence \(M(j,\xi)\sE M(i,\nu)\). By \cref{PointedTransitive},
\(M(i,\nu)\sE M(i,\nu)\), contradicting \cref{GenWellfounded}.
\end{proof}
\end{lma}

Unlike their restrictions to ultrafilters, the relations \(\KE\) and \(=_S\) are
far from trivial on pointed ultrapowers and embeddings. When one of the pointed
ultrapower embeddings involved represents an ultrafilter, \(=_S\) is closely
related to the Rudin-Frol\'ik order (\cref{RFChapter}):

\begin{lma}\label{RFSEquiv}\index{Seed order!seed equivalence!vs. the Rudin-Frol\'ik order}
Suppose \((i,\nu)\) and \((j,\xi)\) are pointed ultrapower embeddings such that
\((i,\nu)\) represents an ultrafilter. Then \((i,\nu) =_S (j,\xi)\) if and only
if there is an internal ultrapower embedding \(e : M(i)\to M(j)\) such that
\(e\circ i = j\) and \(e(\nu) = \xi\).
\begin{proof}
We prove the forwards direction, since the converse is trivial. Let \(M = M(i)\)
and \(N = M(j)\). Fix an internal ultrapower comparison \((k,h) : (M,N)\to P\)
of \((i,j)\) with \(k(\nu) = h(\xi)\). 

We claim that \(k[M]\subseteq h[N]\). Since \(i\) represents an ultrafilter, \(M
= H^M(i[V]\cup \{\nu\}),\) and hence \(k[M] = H^P(k[i[V]\cup \{\nu\}])\). But
\(k(\nu) = h(\xi)\in h[N]\) and \(k[i[V]] = h[j[V]]\subseteq h[N]\). Thus
\(k[i[V]\cup \{\nu\}]\subseteq h[N]\), so that \(k[M] = H^P(k\circ i[V]\cup
\{k(\nu)\})\subseteq h[N]\), as claimed.

Let \(e = h^{-1}\circ k\). Then \(e : M\to N\) is an elementary embedding and
since \(e\circ i = h^{-1} \circ k \circ i = h^{-1} \circ h\circ j\), we have
\(e\circ i = j\). It follows that \(e\) is an ultrapower embedding of \(M\).
Since \(h \circ e = k\) and \(k\) is close to \(M\), \(e\) is an internal
ultrapower embedding of \(M\) by \cref{CloseLemma}. Finally, \(e(\nu) =
h^{-1}(k(\nu)) = h^{-1}(h(\xi)) = \xi\). 
\end{proof}
\end{lma}

\begin{cor}[UA]\label{UniqueUltrapowerEmbedding}\index{Ultrapower!ultrapower embedding!uniqueness}
Suppose \(i : V\to M\) and \(j : V\to M\) are ultrapower embeddings with the
same target model. Then \(i = j\).
\begin{proof}
Fix \(\xi\) such that \(M = H^M(i[V]\cup \{\xi\})\). Since \(M(i,\xi) = (M,\xi)
= M(j,\xi)\), we must have \((i,\xi) =_S (j,\xi)\) by \cref{GenSE}. By
\cref{RFSEquiv}, it follows that there is an internal ultrapower embedding \(k:
M\to M\) such that \(k\circ i = j\) and \(k(\xi) = \xi\). Since \(k\) is
internal to \(M\), \(k\) is the identity, and therefore \(i = j\).
\end{proof}
\end{cor}

The structure of the equivalence relation \(=_S\) on pointed models under UA
seems quite interesting. For example, for all we know, if \(M =_S N\) then there
some \(H\) of which both \(M\) and \(N\) are ultrapowers such that \(M =_S H =_S
N\).

\subsection{The width of an embedding}
As a brief digression, we make some general remarks about the size of
ultrafilters necessary to realize compositions of ultrapower embeddings. It
turns out to be easier to work in a bit more generality, using a definition due
to Cummings \cite{CummingsHandbook}:

\begin{defn}\index{Width of an elementary embedding}
	Suppose \(M\) and \(N\) are transitive models of ZFC and \(j : M\to N\) is
	an extender embedding. The {\it width} of \(j\), denoted
	\(\textsc{width}(j)\), is the least \(M\)-ordinal \(\iota\) such that \(N =
	H^N(j[M]\cup \sup j[\iota])\).
\end{defn}
Note that an embedding is an extender embedding if and only if its width is
well-defined. For ultrapower embeddings, there is a simple relationship between
width and size. Recall from \cref{RelativizedLambda}, which generalizes the
notion of size (i.e., \(\lambda_U\)) to \(M\)-ultrafilters \(U\).

\begin{prp}
	If \(j : M\to N\) is the ultrapower embedding associated to an
	\(M\)-ultrafilter \(U\), then \(\textsc{width}(j) = \lambda_U+1\).\qed
\end{prp}

There are really two key facts about width, both of which are generalized by the
theory of generators (\cref{CompositionGenerators}). The first can be summarized
that narrow embeddings are continuous\index{Continuous embedding} at large
regular cardinals \(\lambda\) (i.e., \(j(\lambda) =\sup j[\lambda]\)):

\begin{lma}\label{NarrowContinuity}
	Suppose \(j: M\to N\) is an extender embedding and \(\lambda\) is an ordinal
	of \(M\)-cofinality at least \(\textsc{width}(j)\). Then \(j(\lambda) = \sup
	j[\lambda]\).
	\begin{proof}
		Suppose \(\alpha\in \text{Ord}^{N}\) and \(\alpha < j(\lambda)\). We
		will show \(\alpha \leq j(\gamma)\) for some \(\gamma < \lambda\). Since
		\(N = H^N(j[M]\cup \sup j[\lambda])\), we can find a function \(f\in M\)
		and an ordinal \(\nu < \lambda\) such that \(\alpha = j(f)(\xi)\) for
		some \(\xi < j(\nu)\). Since the \(M\)-cofinality of \(\lambda\) is
		above \(\nu\), \(f[\nu]\cap \lambda\) is bounded by some \(\gamma <
		\lambda\). Hence \(\alpha = j(f)(\xi) \leq \sup j(f)[j(\nu)]\cap
		j(\lambda) = j(\sup f[\nu]\cap \lambda) = j(\gamma)\), as desired.
	\end{proof}
\end{lma}

This has a useful consequence for ultrapower embeddings (which is essentially
equivalent):

\begin{lma}\label{MUFContinuity}\label{UFContinuity}
	Suppose \(M\) is a transitive model of \textnormal{ZFC} and \(U\) is an
	\(M\)-ultrafilter on a set \(X\in M\). Then for any ordinal \(\delta\) such
	that \(\text{cf}^M(\delta) > \lambda_U\), \(j_U^M\) is continuous at
	\(\delta\).\qed
\end{lma}

It is worth mentioning a related fact here:
\begin{lma}\label{UFCardinality}
	If \(U\) is an ultrafilter on \(X\), then for any cardinal \(\gamma\),
	\(|j_U(\gamma)| \leq \gamma^{|X|}\). Thus if \(\lambda\) is a strong limit
	cardinal above \(|X|\), \(j_U[\lambda]\subseteq \lambda\). If moreover
	\(\textnormal{cf}(\lambda) > |X|\), then \(j_U(\lambda) = \lambda\).\qed
\end{lma}

The second provides a computation of the width of a composition in terms of the
width of the factors:

\begin{lma}\label{WidthLemma}
	Suppose \(M\stackrel{i}{\longrightarrow}N\stackrel{j}{\longrightarrow}P\)
	are elementary embeddings. Then 
	\[\textsc{width}(j\circ i) = \max \{\textsc{width}(i),\gamma\}\] where
	\(\gamma\) is the least ordinal such that \(\textsc{width}(j) \leq \sup
	i[\gamma]\).
	\begin{proof}
		Let \(\iota = \max \{\textsc{width}(i),\gamma\}\).
		
		We first show \(\textsc{width}(j\circ i) \leq \iota\). Since
		\(\textsc{width}(i) \leq \iota\), \(N = H^N(i[M]\cup \sup i[\iota])\).
		Since \(\textsc{width}(j) \leq \sup i[\iota]\), \(P = H^P(j[N]\cup \sup
		j[\sup i[\iota]]) = H^P(j[N]\cup \sup j\circ i[\iota])\). Putting these
		calculations together, \[P = H^P(j[i[M]\cup \sup i[\iota]]\cup \sup
		j\circ i[\iota]) = H^{P}(j\circ i[M]\cup \sup j\circ i[\iota])\] It
		follows that \(\textsc{width}(j\circ i) \leq \iota\).
		
		We now show \(\iota \leq \textsc{width}(j\circ i)\). First, we show
		\(\gamma \leq \textsc{width}(j\circ i)\). Fix \(\eta < \gamma\), and we
		will show \(\eta < \textsc{width}(j\circ i)\). This follows from the
		fact that \[H^P(j\circ i[M]\cup \sup j\circ i[\eta]) \subseteq
		H^P(j[N]\cup \sup j[\sup i[\eta]]) \subsetneq P\] The final inequality
		uses that \(\sup i[\eta] < \textsc{width}(j)\).
		
		We finish by showing \(\textsc{width}(i) \leq\textsc{width}(j\circ i)\).
		This uses the argument from \cref{IncomRK}. Suppose \(\eta <
		\textsc{width}(i)\) is an \(M\)-cardinal, and we will show  \(\eta <
		\textsc{width}(j\circ i)\). Fix \(a\in N\) such that \(a\notin
		H^N(i[M]\cup \sup i[\eta])\). Suppose towards a contradiction that
		\(j(a)\in H^P(j\circ i[M]\cup \sup j\circ i[\eta])\). Fix \(\xi < \eta\)
		and \(f\in M\) such that \(j(a) = j(i(f))(\alpha)\) for some \(\alpha
		\leq j(i(\xi))\). Then by the elementarity of \(j\), \(N\) satisfies
		that \(a = i(f)(\alpha)\) for some \(\alpha \leq i(\xi)\). This
		contradicts our assumption that \(a\notin H^N(i[M]\cup \sup i[\eta])\).
		Therefore \(j(a)\notin H^P(j\circ i[M]\cup \sup j\circ i[\eta])\), so
		\(\eta < \textsc{width}(j\circ i)\), as desired.
	\end{proof}
\end{lma}

\subsection{The direct limit of all ultrapowers}\label{MInftySection}
Under the Ultrapower Axiom, it is possible to take the direct limit of all
ultrapower embeddings. The properties of this structure, which is denoted
\(M_\infty\), turn out to be closely related both to the Ketonen order on
pointed ultrapowers and to the theory of supercompact cardinals. 

To save ink, it is convenient to let \(\infty\) be a formal symbol such that (by
definition) every ordinal \(\lambda\) satisfies \(\lambda < \infty\).

\begin{defn}
	Suppose \(M\) is an ultrapower of the universe and \(\lambda\) is an
	ordinal. Then \(\nu(M,\lambda) = \sup i[\lambda]\) where \(i : V\to M\) is
	any ultrapower embedding. 
\end{defn}

By \cref{MinDefEmb}, \(\nu(M,\lambda)\) does not depend on the choice of \(i\).
We also set \(\nu(M,\infty) = \infty\).

\begin{defn}
If \(\lambda\) is is a cardinal or \(\lambda = \infty\), then \(\mathcal
D_\lambda\) denotes the following category:
\begin{itemize}
	\item An inner model \(M\) is an object of \(\mathcal D_\lambda\) if there is an ultrapower embedding  \(i : V\to M\) with \(\textsc{width}(i) \leq \lambda\).\index{Ultrapower!category of}
	\item If \(M,N\in \mathcal D_\lambda\), an internal ultrapower embedding \(j
	: M\to N\) is a morphism of \(\mathcal D_\lambda\) if \(\textsc{width}(j)
	\leq \nu(M,\lambda)\).
\end{itemize}
\end{defn}

Thus \(\mathcal D_\infty\) is the category of all ultrapowers of the universe
equipped with their internal ultrapower embeddings.  As an immediate consequence
of \cref{UniqueUltrapowerEmbedding}, UA implies that \(\mathcal D_\lambda\) is a
full subcategory of \(\mathcal D_\infty\); that is, it contains all morphisms
between the objects it sees. It is not clear whether this is the case in ZFC
(and it is not even clear whether this subcategory is locally small in the
natural sense).\footnote{This raises an interesting question in the general
theory of elementary embeddings:
\begin{qst}[ZFC] Can two ultrapower embeddings of the universe have the same
	target model but different widths?
\end{qst}
The question has something to do with uniform ultrafilters on singular
cardinals.}

\begin{defn}
	A category \(\mathcal C\) is a {\it partial order} if every pair of objects
	\(a,b\in \mathcal C\) there is at most one morphism from \(a\) to \(b\) in
	\(\mathcal C\). A category \(\mathcal C\) is {\it directed} if for every
	pair of objects \(a,b\in \mathcal C\), there is a further object \(c\in
	\mathcal C\) admitting morphisms \(a\to c\) and \(b \to c\).
\end{defn}

We have the following equivalences:

\begin{lma}
The following are equivalent:
\begin{enumerate}[(1)]
\item The Ultrapower Axiom.
\item For all cardinals \(\lambda\), \(\mathcal D_\lambda\) is a directed
partial order.
\item \(\mathcal D_\infty\) is a directed partial order.
\end{enumerate}
\begin{proof} 	
 	{\it (1) implies (2):} The fact that \(\mathcal D_\lambda\) is a partial
 	order follows immediately from \cref{UniqueUltrapowerEmbedding}. The
 	directedness of \(\mathcal D_\lambda\) follows from an easy localization of
 	UA (\cref{UpperBoundBound}), which states that if \(U\) and \(W\) are
 	countably complete ultrafilters on a cardinal \(\gamma\), then there is a
 	countably complete ultrafilter \(Z\) on \(\gamma\) such that there are
 	internal ultrapower embeddings \(k : M_U\to M_Z\) and \(h : M_W\to M_Z\).
 	
 	{\it (2) implies (3):} Immediate.
 	
  	{\it (3) implies (1):} Immediate.
\end{proof}

\end{lma}

\begin{defn}[UA]\index{Direct limit of all ultrapowers}
If \(\lambda\) is a cardinal or \(\lambda = \infty\), let \[M_\lambda = \lim
\mathcal D_\lambda\] For all \(M\in \mathcal D_\lambda\), \[j_{M,\lambda} : M\to
M_\lambda\] denotes the direct limit embedding.
\end{defn}

The models \(M_\lambda\) are wellfounded by a standard application of the linear
iterability of the universe (\cref{Iterability}). We will see that \(M_\infty\)
need not be set-like. By convention we identify its set-like part with an inner
model. 

The following lemma is the key to the analysis of the models \(M_\lambda\) for
\(\lambda\) a regular cardinal:

\begin{lma}[UA]\label{InftyAbs}
Suppose \(\lambda\) is a regular cardinal or \(\lambda = \infty\). For any
ultrapower embedding \(i :V\to N\) with \(i\in \mathcal D_\lambda\),
\(i(M_\lambda) = M_\lambda\) and \(i(j_{V,\lambda}) = j_{N,\lambda}\).
\begin{proof}
The key point is that since \(\lambda\) is regular and \(\textsc{width}(i) \leq
 \lambda\), \(i(\lambda)  = \sup i[\lambda] = \nu(N,\lambda)\). Thus
 \(i(\mathcal D_\lambda) = \mathcal D_{i(\lambda)}^N =\mathcal
 D_{\nu(N,\lambda)}^N\), which is equal to the cone above \(N\) in \(\mathcal
 D_\lambda\). Since \(\mathcal D_\lambda\) is a directed partial order, this
 cone is cofinal in \(\mathcal D_\lambda\), and thus its direct limit is equal
 to that of \(\mathcal D_\lambda\). In other words, \(i(M_\lambda) = \lim
 i(\mathcal D_\lambda) = \lim \mathcal D_\lambda = M_\lambda\), and similarly
 \(i(j_{V,\lambda}) = j_{N,\lambda}\). 
\end{proof}
\end{lma}
In the case \(\lambda = \infty\) above, we are heavily abusing notation: when
\(M_\infty\) is not set-like (and so cannot be identified with a transitive
class), a more careful statement would involve isomorphism rather than equality.

We now further explore the set-likeness of \(M_\infty\).

\begin{defn}
	We say an ordinal \(\kappa\) {\it can be mapped arbitrarily high by
	ultrapower embeddings} if for all \(\alpha > \kappa\), there is an
	ultrapower embedding \(j : V\to M\) such that \(j(\kappa) > \alpha\). The
	{\it ultrapower threshold}\index{Ultrapower threshold} is the least ordinal
	\(\kappa\) that can be mapped arbitrarily high by ultrapower embeddings.
\end{defn}

The existence of the ultrapower threshold is a large cardinal principle closely
related to two recently popularized weakenings of strong compactness: the
strongly tall cardinals of Hamkins \cite{Hamkins} and the \(\omega_1\)-strongly
compact cardinals of Bagaria-Magidor \cite{MagidorBagaria}. Certainly a strongly
tall cardinal or an \(\omega_1\)-strongly compact cardinal is greater than or
equal to the ultrapower threshold. (By theorems of Gitik \cite{Gitik}, it is
consistent with ZFC that these inequalities are strict.) The ultrapower
threshold is in a sense a hybrid of these notions in the sense that it weakens
strong compactness in the Hamkins and Bagaria-Magidor directions simultaneously,
producing a super-weakening of strong compactness. 

If it exists, the ultrapower threshold is a Beth fixed point, but by the
arguments of Bagaria-Magidor \cite{MagidorBagaria}, it cannot be proved to be
inaccessible in ZFC. The nonexistence of the ultrapower threshold has the
following structural consequence for ultrapower embeddings:

\begin{lma}\label{ArbHigh}
If the ultrapower threshold does not exist, then unboundedly many ordinals are
fixed by all ultrapower embeddings.
\begin{proof}
Fix an ordinal \(\xi\). Let \[T_\xi = \{j(\xi) : j\text{ is an ultrapower
embedding of }V\}\]  and let \(C\) be the class of ordinals fixed by all
ultrapower embeddings. Note that if \(\xi\) is an ordinal and \(i: V\to N\) is
an ultrapower embedding, \[i(T_\xi) = \{j(i(\xi)) : j\text{ is an ultrapower
embedding of }N\}\subseteq T_\xi\] This is because the composition of ultrapower
embeddings is an ultrapower embedding.

Since the ultrapower threshold does not exist, \(T_\xi\) is a set for all
ordinals \(\xi\). So assume \(T_\xi\) is a set, and we will show \(C\setminus
\xi\) is nonempty. Let \(\alpha = \sup(T_\xi)\). Obviously \(\alpha \geq \xi\).
Suppose \(i\) is an ultrapower embedding. Then \[i(\alpha) = \sup i(T_\xi) \leq
\sup(T_\xi) = \alpha\] Thus \(i(\alpha) = \alpha\). It follows that \(\alpha\in
C\setminus \xi\), as desired. Thus if \(T_\xi\) is a set, then there is an
ordinal above \(\xi\) fixed by all ultrapower embeddings. It follows that if
\(T_\xi\) is a set for all \(\xi\), then \(C\) is a proper class.
\end{proof}
\end{lma}

The embedding \(j_{V,\infty}\) actually encodes the class of common fixed points
of ultrapower embeddings by a standard argument: 

\begin{lma}\label{FixClm} An ordinal belongs to the range of \(j_{V,\infty}\) if and only if it is fixed by all ultrapower embeddings.\end{lma}
\begin{proof}
	Suppose first that \(\beta\) is an ordinal in the range of \(j_{V,\infty}\).
	Fix an ordinal \(\alpha\) such that \(\beta = j_{V,\infty}(\alpha)\).
	Suppose \(i: V\to N\) is an ultrapower embedding. Then \(i(\beta) =
	i(j_{V,\infty}(\alpha)) = i(j_{V,\infty})(i(\alpha)) = j_{N,\infty}\circ
	i(\alpha) = j_{V,\infty}(\alpha) = \beta\). Thus \(\beta\) is fixed by all
	ultrapower embeddings.
	
	Conversely, suppose \(\beta\) is fixed by all ultrapower embeddings. Let
	\(\alpha\) be the least ordinal such that \(j_{V,\infty}(\alpha) \geq
	\beta\). Suppose towards a contradiction that \(j_{V,\infty}(\alpha) >
	\beta\). Then there is an ultrapower embedding \(i : V\to N\) and some
	ordinal \(\xi< i(\alpha) \) and \(j_{N,\infty}(\xi) \geq \beta\). But by the
	elementarity of \(i\), \(i(\alpha)\) is the least ordinal \(\alpha'\) such
	that \(i(j_{V,\infty})(\alpha') \geq i(\beta)\). Since \(i(\beta) = \beta\),
	this means that \(i(\alpha)\) is the least ordinal \(\alpha'\) such that
	\(j_{N,\infty}(\alpha')\geq \beta\). This contradicts the existence of \(\xi
	< i(\alpha)\) such that \(j_{N,\infty}(\xi) \geq \beta\).\end{proof}

\begin{thm}[UA]\label{MInftyDichotomy1}\index{Direct limit of all ultrapowers!set-likeness}
	Exactly one of the following holds:
	\begin{enumerate}[(1)]
		\item The ultrapower threshold exists.
		\item \(M_\infty\) is set-like.
	\end{enumerate}
\begin{proof}
Suppose (1) holds. For any ordinal \(\xi\), the images of \(\xi\) under
ultrapower embeddings are bounded above by \(j_{V,\infty}(\xi)\), and so \(\xi\)
is not the ultrapower threshold. Thus (2) fails.
	
Suppose (2) fails. By \cref{ArbHigh}, the class \(C\) of ordinals fixed by all
ultrapower embeddings is unbounded in the ordinals. By \cref{FixClm}, \(C =
j_{V,\infty}[\text{Ord}]\cap \text{Ord}\). The function
\(j_{V,\infty}\restriction \text{Ord}\) is therefore equal to the increasing
enumeration of \(C\). It follows that for every ordinal \(\alpha\),
\(j_{V,\infty}(\alpha)\) is an ordinal. In other words, \(M_\infty\) is
set-like. Thus (1) holds.
\end{proof}
\end{thm}

The analysis of supercompactness under UA has the following surprising
consequence (\cref{ThreshThm}): the ultrapower threshold is supercompact. 
\begin{thm}[UA]\label{SetLikeSupercompact}
Exactly one of the following holds:
\begin{enumerate}[(1)]
\item \(M_\infty\) is set-like.
\item There is a supercompact cardinal.
\end{enumerate}
\begin{proof}[Proof given \cref{ThreshThm}] Suppose \(M_\infty\) is not
	set-like. By \cref{MInftyDichotomy1}, the ultrapower threshold exists. By
	\cref{ThreshThm}, the ultrapower threshold is supercompact.
\end{proof}
\end{thm}
In fact, if the ultrapower threshold \(\kappa\) is supercompact, then
\(j_{V,\infty}(\kappa)\) is isomorphic to \(\text{Ord}\) while
\(j_{V,\infty}\restriction V_\kappa = (j_{V,\infty})^{V_\kappa}\).

We now explain the connection between the models \(M_\lambda\) and the Ketonen
order on pointed ultrapower embeddings.
\begin{defn}\label{ERank}\index{Ketonen order!on pointed ultrapowers!rank (\(o_\lambda(M)\))}\index{\(o_\lambda(M)\) (rank in the Ketonen order on pointed ultrapowers)}
Let \(\mathcal P_\lambda\) denote the collection of pointed ultrapowers
\((M,\xi)\) such that \(M\in \mathcal D_\lambda\). For any \(M\in\mathcal
P_\lambda\), \(o_\lambda(M)\) denotes the rank of \(M\) in the Ketonen order
restricted to \(\mathcal P_\lambda\). (If \(\lambda = \infty\), this rank may
not exist.) If \(W\) is a countably complete ultrafilter then \(o_\lambda(W) =
o_\lambda(M_W,\id_W)\).
\end{defn}

The following theorem shows that the ordinals \(o_\lambda(M)\) are highly
structured under UA:

\begin{thm}[UA]\label{RankEmb}
Assume \(\lambda\) is regular or \(\lambda = \infty\). For any \(M\in \mathcal
P_\lambda\), \(o_\lambda(M) = j_{M,\lambda}(\xi_M)\).
\begin{proof}
Consider the partial function \(\Phi_\lambda : \mathcal P_\lambda\to
\text{Ord}\) defined by \[\Phi_\lambda(M) = j_{M,\lambda}(\xi_M)\] If \(\lambda
= \infty\), we leave \(\Phi_\lambda(M)\) undefined if \(j_{M,\lambda}(\xi_M)\)
is not in the set-like part of \(M_\infty\).

For \(M,N\in \mathcal P_\lambda\), \(M \sE N\) implies \(\Phi_\lambda(M)\sE
\Phi_\lambda(N)\) and \(M \KE N\) implies \(\Phi_\lambda(M) = \Phi_\lambda(N)\).
Moreover the image of \(\Phi_\lambda\) is the set-like initial segment of
\(\text{Ord}^{M_\lambda}\). Therefore \(\Phi_\lambda\) is equal to the rank
function of \((\mathcal P_\lambda,\sE)\). That is, for all \(M\in \mathcal
P_\lambda\), \(o_\lambda(M) = \Phi_\lambda(M) = j_{M,\lambda}(\xi_M)\).
\end{proof}
\end{thm}

Combining this with \cref{SetLikeSupercompact}, we obtain:
\begin{cor}
Exactly one of the following holds:
\begin{enumerate}[(1)]
\item There is a supercompact cardinal.
\item For any pointed ultrapower \(M\), \(o_\infty(M)\) exists.\qed
\end{enumerate}
\end{cor}
In conclusion, the Ketonen order, which is in a sense the simplest structure
associated with the Ultrapower Axiom, bears a deep relationship to
supercompactness under UA. This relationship is one of the topics of
\cref{SCChapter1} and \cref{SCChapter2}.
\subsection{The Lipschitz order}\label{LipschitzSection}
In this short subsection, we describe a generalization of the Ketonen order that
raises an interesting philosophical question. Throughout the section, we fix an
infinite ordinal \(\delta\).
\begin{defn}
Suppose \(f : P(\delta)\to P(\delta)\). Then \(f\) is:
\begin{itemize}
\item a {\it reduction} if for \(A\subseteq \delta\) and \(\alpha < \delta\),
\(f(A)\cap \alpha\) depends only on \(A\cap \alpha\). 
\item a {\it contraction} if for \(A\subseteq\delta\) and \(\alpha < \delta\),
\(f(A)\cap (\alpha+1)\) depends only on \(A\cap \alpha\).
\end{itemize}
We say \(X\) {\it reduces} ({\it contracts}) to \(Y\) if there is a reduction
(contraction) \(f : P(\delta)\to P(\delta)\) such that \(f^{-1}[Y] = X\). In
this case we say \(f\) is a reduction (contraction) from \(X\) to \(Y\).
\end{defn}

These concepts can be formulated in terms of long games:
\begin{defn}\index{Lipschitz game}
In the {\it Lipschitz game of length \(\delta\)} associated to sets
\(X,Y\subseteq P(\delta)\), denoted \(G_\delta(X,Y)\), two players I and II
alternate playing 0s or 1s. I plays at limit stages. The play lasts for
\(\delta\cdot 2\) moves, so that I and II produce sequences
\(x_\text{I},x_\text{II}\in {}^\delta2\). Let \(A_\text{I} = \{\alpha < \delta :
x_\text{I}(\alpha) = 1\}\) and \(A_\text{II}  = \{\alpha < \delta :
x_\text{II}(\alpha) = 1\}\). Then II wins if \(A_\text{I}\in X\iff
A_\text{II}\in Y\).
\end{defn}

Player II has a winning strategy if and only if \(X\) reduces to \(Y\), and
Player I has a winning strategy if and only if \(Y\) contracts to
\(P(\delta)\setminus X\).

\begin{defn}\index{Lipschitz order}
The {\it Lipschitz order} is defined on \(X,Y\subseteq P(\delta)\) by setting
\(X\sLi Y\) if \(X\) and \(P(\delta)\setminus X\) contract to \(Y\). The {\it
nonstrict Lipschitz order} is defined on \(X,Y\subseteq P(\delta)\) by setting
\(X\Li Y\) if \(X\) reduces to \(Y\).
\end{defn}

This notation is perhaps misleading since it might suggest that \(X\sLi Y\) if
and only if \(X\Li Y\) and \(Y\not \Li X\). Under the Axiom of Determinacy, this
is true when \(X\) and \(Y\) are contained in \(P(\omega)\). 

The Lipschitz order is transitive in the following strong sense:

\begin{lma}
The composition of a contraction and a reduction is a contraction. Therefore is
\(X\) contracts to \(Y\) and \(Y\) reduces to \(Z\), then \(X\) contracts to
\(Z\). In particular, if \(X \sLi Y \Li Z\) then \(X\sLi Z\).\qed
\end{lma}

A generalization of the proof of \cref{KOIrreflexive} shows that the Lipschitz
order is irreflexive:
\begin{lma}\label{LipschitzIrreflexive}
Suppose \(X\subseteq P(\delta)\). Then \(X\) does not contract to
\(P(\delta)\setminus X\).
\begin{proof}
It suffices to show that every contraction \(f : P(\delta)\to P(\delta)\) has a
fixed point \(A\): then \(A\in X\) if and only if \(f(A)\in X\) so \(f\) is not
a contraction from \(X\) to \(P(\delta)\setminus X\). 

We define \(A\) by recursion. Suppose \(\alpha < \delta\) and we have defined
\(A\cap \alpha\). We then put \(\alpha\in A\) if and only if \(\alpha\in f(A\cap
\alpha)\). Then for any \(\alpha < \delta\), \begin{align*}\alpha\in A&\iff
\alpha\in f(A\cap \alpha)\\&\iff\alpha\in f(A)\end{align*} The final equivalence
follows from the fact that \(f\) is a contraction. Thus \(f(A) = A\), as
desired.
\end{proof}
\end{lma}

\begin{cor}
The Lipschitz order is a strict partial order.\qed
\end{cor}

By the proof of the Martin-Monk theorem (see \cite{VanWesep}) descending
sequences in the Lipschitz order give rise to pathological subsets of Cantor
space:
\begin{thm}[ZF + DC]\label{MartinMonk}
The following are equivalent:
\begin{enumerate}[(1)]
\item There is a flip set.
\item The Lipschitz order on \(P(\omega)\) is illfounded.
\item The Lipschitz order on \(P(\delta)\) is illfounded.
\end{enumerate}
\begin{proof}
To see (1) implies (2), suppose \(F\subseteq 2^\omega\) is a flip set. Define
\((E_n)_{n < \omega}\) by recursion, setting \(E_0 = F\) and \(E_{n+1} = \{s\in
2^\omega : 0s \in E_n\}\). It is easy to see that \(E_{n+1}\) and \(2^\omega
\setminus E_{n+1}\) both contract to \(E_n\), via the contractions \(s\mapsto
0s\) and \(s\mapsto1s\) respectively.

(2) trivially implies (3).

We finally show (3) implies (1). Fix \(X_0\sgLi X_1 \sgLi X_2\sgLi\cdots\) a
descending sequence of subsets of \(P(\delta)\). For \(n < \omega\), fix
contractions \(f_n^0\) from \(X_{n+1}\) to \(X_n\) and \(f_n^1\) from
\(X_{n+1}\) to \(P(\delta)\setminus X_n\). For each \(s\in 2^\omega\), we define
sets \(A_n^s\subseteq \delta\) such that \[A_n^s = f_n^{s(n)}(A_{n+1}^s)\]
Suppose \(A_n^s\cap\alpha\) has been defined for all \(n < \omega\). Then
\[A_n^s\cap (\alpha + 1) = f_n^{s(n)}(A_{n+1}^s\cap \alpha)\cap (\alpha + 1)\]
Since \(f_n^i\) is a contraction for all \(n < \omega\) and \(i \in \{0,1\}\),
\(A_n^s\) is well-defined and \(A_n^s = f^{s(n)}_n(A^s_{n+1})\). 

Define \(F_n\subseteq 2^\omega\) by putting \(s\in F_n\) if and only if
\(A_n^s\in X_n\). Whether \(s\in F_n\) depends only on \(s\restriction
(\omega\setminus n)\). Moreover, if \(s\in F_{n+1}\) then \(s\in F_n\) if and
only if \(s(n) = 0\). It is easy to show by induction that if \(s\) and \(s'\)
agree on \(\omega\setminus n\) and \(\sum_{k < n} s(k) = \sum_{k < n} s'(k)\mod
2\), then \(s\in F_0\) if and only if \(s'\in F_0\). Similarly, if \(s\) and
\(s'\) agree on \(\omega\setminus n\) and \(\sum_{k < n} s(k) \neq \sum_{k < n}
s'(k)\mod 2\), then \(s\in F_0\) if and only if \(s'\notin F_0\). It follows
that \(F_0\) is a flip set.
\end{proof}
\end{thm}
Of course, (1), (2), and (3) are all provable in ZFC. In the choiceless context
of ZF + DC, however, there may be no flip sets (for example, if every subset of
Cantor space has the Baire property or is Lebesgue measurable). In this case,
\cref{MartinMonk} shows that the Lipschitz order is wellfounded not only on
subsets of Cantor space but also on subsets of \(P(\delta)\).\footnote{Under the
same hypotheses, one can show that the Lipschitz order on \({}^\delta S\) is
wellfounded for any set \(S\) after generalizing the definition of the Lipschitz
order in the natural way.} The proof also shows that the wellfounded part of the
Lipschitz order is equal to the collection of sets that do not lie above a flip
set. 

We turn now to the relationship between the Lipschitz order and the Ketonen
order. 

\begin{defn}\label{SetConcentration}
A set \(Z\subseteq P(\delta)\) {\it concentrates} on a set \(S\) if for all
\(A,B \subseteq \delta\) with \(A\cap S = B\cap S\), \(A\in Z\) if and only if
\(B\in Z\).
\end{defn}

Note that if \(Z\) is an ultrafilter that concentrates on a class \(S\) in the
sense of \cref{UltrafilterConcentration}, then \(Z\) concentrates on \(S\) is
the sense of \cref{SetConcentration}.

\begin{lma}\label{LipChar}
Suppose \(X\subseteq P(\delta)\) and \(W\) is an ultrafilter on \(\delta\). Then
the following are equivalent:
\begin{enumerate}[(1)]
\item \(X\sLi W\).
\item \(X\) contracts to \(W\).
\item For some \(Z\in M_W\) that concentrates on \(\id_W\), \(X = j_W^{-1}[Z]\).
\end{enumerate}
\begin{proof}
{\it (1) implies (2)}: Trivial.

{\it (2) implies (1)}: Assume \(X\) contracts to \(W\). Since \(W\) is an
ultrafilter, \(W\) reduces to \(P(\delta)\setminus W\). Since \(X\) contracts to
\(W\) and \(W\) reduces to \(P(\delta)\setminus W\), \(X\) contracts to
\(P(\delta)\setminus W\). Therefore \(X\sLi W\).

{\it (1) implies (3)}: Let \(f : P(\delta)\to P(\delta)\) be a contraction from
\(X\) to \(W\). For each \(\alpha\), let \(X_\alpha = \{A\subseteq \delta :
\alpha\in f(A)\}\). Since \(f\) is a contraction, \(X_\alpha\) concentrates on
\(\alpha\). Let \(Z = [\langle X_\alpha : \alpha < \delta\rangle]_W\). By \L
o\'s's Theorem, \(Z\) concentrates on \(\id_W\). Then 
\begin{align*}
A\in X&\iff f(A)\in W\\
&\iff \{\alpha < \delta : A\in X_\alpha\}\in W\\
&\iff j_W(A)\in Z
\end{align*}
Thus \(j^{-1}_W[Z] = X\).

{\it (3) implies (1)}: Fix \(Z\in M_W\) concentrating on \(\id_W\) such that \(X
= j_W^{-1}[Z]\). Let \(\langle X_\alpha : \alpha \in I\rangle\) be such that \(Z
= [\langle X_\alpha : \alpha \in I\rangle]_W\) and \(X_\alpha\) concentrates on
\(\alpha\) for all \(\alpha\in I\). Define \(f : P(\delta)\to P(\delta)\) by
setting \(f(X) = \{\alpha \in I : X\in X_\alpha\}\). Then \(f\) is a contraction
since \(X_\alpha\) concentrates on \(\alpha\) for all \(\alpha\in I\). Moreover,
\begin{align*}
A\in X&\iff j_W(A)\in Z\\
&\iff \{\alpha < \delta : A\in X_\alpha\}\in W\\
&\iff f(X)\in W\qedhere
\end{align*}
\end{proof}
\end{lma}

Using \cref{KOChars}, this has the following corollary:

\begin{cor}\label{LipKet}\index{Ketonen order!vs. the Lipschitz order}
The Lipschitz order extends the Ketonen order on \(\mathscr B(\delta)\).\qed
\end{cor}

Under UA, it follows that the two orders coincide:

\begin{cor}[UA] The Lipschitz order and the Ketonen order coincide on \(\mathscr
B(\delta)\). In particular, the Lipschitz order linearly orders \(\mathscr
B(\delta)\).
\begin{proof}
Since \(\sLi\) is a strict partial order extending the total relation \(\sE\)
(\cref{Totality}), the two orders must be equal.
\end{proof} 
\end{cor}

Another way to state this is as a determinacy consequence of UA:

\begin{cor}[UA]\index{Ultrapower Axiom!vs. long determinacy}
	For all ordinals \(\delta\), for any \(U,W\in \mathscr B(\delta)\), the game
	\(G_\delta(U,W)\) is determined.
\end{cor}

We conclude this section with a question that is perhaps of some philosophical
significance:
\begin{qst}\label{LipUA}
Assume that for any ordinal \(\delta\), for any \(U,W\in \mathscr B(\delta)\),
the game \(G_\delta(U,W)\) is determined. Does the Ultrapower Axiom hold?  
\end{qst}
If this were true then the Ultrapower Axiom would be a long determinacy
principle. In \cref{LinKetSection}, we give partial positive answer.

\subsection{The Ketonen order on filters}\label{FilterSection}
We briefly discuss a generalization of the Ketonen order to a wellfounded
partial order on arbitrary countably complete filters that is suggested by the
proof of \cref{KOWellfounded}. This order will not appear elsewhere in this
dissertation, but it seems potentially quite interesting since it identifies a
connection between the Ketonen order and stationary reflection.

\begin{defn}
	Suppose \(F\) is a filter, \(I\in F\), and \(\langle G_i : i\in I\rangle\)
	is a sequence of filters on a fixed set \(Y\). The {\it \(F\)-limit of
	\(\langle G_i : i\in I\rangle\)} is the filter \[F\text{-}\lim_{i\in I} G_i
	= \{A\subseteq Y : \{i \in I : A\in G_i\}\in F\}\]
\end{defn}

\begin{defn}
	If \(F\) is a filter on  a set \(X\) and \(C\) is a class, then \(F\) {\it
	concentrates} on \(C\) if \(C\cap X\in F\).
\end{defn}

\begin{defn}
	Suppose \(X\) is a set and \(C\) is a class. Let \(\mathscr F(X)\) denote
	the set of countably complete filters on \(X\) and let \(\mathscr F(X,C)\)
	denote the set of filters on \(X\) that concentrate on \(C\).
\end{defn}
\begin{defn}\index{Ketonen order!on filters}
	Suppose \(\epsilon\) and \(\delta\) are ordinals, \(F\in\mathscr
	F(\epsilon)\), and \(G\in \mathscr F(\delta)\). The {\it Ketonen order on
	filters} is defined on by setting \(F\sE G\) if there is a set \(I\in G\)
	and a sequence \(\langle F_\alpha : \alpha \in I\rangle\in \prod_{\alpha\in
	I}\mathscr F(\epsilon,\alpha)\) such that \(F \subseteq
	G\text{-}\lim_{\alpha\in I} F_\alpha\).
\end{defn}
Under the Ultrapower Axiom, the restriction to ultrafilters of the Ketonen order
on filters coincides with the Ketonen order as it is defined in
\cref{KOCharsSection}. We do not know whether this is provable in ZFC.

Note that the proof of \cref{KOIrreflexive} breaks down when we consider filters
instead of ultrafilters. In fact, in a sense this simple proof cannot be
remedied, since irreflexivity fails if we allow filters that are countably
incomplete, and it is not clear how countable completeness could come in to the
argument of \cref{KOIrreflexive}.  It is somewhat surprising that one can in
fact prove the irreflexivity of the Ketonen order by instead using countable
completeness and the argument of \cref{KOWellfounded}:

\begin{thm}\label{FilterWellfounded}
	The Ketonen order on filters is wellfounded.
\end{thm}

We include the proof, which is closely analogous to that of
\cref{KOWellfounded}.
\begin{lma}\label{FilterResemblance}
	Suppose \(H\) is a filter and \(F\sE G\) are countably complete filters on
	ordinals \(\epsilon\) and \(\delta\). Suppose \(J\in H\) and \(\langle G_x :
	x \in J\rangle\) is a sequence of countably complete filters such that \(G
	\subseteq H\text{-}\lim_{x \in J} G_x\). Then there is a set \(J'\subseteq
	J\) in \(H\) and a sequence \(\langle F_x: x \in J'\rangle\) of countably
	complete filters such that \(F_x \sE G_x\) for all \(x\in K\) and \(F
	\subseteq H\text{-}\lim_{x \in K} F_x\).
	\begin{proof}
		Since \(F\sE G\), we can fix \(I\in G\) and countably complete filters
		\(\langle D_\alpha : \alpha\in I\rangle\in \prod_{\alpha\in I} \mathscr
		F(\epsilon,\alpha)\) such that \(F \subseteq G\text{-}\lim_{\alpha\in I}
		D_\alpha\).
		
		Let \(J' = \{b \in J : I\in G_x\}\). Since \(I\in G \subseteq
		H\text{-}\lim_{x\in J} G_x\), we have \(J'\in H\) by the definition of a
		limit. For each \(x\in J'\), let \[F_x = G_x\text{-}\lim_{\alpha\in I}
		D_\alpha\] Then \(F_x\in \mathscr B(\epsilon)\), and the sequence
		\(\langle D_\alpha : \alpha\in I\rangle\) witnesses \(F_b \sE G_b\).
		
		Finally,
		\begin{align*}
			F &\subseteq G\text{-}\lim_{\alpha\in I} D_\alpha\\
			&\subseteq (H\text{-}\lim_{x\in J} G_x)\text{-}\lim_{\alpha\in I} D_\alpha\\
			&= H\text{-}\lim_{x\in J'} (G_x\text{-}\lim_{\alpha\in I} D_\alpha)\\
			&= H\text{-}\lim_{x\in J'} F_x
		\end{align*}
		Thus \(F\subseteq H\text{-}\lim_{x\in K} F_x\), as desired.
	\end{proof}
\end{lma}

\begin{proof}[Proof of \cref{FilterWellfounded}] Suppose towards a contradiction
	that \(\delta\) is the least ordinal such that the Ketonen order is
	illfounded below a countably complete filter that concentrates on
	\(\delta\). Fix a descending sequence \(F_0\sgE F_1\sgE F_2\sgE\cdots\) such
	that \(F_0\) concentrates on \(\delta\).
	
	We will define sets of ordinals \(I_1\supseteq I_2\supseteq \cdots\) in
	\(F\) and sequences \(\langle F^{m}_\alpha : \alpha \in I_m\rangle\) of
	countably complete filters such that \[F_m \subseteq F\text{-}\lim_{\alpha
	\in I_m} F^{m}_\alpha\] for all \(1\leq m< \omega\). We will have:
	\begin{itemize}
		\item For all \(\alpha\in I_1\), \(F_\alpha^1\) concentrates on
		\(\alpha\).
		\item For all \(1\leq m< \omega\), for all \(\alpha \in I_{m+1}\),
		\(F^{m+1}_\alpha \sE F^{m}_\alpha\).
	\end{itemize}
	
	Since \(F_1\sE F\), there is a set of ordinals \(I_1\in F\) and a sequence
	\(\langle F^1_\alpha : \alpha\in I_1\rangle\) of countably complete
	ultrafilters such that \(F_1 \subseteq F\text{-}\lim_{\alpha \in I_1}
	F^1_\alpha\) and \(F^1_\alpha\) concentrates on \(\alpha\) for all
	\(\alpha\in I_1\).
	
	Suppose \(1\leq m < \omega\) and \(\langle F^{m}_\alpha : \alpha \in
	I_m\rangle\) has been defined. We now apply \cref{FilterResemblance} with
	\(H = F\), \(F = F_{m+1}\), and \(G =F_m\). This yields a set
	\(I_{m+1}\subseteq I_m\) in \(F\) and a sequence \(\langle F^{m+1}_\alpha :
	\alpha \in I_m\rangle\) of countably complete filters on \(\delta\) such
	that \(F^{m+1}_\alpha \sE F^{m}_\alpha\) for all \(\alpha\in I_{m+1}\) and
	\[F_{m+1} \subseteq F\text{-}\lim_{\alpha\in I_{m+1}}F^{m+1}_\alpha\]
	
	This completes the definition of the sets \(I_1\supseteq I_2\supseteq
	\cdots\) and sequences \(\langle F^{m}_\alpha : \alpha \in I_m\rangle\)  for
	\(1\leq m < \omega\).
	
	Now let \(I = \bigcap_{1\leq m < \omega} I_m\). Since \(F_0\) is countably
	complete, \(I\) is nonempty, so we can fix an ordinal \(\alpha\in I\). Then
	since \(\alpha\in I_m\) for all \(1\leq m < \omega\), \[F^1_\alpha \sgE
	F^2_\alpha \sgE F^3_\alpha \sgE\cdots\] Since \(F^1_\alpha\) concentrates on
	\(\alpha < \delta\), this contradicts the minimality of \(\delta\).
\end{proof}

Recall the following definition, due to Jech \cite{JechStationary}:
\begin{defn}
	Assume \(\delta\) is a regular cardinal. The {\it canonical order on
	stationary sets} is defined on stationary sets \(S,T\subseteq \delta\) by
	setting \(S < T\) if there is a closed unbounded set \(C\subseteq \delta\)
	such that \(S\cap \alpha\) is stationary in \(\alpha\) for all \(\alpha\in
	C\cap T\).
\end{defn}

\begin{defn}
	For any ordinal \(\alpha\), let \(\mathcal C_\alpha\) denote the filter of
	closed cofinal subsets of \(\alpha\).
\end{defn}

\begin{defn}
	Suppose \(F\) is a filter on a set \(X\) and \(S\) is a set such that \(F\)
	does not concentrate on the complement of \(S\). The {\it projection of
	\(F\) on \(S\)} is the filter defined by \[F\mid S = \{A\cap S : A\in F\}\]
\end{defn}

The following proposition connects the canonical order on stationary sets and
the Ketonen order on filters:
\begin{prp}\label{StationaryKetonen} Suppose \(\delta\) is a regular cardinal
and \(S\) and \(T\) are stationary subsets of \(\delta\). Then \(S < T\) implies
\(\mathcal C_\delta\mid S\sE \mathcal C_\delta\mid T\).
	\begin{proof}
		Fix a closed unbounded set \(C\subseteq \delta\) such that \(S\cap
		\alpha\) is stationary in \(\alpha\) for all \(\alpha\in C\cap T\). Note
		that \(C\cap T\in \mathcal C_\delta\mid T\), and for all \(\alpha\in
		C\cap T\), \(\mathcal C_\alpha\mid S\) is a countably complete filter
		concentrating on ordinals less than \(\alpha\).
		\begin{clm}
			\(C_\delta\mid S \subseteq (\mathcal C_\delta\mid
			T)\text{-}\lim_{\alpha\in C\cap T} \mathcal C_\alpha\mid S\). 
		\end{clm}
		\begin{proof}
			Suppose \(A\in \mathcal C_\delta\mid S\). We will show that \(A\in
			(\mathcal C_\delta\mid T)\text{-}\lim_{\alpha\in C\cap T} \mathcal
			C_\alpha\mid S\). Fix \(E\in C_\delta\) such that \(S\cap E\subseteq
			A\). Let \(E'\) be the set of accumulation points of \(E\). Then for
			any \(\alpha\in E'\), \(S\cap (E\cap \alpha)\subseteq A\) and
			\(E\cap \alpha\in \mathcal C_\alpha\), so \(A\in \mathcal
			C_\alpha\mid S\). Thus \[E'\cap C\cap T\subseteq \{\alpha \in C\cap
			T : A\in \mathcal C_\alpha\mid S\}\] Since \(E'\cap C\in \mathcal
			C_\delta\), \(E'\cap C\cap T\in \mathcal C_\delta\mid T\), and
			therefore \(\{\alpha \in C\cap T : A\in \mathcal C_\alpha\mid S\}\in
			\mathcal C_\delta\mid T\). It follows that \(A\in (\mathcal
			C_\delta\mid T)\text{-}\lim_{\alpha\in C\cap T} \mathcal
			C_\alpha\mid S\), as desired.
		\end{proof}
		The claim implies \(C_\delta\mid S\sE \mathcal C_\delta\mid T\), as
		desired.
	\end{proof}
\end{prp}
As a corollary of \cref{FilterWellfounded} and \cref{StationaryKetonen}, we have
the following theorem of Jech:
\begin{cor}
	The canonical order on stationary sets is wellfounded.\qed
\end{cor}

\section{The linearity of the Ketonen order}\label{LinKetSection}
In this final section, we prove a converse to \cref{SeedKet}, which can also be
seen as a partial positive answer to \cref{LipUA}. We say that {\it the Ketonen
order is linear} if  for all ordinals \(\delta\), the Ketonen order on
\(\mathscr B(\delta)\) is a linear order. The Ketonen order is linear if and
only if its restriction to \(\Un\) is a linear order.

\begin{thm}\label{LinKetThm}\index{Ketonen order!linearity}
The Ketonen order is linear if and only if the Ultrapower Axiom holds.
\end{thm}
Given \cref{KOChars}, the linearity of the Ketonen order would appear to be a
much weaker assumption than UA: the linearity of the Ketonen order only
guarantees comparisons \((k,h): (M_U,M_W)\to N\) such that \(h\) is an internal
ultrapower embedding of \(M_W\), while UA asserts the existence of comparisons
with both \(k\) and \(h\) internal. How can one transform partially internal
comparisons into the fully internal comparisons required by UA? 

To properly describe it, let us make some definitions:
\begin{defn}
Suppose \(M_0,\) \(M_1,\) and \(N\) are transitive models of ZFC and
\[(k_0,k_1): (M_0,M_1)\to N\] are elementary embeddings.
\begin{itemize}
\item \((k_0,k_1)\) is {\it \(0\)-internal} if \(k_0\) is definable over \(M_0\).\index{Comparison!\(0\)-internal, \(1\)-internal}
\item \((k_0,k_1)\) is {\it \(1\)-internal} if \(k_1\) is definable over
\(M_1\).
\item \((k_0,k_1)\) is {\it internal} if it is both \(0\)-internal and
\(1\)-internal.
\end{itemize}
\end{defn}
We indicated above that the difficulty in proving \cref{LinKetThm} is that it is
not clear how to transform the \(1\)-internal comparisons given by the linearity
of the Ketonen order into the internal comparisons required to witness UA. In
fact, it is simply impossible to do this in general, since as a consequence of
the proof of \cref{ZFCComparison}, \(1\)-internal comparisons can be proved to
exist in ZFC alone.
\begin{prp}
Any pair of ultrapower embeddings has a \(0\)-internal ultrapower comparison and
a \(1\)-internal ultrapower comparison.\qed
\end{prp}
Thus the true power of the linearity of the Ketonen order lies not in the mere
existence of \(1\)-internal comparisons \((k,h)\) but rather in the existence of
\((k,h)\) witnessing \(U\sE W\) (or \(W\E U\)); that is, with the additional
property \(k(\id_U) < h(\id_W)\).

\cref{LinKetThm} is an immediate consequence of our next theorem, which shows
how to explicitly define a comparison of a pair of ultrafilters:

\begin{thm}\label{Reciprocity1}
	Assume the Ketonen order is linear. Suppose \(\epsilon\) and \(\delta\) are
	ordinals. Suppose \(U\in \mathscr B(\epsilon)\) and \(W\in \mathscr
	B(\delta)\).
	\begin{itemize}
		\item Let \(W_*\) be the least element of \(j_U(\mathscr
		B(\delta),\sE)\) extending \(j_U[W]\). 
		\item Let \(U_*\) be the least element of \(j_W(\mathscr
		B(\epsilon),\sE)\) extending \(j_W[U]\).
	\end{itemize}
	Then \((j_{W_*}^{M_U},j_{U_*}^{M_W})\) is a comparison of \((j_U,j_W)\).
\end{thm}

The definitions of \(W_*\) and \(U_*\) rely on the fact that \(j_U(\mathscr
B(\delta),\sE)\) and \(j_W(\mathscr B(\epsilon),\sE)\) are wellorders, not only
in \(M_U\) and \(M_W\) but also by absoluteness in the true universe \(V\).
This, however, is not the main use of the linearity of the Ketonen order in the
proof. Indeed, it is consistent that there is a pair of countably complete
ultrafilters \(U\) and \(W\) such that the minimum extensions \(W_*\) and
\(U_*\) are well-defined yet \((j_U,j_W)\) admits no comparison.\footnote{Take
\(U\) and \(W\) to be Mitchell incomparable normal ultrafilters. Apply
\cref{KetMO} and \cref{NormalGeneration} to see that \(j_U(W)\) and \(j_W(U)\)
are the {\it only} extensions of \(j_U[W]\) and \(j_W[U]\) in \(M_U\) and
\(M_W\) respectively.} Instead we will use the linearity of the Ketonen order to
compare \((j_{W_*}^{M_U}\circ j_U,j_{U_*}^{M_W}\circ j_W)\):
\begin{lma}\label{SumComparison}
Suppose \(\epsilon\) and \(\delta\) are ordinals. Suppose \(U\in \mathscr B(\epsilon)\) and \(W\in \mathscr B(\delta)\). 
\begin{itemize}
\item Let \(W_*\) be an element of \(j_U(\mathscr B(\delta))\) extending \(j_U[W]\). 
\item Let \(U_*\) be a minimal element of \(j_W(\mathscr B(\epsilon),\sE)\) extending \(j_W[U]\).
\end{itemize}
For any \(1\)-internal ultrapower comparison \[(k,h) : (M^{M_U}_{W_*},M^{M_W}_{U_*})\to P\] of \((j_{W_*}^{M_U}\circ j_U ,j_{U_*}^{M_W}\circ j_W)\), the following hold:
\begin{align}
	h(j_{U_*}^{M_W}(\id_W))&\leq k(\id_{W_*})\label{bs}\\
	h(\id_{U_*})&\leq k(j_{W_*}^{M_U}(\id_U))\label{as}
\end{align}
\begin{proof}
%Let \(W_* = [\langle W_\alpha : \alpha\in I\rangle]_U\) and let \(U_* = [\langle U_\beta : \beta\in J\rangle]_W\). We put down some information for the reader's convenience. By \cref{SumPower},
%\begin{align*}
%	j_F &= j_{W_*}^{M_U}\circ j_U&  j_G &= j_{U_*}^{M_W}\circ j_W\\
%	\alpha_F &= j_{W_*}^{M_U}(\id_U)& \beta_G &= j_{U_*}^{M_W}(\id_W)\\
%	\beta_F &= \id_{W_*}& \alpha_G &= \id_{U_*}
%\end{align*}
Let us direct the reader's attention to the key diagram, \cref{LinKetFig}.
\begin{figure}
	\center
	\includegraphics[scale=.6]{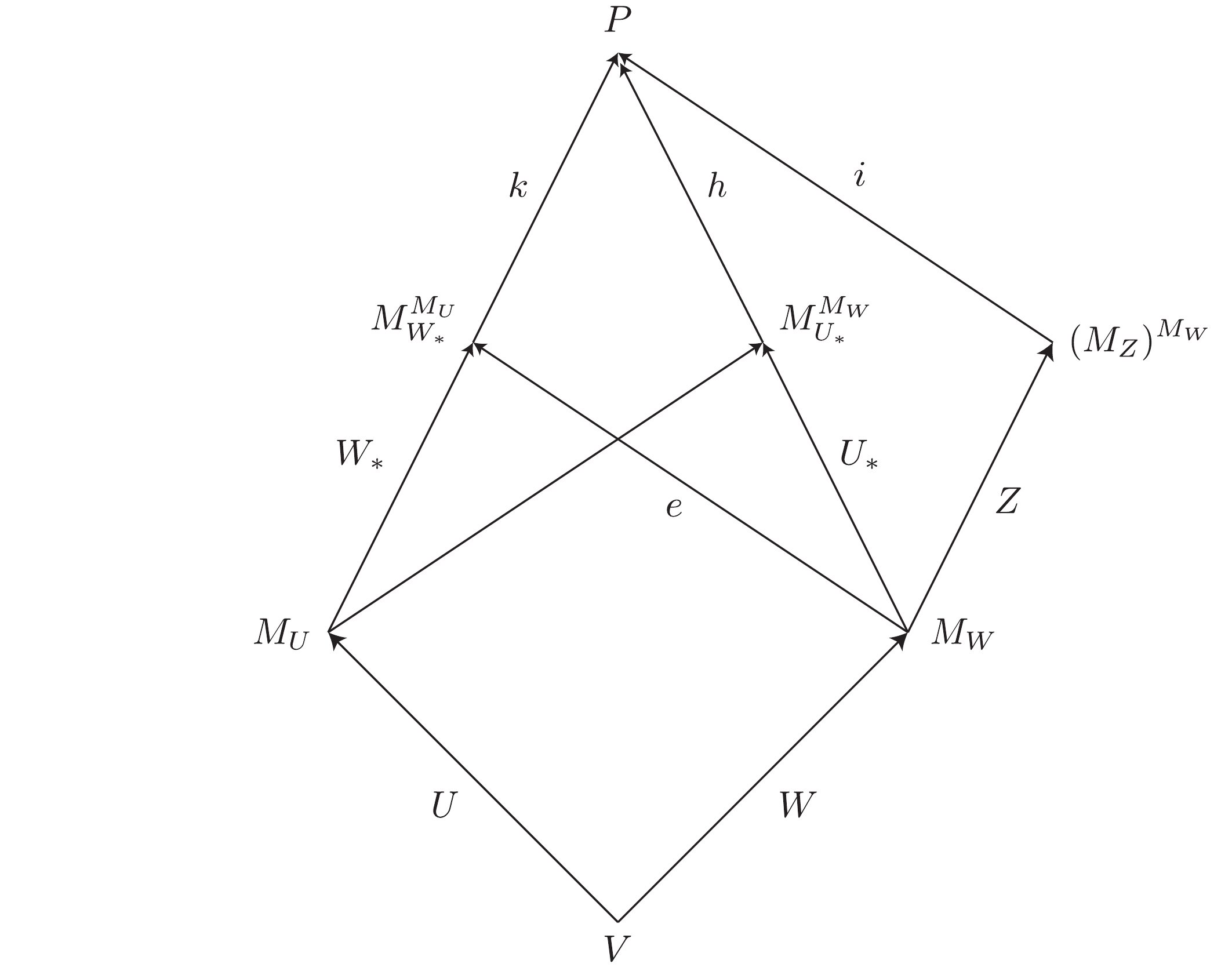}
	\caption{The proof of \cref{SumComparison}.}\label{LinKetFig}
\end{figure}

We first prove \cref{bs}. By \cref{LimitFactor}, there is an elementary embedding \(e : M_W\to M_{W_*}^{M_U}\) such that \(e\circ j_W = j^{M_U}_{W_*}\circ j_U\) and \(e(\id_W) = \id_{W_*}\).
We now apply the minimality of internal embeddings (\cref{MinDefEmb}). Note that \(k\circ e\) and \(h\circ j^{M_W}_{U_*}\) are both elementary embeddings from \(M_W\) to \(P\), but \(h\circ j^{M_W}_{U_*}\) is an internal ultrapower embedding. Thus by \cref{MinDefEmb}, \(h(j^{M_W}_{U_*}(\alpha)) \leq k(e(\alpha))\) for all ordinals \(\alpha\). It follows in particular that \(h(j_{U_*}^{M_W}(\id_W))\leq k(e(\id_W)) = k(\id_{W_{*}})\), proving \cref{bs}.

We now prove \cref{as}. To reduce subscripts, we define: \[\alpha = j_{W_*}^{M_U}(\id_U)\] Let \(Z\) be the \(M_W\)-ultrafilter on \(j_W(\epsilon)\) derived from \(h\circ j_{U_*}^{M_W}\) using \(k(\alpha)\), so \[Z = (h\circ j^{M_W}_{U_*})^{-1}\left[\pr {k(\alpha)} {}\right]\] Since \(h\circ j_{U_*}^{M_W}\) is an internal ultrapower embedding of \(M_W\), \(Z\) is a countably complete ultrafilter of \(M_W\); in other words, \(Z\in j_W(\mathscr B(\epsilon))\). Moreover, it is not hard to compute that \(Z\) extends \(j_W[U]\), or equivalently \(j_W^{-1}[Z] = U\):
\begin{align}
j_W^{-1}[Z] &= j_W^{-1}[(h\circ j^{M_W}_{U_*})^{-1}[\pr {k(\alpha)} {}]]\nonumber\\
&= (h\circ j_{U_*}^{M_W}\circ j_W)^{-1}[\pr {k(\alpha)} {}]\nonumber\\
&= (k\circ j_{W_*}^{M_U}\circ j_U)^{-1}[\pr {k(\alpha)} {}]\nonumber\\
&= (j_{W_*}^{M_U}\circ j_U)^{-1}[k^{-1}[\pr {k(\alpha)} {}]]\nonumber\\
&= (j_{W_*}^{M_U}\circ j_U)^{-1}[\pr {\alpha} {}]   \nonumber  \\
&= j_U^{-1}[(j_{W_*}^{M_U})^{-1}[\pr {j_{W_*}^{M_U}(\id_U)} {}]]\nonumber\\
&= j_U^{-1}[\pr {\id_U} {}] = U\nonumber
\end{align}
Since \(U_*\) is a minimal element of \(j_W(\mathscr B(\epsilon),\sE)\) extending  \(j_W[U]\), \(M_W\) satisfies \(Z\not \sE U_*\).

Since \(Z\) is derived from \(h\circ j_{U_*}^{M_W}\) using \(k(\alpha)\), there is a factor embedding \(i : (M_{Z})^{M_W}\to P\) specified by the following properties:
\begin{align}i\circ j^{M_W}_{Z} &= h\circ j_{U_*}^{M_W}\label{Compy}\\
	i(\id_{Z}) &= k(\alpha)\label{Motion}\end{align}
Note that these properties define \(i\) over \(M_W\). Therefore by \cref{Compy}, \((i,h)\) is a \(1\)-internal ultrapower comparison of \((j_{Z}^{M_W},j_{U_*}^{M_W})\) in \(M_W\). 
The fact that \(Z\not \sE U_*\) in \(M_W\) implies \[h(\id_{U_*})\leq i(\id_Z) = k(\alpha) = k(j_{W_*}^{M_U}(\id_U))\]
proving \cref{as}.
\end{proof}
\end{lma}

\cref{SumComparison} can be read as asserting that the natural ultrafilter representing the embedding \(j_{U_*}^{M_W}\circ j_W\) is not strictly above the one representing \(j_{W_*}^{M_U}\circ j_U\) in the Ketonen order. To make this precise, we need to define what these natural ultrafilters. This is related to the well-known notion of an ultrafilter sum:

\begin{defn}
	Suppose \(U\) is an ultrafilter on \(X\), \(I\) is a set in \(U\), and \(\langle W_i : i \in I\rangle\) is a sequence of ultrafilters on \(Y\). The {\it \(U\)-sum of  \(\langle W_i : i \in I\rangle\)} is the ultrafilter defined by
	\[U\text{-}\sum_{i\in I}W_i = \{A\subseteq X\times Y : \{i\in I : A_i \in W_i\}\in U\}\]
\end{defn}
In the definition above, if \(A\subseteq X\times Y\) and \(i\in X\), then \(A_i = \{j\in Y : (i,j)\in A\}\).

There is an obvious connection between sums and limits: the projection of a sum of ultrafilters onto its second coordinate is precisely equal to the limit of those ultrafilters.

\begin{lma}
	Suppose \(U\) is an ultrafilter, \(I\) is a set in \(U\), and \(\langle W_i : i \in I\rangle\) is a sequence of ultrafilters on \(Y\). Let \(Z = [\langle W_i : i \in I\rangle]_U\) and let \(D =U\text{-}\sum_{i\in I}W_i\). Then \(M_D = M^{M_U}_Z\),  \(j_D = j_Z^{M_U}\circ j_U\), and \(\id_D = (j_Z^{M_U}(\id_U),\id_Z)\).\qed
\end{lma}

Motivated this lemma, we introduce the following nonstandard notation.

\begin{defn}\index{Sum of ultrafilters}\index{\(U\text{-}\sum\) (\(U\)-sum)}
	Suppose \(U\) is an ultrafilter on \(X\), and \(W_*\) is an \(M_U\)-ultrafilter on \(j_U(Y)\). Then \(U\text{-}\sum W_*\) denotes the ultrafilter on \(X\times Y\) derived from \(j_{W_*}^{M_U}\circ j_U\) using \((j_{W_*}^{M_U}(\id_U),\id_{W_*})\).
\end{defn}

In this section we will only require sums of ultrafilters where \(W_*\in M_U\), but it is just more convenient not to choose a representative for \(W_*\).

\begin{lma}\label{SumPower}
	Suppose \(U\) is an ultrafilter and \(W_*\) is an \(M_U\)-ultrafilter on \(j_U(Y)\). Then  \(j_{U\text{-}\sum W_*} = j_{W_*}^{M_U}\circ j_U\), and \(\id_{U\text{-}\sum W_*} = (j_{W_*}^{M_U}(\id_U),\id_{W_*})\).\qed
\end{lma}

In the context of \cref{Reciprocity1}, we would like to use \cref{SumComparison} to conclude that the ultrafilters \(U\text{-}\sum W_*\) and \(W\text{-}\sum U_*\) are either equal or incomparable in the Ketonen order, and thus conclude by the linearity of the Ketonen order that \(U\text{-}\sum W_* = W\text{-}\sum U_*\). The only remaining problem is that \(U\text{-}\sum W_*\) and \(W\text{-}\sum U_*\) are not ultrafilters on ordinals. But obviously we can associate Ketonen orders to an arbitrary wellorder:

\begin{defn}\index{Ketonen order!associated to a wellorder}
	Suppose \((X,\prec)\) is a wellorder. The {\it Ketonen order associated to \((X,\prec)\)} is the order \((\mathscr B(X),\prec^{\Bbbk})\) defined on \(U,W\in \mathscr B(X)\) by setting \(U \prec^{\Bbbk} W\) if there exist \(I\in W\) and \(\langle U_x : x\in I\rangle\in \prod_{x\in I} \mathscr B(X,X_{\prec x})\) such that \(U = W\text{-}\lim_{x\in I} U_x\).
\end{defn}

If \((X,\prec)\) and \((X',\prec')\) are isomorphic wellorders, then the associated Ketonen orders are also isomorphic, so in particular all the characterizations of the Ketonen order generalize to arbitrary wellorders:

\begin{lma}\label{KOWOChar}
	Suppose \((X,\prec)\) is a wellorder and \(U,W\in \mathscr B(X)\). Then the following are equivalent:
	\begin{enumerate}[(1)]
		\item \(U\prec^\Bbbk W\).
		\item There is a \(1\)-internal ultrapower comparison \((k,h) : (M_U,M_W)\to N\) of \((j_U,j_W)\) such that \(k(\id_U)\prec^* h(\id_W)\) where \({\prec}^* =k(j_U(\prec)) = h(j_W(\prec))\).\qed
	\end{enumerate}
\end{lma}

It is convenient to introduce some notation for the statement of \cref{PreciseReciprocity}:

\begin{defn}\label{GodelDef}
	Let \(\text{flip} : \text{Ord}\times  \text{Ord} \to  \text{Ord}\times  \text{Ord}\) be defined by \(\text{flip}(\alpha,\beta) = (\beta,\alpha)\). Let \(\prec\) denote the G\"odel order on \(\text{Ord}\times \text{Ord}\). 
\end{defn}

The only property of the G\"odel order that we need is that \((\alpha_0,\beta_0)\prec (\alpha_1,\beta_1)\) implies that either \(\alpha_0 < \alpha_1\) or \(\beta_0 < \beta_1\).

\begin{lma}\label{PreciseReciprocity}
	Suppose \(\epsilon\) and \(\delta\) are ordinals. Suppose \(U\in \mathscr B(\epsilon)\) and \(W\in \mathscr B(\delta)\). Assume the Ketonen order \((\mathscr B(\epsilon\times \delta),\prec^{\Bbbk})\) is linear.
	\begin{itemize}
		\item Let \(W_*\) be the least element of \(j_U(\mathscr B(\delta),\sE)\) extending \(j_U[W]\). 
		\item Let \(U_*\) be the least element of \(j_W(\mathscr B(\epsilon),\sE)\) extending \(j_W[U]\).
	\end{itemize}
	Then \(U\text{-}\sum W_* = \textnormal{flip}_*(W\text{-}\sum U_*)\).
	\begin{proof}
		Assume towards a contradiction that \(U\text{-}\sum W_* \prec^{\Bbbk} \textnormal{flip}_*(W\text{-}\sum U_*)\). 
		The following identities are easily verified using \cref{SumPower}: 
		 \begin{align*}
		 	j_{U\text{-}\sum W_*} &= j_{W_*}^{M_U} \circ j_U & 	j_{\textnormal{flip}_*(W\text{-}\sum U_*)} &= j_{U_*}^{M_W} \circ j_W\\
		 	\id_{U\text{-}\sum W_*} &= (j_{W_*}^{M_U}(\id_U),\id_{W_*})& \id_{\textnormal{flip}_*(W\text{-}\sum U_*)} &= (\id_{U_*}, j_{U_*}^{M_W}(\id_W))
		 \end{align*}
		By  \cref{KOWOChar}, the assumption that \(U\text{-}\sum W_* \prec^{\Bbbk} \textnormal{flip}_*(W\text{-}\sum U_*)\) is equivalent to the existence of a \(1\)-internal comparison \[(k,h) : (M_{W_*}^{M_U},M_{U_*}^{M_W})\to N\] of \((j_{W_*}^{M_U} \circ j_U,  j_{U_*}^{M_W} \circ j_W)\) such that 
		\[k(j_{W_*}^{M_U}(\id_U),\id_{W_*}) \prec h(\id_{U_*}, j_{U_*}^{M_W}(\id_W))\]
		Therefore either \(k(j_{W_*}^{M_U}(\id_U)) < h(\id_{U_*})\) or \(k(\id_{W_*}) < h(j_{U_*}^{M_W}(\id_W))\), contradicting \cref{SumComparison}.
		
		A symmetric argument shows that we cannot have \(\textnormal{flip}_*(W\text{-}\sum U_*) \prec^{\Bbbk} U\text{-}\sum W_*\) either. Thus by the linearity of \((\mathscr B(\epsilon\times\delta),\prec^\Bbbk)\), we must have \(U\text{-}\sum W_* = \textnormal{flip}_*(W\text{-}\sum U_*)\), which proves the theorem.
	\end{proof}
\end{lma}

As an immediate consequence, we can prove \cref{Reciprocity1}:
\begin{proof}[Proof of \cref{Reciprocity1}]
	Let \(\alpha\) be the ordertype of the G\"odel order on \(\epsilon\times \delta\). Since the Ketonen order is linear on \(\mathscr B(\alpha)\), the isomorphic order \((\mathscr B(\epsilon\times \delta),\prec^\Bbbk)\) is also linear. Thus we can apply \cref{PreciseReciprocity} to conclude that \(U\text{-}\sum W_* = \text{flip}_*(W\text{-}\sum U_*)\). In particular, \(U\text{-}\sum W_*\cong W\text{-}\sum U_*\), so applying \cref{SumPower},
	\[j_{W_*}^{M_U}\circ j_U = j_{U\text{-}\sum W_*} =j_{W\text{-}\sum U_*} = j_{U_*}^{M_W} \circ j_W\]
	Thus \((j_{W_*}^{M_U},j_{U_*}^{M_W})\) is a comparison of \((j_U,j_W)\), as desired.
\end{proof}

Let us make some comments on this theorem. It is not immediately obvious from the definition that the linearity of the Ketonen order on \(\mathscr B(\lambda)\) implies the linearity of the Ketonen order on \(\mathscr B(\delta)\) for all ordinals \(\delta < \lambda^+\).\footnote{Note that if \(\kappa\) is regular, then for any \(n < \omega\), the collection of subsets of \(\kappa^n\) of ordertype less than \(\kappa^n\) forms a  \(\kappa\)-complete ideal; this is closely related to the Milner-Rado Paradox. Therefore for example if \(\kappa\) is \(2^\kappa\)-strongly compact, there is a \(\kappa\)-complete ultrafilter on \(\kappa^n\) that does not concentrate on a set of ordertype less than \(\kappa^n\). (It suffices that \(\kappa\) is measurable.) This suggests it may be nontrivial to reduce the linearity of \((\mathscr B(\kappa^2),\sE)\) to that of \((\mathscr B(\kappa),\sE)\) by a direct combinatorial argument.}

\begin{defn}
	Suppose \(\lambda\) is a cardinal.
	\begin{itemize}
		\item \(\text{UA}_{<\lambda}\) is the assertion that any pair of ultrapower embeddings of width less than \(\lambda\) have an internal ultrapower comparison. 
		\item \(\text{UA}_{\leq\lambda}\) is another way of writing \(\text{UA}_{<\lambda^+}\).
	\end{itemize}
\end{defn}

\begin{cor}
	Suppose \(\lambda\) is an infinite cardinal and the Ketonen order is linear on \(\mathscr B(\lambda)\). Then \(\textnormal{UA}_{\leq\lambda}\) holds. In particular, the Ketonen order is linear on \(\mathscr B(\delta)\) for all \(\delta < \lambda^+\).
	\begin{proof}
		Suppose \(U\) and \(W\) are ultrafilters on \(\lambda\). To see \(\textnormal{UA}_{\leq\lambda}\), it suffices to show that \((j_U,j_W)\) has a comparison. Since the Ketonen order \((\mathscr B(\lambda),\sE)\) is linear, so is \((\mathscr B(X),\prec^{\Bbbk})\) whenever \((X,\prec)\) is a wellorder of ordertype \(\lambda\).
		Since \(\lambda\) is an infinite cardinal, the G\"odel order on \(\lambda\times \lambda\) has ordertype \(\lambda\). Thus \((\mathscr B(\lambda\times\lambda),\prec^{\Bbbk})\) is linear, and so we can apply \cref{PreciseReciprocity} and the proof of \cref{Reciprocity1} to conclude that \((j_U,j_W)\) has a comparison.
	\end{proof}
\end{cor}

Surely with some extra work one can prove the following conjecture:
\begin{conj}
	If the Ketonen order is linear on countably complete incompressible ultrafilters, then the Ultrapower Axiom holds.
\end{conj}

The proof of \cref{LinKetThm} that we have given here uses \L o\'s's Theorem, which makes significant use of the Axiom of Choice. With care, however, the combinatorial content of \cref{LinKetThm}, namely \cref{PreciseReciprocity}, can actually be established in ZF + DC alone. This makes the following question seem interesting:
\begin{qst}
Assume \(\textnormal{AD} + V= L(\mathbb R)\). Is the Ketonen order linear?
\end{qst}
\chapter{The Generalized Mitchell Order}\label{GMOChapter}
\section{Introduction}
\subsection{The linearity of the generalized Mitchell order}
The topic of this section is the generalized Mitchell order, which is defined by
extending the definition of the Mitchell order to a broader class of objects:
\begin{defn}\index{Generalized Mitchell order}\index{Mitchell order!generalized}
The {\it generalized Mitchell order} is defined on countably complete
ultrafilters \(U\) and \(W\) by setting \(U\mo W\) if \(U\in M_{W}\).
\end{defn}
The main question we investigate here is to what extent this generalized order
is linear assuming the Ultrapower Axiom. Recall that UA implies the linearity of
the Mitchell order on normal ultrafilters (\cref{UAMO}). On the other hand the
generalized Mitchell order is obviously not a linear order on arbitrary
countably complete ultrafilters (\cref{NonlinearitySection}). The main theorem
of this chapter is the generalization of the linearity of the Mitchell order on
normal ultrafilters to normal fine ultrafilters:

\begin{defn}\index{\(P_\textnormal{bd}(\lambda)\) (bounded powerset)}
	For any ordinal \(\lambda\), the {\it bounded powerset of \(\lambda\)} is
	the set \(P_\text{bd}(\lambda) = \bigcup_{\alpha < \lambda} P(\alpha)\).
\end{defn}

\begin{repthm}{GCHLinear}[UA] Suppose \(\lambda\) is a cardinal such that
	\(2^{<\lambda} = \lambda\). Then the Mitchell order is linear on normal fine
	ultrafilters on \(P_\textnormal{bd}(\lambda)\).
\end{repthm}
This amounts to the most general form of the linearity of the Mitchell order on
normal fine ultrafilters that one could hope for (\cref{NIso}), except for the
cardinal arithmetic assumption on \(\lambda\) (which we dispense with much later
in \cref{ULinearity}).
\subsection{Outline of \cref{GMOChapter}}
We now outline the rest of this chapter.\\

\noindent {\sc\cref{FolkloreSection}.} This contains various folklore facts about large
cardinals and the generalized Mitchell order. None of the results here are due
to the author. We give a brief exposition of the theory of strong embeddings
(\cref{StrengthSection}) and supercompact embeddings (\cref{SCSection}) centered
around the relationship between these concepts and the generalized Mitchell
order. We also exposit the Kunen Inconsistency Theorem, which is closely related
to the wellfoundedness properties of the Mitchell order. Finally we establish
the basic order theoretic properties of the generalized Mitchell order,
especially its transitivity, wellfoundedness (\cref{MOWF}), and nonlinearity
(\cref{NonlinearitySection}).\\

\noindent {\sc\cref{DoddSection}.} This section introduces the notion of Dodd soundness.
This concept first arose in inner model theory, and our exposition is the first
to put it into a general context. We begin by giving a very simple definition of
Dodd soundness that will hopefully help the reader view it as a natural
refinement of supercompactness. We then prove the equivalence of this notion
with the definition of Dodd soundness from fine structure theory
(\cref{DoddEquiv}).

A theorem of Schlutzenberg \cite{Schlutzenberg} (stated as \cref{Schlutz} below)
shows that the Mitchell order is linear on Dodd sound ultrafilters in the
canonical inner models. We prove this theorem (\cref{DoddMO}) here under the
much weaker assumption of UA and by a completely different and much simpler
argument directly generalizing the proof of the linearity of the Mitchell order
on normal ultrafilters.\\

\noindent {\sc\cref{GeneralizedNormalSection}.} We finally turn to the Mitchell order on
normal fine ultrafilters, the natural generalization of normal ultrafilters
associated with supercompact cardinals. Our analysis proceeds by showing that
normal fine ultrafilters are isomorphic to Dodd sound ultrafilters, and then
citing the linearity of the Mitchell order on Dodd sound ultrafilters. To do
this, we introduce the notion of an {\it isonormal ultrafilter} and prove that
every normal fine ultrafilter is isomorphic to an isonormal ultrafilter
(\cref{IsoNormalThm}). The main difficulty is the ``singular case"
(\cref{SingularSolovaySection}) which amounts to generalizing Solovay's Lemma
\cite{Solovay} (proved as \cref{SolovayLemma}) to singular cardinals.
\cref{IsoSound} states that if \(2^{<\lambda} = \lambda\), then isonormal
ultrafilters on \(\lambda\) are Dodd sound. Putting these theorems together, we
obtain that under the Generalized Continuum Hypothesis, normal fine ultrafilters
are isomorphic to Dodd sound ultrafilters, yielding the main theorem of the
chapter (\cref{GCHLinear}), the linearity of the Mitchell order on normal fine
ultrafilters.
\section{Folklore of the generalized Mitchell order}\label{FolkloreSection}

\subsection{Strength and the Mitchell order}\label{StrengthSection}
The generalized Mitchell order is often viewed as a more finely calibrated
generalization of the concept of the strength of an elementary embedding. In
this subsection, we set down the basic theory of strength and discuss its
relationship with the Mitchell order.

\begin{defn}
Suppose \(M\) is a transitive class and \(\lambda\) is a cardinal. 
\begin{itemize}
\item An elementary embedding \(j: V\to M\) is {\it \(\lambda\)-strong} if \(P(\lambda)\subseteq M\).\index{Strength!\(\lambda\)-strong embedding}
\item An elementary embedding is {\it \({<}\lambda\)-strong} if
\(P_\text{bd}(\lambda)\subseteq M\).
\end{itemize}
\end{defn}

Notice that the property of being \(\lambda\)-strong depends only on \(M\). The
basic lemmas we prove about \(\lambda\)-strong embeddings almost all apply to
arbitrary inner models containing \(P(\lambda)\). (The embedding \(j\) just
comes along for the ride.)

Most authors define \(j\) to be \(\alpha\)-strong if \(V_\alpha\subseteq M\).
The definition used here is arguably preferable (if one is assuming the Axiom of
Choice and not assuming the Generalized Continuum Hypothesis). This is because
it is more expressive:

\begin{lma} Suppose \(j : V\to M\) is an elementary embedding. 
	\begin{itemize}
		\item If \(\alpha\) is an ordinal, then \(V_{\alpha + 1}\subseteq M\) if
		and only if \(j\) is \(\beth_\alpha\)-strong. 
		\item If \(\gamma\) is a limit ordinal, then \(V_{\gamma}\subseteq M\)
		if and only if \(j\) is \({<}\beth_\gamma\)-strong.\qed
	\end{itemize}
\end{lma}

It would be strange to define strong embeddings without defining strong
cardinals, so let us include the definition even though we will have little to
say about the concept:
\begin{defn}
	Suppose \(\kappa \leq\lambda\) are cardinals. Then \(\kappa\) is
	{\it\(\lambda\)-strong} if there is an inner model \(M\) and a
	\(\lambda\)-strong elementary embedding \(j : V\to M\) such that
	\(\textsc{crt}(j) = \kappa\) and \(j(\kappa) > \lambda\). \(\kappa\) is {\it
	strong} if \(\kappa\) is \(\lambda\)-strong for all \(\lambda\).
\end{defn}

The requirement that \(j(\kappa) > \lambda\) above is not actually necessary as
a consequence of \cref{KunenInconsistency}. We use standard notation for
hereditary cardinality:
\begin{defn}
	If \(x\) is a set, \(\text{tc}(x)\) denotes the smallest transitive set
	\(y\) with \(x\subseteq y\). The {\it hereditary cardinality} of \(x\) is
	the cardinality of \(\text{tc}(x)\). For any cardinal \(\lambda\),
	\(H(\lambda)\) denotes the collection of sets of hereditary cardinality less
	than \(\lambda\).\index{\(H(\lambda)\)}
\end{defn}

\begin{lma}
	For any infinite cardinal \(\lambda\), 
	\begin{itemize}
		\item \(H(\lambda^+)\) is a transitive set.
		\item \(H(\lambda^+)\) is bi-interpretable with \(P(\lambda)\).
		\item \(H(\lambda)\) is bi-interpretable with
		\(P_\textnormal{bd}(\lambda)\).\qed
	\end{itemize}
\end{lma}
The bi-interpretability of \(H(\lambda^+)\) and \(P(\lambda)\) yields the
following lemma:
\begin{lma}
An embedding \(j: V\to M\) is \({<}\lambda\)-strong if and only if
\(H(\lambda)\subseteq M\) and \(\lambda\)-strong if and only if
\(H(\lambda^+)\subseteq M\).\qed
\end{lma}

\begin{defn}
The {\it strength} of an elementary embedding \(j : V\to M\), denoted
\(\textsc{str}(j)\), is the largest cardinal \(\lambda\) such that \(j\) is
\({<}\lambda\)-strong.\index{Strength!of an embedding}
\end{defn}

The following fact specifies exactly which powersets are contained in the target
model of an elementary embedding in terms of its strength:
\begin{lma}\label{StrengthInvariance}
Suppose \(j : V\to M\) is an elementary embedding and \(\lambda\) is a cardinal.
Then the following are equivalent:
\begin{enumerate}[(1)]
\item \(\textsc{str}(j) = \lambda\).
\item For all \(X\in M\), \(P(X)\subseteq M\) if and only if \(|X|^M <
\lambda\).\qed
\end{enumerate}
\end{lma}

The main limitation on the strength of an elementary embedding is known as the
{\it Kunen Inconsistency Theorem} \cite{Kunen}:
\begin{thm}[Kunen] Suppose \(j : V\to M\) is a nontrivial elementary embedding
and \(\lambda\) is the first fixed point of \(j\) above \(\textsc{crt}(j)\).
Then \(\textsc{str}(j)\leq\lambda\).\qed
\end{thm}
We prove this and other related facts in \cref{KunenSection}.

The basic relationship between strength and the Mitchell order is given by the
following two lemmas:
\begin{lma}\label{MOStrength}
Suppose \(U\) and \(W\) are countably complete ultrafilters and \(U\mo W\). Then
\(M_W\) is \(\lambda\)-strong where \(\lambda\) is the cardinality of the
underlying set \(X\) of \(U\). In fact, \(P(X)\subseteq M_W\).
\begin{proof}
Clearly \(X\in M_W\) since \(X\in U\in M_W\). It suffices to show that
\(P(X)\subseteq M_W\). Fix \(A\subseteq X\), and we will show \(A\in M_W\).
Since \(U\) is an ultrafilter, either \(A\in U\) or \(X\setminus A\in U\). If
\(A\in U\), then \(A\in U\in M_W\), so \(A\in M_W\). If \(X\setminus A\in U\),
then similarly \(X\setminus A\in M_W\), and since \(X\in M_W\), it follows that
\(A = X\setminus (X\setminus A)\in M_W\). Therefore in either case, \(A\in
M_W\).
\end{proof}
\end{lma}

\begin{lma}
Suppose \(W\) is a countably complete ultrafilter and \(j_W\) is
\(2^\lambda\)-strong. Then for any countably complete ultrafilter \(U\) on
\(\lambda\),  \(U\mo W\).
\begin{proof}
Since \(U\subseteq P(\lambda)\), \(U\in H((2^\lambda)^+)\subseteq M_W\).
\end{proof}
\end{lma}
This strength requirement implicit in the definition of the generalized Mitchell
order may seem somewhat unnatural. What if one modified the Mitchell order,
considering for example the {\it amenability relation} defined on countably
complete ultrafilters by setting \(U\A W\) if and only if \(U\) concentrates on
\(M_W\) and \(U\cap M_W\in M_W\)? Such modified Mitchell orders are the subject
of \cref{InternalSection}.

For the time being, we must point out some irritating properties of the
generalized Mitchell order that suggest that in some sense it may be a little
bit {\it too} general. The issue is that the definition of \(U\mo W\) above has
a strong dependence on the choice of the underlying set of \(U\). For example,
if \(W\) is nonprincipal, then the following hold:
\begin{itemize}
\item There is a principal ultrafilter \(D\) on an ordinal such that \(D\not\mo
W\).
\item There is a set \(x\) such that the principal ultrafilter \(\{\{ x
\}\}\not\mo W\).
\end{itemize}
For the first bullet point, let \(\lambda\) be the strength of \(j_W\), and let
\(D\) be any principal ultrafilter on \(\lambda\). For the second bullet point,
let \(x\) be any set that does not belong to \(M_W\). 

These silly counterexamples suggest that the generalized Mitchell order is only
a well-behaved relation on a restricted class of ultrafilters. Recall that for
any ultrafilter \(U\) on a set \(X\), \(\lambda_U\) is defined to be the least
cardinality of a set in \(U\), and \(U\) is said to be uniform if \(|X| =
\lambda_U\). Hereditary uniformity is a strengthening of uniformity:
\begin{defn}
An ultrafilter \(U\) on a set \(X\) is {\it hereditarily uniform} if
\(\lambda_U\) is the hereditary cardinality of \(X\).\index{Hereditarily uniform
ultrafilter}
\end{defn}
Any ultrafilter \(U\) is isomorphic to a hereditarily uniform ultrafilter since
in fact \(U\) is isomorphic to an ultrafilter on \(\lambda_U\)
(\cref{UniformIso}). The following lemma argues that the generalized Mitchell
order is a reasonable relation on the class of hereditarily uniform
ultrafilters:

\begin{lma}\label{HeredLemma0}
Suppose \(U'\RK U\mo W\) are countably complete ultrafilters. Let \(X\) and
\(X'\) be the underlying sets of \(U\) and \(U'\), and assume \(X'\in M_W\) and
\(M_W\) satisfies \(|X'|\leq|X|\). Then \(U'\mo W\) and \(M_W\) satisfies
\(U'\RK U\). If \(U'\cong U\), then \(M_W\) satisfies \(U'\cong U\).
\end{lma}

\begin{lma}\label{HeredLemma}
Suppose \(U'\RK U\mo W\) are countably complete ultrafilters and \(U'\) is
hereditarily uniform. Then \(U'\mo W\) and \(M_W\) satisfies \(U'\RK U\).  If
\(U'\cong U\), then \(M_W\) satisfies \(U'\cong U\). In particular, the
restriction of the generalized Mitchell order to hereditarily uniform
ultrafilters is isomorphism invariant.
\end{lma}

\cref{HeredLemma0} and \cref{HeredLemma} follow from a fact that is both more
general and easier to prove:
\begin{lma}\label{ZFCStrength}
	Suppose \(M\) is an inner model of \textnormal{ZFC}, \(\lambda\) is a
	cardinal, and \(X\in M\) is a set of cardinality \(\lambda\) such that
	\(P(X)\subseteq M\). 
	\begin{itemize}
		\item For any set \(Y\in M\) such that \(M\vDash |Y| \leq |X|\),
		\(P(Y)\subseteq M\).
		\item For any set \(Y\in M\) such that \(M\vDash  |Y| \leq |X|\),
		\(P(X\times Y)\subseteq M\).
		\item For any set \(Y\in M\) such that \(M\vDash  |Y| \leq |X|\), every
		function from \(X\) to \(P(Y)\) belongs to \(M\).
		\item \(P(\lambda)\subseteq M\).
		\item Every set of hereditary cardinality at most \(\lambda\) belongs to
		\(M\) and has hereditary cardinality at most \(\lambda\) in \(M\).\qed
	\end{itemize}
\end{lma}
The bullet points are arranged in such a way that the reader should have no
trouble proving each one in turn.\footnote{It is likely, however, that the
second bullet-point cannot be established if \(M\) is not assumed to satisfy the
Axiom of Choice.}

\begin{proof}[Proof of \cref{HeredLemma0}] Fix \(f : X\to X'\) such that
	\(f_*(U) = U'\). By \cref{ZFCStrength}, \(f\in M_W\), and hence \(U'=
	f_*(U)\in M_W\). Moreover \(f\) witnesses \(U'\RK U\) in \(M_W\). Finally if
	\(U'\cong U\), then this is also witnessed by some \(g\in M_W\).
\end{proof}
\begin{proof}[Proof of \cref{HeredLemma}] By \cref{ZFCStrength}, the underlying
	set of \(U'\) belongs to \(M_W\) and has hereditary cardinality at most
	\(\lambda_{U'}\leq \lambda_U \leq |X|\) in \(M_W\), so the lemma follows
	from \cref{HeredLemma0}.
\end{proof}

\subsection{Supercompactness and the Mitchell order}\label{SCSection}
We now turn to a concept that is more pertinent to this dissertation than
strength: supercompactness.

\begin{defn}\index{Supercompactness!\(X\)-supercompact embedding}
Suppose \(M\) is a transitive class and \(X\) is a set. An elementary embedding
\(j: V\to M\) is {\it \(X\)-supercompact} if \(j[X]\in M\).
\end{defn}

The following lemma allows us to focus solely on the case of
\(\lambda\)-supercompact embeddings for \(\lambda\) a cardinal:
\begin{lma}\label{XSC}
Suppose \(X\) and \(Y\) are sets such that \(|X| = |Y|\). Then an elementary
embedding \(j : V\to M\) is \(X\)-supercompact if and only if it is
\(Y\)-supercompact. In particular, \(j\) is \(X\)-supercompact if and only if
\(j\) is \(|X|\)-supercompact.
\begin{proof}
Suppose \(j\) is \(X\)-supercompact and \(f : X\to Y\) is a surjection. Then
\[j(f)[j[X]] = j[Y]\] so \(j\) is \(Y\)-supercompact.
\end{proof}
\end{lma}

\begin{defn}\index{Supercompactness!\(\lambda\)-supercompact cardinal}\index{Supercompact cardinal}
	Suppose \(\kappa\leq \lambda\) are cardinals. Then \(\kappa\) is {\it
	\(\lambda\)-supercompact} if there is a \(\lambda\)-supercompact embedding
	\(j : V\to M\) such that \(\textsc{crt}(j) = \kappa\) and \(j(\kappa) >
	\lambda\); \(\kappa\) is supercompact if \(\kappa\) is
	\(\lambda\)-supercompact for all cardinals \(\lambda \geq \kappa\).  
\end{defn}

The results of this dissertation (\cref{PathologicalSection}) single out a class
of ultrapower embeddings that are just shy of \(\lambda\)-supercompact, so the
following is an important definition:

\begin{defn}\index{Supercompactness!\({<}\lambda\)-supercompact embedding}
Suppose \(\lambda\) is a cardinal. An elementary embedding \(j: V\to M\) is {\it
\({<}\lambda\)-supercompact} if \(j\) is \(\delta\)-supercompact for all
cardinals \(\delta < \lambda\).
\end{defn}

The definition of supercompactness is motivated by its relationship with the
closure of \(M\) under \(\lambda\)-sequences:

\begin{lma}\label{SupercompactClosure}
Suppose \(j : V\to M\) is an elementary embedding and \(\lambda\) is a cardinal.
\begin{enumerate}[(1)]
\item \(j\) is \(\lambda\)-supercompact if and only if \(j\restriction
\lambda\in M\).
\item If \(j\) is \(\lambda\)-supercompact, then \(j\) is \(\lambda\)-strong.
\item If \(j\) is \(\lambda\)-supercompact, then \(j[X]\in M\) for all \(X\) of
cardinality \(\lambda\).
\item If \(j\) is \(\lambda\)-supercompact and \(M = H^M(j[V]\cup S)\) for some
\(S\subseteq M\) such that \(S^\lambda\subseteq M\), then \(M^\lambda\subseteq
M\).
\end{enumerate}
\begin{proof}
For (1), note that \(j\restriction \lambda\) is the inverse of the transitive
collapse of \(j[\lambda]\). 

For (2), suppose \(A\subseteq \lambda\). Then \(A = (j\restriction
\lambda)^{-1}[j(A)]\), so since \(j\restriction \lambda\) and \(j(A)\) both
belong to \(M\), so does \(A\). 

(3) is immediate from \cref{XSC}. 

For (4), fix \(\langle x_\alpha : \alpha < \lambda\rangle\in M^\lambda\). Fix
\(\langle f_\alpha : \alpha < \lambda\rangle\) and \(\langle a_\alpha : \alpha <
\lambda\rangle\in S^\lambda\) such that \(x_\alpha = j(f_\alpha)(a_\alpha)\) for
all \(\alpha < \lambda\). The function \(G : j[\lambda]\to M\) defined by
\(G(j(\alpha)) = j(f_\alpha)\) belongs to \(M\) by (3), since \[G =
j[\{(\alpha,f_\alpha) : \alpha < \lambda\}]\] Therefore the sequence \(\langle
j(f_\alpha) : \alpha < \lambda\rangle\) can be computed from \(G\) and
\(j\restriction \lambda\): \[j(f_\alpha) = G\circ(j\restriction
\lambda)(\alpha)\] Since both \(G\) and \(j\restriction \lambda\) belong to
\(M\) by (1), \(\langle j(f_\alpha) : \alpha < \lambda\rangle\in M\). Finally,
\[\langle x_\alpha : \alpha < \lambda\rangle = \langle j(f_\alpha)(a_\alpha) :
\alpha < \lambda\rangle\] can be computed from \(\langle j(f_\alpha) : \alpha <
\lambda\rangle\) and \(\langle a_\alpha : \alpha < \lambda\rangle\). Both these
sequences belong to \(M\), since \(\langle a_\alpha : \alpha < \lambda\rangle\in
S^\lambda\subseteq M\), so \(\langle x_\alpha : \alpha < \lambda\rangle\in M\),
as desired.
\end{proof}
\end{lma}

For the purposes of this dissertation, the most relevant corollary of
\cref{SupercompactClosure} is its application to ultrapower embeddings:

\begin{cor}\label{UltrapowerSC}
An ultrapower embedding \(j : V\to M\) is \(\lambda\)-supercompact if and only
if \(M^\lambda\subseteq M\).
\begin{proof}
Fix \(a\in M\) such that \(M = H^M(j[V]\cup \{a\})\). The corollary follows from
applying \cref{SupercompactClosure} (4) in the case \(S = \{a\}\).
\end{proof}
\end{cor}

We can make good use of \cref{UltrapowerSC} since it is always possible to
derive a \(\lambda\)-supercompact ultrapower embeddings from a
\(\lambda\)-supercompact embedding:
\begin{lma}\label{DerivedSC}
	Suppose \(j : V\to M\) is an \(X\)-supercompact embedding,
	\(V\stackrel{i}\longrightarrow N\stackrel{k}\longrightarrow M\) are
	elementary embeddings, \(k\circ i = j\), and \(j[X]\in k[N]\). Then \(i\) is
	\(X\)-supercompact and \(k(i[X]) = j[X]\). In particular, letting \(\lambda
	= |X|\), \(k\restriction \lambda + 1\) is the identity.
	\begin{proof}
		Fix \(S\in M\) such that \(k(S) = j[X]\). Then \[S = k^{-1}[k(S)] =
		k^{-1}[j[X]] = k^{-1}\circ j[X] = i[X]\] Thus \(i[X]= S \in M\), so
		\(i\) is \(X\)-supercompact, and moreover, \(k(i[X])= k(S) = j[X]\).
		
		Since \(k(i[X]) = j[X]\), the argument of \cref{XSC} shows
		\(k(i[\lambda]) = j[\lambda]\). But then if \(\alpha \leq \lambda\),
		\(k(\alpha) = k(\text{ot}(i[\lambda]\cap i(\alpha)))
		=\text{ot}(k(i[\lambda])\cap k(i(\alpha))) = \text{ot}(j[\lambda]\cap
		j(\alpha)) = \alpha\).
	\end{proof}
\end{lma}

\begin{defn}
	The {\it supercompactness} of an elementary embedding \(j : V\to M\) is the
	least cardinal \(\lambda\) such that \(j\) is not \(\lambda\)-supercompact.
\end{defn}
Which cardinals are the supercompactness of an elementary embedding? Which are
the supercompactness of an ultrapower embedding? This turns out to be a major
distinction:
\begin{prp}
	Suppose \(\lambda\) is a singular cardinal and \(j : V\to M\) is an
	elementary embedding such that \(M^{<\lambda}\subseteq M\). Then
	\(M^\lambda\subseteq M\).\qed
\end{prp}
Thus the supercompactness of an ultrapower embedding is always regular, while it
is easy to see this must fail for arbitrary embeddings if there is a
\(\kappa^+\omega\)-supercompact cardinal. An important point is that if the
cofinality of \(\lambda\) is small, \(\lambda\)-supercompactness is equivalent
to \(\lambda^+\)-supercompactness:
\begin{lma}\label{SmallCfCompact}\label{SmallCfSC}
	Suppose \(\lambda\) is a cardinal, \(j : V\to M\) is elementary embedding,
	and \(\kappa = \textsc{crt}(j)\). If \(j\) is \(\lambda\)-supercompact, then
	\(j\) is \(\lambda^{<\kappa}\)-supercompact.
	\begin{proof}
		Assume \(j[\lambda]\in M\), and we will show that
		\(j[P_\kappa(\lambda)]\in M\).  Note that for \(\sigma\in
		P_\kappa(\lambda)\), \(j(\sigma) = j[\sigma]\). Thus
		\[j[P_\kappa(\lambda)] = \{ j[\sigma] : \sigma\in P_\kappa(\lambda)\} =
		P_\kappa(j[\lambda])\] One consequence of this is that
		\(P_\kappa(j[\lambda])\subseteq M\), since
		\(j[P_\kappa(\lambda)]\subseteq M\), and therefore
		\(P_\kappa(j[\lambda]) = (P_\kappa(j[\lambda]))^M\in M\). It follows
		that \(j[P_\kappa(\lambda)] \in M\), as desired.
	\end{proof}
\end{lma}
It follows for example that the supercompactness of an elementary embedding is
never the successor of a singular cardinal \(\gamma\) of countable cofinality,
since \(\gamma^\omega = \gamma^+\). This is an important component in the proof
of Kunen's Inconsistency Theorem (\cref{KunenInconsistency}).

We now begin to examine the relationship between supercompactness and the
Mitchell order, which turns out to be central to the rest of this dissertation.
The key point is that if \(U\mo W\), then the supercompactness of \(M_W\)
determines the extent to which the ultrapower of \(M_W\) by \(U\) is correctly
computed by \(M_W\).
\begin{lma}\label{MOFactor}
	Suppose \(U\mo W\) are countably complete ultrafilters.  Then there is a
	unique elementary embedding \(k : (M_U)^{M_W}\to j_U(M_W)\) such that
	\(k\circ (j_U)^{M_W} = j_U\restriction M_W\) and \(k(\id_U^{M_W}) = \id_U\).
	Let \(X\) be the underlying set of \(U\). Then \(k \restriction
	j_U((2^\lambda)^{M_W}) + 1\) where \(\lambda = |X|\).
	\begin{proof}
		Since \(P(X)\subseteq M_W\), \(U\) is the ultrafilter derived from
		\(j_U\restriction M_W\) using \(\id_U\). Thus there is a unique factor
		embedding \(k : (M_U)^{M_W}\to j_U(M_W)\) such that \(k\circ (j_U)^{M_W}
		= j_U\restriction M_W\) and \(k(\id_U^{M_W}) = \id_U\). This establishes
		the first part of the lemma.
		
		As for the second part, since \(U\mo W\), we have \(P(X)\subseteq M_W\)
		and hence by \cref{ZFCStrength}, \(P(\lambda)\subseteq M_W\) and every
		function from \(X\) to \(P(\lambda)\) belongs to \(M_W\). It follows
		that \(j_U(P(\lambda))\subseteq \text{ran}(k)\): if \(A\in
		j_U(P(\lambda))\), then \(A = j_U(f)(\id_U)\) for some \(f : X\to
		P(\lambda)\), and therefore \[A =
		k(j_U^{M_W}(f)(\id_U^{M_W}))\in\text{ran}(k)\] Since there is a
		surjection \(g: P(\lambda)\to (2^\lambda)^{M_W}\) in \(M_W\),
		\[j_U(g)[j_U(P(\lambda))] =  j_U((2^\lambda)^{M_W})\subseteq
		\text{ran}(k)\] Moreover \(j_U((2^\lambda)^{M_W})\in j_U[M_W]\subseteq
		\text{ran}(k)\). Thus \(j_U((2^\lambda)^{M_W}) + 1 \subseteq
		\text{ran}(k)\), or in other words, \(k \restriction
		j_U((2^\lambda)^{M_W}) + 1 = \text{id}\).
	\end{proof}
\end{lma}
We will refer to the embedding of \cref{MOFactor} as a {\it factor
embedding}.\index{Factor embedding}

\begin{lma}\label{MOSuper2}
	Suppose \(U\) and \(W\) are countably complete ultrafilters with \(U\mo W\).
	Let \(X\) be the underlying set of \(U\), let \(\lambda = |X|\) and let
	\(\delta = ((2^{\lambda})^+)^{M_W}\). Then \[j^{M_W}_U\restriction
	H^{M_W}(\delta) = j_U\restriction H^{M_W}(\delta)\]
	\begin{proof}
		Let \(k : (M_U)^{M_W}\to j_U(M_W)\) be the factor embedding with
		\(k\circ (j_U)^{M_W} = j_U\restriction M_W\) and \(k(\id_U^{M_W}) =
		\id_U\). Then \cref{MOFactor} implies \(k\restriction
		j_U^{M_W}(\delta)\) is the identity, and therefore \(k\restriction
		j_U^{M_W}(H^{M_W}(\delta))\) is the identity. Now
		\[j^{M_W}_U\restriction H^{M_W}(\delta) = (k\restriction
		j_U^{M_W}(H^{M_W}(\delta)))\circ (j^{M_W}_U\restriction H^{M_W}(\delta))
		= j_U\restriction H^{M_W}(\delta)\qedhere\]
	\end{proof}
\end{lma}

Our next proposition, \cref{MOSuper}, suggests that the Mitchell order on
ultrafilters be seen as a generalization of supercompactness that asks for one
ultrapower \(M_W\) how much it can see of another embedding \(j_U\). (On this
view supercompactness is the special case in which we ask how much of \(j_U\) is
seen by \(M_U\) itself.)

\begin{prp}\label{MOSuper}
Suppose \(U\) and \(W\) are countably complete ultrafilters. Let \(X\) be the
underlying set of \(U\), let \(\lambda = |X|\) and let \(\delta =
((2^{\lambda})^+)^{M_W}\). Then the following are equivalent:
\begin{enumerate}[(1)]
	\item \(U\mo W\).
	\item \(j_U\restriction H^{M_W}(\delta)\in M_W\).
	\item \(j_U\restriction P(\lambda)\in M_W\).
	\item  \(j_U\restriction P(X) \in M_W\).
\end{enumerate}
\begin{proof}
{\it (1) implies (2).} Immediate from \cref{MOSuper2}.

{\it (2) implies (3).} Immediate since \(P(\lambda)\subseteq H^{M_W}(\delta)\).

{\it (2) implies (3).} This is probably clear enough (and in any case, (1)
implies (4) is easy), but let us just make sure. By \cref{ZFCStrength}, \(|X|^M
= \lambda\). Let \(\rho : \lambda\to X\) be a surjection in \(M_W\). For \(A\in
P(X)\), \[j_U(A) = j_U(\rho)[j_U(\rho^{-1}[A])]\]

{\it (3) implies (1).} If \(j_U\restriction P(X)\) belongs to \(M_W\), then \(U
= \{A\subseteq X : \id_U\in j_U(A)\}\) belongs to \(M_W\) as well.
\end{proof}
\end{prp}

Given \cref{MOSuper2}, it is reasonable to wonder whether the whole embedding
\(j_U\restriction M_W\) might be correctly computed by \(M_W\) as well; that is,
perhaps the factor embedding \(k\) is always trivial. We provide a
counterexample in \cref{CummingsExample}.\footnote{This counterexample also
shows that in the context of  \cref{MOFactor}, the lower bound given there on
\(\textsc{crt}(k)\) can be tight in the sense that (consistently) one can have
\[\textsc{crt}(k) = j_U\left((2^\lambda)^{M_W}\right)^{+(M_U)^{M_W}}\]} This is
equivalent to the supercompactness of \(j_W\), a phenomenon we exploit later:

\begin{prp}\label{USupercompact}
	Suppose \(U\mo W\) are countably complete ultrafilters. Then the following
	are equivalent:
	\begin{enumerate}[(1)]
		\item \((j_U)^{M_W} = j_U\restriction M_W\).
		\item \(j_W\) is \(\lambda_U\)-supercompact.
	\end{enumerate}
	\begin{proof}
		{\it (1) implies (2):} Let \(k : (M_U)^{M_W}\to j_U(M_W)\) be the factor
		embedding of \cref{MOFactor}, with \(k\circ j_U^{M_W} = j_U\restriction
		M_W\) and \(k(\id_U^{M_W}) = \id_U\). Since \((j_U)^{M_W} =
		j_U\restriction M_W\), we have that \(k : j_U(M_W)\to j_U(M_W)\) and
		\(k\circ j_U\circ j_W = j_U\circ j_W\). Hence by the basic theory of the
		Rudin-Keisler order (\cref{RKRigid}), \(k\) is the identity.
		
		It follows in particular that \(j_U(j_W)(\id_U)\in \text{ran}(k)\). Fix
		\(f:X\to M_W\) in \(M_W\) such that \[k(j_U^{M_W}(f)(\id^{M_W}_U)) =
		j_U(j_W)(\id_U)\] Thus \(j_U(f)(\id_U) = j_U(j_W)(\id_U)\), so by \L
		o\'s's Theorem, there is a set \(A\in U\) such that \(f\restriction A =
		j_W\restriction A\). Since \(P(X)\subseteq M_W\), \(A\in M_W\), and
		hence \(j_W\restriction A = f\restriction A\in M_W\). In particular,
		\(j_W[A]\in M_W\), so \(j_W\) is \(A\)-supercompact. By \cref{XSC},
		\(j_W\) is \(|A|\)-supercompact, and since \(\lambda_U\leq |A|\), it
		follows that \(j_W\) is \(\lambda_U\)-supercompact.
		
		{\it (2) implies (1):} Obvious.
	\end{proof}
\end{prp}

Down the line (\cref{UltrapowerCorrectness}) we will show that under UA,
whenever \(U\mo W\), in fact \(j_W\) is \(\lambda_U\)-supercompact (and in fact
it suffices that \(P(\lambda_U)\subseteq M_W\)), and thus \(j_W^{M_U} =
j_W\restriction M_U\). For now, let us mention a generalization of
\cref{USupercompact}, whose proof we omit:
\begin{prp}
	Suppose \(U\) and \(W\) are countably complete ultrafilters such that \(U\)
	concentrates on a set in \(M_W\). The following are equivalent:
	\begin{enumerate}[(1)]
		\item \(j_{U\cap M_W}^{M_W} = j_U\restriction M_W\)
		\item There is a function \(f\in M_W\) such that \(f \restriction A =
		j_W\restriction A\) for some \(A\in U\).
		\item For all \(f: I\to M_W\) where \(I\in U\), there is some \(g\in
		M_W\) such that \(g\restriction A = f\restriction A\) for some \(A\in
		U\).\qed
	\end{enumerate}
\end{prp}

We finish this section with a restriction on the supercompactness of an
ultrafilter:
\begin{prp}\label{UFSuperBound}
Suppose \(U\) is an ultrafilter and \(j_U\) is \(\lambda_U^+\)-supercompact.
Then \(U\) is principal.
\end{prp}

We use the following lemma:

\begin{lma}\label{Ineqs}
Suppose \(j : V\to M\) is an elementary embedding that is discontinuous at the
infinite cardinal \(\lambda\). Let \(\lambda_* = \sup j[\lambda]\). Then
\[\lambda^+ \leq \lambda_*^{+{M}} < j(\lambda)^{+M} = j(\lambda^+)\] If \(j\) is
continuous at \(\lambda^+\), then \(j(\lambda^+)\) is a singular ordinal of
cofinality \(\lambda^+\), so \(j(\lambda^+) < j(\lambda)^+\).
\begin{proof}
We first show that \(\lambda^+ \leq \lambda_*^{+M}\). Suppose \(\alpha <
\lambda^+\). Let \(\prec\) be a wellorder of \(\lambda\) such that
\(\text{ot}(\prec) = \alpha\). Then \({\prec}_* = j(\prec)\restriction
\lambda_*\) is a wellorder of \(\lambda_*\) and \(j\) restricts to an
order-preserving embedding from \((\lambda,\prec)\) into
\((\lambda_*,{\prec}_*)\). Therefore \[\alpha \leq
\text{ot}(\lambda_*,{\prec}_*) < \lambda_*^{+{M}}\] The final inequality follows
from the fact that \((\lambda_*,{\prec}_*)\) belongs to \(M\). Since \(\alpha <
\lambda^+\) was arbitrary, it follows that \(\lambda^+ \leq \lambda_*^{+{M}}\).

To prove \(\lambda_*^{+{M}} < j(\lambda)^{+M}\), it is of course enough to show
\(\lambda_*^{+{M}}\leq j(\lambda)\). But \(j(\lambda)\) is a cardinal of \(M\)
that is greater than \(\lambda_*\), and hence \(\lambda_*^{+{M}}\leq
j(\lambda)\).

Finally, assume that \(j\) is continuous at \(\lambda^+\). Obviously
\(j(\lambda^+)\) has cofinality \(\lambda^+\), but the point is that this
implies \(j(\lambda^+)\) is {\it singular}, since the inequalities above show
\(\lambda^+ < j(\lambda^+)\). We can therefore conclude \(j(\lambda)^{+M} <
j(\lambda)^+\): obviously \(j(\lambda)^{+M} \leq j(\lambda)^+\), but the point
is that equality cannot hold since \(j(\lambda^+)\) is singular and
\(j(\lambda)^+\) is regular. 
\end{proof}
\end{lma}

\begin{proof}[Proof of \cref{UFSuperBound}] Let \(\lambda = \lambda_U\). Without
loss of generality, we may assume that \(U\) is a uniform ultrafilter on
\(\lambda\) and \(\lambda\) is infinite. Thus \(j_U\) is discontinuous at
\(\lambda\). Assume towards a contradiction that \(j_U[\lambda^+]\in M_U\). By
\cref{Ineqs}, \(j_U(\lambda^+) > \lambda^+\). But by \cref{UFContinuity},
\(j_U\) is continuous at \(\lambda^+\), and therefore \(j_U[\lambda^+]\in M_U\)
is a cofinal subset of \(j_U(\lambda^+)\) of ordertype \(\lambda^+\). Hence
\(\textnormal{cf}^{M_U}(j_U(\lambda)^{+M_U}) = \lambda^+< j_U(\lambda^+)\), and
this contradicts that \(j_U(\lambda^+)\) is regular in \(M_U\).
\end{proof}

\subsection{The Kunen Inconsistency}\label{KunenSection}
The story of the Kunen Inconsistency Theorem is often cast as a cautionary tale
with the moral that a large cardinal hypothesis may turn out to be false for
nontrivial combinatorial reasons:
\begin{thm}[Kunen]\index{Kunen Inconsistency Theorem}
	There is no nontrivial elementary embedding from the universe to itself.\qed
\end{thm}
A more pragmatic perspective is to view the Kunen Inconsistency as a proof
technique, providing at least some constraint on the elementary embeddings a
large cardinal theorist is bound to analyze. Examples pervade this work, but to
pick the closest one, the Kunen Inconsistency will form a key component of the
proof of the wellfoundedness of the Mitchell order in \cref{MOWFSection}.  Since
our applications of Kunen's theorem will require the basic concepts from the
proof (especially the notion of a critical sequence), we devote this subsection
for a brief exposition of this topic.

We first give a proof of a version of Kunen's inconsistency Theorem that is due
to Harada. (Another writeup of this proof appears in Kanamori's textbook
\cite{Kanamori2}.) The methods are purely ultrafilter-theoretic and very
much in the spirit of this dissertation:
\begin{prp}[Kunen]\label{KunenPrp}
	Suppose \(j : V\to M\) is an elementary embedding and \(\eta\) is a strong
	limit cardinal such that  \(j\) is \(\eta\)-supercompact and \(j[\eta]
	\subseteq \eta\). Then \(j \restriction \eta = \textnormal{id}\).
	\begin{proof}
		Assume the proposition holds for all \(\bar \eta <\eta\). 
		
		If \(\eta\) is has uncountable cofinality, then there is an
		\(\omega\)-closed unbounded set of \(\bar \eta < \eta\) such that
		\(j[\bar \eta]\subseteq \bar \eta\). Therefore \(j \restriction \bar
		\eta=\text{id}\) for unboundedly many \(\bar \eta < \eta\), so
		\(j\restriction \eta = \text{id}\). 
		
		Assume instead that \(\eta\) has countable cofinality. Then \(j\) is
		continuous at \(\eta\), so since \(j[\eta]\subseteq \eta\), we have
		\(j(\eta) = \eta\). We essentially reduce to the case that \(j\) is the
		ultrapower of the universe by an ultrafilter \(U\) on \(P(\eta)\). Let
		\(U\) be the ultrafilter on \(P(\eta)\) derived from \(j\) using
		\(j[\eta]\). Let \(k : M_U\to M\) be the factor embedding. Then
		\(j[\eta]\in k[M_U]\), so by our analysis of the supercompactness of
		derived embeddings (\cref{DerivedSC}), \(j_U\) is \(\eta\)-supercompact
		and \(k\restriction \eta + 1 = \text{id}\). If follows that \(j_U(\eta)
		= \eta\). Moreover, if we show \(j_U\restriction \eta = \text{id}\) then
		we can conclude \(j\restriction \eta = \text{id}\). In fact, we will
		show that \(U\) is principal.
		
		By \cref{SmallCfCompact}, \(j_U\) is actually
		\(\eta^\omega\)-supercompact. Since \(\eta\) is a strong limit cardinal
		of countable cofinality, \(\eta^\omega = 2^\eta\). Thus \(j_U\) is
		\(2^\eta\)-supercompact. Recall that \cref{UFSuperBound} states that if
		\(j_U\) is \(\lambda_U^+\)-supercompact, then \(U\) is principal. Thus
		to show \(U\) is principal, it suffices to show that \(\lambda_U <
		2^\eta\).
		
		Since \(U\) is an ultrafilter on \(P(\eta)\), \(\lambda_U \leq |P(\eta)|
		\leq 2^\eta\), so in fact, we need only show \(\lambda_U \neq 2^\eta\).
		Since \(U\) is isomorphic to a uniform ultrafilter on \(\lambda_U\),
		\(j_U\) is discontinuous at \(\lambda_U\), and in particular
		\(j_U(\lambda_U) \neq \lambda_U\). On the other hand, since \(M_U\) is
		closed under \(2^\eta\)-sequences, \((2^\eta)^{M_U} = 2^\eta\), and
		hence \[j_U(2^\eta) =  (2^{j_U(\eta)})^{M_U} = (2^\eta)^{M_U} = 2^\eta\]
		Since \(\lambda_U\) is moved by \(j_U\) while \(2^\eta\) is fixed,
		\(\lambda_U \neq 2^\eta\), and hence \(\lambda_U < 2^\eta\), as desired.
	\end{proof}
\end{prp}

Suppose \(j : V\to M\) is an elementary embedding with critical point
\(\kappa\). Let \(\lambda\) be the first ordinal above \(\kappa\) such that
\(j[\lambda]\subseteq \lambda\), the first ordinal at which one might be able to
apply the Kunen argument. \cref{KunenPrp} tells us that \(j[\lambda]\notin M\)
if \(\lambda\) is a strong limit; we would like to see that in fact
\(j[\lambda]\) never belongs to \(M\). This follows from the critical sequence
analysis of \(\lambda\):

\begin{defn}\label{CriticalSequence}\index{Critical sequence}
	Suppose \(N\) and \(P\) are transitive models of ZFC with the same ordinals
	and \(j : N\to P\) is a nontrivial elementary embedding. The {\it critical
	sequence} of \(j\) is the sequence \(\langle \kappa_n : n < \omega\rangle\)
	defined by recursion: \(\kappa_0 = \textsc{crt}(j)\) and for all \(n <
	\omega\), \(\kappa_{n+1} = j(\kappa_n)\). 
\end{defn}

In the context of \cref{CriticalSequence}, let \(\lambda = \sup_{n < \omega}
\kappa_n\). Clearly \(\lambda\) is the least ordinal such that
\(j[\lambda]\subseteq \lambda\). If \(\text{cf}^M(\lambda) = \omega\), \(j\) is
continuous at \(\lambda\), so \(j(\lambda) = \lambda\). In particular, if \(N =
V\), which is the case of interest in this section, then \(\lambda\) is the
first fixed point of \(j\) above \(\textsc{crt}(j)\).

In the case \(n > 1\), the conclusion of the following lemma is a considerable
understatement:
\begin{lma}\label{CrtSeqMeas}
	Suppose \(j : V\to M\) is a nontrivial elementary embedding and \(\langle
	\kappa_n :  n < \omega\rangle\) is its critical sequence. For any \(n <
	\omega\), if \(j\) is \(\kappa_n\)-strong then \(\kappa_n\) is measurable.
	\begin{proof}
		The proof is by induction on \(n\). Certainly \(\kappa_0 =
		\textsc{crt}(j)\) is measurable. Assume the lemma is true for \(n = m\),
		and we will show it is true for \(n = m+1\). Therefore assume \(j\) is
		\(\kappa_{m+1}\)-strong. In particular, \(j\) is \(\kappa_m\)-strong, so
		by our induction hypothesis, \(\kappa_m\) is measurable. By
		elementarity, \(\kappa_{m+1} = j(\kappa_m)\) is measurable in \(M\).
		Since \(j\) is \(\kappa_{m+1}\)-strong, \(P(\kappa_{m+1})\subseteq M\).
		Thus the measurability of \(\kappa_{m+1}\) in \(M\) is upwards absolute
		to \(V\), so \(\kappa_{m+1}\) is measurable.
	\end{proof}
\end{lma}

\begin{thm}[Kunen]\label{KunenInconsistency}
	Suppose \(j : V\to M\) is a nontrivial elementary embedding and \(\lambda\)
	is the least ordinal above \(\textsc{crt}(j)\) with \(j[\lambda]\subseteq
	\lambda\). Then \(j[\lambda]\notin M\).
	\begin{proof}
		Assume towards a contradiction that \(j[\lambda]\in M\). By
		\cref{SupercompactClosure}, \(j\) is \(\lambda\)-strong. Therefore
		\(\kappa_n\) is measurable for all \(n < \omega\) by \cref{CrtSeqMeas},
		and so \(\lambda\) is a strong limit cardinal. Since \(j[\lambda]
		\subseteq \lambda\) and \(P(\lambda)\subseteq M\), \(j[\lambda]\in M\).
		Thus \(\lambda\) is a strong limit cardinal, \(j[\lambda]\subseteq
		\lambda\), and \(j\) is \(\lambda\)-supercompact. From \cref{KunenPrp}
		we can therefore conclude that \(\textsc{crt}(j) \geq \lambda\),
		contradicting that \(\textsc{crt}(j) < \lambda\).
	\end{proof}
\end{thm}

A useful structural consequence of Kunen's Inconsistency Theorem is the
following lemma:

\begin{lma}\label{Kunen}
	Suppose \(\gamma\) is a cardinal, \(j : V\to M\) is a nontrivial elementary
	embedding, \( \textsc{crt}(j)\leq \gamma\), and \(P(\gamma)\subseteq M\).
	Then there is a measurable cardinal \(\kappa \leq\gamma\) such that
	\(j(\kappa) > \gamma\).
	\begin{proof}
		Let \(\langle \kappa_n : n < \omega\rangle\) be the critical sequence of
		\(j\) and  \(\lambda = \sup_{n <\omega} \kappa_n\). Thus \(\lambda\) is
		the least ordinal with \(j[\lambda]\subseteq \lambda\). By
		\cref{KunenInconsistency}, \(P(\lambda)\not\subseteq M\), so since
		\(P(\gamma)\subseteq M\), we have \(\gamma < \lambda\). Let \(n <
		\omega\) be least such that \(\kappa_n \leq \gamma < \kappa_{n+1}\).
		\cref{CrtSeqMeas} implies \(\kappa_n\) is measurable, and \(j(\kappa_n)
		= \kappa_{n+1} > \gamma\). Thus taking \(\kappa= \kappa_n\) proves the
		lemma.
	\end{proof}
\end{lma}

In one instance (\cref{GeneralizedSolovay}), we will need a strengthening of
\cref{Kunen} which has essentially the same proof:
\begin{lma}\label{Kunen2}
	Suppose \(\gamma\leq \lambda\) are cardinals and \(j : V\to M\) is a
	nontrivial \(\lambda\)-supercompact elementary embedding with
	\(\textsc{crt}(j)\leq \gamma\). Then there is a \(\lambda\)-supercompact
	cardinal \(\kappa \leq\gamma\) such that \(j(\kappa) > \gamma\).\qed
\end{lma}

\subsection{The wellfoundedness of the generalized Mitchell order}\label{MOWFSection}
The main theorem of this subsection states that the generalized Mitchell order
is a wellfounded partial order when restricted to a reasonable class of
countably complete ultrafilters. In fact, the wellfoundedness of the generalized
Mitchell order on countably complete ultrafilters is a special case of Steel's
wellfoundedness theorem for the Mitchell order on extenders
\cite{SteelWellfounded}, since countably complete ultrafilters are amenable
extenders in the sense of \cite{SteelWellfounded}, but we will give a much
simpler proof here.

We start with the fundamental fact that the Mitchell order is irreflexive:
\begin{lma}\label{MOStrict}
Suppose \(U\) is a countably complete nonprincipal ultrafilter. Then \(U\not\mo
U\).
\begin{proof}
Suppose towards a contradiction that \(U\mo U\). By \cref{HeredLemma}, if \(U'
\cong U\) is a uniform ultrafilter on a cardinal (as given by \cref{UniformIso})
then  \(U'\mo U'\) as well. We can therefore assume without loss of generality
that \(U\) is a uniform ultrafilter on a cardinal \(\lambda\).  By
\cref{MOSuper}, \(j_U\restriction P(\lambda)\in M_U\). In particular,
\(j_U\restriction \lambda\in M_U\), so \(M_U^\lambda\subseteq M_U\) by
\cref{SupercompactClosure}. Therefore \(j^{M_U}_U = j_U\restriction M_U\), for
example as a consequence of \cref{USupercompact}. Thus \(j_U\) is
\(\delta\)-supercompact for all cardinals \(\delta\). This contradicts
\cref{UFSuperBound}.
\end{proof}
\end{lma}

We now turn to the transitivity and wellfoundedness of the generalized Mitchell
order. The following lemma (which in the language of \cite{SteelWellfounded}
states that countably complete ultrafilters are {\it amenable}), is the key to
the proof. 
\begin{lma}\label{AmenabilityLemma}
Suppose \(U\) is a nonprincipal countably complete ultrafilter on a set \(X\).
Suppose \(\lambda\) is a cardinal such that \(P(\lambda)\subseteq M_U\). Then
\(M_U\vDash 2^{\lambda} < j_U(|X|)\).
\begin{proof}
	The proof proceeds by finding a measurable cardinal \(\kappa\leq |X|\) such
	that \(2^\lambda < j_U(\kappa)\).
	
	If \(\lambda < \textsc{crt}(j_U)\), then \(\kappa = \textsc{crt}(j_U)\)
	works. Therefore assume \(\textsc{crt}(j_U)\leq\lambda\). By \cref{Kunen},
	there is an measurable cardinal \(\kappa \leq \lambda\) such that
	\(j_U(\kappa) > \lambda\). We claim that \(\kappa \leq |X|\), which
	completes the proof. Assume not. Then \(|X|\) is smaller than the
	inaccessible cardinal \(\kappa\), and hence \(j_U(\kappa) = \kappa \leq
	\lambda\) (\cref{UFCardinality}), a contradiction.
\end{proof}
\end{lma}

We really only use the following consequence of \cref{AmenabilityLemma}:

\begin{cor}\label{Hered0}
Suppose \(U_0\mo U_1\) are countably complete nonprincipal hereditarily uniform
ultrafilters. Then \(M_{U_1}\vDash 2^{\lambda_{U_0}} < j_{U_1}(\lambda_{U_1})\).
\begin{proof}
This is immediate from \cref{AmenabilityLemma}, using the fact
(\cref{MOStrength}) that if \(U_0\mo U_1\) then \(P(\lambda_{U_0})\subseteq
M_{U_1}\).
\end{proof}
\end{cor}

\begin{cor}
Suppose \(U_0\mo U_1\) are countably complete nonprincipal hereditarily uniform
ultrafilters. Let \(\lambda = \lambda_{U_1}\). Then \(U_0\in
j_{U_1}(H(\lambda))\).
\begin{proof}
Since \(U_0\) is hereditarily uniform, \(M_{U_1}\vDash |\text{tc}(U_0)| =
2^{\lambda_{U_0}}\) By \cref{Hered0}, \(M_{U_1}\vDash 2^{\lambda_{U_0}} <
j_{U_1}(\lambda)\). Therefore \(U_0\in H^{M_{U_1}}(j_{U_1}(\lambda)) =
j_{U_1}(H(\lambda))\).
\end{proof}
\end{cor}

\begin{prp}\label{StrongTransitivity}
Suppose \(U_0\mo U_1\mo U_2\) are countably complete nonprincipal hereditarily
uniform  ultrafilters. Then \(U_0\mo U_2\) and \(M_{U_2}\vDash U_0\mo U_1\).
\begin{proof}
Let \(\lambda =\lambda_{U_1}\). Then \(U_0\in j_{U_1}(H(\lambda))\). By
\cref{MOSuper2}, \(M_{U_2}\) contains \(j_{U_1}(H(\lambda))\), so \(U_0\in
M_{U_2}\), which yields \(U_0\mo U_2\). In fact, by \cref{MOSuper2},
\(j_{U_1}(H(\lambda))= j_{U_1}^{M_{U_2}}(H(\lambda))\), and so \(U_0\in
j_{U_1}^{M_{U_2}}(H(\lambda))\subseteq M^{M_{U_2}}_{U_1}\). Thus \(U_0\in
M^{M_{U_2}}_{U_1}\), or other words, \(M_{U_2}\vDash U_0\mo U_1\).
\end{proof}
\end{prp}

\begin{cor}
The generalized Mitchell order is transitive on countably complete nonprincipal
hereditarily uniform ultrafilters.\qed
\end{cor}

It is worth pointing out that the generalized Mitchell order on extenders is
{\it not} transitive if there is a cardinal that is \(\kappa\)-strong where
\(\kappa\) is a measurable cardinal. The counterexample is described in
\cite{SteelWellfounded}. (The generalized Mitchell order is not  transitive on
arbitrary countably complete ultrafilters either as a consequence of the silly
counterexamples in \cref{StrengthSection}.) 

\begin{prp}\label{MOUniformWF}
The generalized Mitchell order is wellfounded on countably complete nonprincipal
hereditarily uniform ultrafilters.
\begin{proof}
Suppose not, and let \(\lambda\) be the least cardinal such that there is a
descending sequence \[U_0\gmo U_1\gmo U_2\gmo \cdots\] of  countably complete
hereditarily uniform ultrafilters with \(\lambda_{U_0} = \lambda\).

 By \cref{StrongTransitivity} and the closure of \(M_{U_0}\) under countable
sequences, the sequence \(\langle U_n : {1\leq n < \omega}\rangle\) belongs to
\(M_{U_0}\) and \[M_{U_0}\vDash U_1\gmo U_2\gmo \cdots\] Note that in
\(M_{U_0}\), \(U_1\) is a countably complete nonprincipal hereditarily uniform
ultrafilter, and by \cref{Hered0}, \(\lambda_{U_1} < j_{U_0}(\lambda)\). 

On the other hand, by the elementarity of \(j_{U_0}\), from the perspective of
\(M_{U_0}\), \(j_{U_0}(\lambda)\) is the least cardinal \(\lambda'\) such that
there is a descending sequence \[W_0\gmo W_1\gmo W_2\gmo \cdots\] of countably
complete hereditarily uniform ultrafilters such that \(\lambda_{W_0} =
\lambda'\). This is a contradiction.
\end{proof}
\end{prp}
One can prove a slightly more general result than \cref{MOUniformWF} although
this generality is never useful.
\begin{thm}\label{MOWF}
The generalized Mitchell order is wellfounded on nonprincipal countably complete
ultrafilters.\index{Generalized Mitchell order!wellfoundedness}
\begin{proof}
Assume towards a contradiction that \(U_0\gmo U_1\gmo \cdots\) are nonprincipal
countably complete ultrafilters. For each \(n < \omega\), let \(U'_n\) be a
hereditarily uniform ultrafilter such that \(U_n'\cong U_n\). Then by
\cref{HeredLemma}, \(U_0'\gmo U_1'\gmo\cdots\). This contradicts
\cref{MOUniformWF}.
\end{proof}
\end{thm}

\subsection{The nonlinearity of the generalized Mitchell order}\label{NonlinearitySection}
Before we discuss the extent to which the generalized Mitchell order is linear
under UA, it is worth pointing out:
\begin{itemize}
\item the obvious counterexamples to linearity
\item the maximal amount of linearity one can reasonably hope for.
\end{itemize}

The fact is that if there is a measurable cardinal, then the generalized
Mitchell order is not linear, even restricting to uniform countably complete
ultrafilters on cardinals. The known counterexamples to the linearity of the
generalized Mitchell order are closely related to the {\it Rudin-Frol\'ik order}
(the subject of \cref{RFChapter}):
\begin{defn}\index{Rudin-Frol\'ik order}
The {\it Rudin-Frol\'ik order} is defined on countably complete ultrafilters
\(U\) and \(W\) by setting \(U\D W\) if there is an internal ultrapower
embedding \(i : M_D\to M_W\) such that \(i\circ j_D = j_W\).
\end{defn}
By \cref{RKChar}, the Rudin-Keisler order can be defined in exactly the same way
except omitting the requirement that \(i\) be internal.
\begin{prp}
If \(U\D W\) are nonprincipal countably complete ultrafilters, then \(U\) and
\(W\) are incomparable in the generalized Mitchell order.
\begin{proof}
We first show \(U\not \mo W\). Since \(U\D W\), \(M_W\subseteq M_U\). Therefore
the fact that \(U\notin M_U\) implies that \(U\notin M_W\), and hence \(U\not
\mo W\).

We now show \(W\not \mo U\). Assume towards a contradiction that \(W\mo U\).
Assume without loss of generality that \(U\) is a uniform ultrafilter on a
cardinal \(\lambda\). (Since the Mitchell order is isomorphism invariant in its
second argument, this does not change our situation.) Since \(U\D W\), we have
\(U\RK W\) by \cref{RKChar}. Since \(U\) is hereditarily uniform and \(U\RK W
\mo U\), our lemma on the invariance of the Mitchell order (\cref{HeredLemma})
yields that \(U\mo U\). This contradicts \cref{MOStrict}.
\end{proof}
\end{prp}
A similar argument shows the following:
\begin{prp}
Suppose \(U\) and \(W\) are countably complete ultrafilters and there is a
nonprincipal \(D\D U,W\). Then \(U\) and \(W\) are incomparable in the
generalized Mitchell order.\qed
\end{prp}
Even this does not exhaust the known counterexamples to the linearity of the
generalized Mitchell order:

\begin{prp}
Suppose \(U_0\mo U_1\mo U_2\). Suppose \(U_0,U_2\D W\). Then \(U_1\) and \(W\)
are incomparable in the Mitchell order.\qed
\end{prp}
We omit the proof. The hypotheses of the proposition are satisfied if
\(U_0,U_1,U_2\) are normal ultrafilters on measurable cardinals
\(\kappa_0<\kappa_1<\kappa_2\) respectively and \(W = U_0\times U_2\).

All known examples of nonlinearity in the generalized Mitchell order are
accompanied by nontrivial relations in the Rudin-Frol\'ik order. A driving
question in this work is whether assuming UA, these are the only
counterexamples.
\begin{defn}\index{Irreducible ultrafilter}
A nonprincipal countably complete ultrafilter \(W\) is {\it irreducible} if for
all \(U\D W\), either \(U\) is principal or \(U\) is isomorphic to \(W\).
\end{defn}
The Irreducible Ultrafilter Hypothesis (IUH) essentially states that the sort of
counterexamples to the linearity of the Mitchell order that we have described
are the only ones.
\begin{iuh}
Suppose \(U\) and \(W\) are hereditarily uniform irreducible ultrafilters.
Either \(U\cong W\), \(U\mo W\), or \(W\mo U\). \index{Irreducible Ultrafilter
Hypothesis}
\end{iuh}
We can now make precise the question of the extent of the linearity of the
Mitchell order under UA:
\begin{qst}
Does UA imply IUH?
\end{qst}
With this in mind, let us turn to the positive results on linearity.
\section{Dodd soundness}\label{DoddSection}
\subsection{Introduction}
{\it Dodd soundness} is a fine-structural generalization of supercompactness,
introduced by Steel \cite{Schimmerling} in the context of inner model theory as
a strengthening of the initial segment condition. The following remarkable
theorem is due to Schlutzenberg \cite{Schlutzenberg}:
\begin{thm}[Schlutzenberg]\label{Schlutz}
Suppose \(L[\mathbb E]\) is an iterable Mitchell-Steel model and \(U\) is a
countably complete ultrafilter of \(L[\mathbb E]\). Then the following are
equivalent:\index{Schlutzenberg}
\begin{enumerate}[(1)]
\item \(U\) is irreducible.
\item \(U\) is isomorphic to a Dodd sound ultrafilter.
\item \(U\) is isomorphic to an extender on the sequence \(\mathbb E\).\qed
\end{enumerate}
\end{thm}
Since the total extenders on \(\mathbb E\) are linearly ordered by the Mitchell
order, this has the following consequence:
\begin{thm}[Schlutzenberg] Suppose \(L[\mathbb E]\) is an iterable
Mitchell-Steel model. Then \(L[\mathbb E]\) satisfies the Irreducible
Ultrafilter Hypothesis.\qed\index{Irreducible Ultrafilter Hypothesis}
\end{thm}
It is open whether this theorem can be extended to the Woodin models at the
finite levels of supercompactness. The main result of this section
(\cref{DoddMO}) states that UA alone suffices to prove the linearity of the
generalized Mitchell order on Dodd sound ultrafilters.

\subsection{Dodd sound embeddings, extenders, and ultrafilters}
In this subsection, we present a definition of Dodd soundness  due to the author
that is much simpler than the one given in \cite{Schimmerling} and
\cite{Schlutzenberg}, and that is easier to use in certain contexts. (The other
definition is also useful.) We then show that the two definitions are
equivalent.

\begin{defn}\label{EmbeddingSolid}\index{Soundness of an embedding}
Suppose \(M\) is a transitive class, \(j : V\to M\) is an elementary embedding,
and \(\alpha\) is an ordinal. Let \(\delta\) be the least ordinal such that
\(j(\delta)\geq \alpha\). Then \[j^\alpha: P(\delta)\to M\] is the function
defined by \(j^\alpha(X) = j(X)\cap \alpha\). The embedding \(j\) is said to be
{\it \(\alpha\)-sound} if \(j^\alpha\) belongs to \(M\).
\end{defn}
Recall that the bounded powerset of an ordinal \(\delta\) is defined by
\(P_\text{bd}(\delta) = \bigcup_{\xi < \delta} P(\xi)\). In the context of
\cref{EmbeddingSolid}, if \(\alpha = \sup j[\delta]\) it have been natural to
define \(j^\alpha = j\restriction P_{\text{bd}}(\delta)\). With this alternate
definition, \(j^\alpha\in M\) is an a priori weaker requirement. The next lemma
shows that this does not actually make a difference:

\begin{lma}\label{SpaceIrrel}
Suppose \(M\) is a transitive class, \(j : V\to M\) is an elementary embedding,
and \(\delta\) is an ordinal. Let \(\delta_* = \sup j[\delta]\). Then the
following are equivalent:
\begin{enumerate}[(1)]
\item \(j\) is \(\delta_*\)-sound.
\item \(j[P_\textnormal{bd}(\delta)]\in M\) or equivalently \(j\) is
\(2^{<\delta}\)-supercompact.
\item \(j\restriction P_\textnormal{bd}(\delta)\in M\).
\item \(j\restriction P^M_\textnormal{bd}(\delta)\in M\).
\end{enumerate}
\begin{proof}
{\it (1) implies (2):} Trivial. (The equivalence of
\(j[P_\textnormal{bd}(\delta)]\in M\) with \(2^{<\delta}\)-supercompactness is
immediate from \cref{XSC}.)

{\it (2) implies (3):} \(j\restriction P_\textnormal{bd}(\delta)\) is the
inverse of the transitive collapse of \(j[P_\textnormal{bd}(\delta)]\).

%{\it (2) implies (3):} Note first that since supercompactness implies strength
%(\cref{SupercompactClosure}), \(P(\delta)\subseteq M\). Note that \(j(X) \cap
%\delta_* = \bigcap_{\alpha < \delta} j(X\cap \alpha) = \bigcap_{\alpha <
%\delta} j(X\cap \alpha)\).

{\it (3) implies (4):} Trivial.

{\it (4) implies (1):} Assume \(j\restriction P^M_\textnormal{bd}(\delta)\in
M\). Since \(\delta\subseteq P^M_\textnormal{bd}(\delta)\), \[j\restriction
\delta = (j\restriction P^M_\textnormal{bd}(\delta))\restriction \delta\in M\]
Therefore \(j\) is \(\delta\)-supercompact. Since supercompactness implies
strength (\cref{SupercompactClosure}), \(P(\delta)\subseteq M\). In particular
\(j\restriction P^M_\textnormal{bd}(\delta) = j\restriction
P_\text{bd}(\delta)\). Finally for \(X\subseteq \delta\), \(j^{\delta_*}(X) =
\bigcup_{\xi < \delta} j(X\cap \xi)\), so \(j^{\delta_*}\) is definable from
\(j\restriction   P_\text{bd}(\delta)\) and hence \(j^{\delta_*}\in M\), which
shows (1).
\end{proof}
\end{lma}

\begin{lma}\label{SetDodd}
Suppose \(M\) is a transitive class, \(j : V\to M\) is an elementary embedding,
and \(\alpha\) is an ordinal. Then \(j\) is \(\alpha\)-sound if and only if
\(\{j(X)\cap \alpha : X\in V\}\in M\).
\begin{proof}
The forward direction is immediate since \(\{j(X)\cap \alpha : X\in V\} =
\text{ran}(j^\alpha)\). The reverse direction follows from the fact that
\(j^\alpha\) is the inverse of the transitive collapse of \(\{j(X)\cap \alpha :
X\in V\}\).
\end{proof}
\end{lma} 

Our next lemma states that the fragments \(j^\alpha\) ``pull back" under
elementary embeddings.

\begin{lma}\label{DoddPullback}
Suppose \(V\stackrel{i}{\longrightarrow} N\stackrel{k}{\longrightarrow} M\) are
elementary embeddings and \(j = k\circ i\). Suppose \(j^\alpha\in
\textnormal{ran}(k)\). Then \(k^{-1}(j^\alpha) = i^{k^{-1}(\alpha)}\).
\begin{proof}
Let \(\delta\) be the least ordinal such that \(j(\delta) \geq \alpha\). Note
that \(j^\alpha[\text{Ord}] = j[\delta]\in \text{ran}(k)\), so by our analysis
of derived embeddings (\cref{DerivedSC}), \(k\restriction \delta + 1\) is the
identity and \(i\) is \(\delta\)-supercompact. In particular,
\(P_\text{bd}(\delta)\subseteq M\) and \(k(P_\text{bd}(\delta)) =
P_\text{bd}(\delta)\).

Let \(h = k^{-1}(j^\alpha)\). Then \(\text{dom}(h) = k^{-1}(P_\text{bd}(\delta))
= P_\text{bd}(\delta)\). Thus for \(X\in \text{dom}(h)\), \(k(X) = X\), and
hence \[k(h(X)) = k(h)(k(X)) = k(h)(X) = j^\alpha(X) = j(X)\cap \alpha =
k(i(X))\cap \alpha\] By the elementarity of \(k\), this implies that \(h(X) =
i(X)\cap k^{-1}(\alpha)\), or in other words \(k^{-1}(j^\alpha) = h =
i^{k^{-1}(\alpha)}\), as desired.
\end{proof}
\end{lma}

We now turn to Dodd soundness.

\begin{defn}
If \(j : V\to M\) is an extender embedding, the {\it Dodd length of \(j\)},
denoted \(\alpha(j)\), is the least ordinal \(\alpha\) such that every element
of \(M\) is of the form \(j(f)(\xi)\) for some \(\xi < \alpha\).\index{Dodd
length}
\end{defn}

On first glance, one might believe that the Dodd length of an elementary
embedding \(j\) is the same as its {\it natural length}, denoted \(\nu(j)\), the
least \(\nu\) such that \(M = H^M(j[V]\cup \nu)\). In fact, equality may fail:
the issue is that \(\nu(j)\) is the least ordinal such that every element of
\(M\) is of the form \(j(f)(p)\) for a finite set \(p\subseteq \nu\), whereas in
the definition of \(\alpha(j)\), one must write every element of \(M\) in the
form \(j(f)(\xi)\) where \(\xi\) is not a finite set but a single ordinal below
\(\nu\).

Our main focus, of course, is on ultrafilters, and in this case the Dodd length
has an obvious characterization:\footnote{This gives us a counterexample to the
equality of Dodd length and natural length. Suppose \(U\) is a normal
ultrafilter on \(\kappa\). Let \(W = U^2\). Then \(\nu(j_W) = j_U(\kappa) + 1\)
but \(\alpha(j_W) = j_U(\kappa) + \kappa + 1\).}

\begin{lma}\label{UltraDodd} If \(j: V\to M\) is an ultrapower embedding, then \(\alpha(j) = \xi + 1\) where \(\xi\) is the least ordinal such that \(M = H^{M}(j[V]\cup \{\xi\})\). Therefore \(U\) is incompressible if and only if \(U\) is tail uniform and \(\alpha(j_U) = \id_U + 1\).\qed\end{lma}

Our next lemma establishes a limit on the solidity of an extender embedding. (It
is equivalent to the statement that no extender belongs to its own ultrapower.)

\begin{lma}\label{SoundnessLimit}
Suppose \(j : V\to M\) is an extender embedding and \(\alpha = \alpha(j)\). Then
\(j\) is not \(\alpha\)-sound.
\begin{proof}
Let us first show that if \(U\) is a countably complete tail uniform ultrafilter
on an ordinal \(\delta\), then \(j_U\) is not \(\id_U + 1\)-sound. Note that \[U
= \{A\subseteq \delta : \id_U\in j(A)\} = \{A\subseteq \delta : \id_U\in
j^{\id_U + 1}_U(A)\}\] so since \(U\notin M_U\), \(j^{\id_U + 1}_U\notin M\).
Thus \(j_U\) is not \(\id_U + 1\)-sound, as claimed.
	
We now handle the case where \(j\) is an arbitrary extender embedding. By the
definition of Dodd length, there is some \(\xi < \alpha\) and some function
\(f\in V\) such that \(j^\alpha = j(f)(\xi)\). Let \(U\) be the tail uniform
ultrafilter derived from \(j\) using \(\xi\), and let \(k : M_U\to M\) be the
factor embedding. Then \(\xi \in \text{ran}(k)\) and so \(j^\alpha\in
\text{ran}(k)\). Applying our lemma on pullbacks of the fragments \(j^\alpha\)
(\cref{DoddPullback}), \(k^{-1}(j^\alpha) = j_U^{k^{-1}(\alpha)}\). Therefore
\(j_U\) is \(k^{-1}(\alpha)\)-sound. But note that \(\id_U = k^{-1}(\xi) <
k^{-1}(\alpha)\). Hence \(j_U\) is \(\id_U + 1\)-sound, and this contradicts the
first paragraph.
\end{proof}
\end{lma}

An embedding is {\it Dodd sound} if it is as sound as it can possibly be: 

\begin{defn}\label{EmbeddingSoundDef}
Suppose \(M\) is a transitive class and \(j : V\to M\) is an elementary
embedding. Then \(j\) is said to be {\it Dodd sound} if \(j\) is \(\beta\)-sound
for all \(\beta < \alpha(j)\).\index{Dodd soundness}
\end{defn}

We now prove the equivalence between the Dodd soundness of an extender \(E\) as
it is defined in \cite{Schimmerling} and the Dodd soundness of its associated
embedding \(j_E\) as it is defined in \cref{EmbeddingSolid}. 

\begin{defn}
\begin{itemize}\index{Parameter}
	\item A {\it parameter} is a finite set of ordinals. 
	\item The {\it parameter order}\index{Parameter!parameter order} is defined
	on parameters \(p\) and \(q\) by \[p < q\iff \max(p\mathrel{\triangle} q)\in
	q\]
	\item If \(p\) is a parameter, then \(\langle p_i : i < |p|\rangle\) denotes
	the {\it descending} enumeration of \(p\). 
	\item For any \(k \leq |p|\), \(p\restriction k\) denotes the parameter
	\(\{p_i : i < k\}\).
\end{itemize}
\end{defn}

The point of enumerating parameters in descending order is that the parameter
order is then transformed into the lexicographic order:

\begin{lma}
Suppose \(p\) and \(q\) are parameters of length \(n\) and \(m\) respectively.
Then \(p < q\) if and only if \(\langle p_0,\dots,p_{n - 1}\rangle
<_\textnormal{lex} \langle q_0,\dots,q_{m - 1}\rangle\).\qed
\end{lma}

\begin{figure}
	\center
	\includegraphics[scale=.6]{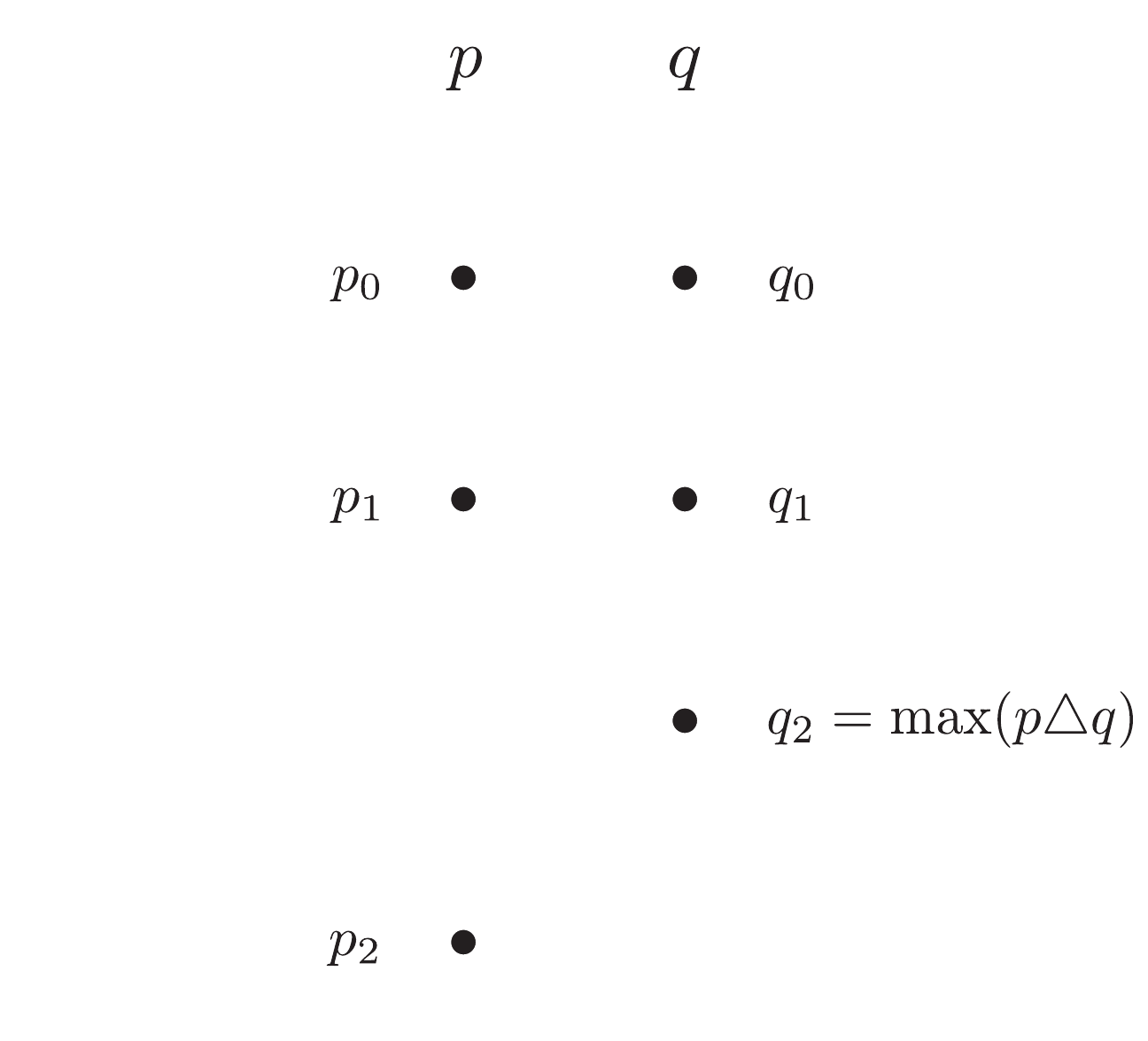}
	\caption{The parameter order}\label{ParamterFig}
\end{figure}

\begin{lma}
The parameter order is a set-like wellorder.\qed
\end{lma}

\begin{defn}
If \(j : V\to M\) is an elementary embedding and \(p\) is a parameter, then
\(\mu_j(p)\) is the least ordinal \(\mu\) such that \(p\subseteq j(\mu)\).
\end{defn}

\begin{defn}\index{Extender}
Suppose \(j : V\to M\) is an elementary embedding, \(p\) is a parameter, and
\(\nu < \min(p)\) is an ordinal. Let \(\delta = \mu_j(p)\). Then the {\it
extender of \(j\) below \((p,\nu)\)} is the set \[E^j\restriction p\cup \nu =
\{(q,A) : q\in [\nu]^{<\omega}\text{, }A\subseteq [\delta]^{<\omega}\text{, and
}p\cup q\in j(A)\}\]
\end{defn}

The restriction \(E^j\restriction p\cup \nu\) can be thought of as an extender
{\it relativized} to the parameter \(p\). It is possible to axiomatize
relativized extenders as directed systems of ultrafilters and associate to them
ultrapower embeddings, namely the direct limit of these systems. Instead we make
the following definition:

\begin{defn}\index{Extender!relativized}
A {\it relativized extender} is a set of the form \(E^j\restriction p\cup \nu\)
for some elementary embedding \(j\). The {\it extender embedding} associated to
a relativized extender \(E\), denoted \[j_E : V\to M_E\] is the unique \(j :
V\to M\) such that \(E = E^j\restriction p\cup \nu\) for some \(p,\nu\) and \(M
= H^M(j[V]\cup p\cup \nu)\).

If \(E\) is a relativized extender, \(\nu\) is an ordinal, and \(p\) is a
parameter, then \[E\restriction p\cup \nu = E^{j}\restriction p\cup \nu\] where
\(j = j_E\).
\end{defn}

The {\it Dodd parameter} of an extender is the key to the fine-structural proofs
of Dodd soundness, which are motivated by the fundamental solidity proofs from
fine structure theory. 

\begin{defn}\label{DoddParamDef1}
Suppose \(j : V\to M\) is an extender embedding. Then \(\eta(j)\) is the least
ordinal \(\eta\) such that for some parameter \(p\), \[M = H^M(j[V]\cup p\cup
\eta)\] The {\it Dodd parameter}\index{Dodd parameter}\index{\(p(j)\) (Dodd
parameter)} of \(j\), denoted \(p(j)\), is the least parameter \(p\) such that
\[M = H^M(j[V]\cup p\cup \eta(j))\]
\end{defn}
Thus if \(j\) is an ultrapower embedding, as it always will be in our
applications, then \(\eta = 0\). More generally, \(\eta\) is obviously always a
limit ordinal.

The Dodd parameter can also be defined recursively using the concept of an
\(x\)-generator of an elementary embedding:
\begin{defn}\index{Generator}
Suppose \(M\) and \(N\) are transitive models of ZFC, \(j : M\to N\) is an
elementary embedding, and \(x\in N\). Then an ordinal \(\xi\in N\) is an {\it
\(x\)-generator of \(j\)} if \(\xi\notin H^{N}(j[M]\cup \xi\cup \{x\})\).
\end{defn}

\begin{lma}
Suppose \(j : V\to M\) is an extender embedding. Let \(q\) be the
\(\subseteq\)-maximum parameter with the property that \(q_k\) is the largest
\(q\restriction k\)-generator of \(j\) for all \(k < |q|\). Then \(p(j) = q\)
and \(\eta(j)\) is the strict supremum of the \(q\)-generators of \(j\).
\begin{proof}
	 Let \(p = p(j)\), \(n = |p|\), and \(\eta = \eta(j)\). Fix \(k < n\). We
	 will show \(p_k\) is the largest \(p\restriction k\)-generator. 
	 
	 Since \(M = H^M(j[V]\cup p\cup \eta)\subseteq H^M(j[V]\cup p\restriction k
	 \cup (p_k+1))\), there are no \(p\restriction k\)-generators strictly above
	 \(p_k\). It therefore suffices to show that \(p_k\) is a \(p\restriction
	 k\)-generator. Assume not. Then \(p_k \in H^M(j[V]\cup p\restriction k \cup
	 p_k)\). Fix \(u\subseteq p_k\) such that \(p_k = j(f)(p\cup r)\) for some
	 function \(f\in V\). Let \(r = p\setminus \{p_k\} \cup u\). Then \(r < p\)
	 in the parameter order, but \(p\subseteq H^M(j[V]\cup r)\), and hence \(M =
	 H^M(j[V]\cup r \cup \eta)\), contrary to the minimality of the Dodd
	 parameter \(p\). 
	 
	 By the maximality of \(q\), this shows that \(p = q\restriction n\). We now
	 show that \(\eta\) is the strict supremum of the \(p\)-generators of \(j\).
	 Since \(M = H^M(j[V]\cup p\cup \eta)\), there are no \(p\)-generators
	 greater than or equal to \(\eta\). It therefore suffices to show that for
	 any \(\alpha < \eta\), there is a \(p\)-generator of \(j\) above
	 \(\alpha\). Suppose \(\alpha < \eta\). By the minimality of \(\eta\), \(M
	 \neq H^M(j[V]\cup p\cup \alpha)\), and so there is a \(p\)-generator of
	 \(j\) above \(\alpha\), as desired.
	 
	 Since \(\eta\) is a limit ordinal, there is no largest \(p\)-generator, and
	 hence \(p = q\).
\end{proof}
\end{lma}

\begin{cor}
Suppose \(j : V\to M\) is an extender embedding and \(p = p(j)\). Then for all
\(i < |p|\), \(p_i\) is a \(\{p_0,\dots,p_{i-1}\}\)-generator.\qed
\end{cor}

The following is Steel's definition of the Dodd soundness of an extender:

\begin{defn}\label{SteelSoundDef}\index{Dodd soundness!Dodd sound extender}\index{Dodd soundness!Dodd solid extender}
Suppose \(E\) is an extender, \(p = p(j_E)\), and \(\eta = \eta(j_E)\).
\begin{itemize}
\item  \(E\) is {\it Dodd solid} if \[E\restriction \{p_0,\dots,p_{i-1}\}\cup
p_i\in M_E\] for all \(i < |p|\).
\item \(E\) is {\it Dodd sound} if \(E\) is Dodd solid and \[E\restriction p\cup
\nu\in M_E\] for all \(\nu < \eta\).
\end{itemize}
\end{defn}

If \(E\) is an extender such that \(j_E\) is an ultrapower embedding, then \(E\)
is Dodd solid if and only if \(E\) is Dodd sound, simply because \(\eta(j_E) =
0\) (so the extra requirement for Dodd soundness holds vacuously). 

The following fact is essentially a matter of rearranging definitions:

\begin{thm}\label{DoddEquiv}
Suppose \(E\) is an extender. Then \(E\) is Dodd sound in the sense of
\cref{SteelSoundDef} if and only if \(j_E\) is Dodd sound in the sense of
\cref{EmbeddingSoundDef}.
\begin{proof}
Before we prove the equivalence, we prove three preliminary claims.

Let \(j = j_E\) and \(M = M_E\). Let \(\eta = \eta(j)\) and let \(p = p(j)\) be
the Dodd parameter of \(j\).

\begin{clm}\label{LeastParamClm}
\(p\cup \{\eta\}\) is the least parameter \(s\) such that every element of \(M\)
is of the form \(j(f)(q)\) for some \(q < s\).
\begin{proof}
Suppose not. Then fix \(s < p\cup \{\eta\}\) such that every element of \(M\) is
of the form \(j(f)(q)\) for some \(q < s\). Fix \(q < s\) such that \(p =
j(f)(q)\) for some \(f\). Then \(M = H^M(j[V]\cup q\cup \eta)\). Since \(p\) is
the least parameter with this property (by the definition of the Dodd
parameter), it follows that \(p\leq q\). In particular \(p < s\). Since \(p < s
< p\cup \{\eta\}\), \(s = p \cup r\) for some \(r \in [\eta]^{<\omega}\). Now
let \(\xi < \eta\) be a \(p\)-generator such that \(r \subseteq \xi\). Then
\(p\cup \{\xi\} = j(f)(u)\) for some \(u < s\). Since \(u\) generates \(p\), we
must have \(p\leq u\). Since \(p \leq u \leq p \cup r\), \(u = p \cup t\) for
some \(t < r\). In particular, since \(r\subseteq \xi\), \(t\subseteq \xi\). Now
\(\xi = j(f)(p\cup r)\) where \(r \in [\xi]^{<\omega}\), contradicting that
\(\xi\) is not a \(p\)-generator. 
\end{proof}
\end{clm}

Let \(\varphi\) be the function that sends a parameter to its rank in the
parameter order.
\begin{clm}\label{ConstClm} Suppose \(x\in M\) and \(q\) is a parameter. Then \(x = j(f)(q)\) for some function \(f\in V\) if and only if \(x = j(g)(\varphi(q))\) for some function \(g\in V\).\end{clm}
\begin{proof}
For the forwards direction, let \(g = f\circ \varphi^{-1}\), and for the reverse
direction, let \(f = g\circ \varphi\).
\end{proof}
From \cref{LeastParamClm} and \cref{ConstClm}, we obtain the following key
identity: 
\begin{equation}\varphi(p\cup \{\eta\}) = \alpha(j)\label{lengthparam}\end{equation}
(Recall that \(\alpha(j)\) denotes the Dodd length of \(j\), the least ordinal \(\alpha\) such that every element of \(M\) is of the form \(j(f)(\xi)\) for some \(\xi < \alpha\).)

\begin{clm}\label{unionclm}
Suppose \(q\) is a parameter and \(m = |q|\). For \(i < m\), let \[F_i = E\restriction \{q_0,\dots,q_{i-1}\}\cup q_i\] Then for any transitive model \(N\) of \textnormal{ZFC}, the following are equivalent:
\begin{enumerate}[(1)]
\item \(F_0,\dots,F_{m-1}\in N\).
\item \(j^{\varphi(q)}\in N\).
\end{enumerate}
\begin{proof}[Sketch]
{\it (1) implies (2):} Let \(\mu = \mu_j(q) = \mu_j(\{q_0\})\). If \(F_0,\dots,F_{m-1}\in N\), then so is the function \(e : P([\mu]^{<\omega})\to M\) defined by \(e(X) = \{r < q : r\in j(X)\}\). (\(e\) is the parameter version of \(j^q\).) This is because \(r\in e(X)\) if and only if \((r,X)\in F_i\) where \(i\) is such that \(\max(q\triangle r) = q_i\). 

Let \(\delta\) be least such that \(j(\delta) \geq \varphi(q)\). Then \(\varphi[\delta]\subseteq\mu\) and for \(A\subseteq\delta\), \(j(A)\cap \varphi(q) = \varphi^{-1}[e(\varphi[A])]\). This shows \(j^{\varphi(q)}\in N\).

{\it (2) implies (1):} Similar. 
\end{proof}
\end{clm}

Having proved the three claims, we finally turn to the equivalence of the two notions of Dodd soundness. (We will leave some of the parameter order theoretic details to the reader.)

Assume first that \(E\) is Dodd sound in the sense of \cref{SteelSoundDef}. Suppose \(\beta < \alpha(j)\), and we will show that \(j\) is \(\beta\)-sound. It suffices to show that \(j\) is \(\beta'\)-sound for some \(\beta' \geq \beta\), which allows us to increase \(\beta\) throughout the argument if necessary. By \cref{lengthparam}, by increasing \(\beta\), we may assume \(\varphi(p)\leq \beta\). Thus \(p\leq \varphi^{-1}(\beta) < \varphi^{-1}(\alpha(j)) = p\cup \{\eta\}\), as a consequence of \cref{lengthparam}. Let \(q = \varphi^{-1}(\beta)\). Then \(p\leq q  < p\cup \{\eta\}\), so \(q = p \cup r\) for some \(r \subseteq \eta\). Since \(\eta\) is a limit ordinal, by increasing \(\beta\) if necessary, we may assume \(|r| \leq 1\). By the Dodd soundness of \(E\), for all \(i < |q|\), \[E \restriction \{q_0,\dots,q_{i-1}\}\cup q_i\in M\]
This is because either \(\{q_0,\dots,q_{i-1}\}\cup q_i = \{p_0,\dots,p_{i-1}\}\cup p_i\) or \(\{q_0,\dots,q_{i-1}\}\cup q_i =p\cup \xi\) for some \(\xi < \eta\). Therefore by \cref{unionclm}, \(j^\beta\in M\) so \(j\) is \(\beta\)-sound.

Conversely, assume that \(j\) is Dodd sound as an elementary embedding. Let \(\beta = \varphi(p)\). Since \(p < p\cup \{\eta\}\), by \cref{lengthparam}, \(\beta < \alpha\). Therefore \(j^\beta\in M\) by the Dodd soundness of \(j\). By \cref{unionclm}, it follows that \(E\restriction \{p_0,\dots,p_{i-1}\}\cup p_i\) for all \(i < |p|\), so \(E\) is Dodd solid. If \(\eta = 0\), it follows that \(E\) is Dodd sound. Assume instead that \(\eta > 0\). Fix \(\xi < \eta\), and we will show \(E\restriction p \cup \xi\in M\). Let \(q = p \cup \{\xi\}\). Then \(q < p\cup \{\eta\}\), so \(\varphi(q) < \alpha\). Therefore by the Dodd soundness of \(j\), \(j^{\varphi(q)}\in M\). Applying \cref{unionclm}, it follows that \(E\restriction p \cup \xi\in M\).
\end{proof}
\end{thm}

It is worth remarking that the proof shows that an extender \(E\) is Dodd solid if and only if \(j_E\) is \(\beta\)-solid where \(\beta\) is the rank of \(p(j_E)\) in the canonical wellorder on parameters.

We now define Dodd sound ultrafilters. One could define an ultrafilter to be Dodd sound if its ultrapower embedding is Dodd sound, but then there would be many isomorphic Dodd sound ultrafilters all with the same associated embedding, which complicates the statements of our theorems and adds no real generality. Instead, we ensure that a Dodd sound ultrafilter is the canonical element of its isomorphism class:  

\begin{defn}\index{Dodd soundness!Dodd sound ultrafilter}
A countably complete ultrafilter is {\it Dodd sound} if it is incompressible and its ultrapower embedding is Dodd sound.
\end{defn}

The following alternate characterization of Dodd soundness for ultrafilters is immediate from \cref{UltraDodd} and \cref{SoundnessLimit}:

\begin{lma}\label{IdDodd}
	A tail uniform ultrafilter \(U\) on an ordinal \(\delta\) is Dodd sound if and only if \(j_U\) is \(\id_U\)-sound. That is, \(U\) is Dodd sound if and only if the function \(h : P(\delta)\to M_U\) defined by \(h(X) = j_U(X)\cap \id_U\) belongs to \(M_U\).\qed
\end{lma}

We finally provide a combinatorial characterization of Dodd soundness for ultrafilters:

\begin{defn}
Suppose \(U\) is an ultrafilter on an ordinal \(\delta\). 
\begin{itemize}
	\item A sequence \(\langle S_\alpha : \alpha < \delta\rangle\) is {\it \(U\)-threadable} if there is a set \(S\subseteq \delta\) such that \(S\cap \alpha = S_\alpha\) for \(U\)-almost all \(\alpha < \delta\).
	\item A {\it soundness sequence} for \(U\) is a sequence \(\langle \mathcal A_\alpha :\alpha < \delta\rangle\) such that for any sequence \(\langle S_\alpha : \alpha <\delta\rangle\), the following are equivalent:
\begin{enumerate}[(1)]
\item \(\langle S_\alpha :\alpha < \delta\rangle\) is \(U\)-threadable.
\item \(S_\alpha\in \mathcal A_\alpha\) for \(U\)-almost all \(\alpha\). 
\end{enumerate}
\end{itemize}
\end{defn}

\begin{thm}
A tail uniform ultrafilter \(U\) is {\it Dodd sound} if and only if it has a soundness sequence.
\end{thm}
\begin{proof}
Note that a sequence \(\langle S_\alpha : \alpha < \delta\rangle\) is \(U\)-threadable if and only if \[[\langle S_\alpha : \alpha < \delta\rangle] = j_U(S)\cap a_U\] some \(S\subseteq \delta\). Thus \(\langle \mathcal A_\alpha :\alpha < \delta\rangle\) is a soundness sequence for \(U\) if and only if \[[\langle \mathcal A_\alpha :\alpha < \delta\rangle]_U = \{j_U(S)\cap \id_U : S\subseteq \delta\}\]
By \cref{SetDodd}, it follows that \(U\) has a soundness sequence if and only if \(j_U\) is \(\id_U\)-sound, or in other words (applying \cref{IdDodd}) \(U\) is Dodd sound.
\end{proof}

\subsection{The generalized Mitchell order on Dodd sound ultrafilters}
In this short section, we prove the linearity of the Mitchell order on Dodd sound ultrafilters. We first prove a stronger statement that characterizes \(P(P(\lambda))\cap M_W\) when \(W\) is Dodd solid in terms of the Lipschitz order on subsets of \(P(\lambda)\).
\begin{prp}\label{LipschitzDodd}
Suppose \(W\) is a Dodd sound ultrafilter on a cardinal \(\lambda\). Then 
\[P(P(\lambda))\cap M_W = \{X\subseteq P(\lambda) : X \sLi W\}\]
\begin{proof}
Suppose \(X\subseteq P(\lambda)\).

Assume first that \(X\sLi W\). By our characterization of the Lipschitz order where the second argument is an ultrafilter (\cref{LipChar}), this means that there is a set \(Z\in M_W\) such that for all \(A\subseteq \delta\), \(A\in X\) if and only if \(j_W(A)\cap \id_W \in Z\). But then \(X = (j^{\id_W})^{-1}[Z]\), so \(X\in M_W\).

Conversely, suppose \(X\in M_W\). Let \(Z = j^{\id_W}[X]\). Then \(Z\in M_W\) and for all \(A\subseteq \delta\), \(A\in X\) if and only if \(j_W(A)\cap \id_W= j^{\id_W}(A) \in Z\). It follows that \(X\sLi W\).
\end{proof}
\end{prp}

\begin{cor}\label{LipMO}
Suppose \(U\) and \(W\) are countably complete ultrafilters on \(\lambda\) and \(W\) is Dodd sound. Then \(U\sLi W\) if and only if \(U\mo W\). In particular, if \(U\E W\) then \(U\mo W\).\qed.
\end{cor}

In particular, the Lipschitz order is wellfounded on Dodd sound ultrafilters.

\begin{thm}[UA]\label{DoddMO}
The generalized Mitchell order is linear on Dodd sound ultrafilters.\index{Generalized Mitchell order!linearity!on Dodd sound ultrafilters}
\begin{proof}
Suppose \(U\) and \(W\) are Dodd sound ultrafilters. By the linearity of the Lipschitz order on \(\Un\), either \(U\sLi W\), \(U = W\), or \(U\sgLi W\). Therefore by \cref{LipschitzDodd}, either \(U\mo W\), \(U = W\), or \(U\gmo W\), as desired.
\end{proof}
\end{thm}
Notice that the linearity of the Mitchell order on Dodd sound ultrafilters actually follows from the linearity of the Lipschitz order, which perhaps is weaker than UA.

As a consequence of \cref{LipMO}, if \(W\) is Dodd sound and \(U\sE W\), then \(U\mo W\). We now prove a strong converse, which is closely related to \cref{USupercompact}:

\begin{prp}\label{SuperMOKet}\index{Seed order!vs. the generalized Mitchell order}\index{Ketonen order!vs. the generalized Mitchell order}
Suppose \(U\) is a countably complete ultrafilter on a cardinal \(\lambda\) and \(W\) is a nonprincipal uniform ultrafilter on a cardinal \(\delta\) such that \(j_W\) is \(\lambda\)-supercompact. If \(U \mo W\), then \(U\swo W\).
\begin{proof}
Note that \((j_U(j_W),j_U\restriction M_W)\) is a \(0\)-internal comparison of \((j_U,j_W)\) by the standard identity: \[j_U(j_W)\circ j_U = j_U\circ j_W\] Since \(j_W\) is \(\lambda\)-supercompact, \(j_U\restriction M_W = j_U^{M_W}\), which is definable over \(M_W\) since \(U\mo W\).

Since \(j_W\) is \(\lambda\)-supercompact, \(\lambda \leq \delta\) by \cref{UFSuperBound}. Therefore for all \(\alpha < \lambda\), \(j_W(\alpha) < \id_W\). Applying \L o\'s's Theorem,
\[j_U(j_W)(\id_U) = [j_W\restriction \lambda]_U < j_U(\id_W)\]
Thus \((j_U(j_W),j_U\restriction M_W)\) witnesses that \(U\swo W\).
\end{proof}
\end{prp}

This raises the question of whether the Ketonen order extends the generalized Mitchell order. One should restrict attention here to countably complete uniform ultrafilters on cardinals, or else there are silly counterexamples. If this were true, it would complete the picture in which the wellfoundedness of the Ketonen order explains that of all the other known wellfounded orders. It is consistently false, however (\cref{CummingsExample}):
\begin{prp}
Suppose \(\kappa\) is \(2^\kappa\)-supercompact and \(2^\kappa = 2^{\kappa^+}\). Then there are \(\kappa\)-complete uniform ultrafilters \(U\) and \(W\) on \(\kappa\) and \(\kappa^+\) respectively such that \(W\mo U\).\qed
\end{prp}

Thus \(W\mo U\) but \(U\sE W\) simply because \(\delta_U < \delta_W\). (This is a consequence of \cref{SpaceLemma}.) By \cref{SuperMOKet}, if \(U\) and \(W\) are uniform ultrafilters on the same cardinal \(\lambda\) and both \(j_U\) and \(j_W\) are \(\lambda\)-supercompact, then \(U\mo W\) implies \(U\sE W\).
 
\begin{lma}\label{MOKetExtend}
Suppose \(\lambda\) is a cardinal, \(W\) is a countably complete ultrafilter on \(\lambda\), and \(Z\) is a countably complete ultrafilter such that \(W\mo Z\). Assume that for all \(\alpha < \lambda\), \(\mathscr B(\lambda,\alpha)\subseteq M_Z\) and \(M_Z\vDash \mathscr B(\lambda,\alpha) \leq 2^\lambda\). Then for any \(U\sE W\), \(U\mo Z\).
\begin{proof}
Since \(W\mo Z\), \(P(\lambda)\subseteq M_Z\) and in fact \(P(\lambda)^\lambda\subseteq M_Z\). Moreover \[M_Z\vDash \left|\bigcup_{\alpha < \lambda} \mathscr B(\lambda,\alpha)\right| \leq 2^\lambda = |P(\lambda)|\] Hence \((\bigcup_{\alpha < \lambda} \mathscr B(\lambda,\alpha))^\lambda\subseteq M_Z\), so \(\prod_{\alpha\in I}\mathscr B(\lambda,\alpha)\in M_Z\) for any set \(I\subseteq \lambda\). 

Now suppose \(U\sE W\). Fix \(I\in W\) and \(\langle U_\alpha : \alpha\in I\rangle\in \prod_{\alpha\in I}\mathscr B(\lambda,\alpha)\) such that \(U = W\text{-}\lim_{\alpha \in I} U_\alpha\). Then the sequence \(\langle U_\alpha : \alpha\in I\rangle\in M_Z\), so \(U\in M_Z\), so \(U\mo Z\), as desired.
\end{proof}
\end{lma}

In fact, this lemma yields the somewhat stronger result that for any \(I\in W\) and sequence \(\langle U_\alpha : \alpha \in I\rangle\) of ultrafilters with \(\delta_{U_\alpha} < \lambda\), \(W\text{-}\lim_{\alpha\in I} U_\alpha \mo Z\).

\begin{cor}[UA]
Assume \(\lambda\) is a cardinal such that \(2^{<\lambda} = \lambda\). If \(W\) and \(Z\) are countably complete ultrafilters on \(\lambda\) such that \(W\mo Z\), then \(W\sE Z\).
\begin{proof}
Given the assumption that \(2^{<\lambda} = \lambda\) and the fact that \(P(\lambda)\subseteq M_Z\), it is not hard to show that \(\Un_\alpha\in M_Z\) and \(M_Z\vDash |\Un_\alpha| \leq 2^\lambda\) for all \(\alpha < \lambda\). Therefore we are in a position to apply \cref{MOKetExtend} to any ultrafilter \(U\sE W\). Assume towards a contradiction that \(W\not \sE Z\). By the linearity of the Ketonen order, \(Z\sE W\). Now \(Z\sE W\mo Z\), so by \cref{MOKetExtend}, \(Z\mo Z\). This contradicts the strictness of the Mitchell order (\cref{MOStrict}).
\end{proof}
\end{cor}

\begin{cor}[UA + GCH]\index{Ketonen order!vs. the generalized Mitchell order}The Ketonen order extends the generalized Mitchell order on countably complete uniform ultrafilters on infinite cardinals.\qed\end{cor}
\cref{MOExtensions} shows that the same conclusion can be deduced from UA alone. This will be achieved by proving from UA that if \(W\mo Z\), then \(Z\) is \(\lambda_W\)-supercompact. The result then follows from \cref{SuperMOKet}.
\section{Generalizations of normality}\label{GeneralizedNormalSection}
In this section, we develop the theory of normal fine ultrafilters, the natural combinatorial generalization of normal ultrafilters, and a central component of the classical theory of supercompact cardinals. The main result of the section (\cref{GCHLinear}) states roughly that UA + GCH implies that all these ultrafilters are linearly ordered by the Mitchell order.

\begin{defn}\index{\(\mathscr N_\lambda\) (normal ultrafilters on \(P_\text{bd}(\lambda)\))}
	For any infinite cardinal \(\lambda\), let \(\mathscr N_\lambda\) be set of normal fine ultrafilters on \(P_\textnormal{bd}(\lambda)\). Let \(\mathscr N = \bigcup_{\lambda }\mathscr N_\lambda\).
\end{defn}

We provide the definitions of normality and fineness in \cref{NFSection}.

\begin{thm}[UA]\label{GCHLinear}\index{Generalized Mitchell order!linearity!on normal fine ultrafilters}
	Suppose \(\lambda\) is a cardinal such that \(2^{<\lambda} = \lambda\). Then \(\mathscr N_\lambda\) is wellordered by the Mitchell order. Therefore assuming the Generalized Continuum Hypothesis, \(\mathscr N\) is linearly ordered by the Mitchell order.
\end{thm}

We asserted that UA + GCH would roughly imply that the Mitchell order is linear on the class of all normal fine ultrafilters, but our theorem only mentions the subclass \(\mathscr N\). In fact, the class of all normal fine ultrafilters is not literally linearly ordered by the Mitchell order for a number of reasons: one reason is that distinct normal fine ultrafilters can be isomorphic and hence Mitchell incomparable. \cref{NIso} below, however, shows that every normal fine ultrafilter is isomorphic to an element of \(\mathscr N\), so \cref{GCHLinear} essentially covers all the bases.

A key concept in the proof of \cref{GCHLinear}, introduced here for the first time, is that of an isonormal ultrafilter.
\begin{defn}\index{Isonormal ultrafilter}
Suppose \(\lambda\) is a cardinal. An ultrafilter \(U\) on \(\lambda\) is {\it isonormal} if \(U\) is weakly normal and \(j_U\) is a \(\lambda\)-supercompact embedding.
\end{defn}

We define weak normality in \cref{WeaklyNormalSection}. The concept dates back to Solovay and Ketonen \cite{Ketonen2}. The other main theorem of this section explains how isonormal ultrafilters get their name:
\begin{repthm}{IsoNormalThm}
Suppose \(U\) is a nonprincipal ultrafilter. Then \(U\) is isonormal if and only if \(U\) is the incompressible ultrafilter isomorphic to a normal fine ultrafilter. In particular, every normal fine ultrafilter is isomorphic to a unique isonormal ultrafilter.
\end{repthm}

The proof appears in \cref{SingularSolovaySection}. The forwards direction is quite easy, but the reverse implication requires quite a bit of work amounting to a generalization of the theorem of \cite{Solovay} known as {\it Solovay's Lemma} to singular cardinals. This generalization constitutes a fundamental and (apparently) new fact about supercompactness whose proof requires some basic notions from PCF theory. 

The investigation of isonormal ultrafilters is related back to the linearity of the Mitchell order by the following proposition:
\begin{repthm}{IsoSound}
Suppose \(2^{<\lambda} = \lambda\). Then every isonormal ultrafilter \(U\) on \(\lambda\) is Dodd sound.
\end{repthm}

This is basically just a matter of defining weakly normal ultrafilters on singular cardinals.  We actually prove our main theorem (\cref{GCHLinear}) right now. But we will need to assume \cref{IsoSound} and \cref{IsoNormalThm}. We also need a lemma that shows \(\mathscr N\) is well-behaved under the Mitchell order assuming GCH:
\begin{lma}\label{NHered}
If \(2^{<\lambda} = \lambda\), then any \(\mathcal U \in \mathscr N_\lambda\) is hereditarily uniform and satisfies \(\lambda_\mathcal U = \lambda\).
\begin{proof}
Since \(P_\text{bd}(\lambda)\) is transitive, \(|\text{tc}(P_\text{bd}(\lambda))| = |P_\text{bd}(\lambda)| = 2^{<\lambda} = \lambda\). On the other hand, since \(j_\mathcal U\) is \(\lambda\)-supercompact, \cref{UFSuperBound} implies \(\lambda_\mathcal U\geq \lambda\). Thus \(|\text{tc}(P_\text{bd}(\lambda))| = \lambda_\mathcal U\), so \(\mathcal U\) is hereditarily uniform.
\end{proof}
\end{lma}

We finally prove \cref{GCHLinear} assuming \cref{IsoSound} and \cref{IsoNormalThm}.

\begin{proof}[Proof of \cref{GCHLinear}]
Suppose \(\mathcal U\) and \(\mathcal W\) are elements of \(\mathscr N_\lambda\). We show that either \(\mathcal U\mo \mathcal W\), \(\mathcal U = \mathcal W\), or \(\mathcal U\gmo\mathcal W\). Applying \cref{IsoNormalThm}, let \(U\) be the isonormal ultrafilter isomorphic to \(\mathcal U\) and let \(W\) be the isonormal ultrafilter isomorphic to \(\mathcal W\). Note that \(U\) and \(W\) are uniform ultrafilters on the cardinal \(\lambda_\mathcal U = \lambda_\mathcal W = \lambda\) (\cref{NHered}). We have \(2^{<\lambda} = \lambda\) by assumption, so \cref{IsoSound} yields that \(U\) and \(W\) are Dodd sound. By the linearity of the Mitchell order on Dodd sound ultrafilters (\cref{DoddMO}), we are in one of the following cases:

\begin{case} \label{EqLin}\(U = W\).\end{case} 
\begin{proof}[Proof in \cref{EqLin}]
Since \(\mathcal U \cong U = W\cong \mathcal W\),  \cref{NormalIso} below implies \(\mathcal U=\mathcal W\).
\end{proof}
\begin{case} \label{moCase}\(U \mo W\).\end{case}
\begin{proof}[Proof in \cref{moCase}]
Since \(W\cong \mathcal W\), we have \(U\mo \mathcal W\). Since \(\mathcal U\) is hereditarily uniform (\cref{NHered}) and isomorphic to \(U\), the isomorphism invariance of the generalized Mitchell order on hereditarily uniform ultrafilters (\cref{HeredLemma}) implies \(\mathcal U\mo \mathcal W\). 
\end{proof}
\begin{case}\label{gmoCase}\(U \gmo W\).\end{case}
\begin{proof}[Proof in \cref{gmoCase}] 
Proceeding as in \cref{moCase}, we obtain \(\mathcal U\gmo \mathcal W\).
\end{proof}
This shows that either \(\mathcal U\mo \mathcal W\), \(\mathcal U = \mathcal W\), or \(\mathcal U\gmo\mathcal W\), as desired.

We finally sketch the proof that \(\mathscr N\) is linearly ordered by the Mitchell order assuming UA + GCH. It suffices to show the following: suppose \(\mathcal U\in \mathscr N_\gamma\), \(\mathcal W\in \mathscr N_\lambda\), and \(2^{<\lambda} = \lambda\). Then \(\mathcal U \mo \mathcal W\). Let \(U\) be the isonormal ultrafilter of \(\mathcal U\), so by the proof of \cref{NHered}, \(U\) is an ultrafilter on \(\gamma\). Since \(2^\gamma \leq 2^{<\lambda} = \lambda\), \(U\in H_{(2^\gamma)^+}\subseteq H_{\lambda^+}\subseteq M_{\mathcal W}\) Since \(P(P_{\text bd}(\gamma))\subseteq M_{\mathcal W}\), this easily implies that \(\mathcal U\mo \mathcal W\).
\end{proof}

\subsection{Normal fine ultrafilters}\label{NFSection}
In this section, we give the general definition of a normal fine ultrafilter, which is the natural combinatorial generalization of the notion of a normal ultrafilter on a cardinal. This begins with the generalized diagonal intersection operation:
\begin{defn}\index{Diagonal intersection}\index{\(\triangle_{\alpha < \delta}\) (diagonal intersection)}
Suppose \(X\) is a set and \(\langle A_x : x\in X\rangle\) is a sequence with \(A_x\subseteq P(X)\) for all \(x\in X\). The {\it diagonal intersection} of \(\langle A_x : x\in X\rangle\) is the set \[\triangle_{x\in X} A_x = \left\{\sigma\in P(X) : \sigma\in\textstyle \bigcap_{x\in \sigma} A_x\right\}\]
\end{defn}

\begin{defn}\index{Family over a set}
	If \(X\) is a set, a {\it family over \(X\)} is a family \(Y\) of subsets of \(X\) such that every element of \(X\) belongs to some element of \(Y\).
\end{defn}

Thus any set \(Y\) is a family on a unique set (namely \(X = \bigcup Y\)).

\begin{defn}\label{NormalFineDef}\index{Normal fine ultrafilter}\index{Fine ultrafilter}
Suppose \(Y\) is a family over \(X\). A filter \(\mathcal F\) on \(Y\) is:
\begin{itemize}
\item {\it fine} if for any \(x\in X\), \(\mathcal F\) concentrates on \(\{\sigma : x\in \sigma\}\).
\item {\it normal} if for any \(\{A_x : x\in X\}\subseteq \mathcal F\), \(\triangle_{x\in X} A_x\in \mathcal F\).
\end{itemize}
\end{defn}

\begin{rmk}
Let us make some remarks regarding this definition. 
\begin{enumerate}[(1)]
	\item It makes sense to discuss normal fine filters on \(Y\) without mention of \(X\), since \(X = \bigcup Y\) is determined from \(Y\).
	\item The structure of the underlying set \(Y\) is usually not that important since a normal fine ultrafilter \(\mathcal U\) on \(Y\) can always be lifted to a normal fine ultrafilter on \(P(X)\) where \(X = \bigcup_{\sigma\in Y} \sigma\). Therefore it is tempting to restrict consideration to normal fine ultrafilters on \(P(X)\) for some \(X\). It is often important for technical reasons, however, that the underlying set \(Y\) be small; usually we want \(|Y| = |\bigcup Y|\).
	\item The structure of the set \(X\) is also usually irrelevant, but sometimes it is useful that \(X\) be transitive or that \(X\) be a cardinal. Suppose \(X\) and \(X'\) are sets and \(f : X\to X'\) is a surjection. If \(Y\) is a family over \(X\), then \(Y' = \{f[\sigma] : \sigma \in Y\}\) is a family over \(X'\) and \(g(\sigma) = f[\sigma]\) defines a surjection from \(Y\) to \(Y'\). If \(\mathcal U\) is an ultrafilter on \(Y\), then \(g_*(\mathcal U)\) is an isomorphic ultrafilter on \(Y'\) and moreover \(\mathcal U'\) is normal (fine) if and only if \(\mathcal U\) is normal (fine). (This is the ultrafilter theoretic analog of \cref{XSC}.) 
	\item An ultrafilter on an ordinal is fine if and only if it is tail uniform. Thus a normal fine ultrafilter on \(\kappa\) is the same thing as a normal ultrafilter on \(\kappa\).
\end{enumerate}
\end{rmk}

The connection between normality and supercompactness is clear from the following lemma:
\begin{lma}\label{NormalFineChar}
Suppose \(Y\) is a family over \(X\) and \(\mathcal U\) is an ultrafilter on \(Y\).
\begin{enumerate}[(1)]
\item \(\mathcal U\) is fine if and only if \(j_\mathcal U[X]\subseteq \id_\mathcal U\).
\item \(\mathcal U\) is normal if and only if \(\id_\mathcal U\subseteq j_\mathcal U[X]\).\qed
\end{enumerate}
Thus \(\mathcal U\) is a normal fine ultrafilter on \(Y\) over \(X\) if and only if \(\id_\mathcal U = j_\mathcal U[X]\), or in other words, \(\id_\mathcal U\) witnesses that \(j_\mathcal U\) is \(X\)-supercompact.
\end{lma}
 \cref{NormalFineChar} yields the main source of normal fine ultrafilters.

\begin{lma}\label{DerivedNF}
Suppose \(j :V\to M\) is an \(X\)-supercompact elementary embedding and \(Y\subseteq P(X)\) is such that \(j[X]\in j(Y)\). 
\begin{itemize}
	\item \(Y\) is a family over \(X\).
	\item The ultrafilter \(\mathcal U\) on \(Y\) derived from \(j\) using \(j[X]\) is a normal fine ultrafilter on \(Y\). 
	\item Let \(k : M_\mathcal U\to M\) be the factor embedding. Then \(k(\alpha) = \alpha\) for all \(\alpha \leq |X|\).
\end{itemize}
\begin{proof}
Immediate from \cref{DerivedSC} and \cref{NormalFineChar}.
\end{proof}
\end{lma}

Another consequence of \cref{NormalFineChar} is the following fact, which does not seem to have a simple combinatorial proof:
\begin{lma}\label{NormalIso}
Suppose \(\mathcal U\) and \(\mathcal W\) are normal fine ultrafilters on \(Y\). If \(\mathcal U\cong \mathcal W\) then \(\mathcal U = \mathcal W\).
\begin{proof}
Let \(X = \bigcup Y\). Since \(\mathcal U\cong\mathcal W\), \(j_\mathcal U = j_\mathcal W\). By \cref{NormalFineChar}, \(\id_\mathcal U = j_\mathcal U[X] =  j_\mathcal W[X]= \id_\mathcal W\). Thus \(\mathcal U = \{A\subseteq Y : \id_\mathcal U \in j_\mathcal U(A)\} =  \{A\subseteq Y : \id_\mathcal W \in j_\mathcal W(A)\}  = \mathcal W\).
\end{proof}
\end{lma}

It also follows that any normal fine ultrafilter is countably complete. This is because the proof that an \(\omega\)-supercompact ultrapower embedding \(j :V\to M\) has the property that \(M^\omega\subseteq M\) does not really require that \(M\) is wellfounded. (The reader will lose nothing by simply appending countable completeness to the definition of normality, rather than proving it from the definition we have given.)

Recall the class \(\mathscr N\) defined in the previous section. We finish this section by proving that every normal fine ultrafilter is isomorphic to a unique element of \(\mathscr N\). 
\begin{prp}\label{NIso}
Any nonprincipal normal fine ultrafilter \(\mathcal D\) is isomorphic to a unique ultrafilter \(\mathcal U\in \mathscr N\).
\end{prp}
For this we will use a basic lemma about supercompactness:
\begin{lma}\label{SmallCfCompact2}
Suppose \(j : V\to M\) is \(\lambda\)-supercompact and \(\sup j[\lambda] = j(\lambda)\). Then \(j\) is \(\lambda^{\iota}\)-supercompact where \(\iota = \textnormal{cf}(\lambda)\). In particular, \(j\) is \(\lambda^+\)-supercompact.
\begin{proof}
Let \(\kappa = \textsc{crt}(j)\). \cref{SmallCfCompact} states that \(j\) is \(\lambda^{<\kappa}\)-supercompact. It suffices to show that \(\iota < \kappa\): then since \(j\) is \(\lambda^{<\kappa}\)-supercompact, \(j\) is \(\lambda^{\iota}\)-supercompact, and so since \(\lambda^\iota>\lambda\), \(j\) is \(\lambda^+\)-supercompact. 

We now show \(\iota < \kappa\). Since \(\sup j[\lambda] = j(\lambda)\) and \(j[\lambda]\in M\), \(\text{cf}^M(j(\lambda)) = \text{cf}(\lambda) = \iota\). On the other hand, by elementarity \(\text{cf}^M(j(\lambda)) = j(\text{cf}(\lambda)) = j(\iota)\). It follows that \(j(\iota) = \iota\). Since \(j\) is \(\iota\)-supercompact, the Kunen Inconsistency Theorem (\cref{KunenInconsistency}) implies \(\iota < \kappa\) where \(\kappa = \textsc{crt}(j)\). 
\end{proof}
\end{lma}
Actually, we always have \(\lambda^+ = \lambda^{<\kappa}\) in the context of \cref{SmallCfCompact2}, and this is how SCH above a supercompact is proved.

\begin{proof}[Proof of \cref{NIso}]
Obviously, any normal fine ultrafilter is isomorphic to a normal fine ultrafilter on \(P(\lambda)\) for some cardinal \(\lambda\). Therefore assume \(\mathcal D\) is a normal ultrafilter on \(P(\lambda)\), and we will show that \(\mathcal D\) is isomorphic to a normal fine ultrafilter on \(P_\text{bd}(\lambda')\) for some cardinal \(\lambda'\).

If \(\mathcal D\) concentrates on \(P_\text{bd}(\lambda)\), we are done, since \(\mathcal D\) is then isomorphic to \(\mathcal D\mid P_\text{bd}(\lambda)\). So assume \(\mathcal D\) does not concentrate on \(P_\text{bd}(\lambda)\). By \L o\'s's Theorem, \(\id_\mathcal D = j_\mathcal D[\lambda]\) is unbounded in \(j_\mathcal D(\lambda)\). In other words, \(j_\mathcal D\) is continuous at \(\lambda\). Therefore by \cref{SmallCfCompact2}, \(j_\mathcal D\) is \(\lambda^+\)-supercompact. Note that \(j_\mathcal D[\lambda^+]\) is not cofinal in \(j_\mathcal D(\lambda^+)\): otherwise \(j_\mathcal D(\lambda^+) =\text{cf}^{M_\mathcal D}(j_\mathcal D(\lambda^+) ) = \lambda^+\), so \(\textsc{crt}(j_\mathcal D) > \lambda^+\) by \cref{KunenInconsistency0}, which implies that \(\mathcal D\) is principal. Therefore let \(\mathcal U\) be the normal fine ultrafilter on \(P_\text{bd}(\lambda^+)\) derived from \(j_\mathcal D\) using \(j_\mathcal D[\lambda^+]\). Then \(\mathcal U\) is isomorphic to \(\mathcal D\): by construction \(\mathcal U\RK \mathcal D\), and on the other hand, the map \(f : P_\text{bd}(\lambda^+)\to Y\) defined by \(f(\sigma) = \sigma\cap \lambda\) pushes \(\mathcal U\) forward to \(\mathcal D\) so \(\mathcal D\RK \mathcal U\).
\end{proof}
\subsection{Weakly normal ultrafilters}\label{WeaklyNormalSection}
Another combinatorial generalization of the notion of a normal ultrafilter, due to Solovay and Ketonen \cite{Ketonen2}, is the notion of a weakly normal ultrafilter.
\begin{defn}\index{Weakly normal ultrafilter}
A uniform ultrafilter \(U\) on a cardinal \(\lambda\) is {\it weakly normal} if for any set \(A\in U\), if \(f : A\to\lambda\) is regressive, then there is some \(B\subseteq A\) such that \(B\in U\) and \(f[B]\) has cardinality less than \(\lambda\).   
\end{defn}
Solovay's definition of a weakly normal ultrafilter applied only to regular cardinals \(\lambda\), asserting that every regressive function on \(\lambda\) is {\it bounded} on a set of full measure. The generalization of the concept of weak normality to singular cardinals is due to Ketonen.
\begin{lma}
Suppose \(U\) is a uniform ultrafilter on a cardinal \(\lambda\). Then the following are equivalent:
\begin{enumerate}[(1)]
\item \(U\) is weakly normal.
\item Suppose \(\langle A_\alpha : \alpha < \lambda\rangle\) is a sequence of subsets of \(\lambda\) such that \(\bigcap_{\alpha\in \sigma} A_\alpha\in U\) for all nonempty \(\sigma\in P_\lambda(\lambda)\). Then \(\triangle_{\alpha <\lambda} A_\alpha\in U\).\qed
\end{enumerate}
\end{lma}
\begin{cor}\index{Weakly normal ultrafilter!on a regular cardinal}A uniform ultrafilter on a regular cardinal is weakly normal if and only if it is closed under decreasing diagonal intersections.\qed\end{cor}

Weakly normal ultrafilters on regular cardinals have a simple characterization in terms of their ultrapowers:
\begin{lma}\label{RegWeaklyNormal}
Suppose \(\lambda\) is a regular cardinal. An ultrafilter \(U\) on \(\lambda\) is weakly normal if and only if \(\id_U = \sup j_U[\lambda]\).
\begin{proof}
Suppose \(U\) is weakly normal. Since \(U\) is a tail uniform ultrafilter on \(\lambda\), \(\id_U > j_U(\alpha)\) for all \(\alpha < \lambda\). We will show that \(j_U[\lambda]\) is cofinal in \(\id_U\), which proves \(\id_U = \sup j_U[\lambda]\). Suppose \(\xi < \id_U\). Then \(\xi = [f]_U\) for some \(f : \lambda\to \lambda\) that is regressive on a set in \(U\). Since \(U\) is weakly normal, there is a set \(A\in U\) such that \(|f[A]| < \lambda\). Since \(\lambda\) is regular, \(f[A]\) is bounded below \(\lambda\). Fix \(\alpha < \lambda\) such that \(f(\xi) < \alpha\) for all \(\xi \in A\). Then \([f]_U < j_U(\alpha)\).

Conversely suppose \(\id_U = \sup j_U[\lambda]\). Since \(\id_U > j_U(\alpha)\) for all \(\alpha <\lambda\), \(\delta_U \geq \lambda\), and hence \(U\) is tail uniform. Since \(\lambda\) is regular, it follows that \(\lambda\) is uniform. Next, suppose \(A\in U\) and \(f :A\to \lambda\) is regressive. Then \([f]_U <\id_U\). Since \(j_U[\lambda]\) is cofinal in \(\id_U\), fix \(\alpha < \lambda\) with \([f]_U < j_U(\alpha)\). Then for a set \(B\in U\) with \(B\subseteq A\), \(f(\beta) < \alpha\) for all \(\beta\in B\). In particular, \(f\) takes fewer than \(\lambda\) values on \(B\). 
\end{proof}
\end{lma}

\cref{RegWeaklyNormal} yields the main source of weakly normal ultrafilters on regular cardinals:

\begin{cor}\label{DerivedWN}
Suppose \(j : V\to M\) is an elementary embedding and \(\lambda\) is a regular cardinal such that \(\sup j[\lambda] < j(\lambda)\). Then the ultrafilter on \(\lambda\) derived from \(j\) using \(\sup j[\lambda]\) weakly normal.
\end{cor}

To help motivate the concept of weak normality on singular cardinals, let us explain its relationship to an isomorphism invariant notion:
\begin{defn}\index{Minimal ultrafilter}
Suppose \(\lambda\) is an infinite cardinal. An ultrafilter \(U\) is {\it \(\lambda\)-minimal} if \(\lambda_U = \lambda\) and for any \(W\sRK U\), \(\lambda_W < \lambda\).
\end{defn}

If \(2^\lambda = \lambda^+\), there is a \(\lambda\)-minimal (countably incomplete) ultrafilter on \(\lambda\), according to a result of Comfort-Negrepontis \cite{Negrepontis}. On the other hand, the existence of a weakly normal ultrafilter (with no completeness assumptions) implies the existence of an inner model with a measurable cardinal \cite{Donder}. Weakly normal ultrafilters, however, are the revised Rudin-Keisler analog (\cref{rRKDef}) of \(\lambda\)-minimal ones:

\begin{lma}\label{WeaklyNormalrRK}
An ultrafilter \(U\) on a cardinal \(\lambda\) is weakly normal if and only if \(\lambda_U = \lambda\) and for all \(W\rRK U\), \(\lambda_W < \lambda\).\qed
\end{lma}

\cref{WeaklyNormalrRK} yields a generalization of Scott's theorem that every countably complete ultrafilter has a derived normal ultrafilter: 
\begin{cor}
If \(Z\) is a countably complete uniform ultrafilter on \(\lambda\), there is a weakly normal ultrafilter \(U\) on \(\lambda\) such that \(U\RK Z\).
\begin{proof}
Since \(\rRK\) is wellfounded on countably complete ultrafilters, there is a countably complete ultrafilter \(U\) that is \(\rRK\)-minimal with the property that \(\lambda_U = \lambda\) and \(U\leq_{\textnormal{RK}} Z\). Then \(U\) satisfies the conditions of \cref{WeaklyNormalrRK}: if \(W\rRK U\), then \(W \leq_{\textnormal{RK}} Z\), so by the \(\rRK\)-minimality of \(U\), it must be the case that \(\lambda_W < \lambda\)
\end{proof}
\end{cor}

The following theorem shows that every countably complete \(\lambda\)-minimal ultrafilter is isomorphic to a weakly normal ultrafilter.
\begin{prp}\label{MinimalIncompressible}
A countably complete uniform ultrafilter \(U\) on a cardinal \(\lambda\) is weakly normal if and only if it is \(\lambda\)-minimal and incompressible.
\begin{proof}
Suppose \(U\) is weakly normal. To see \(U\) is incompressible, note that any function that is regressive on a set in \(U\) takes less than \(\lambda\)-many values on a set in \(U\), and hence is not one-to-one. To see \(U\) is \(\lambda\)-minimal, suppose \(W\sRK U\) and we will show that \(\lambda_W < \lambda\). Since \(W\RK U\), \(W\) is countably complete, and hence \(W\) is isomorphic to an incompressible ultrafilter. We can therefore assume without loss of generality that \(W\) is incompressible. Then by the key lemma about the strict Rudin-Keisler order on incompressible ultrafilters (\cref{IncomRK}) the fact that \(W\sRK U\) implies \(W \rRK U\). Now by \cref{WeaklyNormalrRK}, \(\lambda_W < \lambda\).

Conversely suppose \(U\) is \(\lambda\)-minimal and incompressible. Suppose \(W\rRK U\), and we will show \(\lambda_W < \lambda\). We can then conclude that \(U\) is weakly normal using \cref{WeaklyNormalrRK}. Since \(U\) is incompressible, \(W\rRK U\) implies \(W\sRK U\) (\cref{IncomRK0}, essentially the definition of incompressibility). Therefore by the definition of \(\lambda\)-minimality, \(\lambda_W < \lambda\), as desired.    
\end{proof}
\end{prp}

It is not clear to us whether \cref{MinimalIncompressible} can be proved without the assumption of countable completeness, though of course countable completeness is not required if \(\lambda\) is regular.

The following characterization of weak normality is the one that is most relevant to our investigations of supercompactness.

\begin{prp}\label{WeaklyNormalGenerator}
Suppose \(\lambda\) is an infinite cardinal. A countably complete ultrafilter \(U\) on \(\lambda\) is weakly normal if and only if \(\id_U\) is the unique generator of \(j_U\) that lies above \(j(\delta)\) for all \(\delta < \lambda\).
\end{prp}

For the proof, we will need an obvious lemma:
\begin{lma}\label{UniformGenerator}
Suppose \(\lambda\) is an infinite cardinal.  An ultrafilter \(U\) on \(\lambda\) is uniform if and only if \(\id_U \notin H^{M_U}(j_U[V]\cup j_U(\delta))\) for any \(\delta < \lambda\). \qed
\end{lma}

\begin{proof}[Proof of \cref{WeaklyNormalGenerator}]
We begin with the forwards direction. Suppose \(U\) is weakly normal. 

We first show that for any ordinal \(\xi\) such that \(\xi < \id_U\), \(\xi\in H^{M_U}(j_U[V]\cup j_U(\delta))\) for some \(\delta < \lambda\). Assume not, towards a contradiction. Let \(W\) be the tail uniform ultrafilter derived from \(j_U\) using \(\xi\). Then \(W\rRK U\), as witnessed by the factor embedding \(k :M_W\to M_U\). By \cref{WeaklyNormalrRK}, it follows that \(W\) is not a uniform ultrafilter on \(\lambda\), and so by \cref{UniformGenerator}, there is some \(\delta < \lambda\) such that \(\xi\in H^{M_W}(j_W[V]\cup j_W(\delta))\). It follows that \(\xi\in H^{M_U}(j_U[V]\cup j_U(\delta))\). 

Next we show that \(\id_U\) is a generator of \(j_U\). Since \(U\) is uniform, \cref{UniformGenerator} implies \(\id_U\notin H^{M_U}(j_U[V]\cup j_U(\delta))\) for any \(\delta < \lambda\). But by the previous paragraph, for all \(\xi < \id_U\), \(\xi \in H^{M_U}(j_U[V]\cup j_U(\delta))\) for some \(\delta < \lambda\). Thus \(\id_U\notin H^{M_U}(j_U[V]\cup\vec \xi)\) for any \(\vec \xi \in [\id_U]^{<\omega}\). In other words, \(\id_U\) is a generator of \(j_U\). By the previous paragraph, \(\id_U\) is clearly the unique generator above \(j_U(\delta)\) for all \(\delta < \lambda\).

We now turn to the converse. Assume \(\id_U\) is the unique generator of \(j_U\) that lies above \(j(\delta)\) for all \(\delta < \lambda\). We will show \(U\) is weakly normal by verifying the conditions of \cref{MinimalIncompressible}. Since \(\id_U\) is a generator, \(U\) is incompressible. Since \(M_U\) is wellfounded, there is a least ordinal that does not belong to \(H^{M_U}(j_U[V]\cup j_U(\delta))\) for any \(\delta\), and clearly this ordinal is a generator of \(j_U\) that lies above \(j_U(\delta)\) for all \(\delta < \lambda\). Thus it must equal \(\id_U\). In other words, for any \(\xi < \id_U\), \(\xi\in H^{M_U}(j_U[V]\cup j_U(\delta))\) for some \(\delta < \lambda\).

Fix an ultrafilter \(W\) on \(\lambda\) such that \(W\rRK U\). We will show \(\lambda_W < \lambda\), verifying the second condition of \cref{MinimalIncompressible}. Let \(k : M_W\to M_U\) be an elementary embedding with \(k\circ j_W = j_U\) and \(k(\id_W) < \id_U\). Then by the previous paragraph, \(k(\id_W)\in H^{M_U}(j_U[V]\cup j_U(\delta))\) for some \(\delta < \lambda\). It follows that \(\id_W\in H^{M_W}(j_W[V]\cup j_W(\delta))\) (by the proof of \cref{WidthLemma}). By \cref{UniformGenerator}, this implies \(W\) is not uniform on \(\lambda\), or in other words, \(\lambda_W < \lambda\).
\end{proof}

Using \cref{WeaklyNormalGenerator}, we can prove the Dodd soundness of isonormal ultrafilters on \(\lambda = 2^{<\lambda}\).\index{Isonormal ultrafilter!Dodd soundness}
\begin{thm}\label{IsoSound}
	Suppose \(2^{<\lambda} = \lambda\). Then every isonormal ultrafilter \(U\) on \(\lambda\) is Dodd sound.
\end{thm}
\begin{proof}
Let \(j : V\to M\) be the ultrapower of the universe by \(U\). Since \(j\) is \(\lambda\)-supercompact, \(j\) is \(2^{<\lambda}\)-supercompact. By \cref{SpaceIrrel}, \(j\) is \(\lambda_*\)-sound where \(\lambda_* = \sup j[\lambda]\). 

We now show that \(j\) is \(\xi\)-sound where \(\xi\) is the least generator of \(j\) such that \(\xi\geq \lambda_*\). Since \(\lambda_*\) is closed under pairing, the \(\lambda_*\)-soundness of \(j\) implies that the extender \[E = E_j\restriction \lambda_* = \{(p,X) : p\in [\lambda_*]^{<\omega}\text{, }X\subseteq [\lambda]^{<\omega}\text{, and }p\in j(X)\}\] belongs to \(M_U\). Let \(j_E : V\to M_E\) be the associated extender embedding and let \(k : M_E\to M\) be the factor embedding. Then \[\textsc{crt}(k) = \min \{\alpha :\alpha\notin H^{M}(j[V]\cup \lambda_*)\} = \xi\] by the definition of a generator. Therefore \(j_E^\xi = j^\xi\). Moreover since \(M\) is closed under \(\lambda\)-sequences by \cref{UltrapowerSC}, \(j_E^{M} = j_E\restriction M\). Therefore \(j^\xi = j_E^\xi = (j_E^{M})^\xi\in M\), so \(j\) is \(\xi\)-sound.

By \cref{WeaklyNormalGenerator}, \(\xi = \id_U\). Therefore \(j\) is \(\id_U\)-sound, which implies that \(U\) is Dodd sound.
\end{proof}

We should point out that the assumption \(\lambda = 2^{<\lambda}\) is {\it necessary}:
\begin{lma}
Suppose \(\lambda\) is a cardinal that carries a Dodd sound ultrafilter \(U\). Then \(2^{<\lambda} = \lambda\).
\begin{proof}
Since \(U\) is Dodd sound, \(j_U\) is \(\id_U\)-sound. In particular, \(j_U\) is \(\sup j_U[\lambda]\)-sound since \(\sup j_U[\lambda] \leq \id_U\). Therefore by \cref{SpaceIrrel}, \(j_U\) is \(2^{<\lambda}\)-supercompact. By \cref{UFSuperBound}, \(j_U\) is not \(\lambda^+\)-supercompact. It follows that \(2^{<\lambda} < \lambda^+\), or in other words \(2^{<\lambda} = \lambda\).
\end{proof}
\end{lma}
\subsection{Solovay's Lemma}
A special case of our main theorem, \cref{IsoNormalThm}, was known long before our work.
\begin{thm}[Solovay's Lemma]\label{SolovayLemma}\index{Solovay's Lemma}
Suppose \(\lambda\) is a regular cardinal. Then there is a set \(B\subseteq P(\lambda)\) such that the following hold:
\begin{itemize}
\item For any family \(Y\) over \(\lambda\), any normal fine ultrafilter \(\mathcal U\) on \(Y\) concentrates on \(B\).
\item If \(\sigma\) and \(\tau\) are elements of \(B\) with the same supremum, then \(\sigma = \tau\).
\end{itemize}
\end{thm}
Before proving Solovay's Lemma, let us explain its relevance to isonormal ultrafilters. Essentially, Solovay's Lemma yields the ``regular case" of the key isomorphism theorem for isonormal ultrafilters (\cref{IsoNormalThm}):
\begin{cor}\label{SolovayCor}
Suppose \(\lambda\) is a regular cardinal, \(Y\) is a family over \(\lambda\), and \(\mathcal U\) is a nonprincipal normal fine ultrafilter on \(Y\). Then \(\mathcal U\) is isomorphic to the ultrafilter \[U = \{A\subseteq \lambda : \{\sigma\in Y : \sup \sigma\in A\}\in \mathcal U\}\]
Moreover, \(U\) is an isonormal ultrafilter.
\begin{proof}
To see \(\mathcal U\cong U\), let \(f : P(\lambda)\to \lambda+1\) be the function \(f(\sigma) = \sup \sigma\). Then \(f_*(\mathcal U) = U\) and by \cref{SolovayLemma}, \(f\) is one-to-one on a set in \(\mathcal U\).

To see \(U\) is isonormal, we must verify that \(U\) is weakly normal and \(j_U\) is \(\lambda\)-supercompact. The latter is trivial: \(j_\mathcal U\) is \(\lambda\)-supercompact by \cref{NormalFineChar}, and \(j_\mathcal U = j_U\) since \(\mathcal U\) and \(U\) are isomorphic. As for weak normality, by \cref{PushDerived}, \(U = f_*(\mathcal U)\) is the ultrafilter on \(\lambda\) derived from \(j_\mathcal U\) using \([f]_{\mathcal U}\) so \(U\) is weakly normal by \cref{DerivedWN}.
\end{proof}
\end{cor}

The proof of Solovay's lemma uses the observation that if \(j :V\to M\) is an elementary embedding, \(j[\lambda]\) is definable from the action of \(j\) on a stationary partition:\footnote{Solovay's published proof \cite{Solovay} uses the combinatorics of \(\omega\)-Jonsson algebras instead of stationary sets. Woodin rediscovered the proof using stationary sets, which was already known to Solovay.}

\begin{lma}\label{StationaryPartition}
Suppose \(\lambda\) is a cardinal, \(j : V\to M\) is an elementary embedding, and \(\mathcal P \subseteq P(\lambda)\) is a partition of \(S^\lambda_\omega = \{\alpha < \lambda:\textnormal{cf}(\alpha) = \omega\}\) into stationary sets. Then \[j[\mathcal P] = \{T \in j(\mathcal P) : T\text{ is stationary in }\sup j[\lambda]\}\]
\end{lma}
It is worth noting that \cref{StationaryPartition} is perfectly general; we really do allow \(j\) to be an arbitrary elementary embedding of \(V\).
\begin{proof}
Let \(\lambda_* = \sup j[\lambda]\). 
\begin{clm}\(j[\mathcal P]\subseteq \{T \in j(\mathcal P) : T\text{ is stationary in }\lambda_*\}\).\end{clm}
\begin{proof}Fix \(S\in \mathcal P\). We will show that \(j(S)\) intersects every closed cofinal subset of \(\lambda_*\). Suppose \(C\subseteq \lambda_*\) is closed cofinal in \(\lambda_*\). Then \(j^{-1}[C]\) is \(\omega\)-closed cofinal in \(\lambda\). Since \(S\) is a stationary subset of \(S^\lambda_\omega\), \(S\cap j^{-1}[C]\neq\emptyset\). But \(j(S)\cap C =j(S)\cap C \supseteq j[S\cap j^{-1}[C]]\neq \emptyset\). So \(j(S)\cap C\neq \emptyset\), as desired.
\end{proof}
\begin{clm} \(\{T \in j(\mathcal P) : T\text{ is stationary in }\lambda_*\}\subseteq j[\mathcal P]\). \end{clm}
\begin{proof}
Fix \(T\in j(\mathcal P)\) such that \(T\) is stationary in \(\lambda_*\). We will show that there is some \(S\in \mathcal P\) such that \(j(S) = T\). Since \(j[\lambda]\) is \(\omega\)-closed cofinal in \(\lambda_*\), \(T\cap j[\lambda]\neq \emptyset\). Take \(\xi < \lambda\) such that \(j(\xi)\in T\). Since \(j(\xi)\in T\subseteq j(S^\lambda_\omega)\), \(\xi\in S^\lambda_\omega\). Therefore \(\xi\in S\) for some \(S\in \mathcal P\), since \(\bigcup \mathcal P = S^\lambda_\omega\). Now \(j(\xi)\in j(S)\cap T\). Therefore \(j(S)\) and \(T\) are not disjoint, so since \(j(\mathcal P)\) is a partition, \(j(S) = T\), as desired.
\end{proof}
Combining the two claims yields the lemma.
\end{proof}

\cref{StationaryPartition} leads to a characterization of supercompactness that looks surprisingly weak:

\begin{cor}\label{SCStatCor}\index{Supercompactness!vs. stationary correctness}
Suppose \(j : V\to M\) is an elementary embedding and \(\lambda\) is a regular cardinal. The following are equivalent:
\begin{enumerate}[(1)]
\item \(j\) is \(\lambda\)-supercompact.
\item \(M\) is correct about stationary subsets of \(\lambda_* = \sup j[\lambda]\).
\end{enumerate}
\begin{proof}

{\it (1) implies (2):} Assume \(j\) is \(\lambda\)-supercompact. Suppose \(M\) satisfies that \(S\) is stationary in \(\lambda_*\), and we will show that \(S\) is truly stationary in \(\lambda_*\). Fix a closed cofinal set \(C\subseteq \lambda_*\). We will show \(S\cap C\neq \emptyset\). Note that \(C\cap j[\lambda]\in M\) by \cref{SupercompactClosure} (3). Let \(E\) be the closure of \(C\cap j[\lambda]\) in \(\lambda_*\). Then \(E\in M\), \(E\subseteq C\), and \(E\) is closed cofinal in \(\lambda_*\). Since \(E\in M\) and \(S\) is stationary from the perspective of \(M\), \(S\cap E\neq \emptyset\). In particular, \(S\cap C\neq \emptyset\).

{\it (2) implies (1):} Since \(\lambda\) is regular, there is a partition \(\mathcal P\) of \(S^\lambda_\omega\) into stationary sets such that \(|\mathcal P| = \lambda\). By \cref{StationaryPartition}, \(j[\mathcal P] = \{T\in j(\mathcal P) : T\text{ is stationary in }\lambda_*\}\), which is definable over \(M\) since \(M\) is correct about stationary subsets of \(\lambda_*\). Thus \(j\) is \(\mathcal P\)-supercompact, so by \cref{XSC}, \(j\) is \(\lambda\)-supercompact, as desired.
\end{proof}
\end{cor}
Of course the implication from (1) to (2) is not very surprising, but it allows us to restate \cref{StationaryPartition} in a useful way:

\begin{cor}\label{MStationaryPartition}
	Suppose \(\lambda\) is a regular cardinal, \(j : V\to M\) is a \(\lambda\)-supercompact elementary embedding, and \(\langle S_\alpha : \alpha < \lambda\rangle\) is a partition of \(S^\lambda_\omega\) into stationary sets. Let \(\langle T_\beta : \beta < j(\lambda)\rangle = j(\langle S_\alpha : \alpha < \lambda\rangle)\). Then \(j[\lambda] = \{\beta < j(\lambda) : M\vDash  T_\beta\text{ is stationary in }\lambda_*\}\).\qed
\end{cor}

We now prove Solovay's Lemma.
\begin{proof}[Proof of \cref{SolovayLemma}]
Let \(\langle S_\alpha : \alpha < \lambda\rangle\) be a partition of \(S^\lambda_\omega = \{\alpha < \lambda:\text{cf}(\alpha) = \omega\}\) into stationary sets. Let \[B = \{\sigma\subseteq \lambda : \sigma = \{\beta < \lambda : S_\beta\text{ is stationary in sup}(\sigma)\}\}\]
By construction, any two elements of \(B\) with the same supremum are equal.

To finish, suppose \(Y\) is a family over \(\lambda\) and \(\mathcal U\) is a normal fine on \(Y\). We must show that \(\mathcal U\) concentrates on \(B\), or equivalently, that \(\id_\mathcal U\in j_\mathcal U(B)\). Since \(\id_\mathcal U = j_\mathcal U[\lambda]\) (\cref{NormalFineChar}), this amounts to showing \[j_\mathcal U[\lambda] = \{\beta < j_\mathcal U(\lambda) : M_\mathcal U\vDash j_\mathcal U(S)_\beta\text{ is stationary in sup} j_\mathcal U[\lambda] \}\] which is of course a consequence of \cref{MStationaryPartition}.
\end{proof}

Another corollary of Solovay's Lemma is Woodin's proof of the Kunen Inconsistency Theorem:
\begin{thm}\label{KunenInconsistency0}\index{Kunen Inconsistency Theorem!Woodin's proof}
Suppose \(j : V\to M\) is an elementary embedding, \(\iota\) is a regular cardinal, \(j\) is \(\iota\)-supercompact, and \(j(\iota) = \sup j[\iota]\). Then \(j\restriction \iota + 1\) is the identity.
\end{thm}
\begin{proof}
Let \(\langle S_\alpha : \alpha <\iota\rangle\) be a partition of \(S^\iota_\omega\) into stationary sets. By \cref{MStationaryPartition}, and using the fact that \(j(\iota) = \sup j[\iota]\), 
\[j[\iota] = \{\beta < j(\iota) : M\vDash j(S)_\beta\text{ is stationary in }j(\iota)\} = j(\iota)\]
But this means \(j\restriction \iota + 1\) is the identity, as desired.
\end{proof}
Applying \cref{KunenInconsistency0} at \(\iota  = \lambda^+\) where \(\lambda\) is the first fixed point of \(j\) above \(\textsc{crt}(j)\) yields a proof of the Kunen Inconsistency  (\cref{KunenInconsistency}).
\subsection{Supercompactness and singular cardinals}\label{SingularSolovaySection}
In this section, we finish the proof of \cref{IsoNormalThm}. We do this by proving an analog of Solovay's Lemma at singular cardinals.\index{Solovay's Lemma!at singular cardinals} One basic issue, however, is that \cref{SolovayLemma} itself cannot generalize: in fact, if \(\lambda\) is a singular cardinal, \(Y\) is a family over \(\lambda\), and \(\mathcal U\) is a normal fine ultrafilter on \(Y\), then the supremum function is {\it not} one-to-one on any set in \(\mathcal U\).
\begin{prp}\label{NonSolovay}
Suppose \(\lambda\) is a cardinal of cofinality \(\iota\), \(Y\) is a family over \(\lambda\), and \(\mathcal U\) is a normal fine ultrafilter on \(Y\). Define \(f : Y\to \lambda + 1\) by
\[f(\sigma) = \sup \sigma\]
Define \(g : Y\to \iota + 1\) by
\[g(\sigma) = \sup (\sigma\cap \iota)\]
Then \(f_*(\mathcal U) \cong g_*(\mathcal U)\).
\end{prp}
It is a bit easier to prove the following equivalent statement first (which in any case turns out to be more useful):
\begin{prp}\label{iota*lemma}
Suppose \(j : V\to M\) is an elementary embedding and \(\lambda\) is a cardinal of cofinality \(\iota\). Then \(\sup j[\lambda]\) and \(\sup j[\iota]\) are interdefinable in \(M\) from parameters in \(j[V]\).
\begin{proof}
Let \(h : \iota\to \lambda\) be an increasing cofinal function. Then \[\sup j[\lambda] = \sup j[h[\iota]] = \sup j(h)\circ j[\iota] = \sup j(h)[\sup j[\iota]]\] Therefore \(\sup j[\lambda]\) is definable in \(M\) from \(j(h)\) and \(\sup j[\iota]\). Moreover, \[\sup j[\iota] = \sup j(h)^{-1}[\sup j[\lambda]]\] so \(\sup j[\iota]\) is definable in \(M\) from \(j(h)\) and \(\sup j[\lambda]\).
\end{proof}
\end{prp}
\begin{proof}[Proof of \cref{NonSolovay}]
Let \(j :V\to M\) be the ultrapower of the universe by \(\mathcal U\). Then (by \cref{PushDerived})
\(f_*(\mathcal U)\) is the ultrafilter on \(\lambda+ 1\) derived from \(j\) using \([f]_\mathcal U = \sup j[\lambda]\) and \(g_*(\mathcal U)\) is the ultrafilter on \(\iota + 1\) derived from \(j\) using \([g]_\mathcal U = \sup j[\iota]\).
By \cref{iota*lemma}, \[H^M(j[V]\cup \{\sup j[\lambda]\}) = H^M(j[V]\cup \{\sup j[\iota]\})\] But \[(M_{f_*(\mathcal U)}, j_{f_*(\mathcal U)})\cong (H^M(j[V]\cup \{\sup j[\lambda]\}),j)\cong (M_{g_*(\mathcal U)}, j_{g_*(\mathcal U)})\]
It follows that \(f_*(\mathcal U) \cong g_*(\mathcal U)\).
\end{proof}

\begin{cor}
Suppose \(\lambda\) is a cardinal of cofinality \(\iota\), \(Y\) is a family over \(\lambda\), and \(\mathcal U\) is a normal fine ultrafilter on \(Y\). Then there is a set \(B\in \mathcal U\) on which the supremum function takes at most \(\iota\)-many values.
\begin{proof}
Let \(f : Y\to \lambda\) be the supremum function. Since \(f_*(\mathcal U)\) is isomorphic to an ultrafilter on \(\iota + 1\), \(f\) takes at most \(\iota\)-many values on a set in \(\mathcal U\).
\end{proof}
\end{cor}

What we show instead is that an analog of \cref{StationaryPartition} holds:
\begin{thm}\label{GeneralizedSolovay}
Suppose \(\lambda\) is a cardinal and \(j : V\to M\) is a \(\lambda\)-supercompact elementary embedding. Let \(\theta\) be the least generator of \(j\) with \(\theta\geq\sup j[\lambda]\). Then \[j[\lambda] \in H^M(j[V]\cup \{\theta\})\] Moreover if \(\sup j[\lambda] < j(\lambda)\), then \(\theta < j(\lambda)\).
\end{thm}
As a corollary, we prove the second of the main theorems of this section:
\begin{thm}\label{IsoNormalThm}
	Suppose \(U\) is a nonprincipal ultrafilter. Then \(U\) is isonormal if and only if \(U\) is the incompressible ultrafilter isomorphic to a normal fine ultrafilter. In particular, every normal fine ultrafilter is isomorphic to a unique isonormal ultrafilter.
\end{thm}
\begin{proof}
We begin with the forward direction, which turns out to follow from \cref{MinimalIncompressible}. Suppose \(U\) is an isonormal ultrafilter on a cardinal \(\lambda\). We will show that \(U\) is incompressible and isomorphic to a normal fine ultrafilter on \(P_\text{bd}(\lambda)\). Since \(U\) is weakly normal, \cref{MinimalIncompressible} implies \(U\) is incompressible. 

Since \(U\) is uniform on \(\lambda\), \(\sup j_U[\lambda] < j_U(\lambda)\) and thus \(j_U\restriction \lambda\in j_U(P_\text{bd}(\lambda))\). Let \(\mathcal U\) be the ultrafilter on \(P_\text{bd}(\lambda)\) derived from \(j_U\) using \(j_U\restriction \lambda\). Then \(\mathcal U\RK U\) and \(\mathcal U\) is a normal fine ultrafilter on \(P_\text{bd}(\lambda)\) by \cref{NormalFineChar}. It follows that \(j_\mathcal U\) is \(\lambda\)-supercompact, and therefore \(\lambda_\mathcal U\geq \lambda\) by \cref{UFSuperBound}. Since \(U\) is weakly normal, \cref{MinimalIncompressible} implies \(U\) is \(\lambda\)-minimal and therefore \(\mathcal U\not \sRK U\). Since \(\mathcal U\RK U\) and \(\mathcal U\not \sRK U\), we must have \(\mathcal U\cong U\) (by definition).

Conversely, suppose \(U\) is incompressible and isomorphic to a normal fine ultrafilter, and we will show that \(U\) is isonormal. Since every normal fine ultrafilter is isomorphic to an element of \(\mathscr N\) (\cref{NIso}), for some cardinal \(\lambda\), \(U\) is isomorphic to a normal fine ultrafilter \(\mathcal U\) on \(P_\text{bd}(\lambda)\). In particular \(j_U = j_\mathcal U\) is \(\lambda\)-supercompact. To show that \(U\) is isonormal, it therefore suffices to show that \(U\) is a weakly normal ultrafilter on \(\lambda\). 

Let \(j: V\to M\) be the ultrapower of the universe by \(U\). Let \(\theta\) be the least generator of \(j\) with \(\theta\geq \sup j[\lambda]\). Since \(P_\text{bd}(\lambda)\in \mathcal U\), \(\sup j[\lambda] < j(\lambda)\), and so by \cref{GeneralizedSolovay}, \(\theta < j(\lambda)\). Since \(\theta\) is a generator of \(j = j_U\), \(\theta \leq \id_U\). In fact, we claim \(\id_U = \theta\). On the other hand, by \cref{GeneralizedSolovay}, \[M = H^M(j[V]\cup \{j[\lambda]\}) = H^M(j[V]\cup \{\theta\})\] The ultrapower theoretic characterization of incompressibility (\cref{IncomPower}) implies that \(\id_U\) is the least ordinal \(\alpha\) such that \(M= H^M(j[V]\cup \{\alpha\})\). Thus \(\id_U\leq\theta\). Hence \(\id_U = \theta\), as desired.

Since \(U\) is tail uniform (by the definition of incompressibility) and \(\id_U < j_U(\lambda)\), \(U\) is an ultrafilter on \(\lambda\). Since \(\id_U\) is the least generator of \(j\) above \(\sup j[\lambda]\), the characterization of weakly normal ultrafilters in terms of generators (\cref{WeaklyNormalGenerator}) implies that \(U\) is a weakly normal ultrafilter on \(\lambda\).
\end{proof}

We conclude this chapter by proving \cref{GeneralizedSolovay}. The proof relies on some basic notions from PCF theory.
\begin{defn}
Suppose \(\iota\) is an ordinal. We denote by \(J^\iota_\text{bd}\)\index{\(J_{\textnormal{bd}}\) (bounded ideal)} the ideal of bounded subsets of \(\iota\), omitting the superscript \(\iota\) when it is clear from context. If \(f\) and \(g\) are functions from \(\iota\) to \(\text{Ord}\), 
\begin{itemize}
	\item \(f <_\text{bd} g\) if \(\{\alpha < \iota : f(\alpha) \geq g(\alpha)\}\in J_\text{bd}\).\index{\(<_\text{bd}\) (domination mod bounded)}
	\item \(f =_\text{bd} g\) if \(\{\alpha < \iota : f(\alpha) \neq g(\alpha)\}\in J_\text{bd}\).
\end{itemize}
\end{defn}
\begin{defn}
Suppose \(C\) is a set of functions from \(\iota\) to \(\text{Ord}\). A function \(s: \iota\to \text{Ord}\) is an {\it exact upper bound}\index{Exact upper bound} of \(C\) if the following hold:
\begin{itemize}
\item For all \(f\in C\), \(f <_\text{bd} s\).
\item For all \(g <_\text{bd} s\), for some \(f\in C\), \(g <_\text{bd} f\).
\end{itemize}
\end{defn}

The following trivial fact plays a key role in the proof of \cref{GeneralizedSolovay}:

\begin{lma}\label{EUBUnique}
Suppose \(C\) is a set of functions from \(\iota\) to \(\textnormal{Ord}\) and \(s\) and \(t\) are exact upper bounds of \(C\). Then \(s =_\textnormal{bd} t\).
\begin{proof}
Suppose \(s\) and \(t\) are exact upper bounds of \(C\). Suppose towards a contradiction that \(s\neq_\text{bd} t\). Without loss of generality, we can assume that there is an unbounded set \(A\subseteq \iota\) such that \(s(\alpha) < t(\alpha)\) for all \(\alpha\in A\). Define \(g : \iota\to \text{Ord}\) by setting 
\[g(\alpha) = \begin{cases}
s(\alpha)&\text{if }\alpha\in A\\ 
0&\text{otherwise}\end{cases}\]
Then \(g < t\), so since \(t\) is an exact upper bound of \(C\), there is some \(f\in C\) such that \(g <_\text{bd} f\). Since \(s\) is an upper bound of \(C\), \(f <_\text{bd} s\). Therefore \(g <_\text{bd} s\). This contradicts that \(A = \{\alpha <\iota : g(\alpha) = s(\alpha)\}\) is unbounded in \(\iota\).
\end{proof}
\end{lma}
\begin{defn}\index{Scale}
If \(s : \iota\to \text{Ord}\) is a function and \(\delta\) is an ordinal, a {\it scale of length \(\delta\) in \(\prod_{\alpha < \iota} s(\alpha)\)} is a \(<_\text{bd}\)-increasing cofinal sequence \(\langle f_\alpha : \alpha < \delta\rangle\subseteq \prod_{\alpha < \iota} s(\alpha)\).
\end{defn}
Shelah's Representation Theorem \cite{AbrahamMagidor} states that if \(\lambda\) is a singular cardinal of cofinality \(\iota\), then there is a cofinal continuous sequence \(u : \iota\to \lambda\) such that \(\prod_{\alpha < \iota} u(\alpha)^+\) has a scale of length \(\lambda^+\). This is a deep theorem in the context of ZFC, but since we are assuming large cardinals, we will have enough SCH to get away with using only the following trivial version of Shelah's theorem:
\begin{lma}\label{ScaleExistence}
Suppose \(\lambda\) is a singular cardinal of cofinality \(\iota\) such that \(\lambda^\iota = \lambda^+\). Suppose \(\langle \delta_\alpha :\alpha <\iota\rangle\) is a sequence of regular cardinals cofinal in \(\lambda\). Then there is a scale of length \(\lambda^+\) in \(\prod_{\alpha < \iota} \delta_\alpha\).
\begin{proof}
We start by proving the standard fact that \(\mathbb P = (\prod_{\alpha < \iota} \delta_\alpha,<_\text{bd})\) is a \({\leq}\lambda\)-directed partial order. The proof proceeds in two steps. 

First, we prove that \(\mathbb P\) is \({<}\lambda\)-directed. Suppose \(\gamma < \lambda\) and \(\{g_i: i < \gamma\}\subseteq \mathbb P\). We will find a \(<_\text{bd}\)-upper bound \(g\) of \(\{g_i: i < \gamma\}\). Fix \(\alpha_0\) such that \(\gamma < \delta_{\alpha_0}\). For \(\alpha  < \iota\), define \[g(\alpha) = \begin{cases}
\sup_{i < \gamma} g_i(\alpha) + 1&\text{if }\alpha_0\leq \alpha\\
0&\text{otherwise}
\end{cases}\]
If \(\alpha_0\leq \alpha\), then \(\delta_\alpha\) is a regular cardinal greater than \(\gamma\), so \(\sup_{i < \gamma} g_i(\alpha) < \delta_\alpha\). Hence \(g\in \prod_{\alpha < \iota} \delta_\alpha\) and \(g\) is a \(<_\text{bd}\)-upper bound of \(\{g_i: i < \gamma\}\).

Second, we prove that \(\mathbb P\) is \(\lambda\)-directed. Fix \(\{g_i: i < \lambda\}\subseteq \mathbb P\). For \(\alpha<\iota\), let \(h_\alpha\in \mathbb P\) be a \(<_\text{bd}\)-upper bound of \(\{g_i : i < \delta_\alpha\}\). Finally let \(g\in \mathbb P\) be a \(<_\text{bd}\)-upper bound of \(\{g_i : i < \iota\}\). Then \(g\) is a \(<_\text{bd}\)-upper bound of \(\{g_i: i < \lambda\}\), as desired.

Enumerate \(\prod_{\alpha < \iota} \delta_\alpha\) as \(\{g_\xi : \xi < \lambda^+\}\). We define \(\langle f_\xi : \xi < \lambda^+\rangle\) recursively. If \(\langle f_\xi : \xi < \theta\rangle\) has been defined, choose a \(<_\text{bd}\)-upper bound \(f_\theta\in \mathbb P\) of \(\{ f_\xi : \xi < \theta\}\cup \{g_\xi\}\). (Such a function exists by the \(\lambda\)-directedness of \(\mathbb P\).) By construction \(\langle f_\xi : \xi < \lambda^+\rangle\) is a scale in \(\prod_{\alpha < \iota} \delta_\alpha\).
\end{proof}
\end{lma}

This concludes our summary of the basic notions from PCF theory used in the proof of \cref{GeneralizedSolovay}, which we now commence.

\begin{proof}[Proof of \cref{GeneralizedSolovay}]\index{Solovay's Lemma!at singular cardinals}
For the purposes of the proof, let us say that \(x\) is {\it weakly definable from \(y\) (in \(M\))} if \(x\) is definable in \(M\) from parameters in \(j[V]\cup \{y\}\), or in other words, \(x\in H^M(j[V]\cup \{y\})\). Note that weak definability is a transitive relation.
	
By \cref{StationaryPartition}, we may assume \(\lambda\) is a singular cardinal. Let \(\iota\) be the cofinality of \(\lambda\). Let \(\lambda_* = \sup j[\lambda]\).

\begin{clm}\label{JbdClm}
Suppose \(\langle \delta_\alpha :\alpha < \iota\rangle\) is an increasing sequence of regular cardinals cofinal in \(\lambda\). Let \(e\) be the equivalence class of \(\langle\sup j[\delta_\alpha] : \alpha < \iota\rangle\) modulo \(J_\textnormal{bd}\). Then \(j[\lambda]\) is weakly definable from \(e\) and \(j\restriction \iota\) in \(M\).
\begin{proof}[Proof of \cref{JbdClm}]
%We first show that \( j\restriction \iota\) is weakly definable from \(e\). To see this, note that \(\lambda_*\) is of course definable from \(e\): \(\lambda_* = \limsup s\) for any \(s\in e\). By \cref{iota*lemma}, \(\sup j[\iota]\) is weakly definable from \(\lambda_*\). Moreover by Solovay's Lemma (\cref{MStationaryPartition}), \(j[\iota]\) is weakly definable from \(\sup j[\iota]\). Finally \(j\restriction \iota\) is definable from \(j[\iota]\) as the inverse of its transitive collapse. By the transitivity of weak definability, \(j\restriction \iota\) is weakly definable from \(e\).
	
Fix a sequence \(\langle S^\alpha : \alpha < \iota\rangle\) such that \(S^\alpha = \{S^\alpha_\beta : \beta < \delta_\alpha\}\) is a partition of \(S^{\delta_\alpha}_\omega\) into stationary sets. Note that \(\langle j(S^\alpha) : \alpha < \iota\rangle= j(\langle S^\alpha : \alpha < \iota\rangle)\circ j\restriction \iota\) is weakly definable from \(j\restriction \iota\).

Solovay's Lemma (\cref{MStationaryPartition}) implies that for all \(\alpha < \iota\), \(j[\delta_\alpha]\) is equal to the set \(\{\beta < j(\delta_\alpha) :M\vDash  j(S^\alpha)_\beta\text{ is stationary in }\sup j[\delta_\alpha]\}\).
It follows that
\begin{align*}\beta\in j[\lambda]&\iff \{\alpha < \iota :  M\vDash j(S^\alpha)_\beta\text{ is stationary in }\sup j[\delta_\alpha]\}\notin J_{\text{bd}}\\
&\iff \exists s\in e\ \{\alpha < \iota : M\vDash  j(S^\alpha)_\beta\text{ is stationary in }s(\alpha)\}\notin J_{\text{bd}}
\end{align*}
Thus \(j[\lambda]\) is weakly definable from \(e\) and \(\langle j(S^\alpha) : \alpha < \iota\rangle\). Since \(\langle j(S^\alpha) : \alpha < \iota\rangle\) is weakly definable from  \(j\restriction \iota\), this proves the claim.
\end{proof}
\end{clm}
It is not hard to see that \(j\restriction \iota\) is itself weakly definable from \(e\), but we will not need this. The following observation, however, will be crucial:
\begin{obs}\label{supobs}\(j\restriction \iota\) and \(j[\iota]\) are weakly definable from \(\sup j[\iota]\).\end{obs}
This is an immediate consequence of Solovay's Lemma (\cref{MStationaryPartition}).
 
Let \(\mathcal D\) be the normal fine ultrafilter on \(P(\iota)\) derived from \(j\) using \(j[\iota]\) and let \(k : M_\mathcal D\to M\) be the factor embedding. Let \(\lambda_\mathcal D = \sup j_\mathcal D[\lambda]\).\footnote{If \(\eta^\iota < \lambda\) for all \(\eta < \lambda\), then \(\lambda_\mathcal D= \lambda\), but we do not assume this.} By \cref{DerivedNF}, \(\textsc{crt}(k) > \iota\) and hence \(k(\lambda_\mathcal D) = \sup k[\lambda_\mathcal D] =\lambda_*\).

\begin{obs}\label{kobs}
	\(k[M_\mathcal D]\) consists of all \(x\in M\) that are weakly definable from \(\sup j[\iota]\).
\end{obs}
\cref{kobs} follows from the fact that \(k[M_\mathcal D] = H^M(j[V]\cup \{j[\iota]\})\) combined with the fact (\cref{supobs}) that \(j[\iota]\) and \(\sup j[\iota]\) are weakly definable from each other.

Let \[\theta = \sup k[\lambda_\mathcal D^{+M_\mathcal D}]\] The ordinal \(\theta\) will turn out to be the least generator of \(j\) above \(\lambda_*\). For now, let us just show that there is no smaller generator:
\begin{clm}\label{ThetaClm}
\(\theta\subseteq H^M(j[V]\cup \lambda_*)\).
\begin{proof}
Suppose \(\alpha < \theta\). The claim amounts to showing that \(\alpha\) is weakly definable from a finite set of ordinals below \(\lambda_*\). By the definition of \(\theta\), \(\alpha < k(\xi)\) for some \(\xi < \lambda_\mathcal D^{+M_\mathcal D}\). Fix a surjection \(p : \lambda_\mathcal D\to \xi\) with \(p\in M_\mathcal D\). \cref{kobs} implies \(k(p)\) is weakly definable from \(\sup j[\iota]\). Since \(k(p)\) is a surjection from \(\lambda_*\) onto \(k(\xi)\), for some \(\nu < \lambda_*\), \(\alpha = k(p)(\nu)\). Thus \(\alpha\) is weakly definable from \(\sup j[\iota]\) and \(\nu\), which both lie below \(\lambda_*\), proving the claim.
\end{proof}
\end{clm}

Fix a sequence \(\langle \delta_\alpha :\alpha < \iota\rangle\) of regular cardinals greater than \(\iota\) that is increasing and cofinal in \(\lambda\).
\begin{clm}
In \(M_\mathcal D\), there is a scale \(\vec f = \langle f_\alpha : \alpha < \lambda_\mathcal D^{+M_\mathcal D}\rangle\) in \(\prod_{\alpha < \iota} j_\mathcal D(\delta_\alpha)\).
\begin{proof}
Applying \cref{ScaleExistence} in \(M_\mathcal D\), it suffices to show that \(M_\mathcal D\) satisfies \(\lambda_\mathcal D^\iota = \lambda_\mathcal D^{+\mathcal D}\). By the critical sequence analysis given by the Kunen Inconsistency Theorem (\cref{Kunen2}), there is a \(\lambda\)-supercompact cardinal \(\kappa \leq \iota\) such that \(j_\mathcal D(\kappa) > \iota\). 
Thus \(\iota < j_\mathcal D(\kappa) < \lambda_*\) and \(j_{\mathcal D}(\kappa)\) is \(\lambda_\mathcal D^{+\mathcal D}\)-supercompact in \(M_\mathcal D\). 
By the local version of Solovay's theorem \cite{Solovay} (which appears as \cref{SolovaySCH}) applied in \(M_\mathcal D\), it follows that in \(M_\mathcal D\), \(\lambda_\mathcal D^\iota \leq (\lambda_\mathcal D^+)^{<j_\mathcal D(\kappa)} = \lambda_\mathcal D^{+M_\mathcal D}\), as desired.
\end{proof}
\end{clm}
\begin{clm}\label{EUBClm}
\(\langle\sup j[\delta_\alpha] : \alpha < \iota\rangle\) is an exact upper bound of \(k(\vec f)\restriction \theta\).
\end{clm}

Before proving \cref{EUBClm}, let us show how it implies the theorem. 

Let \(e\) be the equivalence class of \(\langle\sup j[\delta_\alpha] : \alpha < \iota\rangle\) modulo the bounded ideal on \(\iota\). Then \cref{EUBClm} and \cref{EUBUnique} imply that \(e\) is definable in \(M\) from the parameters \(\theta\) and \(k(\vec f)\). Thus by \cref{JbdClm}, \(j[\lambda]\) is weakly definable from \(\theta\), \(k(\vec f)\), and \(j\restriction \iota\).

Note that \(\lambda_*\) is definable in \(M\) from \(\theta\): \(\lambda_*\) is the largest \(M\)-cardinal below \(\theta\). By \cref{iota*lemma}, \(\sup j[\iota]\) is weakly definable from \(\lambda_*\) and hence from \(\theta\). Thus by \cref{supobs} \(j\restriction \iota\) is weakly definable from \(\theta\), and by by \cref{kobs}, \(k(\vec f)\) is weakly definable from \(\theta\). Combining this with the previous paragraph, \(j[\lambda]\) is weakly definable from \(\theta\) alone. This yields:
\[j[\lambda] \in H^M(j[V]\cup \{\theta\})\]

We now show \(\theta\) is the least generator of \(j\) above \(\lambda_*\). It suffices by \cref{ThetaClm} to show that \(\theta\) is a generator of \(j\). Assume towards a contradiction that this fails. Then \(\theta\in H^M(j[V]\cup \theta) = H^M(j[V]\cup \lambda_*)\) by \cref{ThetaClm}. Thus \(j[\lambda]\in H^M(j[V]\cup \lambda_*)\). Fix \(\xi < \lambda_*\) such that \(j[\lambda]\in H^M(j[V]\cup \{\xi\})\). Let \(W\) be the ultrafilter derived from \(j\) using \(\xi\). Then by \cref{DerivedSC}, \(j_W\) is \(\lambda\)-supercompact, yet \(\lambda_W < \lambda\), and this contradicts \cref{UFSuperBound}. Thus our assumption was false, and in fact \(\theta\) is a generator of \(j\).

Thus \(j[\lambda]\in H^M(j[V]\cup \{\theta\})\) where \(\theta\) is the least generator of \(j\) greater than or equal to \(\lambda_*\). To finish, we must show that if \(\lambda_* < j(\lambda)\) then \(\theta < j(\lambda)\). But \(\theta \leq \lambda_*^{+M}\) while \(j(\lambda)\) is a limit cardinal of \(M\) above \(\lambda_*\). Hence \(\lambda_*^{+M} < j(\lambda)\), as desired.

We now turn to the proof of \cref{EUBClm}. It will be important here that for any \(s:\iota\to \text{Ord}\), \(k(s) = k\circ s\) since \(\textsc{crt}(k) > \iota\).
\begin{proof}[Proof of \cref{EUBClm}]
We first show that for all \(\nu < \theta\), \[k(\vec f)_\nu <_\text{bd} \langle \sup j[\delta_\alpha] : \alpha < \iota\rangle\] For any \(\nu < \theta\), there is some \(\xi < \lambda_\mathcal D^{+M_\mathcal D}\) such that \(\nu < k(\xi)\). Therefore \[k(\vec f)_\nu <_\text{bd} k(\vec f)_{k(\xi)} = k(f_\xi)\] Hence it suffices to show that for any \(\xi < \lambda_\mathcal D^{+M_\mathcal D}\), \(k(f_\xi) < \langle  \sup j[\delta_\alpha] : \alpha < \iota\rangle\). For all \(\alpha < \iota\), we have that \(\delta_\alpha \) is a regular cardinal above \(\iota\). By \cref{SolovayCor}, \(\lambda_\mathcal D = \iota\), so since ultrapower embeddings are continuous at regular cardinals above their size (\cref{UFContinuity}), \[j_\mathcal D(\delta_\alpha) = \sup j_\mathcal D[\delta_\alpha]\]  Since \(f_\xi\in \prod_{\alpha < \iota}j_\mathcal D(\delta_\alpha)\), we therefore have \(f_\xi(\alpha) < \sup j_\mathcal D[\delta_\alpha]\) and hence \(k(f_\xi)(\alpha) = k(f_\xi(\alpha)) < \sup j[\delta_\alpha]\) for all \(\alpha < \iota\), as desired.

We finish by showing that for any \(g :\iota\to \text{Ord}\) such that \(g<_\text{bd}\langle \sup j[\delta_\alpha] : \alpha < \iota\rangle\), there is some \(\xi < \lambda_\mathcal D^{+M_\mathcal D}\) such that \(g <_\text{bd} k(f_\xi)\). For \(\alpha < \iota\), let \(h(\alpha) < \delta_\alpha\) be least such that \(g(\alpha) \leq j(h(\alpha))\). Then \(j_\mathcal D\circ h\in M_\mathcal D\) (since \(M_\mathcal D\) is closed under \(\iota\)-sequences by \cref{UltrapowerSC}). Since \(\langle f_\xi : \xi < \lambda^+\rangle\) is cofinal, in \(\prod_{\alpha < \iota}j_\mathcal D(\delta_\alpha)\), there is some \(\xi < \lambda_\mathcal D^{+M_\mathcal D}\) such that \(j_\mathcal D\circ h <_\text{bd} f_\xi\) . It follows that \[g \leq j\circ h = k\circ j_\mathcal D\circ h = k(j_\mathcal D\circ h) <_\text{bd} k(f_\xi)\] as desired.
\end{proof}
This completes the proof of \cref{GeneralizedSolovay}.
\end{proof}
\chapter{The Rudin-Frol\'ik Order}\label{RFChapter}
\section{Introduction}
\subsection{Ultrafilters on the least measurable cardinal}
This chapter is motivated by a single simple question. \cref{MOChapter}
established the linearity of the Mitchell order on normal ultrafilters assuming
UA. As a consequence, the least measurable cardinal \(\kappa\) carries a unique
normal ultrafilter. But what are the other countably complete ultrafilters on
\(\kappa\)? The following theorem of Kunen \cite{Kunen} answers this question
under a hypothesis that is much more restrictive than UA:
\begin{thm}[Kunen]\label{KunenUlts}
Suppose \(U\) is a normal ultrafilter on \(\kappa\) and \(V = L[U]\). Then every
countably complete ultrafilter is isomorphic to \(U^n\) for some \(n < \omega\).
\end{thm}
Here \(U^n\) is the ultrafilter on \([\kappa]^n\) generated by sets of the form
\([A]^n\) where \(A\in U\). An even stronger theorem of Kunen characterizes
every elementary embedding of the universe when \(V = L[U]\):
\begin{thm}[Kunen] Suppose \(V = L[U]\) for some normal ultrafilter \(U\). Then
any elementary embedding \(j :V\to M\) is an iterated ultrapower of \(U\).
\end{thm}
Kunen's proofs of these theorems rely heavily on the structure of \(L[U]\), so
much so that it might seem unlikely UA alone could imply analogous results. The
results of this chapter, however, show that UA does just as well:
\begin{repthm}{KunenUA}[UA] Let \(\kappa\) be the least measurable cardinal.
Then there is a unique normal ultrafilter \(U\) on \(\kappa\), and every
countably complete ultrafilter is isomorphic to \(U^n\) for some \(n < \omega\).

\end{repthm}

\begin{repthm}{TransfiniteKunen}[UA] Let \(\kappa\) be the least measurable
cardinal. Let \(U\) be the unique normal ultrafilter on \(\kappa\), Then any
elementary embedding \(j : V\to M\) such that \(M = H^M(j[V]\cup j(\kappa))\) is
an iterated ultrapower of \(U\).
\end{repthm}
The requirement that \(M = H^M(j[V]\cup j(\kappa))\) is necessary because for
example there could be two measurable cardinals. (Actually one could make do
with the requirement that \(M = H^M(j[V]\cup j(\alpha))\) and there are no
measurable cardinals in the interval \((\kappa,\alpha]\).)

Thus there is an abstract generalization of Kunen's analysis of \(L[U]\) to
arbitrary models of UA. Far more interesting, however, is that this
generalization leads to the discovery of new structure high above the least
measurable cardinal.
\begin{defn}\index{Irreducible ultrafilter}
A nonprincipal countably complete ultrafilter \(U\) is {\it irreducible} if its
ultrapower embedding cannot be written nontrivially as a linear iterated
ultrapower.
\end{defn}
Irreducible ultrafilters arise in the generalization of Kunen's theorem, which
really factors into the following two theorems:
\begin{repthm}{LeastIrreds}[UA] Every irreducible ultrafilter on the least
measurable cardinal \(\kappa\) is isomorphic to the unique normal ultrafilter on
\(\kappa\).
\end{repthm}
\begin{repthm}{FactorizationThm}[UA] Every ultrapower embedding can be written
as a finite linear iterated ultrapower of irreducible ultrafilters.
\end{repthm}
The first of these theorems is highly specific to the least measurable cardinal,
but the second is a perfectly general fact: under UA, the structure of countably
complete ultrafilters in general can be reduced to the structure of irreducible
ultrafilters. The nature of irreducible ultrafilters in general is arguably the
most interesting problem raised by this dissertation, intimately related to the
theory of supercompactness and strong compactness under UA. %For example, let us
mention the following theorem:
%\begin{thm}[UA] Suppose \(\lambda\) is either a successor cardinal or a strong
%limit singular cardinal. Suppose \(U\) is a countably complete uniform
%ultrafilter on \(\lambda\). Then the following are equivalent:
%\begin{enumerate}[(1)] \item \(U\) is isomorphic to a normal fine ultrafilter.
%\item \(U\) is weakly normal and irreducible. \end{enumerate} \end{thm} This
%theorem appears as \cref{WNIrred} and requires the full strength of the
%supercompactness analysis of \cref{SCChapter1} and \cref{SCChapter2}. It
%remains open whether it can be generalized to the case that \(\lambda\) is
%inaccessible.
\subsection{Outline of \cref{RFChapter}}
We now outline the rest of this chapter.\\

\noindent {\sc\cref{RFSection}.} We introduce the fundamental Rudin-Frol\'ik
order, which measures how an ultrapower embedding can be factored as a finite
iterated ultrapower. We explain how the topological definition of the
Rudin-Frol\'ik order is related to the concept of an internal ultrapower
embedding (\cref{RFChar0}). We show that the Ultrapower Axiom is equivalent to
the directedness of the Rudin-Frol\'ik order on countably complete ultrafilters,
and we show that the Rudin-Frol\'ik order is not directed on ultrafilters on
\(\omega\).\\

\noindent {\sc\cref{MuSection}.} In this section, we answer the basic question,
characterizing the ultrafilters on the least measurable cardinal up to
isomorphism. It turns out that this can be done for all ultrafilters below the
least \(\mu\)-measurable cardinal. (In fact, the analysis extends quite a bit
further, but we have omitted this work from this dissertation.) Towards this, in
\cref{IrredMuSection}, we introduce irreducible ultrafilters and analyze the
irreducible ultrafilters up to isomorphism. We then prove that every ultrafilter
can be factored into finitely many irreducible ultrafilters in
\cref{FactorizationSection}.\\

\noindent {\sc\cref{RFStructureSection}.} In this section, we investigate the deeper
structural properties of the Rudin-Frol\'ik order assuming UA. We show that the
Rudin-Frol\'ik order satisfies the local ascending chain condition (\cref{ACC}),
which was actually required as a step in the irreducible factorization theorem.
We show that the Rudin-Frol\'ik order induces a lattice on the isomorphism types
of countably complete ultrafilters. This involves showing that every pair of
ultrapower embeddings has a minimum comparison, which we call a {\it pushout}.
In \cref{FiniteSection}, we use pushouts to prove the local finiteness of the
Rudin-Frol\'ik order: a countably complete ultrafilter has at most finitely many
Rudin-Frol\'ik predecessors assuming UA. Finally, in \cref{TranslationSection},
we study the structure of pushouts and their relationship to the minimal covers
of \cref{LinKetSection}. This involves the key notion of a {\it translation} of
ultrafilters.\\

\noindent {\sc\cref{InternalSection}.} In this section, we use the theory of comparisons
developed in \cref{RFStructureSection} to investigate a variant of the
generalized Mitchell order called the internal relation.
\section{The Rudin-Frol\'ik order}\label{RFSection}
Irreducible ultrafilters are most naturally studied in the setting of the {\it
Rudin-Frol\'ik order}, an order on ultrafilters due to Rudin and Frol\'ik
\cite{Frolik} that dates back to the study of ultrafilters by Mary Ellen Rudin's
school in the late 1960s. The structure of the Rudin-Frol\'ik order on countably
complete ultrafilters turns out to encapsulate many of the phenomena we have
been studying so far. For example, the Ultrapower Axiom is equivalent to the
statement that the Rudin-Frol\'ik order is directed, while irreducible
ultrafilters are simply the minimal elements of the Rudin-Frol\'ik order. The
deeper properties of this order developed in this chapter (especially the
existence of least upper bounds) will provide some of the key tools of the
supercompactness analysis.

In this section, we discuss the theory of the Rudin-Frol\'ik order without yet
restricting to countably complete ultrafilters. For this reason, this subsection
is a bit out of step with the rest of this dissertation, and the only fact that
will be truly essential going forward is the characterization of the
Rudin-Frol\'ik order on countably complete ultrafilters given by \cref{RFChar},
which the reader who is not interested in ultrafilter combinatorics can take as
the definition of the Rudin-Frol\'ik order on countably complete ultrafilters.
\begin{defn}\index{Discrete sequence of ultrafilters}
A sequence of ultrafilters \(\langle W_i : i\in I\rangle\) is {\it discrete} if
there is a sequence of pairwise disjoint sets \(\langle Y_i : i\in I\rangle\)
such that \(Y_i\in W_i\) for all \(i\in I\).
\end{defn}
Typically (for example, in \cref{RFDef}) we will consider discrete sequences of
ultrafilters that all lie on the same set \(X\). Then discreteness says these
ultrafilters can be simultaneously separated from each other.

\begin{defn}\label{RFDef}\index{Rudin-Frol\'ik order}
Suppose \(U\) is an ultrafilter on \(X\) and \(W\) is an ultrafilter on \(Y\).
The {\it Rudin-Frol\'ik order} is defined by setting \(U\D W\) if there is a set
\(I\in U\) and a discrete sequence of ultrafilters \(\langle W_i : i\in
I\rangle\) on \(Y\) such that \(W = U\text{-}\lim_{i\in I}W_i\).
\end{defn}

Recall that if \(U\) is an ultrafilter on \(X\), \(I\) is a set in \(U\), and
\(\langle W_i : i\in I\rangle\) is a sequence of ultrafilters on \(Y\), then the
\(U\)-sum of \(\langle W_i : i \in I\rangle\) is defined by
\[U\text{-}\sum_{i\in I} W_i = \{A \subseteq X\times Y : \{i \in I : A_i\in
W_i\}\in U\}\] The projection \(\pi^0 : X\times Y\to X\) defined by \(\pi^0(i,j)
= i\) satisfies \(\pi^0_*\left(U\text{-}\sum_{i\in I} W_i\right) = U\), and the
projection \(\pi^1 : X\times Y\to Y\) defined by \(\pi^1(i,j) = j\) satisfies
\[\pi^1_*\left(U\text{-}\sum_{i\in I} W_i\right) = U\text{-}\lim_{i\in I} W_i \]

The model-theoretic characterization of the Rudin-Frol\'ik order uses the
following lemma:

\begin{lma}\label{DiscreteSum}
Suppose \(U\) is an ultrafilter, \(I\in U\), and \(\langle W_i : i\in I\rangle\)
is a sequence of ultrafilters on \(Y\). Then the following are equivalent:
\begin{enumerate}[(1)]
\item There is a \(U\)-large set \(J\subseteq I\) such that \(\langle W_i : i\in
J\rangle\) is discrete.
\item \(\pi^1\) is one-to-one on a set in  \(U\text{-}\sum_{i\in I}W_i\).
\item \(U\text{-}\sum_{i\in I}W_i\cong U\text{-}\lim_{i\in I} W_i\).
\end{enumerate}
\begin{proof}
{\it (1) implies (2):} Fix \(J\in U\) contained in \(I\) and pairwise disjoint
sets \(\langle Y_i : i \in J\rangle\) with \(Y_i\in W_i\) for all \(i\in J\).
We will show \(\pi^1\) is one-to-one on a set in \(U\text{-}\sum_{i\in I}W_i\).
Let \[A = \{(i,j) : i\in J\text{ and }j\in Y_i\}\] Then \(A\in
U\text{-}\sum_{i\in I}W_i\) and \(\pi^1\) is one-to-one on \(A\). 

{\it (2) implies (1):} Fix \(A\in U\text{-}\sum_{i\in I}W_i\) on which \(\pi^1\)
is one-to-one. For each \(i\in I\), let \(Y_i = \{j\in Y : (i,j)\in A\}\). Since
\(\pi^1\) is one-to-one on \(A\), the sets \(Y_i\) are disjoint. Since \(A\in
U\text{-}\sum_{i\in I}W_i\), the set \(J = \{i\in I : Y_i\in W_i\}\) belongs to
\(U\). Thus \(J\in U\), \(J\subseteq I\), and \(\langle W_i : i\in J\rangle\) is
witnessed to be discrete by \(\langle Y_i : i\in J\rangle\), as desired.

{\it (2) implies (3):} Trivial.

{\it (3) implies (2):} By \cref{RKRigid}, if \(Z \cong Z'\) and \(f\) is such
that \(f_*(Z) =Z'\), then \(f\) is one-to-one on a set in \(Z\). Therefore since
\(\pi^1_*(U\text{-}\sum_{i\in I}W_i) = U\text{-}\lim_{i\in I} W_i\), \(\pi^1\)
is one-to-one on a set in \(U\text{-}\sum_{i\in I}W_i\).
\end{proof}
\end{lma}

\begin{cor}\label{RFSum}
	If \(U\) and \(W\) are ultrafilters, the following are equivalent:
	\begin{enumerate}[(1)]
		\item \(U\D W\).
		\item There exist \(I\in U\) and ultrafilters \(\langle W_i : i\in
		I\rangle\) on a set \(Y\) such that \(W \cong U\text{-}\sum_{i\in I}
		W_i\).
	\end{enumerate}
\begin{proof}
	{\it (1) implies (2):} Obvious from \cref{DiscreteSum}.
	
	{\it (2) implies (1):} The proof uses the fact that the Rudin-Frol\'ik order
	is isomorphism invariant, which should be easy enough to see from the
	definition.
	
	Let \(Y' = I\times Y\). Let \(f^i : Y\to Y'\) be the embedding defined by
	\(f^i(y) = (i,y)\), and let \(W'_i = f^i_*(W_i)\). Then \(W_i' \cong W_i\)
	and \(\langle W_i' : i \in I\rangle\) is discrete. We have \[W \cong
	U\text{-}\sum_{i\in I}W_i\cong U\text{-}\sum_{i\in I}W'_i\cong
	U\text{-}\lim_{i\in I}W_i'\] where the last isomorphism follows from
	\cref{DiscreteSum}. By the definition of the Rudin-Frol\'ik order \(U\D
	U\text{-}\lim_{i\in I} W_i'\), so by the isomorphism invariance of the
	Rudin-Frol\'ik order, \(U\D W\).
\end{proof}
\end{cor}

The following generalization of closeness to possibly illfounded models in our
view simplifies the theory of the Rudin-Frol\'ik order on countably incomplete
ultrafilters:
\begin{defn}
Suppose \(N\) and \(M\) are models of ZFC. A cofinal elementary embedding \(h :
N\to M\) is {\it close} to \(N\) if for all \(X\in N\) and all \(a\in M\) such
that \(M\vDash a\in h(X)\), the \(N\)-ultrafilter on \(X\) derived from \(h\)
using \(a\) belongs to \(N\).
\end{defn}
It is really not quite accurate to say that this \(N\)-ultrafilter belongs to
\(N\); we really mean that it is the extension of a point in \(N\).

\begin{lma}\label{CloseChar}
If \(h : N\to M\) is close and \(M = H^M(h[N]\cup \{a\})\) for some \(a\in M\),
then there is an ultrafilter \(Z\) of \(N\) and an isomorphism \(k : M_Z^N\to
M\) such that \(k\circ j_Z^N = h\).\qed
\end{lma}

\begin{cor}\label{RFChar0}
If \(U\) and \(W\) are ultrafilters, the following are equivalent:
\begin{enumerate}[(1)]
\item \(U\D W\).
\item There is a close embedding \(h : M_U\to M_W\) such that \(h\circ j_U =
j_W\).
\end{enumerate}
\begin{proof}[Sketch] {\it (1) implies (2):} By \cref{RFSum}, fix \(I\in U\) and
	a sequence of ultrafilters \(\langle W_i : i \in I\rangle\) such that \(W
	\cong U\text{-}\sum_{i\in I} W_i\). Let \(D = U\text{-}\sum_{i\in I} W_i\)
	and let \(Z = [\langle W_i : i\in I\rangle]_U\). We have
	\((M_Z^{M_U},j_Z^{M_U}\circ j_U) \cong (M_{D},j_D)\cong (M_W,j_W)\), so fix
	an isomorphism \(k : M_Z^{M_U}\to M_W\) such that \(k\circ j_Z^{M_U}\circ
	j_U = j_W\). It is easy to see that \(k\circ j_Z^{M_U}\) is close to
	\(M_U\).
	
	{\it (2) implies (1):} Let \(Y\) be the underlying set of \(W\) and let
	\(Z\) be the \(M_U\)-ultrafilter on \(j_U(Y)\) derived from \(h\) using
	\(\id_W\). Let \(k : M_{Z}^{M_U}\to M_W\) be the factor embedding. It is
	easy to see that \(k\) is surjective. Thus \(k\) is an isomorphism. It
	follows that \(U\text{-}\sum Z \cong W\), so by \cref{RFSum}, \(U\D W\).
	\end{proof}
\end{cor}

Note that the close embedding given by \cref{RFChar0} is ``isomorphic to" a
(possibly illfounded) internal ultrapower embedding of \(M_U\). But the language
of close embeddings makes it possible to work with the Rudin-Frol\'ik order in
fairly simple model theoretic terms while keeping our language precise.

In the countably complete case, \cref{RFChar0} really does imply that there is
an internal ultrapower embedding from \(M_U\) to \(M_W\):

\begin{cor}\label{RFChar}
	If \(U\) and \(W\) are countably complete ultrafilters, then \(U\D W\) if
	and only if there is an internal ultrapower embedding \(h : M_U\to M_W\)
	such that \(h\circ j_U = j_W\).\qed
\end{cor}

\begin{cor}\label{RFDirected}\index{Rudin-Frol\'ik order!directedness}
	The Ultrapower Axiom holds if and only if the Rudin-Frol\'ik order is
	directed on countably complete ultrafilters.
	\begin{proof}
		Assume the Ultrapower Axiom. Suppose \(U\) and \(W\) are countably
		complete ultrafilters. Let \(j :V\to M\) and \(i : V\to N\) be their
		respective ultrapower embeddings. Using UA, fix an internal ultrapower
		comparison \((k,h) : (M,N)\to P\). Then the composition \(k\circ j =
		h\circ i\) is an ultrapower embedding of \(V\), associated say to the
		countably complete ultrafilter \(D\). Then \(U\D D\) since \(k :M_U\to
		M_D\) is an internal ultrapower embedding such that \(k\circ j_U =
		k\circ j = j_D\). Similarly, \(W\D D\). Thus the Rudin-Frol\'ik order is
		directed on countably complete ultrafilters. The converse is similar.
	\end{proof}
\end{cor}

\cref{RFChar0} makes the relationship between the Rudin-Frol\'ik order and the
Rudin-Keisler order clear:
\begin{cor}\label{RKExtend}
The Rudin-Keisler order extends the Rudin-Frol\'ik order.\index{Rudin-Frol\'ik
order!vs. the Rudin-Keisler order}
\begin{proof}
Suppose \(U\D W\). Then by \cref{RFChar0}, there is an elementary embedding \(h
: M_U\to M_W\) such that \(h\circ j_U = j_W\). By \cref{RKChar}, \(U\RK W\).
\end{proof}
\end{cor}

Thus by \cref{RKThm}, if \(U\D W\) and \(W\D U\), then \(U\cong W\). This
motivates the following definition:
\begin{defn}
The {\it strict Rudin-Frol\'ik order} is defined on ultrafilters \(U\) and \(W\)
by setting \(U\sD W\) if \(U\D W\) but \(U\not \cong W\).
\end{defn}

\begin{lma}\label{RFWF}
The strict Rudin-Frol\'ik order is wellfounded on countably complete
ultrafilters.
\begin{proof}
This follows from the fact that the strict Rudin-Keisler order extends the
strict Rudin-Frol\'ik order (\cref{RKExtend}) and is wellfounded on countably
complete ultrafilters (\cref{RKWF}). 
\end{proof}
\end{lma}

The Rudin-Frol\'ik order is {\it not} directed on arbitrary ultrafilters. In
fact, the Rudin-Frol\'ik order restricted to ultrafilters on \(\omega\) already
fails to be directed. We sketch a proof of this fact that bears a striking
resemblance to many of the comparison arguments used throughout this
dissertation, especially \cref{MuMeasure} below. We hope it demonstrates that
the close embedding approach to the Rudin-Frol\'ik order really yields some
simplifications.
\begin{thm}[Rudin] If \(U\) and \(W\) are ultrafilters on \(\omega\) that have
an upper bound in the Rudin-Frol\'ik order, then either \(U\D W\) or \(W\D U\).
\begin{proof}[Sketch] By \cref{RFChar0} (3), the existence of an \(\D\)-upper
bound of \(U\) and \(W\) implies the existence of close embeddings \((k,h) :
(M_U,M_W)\to N\) such that \(k\circ j_U = h\circ j_W\). Assume without loss of
generality that \(k(\id_U) < h(\id_W)\). Let \(Z\) be the \(M_W\)-ultrafilter on
\(j_W(\omega)\) derived from \(h\) using \(k(\id_U)\). Then \(Z\) concentrates
on \(\id_W < j_W(\omega)\). Since \(Z\) belongs to \(M_W\) and concentrates on
an \(M_W\)-finite set, \(Z\) is principal. Since \(Z\) is derived from \(h\)
using \(k(\id_U)\), we must in fact have \(h(\id_Z) = k(\id_U)\). 

We now follow the argument of \cref{RFSEquiv}. Since \(k(\id_U)\in h[M_W]\), it
is easy to see that \(k[M_U] = H^N(k\circ j_U[V]\cup \{k(\id_U)\})\subseteq
h[M_W]\). Define \(e : M_U\to M_W\) by \(e = h^{-1}\circ k\). Then \(e\) is an
elementary embedding and \(h \circ e = k\), so since \(k\) is close to \(M_U\),
\(e\) is close to \(M_U\). Thus there is a close embedding \(e : M_U\to M_W\),
and it follows that \(U\D W\).
\end{proof}
\end{thm}
This theorem is often summarized by the statement that ``the Rudin-Frol\'ik
order forms a tree," but this is only true of the Rudin-Frol\'ik order on
\(\omega\). The reader should note that this proof is very similar to the proof
of the linearity of the Mitchell order from UA. The argument shows that natural
generalization of the seed order to \(\beta(\omega)\) is equal to the
Rudin-Frol\'ik order, while the natural generalization of the Ketonen order is
equal to the revised Rudin-Keisler order.
\begin{cor}\index{Rudin-Frol\'ik order!directedness}
The Rudin-Frol\'ik order on \(\beta(\omega)\) is not directed.
\begin{proof}
Assume towards a contradiction that the Rudin-Frol\'ik order on
\(\beta(\omega)\) is directed. Then it is linear. This contradicts the
well-known theorem of Kunen \cite{KunenRK} that the Rudin-Keisler order is not
linear on ultrafilters on \(\omega\).
\end{proof}
\end{cor}
Thus, unsurprisingly, the analog of the Ultrapower Axiom for countably
incomplete ultrafilters is false.

We conclude this section with a basic rigidity lemma for the Rudin-Frol\'ik
order on countably complete ultrafilters that apparently had not been noticed:
\begin{thm}\label{RFRigid}
Suppose \(U\) is a countably complete ultrafilter. Suppose \(I\in U\) and
\(\langle W_i : i\in I\rangle\) and \(\langle W'_i : i\in I\rangle\) are
discrete sequences of countably complete ultrafilters such that 
\[U\text{-}\lim_{i\in I} W_i = U\text{-}\lim_{i\in I} W'_i\] Then for
\(U\)-almost all \(i\in I\), \(W_i = W_i'\).
\end{thm}
In other words, there is at most one way to realize one countably complete
ultrafilter as a discrete limit with respect to another.
\begin{lma}\label{UniqueFactor}
Suppose \(U\) and \(W\) are countably complete ultrafilters. Then there is at
most one internal ultrapower embedding \(h : M_U\to M_W\) such that \(h\circ j_U
= j_W\).
\begin{proof}
Suppose \(h,k : M_U\to M_W\) are internal ultrapower embeddings such that
\(h\circ j_U = k\circ j_U\). In other words, \(h\restriction j_U[V] =
k\restriction j_U[V]\). Moreover \(h\restriction \text{Ord} = k\restriction
\text{Ord}\) by \cref{MinDefEmb}. Since \(M_U = H^{M_U}(j_U[V]\cup
\text{Ord})\), it follows that \(h = k\). 
\end{proof}
\end{lma}

\begin{proof}[Proof of \cref{RFRigid}] Let \(Z = [\langle W_i : i\in
I\rangle]_U\) and let \(Z' = [\langle W'_i : i\in I\rangle]_U\). By
\cref{DiscreteSum},  \(U\text{-}\sum_{i\in I} Z_i \cong U\text{-}\sum_{i\in I}
Z'_i\) and their projections to the second coordinate are equal. Using the
ultrapower theoretic characterization of sums (\cref{SumPower}), this means: 
\begin{align*}
	j_Z^{M_U}\circ j_U &= j_{Z'}^{M_U}\circ j_U\\
	\id_Z &= \id_{Z'}
\end{align*}\cref{UniqueFactor} now implies \(j_Z^{M_U} = j_{Z'}^{M_U}\). But \(Z\) and \(Z'\) are derived from \(j_Z^{M_U} = j_{Z'}^{M_U}\) using \(	\id_Z = \id_{Z'}\). Thus \(Z = Z'\). Finally, by \L o\'s's Theorem we have that \(W_i = W_i'\) for \(U\)-almost all \(i\in I\).
\end{proof}
The author's intuition is that \cref{RFRigid} should be true for countably
incomplete ultrafilters as well, and the fact that the proof does not just
generalize is a bit of a subtle point.
\section{Below the first \(\mu\)-measurable cardinal}\label{MuSection}
\subsection{Introduction}
In a sense, the first large cardinal axiom that is significantly beyond any
``normal ultrafilter axiom" is the existence of a \(\mu\)-measurable cardinal:
\begin{defn}\index{\(\mu\)-measurable cardinal}
A cardinal \(\kappa\) is said to be {\it \(\mu\)-measurable} if there is an
elementary embedding \(j : V\to M\) with critical point \(\kappa\) such that the
ultrafilter on \(\kappa\) derived from \(j\) using \(\kappa\) belongs to \(M\).
\end{defn}

The existence of a \(\mu\)-measurable cardinal is a large cardinal axiom that is
stronger than the existence of a measurable cardinal \(\kappa\) such that
\(o(\kappa) = 2^{2^\kappa}\), but weaker than the existence of a cardinal
\(\kappa\) that is \(2^\kappa\)-strong. 

As an example of the strength of \(\mu\)-measurable cardinals, let us show the
following fact:
\begin{prp}
Suppose \(\kappa\) is a \(\mu\)-measurable cardinal. Then there is a normal
ultrafilter on \(\kappa\) that concentrates on cardinals \(\delta\) such that
for any \(A\subseteq P(\delta)\), there is a normal ultrafilter \(D\) on
\(\delta\) such that \(A\in M_D\).
\begin{proof}
Let \(j : V\to M\) witness that \(\kappa\) is \(\mu\)-measurable and let \(U\)
be the normal ultrafilter on \(\kappa\) derived from \(j\) using \(\kappa\).
Thus \(U\in M\).
\begin{clm} \(M_U\) satisfies the statement that for all \(A\subseteq
P(\kappa)\), there is a normal ultrafilter \(D\) on \(\kappa\) such that \(A\in
M_D\).
\begin{proof}
Suppose not, and fix \(A\subseteq P(\kappa)\) such that \(M_U\) satisfies that
there is no normal ultrafilter \(D\) on \(\kappa\) with \(A\in (M_D)^{M_U}\).
Let \(k : M_U\to M\) be the factor embedding. By \cref{DerivedNF},
\(\textsc{crt}(k) > \kappa\) and \(P(\kappa)\cap M_U = P(\kappa) = P(\kappa)\cap
M\), so \(k(A) = A\). Therefore since \(k\) is elementary,  \(M\) satisfies that
there is no normal ultrafilter \(D\) on \(\kappa\) with \(A\in (M_D)^M\). But
\(A\in j_U(V_\kappa) \subseteq (M_U)^M\), and \(U\in M\) is a normal
ultrafilter. This is a contradiction.
\end{proof}
\end{clm}
By \L o\'s's Theorem, \(U\) concentrates on cardinals \(\delta\) such that for
all \(A\subseteq P(\delta)\), there is a normal ultrafilter \(D\) on \(\delta\)
such that \(A\in M_D\).
\end{proof}
\end{prp}
Thus a \(\mu\)-measurable cardinal is a limit of cardinals \(\delta\) such that
\(o(\delta) = 2^{2^\delta}\).

\subsection{Irreducible ultrafilters and \(\mu\)-measurability}\label{IrredMuSection}
The goal of the next few subsections is to analyze the countably complete
ultrafilters in \(V_\kappa\) where \(\kappa\) is the least \(\mu\)-measurable
cardinal. We first analyze simpler ultrafilters called {\it irreducible
ultrafilters} and then we reduce the general case to the irreducible case.
\begin{defn}\index{Irreducible ultrafilter}
An a nonprincipal countably complete ultrafilter \(U\) is {\it irreducible} if
every ultrafilter \(D\sD U\) is principal.
\end{defn}
Let us give some examples of irreducible ultrafilters.
\begin{prp}
If \(U\) is a normal ultrafilter on a cardinal \(\kappa\), then \(U\) is
irreducible.
\begin{proof}
Suppose \(D\sD U\). By \cref{RKExtend}, \(D\sRK U\), and therefore by
\cref{MinimalIncompressible}, \(\lambda_D < \kappa\). But since \(D\RK U\),
\(D\) is \(\kappa\)-complete. Since \(D\) is \(\kappa\)-complete and \(\lambda_D
< \kappa\), \(D\) is principal.
\end{proof}
\end{prp}
A direct generalization of this yields:
\begin{prp}
Normal fine ultrafilters are irreducible.
\begin{proof}
Suppose \(\mathcal U\) is a normal fine ultrafilter. By \cref{IsoNormalThm},
\(\mathcal U\) is isomorphic to an isonormal ultrafilter \(U\) on a cardinal
\(\lambda\). It suffices to show that \(U\) is irreducible. Suppose \(D\sD U\),
and we will show \(D\) is principal. By \cref{RKExtend}, \(D\sRK U\), and
therefore by \cref{MinimalIncompressible}, \(\lambda_D < \lambda\). Since \(D \D
U\), \(M_U\) is contained in \(M_D\), and so because \(j_U\) is
\(\lambda\)-supercompact, using \cref{UltrapowerSC}, \(\text{Ord}^\lambda
\subseteq M_U\subseteq M_D\). In particular, \(j_D\restriction \lambda\in M_D\),
so \(j_D\) is \(\lambda\)-supercompact. Since \(\lambda_D < \lambda\) and
\(j_D\) is \(\lambda\)-supercompact, \(D\) is principal by \cref{UFSuperBound}.
\end{proof}
\end{prp}
Dodd sound ultrafilters are also irreducible:
\begin{prp}
If \(U\) is a Dodd sound ultrafilter, then \(U\) is irreducible.
\begin{proof}
Suppose \(D\sD U\), and we will show \(D\) is principal. We may assume without
loss of generality that \(D\) is incompressible. Then since \(D\sRK U\), in fact
\(D\sE U\) by \cref{RKKet}. Since the Lipschitz order extends the Ketonen order,
\(D\sLi U\), so by \cref{LipMO}, \(D\mo U\). But then \(D\in M_U\subseteq M_D\),
so \(D\mo D\), which implies \(D\) is principal by \cref{MOStrict}.
\end{proof}
\end{prp}
Finally returning to \(\mu\)-measurable cardinals, we have the following fact:
\begin{prp}\label{muIrred}
Suppose \(j : V\to M\) is such that \(\textsc{crt}(j) = \kappa\) and \(U_0\in
M\) where \(U_0\) is the normal ultrafilter on \(\kappa\) derived from \(j\).
Let \(U_1\) be the ultrafilter on \(V_\kappa\) derived from \(j\) using \(U_0\).
Then \(U_1\) is irreducible and \(U_1\) is not isomorphic to a normal
ultrafilter.
\begin{proof}
Let \(j_1 : V\to M_1\) be the ultrapower of \(V\) by \(U_1\). The key point,
which is easily verified, is that \(\id_{U_1} = U_0\). Also note that \[M_{1} =
H^{M_{1}}(j_{1}[V]\cup (2^{2^\kappa})^{M_1})\] since \(\id_{U_1} = U_0\in
H^{M_{1}}(j_{1}[V]\cup (2^{2^\kappa})^{M_1})\), \(U_0\) being a subset of
\(P(\kappa)\).

We now show that \(U_1\) is irreducible. Suppose \(D\D U_1\) and \(D\) is
nonprincipal. We must show \(D \cong U_1\). Since \(\lambda_D = \kappa\), we
have \(\textsc{crt}(j_D) = \kappa\). Let \(k : M_D\to M_{1}\) be the unique
internal ultrapower embedding with \(k\circ j_D = j_{1}\).  We claim \(k(\kappa)
= \kappa\). Supposing the contrary, we have that \(k(\kappa) >\kappa\) is an
inaccessible cardinal that is a generator of \(j_1\), contradicting that \(M_1 =
H^{M_{1}}(j_{1}[V]\cup (2^{2^\kappa})^{M_1})\). Thus \(k(\kappa) = \kappa\).
Since \(M_1\subseteq M_D\), \(U_0\in M_D\), and since \(k(\kappa) = \kappa\),
\(k(U_0) = U_0\). Since \(U_0 = \id_{U_1}\), it follows that \(k\) is
surjective. Thus \(k\) is an isomorphism, and it follows that \(D\cong U_1\).

Finally we show that \(U_1\) is not isomorphic to a normal ultrafilter. Suppose
towards a contradiction that it is. Then in fact, \(U_1\) is isomorphic to the
ultrafilter on \(\kappa\) derived from \(j_{U_1}\) using \(\kappa\), namely
\(U_0\). In particular, \(M_{U_0} = M_{U_1}\), so since \(U_0\in M_{U_1}\), in
fact \(U_0\in M_{U_0}\). This contradicts the fact that the Mitchell order is
irreflexive (\cref{MOStrict}).
\end{proof}
\end{prp}
Under UA, \cref{muIrred} has a converse:
\begin{thm}[UA]\label{MuDichotomy}
Suppose \(\kappa\) is a measurable cardinal. Exactly one of the following holds:
\begin{enumerate}[(1)]
\item \(\kappa\) is \(\mu\)-measurable.
\item Every irreducible ultrafilter \(U\) of completeness \(\kappa\) is
isomorphic to a normal ultrafilter.
\end{enumerate}
\end{thm}

The proof will use some of the machinery from \cref{KetonenChapter}. Recall that
a {\it pointed ultrapower embedding} is a pair \((j,\alpha)\) such that \(j
:V\to M\) is an ultrapower embedding and \(\alpha\) is an ordinal. For the
reader's convenience we restate here the definition of the Ketonen order and the
seed order to pointed ultrapower embeddings:
\begin{defn}
Suppose \((j,\alpha)\) and \((i,\beta)\) are pointed ultrapower embeddings.
\begin{itemize}
\item \((j,\alpha)\E (i,\beta)\) (resp.  \((j,\alpha)\sE (i,\beta)\)) if there
is a \(1\)-internal comparison \((k,h)\) of \((j,i)\) such that \(k(\alpha) \leq
h(\beta)\) (resp. \(k(\alpha) < h(\beta)\)).
\item \((j,\alpha)\KE (i,\beta)\) if \((j,\alpha)\E (i,\beta)\) and
\((i,\beta)\E (j,\alpha)\).
\item \((j,\alpha)\wo (i,\beta)\) (resp.  \((j,\alpha)\swo (i,\beta)\)) if there
is an internal comparison \((k,h)\) of \((j,i)\) such that \(k(\alpha) \leq
h(\beta)\) (resp.  \(k(\alpha) < h(\beta)\)).
\item \((j,\alpha)=_S (i,\beta)\) if \((j,\alpha)\wo (i,\beta)\) and
\((i,\beta)\wo (j,\alpha)\).
\end{itemize}
\end{defn}
Equivalently \((j,\alpha)=_S (i,\beta)\) if there is an internal comparison
\((k,h)\) of \((j,i)\) such that \(k(\alpha) = h(\beta)\). Two fundamental
consequences of UA are that \(\E\) and \(\wo\) coincide on pointed ultrapower
embeddings (\cref{GenSE}) and are prewellorders (\cref{GenTrichotomy}).

The following lemma is an immediate consequence of \cref{RFSEquiv}:

\begin{lma}\label{SeedEq}
Suppose \(U\) and \(W\) are countably complete ultrafilters concentrating on
ordinals. Then \(U\D W\) if and only if for some ordinal \(\alpha\),
\((j_U,\id_U) =_S (j_W,\alpha)\).
\end{lma}

The following theorem can be viewed as yet another generalization of the proof
that the Mitchell order is linear under UA.

\begin{thm}[UA]\label{MuMeasure}
Suppose \(U\) is a countably complete ultrafilter. Let \(D\) be the normal
ultrafilter on \(\kappa = \textsc{crt}(j_U)\) derived from \(j_U\) using
\(\kappa\). Then either \(D\D U\) or \(D\mo U\).
\begin{proof}
Let \(i : M_D\to M_U\) be the factor embedding. Then \((i,\text{id}) :
(M_D,M_U)\to M_U\) witnesses that \((j_D,\kappa) \E (j_U,\kappa)\). Thus
\((j_D,\kappa)\wo (j_U,\kappa)\), so let \((k,h) : (M_D,M_U)\to N\) be an
internal ultrapower comparison of \((j_D,j_U)\) witnessing this. In other words,
\(k(\kappa) \leq h(\kappa)\). The proof now breaks into two cases:

\begin{case}\label{RFCase}
\(k(\kappa) = h(\kappa)\)
\begin{proof}[Proof in \cref{RFCase}] Then \((k,h)\) witnesses \((j_D,\kappa)=_S
(j_U,\kappa)\). \cref{SeedEq} therefore implies that \(D\D U\).
\end{proof}
\end{case}

\begin{case}\label{MOCase}
\(k(\kappa) < h(\kappa)\)
\begin{proof}[Proof in \cref{MOCase}] We will show that \(D\in M_U\). The key
point is that for any \(A\subseteq \kappa\), \[h(j_U(A))\cap h(\kappa) =
h(A)\cap h(\kappa)\] and therefore
\begin{align*}A\in D&\iff \id_D\in j_D(A)\\
&\iff k(\id_D)\in k(j_D(A))\\
&\iff k(\id_D)\in h(j_U(A))\\
&\iff k(\id_D)\in h(A)
\end{align*}
Since \(h\) is definable over \(M_U\) and \(P(\kappa)\subseteq M_U\), it follows
that \(D\) is a definable over \(M_U\), and hence \(D\in M_U\).
\end{proof}
\end{case}
Thus in \cref{RFCase}, \(D\D U\), and in \cref{MOCase}, \(D\mo U\). This proves
the theorem.
\end{proof}
\end{thm}

A more abstract perspective on this argument is that it generalizes the
linearity of the Mitchell order on Dodd sound ultrafilters (\cref{DoddMO}):

\begin{prp}
	Suppose \(\alpha\) is an ordinal and \(i : V\to N\) is an \(\alpha\)-sound
	elementary embedding. Suppose \(U\) is a countably complete tail uniform
	ultrafilter on an ordinal \(\eta\) such that \[(j_U,\id_U)\sE (i,\alpha)\]
	Then \(U\in N\).
	\begin{proof}[Sketch] Recall that \(i^\alpha : P(\delta) \to N\) is the
		function \(i^\alpha(A) = i(A)\cap \alpha\) where \(\delta\) is least
		such that \(i(\delta)\geq \alpha\). The \(\alpha\)-soundness of \(i\)
		amounts to the fact that \(i^\alpha\in N\).
		
		Fix a \(1\)-internal comparison \((k,h) : (M_U,N)\to P\). We have \(U =
		i^{-1}[h^{-1}[\pr {k(\id_U)} {h(i(\eta))}]]\) by the usual argument, so
		since \(k(\id_U) < h(\alpha)\), 
		\[U = (i^\alpha)^{-1}[h^{-1}[\pr {k(\id_U)} {h(\alpha)}]]\cap P(\eta)\]
		Since \(i^\alpha\in N\) and \(h\) is definable over \(N\), it follows
		that \(U\in N\).
	\end{proof}
\end{prp}

\cref{MuMeasure} leads to the proof of \cref{MuDichotomy}.

\begin{proof}[Proof of \cref{MuDichotomy}] Assume (1) fails, and we will show
(2) holds. Let \(U\) be an irreducible ultrafilter such that \(\textsc{crt}(j_U)
= \kappa\) but \(U\) is not isomorphic to a normal ultrafilter. Let \(D\) be the
normal ultrafilter on \(\kappa\) derived from \(j_U\) using \(\kappa\). By
\cref{MuMeasure}, either \(D\D U\) or \(D\mo U\). If \(D\D U\) then since \(D\)
is nonprincipal and \(U\) is irreducible, \(D\cong U\), contrary to our
hypothesis that \(U\) is not isomorphic to a normal ultrafilter. Therefore
\(D\mo U\). Then \(j_U : V\to M_U\) has critical point \(\kappa\) and the normal
ultrafilter on \(\kappa\) derived from \(j_U\) using \(\kappa\) belongs to
\(M_U\), so \(\kappa\) is a \(\mu\)-measurable cardinal. Therefore (2) holds.

If (2) holds, then (1) fails as a consequence of \cref{muIrred}.
\end{proof}

\begin{cor}[UA] Suppose \(\kappa\) is the least \(\mu\)-measurable cardinal.
Then every irreducible ultrafilter in \(V_\kappa\) is isomorphic to a normal
ultrafilter.
\begin{proof}
This follows from \cref{MuDichotomy} applied in \(V_\kappa\), which is a model
of ZFC + UA that also satisfies the statement that there are no
\(\mu\)-measurable cardinals.
\end{proof}
\end{cor}

\begin{cor}[UA]\label{LeastIrreds}
Let \(\kappa\) be the least measurable cardinal. Then \(\kappa\) carries a
unique irreducible ultrafilter up to isomorphism.\qed
\end{cor}
\subsection{Factorization into irreducibles}\label{FactorizationSection}
The results of the previous section motivate understanding how arbitrary
countably complete ultrafilters relate to irreducible ultrafilters. The main
theorem of this subsection answers the question in complete generality: every
ultrapower embedding can be written as a finite iterated ultrapower of
irreducible ultrafilters. To be perfectly precise, let us introduce some
notation for iterated ultrapowers.
\begin{defn}\index{Ultrapower!iterated ultrapower}\index{Iterated ultrapower}
Suppose \(\nu\) is an ordinal. An {\it iterated ultrapower of length \(\nu\)} is
a sequence \(\langle M_\beta,U_\alpha, j_{\alpha,\beta} : \alpha < \beta <
\nu\rangle\) such that the following hold:
\begin{itemize}
\item For all \(\alpha\) with \(\alpha + 1 < \nu\), \(U_\alpha\) is a countably
complete ultrafilter of \(M_\alpha\) and \(j_{\alpha,\alpha+1} : M_\alpha\to
M_{\alpha+1}\) is the ultrapower of \(M_\alpha\) by \(U_\alpha\).
\item For all \(\alpha < \beta < \gamma < \nu\), \(j_{\alpha,\gamma} = j_{\beta,
\gamma}\circ j_{\alpha,\beta}\).
\item For all limit ordinals \(\gamma < \nu\), \(M_\gamma\) is the direct limit
of \(\langle M_\alpha, j_{\alpha,\beta} : \alpha < \beta < \gamma\rangle\) and
for all \(\alpha < \gamma\), \(j_{\alpha,\gamma} : M_\alpha\to M_\gamma\) is the
direct limit embedding.
\end{itemize}
\end{defn}
Note that the iterated ultrapower \(\langle M_\beta,U_\alpha, j_{\alpha,\beta} :
\alpha < \beta < \nu\rangle\) is actually completely determined by the sequence
\(\langle U_\alpha : \alpha + 1 < \nu\rangle\). We make the convention that for
\(\beta < \nu\), \(j_{\beta,\beta}\) is the identity function on \(M_\beta\).

\begin{thm}[UA]\label{FactorizationThm}
Suppose \(W\) is a countably complete ultrafilter. Then there is a finite linear
iterated ultrapower \(\langle M_n,U_m,j_{m,n} : m < n \leq \ell\rangle\) such
that \(M_0 = V\), \(M_\ell = M_W\), \(U_m\) is an irreducible ultrafilter of
\(M_m\) for all \(m <\ell\), and \(j_W = j_{0,\ell}\).\index{Irreducible
ultrafilter!Factorization Theorem}
\end{thm}

The proof of this theorem relies on a stronger structural property of the
Rudin-Frol\'ik order:

\begin{thm}[UA]\label{ACC}\index{Rudin-Frol\'ik order!local ascending chain condition}
Suppose \(W\) is a countably complete ultrafilter. Then there is no ascending
chain \(D_0\sD D_1\sD D_2\sD\cdots\) such that \(D_n\D W\) for all \(n <
\omega\).
\end{thm}

More succinctly, the Rudin-Frol\'ik order satisfies the {\it local ascending
chain condition}.  Later we will give a deeper explanation of why this is true
(\cref{RFFinite}): a countably complete ultrafilter has only finitely many
Rudin-Frol\'ik predecessors up to isomorphism. 

We defer the proof of \cref{ACC} to the next section. In this section we will
derive \cref{FactorizationThm} from \cref{ACC}, and show how this can be used to
analyze ultrafilters on the least measurable cardinal. 

Before we can proceed, we need a simple lemma about the pervasiveness of
irreducible ultrafilters:
\begin{lma}\label{Atomic}
Suppose \(D\sD W\) are countably complete ultrafilters. Then there is a
countably complete ultrafilter \(F\) with \(D\sD F\D W\) and an irreducible
ultrafilter \(U\) of \(M_D\) such that \(j_{F} = j_U^{M_D}\circ j_D\).
\begin{proof}
By the wellfoundedness of the Rudin-Frol\'ik order on countably complete
ultrafilters (\cref{RFWF}), let \(F\) be \(\sD\)-minimal among ultrafilters
\(Z\) such that \(D\sD Z \D W\). By \cref{RFChar}, fix a countably complete
ultrafilter \(U\) of \(M_Z\) such that \(j_{F} = j_U^{M_Z}\circ j_D\). 

We claim \(U\) is an irreducible ultrafilter of \(M_D\). Suppose \(\bar U\sD U\)
in \(M_D\), and we will show that \(\bar U\) is principal in \(M_D\). Let \(\bar
F\) be a countably complete ultrafilter such that \(j_{\bar F} = j_{\bar
U}^{M_D}\circ j_D\). One easily computes: \[D\D \bar F\sD F\D W\] Assume towards
a contradiction \(D\sD \bar F\); then \(D\sD \bar F \D W\) and \(\bar F\sD F\),
contradicting that \(F\) is \(\sD\)-minimal among ultrafilters \(Z\) such that
\(D\sD Z \D W\). Therefore \(D\not \sD \bar F\), or in other words \(D\cong \bar
F\). Now \[M_D = M_{\bar F} = M_{\bar U}^{M_D}\] It follows that \(\bar U\) is
principal in \(M_D\).
\end{proof}
\end{lma}

We now deduce \cref{FactorizationThm} from \cref{ACC}.
\begin{proof}[Proof of \cref{FactorizationThm} assuming \cref{ACC}] By
recursion, we construct a finite sequence of countably complete ultrafilters
\(D_0\sD D_1\sD\cdots \sD D_\ell \cong W\) and an iterated ultrapower \(\langle
M_n,U_m,j_{m,n} : m < n  \leq \ell\rangle\) such that \(M_0 = V\), \(U_m\) is an
irreducible ultrafilter of \(M_m\) for all \(m < \ell\), and \(j_{0,n} =
j_{D_n}\) for all \(n \leq \ell\).

To begin, let \(M_0 = V\) and let \(D_0\) be principal.

Suppose \(D_0\D D_1\D\cdots \D D_k \D W\) and \(\langle M_n,U_m,j_{m,n} : m < n
\leq k\rangle\) have been constructed. If \(D_k\cong W\), we set \(\ell = k\)
and terminate the construction. Otherwise, \(D_k\sD W\). Using \cref{Atomic},
fix \(D_{k+1}\) with \(D_k\sD D_{k+1}\D W\) and an irreducible ultrafilter
\(U_k\) of \(M_{D_k} = M_k\) such that \(j_{D_{k+1}} = j_{U_k}^{M_k}\circ
j_{D_k}\). Let \(\langle M_n,U_m,j_{m,n} : m < n \leq k + 1\rangle\) be the
iterated ultrapower given by the sequence \(\langle U_m : m \leq k\rangle\).

This recursion must terminate in finitely many steps, since otherwise we will
produce \(D_0\sD D_1\sD \cdots\) with \(D_n\D W\) for all \(n < \omega\),
contradicting the local ascending chain condition (\cref{ACC}). When the process
terminates, we have \(D_\ell \cong W\). This yields the objects promised in the
first paragraph.

In particular, we have produced an iterated ultrapower \(\langle M_n,U_m,j_{m,n}
: m < n  \leq \ell\rangle\) such that \(U_m\) is an irreducible ultrafilter of
\(M_m\) for all \(m < \ell\) and \(j_{0,\ell} = j_{D_\ell} = j_W\), as desired.
\end{proof}

We now turn our sights back to the countably complete ultrafilters below the
least \(\mu\)-measurable cardinal.

\begin{thm}[UA]\label{NormalIteration}
Assume that there are no \(\mu\)-measurable cardinals. Suppose \(W\) is a
countably complete ultrafilter. Then there is a finite iterated ultrapower
\(\langle M_n,U_m,j_{m,n} : m < n \leq \ell\rangle\) such that \(M_0 = V\),
\(M_\ell = M_W\), \(U_m\) is a normal ultrafilter of \(M_m\) for all \(m
<\ell\), and \(j_W = j_{0,\ell}\).
\begin{proof}
This is immediate from \cref{MuDichotomy} and \cref{FactorizationThm}.
\end{proof}
\end{thm}

Stated more succinctly, if there are no \(\mu\)-measurable cardinals and the
Ultrapower Axiom holds, then every ultrapower embedding is given by a finite
iteration of normal ultrafilters. Combining this with the linearity of the
Mitchell order on normal ultrafilters, \cref{NormalIteration} comes very close
to a complete analysis of all countably complete ultrafilters below the least
\(\mu\)-measurable cardinal on the assumption of the Ultrapower Axiom alone. In
any case, it gives as complete an analysis as the Ultrapower Axiom ever will:

\begin{prp}
The following are equivalent:
\begin{enumerate}[(1)]
\item The Mitchell order is linear and every ultrapower embedding is given by a
finite iteration of normal ultrafilters.
\item The Ultrapower Axiom holds and there are no \(\mu\)-measurable
cardinals.\qed
\end{enumerate}
\end{prp}
The proof is as obvious as it is tedious, and it is omitted.

We now derive the analog of Kunen's theorem (\cref{KunenUlts} above):

\begin{thm}[UA]\label{KunenUA}
Suppose \(\kappa\) is the least measurable cardinal. Let \(U\) be the unique
normal ultrafilter on \(\kappa\). Then every countably complete ultrafilter on
\(\kappa\) is isomorphic to \(U^n\) for some \(n < \omega\).
\begin{proof}
We first prove the theorem assuming \(\kappa\) is the only measurable cardinal.
Then \(U\) is the only normal ultrafilter. Thus by \cref{FactorizationThm},
every ultrapower embedding is given by a finite iterated ultrapower of \(U\). In
other words, every countably complete ultrafilter is isomorphic to \(U^n\) for
some \(n < \omega\).

We now prove the theorem assuming there are two measurable cardinals. Let
\(\delta\) be the second one. Since \(V_\delta\) is a model of UA and satisfies
that \(\kappa\) is the only measurable cardinal, by the previous paragraph
\(V_\delta\) satisfies that every countably complete ultrafilter is isomorphic
to \(U^n\) for some \(n < \omega\). Since every countably complete ultrafilter
on \(\kappa\) belongs to \(V_\delta\), it follows that (in \(V\)) every
countably complete ultrafilter on \(\kappa\) is isomorphic to \(U^n\) for some
\(n < \omega\).
\end{proof}
\end{thm}
We sketch how this implies the transfinite version of Kunen's theorem.
\begin{defn}
Suppose \(U\) is a countably complete ultrafilter and \(\nu\) is an ordinal.
Then \(j_{U^\nu} : V\to M_{U^\nu}\) denotes the elementary embedding \(j_{0,\nu}
: V \to M_\nu\) where \(\langle M_\beta,U_\alpha,j_{\alpha,\beta} : \alpha <
\beta \leq\nu\rangle\) is the iterated ultrapower given by setting \(U_\alpha =
j_{0,\alpha}(U)\) for all \(\alpha < \nu\).
\end{defn}

\begin{thm}[UA]\label{TransfiniteKunen}
Let \(\kappa\) be the least measurable cardinal. Let \(U\) be the unique normal
ultrafilter on \(\kappa\). Suppose \(M\) is an inner model and \(j : V\to M\) is
an elementary embedding such that \(M = H^M(j[V]\cup j(\kappa))\). Then \(j =
j_{U^\nu}\) for some ordinal \(\nu\).
\end{thm}

\begin{lma}\label{GeneratorDetermine}
Suppose \(M\) is an inner model, \(j : V\to M\) is an elementary embedding, and
\(\langle \xi_\alpha : \alpha < \nu\rangle\) is the increasing enumeration of
the generators of \(j\). For any \(p\in [\nu]^{<\omega}\), let \(U_p\) be the
ultrafilter on \([\mu_{j}(p)]^{|p|}\) derived from \(j\) using \(\{\xi_\alpha :
\alpha\in p\}\). Then \(j\) is uniquely determined by the sequence \(\langle U_p
: p\in [\nu]^{<\omega}\rangle\).
\begin{proof}[Sketch] This follows from the usual extender ultrapower
construction. This proof is not intended as an exposition of this construction;
we are merely checking, for the sake of the reader already familiar with this
construction, that a slightly modified version (i.e., using only generators)
works just as well. 

For \(p\in [\nu]^{<\omega}\), let \(j_p : V\to M_p\) be the ultrapower of the
universe by \(U_p\) and let \(k_p : M_p\to M\) be the unique elementary
embedding such that \(k_p\circ j_p = j\) and \(k_p(\id_{U_p}) = \{\xi_\alpha :
\alpha\in p\}\).

For \(p\subseteq q\in [\nu]^{<\omega}\), define \(k_{p,q} : M_p\to M_q\) by
setting \[k_{p,q}([f]_{U_p}) = [f']_{U_q}\] where, letting \(e : |p|\to |q|\) be
the unique function such that \(q_{e(n)} = p_n\), \(f'\) is defined for
\(U_q\)-almost every \(s\) by \[f'(s) = f(\{s_{e(n)} : n < |p|\})\] Then 
\[\langle M_{U_q},k_{p,q} : p\subseteq q\in [\nu]^{<\omega}\rangle\] is a
directed system. Let \(N\) be its direct limit and let \(j_{p,\infty} : M_p\to
N\) be the direct limit map. 

For any \(p\subseteq q\in [\nu]^{<\omega}\), it is easy to check that \(k_q\circ
k_{p,q} = k_p\). Therefore by the universal property of the direct limit, there
is a map \(k : N\to M\) such that \(k\circ j_{p,\infty}\) is equal to \(k_p :
M_p\to M\). 

We claim \(k\) is the identity. Towards a contradiction, suppose not. Let \(\xi
= \textsc{crt}(k)\). Then \(\xi\) is a generator of \(j\), so \(\xi =
\xi_\alpha\) for some \(\alpha < \nu\). But then letting \(a=
\id_{U_{\{\alpha\}}}\), we have \(\{\xi\} = k_{\{\alpha\}}(a) = k\circ
j_{p,\infty}(a) \in \text{ran}(k)\), so \(\xi\in \text{ran}(k)\), contradicting
that \(\xi = \textsc{crt}(k)\).

Since \(k\) is the identity, \(j_{0,\infty} = j\). Since the directed system
\(\langle M_{U_q},k_{p,q} : p\subseteq q\in [\nu]^{<\omega}\rangle\), and thus
the embedding \(j_{0,\infty}\), were defined only with reference to the sequence
\(\langle U_p : p\in [\nu]^{<\omega}\rangle\), the lemma follows.
\end{proof}
\end{lma}

\begin{lma}\label{TransfiniteReduction}
Suppose \(U\) is a normal ultrafilter, \(M\) is an inner model, and \(j : V\to
M\) is an elementary embedding such that for any \(a\in M\), the ultrafilter
derived from \(j\) using \(a\) is isomorphic to \(U^n\) for some \(n < \omega\).
Then \(M = M_{U^\nu}\) and \(j = j_{U^\nu}\) for some ordinal \(\nu\).
\begin{proof}[Sketch] For all \(m < \omega\), let \(\kappa_m =
j_{U^m}(\kappa)\), so \(\kappa_m\) is the \(m\)-th generator of \(j_{U^n}\) for
any \(n > m\). Let \(W_n\) be the ultrafilter on \([\kappa]^n\) derived from
\(j_{U^n}\) using \(\{\kappa_{n-1},\dots,\kappa_0\}\). Thus \(W_n\) is the
unique ultrafilter with the following properties:
\begin{itemize}
\item \(W_n\cong U^m\) for some \(m < \omega\).
\item The underlying set of \(W_n\) is \([\kappa]^n\).
\item Every element of \(\id_{W_m}\) is a generator of \(j_{W_m}\).
\end{itemize}

Since every ultrafilter \(Z\) derived from \(j\) is isomorphic to an ultrafilter
on \(\kappa\), the class of generators of \(j\) is contained in \(j(\kappa)\),
and in particular it forms a set. Let \(\langle \xi_\alpha : \alpha <
\nu\rangle\) enumerate this set in increasing order.  For any finite set
\(p\subseteq \nu\), the ultrafilter on \([\kappa]^n\) derived from \(j\) using
\(\{\xi_\alpha : \alpha\in p\}\) has the properties enumerated above, and hence
is equal to \(W_n\).

Let \(\langle \xi'_\alpha : \alpha < \nu\rangle\) denote the sequence of
generators of \(j_{U^\nu}\). Then for any finite set \(p\subseteq \nu\), the
ultrafilter on \([\kappa]^n\) derived from \(j_{U^\nu}\) using \(\{\xi'_\alpha :
\alpha\in p\}\) is equal to \(W_n\).

By \cref{GeneratorDetermine}, it follows that \(j = j_{U^\nu}\).
\end{proof}
\end{lma}

\begin{proof}[Proof of \cref{TransfiniteKunen}] The assumption that \(M =
H^M(j[V]\cup j(\kappa))\) implies that every ultrafilter derived from \(j\) is
isomorphic to an ultrafilter on \(\kappa\). By \cref{KunenUlts}, it follows that
every ultrafilter derived from \(j\) is isomorphic to \(U^n\) for some \(n\). By
\cref{TransfiniteReduction}, \(j = j_{U^\nu}\) for some ordinal \(\nu\).
\end{proof}

\section{The structure of the Rudin-Frol\'ik order}\label{RFStructureSection}
\subsection{The local ascending chain condition}\label{ACCSection}
The goal of this subsection is to prove \cref{ACC}, the local ascending chain
condition for the Rudin-Frol\'ik order. This uses two lemmas, the first of which
is often useful in the context of UA. The approach taken here uses the following
concept:
\begin{defn}\label{RFTrans}\index{Translation of an ultrafilter (\(\tr U W\))!when \(U\D W\)}\index{\(\tr U W\) (translation)}
Suppose \(Y\) is a set, \(W \in \mathscr B(Y)\), and \(U\D W\). Then the {\it
translation of \(U\) by \(W\),} denoted \(\tr U W\), is the unique
\(M_U\)-ultrafilter \(Z\in j_U(\mathscr B(Y))\) such that \(j_Z^{M_U}\circ j_U =
j_W\) and \(\id_{Z} = \id_W\).
\end{defn}

The uniqueness of \(Z\) follows from the fact (\cref{UniqueFactor}) that there
is at most one internal ultrapower embedding \(k : M_U\to M_W\) such that
\(k\circ j_U = j_W\). Then \(\tr U W\) must be the \(M_U\)-ultrafilter on
\(j_U(Y)\) derived from \(k\) using \(\id_W\). We view \(\tr U W\) as a version
of \(W\) inside \(M_U\). 

A more elegant, less comprehensible characterization of \(\tr U W\) is immediate
from the proof of \cref{RFRigid}:

\begin{lma}\label{TCombo}
Suppose \(U\) and \(W\) are countably complete ultrafilters.  Suppose \(I\) is a
set in \(U\) and \(\langle W_i : i\in I\rangle\) is a discrete sequence of
ultrafilters such that \(W = U\text{-}\lim_{i\in I} W_i\). Then \(\tr U W =
[\langle W_i : i\in I\rangle]_U\).\qed
\end{lma}

The following lemma links translations to the minimal covers from the proof of
\cref{LinKetThm}:

\begin{lma}\label{TMinimal}
	Suppose \(\delta\) is an ordinal, \(W\in \mathscr B(\delta)\), and \(U\D W\)
	is a countably complete ultrafilter. Then \(\tr U W\) is
	\(\sE^{M_U}\)-minimal among all \(Z\in j_U(\mathscr B(\delta))\) such that
	\(j_U[W]\subseteq Z\).
	\begin{proof}
		Fix \(Z\in j_U(\mathscr B(\delta))\) with \(j_U[W]\subseteq Z\). For
		ease of notation, let \(N = M_Z^{M_U}\). Then by \cref{LimitFactor},
		there is a unique embedding \(e : M_W \to N\) such that \(e\circ j_W =
		j_Z^{M_U} \circ j_U\) and \(e(\id_W) = \id_Z\). Thus the  \(1\)-internal
		comparison \((e,\text{id}) : (M_W, N)\to N\) witnesses \[(M_W,\id_W) \E
		(N, \id_Z)\] Suppose now towards a contradiction that \(Z \sE \tr U W\)
		in \(M_U\). Let \((k,h)\) be a \(1\)-internal comparison of
		\((j_Z^{M_U},j_{\tr U W}^{M_U})\) such that \(k(\id_Z) < h(\id_{\tr U
		W})\). Since \(M_Z^{M_U} = N\), \(M_{\tr U W}^{M_U} = M_W\), and
		\(\id_{\tr U W} = \id_W\), \((k,h) : (N,M_W)\to P\) is a \(1\)-internal
		comparison witnessing \[(N,\id_Z) \sE (M_W,\id_W)\] This contradicts the
		wellfoundedness of the Ketonen order on pointed ultrapowers
		(\cref{GenWellfounded}).
	\end{proof}
\end{lma}

\begin{lma}\label{TnotJ}
	Suppose \(U\D W\) are countably complete ultrafilters. If \(U\) is
	nonprincipal, then \(\tr U W \neq j_U(W)\).
	\begin{proof}
		Assume \(\tr U W = j_U(W)\), and we will show that \(U\) is principal.
		By \cref{TCombo}, fix a set \(I\in U\) and a discrete sequence \(\langle
		W_i : i\in I\rangle\) such that \([\langle W_i : i\in I\rangle]_U = \tr
		U W\). Since \(\langle W_i : i\in I\rangle\) is discrete, in particular
		the \(W_i\) are pairwise distinct. Since \(\tr U W = j_U(W)\), \L o\'s's
		Theorem implies that there is a \(U\)-large set \(J\subseteq I\) such
		that \(W_i = W\) for all \(i\in J\). Since the \(W_i\) are pairwise
		distinct, it follows that \(|J| = 1\). Thus \(U\) contains a set of size
		\(1\), so \(U\) is principal.
	\end{proof}
\end{lma}

\begin{prp}[UA]\label{Pushdown}
Suppose \(\delta\) is an ordinal, \(W\in \mathscr B(\delta)\), and \(U\D W\) is
a nonprincipal ultrafilter. Then \( \tr U W\sE j_U(W)\) in \(M_U\).
\begin{proof}
	By \cref{TMinimal} and the linearity of the Ketonen order, \(\tr U W\E
	j_U(W)\). By \cref{TnotJ}, \(\tr U W\neq j_U(W)\). It follows that \(\tr U
	W\sE j_U(W)\).
\end{proof}
\end{prp}

The following simple lemma on the preservation of the Rudin-Frol\'ik order under
translation functions will be used in the proof of \cref{ACC}:
\begin{lma}\label{TranslationAssociative}
Suppose \(U\), \(W\), and \(Z\) are countably complete ultrafilters with \(U\D
W,Z\). 
\begin{itemize}
	\item \(W \D Z\) if and only if \(\tr U W\D \tr U Z\) in \(M_U\).
	\item \(W \sD Z\) if and only if \(\tr U W\sD \tr U Z\) in \(M_U\).\qed
\end{itemize}
\end{lma}

We finally prove the local ascending chain condition.
\begin{proof}[Proof of \cref{ACC}]\index{Rudin-Frol\'ik order!local ascending chain condition}
Assume towards a contradiction that the theorem is false. Let \(C\) be the class
of countably complete tail uniform ultrafilters \(Z\) such that there is an
infinite \(\sD\)-ascending sequence \(\langle U_n : n < \omega\rangle\) sequence
\(\D\)-bounded above by \(Z\). Let \(W\) be a \(\sE\)-minimal element of \(C\),
and fix \(U_0\sD U_1\sD\cdots\) such that \(U_n \D W\) for all \(n < \omega\).
We may assume without loss of generality that \(U_0\) is nonprincipal. By
elementarity, \(j_{U_0}(W)\) is a \(\sE^{M_{U_0}}\)-minimal element of
\(j_{U_0}(C)\). 

Since translation functions preserve the Rudin-Frol\'ik order
(\cref{TranslationAssociative}), \(M_{U_0}\) satisfies \(\tr {U_0} {U_0}\sD \tr
{U_0} {U_1}\sD \tr {U_0} {U_2}\sD\cdots\) and \(\tr {U_0} {U_n}\D \tr {U_0} W\)
for all \(n \leq \omega\). Since \(M_{U_0}\) is closed under countable
sequences, it follows that \(\tr {U_0} W\in j_{U_0}(C)\). But by
\cref{Pushdown}, \(\tr {U_0} W \sE j_{U_0}(W)\). This contradicts that
\(j_{U_0}(W)\) is a \(\sE^{M_{U_0}}\)-minimal element of \(j_{U_0}(C)\).
\end{proof}

\subsection{Pushouts and the Rudin-Frol\'ik lattice}\label{PushoutSection}
In this section we establish that the Rudin-Frol\'ik order is a lattice:
\begin{thm*}[UA]\index{Rudin-Frol\'ik order!as a lattice}
The Rudin-Frol\'ik order is a lattice in the following sense:
\begin{itemize}
\item If \(U_0\) and \(U_1\) are countably complete ultrafilters, there is an
\(\D\)-minimum countably complete ultrafilter \(W\gD U_0,U_1\).
\item If \(U_0\) and \(U_1\) are countably complete ultrafilters, there is an
\(\D\)-maximum countably complete ultrafilter \(D\D U_0,U_1\).
\end{itemize}
\end{thm*}
The Rudin-Frol\'ik order is not literally a lattice because it is only a
preorder, but the theorem above shows that it induces a lattice structure on
isomorphism types of countably complete ultrafilters. The two parts will be
proved as \cref{Join} and \cref{Meet} below.

We begin by establishing the existence of least upper bounds in the
Rudin-Frol\'ik order, which is by far the most important part of the theorem.
Here it is cleaner to work with the elementary embeddings rather than the
ultrafilters:
\begin{defn}\index{Pushout}\index{Comparison!pushout}
Suppose \(j_0 : V\to M_0\) and \(j_1 : V\to M_1\) are ultrapower embeddings. An
internal ultrapower comparison \((i_0,i_1) : (M_0,M_1)\to N\) is a {\it pushout}
of \((j_0,j_1)\) if for any internal ultrapower comparison \((k_0,k_1) :
(M_0,M_1)\to P\), there is a unique internal ultrapower embedding \(h : N\to P\)
such that \(h \circ i_0 = k_0\) and \(h\circ i_1 = k_1\).
\end{defn}
\begin{figure}
	\center
	\includegraphics[scale=.65]{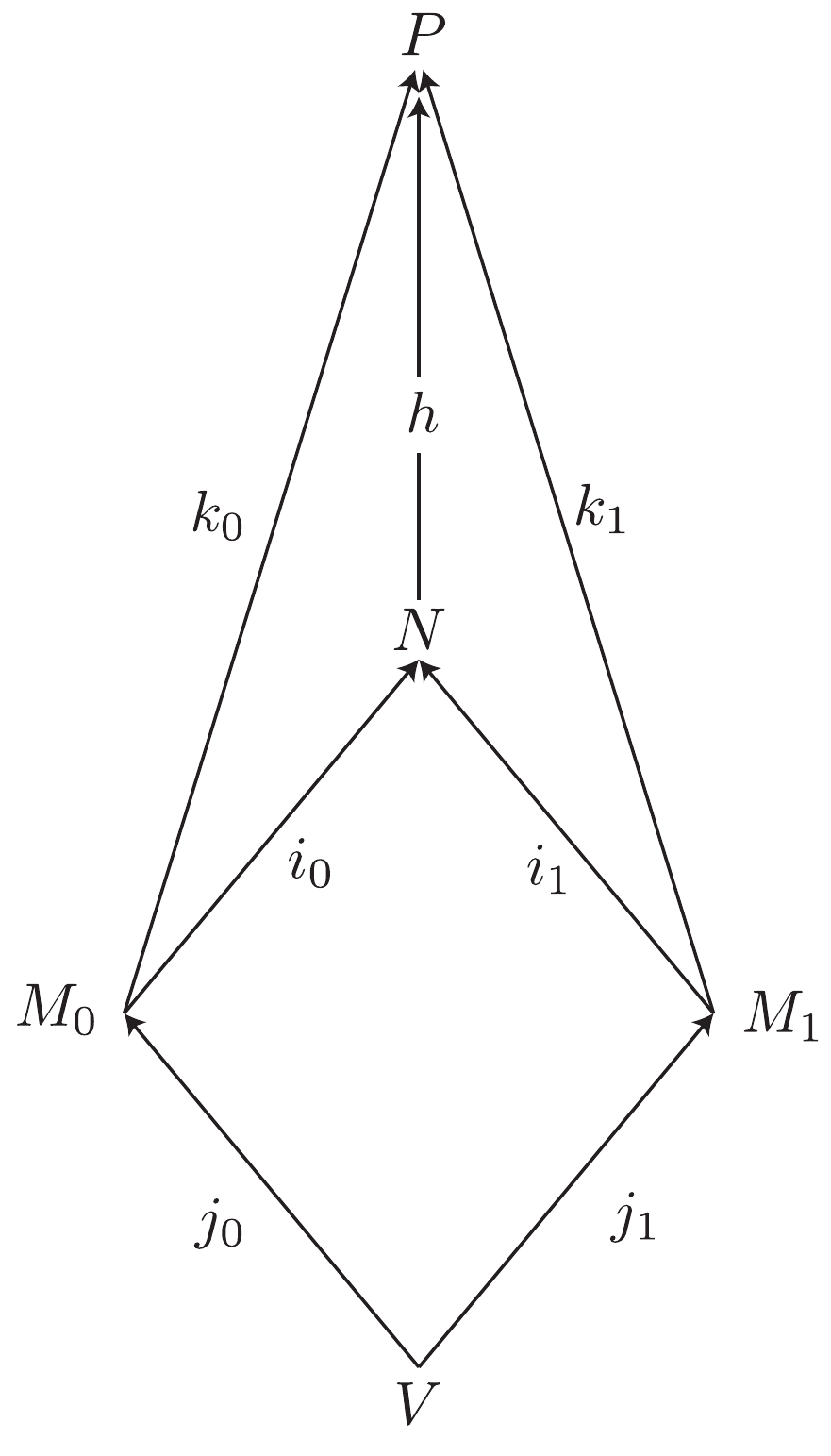}
	\caption{The pushout of \((j_0,j_1)\).}\label{PushoutFig}
\end{figure}
Pushout comparisons are simply the model theoretic manifestation of least upper
bounds in the Rudin-Frol\'ik order. We will prove:
\begin{thm}[UA]\label{Pushout}
Every pair of ultrapower embeddings has a unique pushout.
\end{thm}

The uniqueness of pushouts is a standard category theoretic fact: the pushout of
a pair of embeddings is what a category theorist would call the pushout of these
morphisms in the category \(\mathcal D_\infty\) of all ultrapower embeddings. In
general, if two morphisms in a category have a pushout, it is unique up to
isomorphism. Since the only isomorphisms in \(\mathcal D_\infty\) are identity
functions, this implies the uniqueness of ultrapower pushouts up to equality.

We now begin the proof of \cref{Pushout}. The proof involves the following
auxiliary concept:
\begin{defn}\label{MinimalDefinition}
Suppose \(M_0\) and \(M_1\) are transitive models of ZFC. A pair of elementary
embeddings \((i_0,i_1) : (M_0,M_1)\to N\) to a transitive model \(N\) is {\it
minimal} if \(N = H^N(i_0[M_0]\cup i_1[M_1])\).\index{Minimal pair of
embeddings}\index{Comparison!minimal}
\end{defn}
In the context of ultrapower embeddings, minimality has the following alternate
characterization:
\begin{lma}\label{MinimalUltrapowers}
Suppose \(j_0 : V\to M_0\) and \(j_1 : V\to M_1\) are elementary embeddings and
\((i_0,i_1) : (M_0,M_1)\to N\) is a comparison of \((j_0,j_1)\). Suppose \(a\in
M_1\) is such that \(M_1 = H^{M_1}(j_1[V]\cup \{a\})\). Then \((i_0,i_1)\) is
minimal if and only if \(N = H^N(i_0[M_0]\cup \{i_1(a)\})\).\qed
\end{lma}
Embedded in any pair \((k_0,k_1) : (M_0,M_1)\to P\), there is a unique minimal
pair \((i_0,i_1) : (M_0,M_1)\to N\). This follows from a trivial hull argument:
\begin{lma}\label{MinimalExistence}
Suppose \((k_0,k_1) : (M_0,M_1)\to P\) is a pair of elementary embeddings. Then
there exists a unique minimal \((i_0,i_1) : (M_0,M_1)\to N\) admitting an
elementary embedding \(h: N\to P\) such that \(h\circ i_0 = k_0\) and \(h\circ
i_1 = k_1\).
\begin{proof}
Let \(H = H^P(k_0[M_0]\cup k_1[M_1])\). Let \(N\) be the transitive collapse of
\(H\). Let \(h : N\to P\) be the inverse of the transitive collapse. Let \(i_0 =
h^{-1}\circ k_0\) and \(i_1 = h^{-1}\circ k_1\). Then  \((i_0,i_1) :
(M_0,M_1)\to N\) and \(h\circ i_0 = k_0\) and \(h\circ i_1 = k_1\). Moreover
 \[h[H^N(i_0[M_0]\cup i_1[M_1])] = H^P(k_0[M_0]\cup k_1[M_1]) = h[N]\] which
 implies \(H^N(i_0[M_0]\cup i_1[M_1]) = N\) since \(h\) is injective. Thus
 \((i_0,i_1)\) is minimal.

Uniqueness is obvious; we omit the proof.
\end{proof}
\end{lma}

\begin{cor}[UA]\label{MinimalInternal}
Every pair of ultrapower embeddings of \(V\) has a minimal internal ultrapower
comparison.
\begin{proof}
	Suppose \(j_0 : V\to M_0\) and \(j_1 : V\to M_1\) are ultrapower embeddings.
	Fix an internal ultrapower comparison \((k_0,k_1) : (M_0,M_1)\to P\) of
	\((j_0,j_1)\). By \cref{MinimalExistence}, there is a minimal pair
	\((i_0,i_1) : (M_0,M_1)\to N\) and an elementary \(h : N\to P\) with
	\(h\circ i_0 = k_0\) and \(h\circ i_1 = k_1\). It follows immediately that
	\((i_0,i_1)\) is a comparison of \((j_0,j_1)\). By
	\cref{MinimalUltrapowers}, \(i_0\) is an ultrapower embedding of \(M_0\).
	Since \(k_0\) is close to \(M_0\) and \(h\circ i_0 = k_0\), \(i_0\) is close
	to \(M_0\). Thus \(i_0\) is a close ultrapower embedding of \(M_0\), so
	\(i_0\) is an internal ultrapower embedding of \(M_0\). Similarly \(i_1\) is
	an internal ultrapower embedding of \(M_1\). Thus \((i_0,i_1)\) is a minimal
	internal ultrapower comparison of \((j_0,j_1)\).
\end{proof}
\end{cor}

\begin{lma}\label{MinimalEmbeddingUnique}
Suppose \((k_0,k_1) : (M_0,M_1)\to P\) is a pair of elementary embeddings and
\((i_0,i_1) : (M_0,M_1)\to N\) is a minimal pair. Then there is at most one
elementary embedding \(h : N\to P\) such that \(h\circ i_0 = k_0\) and \(h\circ
i_1 = k_1\).
\begin{proof}
Suppose \(h,h': N\to P\) satisfy \(h\circ i_0 = h'\circ i_0 = k_0\) and \(h\circ
i_1 = h'\circ i_1 = k_1\). Then \(h\restriction i_0[M_0] = h'\restriction
i_0[M_0]\) and \(h\restriction i_1[M_1] = h'\restriction i_1[M_1]\). Since \(N =
H^N(i_0[M_0]\cup i_1[M_1])\), it follows that \(h = h'\).
\end{proof}
\end{lma}

\begin{lma}[UA]\label{MinimalPushout}
Suppose \(j_0: V\to M_0\) and \(j_1 : V\to M_1\) are ultrapower embeddings and
\((i_0,i_1) : (M_0,M_1)\to N\) is a minimal comparison of \((j_0,j_1)\). Then
\((i_0,i_1)\) is the pushout of \((j_0,j_1)\).
\begin{proof}
Suppose \(j_0 : V\to M_0\) and \(j_1 : V\to M_1\) are ultrapower embeddings.
Suppose \((k_0,k_1) : (M_0,M_1)\to P\) is a comparison of \((j_0,j_1)\). It
suffices to find an internal ultrapower embedding \(h : N\to P\) such that
\(h\circ i_0 = k_0\) and \(h\circ i_1 = k_1\); uniqueness is then immediate from
\cref{MinimalEmbeddingUnique}.

Fix  \(a\in M_1\) such that \(M_1 = H^{M_1}(j_1[V]\cup \{a\})\). By
\cref{MinimalUltrapowers},
\[N = H^N(i_0[M_0]\cup i_1(a))\] By the definition of \(=_S\), we have:
\[(N,i_1(a))  =_S (M_1,a) =_S (P,k_1(a))\] Thus by the transitivity of the seed
order, \((N,i_1(a)) =_S (P,k_1(a))\). Since the objects witnessing \((N,i_1(a))
=_S (P,k_1(a))\) are internal ultrapower embeddings of \(N\) and \(P\), which
are themselves internal ultrapowers of \(M_0\), it follows that \(M_0\)
satisfies \((N,i_1(a)) =_S (P,k_1(a))\). By the equivalence between the seed
order on models and the seed order on embeddings (\cref{GenSE}), \(M_0\)
satisfies \((i_0,i_1(a)) =_S (k_0,k_1(a))\). Applying \cref{SeedEq} in \(M_0\),
it follows that there is an internal ultrapower embedding \(h : N\to P\) such
that \(h\circ i_0 = k_0\) and \(h(i_1(a)) = k_1(a)\).

We claim \(h\circ i_1 = k_1\). Note that \(h\circ i_1\circ j_1 = h\circ i_0
\circ j_0 = k_0\circ j_0 = k_1\circ j_1\), so \(h\circ i_1\restriction j_1[V] =
k_1 \restriction j_1[V]\). Moreover \(h(i_1(a)) = k_1(a)\). Thus \[h\circ
i_1\restriction j_1[V]\cup \{a\} = k_1\restriction j_1[V]\cup \{a\}\] Since
\(M_1 = H^{M_1}(j_1[V]\cup \{a\})\), it follows that \(h \circ i_1 = k_1\), as
desired.

Thus \(h : N\to P\) is an internal ultrapower embedding with \(h\circ i_0 =
k_0\) and \(h\circ i_1 = k_1\). 
\end{proof}
\end{lma}

\begin{proof}[Proof of \cref{Pushout}] The existence of pushouts is an immediate
	consequence of \cref{MinimalInternal} and \cref{MinimalPushout}.
\end{proof}

The existence of least upper bounds in the Rudin-Frol\'ik order is a trivial
restatement of \cref{Pushout}:
\begin{cor}\label{Join}
Suppose \(U_0\) and \(U_1\) are countably complete ultrafilters. Suppose
\((i_0,i_1): (M_{U_0},M_{U_1})\to N\) is the pushout of \((j_{U_0},j_{U_1})\).
Suppose \(W\) is a countably complete ultrafilter such that \(j_W = i_0\circ j_0
= i_1\circ j_1\). Then \(W\) is the \(\D\)-minimum countably complete
ultrafilter \(W\gD U_0,U_1.\) 
\begin{proof}
The internal ultrapower embeddings \(i_0\) and \(i_1\) witness that \(U_0\D W\)
and \(U_1\D W\). Suppose \(U_0\D Z\) and \(U_1\D Z\). We will show \(W\D Z\).
Let \(k_0 : M_{U_0}\to M_{Z}\) and \(k_1 : M_{U_1}\to M_{Z}\) witness  \(U_0\D
Z\) and \(U_1\D Z\). Then since \((i_0,i_1)\) is a pushout and \((k_0,k_1):
(M_{U_0},M_{U_1})\to M_Z\), there is an internal ultrapower embedding \(h :
M_W\to M_Z\) such that \(h\circ i_0 = k_0\) and \(h\circ i_1 = k_1\). In
particular \(h\circ j_W = h\circ i_0\circ j_{U_0} = k_0\circ j_{U_0} = j_Z\), so
\(h\) witnesses that \(W\D Z\).
\end{proof}
\end{cor}

It is worth noting the following bound here:
\begin{prp}\label{UpperBoundBound}
	Suppose \(U_0\) and \(U_1\) are countably complete ultrafilters. If \(W\) is
	a minimal upper bound of \(U_0\) and \(U_1\) in the Rudin-Frol\'ik order,
	then \(\lambda_W = \max\{\lambda_{U_0},\lambda_{U_1}\}\).
	\begin{proof}
		Let \(\lambda = \max\{\lambda_{U_0},\lambda_{U_1}\}\). Let \(j_0 :V\to
		M_0\) and \(j_1 : V\to M_1\) be the ultrapowers by \(U_0\) and \(U_1\)
		respectively. There is a minimal comparison \((i_0,i_1) : (M_0,M_1)\to
		N\) of \((j_0,j_1)\) such that \(i_0\circ j_0 = i_1\circ j_1 = j_W\).
		Fix \(\alpha < j_0(\lambda)\) such that \(M_0 = H^{M_0}(j_0[V]\cup
		\{\alpha\})\). By \cref{MinimalUltrapowers}, \(N = H^N(i_1[M_1]\cup
		\{i_0(\alpha)\})\subseteq H^N(i_1[M_1]\cup i_1(j_1(\lambda)))\). It
		follows that \(\textsc{width}(i_1) \leq j_1(\lambda) + 1\). Therefore by
		our lemma on the width of the composition of two elementary embeddings
		(\cref{WidthLemma}), \(\textsc{width}(j_W) = \textsc{width}(i_1\circ
		j_1) = \lambda+1\). In other words, \(\lambda_W = \lambda\).
	\end{proof}
\end{prp}

We now show the existence of greatest lower bounds in the Rudin-Frol\'ik order.
In fact we do a bit better:

\begin{prp}\label{Meet}
Suppose \(A\) is a nonempty class of countably complete ultrafilters. Then \(A\)
has a greatest lower bound in the Rudin-Frol\'ik order.
\end{prp}

This follows purely abstractly from what we have proved so far. Recall that a
partial order \((P,\leq)\) has the local ascending chain condition if for any
\(p\in P\), there is no ascending sequence \(a_0 < a_1 < \cdots\) in \(P\) with
\(a_n \leq p\) for all \(n  < \omega\). 
\begin{lma}\label{Algebra}
Suppose \((P,\leq)\) is a join semi-lattice with a minimum element that
satisfies the local ascending chain condition. For any nonempty set \(A\subseteq
P\), \(A\) has a greatest lower bound in \(P\).
\begin{proof}
Consider the set \(B\subseteq P\) of lower bounds of \(A\). In other words, \[B
= \{b\in P : \forall a\in A\ b \leq a\}\] Since \(P\) has a minimum element,
\(B\) is nonempty. Since \(A\) is nonempty, fixing \(p\in A\), every element of
\(B\) lies below \(p\). Therefore by the local ascending chain condition, \(B\)
has a maximal element \(b_0\). (The ascending chain condition says that the
relation \(>\) is wellfounded on \(\{c \in P : c\leq p\}\), so the nonempty set
\(B\) has a \(>\)-minimal element, or equivalently a \(<\)-maximal element.)

We claim \(B\) is a directed subset of \((P,\leq)\). Suppose \(b,c\in B\). For
any \(a\in A\), by the definition of \(B\), \(b,c\leq a\), and therefore their
least upper bound \(b\vee c \leq a\). In other words, \(b\vee c\leq a\) for all
\(a\in A\), so \(b\vee c\in B\). This shows that \(B\) is directed.

Finally since \(b_0\) is a maximal element of the directed set \(B\), in fact
\(b_0\) is its maximum element.
\end{proof}
\end{lma}

\begin{proof}[Proof of \cref{Meet}] The Rudin-Frol\'ik order induces a partial
order on the isomorphism types of countably complete ultrafilters. This partial
order is a join semi-lattice by \cref{Join}, and it has the local ascending
chain condition by \cref{ACC}. It has a minimum element, namely the isomorphism
type of the principal ultrafilters. Therefore the conditions of \cref{Algebra}
are met (except that we are considering a set-like partial order instead of a
set, which makes no difference). This implies the proposition.
\end{proof}

Let us give another application of pushouts to the Rudin-Frol\'ik order. The
following characterization of the internal ultrapower embeddings of a pushout is
remarkably easy to prove:
\begin{thm}\label{PushoutInternal}\index{Pushout!internal ultrapower embeddings}
	Suppose \(j_0 : V\to M_0\) and \(j_1 : V\to M_1\) are ultrapower embeddings
	and \((i_0,i_1) : (M_0,M_1)\to N\) is their pushout. Suppose \(h: N \to P\)
	is an ultrapower embedding. Then the following are equivalent:
	\begin{enumerate}[(1)]
		\item \(h\) is amenable to both \(M_0\) and \(M_1\).
		\item \(h\) is an internal ultrapower embedding of \(N\).
	\end{enumerate}
	\begin{proof}
		{\it (1) implies (2):} Let \(k_0 = h\circ i_0\) and \(k_1 = h\circ
		i_1\). Since \(h\) is an ultrapower embedding of \(N\), \(k_0\) is an
		ultrapower embedding of \(M_0\). Since \(h\) is amenable to \(M_0\),
		\(k_0\) is amenable to \(M_0\), and hence \(k_0\) is close to \(M_0\).
		Since \(k_0\) is a close ultrapower embedding of \(M_0\), in fact
		\(k_0\) is an internal ultrapower embedding of \(M_0\). Similarly
		\(k_1\) is an internal ultrapower embedding of \(M_1\). Thus
		\((k_0,k_1)\) is an internal ultrapower comparison of \((j_0,j_1)\).
		Since \((i_0,i_1)\) is a pushout, there is an internal ultrapower
		embedding \(h' : N\to P\) such that \(h' \circ i_0 = k_0\) and \(h'\circ
		i_1 = k_1\). By \cref{MinimalEmbeddingUnique}, however, \(h\) is the
		unique elementary embedding from \(N\) to \(P\) such that \(h\circ i_0 =
		k_0\) and \(h\circ i_1 = k_1\). Thus \(h = h'\), so \(h\) is an internal
		ultrapower embedding, as desired.
		
		{\it (2) implies (1):} Trivial.
	\end{proof}
\end{thm}

An elegant way to restate this is in terms of the ultrafilters amenable to a
pushout:

\begin{cor}\label{PushoutAmenable}
	Suppose \(j_0 : V\to M_0\) and \(j_1 : V\to M_1\) are ultrapower embeddings
	and \((i_0,i_1) : (M_0,M_1)\to N\) is their pushout. Suppose \(W\) is a
	countably complete \(N\)-ultrafilter. Then \(W\in N\) if and only if \(W\in
	M_0\cap M_1\).\qed
\end{cor}

\cref{PushoutInternal} permits an interesting generalization of the uniqueness
of ultrapower embeddings:

\begin{cor}[UA]\index{Rudin-Frol\'ik order!vs. inclusion of ultrapowers}
	 Suppose \(U\) and \(W\) are countably complete ultrafilters. Then the
	 following are equivalent:
	 \begin{enumerate}[(1)]
	 	\item \(U\D W\).
	 	\item \(M_W\subseteq M_U\).
	 \end{enumerate}
 \begin{proof}
 	{\it (1) implies (2):} Trivial.
 	
 	{\it (2) implies (1):} Let \((h.k): (M_U,M_W)\to N\) be the pushout of
 	\((j_U,j_W)\). Since \(M_W\subseteq M_U\) and \(k\) is an internal
 	ultrapower of \(M_U\), \(k\) is amenable to \(M_U\). In particular,
 	\(k\restriction N\) is amenable to both \(M_U\) and \(M_W\). Therefore
 	\(k\restriction N\) is an internal ultrapower of \(N\). Thus \(k\) is
 	\(\gamma\)-supercompact for all ordinals \(\gamma\). It follows from
 	\cref{UFSuperBound} that \(k\) is the identity. Hence \(h : M_U\to M_W\) is
 	an internal ultrapower embedding with \(h\circ j_U = j_W\), so \(U\D W\).
 \end{proof}
\end{cor}
\subsection{The finiteness of the Rudin-Frol\'ik order}\label{FiniteSection}
The goal of this subsection is to prove the central structural fact about the
Rudin-Frol\'ik order under UA: any countably complete ultrafilter has at most
finitely many predecessors in the Rudin-Frol\'ik order up to isomorphism. The
following terminology allows us to state this more precisely:
\begin{defn}
The {\it type} of an ultrafilter \(U\) is the class \(\{U' : U'\cong U\}\).
\end{defn}

\begin{thm}[UA]\label{RFFinite}\index{Rudin-Frol\'ik order!local finiteness}
If \(W\) is a countably complete ultrafilter, then \(\{U : U \D W\}\) is the
union of finitely many types.
\end{thm}

The proof heavily uses the concept of a Dodd parameter, introduced in
\cref{DoddSection} in a slightly more general context. Let us just remind the
reader what this is in the special case of ultrapower embeddings. We identify
finite sets of ordinals with their {\it decreasing} enumerations: if
\(p\subseteq\text{Ord}\) and \(|p| = \ell\), then \(\langle p_n : n <
\ell\rangle\) denotes the unique decreasing sequence such that \(p =
\{p_0,\dots,p_{\ell-1}\}\). The canonical order on finite sets of ordinals is
then the lexicographic order on their decreasing enumerations. 

\begin{defn}\label{DoddParamDef2}\index{Dodd parameter!of an ultrapower embedding}\index{\(p(j)\) (Dodd parameter)}
Suppose \(j : V\to M \) is an ultrapower embedding. The {\it Dodd parameter} of
\(j\), denoted \(p(j)\), is the least finite set of ordinals \(p\) such that
\(H^M(j[V]\cup p) = M\).
\end{defn}
Note that since \(j\) is an ultrapower embedding, \(M = H^M(j[V]\cup
\{\alpha\})\) for some ordinal \(\alpha\), so \(p(j)\) certainly exists. 

Recall the notion of an \(x\)-generator of an elementary embedding: if \(j :
M\to N\) is an elementary embedding between transitive models of ZFC and \(x\in
N\), then an ordinal \(\xi\in N\) is an \(x\)-generator of \(j\) if \(\xi\notin
H^N(j[V]\cup \xi\cup \{x\})\). We need a basic but not completely trivial lemma
about \(x\)-generators:
\begin{lma}\label{CompositionGenerators}
Suppose \(M\stackrel{j}{\longrightarrow} N\stackrel{i}{\longrightarrow} P\) are
elementary embeddings between transitive models and \(\xi\) is an
\(x\)-generator of \(j\). Then \(i(\xi)\) is an \(i(x)\)-generator of \(i\circ
j\).
\begin{proof}
Suppose not, and fix a function \(f\) and a finite set \(p\subseteq i(\xi)\)
such that \[i(\xi) = i(j(f))(p,i(x))\] Then \(P\) satisfies the statement that
for some finite set \(q\subseteq i(\xi)\), \(i(\xi) = i(j(f))(q,i(x))\). Since
\(i\) is elementary, \(N\) satisfies that for some finite set \(q\subseteq
\xi\), \(\xi = j(f)(q,x)\), and this contradicts that \(\xi\) is an
\(x\)-generator of \(j\).
\end{proof}
\end{lma}

The key lemma regarding the Dodd parameter is that each of its elements is a
generator in a strong sense:
\begin{lma}\label{DoddGenerators}
Suppose \(j : V\to M\) is an ultrapower embedding. Let \(p = p(j)\). Let \(\ell
= |p|\). Then for any \(n < \ell\), \(p_n\) is the largest \(p\restriction
n\)-generator of \(j\).
\begin{proof}
We first show that \(p_n\) is a \(p\restriction n\)-generator of \(j\). Suppose
not, towards a contradiction. Fix a finite set \(q \subseteq p_n\) such that
\(p_n \in H^M(j[V]\cup p\restriction n\cup q)\). Let \(r = (p\setminus \{p_n\})
\cup q\). Then \(r < p\) but \(p\in  H^M(j[V]\cup r)\). Therefore \[M =
H^M(j[V]\cup p) \subseteq H^M(j[V]\cup r)\] so \(H^M(j[V]\cup r) = M\), contrary
to the minimality of the Dodd parameter \(p\).

Now let \(\xi\) be the largest \(p\restriction n\)-generator of \(j\). Suppose
towards a contradiction that \(p_n < \xi\). Then \(p \subseteq \xi\cup
\{p_0,\dots,p_{n-1}\}\), so since \(\xi\notin H^M(j[V]\cup \xi\cup p\restriction
n)\), in fact \(\xi\notin H^M(j[V]\cup p)\). This contradicts the definition of
\(p(j)\).
\end{proof}
\end{lma}

The key to the proof of the finiteness of the Rudin-Frol\'ik order is to
partition the Rudin-Frol\'ik predecessors of a countably complete ultrafilter
according to their relationship with its Dodd parameter.
\begin{defn}\label{nUW}
Suppose \(U\sD W\) are countably complete ultrafilters. Let \(p = p(j_W)\). Let
\(i : M_U\to M_W\) be the unique internal ultrapower embedding such that
\(i\circ j_U = j_W\). Then \(n(U,W)\) is the least number \(n\) such that
\(p_n\notin i[M_U]\).
\end{defn}
Note that \(n(U,W)\) depends only on the types of \(U\) and \(W\). Note moreover
that \(n(U,W)\) exists whenever \(U\sD W\): otherwise \(p\subseteq i[M_U]\), so
\(M_W = H^{M_W}(j_W[V]\cup p)\subseteq i[M_U]\), which implies that \(i\) is
surjective; thus \(i\) is an isomorphism, so \(U\cong W\), contrary to the
assumption that \(U\sD W\).
\begin{lma}\label{DChar}
Suppose \(U\sD W\) are countably complete ultrafilters. Let \(p = p(j_W)\). Let
\(i : M_U\to M_W\) be the unique internal ultrapower embedding such that
\(i\circ j_U = j_W\). Let \(n = n(U,W)\). Then 
\[i[M_U]\subseteq H^{M_W}(j_W[V]\cup p\restriction n\cup p_n)\]
\begin{proof}
Suppose towards a contradiction that the lemma fails. Let \(\xi\) be the least
ordinal such that \(i(\xi)\notin H^{M_W}(j_W[V]\cup p\restriction n\cup p_n)\).
Then \(i[\xi]\subseteq H^{M_W}(j_W[V]\cup p\restriction n\cup p_n)\).

By the minimality of \(n\), \(p\restriction n\in i[M_U]\). Therefore let \(q\in
M_U\) be such that \(i(q) = p\restriction n\). We claim \(\xi\) is a
\(q\)-generator of \(j_U\). Supposing the contrary, we have \(\xi\in
H^{M_U}(j_U[V]\cup \xi\cup q)\), so \[i(\xi)\in i[H^{M_U}(j_U[V]\cup \xi\cup
q)]\subseteq H^{M_W}(j_W[V]\cup p\restriction n\cup p_n)\] which contradicts the
definition of \(\xi\).

Since \(\xi\) is a \(q\) generator of \(j_U\), \(i(\xi)\) is an
\(i(q)\)-generator of \(i\circ j_U\) by \cref{CompositionGenerators}. In other
words, \(i(\xi)\) is a \(p\restriction n\)-generator of \(j_W\). By
\cref{DoddGenerators}, \(p_n\) is the largest \(p\restriction n\)-generator of
\(j_W\), so \(i(\xi) \leq p_n\). This contradicts that \(i(\xi)\notin
H^{M_W}(j_W[V]\cup p\restriction n\cup p_n)\).
\end{proof}
\end{lma}

\begin{defn}
Suppose \(W\) is a  countably complete ultrafilter and \(p = p(j_W)\). For any
\(n < |p|\), \(\mathscr D_n(W) = \{U\sD W : n(U,W) = n\}\).
\end{defn}

\begin{lma}\label{DPartition}
For any countably complete ultrafilter \(W\), \[\{U: U\sD W\} = \bigcup_{n <
|p(j_W)|} \mathscr D_n(W)\]
\begin{proof}
See the remarks following \cref{nUW}.
\end{proof}
\end{lma}
The following fact is the key to the proof of the finiteness of the
Rudin-Frol\'ik order:
\begin{lma}\label{DDirected}
Suppose \(U_0,U_1\in \mathscr D_n(W)\) and \(D\) is the \(\D\)-minimum countably
complete ultrafilter such that \(U_0,U_1\D D\). Then \(D\in \mathscr D_n(W)\).
\begin{proof}
Let \(M_0 = M_{U_0}\) and let \(M_1 = M_{U_1}\). Let \((i_0,i_1) : (M_0,M_1)\to
M_D\) be internal ultrapower embeddings witnessing that \(U_0,U_1\D D\) and let
\((k_0,k_1) : (M_0,M_1)\to M_W\) be internal ultrapower embeddings witnessing
that \(U_0,U_1\D W\). 

Since \(D\) is the \(\D\)-minimum countably complete ultrafilter with
\(U_0,U_1\D D\), in fact \(D\D W\). Let \(h : M_D\to M_W\) be the unique
internal ultrapower embedding such that \(h\circ j_D = j_W\). Notice that 
\begin{align*}h\circ i_0 &= k_0\\ 
h\circ i_1 &= k_1\end{align*}
by \cref{UniqueFactor}.

Since \(D\) is the \(\D\)-minimum ultrafilter with \(U_0,U_1\D D\), \((i_0,i_1)
: (M_0,M_1)\to M_D\) must be minimal in the sense of \cref{MinimalDefinition}:
\[M_D = H^{M_D}(i_0[M_0]\cup i_1[M_1])\] (The proof is a trivial diagram chase.
Let \((\bar i_0,\bar i_1) :(M_0,M_1)\to N\) be the unique minimal pair admitting
\(e : N \to M_D\) such that \(e\circ \bar i_0 = i_0\) and \(e\circ \bar i_1 =
i_1\). By the proof of \cref{MinimalInternal}, \(N\) is an internal ultrapower
of \(M_0\) and \(M_1\), so since \(D\) is a least upper bound of \(U_0,U_1\),
there is an internal ultrapower embedding \(d: M_D\to N\) such that \(d \circ
i_0 = \bar i_0\) and \(d\circ i_1 = \bar i_1\). Then \(d\circ e : N\to N\)
satisfies \(d\circ e \circ \bar i_0 = \bar i_0\) and \(d\circ e \circ \bar i_1 =
\bar i_1\), and hence by \cref{MinimalEmbeddingUnique}, \(d\circ e\) must be the
identity map. Hence \(d\) and \(e\) are inverses, so by transitivity \(N = M_D\)
and \(e\) is the identity. Now \(\bar i_0 = e\circ \bar i_0   = i_0\) and \(\bar
i_1 = e \circ \bar i_1 = i_1\) so \((\bar i_0,\bar i_1) = (i_0,i_1)\). Since
\((\bar i_0,\bar i_1)\) is minimal, so is \((i_0,i_1)\).) Therefore \[h[M_D] =
h[H^{M_D}(i_0[M_0]\cup i_1[M_1])] = H^{M_W}(k_0[M_0]\cup k_1[M_1])\]

Let \(p = p(j_W)\). Since \(U_0\in \mathscr D_n(W)\), \(k_0[M_0]\subseteq
H^{M_W}(j_W[V]\cup p\restriction n\cup p_n)\) by \cref{DChar}. Similarly,
\(k_1[M_1]\subseteq H^{M_W}(j_W[V]\cup p\restriction n\cup p_n)\). Thus
\[k_0[M_0]\cup k_1[M_1]\subseteq H^{M_W}(j_W[V]\cup p\restriction n\cup p_n)\]
It follows that \(h[M_D] = H^{M_W}(k_0[M_0]\cup k_1[M_1])\subseteq
H^{M_W}(j_W[V]\cup p\restriction n\cup p_n)\). In particular, since \(p_n\) is a
\(p\restriction n\)-generator of \(j_W\) by \cref{DoddGenerators}, \(p_n\notin
h[M_D]\). Clearly \[p\restriction n\in k_0[M_0]\subseteq h[M_D]\] so \(n\) is
the least number such that \(p_n\notin h[M_D]\). It follows that \(n(D,W) = n\).
In other words, \(D\in \mathscr D_n(W)\).
\end{proof}
\end{lma}

The point now is that by \cref{ACC} and \cref{Join}, we can find a maximum
element of \(\mathscr D_n\):
\begin{prp}[UA]\label{DMax}
Suppose \(W\) is a countably complete ultrafilter and \(n < |p(j_W)|\). If
\(\mathscr D_n(W)\) is nonempty, then \(\mathscr D_n(W)\) has a \(\D\)-maximum
element.
\begin{proof}
By \cref{Join}, every pair of countably complete ultrafilters has a least upper
bound in the Rudin-Frol\'ik order. Combining this with \cref{DDirected}, the
class \(\mathscr D_n(W)\) is directed under \(\D\). Moreover it is bounded below
\(W\) in \(\D\). Therefore by \cref{ACC}, it has a maximal element \(U\). By the
\(\D\)-directedness of \(\mathscr D_n(W)\), this maximal element is a maximum
element.
\end{proof}
\end{prp}
We finally prove \cref{RFFinite}. 

\begin{proof}[Proof of \cref{RFFinite}] The proof is by induction on the
wellfounded relation \(\sD\). (See \cref{RFWF}.) Assume \(W\) is a countably
complete ultrafilter. Our induction hypothesis is that for all \(U\sD W\), \(\{D
: D\D U\}\) is the union of finitely many types. We aim to show that \(\{U : U\D
W\}\) is the union of finitely many types.

Let \(p = p(j_W)\) and let \(\ell = |p(j_W)|\). By \cref{DPartition},
\[\{U : U\sD W\} = \bigcup_{n < \ell} \mathscr D_n(W)\] We claim that for any
\(n < \ell\), \(\mathscr D_n(W)\) is the union of finitely many types. If
\(\mathscr D_n(W)\) is empty, this is certainly true. If \(\mathscr D_n(W)\) is
nonempty, then by \cref{DMax}, there is a \(\D\)-maximum element \(U\) of
\(\mathscr D_n(W)\). Since \(U\in \mathscr D_n(W)\), \(U\sD W\) so by our
induction hypothesis \(\{D : D\D U\}\) is the union of finitely many types. But
since \(U\) is a \(\D\)-maximum element of \(\mathscr D_n(W)\), \(\mathscr
D_n(W) \subseteq \{D : D\D U\}\). Thus \(\mathscr D_n(W)\) is the union of
finitely many types.

Since \(\{U : U\sD W\} = \bigcup_{n < \ell} \mathscr D_n(W)\) is a finite union
of classes \(\mathscr D_n(W)\) each of which is itself the union of finitely
many types,  \(\{U : U\sD W\}\) is the union of finitely many types. The
collection \(\{U : U\D W\}\) contains just one more type than  \(\{U : U\sD
W\}\), namely that of \(W\). So \(\{U : U\D W\}\) is the union of finitely many
types, completing the induction.
\end{proof}
\subsection{Translations and limits}\label{TranslationSection}
In this section we explain the relationship between pushouts, ultrafilter
translations, and the minimal covers defined for the proof of UA from the
linearity of the Ketonen order in \cref{LinKetSection}.

Recall \cref{RFTrans}, which defined for any countably complete ultrafilters
\(U\D W\) the translation of \(W\) by \(U\), the canonical countably complete
ultrafilter of \(M_U\) that leads from \(M_U\) into \(M_W\). It turns out that
there is a natural generalization of \(\tr U W\) for any ultrafilters that admit
a pushout:
\begin{defn}\index{\(\tr U W\) (translation)}\index{Translation of an ultrafilter (\(\tr U W\))!associated to the pushout}
Suppose \(U\) and \(W\in \mathscr B(Y)\) are countably complete ultrafilters.
Suppose \((k,h) : (M_U,M_W)\to N\) is the pushout of \((j_U,j_W)\). Then \(\tr U
W\) denotes the \(M_U\)-ultrafilter on \(j_U(Y)\) derived from \(k\) using
\(h(\id_W)\).
\end{defn}
The point of this definition is that \(\tr U W\) is the canonical ultrafilter of
\(M_U\) giving rise to the \(M_U\)-side of the pushout of \((j_U,j_W)\):

\begin{lma}\label{TransChar}
	Suppose \(U\) and \(W\in \mathscr B(Y)\) are countably complete
	ultrafilters. Suppose \((k,h) : (M_U,M_W)\to N\) is the pushout of
	\((j_U,j_W)\). Then \(\tr U W\) is the unique ultrafilter \(Z\in
	j_U(\mathscr B(Y))\) such that \(j_Z^{M_U} = k\) and \(\id_Z =
	h(\id_W)\).\qed
\end{lma}

We will try to omit superscripts when we can:

\begin{cor}\label{TransCor}
	If \(U\) and \(W\) are countably complete ultrafilters, then \((j_{\tr U W},
	j_{\tr W U})\) is the pushout of \((j_U,j_W)\) if it exists.
\end{cor}

The notation \(\tr {U} {W}\) generalizes the notation \(\tr U W\) introduced in
\cref{RFTrans} when \(U\D W\). To see this, assume \(U\D W\) and let \(k :
M_U\to M_W\) be the unique internal ultrapower embedding of \(M_U\) such that
\(k\circ j_U = j_W\). Then \((k,\text{id}): (M_U,M_W)\to M_W\) is the pushout of
\((j_U,j_W)\), and hence \(\tr U W\) as we have defined it here is just the
\(M_U\)-ultrafilter derived from \(k\) using \(\id_W\), which is precisely \(\tr
U W\) as defined in \cref{RFTrans}.

It turns out that in the definition of a translation, one does not need to use
the pushout (as long as the pushout exists):
\begin{lma}
Suppose \(U\) and \(W\in \mathscr B(Y)\) are countably complete ultrafilters
such that the pair \((j_{U},j_{W})\) has a pushout. Let \((k,h) :
(M_{U},M_{W})\to P\) be a close comparison of \((j_{U},j_{W})\). Then \(\tr {U}
{W}\) is the \(M_{U}\)-ultrafilter on \(j_{U}(Y)\) derived from \(k\) using
\(h(\id_W)\).
\qed
\end{lma}

It is not hard to see that translations are isomorphism invariant:
\begin{lma}
Suppose \(U\cong U'\) and \(W\cong W'\). Then \(\tr {U} {W} \cong \tr {U'}
{W'}\) in \(M_U\).
\end{lma}
In fact, we can do quite a bit better than this: translation functions preserve
the Rudin-Frol\'ik order.
\begin{lma}\label{TransRF}
Suppose \(U\), \(W\), and \(Z\) are countably complete ultrafilters. If \(W\D
Z\), then \(\tr {U} {W}\D \tr {U} {Z}\) in \(M_{U}\).
\begin{proof}
	Let \(N = M_{\tr U W}^{M_U}\) and let \(P= M_{\tr U Z}^{M_U}\). The proof is
	contained in \cref{TransRFFig}.
	\begin{figure}
		\center
		\includegraphics[scale=.65]{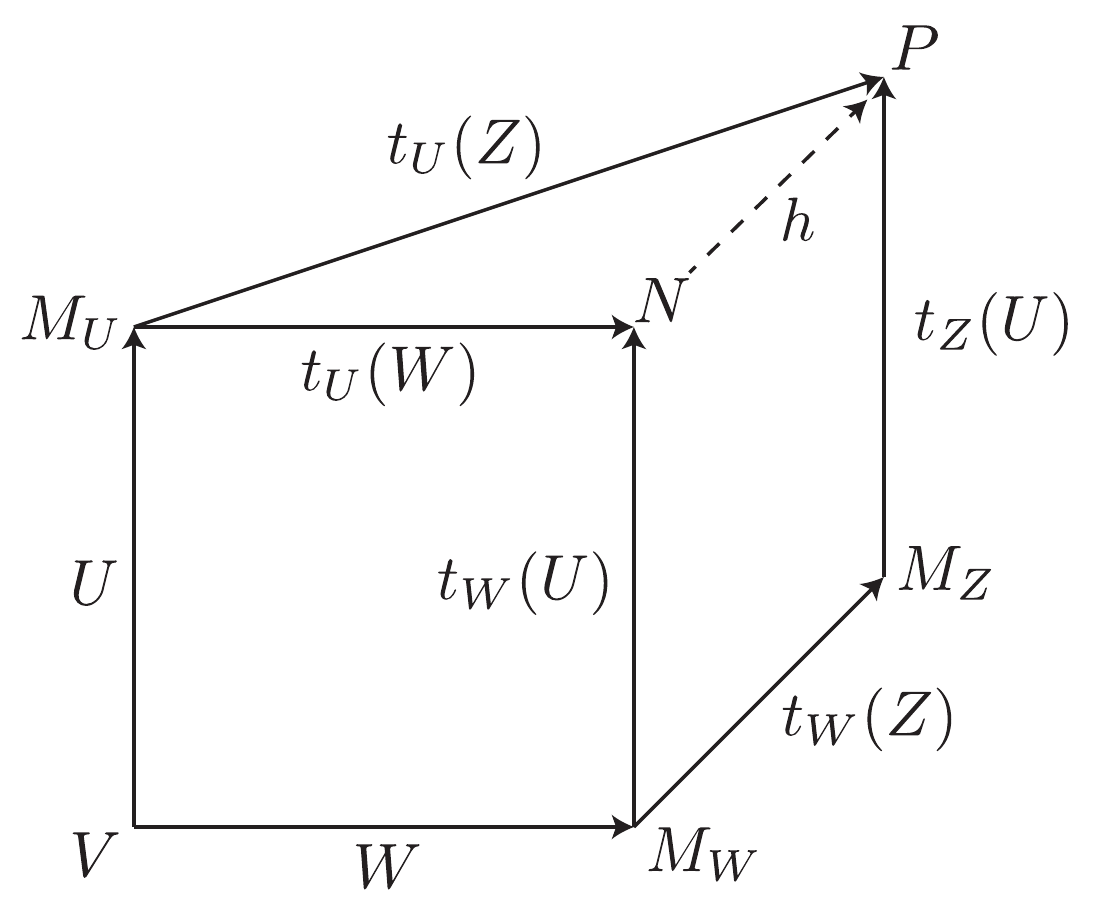}
		\caption{The proof of \cref{TransRF}.}\label{TransRFFig}
	\end{figure}
	By \cref{TransCor}:
	\begin{itemize}
		\item \((j_{\tr U W}, j_{\tr W U}) : (M_U,M_W)\to N\) is the pushout of
		\((j_U,j_W)\). 
		\item \((j_{\tr U Z}, j_{\tr Z U}\circ j_{\tr W Z}): (M_U,M_W)\to P\) is
		an internal ultrapower comparison of \((j_U,j_W)\).
	\end{itemize}
	Since \((j_{\tr U W}, j_{\tr W U})\) is a pushout, there is an internal
	ultrapower embedding \(h : N\to P\) such that \(h\circ j_{\tr U W} = j_{\tr
	U Z}\) and \(h \circ j_{\tr W U} =  j_{\tr Z U}\circ j_{\tr W Z}\). In
	particular, the first of these equations says that \(h\) witnesses \(\tr U
	W\D \tr U Z\) in \(M_U\).
\end{proof}
\end{lma}

We occasionally use the following fact, which is an immediate consequence of
\cref{TransChar}:
\begin{lma}\label{TransRF0}
Suppose \(U\) and \(W\) are countably complete ultrafilters on \(X\) and \(Y\).
Then the following are equivalent:
\begin{enumerate}[(1)]
	\item \(U\D W\).
	\item For some \(I\in U\) and some discrete sequence \(\langle W_i : i\in
	I\rangle\) of countably complete ultrafilters on \(Y\), \(\tr U W = [\langle
	W_i : i\in I\rangle]_U\).
	\item \(j_{\tr U W}\circ j_U = j_W\).
	\item \(\tr W U\) is a principal ultrafilter of \(M_W\).
	\item \(\tr W U = \pr {h(\id_U)} {j_W(X)}\) where  \(h : M_U\to M_W\) is the
	unique internal ultrapower embedding such that \(h\circ j_U = j_W\).\qed
\end{enumerate}
\end{lma}

The following fundamental fact connects translations back to the minimal covers
of \cref{LinKetSection}:
\begin{thm}[UA]\label{TranslationExtension}\label{Reciprocity}\index{Translation of an ultrafilter (\(\tr U W\))! as the minimum extension of \(j_U[W]\)}
	Suppose \(\delta\) is an ordinal, \(U\) is a countably complete ultrafilter,
	and \(W\in \mathscr B(\delta)\). Then \(\tr U W\) is the least element of
	\(j_U(\mathscr B(\delta),\sE)\) that extends \(j_U[W]\).
	\begin{proof}
		By replacing \(U\) with an isomorphic ultrafilter, we may assume that
		for some ordinal \(\epsilon\), \(U\in \mathscr B(\epsilon)\), putting us
		in a position to apply the results of \cref{LinKetSection}. 
		
		Let \(W_*\) be the least element of \(j_U(\mathscr B(\delta),\sE)\) that
		extends \(j_U[W]\) and let \(U_*\) be the least element of
		\(j_W(\mathscr B(\epsilon),\sE)\) that extends \(j_W[U]\). By
		\cref{Reciprocity1}, \[(j_{W_*}^{M_U},j_{U_*}^{M_W}) : (M_U,M_W)\to N\]
		is a comparison of \((j_U,j_W)\). Moreover, as a consequence of
		\cref{PreciseReciprocity}, \(\id_{W_*} = j_{U_*}^{M_W}(\id_W)\). In
		particular, \[N = H^N(j_{W_*}^{M_U}[M_U]\cup \{\id_{W_*}\}) =
		H^N(j_{W_*}^{M_U}[M_U]\cup \{j_{U_*}^{M_W}(\id_W)\})\] It follows from
		\cref{MinimalUltrapowers} that \((j_{W_*}^{M_U},j_{U_*}^{M_W})\) is
		minimal. Therefore by \cref{MinimalPushout},
		\((j_{W_*}^{M_U},j_{U_*}^{M_W})\) is the pushout of \((j_U,j_W)\). Since
		\(W_*\) is the \(M_U\)-ultrafilter on \(j_U(\delta)\) derived from
		\(j_{W_*}^{M_U}\) using \(j_{U_*}^{M_W}(\id_W)\), by definition \(W_* =
		\tr U W\).
	\end{proof}
\end{thm}

This yields the following bound on \(\tr U W\) that is not a priori obvious:
\begin{cor}[UA]\label{Bounding}
	Suppose \(U\) is a countably complete ultrafilter and \(W\) is a countably
	complete ultrafilter on an ordinal. Then \(\tr U W\E j_U(W)\) in \(M_U\).
	\begin{proof}
		Let \(\delta\) be the underlying ordinal of \(W\). Then \(j_U(W)\in
		j_U(\mathscr B(\delta))\) and \(j_U[W]\subseteq j_U(W)\). Thus \(\tr U W
		\E j_U(W)\) in \(M_U\) by \cref{TranslationExtension}.
	\end{proof}
\end{cor}

We finally show that translation functions preserve the Ketonen order:
\begin{thm}[UA]\label{OrderPreserving}
Translation functions preserve the Ketonen order. More precisely, suppose \(Z\)
is a countably complete ultrafilter and \(U\) and \(W\) are countably complete
ultrafilters on ordinals. Then \(U\sE W\) if and only if \(M_Z\vDash \tr Z {U}
\sE \tr Z {W}\).
\end{thm}

For this we need the strong transitivity of the Ketonen order
(\cref{StrongTrans}). We actually use the following immediate corollary of
\cref{StrongTrans} and the characterization of limits in terms of inverse images
(\cref{LimitInverse}):  
\begin{lma}\label{StrongTransCor}
Suppose \(Z\) is an ultrafilter, \(\delta\) is an ordinal, and \(U,W\in \mathscr
B(\delta)\) satisfy \(U\sE W\). For any \(W_*\in j_Z(\mathscr B(\delta))\) with
\(j_Z[W]\subseteq W_*\), there is some \(U_*\in j_Z(\mathscr B(\delta))\) with
\(U_*\sE^{M_Z} W_*\) and \(j_Z[U]\subseteq U_*\).\qed
\end{lma}

With \cref{TranslationExtension} and \cref{StrongTransCor} in hand, we can prove
\cref{OrderPreserving}.
\begin{proof}[Proof of \cref{OrderPreserving}] Assume that \(U \sE W\) are
countably complete ultrafilters on ordinals. We will show \(\tr Z U \sE^{M_Z}
\tr Z W\). For ease of notation, we will assume (without real loss of
generality) that \(U, W\in \mathscr B(\delta)\) for a fixed ordinal \(\delta\). 

Let \(W_* = \tr Z W\). \cref{TranslationExtension} implies that
\(j_Z[W]\subseteq W_*\). (This is actually much easier to prove that
\cref{TranslationExtension}.) By \cref{StrongTransCor}, it follows that there is
some \(U_*\in j_Z(\mathscr B(\delta))\) with \[U_*\sE^{M_Z} W_*\] and
\(j_Z[U]\subseteq U_*\). Since \(\tr Z U\) is the minimal extension of
\(j_Z[U]\) by \cref{TranslationExtension}, we have \[\tr Z U \E^{M_Z} U_*\] By
the transitivity of the Ketonen order, \(\tr Z U\E^{M_Z} \tr Z W\), as desired.
\end{proof}

\section{The internal relation}\label{InternalSection}
\subsection{A generalized Mitchell order}
In this section, we introduce a variant of the generalized Mitchell order that
will serve as a powerful tool in the theory of countably complete ultrafilters.
The trouble with using the Mitchell order itself to prove general theorems about
countably complete ultrafilters is that the Mitchell order is only meaningful
for ultrafilters that have a certain amount of strength: a precondition for
\(U\mo W\) is that \(P(\lambda_U)\subseteq M_W\). In order to analyze a more
general class of ultrafilters, we need a way to talk about the Mitchell order on
ultrafilters that are not assumed to be strong. 

There are a number of possible approaches, but the one that has proved most
successful is called the {\it internal relation:}
\begin{defn}\index{Internal relation}
	The {\it internal relation} is defined on countably complete ultrafilters
	\(U\) and \(W\) by setting \(U\I W\) if \(j_U\restriction M_W\) is an
	internal ultrapower embedding of \(M_W\).
\end{defn}
The topic of this section is the theory of the internal relation under UA. The
reason that we have saved it for this chapter is that it is closely related to
the theory of pushouts from \cref{PushoutSection}. 

Before we proceed through the basic theory below, let us mention that the
supercompactness analysis of \cref{SCChapter1} and \cref{SCChapter2} yields a
set theoretically simpler description of the internal relation on a very large
class of ultrafilters. In fact, the internal relation and the Mitchell order are
essentially one and the same:
\begin{repthm}{InternalChar}[UA] Suppose \(U\) and \(W\) are hereditarily
	uniform irreducible ultrafilters. Then the following are equivalent:
	\begin{enumerate}[(1)]
		\item \(U\I W\).
		\item Either \(U\mo W\) or \(W\in V_\kappa\) where \(\kappa =
		\textsc{crt}(j_U)\).
	\end{enumerate}
\end{repthm}
The second part of condition (2) should be compared with Kunen's commuting
ultrapowers lemma (\cref{KunenCommute}). 
\subsection{The Mitchell order versus the internal relation}
To understand the nature of the internal relation, it helps to consider its
relationship with the Mitchell order.
\begin{prp}\label{InternalToMO}\index{Internal relation!vs. the generalized Mitchell order}
	Suppose \(U\) is a countably complete ultrafilter on a set \(X\) and \(W\)
	is a countably complete ultrafilter such that \(X \in M_W\) and \(U\I W\).
	Then the \(M_W\)-ultrafilter \(U\cap M_W\) belongs to \(M_W\). In
	particular, if \(P(X)\subseteq M_W\), then \(U\mo W\).\qed
\end{prp}

In general, however, \(U\I W\) does not imply \(U\mo W\). This is a consequence
of Kunen's commuting ultrapowers lemma (\cref{KunenCommute}):
\begin{prp}
	Suppose \(\kappa\) is a measurable cardinal, \(U\in V_\kappa\) is a
	countably complete ultrafilter and \(W\) is a \(\kappa\)-complete
	ultrafilter. Then \(W\I U\).\qed
\end{prp}
Note that in the situation above, if \(W\) is nonprincipal, then \(\lambda_W
\geq \kappa\), and in particular \(W\not \mo U\) since \(P(\kappa)\nsubseteq
M_U\). 

Whether \(U\mo W\) always implies \(U\I W\) is a considerably subtler question.
This implication is consistently false. (This is closely related to
\cref{USupercompact}.) We begin with the following fact:
\begin{prp}\label{Cummings0}
	Suppose \(\kappa\) is \(2^\kappa\)-supercompact and \(2^{\kappa} =
	2^{(\kappa^+)}\). Then there is a normal ultrafilter \(D\) on \(\kappa\) and
	a \(\kappa\)-complete normal fine ultrafilter \(\mathcal U\) on
	\(P_\kappa(\kappa^+)\) such that \(\mathcal U\mo D\).
	\begin{proof}[Sketch] Since \(\kappa\) is \(\kappa^+\)-supercompact, there
		is a \(\kappa\)-complete  normal fine ultrafilter \(\mathcal U\) on
		\(P_\kappa(\kappa^+)\). By Solovay's theorem on SCH above a strongly
		compact cardinal (\cref{SolovayThm}), \(|P_\kappa(\kappa^+)| =
		\kappa^+\). By Solovay's ultrafilter-capturing  theorem
		(\cref{SolovayUltrafilters}), for any set \(A\) of hereditary
		cardinality at most \(2^\kappa\), there is a normal ultrafilter \(D\) on
		\(\kappa\) such that \(A\in M_D\). But \(\mathcal U\subseteq
		P(P_\kappa(\kappa^+))\) has hereditary cardinality \(2^{\kappa^+}
		=2^\kappa\). Thus there is a normal ultrafilter \(D\) on \(\kappa\) such
		that \(\mathcal U\in M_D\), or in other words, \(\mathcal U\mo D\).
	\end{proof}
\end{prp}
Thus given a failure of the weak GCH at a supercompact, one must have a rather
unusual situation in which \(\mathcal U\mo D\) even though \(\lambda_\mathcal U
> \lambda_D\). On the other hand, the internal relation does not hold between
these ultrafilters:
\begin{prp}\label{CummingsExample} Assume \(D\) is a \(\kappa\)-complete uniform
ultrafilter on \(\kappa\) and \(\mathcal U\) is a \(\kappa\)-complete normal
fine ultrafilter on \(P_\kappa(\kappa^+)\).\footnote{The proof only requires
that \(\mathcal U\) is a \(\kappa\)-complete fine ultrafilter on
\(P_\kappa(\kappa^+)\).} Then \(\mathcal U \not \I D\). 
	\begin{proof}
		Suppose towards a contradiction that \(\mathcal U\I D\). Then
		\(j_\mathcal U(M_D) \subseteq M_D\) since \(j_\mathcal U\restriction
		M_D\) is an internal ultrapower embedding of \(M_D\). But \(j_\mathcal
		U(M_D) = (M_{j_{\mathcal U}(D)})^{M_\mathcal U}\). Since \(j_\mathcal
		U(D)\) is \(j_\mathcal U(\kappa)\)-complete in \(M_\mathcal U\), 
		\[\text{Ord}^{(\kappa^+)}\subseteq \text{Ord}^{j_D(\kappa)}\cap
		M_\mathcal U\subseteq M_{j_{\mathcal U}(D)})^{M_\mathcal U} = j_\mathcal
		U(M_D) \subseteq M_D\] It follows that \(j_D\) is
		\(\kappa^+\)-supercompact, and this contradicts the bound on the
		supercompactness of the ultrapower by an ultrafilter on \(\kappa\)
		(\cref{UFSuperBound}).
	\end{proof}
\end{prp}
We have not checked that the implication from \(U\mo W\) to \(U\I W\) can fail
under the Generalized Continuum Hypothesis, but we are confident that it can.
Under UA, however, this implication is a theorem:
\begin{repthm}{UltrapowerCorrectness}[UA] Suppose \(U\) and \(W\) are countably
	complete ultrafilters. If \(U\mo W\), then \(U\I W\).\qed
\end{repthm}
\subsection{Basic theory of the internal relation}
The true motivation for the definition of the internal relation comes from the
theory of ultrapower comparisons:

\begin{lma}\label{InternalComparison}
	Suppose \(U\) and \(W\) are countably complete ultrafilters. Then
	\[(j_U(j_W),j_U\restriction M_W) : (M_U,M_W)\to j_U(M_W)\] is a
	\(0\)-internal minimal comparison of \((j_U,j_W)\). It is an internal
	ultrapower comparison if and only if \(U\I W\).
	\begin{proof}
		The fact that \((j_U(j_W),j_U\restriction M_W)\) is a comparison of
		\((j_U,j_W)\) is immediate from the standard application-composition
		identity:
		\[j_U(j_W)\circ j_U = (j_U \restriction M_W)\circ j_W\] Since \(j_W\) is
		an internal ultrapower embedding of \(V\), \(j_U(j_W)\) is an internal
		ultrapower embedding of \(M_U\) by the elementarity of \(j_U\). In
		particular, \((j_U(j_W),j_U\restriction M_W)\) is \(0\)-internal.
		Moreover, if \(U\I W\) then \(j_U\restriction M_W\) is an internal
		ultrapower embedding of \(M_W\), and hence  \((j_U(j_W),j_U\restriction
		M_W)\) is an internal ultrapower comparison.
		
		Let us finally show that  \((j_U(j_W),j_U\restriction M_W)\) is a
		minimal comparison of \((j_U,j_W)\), or in other words that \[j_U(M_W) =
		H^{j_U(M_W)}(j_U(j_W)[M_U]\cup j_U[M_W])\] The proof begins with the
		standard fact that \(M_W = H^{M_W}(j_W[V]\cup \{\id_W\})\). Applying
		\(j_U\) to both sides of the equation, we obtain: \[j_U(M_W) =
		H^{j_U(M_W)}(j_U(j_W)[M_U]\cup \{j_U(\id_W)\})\] Since \(j_U(\id_W)\in
		j_U[M_W]\), \[H^{j_U(M_W)}(j_U(j_W)[M_U]\cup \{j_U(\id_W)\})\subseteq
		H^{j_U(M_W)}(j_U(j_W)[M_U]\cup j_U[M_W])\] This yields that
		\(j_U(M_W)\subseteq H^{j_U(M_W)}(j_U(j_W)[M_U]\cup j_U[M_W])\), which of
		course implies that equality holds, as desired.
\end{proof}
\end{lma}

Combining \cref{InternalComparison} with the fact that minimal comparisons of
ultrapowers are ultrapower comparisons (\cref{MinimalUltrapowers}), we obtain
the following lemma:
\begin{lma} Suppose \(U\) and \(W\) are countably complete ultrafilters. Then \(j_U\restriction M_W\) is an ultrapower embedding of \(M_W\).\end{lma} 
Of course, we do not mean that \(j_U\restriction M_W\) is necessarily an {\it
internal} ultrapower embedding of \(M_W\), just that there is a point \(a\in
j_U(M_W)\) such that \(j_U(M_W) = H^{j_U(M_W}(j_U[M_W]\cup \{a\})\). An
important point is that this point \(a\) need not be \(\id_U\) itself.

Applying the proof of \cref{MinimalUltrapowers} in to the minimal comparison
\((j_U(j_W), j_U\restriction M_W)\) identifies a specific \(M_W\)-ultrafilter
giving rise to the embedding \(j_U\restriction M_W\):

\begin{defn}\index{\(s_W(U)\) (pushforward)}\index{Pushforward of an ultrafilter into an ultrapower \(s_W(U)\)}\index{Internal relation!\(s_W(U)\)}
	Suppose \(U\) and \(W\) are countably complete ultrafilters. Let \(X\) be
	the underlying set of \(U\). Then the {\it pushforward of \(U\) into
	\(M_W\)} is the \(M_W\)-ultrafilter \(s_W(U)\) on \(j_W(X)\) defined as
	follows: if \(A\subseteq j_W(X)\) and \(A\in M_W\), \[A\in s_W(U) \iff
	j_W^{-1}[A]\in U\]
\end{defn}

The reason we call \(s_W(U)\) a pushforward is that it is literally equal to the
pushforward \(f_*(U)\cap M_W\) where \(f : X\to j_W(X)\) is the restriction \(f
= j_W\restriction X\).

For the reader's convenience, let us chase through all the lemmas and prove that
\(s_W(U)\) behaves as it should:

\begin{lma}\label{PushUlt}
	Suppose \(U\) and \(W\) are countably complete ultrafilters on \(X\) and
	\(Y\) Then \(s_W(U)\) is the \(M_W\)-ultrafilter on \(j_W(X)\) derived from
	\(j_U\restriction M_W\) using \(j_U(j_W)(\id_U)\). Moreover,
	\[j_{s_W(U)}^{M_W} = j_U\restriction M_W\] Thus \(U\I W\) if and only if
	\(s_W(U)\in M_W\).
	\begin{proof}
		Let \(f = j_W\restriction X\). Then \(f_*(U)\) is the ultrafilter
		derived from \(j_U\) using \(j_U(f)(\id_U)\) by the basic theory of the
		Rudin-Keisler order (\cref{PushDerived}). Thus \(f_*(U)\cap M_W\) is the
		\(M_W\)-ultrafilter derived from \(j_U\restriction M_W\) using
		\(j_U(f)(\id_U) = j_U(j_W)(\id_U)\). But \(f_*(U)\cap M_W = s_W(U)\), so
		\(s_W(U)\) is the \(M_W\)-ultrafilter on \(j_W(X)\) derived from
		\(j_U\restriction M_W\) using \(j_U(j_W)(\id_U)\).
		
		We finish by proving \(j_{s_W(U)}^{M_W} = j_U\restriction M_W\). Since
		\(s_U(W)\) is derived from \(j_U\restriction M_W\) using
		\(j_U(j_W)(\id_U)\), there is a factor embedding \(k : M_{s_W(U)}^{M_W}
		\to j_U(M_W)\) with \(k\circ j_{s_W(U)}^{M_W} = j_U\restriction M_W\)
		and \(k(\id_{s_W(U)}) = j_U(j_W)(\id_U)\). Since
		\((j_U(j_W),j_U\restriction M_W) : (M_U,M_W)\to j_U(M_W)\) is a minimal
		comparison of \((j_U,j_W)\), \cref{MinimalUltrapowers} yields:
		\[j_U(M_W) = H^{j_U(M_W)}(j_U[M_W]\cup \{j_U(j_W)(\id_U)\})\] But \(
		H^{j_U(M_W)}(j_U[M_W]\cup \{j_U(j_W)(\id_U)\})\subseteq
		k[M_{s_W(U)}^{M_W}]\). In other words, \(k\) is a surjection. It follows
		that \(M_{s_W(U)}^{M_W} = j_U(M_W)\) and \(k\) is the identity.
		Therefore \(j_{s_W(U)}^{M_W} = k\circ j_{s_W(U)}^{M_W} = j_U\restriction
		M_W\) as desired.
	\end{proof}
\end{lma}

As a corollary, one can characterize the internal relation in terms of
amenability of ultrafilters.
\begin{lma}\label{InternalAmenable}
	Suppose \(U\) and \(W\) are countably complete ultrafilters. Then the
	following are equivalent:
	\begin{enumerate}[(1)]
		\item \(U\I W\).
		\item For all \(U'\RK U\), \(U'\cap M_W\in M_W\).
		\item For all \(U'\cong U\), \(U'\cap M_W\in M_W\).
	\end{enumerate}
	\begin{proof}
		{\it (1) implies (2):} Suppose \(U' \RK U\I W\). Fix a set \(X\) and a
		point \(a\in M_U\) such that \(U'\) is the ultrafilter on \(X\) derived
		from \(j_U\) using \(a\). If \(X\cap M_W\notin U'\), then \(U'\cap M_W =
		\emptyset\), and so \(U'\cap M_W\in M_W\) vacuously. Therefore assume
		\(X\cap M_W\in U'\). In other words, \(a\in j_U(X\cap M_W)\), so \(a\in
		j_U(M_W)\). Then \(U'\cap M_W\) is the ultrafilter derived from
		\(j_U\restriction M_W\) using \(a\), so since \(j_U\restriction M_W\) is
		an internal ultrapower embedding of \(M_W\), \(U'\cap M_W\in M_W\).
		
		{\it (2) implies (3):} Trivial.
		
		{\it (3) implies (1):} Let \(X\) be the underlying set of \(U\). Let \(f
		:X \to j_W(X)\) be the restriction \(f = j_W\restriction X\). Since
		\(j_W\) is injective, \(f_*(U) \cong U\). Moreover \(f_*(U)\cap M_W =
		s_W(U)\), so if \(f_*(U)\cap M_W\in M_W\), then \(U\I W\) by
		\cref{PushUlt}.
	\end{proof}
\end{lma}

This has the following corollary, which is perhaps not immediately obvious:

\begin{cor}\label{RKInternal}
	Suppose \(U\), \(W\), and \(Z\) are countably complete ultrafilters and \[Z
	\RK U\I W\] Then \(Z\I W\).
	\begin{proof}
		By \cref{InternalAmenable}, for all \(U'\RK U\), \(U'\I W\). In
		particular (by the transitivity of the Rudin-Keisler order), for all
		\(U'\RK Z\), \(U'\I W\). Applying \cref{InternalAmenable} again, \(Z\I
		W\), as desired. 
	\end{proof}
\end{cor}

There is also an obvious relationship in the other direction between the
Rudin-Frol\'ik order and the internal relation:
\begin{prp}\label{RFInternal}
 Suppose \(U\), \(W\), and \(Z\) are countably complete ultrafilters and \[U \D
 W\gI Z\] Then \(Z\I U\).
 \begin{proof}
 	Since \(Z\I W\), \cref{PushUlt} implies \(s_W(Z)\in M_W\). Since \(U\D W\),
 	there is an internal ultrapower embedding \(h : M_U\to M_W\). We claim that
 	\(h^{-1}[s_W(Z)] = s_U(Z)\). Let \(X\) be the underlying set of \(Z\). If
 	\(A\in j_U(P(X))\),
 	\begin{align*}A\in h^{-1}[s_W(Z)]&\iff h(A)\in s_W(Z)\\
 		 									&\iff j_W^{-1}[h(A)]\in Z\\
 											&\iff (h\circ j_U)^{-1}[h(A)]\in Z\\
 											&\iff j_U^{-1}[A]\in Z\\
 											&\iff A\in s_U(Z)
 	\end{align*}
 	Since \(h\) is definable over \(M_U\) and \(s_W(Z)\in M_W\subseteq M_U\),
 	\( s_U(Z)= h^{-1}[s_W(Z)]\in M_U\). Hence \(Z\I U\) by \cref{PushUlt}, as
 	desired.
 \end{proof}
\end{prp}

The key to understanding the internal relation under UA is the following
theorem, which takes advantage of the theory of pushouts and translations
(\cref{PushoutSection} and \cref{TranslationSection}):

\begin{lma}[UA]\label{IChar}\index{Internal relation!and ultrafilter translations}\index{Translation of an ultrafilter (\(\tr U W\))!and the internal relation}
	Suppose \(U\) and \(W\) are countably complete ultrafilters. Then the
	following are equivalent:
	\begin{enumerate}[(1)]
		\item \(U\I W\).
		\item \((j_U(j_W),j_U\restriction M_W)\) is the pushout of
		\((j_W,j_U)\).
		\item \(\tr U W  = j_U(W)\).
		\item \(\tr W U = s_W(U)\).
	\end{enumerate}
	If the underlying set of \(W\) is an ordinal, we can add to the list:
	\begin{enumerate}[(5)]
		\item \(M_U\vDash j_U(W)\E \tr U W\).	
	\end{enumerate}
\begin{proof}
	{\it (1) implies (2):} Since \(U\I W\),  \((j_U(j_W),j_U\restriction M_W)\)
	is a minimal internal ultrapower comparison of \((j_U,j_W)\). Therefore by
	\cref{MinimalPushout}, \((j_U(j_W),j_U\restriction M_W)\) is the pushout of
	\((j_U,j_W)\), so (2) holds.
	
	{\it (2) implies (3):} Let \(X\) be the underlying set of \(W\). By the
	definition of \(\tr U W\), \(\tr U W\) is the \(M_U\)-ultrafilter on
	\(j_U(X)\) derived from \(k\) using \(h(\id_W)\) where \((k,h) :
	(M_U,M_W)\to N\) is the pushout of \((j_U,j_W)\). By (2), \((k,h) =
	(j_U(j_W),j_U\restriction M_W)\), and hence \(\tr U W\) is the
	\(M_U\)-ultrafilter on \(j_U(X)\) derived from \(j_U(j_W)\) using
	\(j_U(\id_W)\). Since \(W\) is the ultrafilter on \(X\) derived from \(j_W\)
	using \(\id_W\), by the elementarity of \(j_U\), \(j_U(W)\) is the
	ultrafilter on \(j_U(X)\) derived from \(j_U(j_W)\) using \(j_U(\id_W)\).
	This yields that \(\tr U W = j_U(W) \), so (3) holds. 
	
	{\it (3) implies (4):} Let \((k,h) : (M_U,M_W)\to N\) be the pushout of
	\((j_U,j_W)\). Since \(\tr U W = j_U(W)\), \cref{TransChar} implies \(k =
	j_U(j_W)\) and \(h(\id_W) = \id_{j_U(W)} = j_U(\id_W)\). 
	
	We claim that \(h = j_U\restriction M_W\). Note that \(h\restriction j_W[V]
	= j_U\restriction j_W[V]\) since \(h\circ j_W = k\circ j_U = j_U(j_W)\circ
	j_U = j_U\circ j_W\). Moreover \(h(\id_W) = j_U(\id_W)\), so \[h\restriction
	j_W[V]\cup \{\id_W\} = j_U\restriction j_W[V]\cup \{\id_W\}\] Since \(M_W =
	H^{M_W}(j_W[V]\cup \{\id_W\})\) it follows that \(h = j_U\restriction M_W\),
	as claimed.
	
	Now \(\tr W U\) is the \(M_W\)-ultrafilter derived from \(h =
	j_U\restriction M_W\) using \(k(\id_U) = j_U(j_W)(\id_U)\). By
	\cref{PushUlt}, \(\tr W U = s_W(U)\).
	
	{\it (4) implies (1):} Since \(\tr W U = s_W(U)\), \(s_W(U)\in M_W\). By
	\cref{PushUlt}, \(U\I W\).
	
	Finally, assume that the underlying set of \(W\) is an ordinal \(\delta\),
	and we will show the equivalence of (3) and (5). Clearly (3) implies (5), so
	let us prove the converse. Assume (5) holds. By \cref{Bounding}, \(\tr U W
	\E j_U(W)\) in \(M_U\). Thus \(\tr U W\E j_U(W)\) and \(j_U(W)\E \tr U W\)
	in \(M_U\), so \(j_U(W) = \tr U W\) since the Ketonen order is
	antisymmetric.
\end{proof}
\end{lma}

\subsection{Commuting ultrapowers and wellfoundedness}\label{CommutingSection}
The comparison characterization of the internal relation
(\cref{InternalComparison}) leads to a connection between the internal relation
and the seed order on pointed ultrapower embeddings, which will give us some
insight into the wellfoundedness of the internal relation:
\begin{lma}\label{InternalSeed}\index{Internal relation!vs. the seed order}\index{Seed order!on pointed ultrapowers!vs. the internal relation}
	Suppose \(\delta\) is a limit ordinal and \(U\I W\) are countably complete
	ultrafilters. Then for any \(\alpha < j_U(\delta)\), \((j_U,\alpha) \swo
	(j_W,\sup j_W[\delta])\).
	\begin{proof}
		Since \(U\I W\), \((j_U(j_W),j_U\restriction M_W)\) is an internal
		ultrapower comparison of \((j_U,j_W)\) by \cref{InternalComparison}. To
		show that \((j_U,\alpha) \swo (j_W,\sup j_W[\delta])\), it therefore
		suffices to show that \(j_U(j_W)(\alpha) < (j_U\restriction M_W)(\sup
		j_W[\delta])\). Note however that \[(j_U\restriction M_W)(\sup
		j_W[\delta])  = j_U(\sup j_W[\delta]) = \sup j_U(j_W)[j_U(\delta)] >
		j_U(j_W)(\alpha)\] since \(\alpha < j_U(\delta)\).
	\end{proof}
\end{lma}

As an immediate corollary, we have that the seed order extends the internal
relation in many cases:

\begin{thm}\label{InternalKet}
	Suppose \(U\I W\) are ultrafilters concentrating on ordinals. Then \(U\swo
	W\) if and only if \(\delta_U \leq \delta_W\).\qed
\end{thm}

We also obtain a wellfoundedness theorem for the internal relation, which
becomes more interesting when one realizes that the internal relation is not in
fact wellfounded.

\begin{thm}\label{InternalWF}
	Suppose \(\delta\) is an ordinal. Then the internal relation is wellfounded
	on the class of countably complete ultrafilters whose ultrapowers are
	discontinuous at \(\delta\).
	\begin{proof}
		Suppose towards a contradiction \(U_0 \gI U_1\gI U_2\gI \cdots\) are all
		discontinuous at \(\delta\). For \(n < \omega\), let \(j_n : V\to M_n\)
		denote the ultrapower of the universe by \(U_n\), and let \(\delta_n =
		\sup j_n[\delta]\). Since \(\delta_{n+1} < j_{n+1}(\delta)\) and
		\(U_{n+1} \I U_n\), \cref{InternalSeed} implies \((j_{n+1},\delta_{n+1})
		\swo (j_n,\delta_n)\). Writing this a different way, we have:
		\[(j_0,\delta_0)\slwo(j_1,\delta_1)\slwo (j_2,\delta_2)\slwo\cdots\]
		This immediately contradicts the wellfoundedness of the Ketonen order on
		pointed models (\cref{GenWellfounded}).
	\end{proof}
\end{thm}

\begin{cor}\index{Internal relation!Irreflexivity}
	If \(U\) is a nonprincipal countably complete ultrafilter, then \(U\not \I
	U\).\qed
\end{cor}

Unlike the Mitchell order, the internal relation is not strict. In fact, it has
2-cycles, which typically come from the phenomenon of {\it commuting
ultrafilters}:

\begin{defn}\index{Commuting ultrafilters|seealso{Kunen's commuting ultrapowers lemma}}
	Suppose \(U\) and \(W\) are countably complete ultrafilters. Then \(U\) and
	\(W\) {\it commute} if \(j_U(j_W) = j_W\restriction M_U\) and \(j_W(j_U) =
	j_U\restriction M_W\).	
\end{defn}

Clearly if \(U\) and \(W\) commute, then \(U\I W\) and \(W\I U\). Let us provide
some obvious combinatorial characterizations of commuting ultrafilters:
\begin{lma}
	Suppose \(U\) and \(W\) are countably complete ultrafilters on sets \(X\)
	and \(Y\). The following are equivalent:
	\begin{enumerate}[(1)]
		\item \(U\) and \(W\) commute.
		\item For all \(A\subseteq X\times Y\), \(\forall^U x\ \forall^W y\
		(x,y)\in A\iff \forall^W y\ \forall^U x\  (x,y)\in A\).
		\item The function \(\textnormal{flip}(x,y) = (y,x)\) satisfies
		\(\textnormal{flip}_*(U\times W) = W\times U\).\qed
	\end{enumerate}
\end{lma}
Somewhat surprisingly, there are nontrivial examples of commuting ultrafilters:

\begin{thm}[Kunen]\label{KunenCommute}\index{Kunen's commuting ultrapowers lemma}
	Suppose \(U\) and \(W\) are countably complete ultrafilters and \(U\in
	V_\kappa\) where \(\kappa = \textsc{crt}(j_W)\). Then \(j_W(j_U) =
	j_U\restriction M_W\) and \(j_U(j_W) = j_W\restriction M_U\).
\end{thm}

Let us give our pet proof of \cref{KunenCommute}, which uses the following
somewhat surprising reformulation of commutativity:
\begin{prp}\label{WeakCommute}
	Suppose \(U\) and \(W\) are countably complete ultrafilters such that
	\(j_W(j_U) = j_U\restriction M_W\). Then \(U\) and \(W\) commute.
	\begin{proof}
		To show \(U\) and \(W\) commute, we must show that \(j_W\restriction M_U
		= j_U(j_W)\). By \cref{InternalComparison}, \((j_W\restriction M_U,
		j_W(j_U))\) and \(( j_U(j_W), j_U\restriction M_W)\) are 0-internal and
		1-internal minimal comparisons of \((j_U,j_W)\). Since \(j_W(j_U) =
		j_U\restriction M_W\), we can conclude that \[(j_W\restriction M_U)\circ
		j_U = j_U(j_W)\circ j_U\]
		
		In particular, \(j_W\restriction M_U\) and \(j_U(j_W)\) are elementary
		embeddings of \(M_U\) with the same target model, which we will denote
		by \[N = j_W(M_U) = j_U(j_W)(M_U) = j_U(M_W) = j_W(j_U)(M_W)\] so since
		\(j_U(j_W)\) is an internal ultrapower embedding of \(M_U\),
		\(j_U(j_W)(\alpha) \leq j_W(\alpha)\) for all ordinals \(\alpha\).
		
		Let \(\xi\) be the least ordinal such that \(M_U = H^{M_U}(j_U[V]\cup
		\{\xi\})\). We claim that \[j_W(\xi) = j_U(j_W)(\xi)\] By the previous
		paragraph, we have \(j_U(j_W)(\xi)\leq j_W(\xi)\), so it suffices to
		prove the reverse inequality.
		
		By elementarity, \(j_W(\xi)\) is the least ordinal \(\alpha\) with \(N =
		H^{N}(j_W(j_U)[M_W]\cup \{\alpha\})\). On the other hand, since
		\((j_U(j_W),j_U\restriction M_W)\) is a minimal comparison of
		\((j_U,j_W)\) (\cref{InternalComparison}),  \(N = H^N(j_U[M_W]\cup
		\{j_U(j_W)(\xi)\})\) (\cref{MinimalUltrapowers}). Since \(j_U
		\restriction M_W = j_W(j_U)\restriction M_W\), this yields \[N =
		H^N(j_W(j_U)[M_W]\cup \{j_U(j_W)(\xi)\})\] By the minimality of
		\(j_W(\xi)\), \(j_W(\xi) \leq j_U(j_W)(\xi)\), as desired.
		
		Thus \(j_U(j_W)\) and \(j_W\restriction M_U\) coincide on \(j_U[V]\cup
		\{\xi\}\). Since \(M_U = H^{M_U}(j_U[V]\cup \{\xi\})\), it follows that
		\( j_U(j_W) = j_W\restriction M_U\), as desired.
	\end{proof}
\end{prp}

\begin{proof}[Proof of \cref{KunenCommute}] It is trivial to see that \(j_W(j_U)
	= j_U\restriction M_W\). Hence by \cref{WeakCommute}, \(U\) and \(W\)
	commute.
\end{proof}

Under UA, the only counterexamples to the strictness of the internal relation
are commuting ultrafilters:

\begin{thm}[UA]\label{Commute}\index{Commuting ultrafilters!vs. the internal relation}
	Suppose \(U\I W\) and \(W\I U\). Then \(U\) and \(W\) commute.
	\begin{proof}
		Since \(U\I W\), \(\tr U W = j_U(W)\). Since \(W \I U\), \(\tr U W =
		s_U(W)\). Therefore \(j_U(W) = s_U(W)\). It follows that \(j_U(j_W) =
		j^{M_U}_{j_U(W)} = j^{M_U}_{s_U(W)} = j_W\restriction M_U\) by
		\cref{PushUlt}. Similarly, \(j_W(j_U) = j_U\restriction M_W\). In other
		words, \(U\) and \(W\) commute, as desired.
	\end{proof}
\end{thm}

This raises an interesting technical question:
\begin{qst}[ZFC] Suppose \(U\) and \(W\) are countably complete ultrafilters
	such that \(U\I W\) and \(W\I U\). Do \(U\) and \(W\) commute?
\end{qst}

\cref{InternalWF} gives some information regarding this question:

\begin{prp}
	If \(U\I W\) and \(W\I U\), then \(U\) and \(W\) have no common points of
	discontinuity.\qed
\end{prp}

The supercompactness analysis of \cref{SCChapter1} occasionally requires a
partial converse to \cref{KunenCommute}: the only way certain nice pairs of
ultrafilters can commute is if one lies below the completeness of the other.

\begin{defn}\label{LambdaInternalDef}\index{Ultrafilter!\(\lambda\)-internal}
	Suppose \(\lambda\) is a cardinal. A countably complete ultrafilter \(W\) is
	{\it \(\lambda\)-internal} if \(U\I W\) for all \(U\) such that \(\lambda_U
	< \lambda\).
\end{defn}

\begin{prp}\label{IUniformPrp}\index{Kunen's commuting ultrapowers lemma!converse}
	Suppose \(U\) and \(W\) are countably complete hereditarily uniform
	ultrafilters such that \(U\) is \(\lambda_U\)-internal and \(W\) is
	\(\lambda_W\)-internal. Let \(\kappa_U = \textsc{crt}(j_U)\) and \(\kappa_W
	= \textsc{crt}(j_W)\). Then the following are equivalent:
	\begin{enumerate}[(1)]
		\item \(U\) and \(W\) commute.
		\item Either \(U\in V_{\kappa_W}\) or \(W\in V_{\kappa_U}\).
	\end{enumerate}
\end{prp}

One can also state \cref{IUniformPrp} avoiding the notion of hereditary
uniformity: if \(U\) is \(\lambda_U\)-internal and \(W\) is
\(\lambda_W\)-internal, then \(U\) and \(W\) commute if and only if \(\lambda_U
< \kappa_W\) or \(\lambda_W < \kappa_U\).

The proof of \cref{IUniformPrp} requires a number of lemmas. The first allows us
to approximate an arbitrary ultrapower embedding by a small ultrafilter:

\begin{lma}\label{Exponential}
	Suppose \(j : V\to M\) is an ultrapower embedding. Then for any cardinal
	\(\lambda\), there is a countably complete ultrafilter \(D\) with
	\(\lambda_D \leq 2^{\lambda}\) such that there is an elementary embedding
	\(k : M_D\to M\) with \(k\circ j_D = j\) and \(\textsc{crt}(k) > \lambda\).
	\begin{proof}
		Suppose \(\gamma\) is an ordinal. We will find an ultrafilter \(D\) on
		\(\gamma^\gamma\) such that there is an elementary embedding \(k :
		M_D\to M\) with \(k\circ j_D = j\) and \(\textsc{crt}(k) \geq \gamma\).
		Taking \(\gamma = \lambda+1\) proves the lemma.
		
		Fix \(a\in M\) such that \(M = H^M(j[V]\cup \{a\})\) and \(X\) such that
		\(a\in j(X)\). Fix functions \(\langle f_\alpha : \alpha <
		\gamma\rangle\) on \(X\) such that \(\alpha = j(f_\alpha)(a)\). Define a
		function \(g: X\to \gamma^\gamma\) by letting \(g(x)\) be the function
		with \(g(x)(\alpha) = f_\alpha(x)\) for all \(\alpha < \gamma\). 
		
		Let \(D\) be the ultrafilter on \(\gamma^\gamma\) derived from \(j\)
		using \(j(g)(a)\). Let \(k : M_D\to M\) be the factor embedding such
		that \(k\circ j_D = j\) and \(k(\id_D) = j(g)(a)\).
		
		We claim that \(\textsc{crt}(k) \geq \gamma\). It suffices to show that
		\(\gamma\subseteq k[M_D] = H^M(j[V]\cup \{j(g)(a)\})\). Fix \(\alpha <
		\gamma\). Then \[\alpha = j(f_\alpha)(a) = j(f)_{j(\alpha)}(a) =
		j(g)(a)(j(\alpha))\] Thus \(\alpha\) is definable in \(M\) from
		\(j(g)(a)\) and \(j(\alpha)\). Thus \(\alpha\in H^M(j[V]\cup
		\{j(g)(a)\})\), as desired.
	\end{proof}
\end{lma}
The coarseness of the bound \(2^{\lambda}\) actually causes a number of problems
down the line. An argument due to Silver (which appears as \cref{Silver})
provides a major improvement in a special case, and is instrumental in our
analysis of the linearity of the Mitchell order on normal fine ultrafilters
under UA without GCH assumptions. Further improvements could potentially solve
the problems concerning so-called isolated cardinals discussed in
\cref{IsolationSection}.

Using \cref{Exponential}, we prove the following lemma, which can be seen as a
version of the Kunen Inconsistency Theorem (\cref{KunenInconsistency}) that
replaces the strength requirement of that theorem with a requirement involving
the internal relation:

\begin{lma}\label{InternalComplete}
	Suppose \(U\) is a countably complete ultrafilter and \(\kappa\) is a strong
	limit cardinal. Then the following are equivalent:
	\begin{enumerate}[(1)]
		\item \(U\) is \(\kappa\)-internal and \(\sup j_U[\kappa]\subseteq
		\kappa\).
		\item \(U\) is \(\kappa\)-complete.
	\end{enumerate}
	\begin{proof}
		{\it (1) implies (2).} Let \(j : V\to M\) be the ultrapower of the
		universe by \(U\). We first show that \(j\) is
		\({<}\kappa\)-supercompact. Fix \(\gamma < \kappa\), and we will prove
		that \(j\restriction \gamma\in M\). Let \(\lambda = j(\gamma)\), so
		\(\lambda < \kappa\) by the assumption that \(j[\kappa]\subseteq
		\kappa\). By \cref{Exponential}, one can find a countably complete
		ultrafilter \(D\) with \(\lambda_D \leq 2^\lambda < \kappa\) and an
		elementary embedding \(k : M_D\to M\) with \(k\circ j_D = j\) and
		\(\textsc{crt}(k) > \lambda = j(\gamma)\). In particular
		\(j_D\restriction \gamma = j\restriction \gamma\). Moreover since
		\(\lambda_D < \kappa\), \(D\I U\). Therefore \(j\restriction \gamma =
		j_D\restriction \gamma\in M\), as desired.
		
		Now \(j\) is \({<}\kappa\)-supercompact and \(j[\kappa]\subseteq
		\kappa\). Since \(j\) is an ultrapower embedding, if \(\kappa\) is
		singular, then \(j\) is \(\kappa\)-supercompact. Therefore the Kunen
		Inconsistency Theorem (\cref{KunenInconsistency} or
		\cref{KunenInconsistency0}) implies \(\textsc{crt}(j)\geq \kappa\), so
		\(U\) is \(\kappa\)-complete.
		
		{\it (2) implies (1).} Trivial.
	\end{proof}
\end{lma}

\begin{lma}\label{NonCommute}
	Suppose \(U\) and \(W\) are nonprincipal countably complete ultrafilters.
	Let \(\kappa_U = \textsc{crt}(j_U)\) and \(\kappa_W  = \textsc{crt}(j_W)\).
	Assume \(U\) is \(\kappa_W\)-internal and \(W\) is \(\kappa_U\)-internal.
	Then either \(j_U(\kappa_W) > \kappa_W\) or \(j_W(\kappa_U) > \kappa_U\).
	\begin{proof}
		Assume towards a contradiction that \(j_U(\kappa_W)= \kappa_W\) and
		\(j_W(\kappa_U)= \kappa_U\). Since \(U\) is \(\kappa_W\)-internal and
		\(j_U[\kappa_W]\subseteq \kappa_W\), \(U\) is \(\kappa_W\)-complete.
		Therefore \(\kappa_U\geq \kappa_W\). By symmetry, \(\kappa_W\geq
		\kappa_U\). Thus \(\kappa_U = \kappa_W\). This contradicts that
		\(j_U(\kappa_W) = \kappa_W\) while \(j_U(\kappa_U) > \kappa_U\) by the
		definition of a critical point.
	\end{proof}
\end{lma}

We can finally prove \cref{IUniformPrp}:

\begin{proof}[Proof of \cref{IUniformPrp}] {\it (1) implies (2):} Since \(U\)
	and \(W\) commute, \(j_U(\kappa_W) = \kappa_W\) and \(j_W(\kappa_U) =
	\kappa_U\). By \cref{NonCommute}, either \(U\) is not \(\kappa_W\)-internal
	or \(W\) is not \(\kappa_U\)-internal. Therefore either \(\lambda_U <
	\kappa_W\) or \(\lambda_W < \kappa_U\). 
	
	Assume first that \(\lambda_U < \kappa_W\). Since \(U\) is hereditarily
	uniform, the underlying set of \(U\) has hereditary cardinality
	\(\lambda_U\), and hence \(U\in V_{\kappa_W}\) since \(\kappa_W\) is
	inaccessible.
	
	If instead \(\lambda_W < \kappa_U\), then \(W\in V_{\kappa_U}\) by a similar
	argument.
	
	{\it (2) implies (1):} Immediate from \cref{KunenCommute}.
\end{proof}
\subsection{\(j\) on the ordinals}
In this section, we briefly survey some results that tie the structure of the
internal relation under UA to the behavior of elementary embeddings on the
ordinals. We only sketch most of the proofs since the material is a bit of a
detour from the main line of this dissertation.

Recall the notion of the rank of a pointed ultrapower in the Ketonen order
(\cref{ERank}): if \(\lambda\) is a cardinal, then \(\mathcal P_\lambda\)
denotes the collection of pointed ultrapowers \((M,\xi)\) such that \(M\) is the
ultrapower by an ultrafilter \(U\) with \(\lambda_U < \lambda\) and \(\xi\) is
an ordinal; \(o_\lambda(M,\xi)\) denotes the rank of \((M,\xi)\) in the
prewellorder \((\mathcal P_\lambda,\sE)\). 

\begin{lma}[UA]\index{Ketonen order!on pointed ultrapowers!rank (\(o_\lambda(M)\))}\index{Internal relation!fixed points}
	Suppose \(\lambda\) is a regular cardinal, \(j : V\to M\) is an ultrapower
	of width less than \(\lambda\), \(i : M \to N\) is an ultrapower embedding
	such that \(i\circ j\) has width less than \(\lambda\), and \(\xi\) is an
	ordinal such that \(M = H^M(j[V]\cup \{\xi\})\). Then the following are
	equivalent:
	\begin{enumerate}[(1)]
		\item \(i\) is an internal ultrapower embedding.
		\item \((M,\xi) =_S (N,i(\xi))\).
		\item \(i(o_\lambda(M,\xi)) = o_\lambda(N,i(\xi))\).
	\end{enumerate}
	\begin{proof}
		The equivalence of (1) and (2) is an immediate consequence of
		\cref{RFSEquiv}. The equivalence of (2) and (3) follows from the fact
		that \((M,\xi), (N,i(\xi))\in \mathcal P_\lambda\) (and does not require
		the assumption that \(M = H^M(j[V]\cup \{\xi\})\)). 
	\end{proof}
\end{lma}

A surprising consequence of this is that under UA, the ultrafilters to which a
countably complete ultrafilter \(U\) is internal are determined solely by the
class of fixed points of \(j_U\). Recall here that if \(W\) is a countably
complete ultrafilter on an ordinal with \(\lambda_W < \lambda\), then
\(o_\lambda(W) = o_\lambda(M_W,\id_W)\).
\begin{thm}[UA]\label{IFix}
	Suppose \(U\) is a countably complete ultrafilter and \(W\) is a countably
	complete ultrafilter on an ordinal. Then the following are equivalent:
	\begin{enumerate}[(1)]
		\item \(U\I W\)
		\item \(j_U(o_\lambda(W)) = o_\lambda(W)\) for all regular cardinals
		\(\lambda > \lambda_U,\lambda_W\).
		\item \(j_U(o_\lambda(W)) = o_\lambda(W)\) for some regular cardinal
		\(\lambda > \lambda_U,\lambda_W\).\qed
	\end{enumerate}
\end{thm}

In particular, if \(j_U\restriction \text{Ord} = j_{U'}\restriction
\text{Ord}\), then \(U\) and \(U'\) are internal to exactly the same
ultrafilters. The following observation shows that this is not vacuous, in that
there are many nonisomorphic ultrafilters that have the same action on the
ordinals:

\begin{thm}\label{RanFix1}
	Suppose \(U\) is a countably complete ultrafilter and \(Z\) is a countably
	complete ultrafilter of \(M_U\) such that \(j_U(Z)\) and \(j_U(j_U)(Z)\)
	commute in \(j_U(M_U)\). Then \(j^{M_U}_Z\) fixes every ordinal in the range
	of \(j_U\).
\end{thm}

The proof uses a lemma due to Kunen (which should be compared with
\cref{InternalWF}):
\begin{lma}[Kunen]\label{KunenFinite}
	Suppose \(\alpha\) is an ordinal and \(S\) is a set of pairwise commuting
	countably complete ultrafilters such that \(j_W(\alpha) > \alpha\) for all
	\(W\in S\). Then \(S\) is finite.
	\begin{proof}
		Suppose towards a contradiction that \(\alpha\) is the least ordinal
		such that there is an infinite set \(S\) of pairwise commuting countably
		complete ultrafilters such that for all \(W\in S\), \(j_W(\alpha) >
		\alpha\). Fix \(W\in S\). Let \(T\) be a countably infinite subset of
		\(S\) such that \(W\notin T\). Then in \(M_W\), \(j_W(T)\) is an
		infinite set of pairwise commuting countably complete ultrafilters. For
		any \(U\in T\), since \(U\) and \(W\) commute, \[j_W(j_U)(\alpha) =
		j_U(\alpha) > \alpha\] Thus for any \(Z\in j_W(T)\), \(j_Z^{M_W}(\alpha)
		> \alpha\). Here we use that \(T\) is countable so \(j_W(T) = j_W[T]\).
		
		In particular, in \(M_W\) there is an infinite set of pairwise commuting
		ultrafilters all of whose associated embeddings move \(\alpha\). But by
		the elementarity of \(j_W\) and the definition of \(\alpha\), \(M_W\)
		satisfies that \(j_W(\alpha)\) is the least ordinal \(\xi\) such that
		there is an infinite set of pairwise commuting ultrafilters all of whose
		associated embeddings move \(\xi\). Since \(j_W(\alpha) > \alpha\), this
		is a contradiction.
	\end{proof}
\end{lma}

\begin{proof}[Proof of \cref{RanFix1}] Let \(X\) be the underlying set of \(U\).
		Choose countably complete ultrafilters \(\langle Z_x : x\in X\rangle\)
		such that \[Z = [\langle Z_x : x\in X\rangle]_U\] The assumption that
		\(j_U(Z)\) and \(j_U(j_U)(Z)\) commute can be reformulated as follows:
		\begin{equation}\label{AlmostAlwaysCommute}\{(x,y) \in X\times X :
		Z_x\text{ and }Z_y\text{ commute}\}\in U\times U\end{equation}
		
		Fix an ordinal \(\alpha\). We must show that \(j_Z^{M_U}(j_U(\alpha)) =
		j_U(\alpha)\). By \L o\'s's Theorem, it suffices to show that for almost
		all \(x\in X\), \(\{x\in X : j_{Z_x}(\alpha) = \alpha\}\in U\). Let \[A
		= \{x\in X : j_{Z_x}(\alpha) > \alpha\}\] and assume towards a
		contradiction that \(A\in U\). 
		
		Fix \(i <\omega\). If \(n < m < i\), the function \(f : X^i\to X^2\)
		defined by \(f(x_0,\dots,x_{i-1})= (x_n,x_m)\) pushes \(U^i\) forward
		onto \(U^2\), so by \cref{AlmostAlwaysCommute}, \[\{(x_0,\dots,x_{i-1})
		: Z_{x_n}\text{ and }Z_{x_m}\text{ commute}\}\in U^i\] Thus \(C_i\in
		U^i\) where \[C_i = \{(x_0,\dots,x_{i-1})\in X^i : \textnormal{for all
		\(n,m < i\), \(Z_{x_n}\) and \(Z_{x_m}\) commute}\}\] Letting \(B_i =
		A^i\cap C_i\), since \(A\in U\) by assumption, we have \(B_i\in U^i\). 
		
		Since \(\langle U^i : i < \omega\rangle\) is a countably complete tower
		of ultrafilters with \(B_i\in U^i\) for all \(i < \omega\), there is a
		sequence \(\langle x_n : n < \omega\rangle\) such that
		\((x_0,\dots,x_{i-1})\in B_i\) for all \(i < \omega\). Now \(\{ Z_{x_i}
		: i < \omega\}\) is an infinite set of pairwise commuting ultrafilters
		whose associated embeddings all move \(\alpha\), contradicting
		\cref{KunenFinite}. Thus our assumption was false, so in fact
		\(j^{M_U}_Z(j_U(\alpha)) = j_U(\alpha)\), as desired.
\end{proof}

A corollary of the proof of \cref{RanFix2} is the following more general fact
about extenders. 
\begin{defn}A pair of extenders \(E\) and \(F\) {\it commute} if \(j_E(j_F) = j_F\restriction M_E\) and \(j_F(j_E) = j_E\restriction M_F\).\end{defn}

\begin{cor}\label{RanFix2}
	Suppose \(E\) is an extender. Suppose \(F\) is an \(M_E\)-extender such that
	\(j_E(F)\) and \(j_E(j_E)(F)\) commute in \(j_E(M_E)\). Then \(j_F^{M_E}\)
	fixes every ordinal in the range of \(j_E\).
	\begin{proof}[Sketch] The first step is to reduce to the case that \(E\) is
		an ultrafilter. Let \(U\) be the ultrafilter derived from \(j_E\) using
		\(F\). Let \(\bar F = \id_U\). A simple diagram chase shows that \(U^2\)
		is the ultrafilter derived from \(j_E(j_E)\circ j_E\) using
		\((j_E(j_E)(F),j_E(F))\). As a consequence of this, \(j_U(\bar F)\) and
		\(j_U(j_U)(\bar F)\) commute in \(j_U(M_U)\). It suffices to show that
		\(j_{\bar F}^{M_U}\) fixes every ordinal in the range of \(j_U\), since
		then by the elementarity of the factor embedding \(k : M_U\to M_E\),
		\(j_F^{M_E}\) fixes every ordinal in the range of \(j_E\). 
		
		Now one generalizes the proof of \cref{RanFix1} to the case where \(Z\)
		is an extender \(F\) rather than an ultrafilter. This presents no real
		difficulties.
	\end{proof}
\end{cor}

Another interesting corollary regards the relationship between the Mitchell
order and pointwise domination of elementary embeddings on the ordinals.

\begin{thm}
	Suppose \(F\mo E\) are extenders. Let \(\kappa = \textsc{crt}(F)\) and
	\(\iota = \textsc{width}(F)\). Assume that the following hold:
	\begin{itemize}
		\item \((M_F)^{<\kappa}\subseteq M_F\) and \((M_E)^{<\iota} \subseteq
		M_E\).
		\item \(j_E(j_E)(F) \in j_E(V_\kappa)\).
	\end{itemize}
	Then for all ordinals \(\alpha\), \(j_F(\alpha) \leq j_E(\alpha)\) with
	equality if and only if \(j_E(\alpha) = \alpha\).
	\begin{proof}[Sketch] Since \((M_E)^{<\iota} \subseteq M_E\), we have
		\(j_F^{M_E} = j_F\restriction M_E\). Since \(j_E(j_E)(F) \in
		j_E(V_\kappa)\) and \(M_F^{<\kappa}\subseteq M_F\), \(j_E(j_E)(F)\) and
		\(j_E(F)\) commute in \(j_E(M_E)\) by a generalization of the proof of
		\cref{KunenCommute}. Therefore applying \cref{RanFix2} yields that
		\(j_F(j_E(\alpha)) = j_E(\alpha)\) for all ordinals \(\alpha\), which
		easily implies the conclusion of the theorem.
	\end{proof}
\end{thm}

The requirement that \(j_E(j_E)(F)\in j_E(V_\kappa)\) may seem ad hoc, but in
fact it is necessary. For example, suppose \(\kappa < \lambda\) are cardinals,
\(F\) is a \((\kappa,\lambda)\)-extender that witnesses that \(\kappa\) is
\(\lambda\)-strong, and \(U\) is a normal ultrafilter on \(\lambda\). Trivially
\(F\mo U\), \((M_U)^{<\kappa}\subseteq M_U\), \((M_F)^{<\kappa}\subseteq M_F\),
and yet \(j_F(\kappa) > j_U(\kappa)\).

\cref{IFix} above implies that under UA, the question of whether \(U\I W\)
depends only on \(j_U\restriction \text{Ord}\) and \(M_W\). The following
theorem explains why:

\begin{thm}[UA] Suppose \(U\) and \(W\) are countably complete ultrafilters such
	that \(j_U\restriction \lambda\in M_W\) for all cardinals \(\lambda\). Then
	\(U\I W\).
	\begin{proof}[Sketch] For the proof, we need the following weak consequence
		of the analysis of directed systems of internal ultrapower embeddings
		(\cref{InftyAbs}): there is an inner model \(N\) such that the following
		hold:\footnote{The inner model \(N\) can be taken to be the direct limit
		\(M_\lambda\) of all ultrapowers of width less than \(\lambda\) for some
		sufficiently large regular cardinal \(\lambda\). Then by
		\cref{InftyAbs}, \(N = j_U(N) = j_W(N)\subseteq M_W\), as desired. The
		embedding \(k\) is obtained by setting \(k=j_W(j_\lambda) =
		j_{M_W,\lambda}\) (by \cref{InftyAbs} again). Then \(k\) is an
		elementary embedding from \(j_W(V) = M_W\) to \(j_W(N) = N\), as
		desired. Note that we cannot assume that \(k\) is an ultrapower
		embedding.}
		\begin{itemize}
			\item \(j_U(N) \subseteq M_W\).
			\item There is an elementary embedding \(k : M_W\to N\) that is
			amenable to \(M_W\).
		\end{itemize}

		We claim that the embedding \(j_U\restriction N\) is amenable to
		\(M_W\). To see this, suppose \(X\in N\) is a transitive set. We will
		show \(j_U\restriction X\in M_W\). Since \(X\) is transitive, it
		suffices to show that \(j_U[X]\in M_W\). Let \(\lambda\) be a cardinal
		and \(p:\lambda\to X\) be a surjection with \(p\in N\). Then \(j_U(p)\in
		j_U(N) \subseteq M_W\), so \(j_U[X] = j_U(p)[j_U[\lambda]] \in M_W\), as
		desired.
		
		Let \(U_* = (k\circ j_W)_*(U)\cap N\). In other words, by the basic
		theory of the Rudin-Keisler order (\cref{PushDerived}), \(U_*\) is the
		\(N\)-ultrafilter derived from \(j_U\restriction N\) using \(j_U(k\circ
		j_W)(\id_U)\). Since \(j_U\restriction N\) is amenable to \(M_W\), it
		follows that \(U_*\in M_W\). Since \(k\) is amenable to \(M_W\),
		\(k^{-1}[U_*]\in M_W\), but \[k^{-1}[U_*] = k^{-1}[(k\circ j_W)_*(U)\cap
		N] = (j_W)_*(U)\cap M_W = s_W(U)\] Thus \(s_W(U)\in M_W\), so \(U \I W\)
		by \cref{PushUlt}.
	\end{proof}
\end{thm}

\chapter{Ordinal Definability and Cardinal Arithmetic under UA}\label{GCHChapter}
\section{Introduction}
\subsection{The universe above a supercompact cardinal}
In this short chapter, we exposit two results that show that something
remarkable happens when UA is combined with very large cardinal hypotheses:
instead of simply proving structural results for countably complete
ultrafilters, the axiom now begins to resolve major questions independent from
the usual axioms of set theory.

Let us describe the main results of this section. Since UA is preserved by
forcing to add a Cohen real, UA does not imply \(V = \textnormal{HOD}\), no
matter what large cardinals one assumes in addition to UA. But it turns out it
is possible to prove that forcing is the only obstruction:
\begin{repthm}{HODGen}[UA] Assume there is a supercompact cardinal. Then \(V\)
is a generic extension of \textnormal{HOD}.
\end{repthm}

Similarly, UA is preserving by forcing to change the value of the continuum, so
UA does not imply the Continuum Hypothesis. But UA for sufficiently large
cardinals \(\lambda\), UA implies \(2^\lambda = \lambda^+\):
\begin{thm}[UA]\label{GCH}\index{Generalized Continuum Hypothesis}
Assume \(\kappa\) is supercompact. Then for all cardinals \(\lambda\geq
\kappa\), \(2^\lambda = \lambda^+\).
\end{thm}
It seems that above a supercompact cardinal, UA imposes incredible structure on
the universe of sets. This is explored further in \cref{SCChapter1} and
\cref{SCChapter2}.
\subsection{Outline of \cref{GCHChapter}}
We now outline the rest of the chapter.\\

\noindent {\sc\cref{ODSection}.} We prove the results on ordinal definability
under UA and large cardinals. This is quite straightforward, but many open
questions remain. For example, we prove that if \(\kappa\) is supercompact and
UA holds, then \(V\) is a generic extension of \(\textnormal{HOD}\). How small
is the forcing? The best upper bound we know is \(\kappa^{++}\), which comes
from \cref{VopenkaSection} below.\\

\noindent {\sc\cref{GCHSection}.} We prove the results on GCH under UA and large
cardinals. We begin by discussing related results in ZFC, especially Solovay's
theorem on SCH above a supercompact cardinal. In \cref{MoreMOSection}, we prove
a result regarding the Mitchell order and supercompactness that shows that under
UA, if \(D\) and \(U\) are ultrafilters with \(\lambda_D\) below the
supercompactness of \(U\), then \(D\mo U\). This is immediate given GCH, but
proving this using UA alone is a little bit subtle. In \cref{GCHProofSection},
we use this result to conclude that GCH holds above a supercompact.
\section{Ordinal definability}\label{ODSection}
\cref{HODGen} is quite easy given what we have shown so far, ultimately relying
on the following simple fact:
\begin{prp}[UA]\label{ODUF}
Every countably complete ultrafilter on an ordinal is ordinal definable.
\begin{proof}
Suppose \(\delta\) is an ordinal. Then the set \(\mathscr B(\delta)\) of all
countably complete ultrafilters on \(\delta\) is wellordered by the Ketonen
order. Thus every element of \(\mathscr B(\delta)\) is ordinal definable from
its rank in the Ketonen order.
\end{proof}
\end{prp}

\begin{cor}[UA]\label{HODAmenable}
For any set of ordinals \(X\) and any ultrapower embedding \(j : V\to M\):
\begin{enumerate}[(1)]
\item \(j(\textnormal{OD}_X)\subseteq \textnormal{OD}_X\).
\item \(j(\textnormal{HOD}_X)\subseteq \textnormal{HOD}_X\).
\item For any \(Y\in \textnormal{HOD}_X\), \(j\restriction Y\in
\textnormal{HOD}_X\).
\end{enumerate} 
\begin{proof}
We first prove (1). We have \(j(\text{OD}_X) = \text{OD}_{j(X)}^{M}\). Fix a
countably complete ultrafilter \(U\) on an ordinal such that \(j = j_U\). Then
since \(M\) is definable from \(U\) and \(U\in \text{OD}\) by \cref{ODUF},
\(\text{OD}_{j(X)}^M\subseteq \text{OD}_{j(X)}\). Moreover \(j(X) = j_U(X)\) is
definable from \(X\) and \(U\), so \(j(X)\in \text{OD}_X\). Hence
\(\text{OD}_{j(X)}^M\subseteq \text{OD}_{j(X)} \subseteq \text{OD}_X\).

For (2), note that \(j(\textnormal{HOD}_X)\) is the class of sets that are
hereditarily \(j(\text{OD}_X)\), and this is contained in the class of sets that
are hereditarily \(\text{OD}_X\) by (1).

For (3), clearly \(j\restriction Y\in \text{OD}_{Y}\subseteq \text{OD}_X\). But
moreover by (2), \(j\restriction Y\subseteq \text{HOD}_X\). Therefore
\(j\restriction Y\in \textnormal{HOD}_X\).
\end{proof}
\end{cor}

The following lemma should be compared with the theorem of Shelah that if
\(\lambda\) is a singular strong limit cardinal of uncountable cofinality, then
for any \(X\) such that \(P(\alpha)\subseteq \text{HOD}_X\) for all \(\alpha <
\lambda\), in fact  \(P(\lambda)\subseteq \text{HOD}_X\).

\begin{lma}[UA]\label{PowerPropagation}
Suppose \(\kappa\) is \(\lambda\)-supercompact and \(X\subseteq \kappa\) is such
that \(V_\kappa\subseteq \textnormal{HOD}_X\). Then \(P(\lambda)\subseteq
\textnormal{HOD}_X\).
\begin{proof}
Fix a \(\lambda\)-supercompact ultrapower embedding \(j :V \to M\) such that
\(\textsc{crt}(j) = \kappa\) and \(j(\kappa) > \lambda\). Then
\[P(\lambda)\subseteq j(V_\kappa)\subseteq j(\text{HOD}_X)\subseteq
\textnormal{HOD}_X\] The final inclusion follows from \cref{HODAmenable}.
\end{proof}
\end{lma}

\begin{thm}[UA]\label{HODofSet}
Suppose \(\kappa\) is supercompact. Then \(V = \textnormal{HOD}_X\) for some
\(X\subseteq \kappa\).
\begin{proof}
Fix \(X\subseteq \kappa\) such that \(V_\kappa\subseteq
\text{HOD}_X\).\footnote{To obtain such a set \(X\), let \(E\) be a binary
relation on \(\kappa\) such that \((V_\kappa,\in) \cong (\kappa,E)\) using the
fact that \(|V_\kappa| = \kappa\). Code \(E\) as a subset of \(\kappa\) using a
pairing function \(\kappa\to \kappa\times \kappa\).} Since \(\kappa\) is
supercompact, \cref{PowerPropagation} implies that for all \(\lambda\geq
\kappa\), \(P(\lambda)\subseteq \textnormal{HOD}_X\), and therefore \(V =
\text{HOD}_X\).
\end{proof}
\end{thm}

To connect this to generic extensions of \(\text{HOD}\), we use Vop\v enka's
Theorem.

\begin{defn}
Suppose \(X\) is a set such that and \(X\cup \{X\}\subseteq \textnormal{OD}\).
The {\it \textnormal{OD}-cardinality of \(X\)}, denoted \(|X|^\textnormal{OD}\),
is the least ordinal \(\lambda\) such that there is an \(\text{OD}\) bijection
between \(\lambda\) and \(X\).
\end{defn}
The OD cardinality of \(X\) is defined for all \(X\) with \(X\cup \{X\}\subseteq
\text{OD}\). It is always a \(\text{HOD}\)-cardinal. In fact OD cardinality
satisfies all the usual properties of cardinality; for example,
\(|X|^\text{OD}\) is the least ordinal that ordinal definably surjects onto
\(X\) and the least ordinal into which \(X\) ordinal definably injects.

\begin{defn}\index{Vopenka algebra}
Suppose \(\kappa\) is an ordinal. Let \(\mathbb A_\kappa\) be the Boolean
algebra \(P(P(\kappa))\cap \textnormal{OD}\) and let \(\lambda = |\mathbb
A_\kappa|^\text{OD}\). Fix an \(\textnormal{OD}\) bijection \(\pi : \lambda\to
\mathbb A\). Then \(\mathbb V_\kappa\) is the Boolean algebra on \(\lambda\)
given by pulling back the operations on \(\mathbb A\) under \(\pi\).
\end{defn}

Note that \(\mathbb V_\kappa\in \textnormal{HOD}\). The Boolean algebra
\(\mathbb V_\kappa\) is called the {\it Vop\v enka algebra} at \(\kappa\).

\begin{thm}[Vop\v enka]\label{Vopenka}
If \(\kappa\) is an ordinal, then \(\mathbb V_\kappa\) is a complete Boolean
algebra and for any \(X\subseteq \kappa\), there is a
\(\textnormal{HOD}\)-generic ultrafilter \(G\subseteq \mathbb V_\kappa\) such
that \(\textnormal{HOD}_X\subseteq \textnormal{HOD}[G]\).\qed
\end{thm}

This yields a proof of our main theorem on HOD:
\begin{thm}[UA]\label{HODGen}
	Assume there is a supercompact cardinal. Then \(V\) is a generic extension
	of \textnormal{HOD}.
\end{thm}
\begin{proof}
Let \(\kappa\) be the least supercompact cardinal. By \cref{HODofSet}, \(V =
\text{HOD}_X\) for some \(X\subseteq \kappa\), so by \cref{Vopenka}, \(V =
\textnormal{HOD}[G]\) for some generic \(G\subseteq \mathbb V_\kappa\).
\end{proof}

\begin{qst}[UA] Let \(\kappa\) be the least supercompact cardinal.
\begin{itemize}
\item Is \(V = \text{HOD}[X]\) for some \(X\subseteq \kappa\)? 
\item Is \(V = \text{HOD}[G]\) for \(G\subseteq \kappa\) generic for a partial
order \(\mathbb P \in \textnormal{HOD}\) such that \(|\mathbb P| \leq \kappa\)?
What about a \(\kappa\)-cc Boolean algebra?
\item Is \(V = \text{HOD}_{V_\kappa}\)?
\end{itemize}
\end{qst}

Assuming UA, one can actually calculate the cardinality of \(\mathbb V_\kappa\)
precisely. For example, in the next section, we will obtain:
\begin{repthm}{VopSizeThm}[UA] If \(\kappa\) is \(\kappa^{++}\)-supercompact
then \(|\mathbb V_\kappa|^\textnormal{HOD} = \kappa^{++}\). 
\end{repthm}

Thus if \(\kappa\) is supercompact, then \(V = \textnormal{HOD}[A]\) for  some
\(A\subseteq\kappa^{++}\). As an immediate consequence, we have that
\(\textnormal{HOD}\) is very close to \(V\):

\begin{cor}[UA]\label{HODClose}
Let \(\kappa\) be the least supercompact cardinal. Then for all cardinals
\(\lambda \geq \kappa^{++}\):
\begin{enumerate}[(1)]
\item \(\lambda^{+\textnormal{HOD}} = \lambda^+\).
\item \((2^\lambda)^{\textnormal{HOD}} = 2^\lambda\).
\end{enumerate}
Moreover if \(\delta > \kappa^{++}\) is regular, then \(\textnormal{HOD}\) is
correct about stationary subsets of \(\delta\).\qed
\end{cor}

Of course by the Levy-Solovay Theorem \cite{LevySolovay}, HOD is also close to
\(V\) in the sense that it absorbs large cardinals above \(\kappa\). The
structure of HOD at \(\kappa\) itself becomes a key question:
\begin{qst}[UA] Assume \(\kappa\) is supercompact. Is \(\kappa^{+\text{HOD}} =
\kappa^+\)?
\end{qst}

Another question in this vein is whether \(\kappa\) is supercompact in
\(\textnormal{HOD}\). Here the answer turns out to be yes:
\begin{defn}
If \(N\) is an inner model and \(S\) is a set, we say \(S\) is {\it amenable to
\(N\)} if \(S\cap N\in N\).
\end{defn}

\begin{defn}\index{Weak extender model}
Suppose \(\kappa\) is supercompact. An inner model \(N\) is a {\it weak extender
model at \(\kappa\)} if for all ordinals \(\lambda\geq \kappa\), there is a
normal fine \(\kappa\)-complete ultrafilter on \(P_\kappa(\lambda)\) that
concentrates on \(N\) and is amenable to \(N\).
\end{defn}

\begin{lma}\label{UnboundedWEM}
Suppose \(N\) is an inner model and \(\kappa\) is supercompact. Then the
following are equivalent:
\begin{enumerate}[(1)]
\item \(N\) is a weak extender model at \(\kappa\).
\item For arbitrarily large \(\delta\geq \kappa\), there is a normal fine
\(\kappa\)-complete ultrafilter on \(P_\kappa(\delta)\) that concentrates on
\(N\) and is amenable to \(N\).
\end{enumerate}
\begin{proof}
{\it (1) implies (2):} Trivial.
	
{\it (2) implies (1):} Fix \(\lambda \geq \kappa\). We will show that there is a
normal fine \(\kappa\)-complete ultrafilter on \(P_\kappa(\lambda)\) that
concentrates on \(N\) and is amenable to \(N\). By (2), there is some \(\delta
\geq\lambda\) such that there is a normal fine \(\kappa\)-complete ultrafilter
\(\mathcal U\) on \(P_\kappa(\delta)\) that concentrates on \(N\) and is
amenable to \(N\). Let \(\mathcal W = f_*(\mathcal U)\) where \(f :
P_\kappa(\delta)\to P_\kappa(\lambda)\) is defined by \(f(\sigma) = \sigma\cap
\lambda\). Easily \(\mathcal W\) is a normal fine ultrafilter. Moreover
\(f^{-1}[P_\kappa(\lambda)\cap M] = P_\kappa(\delta)\cap M\in \mathcal U\), so
\(P_\kappa(\lambda)\cap M\in \mathcal W\). Thus \(\mathcal W\) concentrates on
\(M\). Finally, letting \(g = f\restriction M\), clearly \(g\in M\) and hence
\(\mathcal W\cap M = f_*(\mathcal U)\cap  M = g_*(\mathcal U\cap M) \in M\)
since \(\mathcal U\cap M\in M\). Thus \(\mathcal W\) is amenable to \(M\).
\end{proof}
\end{lma}

\begin{thm}[UA]\label{WEM}
Let \(\kappa\) be the least supercompact cardinal. Then \(\textnormal{HOD}\) is
a weak extender model at \(\kappa\).
\begin{proof}
First note that every normal fine ultrafilter on an ordinal definable set is
ordinal definable. We will prove this using the fact that isomorphic normal fine
ultrafilters on the same set are equal (\cref{NormalIso}). Suppose \(\mathcal
U\) is a normal fine ultrafilter on \(Y\in \text{OD}\), and let \(U\) be a
countably complete ultrafilter on an ordinal isomorphic to \(\mathcal U\); then
by \cref{NormalIso}, \(\mathcal U\) is the {\it unique} normal fine ultrafilter
on \(Y\) isomorphic to \(U\), and hence \(\mathcal U\in \text{OD}_{Z,U} =
\text{OD}_U = \text{OD}\), with the final equality coming from \cref{ODUF}.

In particular, for all \(\lambda\geq \kappa\), every normal fine ultrafilter
\(\mathcal U\) on \(P_\kappa(\lambda)\) is amenable to \(\textnormal{HOD}\). The
issue is to show that there are such \(\mathcal U\) concentrating on
\(\text{HOD}\).

Fix a regular cardinal \(\delta > \kappa^{++}\). Then by \cref{HODClose},
\(\text{HOD}\) is correct about stationary subsets of \(\delta\). Let \(\langle
S_\alpha : \alpha < \delta\rangle\in \text{HOD}\) be a partition of
\(S^\delta_\omega\) into stationary subsets. Let \(j : V\to M\) be an elementary
embedding with critical point \(\kappa\) such that \(j(\kappa) > \delta\) and
\(j[\delta]\in M\). We claim that \(j[\delta]\in \text{HOD}^M\).

By \cref{MStationaryPartition}, \[j[\delta] = \{\alpha < j(\delta) : M\vDash
j(S)_\alpha \text{ is stationary in }\sup j[\delta]\}\] Thus \[j[\delta]\in
\text{HOD}^M_{j(\langle S_\alpha :\alpha < \delta\rangle)}\] But since \(\langle
S_\alpha :\alpha < \delta\rangle\) is in \(\textnormal{HOD}\), \(j(\langle
S_\alpha :\alpha < \delta\rangle)\in \textnormal{HOD}^M\). Thus \(j[\delta]\in
\text{HOD}^M\).

Let \(\mathcal U\) be the ultrafilter on \(P_\kappa(\delta)\) derived from \(j\)
using \(j[\delta]\). Since \(j[\delta]\in \text{HOD}^M = j(\textnormal{HOD})\),
\(\mathcal U\) concentrates on \(\textnormal{HOD}\) by \L o\'s's Theorem. Thus
\(\mathcal U\) is a normal fine \(\kappa\)-complete ultrafilter on
\(P_\kappa(\delta)\) that concentrates on \(\textnormal{HOD}\) and is amenable
to \(\textnormal{HOD}\).

This shows that for unboundedly many cardinals \(\delta\), there is a normal
fine \(\kappa\)-complete ultrafilter on \(P_\kappa(\delta)\) that concentrates
on \(\textnormal{HOD}\) and is amenable to \(\textnormal{HOD}\). Therefore by
\cref{UnboundedWEM}, HOD is a weak extender model at \(\kappa\).
\end{proof}
\end{thm}

As a consequence of theorems of Woodin \cite{Woodin}, this implies that a
version of Jensen's Covering Lemma is true for HOD:

\begin{cor}[UA] Any set \(A\subseteq \textnormal{HOD}\) is contained in a set
\(B\in \textnormal{HOD}\) such that \(|B| \leq |A|+\gamma\) for some \(\gamma <
\kappa\).\qed
\end{cor}
We omit the proof. Of course one has a much stronger covering results above
\(\kappa^{+++}\) as a consequence of \cref{VopSizeThm}.

\section{The Generalized Continuum Hypothesis}\label{GCHSection}
\subsection{Introduction}
In this section, we prove that GCH holds above the least supercompact assuming
UA. We actually prove a more local version of this theorem. We stress that
proving this local version requires some far from obvious tricks that are not
actually necessary for the global result. We need the local result at various
points in \cref{SCChapter1} and \cref{SCChapter2}.
\subsection{Two theorems of Solovay}
Let us begin by explaining the intuition that led to the expectation that UA
might imply the eventual Generalized Continuum Hypothesis. This begins with two
remarkable theorems of Solovay. First, of course, is his theorem on the Singular
Cardinals Hypothesis:
\begin{thm}[Solovay]\label{SolovaySCHLocal}
Suppose \(\kappa\) and \(\delta\) are cardinals with
\(\textnormal{cf}(\delta)\geq \kappa\). Suppose \(\kappa\) is
\(\delta\)-strongly compact. Then \(\delta^{<\kappa} = \delta\).\qed
\end{thm}
We will give a proof in \cref{SolovaySCH}. As a corollary, the Singular
Cardinals Hypothesis holds above a strongly compact cardinal:\index{Singular
Cardinals Hypothesis!above a strongly compact cardinal}
\begin{cor}\label{SolovaySCH}
Suppose \(\kappa \leq \lambda\) are cardinals, \(\lambda\) is a strong limit
singular cardinal, and \(\kappa\) is \(\lambda\)-strongly compact. Then
\(2^\lambda = \lambda^+\).
\begin{proof}
Note that \(2^\lambda = \lambda^{\text{cf}(\lambda)}\) since \(\lambda\) is a
strong limit cardinal and in general \(2^\lambda =
(2^{<\lambda})^{\text{cf}(\lambda)}\).

First assume \(\text{cf}(\lambda) < \kappa\). Then \[2^\lambda =
\lambda^{\text{cf}(\lambda)} \leq \lambda^{<\kappa} \leq (\lambda^+)^{<\kappa} =
\lambda^+\] by \cref{SolovaySCHLocal}.

Assume instead that \(\kappa\leq\text{cf}(\lambda)\).
\begin{comment}
	\footnote{Many texts
finish the theorem here by citing Silver's Theorem that SCH cannot first fail at
a cardinal of uncountable cofinality. This is anachronistic since Solovay's
Theorem preceded Silver's by several years. In fact, the argument that Solovay
originally used in this case must have inspired Magidor's proof that SCH cannot
first fail at \(\aleph_{\omega_1}\) assuming a precipitous ideal on
\(\omega_1\). It was Magidor's result that inspired Silver to prove his theorem
in the first place.}
\end{comment}
Assume by induction that \(2^{\gamma} = \gamma^+\) for all
strong limit singular cardinals in the interval \((\kappa,\lambda)\). Let
\(\iota = \text{cf}(\lambda)\). Let \(j:V\to M\) be an elementary embedding such
that \(\text{cf}^M(\sup j[\iota]) < j(\kappa)\), which exists since \(\kappa\)
is \(\text{cf}(\lambda)\)-strongly compact. Then \(\lambda_*= \sup j[\lambda]\)
is a strong limit singular cardinal of \(M\) and \[\text{cf}^M(\lambda_*) =
\text{cf}^M(\sup j[\iota]) < j(\kappa)\] Therefore \((2^{\lambda_*})^M =
\lambda_*^{+M}\). But there is an injection from \(P(\lambda)\) to
\(P^M(\lambda_*)\), namely the map \(A\mapsto j(A)\cap \lambda_*\). Therefore
\[2^\lambda \leq | \lambda_*^{+M}|\]

By the usual computations of cardinalities of ultrapowers
(\cref{UFCardinality}), the fact that \(\lambda_*\) is a strong limit implies
\(\lambda_* = \lambda\). Thus \(2^\lambda\leq \lambda_*^{+M} = \lambda^{+M} \leq
\lambda^+\), as desired.
\end{proof}
\end{cor}

The second of Solovay's theorems regards the number of normal fine ultrafilters
generated by supercompactness assumptions:
\begin{thm}[Solovay]\label{SolovayUltrafilters}
Suppose \(\kappa\) and \(\delta\) are cardinals, \(\kappa \leq
\textnormal{cf}(\delta)\), and \(\kappa\) is \(2^{\delta}\)-supercompact. Then
for all \(A\subseteq P(\delta)\), there is a normal fine \(\kappa\)-complete
ultrafilter \(\mathcal U\) on \(P_\kappa(\delta)\) such that \(A\in M_{\mathcal
U}\).
\end{thm}

Since this argument will be used repeatedly, it is worth working in a slightly
more general context.
\begin{lma}\label{MuGen}
	Suppose \(\lambda\) is a cardinal and \(j : V\to M\) is a
	\(\lambda\)-supercompact elementary embedding. Suppose \(Y\subseteq
	P(\lambda)\) is a set such that \(j[\lambda]\in j(Y)\). Let \(\mathcal U\)
	be the normal fine ultrafilter on \(Y\) derived from \(j\) using
	\(j[\lambda]\). Assume \(Y\in M_\mathcal U\) and \(\mathcal U\in M\). Then
	both \(M_\mathcal U\) and \(M\) satisfy the following statement: for any
	\(A\subseteq P(\lambda)\), there is a normal fine ultrafilter \(\mathcal W\)
	on \(Y\) such that \(A\in M_\mathcal W\).
	\begin{proof}
		Let \(k : M_\mathcal U\to M\) be the factor embedding with
		\(k(\id_\mathcal U) = j[\lambda]\) and \(k\circ j_\mathcal U =j\). Then
		by \cref{DerivedNF}, \(\textsc{crt}(k) > \lambda\) if \(k\) is
		nontrivial. In particular, \(k(\lambda) = \lambda\) and \(k(Y) = Y\).
		Therefore by the elementarity of \(k\), if \(M_\mathcal U\) satisfies
		the statement that for any \(A\subseteq P(\lambda)\), there is a normal
		fine ultrafilter \(\mathcal W\) on \(Y\) such that \(A\in M_\mathcal
		W\), then so does \(M\). Therefore it suffices to show that this
		statement is true in \(M_\mathcal U\).
		
		Let \(D\) be the ultrafilter derived from \(j\) using \(\langle \mathcal
		U,j[\lambda]\rangle\). Note that \(j_D: V\to M_D\) is a
		\(\lambda\)-supercompact elementary embedding, \(\mathcal U\) is the
		normal fine ultrafilter on \(Y\) derived from \(j_D\) using
		\(j[\lambda]\), and \(\mathcal U\in M_D\). By replacing \(j\) with
		\(j_D\), we may therefore assume that \(j\) is an ultrapower embedding.
		In particular, by \cref{UltrapowerSC}, \(M^\lambda \subseteq M\).
		
		We claim that \(P(P(\lambda)) \cap (M_\mathcal U)^M= P(P(\lambda))\cap
		M_\mathcal U\). (This is a consequence of \cref{MOSuper2}, but we give a
		proof here.) Since \(M\) is closed under \(\lambda\)-sequences,
		\((M_\mathcal U)^M = j_\mathcal U(M)\) by \cref{USupercompact}. In
		particular, \((M_\mathcal U)^M\subseteq M_\mathcal U\), so
		\(P(P(\lambda)) \cap (M_\mathcal U)^M\subseteq P(P(\lambda))\cap
		M_\mathcal U\). We now show the reverse inclusion. By \cref{Kunen},
		there is an inaccessible cardinal \(\kappa\leq \lambda\) such that
		\(j(\kappa) > \lambda\). Since \(j\) is \(\lambda\)-supercompact, \(j\)
		is \(\lambda\)-strong, and so in particular \(V_\kappa\subseteq M\).
		Therefore \[P(P(\lambda)) \cap (M_\mathcal U)^M \subseteq V_{j_\mathcal
		U(\kappa)}\cap M_\mathcal U = j_\mathcal U(V_\kappa)\subseteq j_\mathcal
		U(M) = (M_\mathcal U)^M\] as desired.

		Suppose \(A\subseteq P(\lambda)\) and \(A\in M_\mathcal U\). Note that
		\(k(A) = A \in  P(P(\lambda))\cap  M_\mathcal U \subseteq (M_\mathcal
		U)^M\). Thus \(M\) satisfies that \(k(A)\) belongs to \(M_\mathcal W\)
		for some normal fine ultrafilter on \(k(Y)\). (Namely, take \(\mathcal W
		= \mathcal U\).) By the elementarity of \(k\), it follows that
		\(M_\mathcal U\) satisfies that \(A\) belongs to \(M_\mathcal W\) for
		some normal fine ultrafilter on \(Y\).
		
		This shows that \(M_\mathcal U\) satisfies the statement that for any
		\(A\subseteq P(\lambda)\), there is a normal fine ultrafilter \(\mathcal
		W\) on \(Y\) such that \(A\in M_\mathcal W\), completing the proof.
	\end{proof}
\end{lma}

\begin{proof}[Proof of \cref{SolovayUltrafilters}] By our large cardinal
assumption, there is an elementary embedding \(j : V\to M\) such that the
following hold:
\begin{itemize}
\item \(\textsc{crt}(j) = \kappa\) and \(j(\kappa) > \delta\).
\item \(j\) is \(\delta\)-supercompact.
\item \(j\) is \(2^\delta\)-strong.
\end{itemize}
 Let \(\mathcal D\) be the normal fine \(\kappa\)-complete ultrafilter on
 \(P_\kappa(\delta)\) derived from \(j\) using \(j[\delta]\), and let \(k :
 M_\mathcal D\to M\) be the factor embedding.
\begin{clm}\label{muClm}
\(\mathcal D\in M\).
\begin{proof}
Since \(j\) is \(2^\delta\)-strong, \(H_{{2^\delta}^+}\subseteq M\). Since
\(|P_\kappa(\delta)| = \delta\) by \cref{SolovaySCH}, \(P(P_\kappa(\delta))\in
H_{{2^\delta}^+}\). Therefore since  \(\mathcal D\in P(P_\kappa(\delta))\),
\(\mathcal D\in M\).
\end{proof}
\end{clm}

Therefore applying \cref{MuGen} yields that \(M\) satisfies that for all
\(A\subseteq P(\delta)\), there is a normal fine ultrafilter \(\mathcal W\) on
\(P_\kappa(\delta)\) such that \(A\in M_\mathcal W\). But every \(A\subseteq
P(\delta)\) belongs to \(M\). Moreover, if \(M\) satisfies that \(\mathcal W\)
is a normal fine ultrafilter on \(P_\kappa(\delta)\), then \(\mathcal W\)
actually is a normal fine ultrafilter on \(P_\kappa(\delta)\) and moreover
\(P(P(\delta))\cap (M_\mathcal W)^M= P(P(\delta))\cap M_\mathcal W\). Thus we
can conclude that \(V\) satisfies that every \(A\subseteq P(\delta)\) belongs to
\(M_\mathcal W\) for some normal fine ultrafilter \(\mathcal W\) on
\(P_\kappa(\delta)\). This proves the theorem.
\end{proof}

As a corollary, Solovay observed that instances of GCH follow from the linearity
of the Mitchell order:
\begin{thm}[Solovay]\label{SolovayCounting}
Suppose \(\kappa\) is \(2^\kappa\)-supercompact and the set of normal
ultrafilters on \(\kappa\) is linearly ordered by the Mitchell order. Then
\(2^{2^\kappa} = (2^\kappa)^+\).
\begin{proof}
First note that for any normal ultrafilter \(U\) on \(\kappa\),
\[|P(P(\kappa))\cap M_U| \leq |j_U(V_\kappa)| \leq |(V_\kappa)^\kappa| =
2^\kappa\] The first inequality follows from the inclusion \(P(P(\kappa))\cap
M_U\subseteq j_U(V_\kappa)\), and the second from the existence of a surjection
\(\pi : (V_\kappa)^\kappa\to j_U(V_\kappa)\), defined by \(\pi(f) =
j_U(f)(\kappa)\).

Let \(\mathcal N\) be the set of normal ultrafilters on \(\kappa\). The result
will follow from counting \(\mathcal N\) in two different ways.

First, note that a normal ultrafilter \(U\) has at most \(2^\kappa\)
predecessors in the Mitchell order. We are assuming the Mitchell order on
\(\mathcal N\) is a wellorder, and therefore the ordertype of \((\mathcal
N,\mo)\) is at most \((2^\kappa)^+\): any proper initial segment of \((\mathcal
N,\mo)\) has cardinality \(2^\kappa\). In particular \(|\mathcal N| \leq
(2^\kappa)^+\).

Second, note that \[P(P(\kappa)) = \bigcup_{U\in \mathcal N} P(P(\kappa))\cap
M_U\] Thus \[2^{2^\kappa} = |P(P(\kappa))| = |\mathcal N|\cdot
\sup_{U\in\mathcal N} |P(P(\kappa))\cap M_U| = |\mathcal N|\cdot 2^\kappa\] It
follows that \(|\mathcal N| = 2^{2^\kappa}\).

Thus we have shown \(2^{2^\kappa} = |\mathcal N|\leq (2^\kappa)^+\). It follows
that \(2^{2^\kappa} = (2^\kappa)^+\), as desired.
\end{proof}
\end{thm}
More generally, and by exactly the same argument, one can show:
\begin{prp}[UA]\label{GeneralizedCounting}
Assume \(\lambda\) is a cardinal such that \(2^{<\lambda} = \lambda\) and every
\(A\subseteq P(\lambda)\) belongs to \(M_\mathcal W\) for some normal fine
ultrafilter \(\mathcal W\) on \(P_{\textnormal{bd}}(\lambda)\). Then
\(2^{2^\lambda} = (2^\lambda)^+\).
\begin{proof}
	Recall that \(\mathscr N_\lambda\) denotes the set of normal fine
	ultrafilters on \(P_{\textnormal{bd}}(\lambda)\).
	
	We claim that for any \(\mathcal W\in \mathscr N_\lambda\),
	\(P(P(\lambda))\cap M_\mathcal W\) has cardinality at most \(2^\lambda\). By
	\cref{IsoNormalThm}, \(\lambda_\mathcal W = \lambda\). By \cref{Kunen},
	there is an inaccessible cardinal \(\kappa\leq \lambda\) such that
	\(j_\mathcal W(\kappa) > \lambda\). Thus \(P(P(\lambda))\cap M_\mathcal
	W\subseteq j_\mathcal W(V_\kappa)\). But \(|j_\mathcal W(V_\kappa)|\leq
	|V_\kappa|^{\lambda_\mathcal W} = \kappa^\lambda = 2^\lambda\). In
	particular \(|P(P(\lambda))\cap M_\mathcal W|\leq 2^\lambda\).
	
	This bound has two consequences. 
	
	First, it follows that any \(\mathcal W\in \mathscr N_\lambda\) has at most
	\(2^\lambda\) predecessors in the Mitchell order. This is because if
	\(\mathcal U\mo \mathcal W\), then \(\mathcal U\in
	P(P(P_\text{bd}(\lambda)))\cap M_\mathcal W\), and
	\[|P(P(P_\text{bd}(\lambda)))\cap M_\mathcal W| = |P(P(\lambda))\cap
	M_\mathcal W|\] since \(|P_\text{bd}(\lambda)|^{M_\mathcal W} =
	(2^{<\lambda})^{M_\mathcal W} = \lambda\). (One does not actually need to
	use \(2^{<\lambda} = \lambda\) here, but it is convenient.) Hence
	\((\mathscr N_\lambda,\mo)\) is a wellfounded partial order of rank at most
	\((2^\lambda)^+\). Since \(2^{<\lambda} = \lambda\), \cref{GCHLinear}
	implies that \((\mathscr N_\lambda,\mo)\) is a wellorder, and hence
	\(|\mathscr N_\lambda| \leq (2^\lambda)^+\).
	
	Second, it follows that \(|\mathscr N_\lambda| = 2^{2^\lambda}\): by our
	assumption that every \(A\subseteq P(\lambda)\) belongs to \(M_\mathcal W\)
	for some \(\mathcal W\in \mathscr N_\lambda\),
	\[P(P(\lambda))  = \bigcup_{\mathcal W\in \mathscr N_\lambda}
	P(P(\lambda))\cap M_\mathcal W\] Thus \[2^{2^\lambda} = |P(P(\lambda))| =
	|\mathscr N_\lambda|\cdot \sup_{\mathcal W\in \mathscr N_\lambda}
	|P(P(\lambda))\cap M_\mathcal W| =|\mathscr N_\lambda|\cdot 2^\lambda =
	|\mathscr N_\lambda|\]
	
	Putting everything together, \(2^{2^\lambda} = |\mathscr N_\lambda| \leq
	(2^\lambda)^+\), which proves the theorem.
\end{proof}
\end{prp}

Regarding this lemma, a much more complicated argument (\cref{UFCounting}) shows
that under UA, a set \(X\) carries at most \((2^{|X|})^+\) countably complete
ultrafilters. 

Let us mention a little fact, proved very early on in this work, that gave the
first indication that GCH  above a supercompact cardinal might be provable from
UA:
\begin{prp}[UA] Suppose \(\kappa\) is supercompact. Let \(\lambda =
\beth_\kappa(\kappa)\). Then \(2^\delta = \delta^+\) for all cardinals \(\delta
\in [\lambda,\lambda^{+\omega}]\).
\begin{proof}
Since \(\lambda\) is a singular strong limit cardinal, \cref{SolovaySCH} implies
\(2^\lambda =\lambda^+\). Since \(2^{<\lambda} = \lambda\),
\cref{GeneralizedCounting} implies \(2^{2^\lambda} = (2^\lambda)^+\). In other
words, \(2^{(\lambda^+)} = \lambda^{++}\). Since \(2^{<\lambda^+} = 2^\lambda =
\lambda^+\), \cref{GeneralizedCounting} implies \(2^{2^{(\lambda^+)}} =
(2^{(\lambda^+)})^+\). In other words, \(2^{(\lambda^{++})} = \lambda^{+++}\).
Continuing this way yields the result for cardinals \(\delta\) such that
\(\lambda \leq \delta < \lambda^{+\omega}\). Then \(\lambda^{+\omega}\) is a
strong limit cardinal, so \(2^{(\lambda^{+\omega})} = \lambda^{+\omega+1}\) by
\cref{SolovaySCH}.
\end{proof}
\end{prp}
This proof breaks down completely at \(\lambda^{+\omega+1}\), and it gives no
hint of whether \(2^\kappa = \kappa^+\) should hold when \(\kappa\) is
supercompact. But the fact that one gets GCH at \(\omega+1\) cardinals in a row
strongly suggests that one should be able to prove the eventual GCH. To handle
the case \(\delta = \lambda^{+\omega+1}\) and the case \(\delta = \kappa\) turns
out to require a completely different argument, which we turn to now.
\subsection{More on the Mitchell order}\label{MoreMOSection}
\begin{defn}\label{LambdaMitchellDef}
	A countably complete ultrafilter \(U\) is {\it
	\(\lambda\)-Mitchell}\index{Ultrafilter!\(\lambda\)-Mitchell} if for all
	hereditarily uniform ultrafilters \(D\) such that \(\lambda_D < \lambda\),
	\(D\mo U\). 
\end{defn}

If \(\lambda = 2^{<\lambda}\) and \(U\) is a countably complete ultrafilter such
that \(j_U\) is \(\lambda\)-strong, then \(U\) is \(\lambda\)-Mitchell. The
first step in the proof of GCH we will give is to prove the same result without
assuming that \(2^{<\lambda} = \lambda\), and instead using UA and a
supercompactness hypothesis.

\begin{prp}[UA]\label{MitchellLemma}\index{Generalized Mitchell order!linearity}
Suppose \(\lambda\) is an infinite cardinal and \(U\) is a countably complete
ultrafilter such that \(j_U\) is \(\lambda\)-supercompact. Then \(U\) is
\(\lambda\)-Mitchell.
\end{prp}

In order to prove \cref{MitchellLemma}, we need two preliminary lemmas. The
first is the obvious attempt to extend the proof of the linearity of the
Mitchell order on normal ultrafilters from \cref{MOChapter} to normal fine
ultrafilters. (This was in fact the first proof we attempted in the very early
days of UA, before realizing that the generalization was related to cardinal
arithmetic.)

\begin{lma}\label{FirstAttemptLemma}
Suppose \(\lambda\) is a cardinal, \(U\) is a countably complete ultrafilter,
and \(j_U\) is \(\lambda\)-supercompact. Suppose \(D\) is a countably complete
ultrafilter on an ordinal \(\gamma \leq \lambda\). Suppose \((k,i) :
(M_D,M_U)\to N\) is a \(1\)-internal comparison of \((j_D,j_U)\) such that
\(k([\textnormal{id}]_D) \in i(j_U[\lambda])\). Then \(D\mo U\).
\begin{proof}
Note that for any \(A\subseteq \gamma\),
\begin{align*}
A\in D&\iff [\text{id}]_D\in j_D(A)\\
&\iff k([\text{id}]_D)\in k(j_D(A))\\
&\iff k([\text{id}]_D)\in i(j_U(A))\\
&\iff k([\text{id}]_D)\in i(j_U(A))\cap i(j_U[\lambda])\\
&\iff k([\text{id}]_D)\in i(j_U(A)\cap j_U[\lambda])\\
&\iff k([\text{id}]_D)\in i(j_U[A])
\end{align*}
Therefore \begin{equation}\label{DInternal}D = \{A\subseteq \gamma :
k([\text{id}]_D)\in i(j[A])\}\end{equation} Since \(j\restriction \gamma\in
M_U\), the function defined on \(P(\gamma)\) by \(A\mapsto j[A]\) belongs to
\(M_U\). Moreover \(i\) is an internal ultrapower embedding of \(M_U\).
Therefore \cref{DInternal} shows that \(D\) is definable over \(M_U\) from
parameters in \(M_U\), and hence \(D\mo U\).  
\end{proof}
\end{lma}

Incidentally, this lemma suggests considering the following generalized Ketonen
order: for \(D\in \mathscr B(X)\) and \(\mathcal U\in \mathscr B(Y)\), set \(D
\mathrel{\in^\Bbbk} \mathcal U\) if there exist \(I\in \mathcal U\) and
\(\langle D_\sigma : \sigma\in I\rangle\in \prod_{\sigma\in I} \mathscr
B(X,\sigma)\) such that \(D = \mathcal U\text{-}\lim_{\sigma\in I} D_\sigma\).
\cref{FirstAttemptLemma} can be restated as follows: if \(\lambda\) is a
cardinal, \(D\) is a countably complete ultrafilter on \(\lambda\), and
\(\mathcal U\) is a normal fine ultrafilter on \(P(\lambda)\), then \(D\in^\Bbbk
\mathcal U\) if and only if \(D\mo \mathcal U\).

Our next lemma puts us in a position to apply \cref{FirstAttemptLemma}. For the
proof of \cref{MitchellLemma}, we will only need the case \(A = j_U[\lambda]\),
but the general statement is used in the proof of level-by-level equivalence at
singular cardinals (\cref{ClubEQ}).

\begin{lma}\label{ClubLemma}
Suppose \(\lambda\) is a cardinal, \(U\) is a countably complete ultrafilter,
and \( A\subseteq j_U(\lambda)\) is a nonempty set that is closed under
\(j_U(f)\) for every \(f : \lambda\to \lambda\). Suppose \(D\) is a countably
complete ultrafilter on an ordinal \(\gamma < \lambda\). Suppose \((k,i) :
(M_D,M_U)\to N\) is a \(0\)-internal comparison of \((j_D,j_U)\). Then
\(k([\textnormal{id}]_D) \in i(A)\).
\begin{proof}
Let \(B = k^{-1}[i(A)]\). By the definition of a \(1\)-internal comparison, \(k
: M_D\to N\) is an internal ultrapower embedding, and therefore \(B\in M_D\). We
must show that \([\text{id}]_D\in B\). 

We first show that \(j_D[\lambda]\subseteq B\). Note that
\(j_U[\lambda]\subseteq A\) since \(A\) is nonempty and closed under
\(j(c_\alpha)\) for any \(\alpha < \lambda\), where \(c_\alpha : \lambda\to
\lambda\) is the constant function with value \(\alpha\). Thus \(i\circ
j_U[\lambda]\subseteq i(A)\). Since \((i,k)\) is a \(1\)-internal comparison,
\(k\circ j_D[\lambda] = i\circ j_U[\lambda]\subseteq i(A)\). So
\(j_D[\lambda]\subseteq k^{-1}[i(A)] = B\). 

We now show that \(B\) is closed under \(j_D(f)\) for any \(f : \lambda\to
\lambda\). Fix \(\xi \in B\) and \(f : \lambda\to \lambda\); we will show
\(j_D(f)(\xi)\in B\). By assumption \(A\) is closed under \(j_U(f)\), and so by
elementarity \(i(A)\) is closed under \(i(j_U(f))\). In particular, since
\(k(\xi)\in i(A)\), \(i(j_U(f))(k(\xi))\in i(A)\). But \(i(j_U(f))(k(\xi)) =
k(j_D(f)(\xi))\). Now \(k(j_D(f)(\xi))\in i(A)\) so \(j_D(f)( \xi)\in
k^{-1}(i(A)) = B\), as desired.

Since \(\gamma < \lambda\) and \(j_D[\lambda]\subseteq B\), in particular
\(j_D[\gamma^+]\subseteq B\). Thus \(B\) is cofinal in the \(M_D\)-regular
cardinal \(j_D(\gamma^+) = \sup j_D[\gamma^+]\). In particular, \(|B|^{M_D}\geq
j_D(\gamma^+)\). Fix \(\langle B_\xi : \xi < \gamma\rangle\) with \(B =
j_D(\langle B_\xi : \xi < \gamma\rangle)_{[\text{id}]_D}\). By \L o\'s's
Theorem, we may assume without loss of generality that \(B_\xi\subseteq
\lambda\) and \(|B_\xi| \geq \gamma^+\) for all \(\xi < \gamma\). Therefore
there is an injective function \(g :\gamma \to \lambda\) such that \(g(\xi)\in
B_\xi\) for all \(\xi < \gamma\). By \L o\'s's Theorem,
\(j_D(g)([\text{id}]_D)\in B\). Since \(g\) is injective, there is a function
\(f : \lambda\to \lambda\) be a function satisfying \(f(g(\xi)) = \xi\) for all
\(\xi < \gamma\). But \(B\) is closed under \(j_D(f)\), and
\(j_D(f)(j_D(g)([\text{id}]_D)) = [\text{id}]_D\), so \([\text{id}]_D\in B\), as
desired.
\end{proof}
\end{lma}

\cref{MitchellLemma} now follows easily.
\begin{proof}[Proof of \cref{MitchellLemma}] Fix a countably complete
hereditarily uniform ultrafilter \(D\) with \(\lambda_D < \lambda\). We must
show \(D\mo U\). By the isomorphism invariance of the Mitchell order on
hereditarily uniform ultrafilters (\cref{HeredLemma}), we may assume \(D\) lies
on an ordinal \(\gamma < \lambda\). By UA, there is an internal ultrapower
comparison \((k,i) : (M_D,M_U)\to N\) of \((j_D,j_U)\). By \cref{ClubLemma} with
\(A = j[\lambda]\), \(k([\text{id}]_D)\in i(j[\lambda])\). Therefore the
hypotheses of \cref{FirstAttemptLemma} are satisfied, so \(D\mo U\), as desired.
\end{proof}

It is natural to hope that the proof of \cref{MitchellLemma} can be generalized
to show the linearity of the Mitchell order on normal fine ultrafilters without
assuming GCH. The trouble of course is removing the assumption \(\gamma <
\lambda\) in \cref{ClubLemma}. If \(\lambda\) is regular (which turns out to be
the hard case), the proof of \cref{ClubLemma} goes through under the assumption
that \(k[j_D(\lambda)] = \sup i\circ j_U[\lambda]\). We know how to prove the
lower bound \(k[j_D(\lambda)] \leq \sup i\circ j_U[\lambda]\) (using
\cref{NormalGeneration}), but we do not know how to prove \(k[j_D(\lambda)] \geq
\sup i\circ j_U[\lambda]\) directly. The only proof we know of the linearity of
the Mitchell order that does not require a GCH assumption (\cref{ULinearity})
requires a good deal of the supercompactness analysis of \cref{SCChapter1}.

\subsection{The proof of GCH}\label{GCHProofSection}
\cref{GCH} above follows immediately from the following statement, which is much
more local (and much harder to prove):

\begin{thm}[UA]\label{MainTheorem}\label{MainThm}
Suppose \(\kappa \leq \delta\) are cardinals with \(\kappa \leq
\textnormal{cf}(\delta)\). If \(\kappa\) is \(\delta^{++}\)-supercompact, then
for any cardinal \(\lambda\) with \(\kappa\leq \lambda \leq\delta^{++}\),
\(2^\lambda = \lambda^+\).\index{Generalized Continuum Hypothesis!local proof}
\end{thm}

Combining \cref{MainTheorem} with the results of \cref{SCChapter1}, the
hypothesis that \(\kappa\) is \(\delta^{++}\)-supercompact can be weakened to
the assumption that \(\kappa\) is \(\delta^{++}\)-strongly compact.

The hard part of the proof is contained in the following theorem:
\begin{thm}[UA]\label{HardPart}
Suppose \(\kappa\) and \(\delta\) are cardinals such that
\(\textnormal{cf}(\delta)\geq \kappa\). Suppose \(\kappa\) is
\(\delta^{++}\)-supercompact. Then \(2^\delta = \delta^+\).
\begin{proof}
Assume towards a contradiction that \(2^\delta > \delta^+\). We use this
assumption to prove the following claim, following the proof of
\cref{SolovayUltrafilters}:
\begin{clm} Every subset of \(\delta^{++}\) belongs to the ultrapower of the universe by a normal fine \(\kappa\)-complete ultrafilter on \(P_\kappa(\delta)\). \end{clm}
\begin{proof}
	Let \(j : V\to M\) be a \(\delta^{++}\)-supercompact ultrapower embedding
	with \(\textsc{crt}(j)=\kappa\) and \(j(\kappa) > \delta^{++}\). Since
	\(P(\delta^{++})\subseteq M\), it suffices to show the claim is true in
	\(M\). Let \(\mathcal U\) be the normal fine ultrafilter on
	\(P_\kappa(\delta)\) derived from \(j\) using \(j[\delta]\). By
	\cref{MitchellLemma}, \(\mathcal U\in M\) (since in fact every countably
	complete ultrafilter on \(P_\kappa(\delta)\) is in \(M\)). Therefore by
	\cref{MuGen}, \(M\) satisfies that every subset of \(P(\delta)\) belongs to
	the ultrapower of the universe by a normal fine \(\kappa\)-complete
	ultrafilter on \(P_\kappa(\delta)\). Since \(\delta^{++}\leq (2^\delta)^M\),
	it follows that \(M\) satisfies that every subset of \(\delta^{++}\) belongs
	to the ultrapower of the universe by a normal fine \(\kappa\)-complete
	ultrafilter on \(P_\kappa(\delta)\), as desired.
\end{proof}

Let \(W\) be a \(\delta^+\)-supercompact ultrafilter on \(\delta^+\) with
\(j_W(\kappa) > \delta^+\). We claim \(P(\delta^{++})\subseteq M_W\). Suppose
\(A\subseteq \delta^{++}\). For some normal fine \(\kappa\)-complete ultrafilter
\(U\) on \(P_\kappa(\delta)\), \(A\in M_U\). But since \(|P_\kappa(\delta)|=
\delta\) (by Solovay's Theorem on SCH above a strongly compact cardinal,
\cref{SolovaySCHLocal}), \cref{MitchellLemma} implies \(U\in M_W\). It is easy
to see that this implies \(A\in M_W\).

Let \(Z\) be a \(\delta^{++}\)-supercompact ultrafilter on \(\delta^{++}\) with
\(j_Z(\kappa) > \delta^{++}\). Let \(k : M_W\to N\) be the ultrapower of \(M_W\)
by \(Z\) using functions in \(M_W\). We have \(\textsc{width}(j_W) =
\delta^{++}\) and \(\textsc{width}(k) = \delta^{++} < j_W(\delta^{++})\), so by
the lemma on the width of compositions (\cref{WidthLemma}),
\(\textsc{width}(k\circ j_W) = \delta^{++}\). In other words, there is an
ultrafilter \(D\) on \(\delta^+\) such that \(M_D = N\) and \(j_D = k\circ
j_W\). 

Since \(P(\delta^{++})\subseteq M_W\), \((V_\kappa)^{\delta^{++}}\subseteq
M_W\). Therefore letting \(\kappa' = j_Z(\kappa) = k(\kappa)\), we have
\(V_{\kappa'}\cap N = V_{\kappa'}\cap M_Z\). By \cref{MitchellLemma}, \(D\mo
Z\), and so since \(\kappa'> \delta^+\), \(D\in V_{\kappa'}\cap M_Z\subseteq
N\). It follows that \(D\in N = M_D\), contradicting the irreflexivity of the
Mitchell order (\cref{MOStrict}).
\end{proof}
\end{thm}

\begin{proof}[Proof of \cref{MainTheorem}] Suppose \(\lambda\) is a cardinal
with \(\kappa \leq \lambda \leq \delta^{++}\). 
\begin{case} \(\lambda\leq \delta\)\end{case}
If \(\lambda\) is regular then by \cref{HardPart}, \(2^\lambda = \lambda^+\). If
\(\lambda\) is singular then \(2^{<\lambda} = \lambda\) by \cref{HardPart}, so
\(2^\lambda = \lambda^+\) by the local version of Solovay's theorem
\cite{Solovay}.
\begin{case}\label{Last} \(\lambda =\delta^+\). \end{case}
Since \(\kappa\) is \(\delta^+\)-supercompact and \(2^\delta = \delta^+\),
\(\kappa\) is \(2^\delta\)-supercompact. Therefore by
\cref{GeneralizedCounting}, \(2^{2^\delta} = (2^\delta)^+\). In other words,
\(2^{(\delta^+)} = \delta^{++}\). 
\begin{case} \(\lambda =\delta^{++}\) \end{case}
Given that \(2^{(\delta^+)} =\delta^{++}\) by \cref{Last}, the case that
\(\lambda = \delta^{++}\) can be handled in the same way as \cref{Last}.
\end{proof}

\begin{cor}[UA] Suppose \(\kappa\leq \delta\) and \(\kappa\) is
\(2^\delta\)-supercompact. Then \(2^\delta = \delta^+\).
\begin{proof}
Assume first that \(\delta\) is singular. Since \(\kappa\) is
\({<}\delta\)-supercompact, \(2^{<\delta} = \delta\) by \cref{MainTheorem}. Now
\(2^\delta = \delta^+\) by \cref{SolovaySCH}.

Assume instead that \(\delta\) is regular. Assume towards a contradiction that
\(2^\delta \geq \delta^{++}\). Then \(\kappa\) is \(\delta^{++}\)-supercompact,
so by \cref{MainTheorem}, \(2^\delta = \delta^+\), a contradiction.
\end{proof}
\end{cor}

Let us point out another consequence that one can obtain using a result in
\cref{SCChapter1}:

\begin{thm}[UA]\label{SuccGCH}
Suppose \(\nu\) is a cardinal and \(\nu^+\) carries a countably complete uniform
ultrafilter. Then \(2^{<\nu} = \nu\).
\begin{proof}
By \cref{SuccessorThm} below, some cardinal \(\kappa \leq\nu\) is
\(\nu^+\)-supercompact. If \(\kappa = \nu\) then obviously \(2^{<\nu} = \nu\).
So assume \(\kappa < \nu\). If \(\nu\) is a limit cardinal, then the hypotheses
of \cref{MainTheorem} hold for all sufficiently large \(\lambda < \nu\) and
hence GCH holds on a tail below \(\nu\), so \(2^{<\nu} = \nu\). So assume \(\nu
= \lambda^+\) is a successor cardinal. If \(\lambda\) is singular, then
\(\lambda\) is a strong limit singular cardinal by \cref{MainTheorem}, so
\(2^\lambda = \lambda^+\) by Solovay's theorem \cref{SolovaySCH}, and hence
\(2^{<\nu} = \nu\). Finally if \(\lambda\) is regular, we can apply
\cref{MainTheorem} directly to conclude that \(2^\lambda = \lambda^+\), so again
\(2^{<\nu} = \nu\). 
\end{proof}
\end{thm}

This leaves open some questions about further localizations of the GCH proof.
\begin{qst}[UA] Suppose \(\kappa\) is \(\delta\)-supercompact. Must \(2^\delta =
\delta^+\)?
\end{qst}
We conjecture that it is consistent with UA that \(\kappa\) is measurable but
\(2^\kappa > \kappa^+\), which would give a negative answer in the case \(\kappa
= \delta\). In certain cases, the question has a positive answer as an
essentially immediate consequence of our main theorem:
\begin{prp}[UA] Suppose \(\kappa \leq \lambda\), \(\textnormal{cf}(\lambda) =
\omega\), and \(\kappa\) is \(\lambda\)-supercompact. Then \(2^\lambda =
\lambda^+\).

Suppose \(\kappa \leq \lambda\), \(\omega_1 \leq \textnormal{cf}(\lambda) <
\lambda\), and \(\kappa\) is \({<}\lambda\)-supercompact. Then \(2^\lambda =
\lambda^+\).

Suppose \(\kappa \leq \lambda\), \(\lambda\) is the double successor of a
cardinal of cofinality at least \(\kappa\), and \(\kappa\) is
\(\lambda\)-supercompact. Then \(2^\lambda = \lambda^+\).\qed
\end{prp}

Another interesting localization question is the following:
\begin{qst}[UA] Suppose \(\kappa\) is the least ordinal \(\alpha\) such that
there is an ultrapower embedding \(j : V\to M\) with \(j(\alpha) >
(2^{\kappa})^+\). Must \(2^\kappa = \kappa^+\)?
\end{qst}

\subsection{\(\diamondsuit\) on the critical cofinality}
We conclude with the observation that stronger combinatorial principles than GCH
follow from UA.

\begin{thm}[UA]\label{Diamond}
Suppose \(\kappa\) is \(\delta^{++}\)-supercompact where
\(\textnormal{cf}(\delta)\geq \kappa\). Then
\(\diamondsuit(S^{\delta^{++}}_{\delta^+})\) holds.
\end{thm}
For the proof, we need a theorem of Kunen.
\begin{defn}
Suppose \(\lambda\) is a regular uncountable cardinal and \(S\subseteq \lambda\)
is a stationary set. Suppose \(\langle \mathcal A_\alpha : \alpha \in S\rangle\)
is a sequence of sets with \(\mathcal A_\alpha \subseteq P(\alpha)\) and
\(|A_\alpha| \leq \alpha\) for all \(\alpha < \lambda\). Then \(\langle \mathcal
A_\alpha : \alpha \in S\rangle\) is a \(\diamondsuit^-(S)\)-sequence if for all
\(X\subseteq \lambda\), \(\{\alpha \in S : X\cap \alpha\in \mathcal A_\alpha\}\)
is stationary.
\end{defn}
\begin{defn}
\(\diamondsuit^-(S)\) is the assertion that there is a
\(\diamondsuit^-(S)\)-sequence.
\end{defn}

\begin{thm}[Kunen, \cite{KunenText}]\label{KunenDiamond}
Suppose \(\lambda\) is a regular uncountable cardinal and \(S\subseteq \lambda\)
is a stationary set. Then \(\diamondsuit^-(S)\) is equivalent to
\(\diamondsuit(S)\).\qed
\end{thm}

\begin{proof}[Proof of \cref{Diamond}] By \cref{MainTheorem}, GCH holds on the
interval \([\kappa,\delta^{++}]\), and we will use this without further comment.

For each \(\alpha < \delta^{++}\), let \(\mathcal U_\alpha\) be the unique
ultrafilter of rank \(\alpha\) in the Mitchell order on normal fine
\(\kappa\)-complete ultrafilters on \(P_\kappa(\delta)\). The uniqueness of
\(\mathcal U_\alpha\) follows from the linearity of the Mitchell order on normal
fine ultrafilters on \(P_\kappa(\delta)\), a consequence of \cref{GCHLinear}
which applies in this context since \(2^{<\delta} = \delta\). Let \(\mathcal
A_\alpha = P(\alpha)\cap M_{\mathcal U_\alpha}\). Note that \(|\mathcal
A_\alpha| \leq \kappa^\delta = \delta^+\). Let \[\vec{\mathcal A} = \langle
\mathcal A_\alpha : \alpha < \delta^{++}\rangle\] Note that \(\vec{\mathcal A}\)
is definable in \(H_{\delta^{++}}\) without parameters.

\begin{clm} \(\vec{\mathcal A}\) is a \(\diamondsuit^-(S^{\delta^{++}}_{\delta^+})\)-sequence.\end{clm}
\begin{proof}
Suppose towards a contradiction that \(\vec{\mathcal A}\) is not a
\(\diamondsuit^-(S^{\delta^{++}}_{\delta^+})\)-sequence. Let \(\mathcal W\) be a
\(\kappa\)-complete normal fine ultrafilter on \(P_\kappa(\delta^{++})\). Then
in \(M_\mathcal W\),  \(\vec{\mathcal A}\) is not a
\(\diamondsuit^-(S^{\delta^{++}}_{\delta^+})\)-sequence. Let \(\mathcal U\) be
the \(\kappa\)-complete normal fine ultrafilter on \(\delta\) derived from
\(\mathcal W\) and let \(k : M_\mathcal U \to M_\mathcal W\) be the factor
embedding. Let \(\gamma = \textsc{crt}(k) = \delta^{++M_\mathcal U}\).

Since  \(\vec{\mathcal A}\) is definable in \(H_{\delta^{++}}\) without
parameters, \(\vec {\mathcal A}\in \text{ran}(k)\). Therefore \(k^{-1}(\vec
{\mathcal A}) = \vec {\mathcal A}\restriction \gamma\) is not a
\(\diamondsuit^-(S^{\gamma}_{\delta^+})\)-sequence in \(M_\mathcal U\). Fix a
witness \(A\in P(\gamma)\cap M_\mathcal U\) and a closed unbounded set \(C\in
P(\gamma)\cap M_\mathcal U\) such that for all \(\alpha \in C\cap
S^{\gamma}_{\delta^+}\), \(A\cap \alpha\notin \mathcal A_\alpha\). By
elementarity, for all \(\alpha\in k(C)\cap S^{\delta^{++}}_{\delta^+}\),
\(k(A)\cap \alpha\notin \mathcal A_\alpha\). Since \(\mathcal U\) is
\(\delta\)-supercompact, \(\text{cf}(\gamma) = \delta^+\), and so in particular
\(k(A)\cap \gamma\notin \mathcal A_\gamma\). Since \(\gamma = \textsc{crt}(k)\),
this means \(A\notin \mathcal A_\gamma\).

Note however that \(\mathcal U\) has Mitchell rank \(\delta^{++M_\mathcal U} =
\gamma\), so \(\mathcal U = \mathcal U_\gamma\). Therefore \(\mathcal A_\gamma =
P(\gamma)\cap M_\mathcal U\), so \(A\in \mathcal A_\gamma\) by choice of \(A\).
This is a contradiction.
\end{proof}
By \cref{KunenDiamond}, this completes the proof.
\end{proof}

\subsection{The size of the Vop\v enka algebra}\label{VopenkaSection}

\begin{thm}[UA]\label{VopenkaSize}
Suppose \(\kappa\) is an inaccessible cardinal such that every \(A \subseteq
P(\kappa)\) belongs to \(M_U\) for some countably complete ultrafilter \(U\) on
\(\kappa\). Then \(|\mathbb V_\kappa|^{\textnormal{HOD}} = (2^\kappa)^+\).
\begin{proof}
Let \(\lambda = |\mathbb V_\kappa|^{\text{HOD}}\). Note that \(\lambda =
|P(P(\kappa))\cap \textnormal{OD}|^\textnormal{OD}\).

Recall that \(\mathscr B(\kappa)\) denotes the set of countably complete
ultrafilters on \(\kappa\). As in \cref{SolovayCounting}, \(|\mathscr B(\kappa)|
= 2^{2^\kappa}\).

We claim that in fact \(|\mathscr B(\kappa)| = (2^\kappa)^+\). It suffices to
show the upper bound \(\mathscr B(\kappa)\leq 2^{2^\kappa}\). For this, we show
that every initial segment of the Ketonen order has cardinality \(2^\kappa\).

 Since \(\kappa\) is inaccessible, for any \(\alpha < \kappa\), the set
 \(\mathscr B(\kappa,\alpha)\) of countably complete ultrafilters on \(\kappa\)
 that concentrate on \(\alpha\) has cardinality less than \(\kappa\). Thus for
 any \(U\in \mathscr B(\kappa)\), \(U\) has at most \(2^\kappa\cdot
 \prod_{\alpha < \kappa} |S_\alpha| = 2^\kappa\) predecessors in the Ketonen
 order, since if \(W\sE U\), then
\[W = U\text{-}\lim_{\alpha\in I} W_\alpha\] for some \(I\in U\) and some
sequence \(\langle W_\alpha : \alpha\in I\rangle\in \prod_{\alpha\in I}\mathscr
B(\kappa,\alpha)\).

Therefore let \(\langle U_\alpha : \alpha < (2^\kappa)^+\rangle\) be the
\(\sE\)-increasing enumeration of \(\mathscr B(\kappa)\).

For the lower bound \( (2^\kappa)^+\leq \lambda\), we apply the fact that every
countably complete ultrafilter on an ordinal is OD (\cref{ODUF}) to obtain
\(\mathscr B(\kappa)\subseteq P(P(\kappa))\cap \textnormal{OD}\), so in fact
\(\lambda \geq |\mathscr B(\kappa)| = 2^{2^\kappa} = (2^\kappa)^+\).

We now turn to the upper bound.

Suppose \(U\in \mathscr B(\kappa)\). Then \(|P(P(\kappa))\cap M_U| \leq
|j_U(V_\kappa)| \leq |(V_\kappa)^\kappa| = 2^\kappa\). Let \[\mathcal A_U =
P(P(\kappa))\cap M_U\cap \text{OD}\] Note that \(P(P(\kappa)\cap M_U\cap
\text{OD}\) is an ordinal definable subset of \(\text{OD}\), so let \(\gamma_U =
|\mathcal A_U|^\text{OD}\) and let \(\pi_U : \gamma_U\to \mathcal A_U\) be the
\(\text{OD}\)-least bijection. Note that \(|\mathcal A_U| \leq 2^\kappa\) so
\(\gamma_U < (2^\kappa)^+\).

Let \(\lambda_0 = \sup\{\gamma_U : U\in S\}\), so \(\lambda_0 \leq
(2^\kappa)^+\). Define \(\pi : (2^\kappa)^+\times \lambda_0\to P(P(\kappa))\cap
\textnormal{OD}\) by
\[\pi(\alpha,\beta) = \pi_{f(\alpha)}(\beta)\] Then our large cardinal
assumption on \(\kappa\) implies that \(\pi\) is a surjection and \(\pi\) is
ordinal definable, so \(\lambda\leq  (2^\kappa)^+\cdot \lambda_0 =
(2^\kappa)^+\).
\end{proof}
\end{thm}
We finally prove \cref{VopSizeThm}, the fact that under UA, if \(\kappa\) is
supercompact then \(V\) is a generic extension of \(\textnormal{HOD}\) for a
forcing of size \(\kappa^{++}\).
\begin{thm}[UA]\label{VopSizeThm}
	If \(\kappa\) is \(\kappa^{++}\)-supercompact then \(|\mathbb
	V_\kappa|^\textnormal{HOD} = \kappa^{++}\). 
\end{thm}
\begin{proof}
Note that since \(\kappa\) is \(\kappa^{++}\)-supercompact, by
\cref{MainTheorem}, \((2^\kappa)^+ =\kappa^{++}\). In particular, \(\kappa\) is
\(2^\kappa\)-supercompact, so the hypotheses of \cref{VopenkaSize} hold by
\cref{SolovayUltrafilters}. Thus \(|\mathbb V_\kappa|^\text{HOD} = (2^\kappa)^+
=  \kappa^{++}\). 
\end{proof}
\chapter{The Least Supercompact Cardinal}\label{SCChapter1}
\section{Introduction}
\subsection{The identity crisis}
How large is the least strongly compact cardinal? This question was first posed
by Tarski in a precise form shortly after his discovery of strong compactness:
{\it is the least strongly compact cardinal larger than the least measurable
cardinal?} About a decade later, Solovay mounted the first serious attack on
this problem. He fused the Scott's elementary embedding analysis of
measurability with the combinatorial properties of strongly compact cardinals to
produce what has become the central large cardinal concept: supercompactness. He
then conjectured that every strongly compact cardinal is supercompact. This is
certainly a natural conjecture to make since supercompact cardinals and strongly
compact cardinals share some rather deep structural similarities. (See
\cref{StrongCSection} and especially \cref{KetonenSection}.) But unlike the
least strongly compact cardinal, the size of the least supercompact cardinal is
no mystery at all: it is upon first glance a staggeringly large object, much
larger than the least measurable cardinal. Thus Solovay's conjecture implies a
positive answer to Tarski's question.

Telis Menas, then a graduate student under Solovay at UC Berkeley, was the first
to realize that Solovay's conjecture is false. Menas climbed up far beyond the
least strongly compact cardinal, and up there he discovered a strongly compact
cardinal that is not supercompact.
\begin{repthm}{Menas}[Menas] The least strongly compact limit of strongly
	compact cardinals is not supercompact.
\end{repthm} 
This theorem closed off Solovay's approach to Tarski's question while leaving
the question itself wide open. The fundamental breakthrough occurred mere months
after Menas's discovery on the other side of the world, with Magidor's landmark
independence result \cite{Magidor}:
\begin{thm*}[Magidor]\index{Identity crisis}
	Suppose \(\kappa\) is a cardinal.
	\begin{itemize}
		\item If \(\kappa\) is strongly compact, then there is a forcing
		extension in which \(\kappa\) remains strongly compact but becomes the
		least measurable cardinal.
		\item If \(\kappa\) is supercompact, then there is a forcing extension
		in which \(\kappa\) remains supercompact but becomes the least strongly
		compact cardinal.\qed
	\end{itemize}
\end{thm*}
Thus the ZFC axioms are insufficient to answer Tarski's question. Magidor
described this peculiar situation as an  ``identity crisis" for the least
strongly compact cardinal. The main result of this chapter is that the
Ultrapower Axiom resolves this crisis:
\begin{repthm}{StrongSuper}[UA] The least strongly compact cardinal is
	supercompact.
\end{repthm}
We will prove much stronger results than this that explain exactly why the least
strongly compact cardinal is supercompact, and that identify much weaker
properties that are sufficient (under UA) for supercompactness. We defer until
the final chapter the analysis of larger strongly compact cardinals.
\subsection{Outline of \cref{SCChapter1}}
We now outline the rest of the chapter.\\

\noindent {\sc \cref{StrongCSection}.} We exposit the basic theory of strong compactness.
We use the theory of the Ketonen order to prove Ketonen's Theorem \cite{Ketonen}
that \(\kappa\) is strongly compact if and only if every regular cardinal
carries a \(\kappa\)-complete ultrafilter (\cref{KetonenThm}). This argument is
the basis for much of the theory of this chapter. We use Ketonen's Theorem to
prove the local version of Solovay's Theorem \cite{Solovay} on SCH above a
strongly compact that we have by now cited several times (\cref{SolovayThm}).\\

\noindent {\sc \cref{KLambdaSection}.} In this section, we introduce the notion of a
Fr\'echet cardinal and its associated Ketonen ultrafilters. Under UA, each
Fr\'echet cardinal \(\lambda\) carries a unique Ketonen ultrafilter \(\mathscr
K_\lambda\). For regular \(\lambda\), we analyze \(\mathscr K_\lambda\) under
the assumption that some \(\kappa \leq \lambda\) is  \(\lambda\)-strongly
compact, showing that its associated embedding is \({<}\lambda\)-supercompact
and \(\lambda\)-tight (\cref{GeneralThm}).\\

\noindent {\sc \cref{FrechetSection}.} Given the analysis of \(\mathscr K_\lambda\) in the
previous section, we would like to show that if \(\lambda\) is a regular
Fr\'echet cardinal, then some cardinal \(\kappa\leq \lambda\) is
\(\lambda\)-strongly compact. In this section, we prove our best result towards
this, showing that this is true unless \(\lambda\) is {\it isolated}
(\cref{NonisolatedCompact}). Isolated cardinals are rare enough that this
implies the supercompactness of the least strongly compact cardinal
(\cref{StrongSuper}).\\

\noindent {\sc \cref{IsolationSection}.} In this section we study the structure of
isolated cardinals, which arose in \cref{FrechetSection} as a pathological case
in our analysis of Fr\'echet cardinals. We rule out pathological isolated
cardinals assuming GCH (\cref{IsolatedStrongLimit}). Without GCH, assuming just
UA, we are still able to fully analyze ultrafilters on an isolated cardinal
(\cref{IsolatedUFSection}), which turn out to look just like ultrafilters on the
least measurable cardinal. We prove that nonmeasurable isolated cardinals are
associated with serious failures of GCH (\cref{LowerBound} and
\cref{UpperBound}).  We leverage these results to prove the linearity of the
Mitchell order on normal fine ultrafilters without assuming GCH
(\cref{ULinearity}).
\section{Strong compactness}\label{StrongCSection}
\subsection{Some characterizations of strong compactness}
After almost 70 years of research, strong compactness remains one of the most
important and mysterious large cardinal notions. Strongly compact cardinals were
first isolated by Tarski in the context of infinitary logic: \(\kappa\) is
strongly compact if the logic \(\mathcal L_{\kappa,\kappa}\) satisfies a
generalized version of the Compactness Theorem. In keeping with modern large
cardinal theory, we will introduce strongly compact cardinals in terms of
elementary embeddings of the universe of sets into inner models with closure
properties. The closure property we have in mind is a two-cardinal version of
the covering property:
\begin{defn}\label{CoveringDef}\index{Covering property}
Suppose \(M\) is an inner model, \(\lambda\) is a cardinal, and \(\delta\) is an
\(M\)-cardinal. Then \(M\) has the {\it \((\lambda,\delta)\)-covering property}
if every set \(A\subseteq M\) such that \(|A| < \lambda\) is contained in a set
\(B\in M\) such that \(|B|^M < \delta\).
\end{defn}

\begin{defn}\index{Strongly compact cardinal}
	A cardinal \(\kappa\) is strongly compact if for any cardinal
	\(\lambda\geq\kappa\), there is an elementary embedding \(j : V\to M\) such
	that \(\textsc{crt}(j) = \kappa\) and \(M\) has the
	\((\lambda,j(\kappa))\)-covering property.
\end{defn}

\begin{defn}
We make the following abbreviations:
\begin{itemize}
\item The {\it \(({\leq}\lambda,\delta)\)-covering property} is the
\((\lambda^+,\delta)\)-covering property. 
\item The {\it \((\lambda,{\leq}\delta)\)-covering property} is the \((\lambda,
\delta^{+M})\)-covering property. 
\item The {\it \(({\leq}\lambda,{\leq}\delta)\)-covering property} is the
\((\lambda^+,\delta^{+M})\)-covering property.
\item The {\it \(\lambda\)-covering property} is the
\((\lambda,\lambda)\)-covering property. 
\item The {\it \({\leq}\lambda\)-covering property} is the
\(({\leq}\lambda,{\leq}\lambda)\)-covering property.
\end{itemize}
\end{defn}

This notation is chosen so that, for example, an inner model \(M\) has the
\(({\leq}\lambda,{\leq}\delta)\)-covering property if every subset \(A\subseteq
M\) such that \(|A|\leq\lambda\) is contained in a set \(B\in M\) such that
\(|B|^M \leq \delta\).

We will be particularly interested in the following local version of strong
compactness (especially when \(\lambda\) is regular):
\begin{defn}
Suppose \(\kappa\leq \lambda\) are cardinals. Then \(\kappa\) is {\it
\(\lambda\)-strongly compact} if there is an inner model \(M\) and an elementary
embedding \(j : V\to M\) with \(\textsc{crt}(j) = \kappa\) such that \(M\) has
the \(({\leq}\lambda,j(\kappa))\)-covering property.
\end{defn}

Note that if \(j : V\to M\) and \(M\) has the
\(({\leq}\lambda,j(\kappa))\)-covering property, then \(j(\kappa) > \lambda\).

\cref{StrongCChar} puts down several equivalent reformulations of strong
compactness. These involve the notions of {\it tightness} and {\it filter
bases}, which we now define.

The concept of {\it tightness} had not been given a name before this
dissertation, but it plays a role analogous to that of supercompactness in the
theory of supercompact cardinals:
\begin{defn}\index{Tightness of an elementary embedding}
Suppose \(M\) is an inner model, \(\lambda\) is a cardinal, and \(\delta\) is an
\(M\)-cardinal. An elementary embedding \(j : V\to M\) is {\it
\((\lambda,\delta)\)-tight} if there is a set \(A\in M\) with \(|A|^M \leq
\delta\) such that \(j[\lambda]\subseteq A\). An elementary embedding is said to
be {\it \(\lambda\)-tight} if it is \((\lambda,\lambda)\)-tight.
\end{defn}

Thus \((\lambda,\delta)\)-tightness is a weakening of
\(\lambda\)-supercompactness. Any \(j :V\to M\) such that \(M\) has the
\(({\leq}\lambda,j(\kappa))\)-covering property is
\(({\leq}\lambda,j(\kappa))\)-tight. Moreover, many of the general theorems
about supercompact embeddings generalize to the context of
\((\lambda,\delta)\)-tight ones. For example, \cref{XSC} generalizes:
\begin{lma}\label{XTight}
Suppose \(j :V\to M\) is an elementary embedding. The following are equivalent:
\begin{enumerate}[(1)]
\item \(j\) is \((\lambda,\delta)\)-tight.
\item For some \(X\) with \(|X| = \lambda\), there is some \(Y\in M\) with
\(|Y|^M \leq \delta\) such that \(j[X]\subseteq Y\).
\item For any \(A\) such that \(|A| \leq \lambda\), there is some \(B\in M\)
with \(|B|^M \leq \delta\) such that \(j[A]\subseteq B\).
\end{enumerate}
\begin{proof}
{\it (1) implies (2):} Trivial.

{\it (2) implies (3):} Suppose \(|A|\leq \lambda\). We will find \(B\in M\) with
\(|B|^M \leq \delta\) such that \(j[A]\subseteq B\). Using (2), fix \(X\) with
\(|X| = \lambda\) such that for some  \(Y\in M\) with \(|Y|^M \leq \delta\) such
that \(j[X]\subseteq Y\). Let \(p : X\to A\) be a surjection. Then \[j[A] =
j(p)[j[X]]\subseteq j(p)[Y]\] Let \(B = j(p)[Y]\). Then \(j[A]\subseteq B\),
\(B\in M\), and \(|B|^{M} \leq |Y|^M \leq \delta\).

{\it (3) implies (1):} Trivial.
\end{proof}
\end{lma}

The relationship between the \(\lambda\)-supercompactness of an embedding and
the closure of its target model under \(\lambda\)-sequences is analogous to the
relationship between the \((\lambda,\delta)\)-tightness of an elementary
embedding and the \(({\leq}\lambda,{\leq}\delta)\)-covering property of its
target model. For example, there is an analog of \cref{UltrapowerSC}:

\begin{lma}\label{UltrapowerStrongC}\index{Tightness of an elementary embedding!vs. the covering property}
Suppose \(j : V\to M\) is a \((\lambda,\delta)\)-tight ultrapower embedding.
Then \(M\) has the \(({\leq}\lambda,{\leq}\delta)\)-covering property.
\begin{proof}
Suppose \(A\subseteq M\) with \(|A|\leq \lambda\), and we will find \(B\in M\)
such that \(|B|^M \leq \delta\) and \(A\subseteq B\). Fix \(a\in M\) such that
\(M = H^M(j[V]\cup \{a\})\). Fix a set of functions \(F\) of cardinality
\(\lambda\) such that \(A = \{j(f)(a) : f\in F\}\). By \cref{XTight}, fix \(G\in
M\) with \(|G|^M \leq \delta\) and such that \(j[F]\subseteq G\). Let \(B =
\{g(a) : g\in G\}\). Then \(B\in M\), \(A\subseteq B\), and \(|B|^M \leq |G|^M
\leq \delta\), as desired.
\end{proof}
\end{lma}

Many filters are most naturally presented in terms of a smaller family of sets
that ``generates" the filter. The notion of a filter base makes this precise:
\begin{defn}\index{Filter base}
A {\it filter base} on \(X\) is a family \(\mathcal B\) of subsets of \(X\) with
the finite intersection property: for all \(A_0,A_1\in \mathcal B\), \(A_0\cap
A_1\neq \emptyset\). If \(\kappa\) is a cardinal, a filter base \(\mathcal B\)
is said to be {\it \(\kappa\)-complete} if for all \(\nu < \kappa\), for all
\(\{A_\alpha : \alpha < \nu\}\subseteq \mathcal B\), \(\bigcap_{\alpha < \nu}
A_\alpha\neq \emptyset\).
\end{defn}
The term ``filter base" is motivated by the fact that every filter base
\(\mathcal B\) on \(X\) generates a filter.

\begin{defn}
	Suppose \(\mathcal B\) is a filter base. The {\it filter generated by
	\(\mathcal B\)}\index{Filter base!filter generated by} is the filter
	\(F(\mathcal B) = \{A\subseteq X : \exists A_0,\dots,A_{n-1}\in \mathcal B\
	A_0\cap\cdots\cap A_{n-1}\subseteq A\}\). If \(\mathcal B\) is a
	\(\kappa\)-complete filter base, the {\it \(\kappa\)-complete filter
	generated by \(\mathcal B\)} is the filter
	\[F_\kappa(\mathcal B) = \left\{A\subseteq X :  \exists S\in
	P_\kappa(\mathcal B)\ \textstyle\bigcap S \subseteq A\right\}\]
\end{defn}

We finally prove our equivalences with strong compactness:

\begin{thm}\label{StrongCChar}\index{Strongly compact cardinal}
Suppose \(\kappa\leq \lambda\) are uncountable cardinals. Then the following are
equivalent:
\begin{enumerate}[(1)]
\item \(\kappa\) is \(\lambda\)-strongly compact.
\item There is an elementary embedding \(j : V\to M\) with critical point
\(\kappa\) that is \((\lambda,\delta)\)-tight for some \(M\)-cardinal \(\delta <
j(\kappa)\).
\item Every \(\kappa\)-complete filter base of cardinality \(\lambda\) extends
to a \(\kappa\)-complete ultrafilter.
\item There is a \(\kappa\)-complete fine ultrafilter on \(P_\kappa(\lambda)\).
\item There is an ultrapower embedding \(j : V\to M\) with critical point
\(\kappa\) that is \((\lambda,\delta)\)-tight for some \(M\)-cardinal \(\delta <
j(\kappa)\).
\item There is an elementary embedding \(j: V\to M\) with critical point
\(\kappa\) such that \(M\) has the \((\lambda,\delta)\)-covering property for
some \(M\)-cardinal \(\delta < j(\kappa)\).
\end{enumerate}
\begin{proof}
{\it (1) implies (2):} Trivial.

{\it (2) implies (3):} Let \(j : V\to M\) be an elementary embedding such that
\(\textsc{crt}(j)= \kappa\) and \(j\) is \((\lambda,\delta)\)-tight for some
\(M\)-cardinal \(\delta < j(\kappa)\). Suppose \(\mathcal B\) is a
\(\kappa\)-complete filter base on \(X\) of cardinality \(\lambda\). By
\cref{XTight}, there is a set \(S\in M\) such that \(j[\mathcal B]\subseteq S\)
and \(|S|^M < j(\kappa)\). By replacing \(S\) with \(S\cap j(\mathcal B)\), we
may assume without loss of generality that \(S\subseteq j(\mathcal B)\). By the
elementarity of \(j\), since \(j(\mathcal B)\) is \(j(\kappa)\)-complete, the
intersection \(\bigcap j(S)\) is nonempty. Fix \(a\in \bigcap j(S)\). Since
\(j[\mathcal B]\subseteq S\), it follows that \(a\in j(A)\) for all \(A\in
\mathcal B\). Let \(U\) be the ultrafilter on \(X\) derived from \(j\) using
\(a\). Then \(U\) extends \(\mathcal B\) and \(U\) is \(\kappa\)-complete since
\(\textsc{crt}(j) = \kappa\).

{\it (3) implies (4):} For any \(\alpha < \lambda\), let \(A_\alpha = \{\sigma
\in P_\kappa(\lambda) : \alpha\in \sigma\}\), and let \({\mathcal B}= \{A_\alpha
: \alpha < \lambda\}\). Then \(\mathcal B\) a \(\kappa\)-complete filter base on
\(P_\kappa(\lambda)\), and any filter on \(P_\kappa(\lambda)\) that extends
\({\mathcal B}\) is fine. By (3), there is a \(\kappa\)-complete ultrafilter
extending \({\mathcal B}\). Thus there is a \(\kappa\)-complete fine ultrafilter
on \(P_\kappa(\lambda)\), as desired.

{\it (4) implies (5):} Suppose \(\mathcal U\) is a \(\kappa\)-complete fine
ultrafilter on \(P_\kappa(\lambda)\). Let \(j : V\to M\) be the ultrapower of
the universe by \(\mathcal U\). The \(\kappa\)-completeness of \(\mathcal U\)
implies that \(\textsc{crt}(j) \geq \kappa\). By \cref{NormalFineChar},
\(j[\lambda]\subseteq \id_\mathcal U\). Moreover \(\id_\mathcal U\in
j(P_\kappa(\lambda))\), so letting \(\delta = |\id_\mathcal U|^M\), \(\delta <
j(\kappa)\). Therefore \(j\) is an ultrapower embedding that is
\((\lambda,\delta)\)-tight for some \(\delta < j(\kappa)\). Since \(\kappa\leq
\lambda\) and \(\lambda \leq \text{ot}(j[\lambda]) \leq \delta^{+M} <
j(\kappa)\), it follows that \(j(\kappa) > \lambda\). In particular,
\(\textsc{crt}(j) = \kappa\).

{\it (5) implies (6):} This is an immediate consequence of the fact that tight
ultrapowers have the covering property (\cref{UltrapowerStrongC}).

{\it (6) implies (1):} Trivial.
\end{proof}
\end{thm}
\subsection{Ketonen's Theorem}\label{KetonenSection}
The main theorem of this subsection is a famous theorem of Ketonen
\cite{Ketonen} that amounts to a deeper ultrafilter theoretic characterization
of strong compactness:
\begin{thm}[Ketonen]\label{KetonenStrongC}\index{Ketonen's Theorem on strongly compact cardinals}
	A cardinal \(\kappa\) is strongly compact if and only if every regular
	cardinal \(\lambda\geq \kappa\) carries a uniform \(\kappa\)-complete
	ultrafilter.
\end{thm}

Part of what is surprising about this theorem is that it does not even require
that the ultrafilters in the hypothesis be \(\kappa^+\)-incomplete. Beyond this,
it is not even obvious at the outset that the existence of \(\kappa\)-complete
ultrafilters on, say, \(\kappa\) and \(\kappa^+\) implies that \(\kappa\) is
\(\kappa^+\)-strongly compact.

We begin, however, with a less famous but no less important theorem of Ketonen,
which is also a key step in the proof of \cref{KetonenStrongC}. This theorem is
in a sense the strongly compact generalization of Solovay's Lemma
\cite{Solovay}. Suppose \(j : V\to M\) is an elementary embedding. For regular
cardinals \(\lambda\), Solovay's Lemma (or more specifically \cref{SCStatCor})
yields a simple criterion for whether \(j\) is \(\lambda\)-supercompact solely
in terms of the inner model \(M\) and the ordinal \(\sup j[\lambda]\):
\begin{repthm}{SCStatCor}[Solovay] Suppose \(j : V\to M\) is an elementary
embedding and \(\lambda\) is a regular cardinal. Then \(j\)  is
\(\lambda\)-supercompact if and only if \(M\) is correct about stationary
subsets of \(\sup j[\lambda]\).\qed
\end{repthm}
Ketonen proved a remarkable analog of this theorem for strongly compact
embeddings:
\begin{thm}[Ketonen]\label{KetonenCov}
Suppose \(j : V\to M\) is an elementary embedding, \(\lambda\) is a regular
uncountable cardinal, and \(\delta\) is an \(M\)-cardinal. Then \(j\) is
\((\lambda,\delta)\)-tight if and only if \(\textnormal{cf}^M(\sup
j[\lambda])\leq \delta\).\index{Tightness of an elementary embedding!and
\(\text{cf}^M(\sup j[\lambda])\)}
\end{thm}

For example, suppose \(j : V\to M\) is an ultrapower embedding.
\cref{KetonenCov} implies that all that is required for \(M\) to have the
\({\leq}\lambda\)-covering property is that \(M\) correctly compute the
cofinality of \(\sup j[\lambda]\).

The proof of \cref{KetonenCov} we give is due to Woodin, and is a bit different
from Ketonen's original proof. The trick is to choose the cover first, and then
choose the set whose image is being covered:
\begin{proof}[Proof of \cref{KetonenCov}] First assume \(j\) is
\((\lambda,\delta)\)-tight. Fix \(A\in M\) with \(j[\lambda]\subseteq A\) such
that \(|A|^M \leq \delta\). Then \(A\cap \sup j[\lambda]\) is cofinal in \(\sup
j[\lambda]\), so \(\sup j[\lambda]\) has cofinality at most \(|A|^M\) in \(M\).

Now we prove the converse. Assume \(\textnormal{cf}^M(\sup j[\lambda])\leq
\delta\). Let \(Y\in M\) be an \(\omega\)-closed cofinal subset of \(\sup
j[\lambda]\) of order type at most \(\delta\). Note that \(j[\lambda]\) is
itself an \(\omega\)-closed cofinal subset of \(\sup j[\lambda]\), so since
\(\sup j[\lambda]\) has uncountable cofinality, \(Y\cap j[\lambda]\) is an
\(\omega\)-closed cofinal subset of \(\lambda\). In particular, since
\(\text{cf}(\sup j[\lambda]) = \lambda\), \(Y\cap j[\lambda]\) has order type at
least \(\lambda\). Let \(X = j^{-1}[Y]\). Then \(j[X] = Y\cap j[\lambda]\), so
\[\text{ot}(X) = \text{ot}(j[X]) = \text{ot}(Y\cap j[\lambda])\geq \lambda\]
Thus \(|X| = \lambda\). Since \(|X| = \lambda\), \(Y\in M\), \(j[X]\subseteq
Y\), and \(|Y|^M \leq \delta\), \cref{XTight} implies that \(j\) is
\((\lambda,\delta)\)-tight.
\end{proof}

With \cref{KetonenCov} in hand, we turn to the proof of Ketonen's
characterization of strong compactness. The key point is that the strong
compactness of an elementary embedding is equivalent to an ultrafilter theoretic
property:
\begin{prp}\label{TightUF}
Suppose \(\kappa \leq\lambda\) are uncountable cardinals and \(\lambda\) is
regular. Suppose \(M\) is an inner model and \(j : V\to M\) is an elementary
embedding. Suppose every regular cardinal in the interval \([\kappa,\lambda]\)
carries a uniform \(\kappa\)-complete ultrafilter. Then the following are
equivalent:
\begin{enumerate}[(1)]
\item \(j\) is \((\lambda,\delta)\)-tight for some \(M\)-cardinal \(\delta <
j(\kappa)\).
\item \(\sup j[\lambda]\) carries no \(j(\kappa)\)-complete tail uniform
ultrafilter.
\end{enumerate}
\begin{proof}
{\it (1) implies (2):} Assume (1). By \cref{KetonenCov}, \(\text{cf}^M(\sup
j[\lambda]) < j(\kappa)\). Therefore the tail filter on \(\sup j[\lambda]\) is
not \(j(\kappa)\)-complete in \(M\), so \(\sup j[\lambda]\) does not carry a
\(j(\kappa)\)-complete tail uniform ultrafilter in \(M\).

{\it (2) implies (1):} Assume (2). Then in particular \(\textnormal{cf}^{M}(\sup
j[\lambda])\) carries no uniform \(j(\kappa)\)-complete ultrafilter in \(M\). By
elementarity, every \(M\)-regular cardinal in the interval
\(j([\kappa,\lambda])\) carries a uniform \(\kappa\)-complete ultrafilter.
Therefore \(\textnormal{cf}^{M}(\sup j[\lambda])\) does not lie in the interval
\(j([\kappa,\lambda])\). Clearly \(\textnormal{cf}^{M}(\sup j[\lambda]) \leq
j(\lambda)\), so it follows that \(\textnormal{cf}^{M}(\sup j[\lambda]) <
j(\kappa)\).
\end{proof}
\end{prp}

Ketonen introduced the Ketonen order as a tool to prove the following theorem,
generalizing a theorem of Solovay that states that any measurable cardinal
carries a normal ultrafilter that concentrates on nonmeasurable cardinals.
\begin{thm}[Ketonen]\label{KetonenExistence}
Suppose \(\lambda\) is a regular cardinal. If \(\lambda\) carries a
\(\kappa\)-complete uniform ultrafilter, then \(\lambda\) carries a
\(\kappa\)-complete uniform ultrafilter \(U\) such that \(\sup j_U[\lambda]\)
carries no tail uniform \(\kappa\)-complete ultrafilter in \(M_U\). Indeed, any
\(\sE\)-minimal \(\kappa\)-complete uniform ultrafilter on \(\lambda\) has this
property.
\begin{proof}
Let \(U\) be a \(\sE\)-minimal element of the set of uniform \(\kappa\)-complete
ultrafilters on \(\lambda\). Suppose towards a contradiction that in \(M_U\),
\(\sup j_U[\lambda]\) carries a tail uniform \(\kappa\)-complete ultrafilter.
Equivalently, there is a \(\kappa\)-complete ultrafilter \(Z\) on
\(j_U(\lambda)\) such that \(\delta_Z = \sup j_U[\lambda]\). Let \(W =
j_U^{-1}[Z]\). Then \(\textsc{crt}(j_W)\geq \textsc{crt}(j_Z^{M_U}\circ j_U)\)
(by \cref{LimitFactor}), so \(W\) is \(\kappa\)-complete. Moreover since
\(\delta_Z = \sup j_U[\lambda]\), \(\delta_W = \lambda\). Thus \(W\) is a
\(\kappa\)-complete uniform ultrafilter on \(\lambda\). Since \(Z\) concentrates
on \(\sup j_U[\lambda] \leq \id_U\), \(W\sE U\) by the definition of the Ketonen
order (\cref{KOChars}). This contradicts the \(\sE\)-minimality of \(U\).
\end{proof}
\end{thm}

We can now prove a local version of Ketonen's theorem, which fits into the list
of reformulations of \(\lambda\)-strong compactness from \cref{StrongCChar}:
\begin{thm}[Ketonen]\label{KetonenThm}
Suppose \(\kappa \leq \lambda\) are regular uncountable cardinals. Then the
following are equivalent: 
\begin{enumerate}[(1)]
\item \(\kappa\) is \(\lambda\)-strongly compact.
\item Every regular cardinal in the interval \([\kappa,\lambda]\) carries a
uniform \(\kappa\)-complete ultrafilter.
\item \(\lambda\) carries a \(\kappa\)-complete ultrafilter \(U\) such that
\(j_U\) is \((\lambda,\delta)\)-tight for some \(\delta < j_U(\kappa)\).
\end{enumerate}
\begin{proof}
{\it (1) implies (2):} Note that the Fr\'echet filter on a regular cardinal
\(\delta\) is \(\delta\)-complete. Thus (2) follows from (1) as an immediate
consequence of the filter extension property of strongly compact cardinals
(\cref{StrongCChar} (3)).

{\it (2) implies (3)}: Assume (2). By \cref{KetonenExistence}, there is a
\(\kappa\)-complete ultrafilter \(U\) on \(\lambda\)  such that \(\sup
j_U[\lambda]\) carries no tail uniform \(\kappa\)-complete ultrafilter in
\(M_U\). Therefore by \cref{TightUF}, \(j_U\) is \((\lambda,\delta)\)-tight for
some \(\delta < j_U(\kappa)\).

{\it (3) implies (1):} See \cref{StrongCChar} (5).
\end{proof}
\end{thm}
\subsection{Solovay's Theorem}
In this section we give a proof of a local version of Solovay's theorem that we
use throughout this dissertation.
\begin{thm}[Solovay]\label{SolovayThm}
	\index{Singular Cardinals Hypothesis!above a strongly compact cardinal}
Suppose \(\kappa\leq\lambda\) are uncountable cardinals, \(\lambda\) is regular,
and \(\kappa\) is \(\lambda\)-strongly compact. Then \(\lambda^{<\kappa} =
\lambda\).
\end{thm}

We need the following lemma, which is in a sense an analog of
\cref{UFSuperBound}, though much easier:
\begin{lma}\label{UFStrongCBound}
Suppose \(U\) is a countably complete ultrafilter. Let \(j : V\to M\) be the
ultrapower of the universe by \(U\). Then for any \(\eta\geq \lambda_U^+\),
\(j\) is not \((\eta,\delta)\)-tight for any \(M\)-cardinal \(\delta <
j(\eta)\).
\begin{proof}
We may assume by induction that \(\eta\) is a successor cardinal. In particular,
\(\eta\) is regular, so by \cref{UFContinuity}, \(j(\eta) = \sup j[\eta]\).
Suppose towards a contradiction that \(\delta < j(\eta)\) is an \(M\)-cardinal
such that \(j\) is \((\eta,\delta)\)-tight. By \cref{KetonenCov},
\(\text{cf}^M(j(\eta)) = \text{cf}^M(\sup j[\eta]) \leq\delta < j(\eta)\). This
contradicts that \(\eta\) is regular in \(M\) by elementarity.
\end{proof}
\end{lma}

\begin{lma}\label{DumbLemma}
Suppose \(\kappa \leq \gamma\) are cardinals. Suppose \(\gamma\) is singular and
\begin{equation}\label{AnnoyingAssumption}\sup_{\eta <\gamma}\eta^{<\kappa}\leq
\gamma\end{equation} Suppose \(\gamma^+\) carries a uniform \(\kappa\)-complete
ultrafilter \(U\). Then \(\gamma^{<\kappa} \leq \gamma^+\).
\begin{proof}
Let \(\lambda = \gamma^+\). We will prove the equivalent statement that
\(\lambda^{<\kappa} = \lambda\).

Let \(j :V\to M\) be the ultrapower of the universe by \(U\).  Let \(\delta =
\text{cf}^M(\sup j[\lambda])\). Note that \(\delta < j(\lambda) \), so
\(\delta\leq j(\gamma)\). In fact, since \(j(\gamma)\) is singular in \(M\),
\(\delta < j(\gamma)\). Therefore by \cref{AnnoyingAssumption} and the
elementarity of \(j\): \begin{equation}\label{ClusterFuck}(\delta^{<\kappa})^M
\leq (\delta^{<j(\kappa)})^M \leq j(\gamma)\end{equation}  By \cref{KetonenCov},
\(j\) is \((\lambda,\delta)\)-tight, so we can  fix \(B\in M\) with
\(j[\lambda]\subseteq B\). Now \(j[P_\kappa(\lambda)]\subseteq B^{<\kappa}\),
and since \(M\) is closed under \(\kappa\)-sequences, \(B^{<\kappa}\in M\).
\cref{XTight} now implies that \(j\) is
\((\lambda^{<\kappa},(\delta^{<\kappa})^M)\)-tight. 

Assume towards a contradiction that \(\lambda^{<\kappa} \geq \lambda^+\). Then
\(j\) is \((\lambda^+,(\delta^{<\kappa})^M)\)-tight. Since \(\lambda_U =
\lambda\), it follows from \cref{UFStrongCBound} that \((\delta^{<\kappa})^M\geq
j(\lambda^+)\), contradicting \cref{ClusterFuck}.
\end{proof}
\end{lma}
We now prove Solovay's theorem:
\begin{proof}[Proof of \cref{SolovayThm}] Suppose \(\kappa\) is
\(\lambda\)-strongly compact. Assume by induction that for all regular \(\iota <
\lambda\), \(\iota^{<\kappa} = \iota\). Since \(\lambda\) is regular, every
element of \(P_\kappa(\lambda)\) is bounded below \(\lambda\), so
\(P_\kappa(\lambda) = \bigcup_{\eta < \lambda} P_\kappa(\eta)\). Thus computing
cardinalities: \[\lambda^{<\kappa} = \sup_{\eta < \lambda} \eta^{<\kappa}\] If
\(\lambda\) is a limit cardinal, it follows immediately from our induction
hypothesis that \(\lambda^{<\kappa} = \lambda\). Therefore assume \(\lambda\) is
a successor cardinal. If the cardinal predecessor of \(\lambda\) is a regular
cardinal \(\iota\), then applying our induction hypothesis we obtain:
\[\lambda^{<\kappa} = \sup_{\eta < \lambda} \eta^{<\kappa} = \lambda\cdot
\iota^{<\kappa} = \lambda\] Therefore assume the cardinal predecessor of
\(\lambda\) is a singular cardinal \(\gamma\). Then \(\sup_{\eta
<\gamma}\eta^{<\kappa}\leq \gamma\). In this case, by \cref{DumbLemma},
\(\lambda^{<\kappa} = \lambda\).
\end{proof}

\section{Fr\'echet cardinals and the least ultrafilter \(\mathscr K_\lambda\)}\label{KLambdaSection}
\subsection{Fr\'echet cardinals}
In this section, we begin our systematic study of strong compactness assuming
UA. We will ultimately prove that UA implies that strong compactness and
supercompactness coincide to the extent that this is possible. (A theorem of
Menas shows that assuming sufficiently large cardinals, not all strongly compact
cardinals are supercompact; see \cref{MenasSection}.) An oddity of the proof is
that it requires a preliminary analysis of the first strongly compact cardinal.
Indeed, to obtain the strongest results, one must enact a hyperlocal analysis of
essentially the weakest ultrafilter-theoretic forms of strong compactness.

With this in mind, we introduce the following central concept:
\begin{defn}\index{Fr\'echet cardinal}
	An uncountable cardinal \(\lambda\) is {\it Fr\'echet} if \(\lambda\)
	carries a countably complete uniform ultrafilter.
\end{defn} 
Fr\'echet cardinals almost certainly do not appear in the work of Fr\'echet.
Their name derives from the fact that \(\lambda\) is Fr\'echet if and only if
the Fr\'echet filter on \(\lambda\) extends to a countably complete ultrafilter.

The following proposition is almost tautological:
\begin{prp}\label{FrechetTaut}
	A cardinal \(\lambda\) is Fr\'echet if and only if \(\lambda = \lambda_U\)
	for some countably complete ultrafilter \(U\).\qed
\end{prp}

For regular cardinals \(\lambda\), we have the following obvious
characterizations of Fr\'echetness:
\begin{prp}\label{FrechetRegChar}
	Suppose \(\lambda\) is a regular uncountable cardinal. The following are
	equivalent:
	\begin{enumerate}[(1)]
		\item \(\lambda\) is Fr\'echet.
		\item There is a countably complete tail uniform ultrafilter on
		\(\lambda\).
		\item Some ordinal of cofinality \(\lambda\) carries a tail uniform
		ultrafilter.
		\item Every ordinal of cofinality \(\lambda\) carries a tail uniform
		ultrafilter.
		\item There is an elementary embedding \(j : V\to M\) that is
		discontinuous at \(\lambda\).
		
	\end{enumerate}
\begin{proof}
	{\it (1) implies (2):} Since \(\lambda\) is a cardinal, any uniform
	ultrafilter on \(\lambda\) is tail uniform. Thus since there is a countably
	complete uniform ultrafilter on \(\lambda\), there is a countably complete
	tail uniform ultrafilter on \(\lambda\).
	
	{\it (2) implies (3):} Trivial.
	
	{\it (3) implies (4):} Recall that two ordinals \(\alpha\) and \(\beta\)
	have the same cofinality if and only if there is a weakly order preserving
	cofinal function \(f: (\alpha,\leq)\to (\beta,\leq)\). In particular,
	\(f_*(T_\alpha) = T_\beta\) where \(T_\alpha\) is the tail filter on
	\(\alpha\).  Thus if \(\alpha\) carries a countably complete tail uniform
	ultrafilter \(U\), then so does \(\beta\), namely \(f_*(U)\).
	
	{\it (4) implies (5):} Suppose \(U\) is a countably complete tail uniform
	ultrafilter on \(\lambda\). Let \(j : V\to M\) be the ultrapower of the
	universe by \(U\). Note that for any \(\alpha < \lambda\), \(j(\alpha) <
	\id_U\) since \(\alpha < \delta_U\). Thus \(\sup j[\lambda] \leq \id_U <
	j(\lambda)\). In other words, \(j\) is discontinuous at \(\lambda\).
\end{proof}
\end{prp}
Singular Fr\'echet cardinals are more subtle, especially when one does not
assume the Generalized Continuum Hypothesis. The following fact gives a sense of
how singular Fr\'echet cardinals should arise:
\begin{prp}\label{FrechetSingularChar}
	Suppose \(\lambda\) is a singular limit of Fr\'echet cardinals. Let
	\(\iota\) be the cofinality of \(\lambda\). Then \(\lambda\) is Fr\'echet if
	and only if \(\iota\) is Fr\'echet.
	\begin{proof}
		If \(\lambda\) is Fr\'echet, then \(\iota\) is Fr\'echet by
		\cref{FrechetRegChar} (4), and this does not require that \(\lambda\) is
		a limit of Fr\'echet cardinals.
		
		We now turn to the converse. Let \(\langle \lambda_\alpha : \alpha <
		\iota\rangle\) be an increasing cofinal sequence of Fr\'echet cardinals
		less than \(\lambda\). Let \(U_\alpha\) be a countably complete
		ultrafilter on \(\lambda\) with \(\lambda_{U_\alpha} = \lambda_\alpha\).
		Let \(D\) be a countably complete uniform ultrafilter on \(\iota\). Let
		\[U = D\text{-}\lim_{\alpha < \iota} U_\alpha\] 
		
		Clearly \(U\) is a countably complete ultrafilter on \(\lambda\). We
		claim that \(U\) is uniform, or in other words that every set \(X\in U\)
		has cardinality \(\lambda\). Suppose \(X\subseteq \lambda\) is such a
		set. By the definition of ultrafilter limits, \(\{\alpha < \iota : X\in
		U_\alpha\}\in D\). Since \(D\) is a uniform ultrafilter, the set
		\(\{\alpha < \iota : X\in U_\alpha\}\) is unbounded in \(\iota\).
		Therefore \(X\in U_\alpha\) for unboundedly many \(\alpha < \iota\), and
		in particular \(|X| \geq \lambda_{U_\alpha} = \lambda_\alpha\) for
		unboundedly many \(\alpha < \iota\). Thus \(|X| \geq \sup_{\alpha <
		\iota} \lambda_\alpha = \lambda\), as desired. Since \(\lambda\) carries
		a countably complete uniform ultrafilter, follows that \(\lambda\) is a
		Fr\'echet cardinal. 
	\end{proof}
\end{prp}
\cref{FrechetSingularChar} tells us that when \(\lambda\) is a singular limit of
Fr\'echet cardinals, whether \(\lambda\) is Fr\'echet depends only on whether
the regular cardinal \(\text{cf}(\lambda)\) is Fr\'echet. One might therefore
hope to reduce problems about Fr\'echet cardinals in general to the regular
case, where we have a bit more information. It is not provable in ZFC, however,
that a singular Fr\'echet cardinal must be a limit of Fr\'echet cardinals. The
Fr\'echet cardinals where this fails are called {\it isolated cardinals}, and
arise as a major issue in our analysis of strong compactness under UA.  Isolated
cardinals are studied in \cref{FrechetSection} and especially
\cref{IsolationSection}.

\subsection{Ketonen ultrafilters}
The following definition is inspired by the proof of \cref{KetonenThm}, which
turned on the existence of a \(\kappa\)-complete ultrafilter \(U\) on
\(\lambda\) such that \(\sup j_U[\lambda]\) carries no \(\kappa\)-complete tail
uniform ultrafilter in \(M_U\).

Recall from \cref{RegWeaklyNormal} that a uniform ultrafilter \(U\) on a regular
cardinal \(\lambda\) is weakly normal if and only if letting \(j :V\to M\) be
the ultrapower of the universe by \(U\), \(\id_U = \sup j[\lambda]\).
Equivalently, \(U\) is weakly normal if it is closed under decreasing diagonal
intersections. 
\begin{defn}\label{KetonenRegDef}
If \(\lambda\) is a regular cardinal, an ultrafilter \(U\) on \(\lambda\) is a
{\it Ketonen ultrafilter} if the following hold:\index{Ketonen ultrafilter!on a
regular cardinal}
\begin{itemize}
\item \(U\) is countably complete and weakly normal.
\item \(U\) concentrates on ordinals that carry no countably complete tail
uniform ultrafilter.
\end{itemize}
\end{defn}

By \cref{RegWeaklyNormal} and \cref{FrechetRegChar}, we have the following
characterization of Ketonen ultrafilters on regular cardinals:

\begin{lma}\label{KetonenChar}
	Suppose \(\lambda\) is a regular cardinal and \(U\) is a countably complete
	ultrafilter on \(\lambda\). Then \(U\) is Ketonen if and only if \(\id_U =
	\sup j_U[\lambda]\) and either of the following equivalent statements holds:
	 \begin{itemize}
	 	\item \(\sup j_U[\lambda]\) carries no countably complete tail uniform
	 	ultrafilter in \(M_U\).
	 	\item \(\textnormal{cf}^{M_U}(\sup j_U[\lambda])\) is not Fr\'echet in
	 	\(M_U\).\qed
	 \end{itemize}
\end{lma}
In this way the key ordinal  \(\textnormal{cf}^{M_U}(\sup j_U[\lambda])\) from
\cref{KetonenCov} arises immediately in the study of Ketonen ultrafilters on
regular cardinals.

The following theorem asserts that Ketonen ultrafilters are analogous to
\(\lambda\)-minimal ultrafilters of \cref{WeaklyNormalSection}, except that
Ketonen ultrafilters are minimal in the Ketonen order rather than merely being
minimal in the Rudin-Keisler order.
\begin{lma}\label{KetonenMinimality}
Suppose \(\lambda\) is a regular cardinal. Then \(U\) is a Ketonen ultrafilter
on \(\lambda\) if and only if \(U\) is a \(\sE\)-minimal element of the set of
countably complete uniform ultrafilters on \(\lambda\).
\begin{proof}
Suppose first that \(U\) is a Ketonen ultrafilter. Let \[\alpha = \id_U = \sup
j_U[\lambda]\] Suppose \(W\sE U\). We will show that \(\lambda_W < \lambda\). By
the definition of the Ketonen order (\cref{KOChars}), there is some \(Z\in
\mathscr B^{M_U}(j_U(\lambda),\sup j_U[\lambda])\) such that \(j_U^{-1}[Z] =
W\). Since \(\sup j_U[\lambda]\) does not carry a countably complete tail
uniform ultrafilter in \(M_U\), there is some \(\beta < \sup j_U[\lambda]\) such
that \(Z\) concentrates on \(\beta\). Fix \(\alpha < \lambda\) such that
\(j_U(\alpha) \geq \beta\). Then \(j_U(\alpha) \in Z\), so \(\alpha\in W\). Thus
\(\lambda_W < \lambda\) as desired. 

Conversely, assume \(U\) is a \(\sE\)-minimal element of the set of uniform
ultrafilters on \(\lambda\). In particular, \(U\) is an \(\rRK\)-minimal element
of the set of uniform ultrafilters on \(\lambda\), which by
\cref{WeaklyNormalrRK} is equivalent to being weakly normal. 

Finally, fix \(Z\in \mathscr B^{M_U}(j_U(\lambda))\), and we will show that
\(\delta_Z < \sup j_U[\lambda]\). Let \(W = j_U^{-1}[Z]\). Then \(W\sE U\) by
the definition of the Ketonen order. It from the minimality of \(U\) that
\(\delta_W < \lambda\), so for some \(\alpha < \lambda\), \(\alpha\in W\). Now
\(j_U(\alpha)\in Z\), so \(\delta_Z \leq j_U(\alpha) < \sup j_U[\lambda]\), as
desired.

It follows that \(\sup j_U[\lambda]\) does not carry a countably complete tail
uniform ultrafilter in \(M_U\), so \(U\) is Ketonen by \cref{KetonenChar}.
\end{proof}
\end{lma}

Reflecting on \cref{KetonenMinimality}, we obtain a definition of Ketonen
ultrafilters on arbitrary cardinals:
\begin{defn}\index{Ketonen ultrafilter}
Suppose \(\lambda\) is a Fr\'echet cardinal. An ultrafilter \(U\) on \(\lambda\)
is {\it Ketonen} if \(U\) is a \(\sE\)-minimal element of the set of countably
complete uniform ultrafilters on \(\lambda\).
\end{defn}

The wellfoundedness of the Ketonen order (\cref{KOWellfounded}) immediately
yields the existence of Ketonen ultrafilters:

\begin{thm}
	Every Fr\'echet cardinal carries a Ketonen ultrafilter.\qed
\end{thm}

When \(\lambda\) is singular, it is important that the definition of a Ketonen
ultrafilter demands minimality only among uniform ultrafilters and not among the
broader class of tail uniform ultrafilters, since an ultrafilter on \(\lambda\)
that is minimal in this stronger sense is essentially the same thing as a
Ketonen ultrafilter on \(\text{cf}(\lambda)\):

\begin{lma}\label{KetonenOrdinals}
	Suppose \(\gamma\) is an ordinal and \(U\) is a \(\sE\)-minimal among
	countably complete ultrafilters \(W\) with \(\delta_W = \gamma\). Let
	\(\lambda = \textnormal{cf}(\gamma)\) and let \(f : \lambda\to \gamma\) be a
	continuous cofinal function. Then \(U = f_*(D)\) for some Ketonen
	ultrafilter \(D\) on \(\lambda\).
	\begin{proof}
		Since \(U\) is \(\sE\)-minimal among countably complete ultrafilters
		\(W\) with \(\delta_W = \gamma\), in particular \(U\) is
		\(\rRK\)-minimal, so every function \(g: \gamma\to \gamma\) that is
		regressive on a set in \(U\) is bounded on a set in \(U\). It follows
		that \(U\) contains every closed cofinal \(C\subseteq\gamma\): letting
		\(A = \gamma\setminus C\) and \(g(\alpha) = \sup (C\cap \alpha)\), \(g\)
		is regressive on \(A\) and unbounded on any cofinal subset of \(A\).
		
		Let \(C = f[\lambda]\). Then \(C \in U\). Let \(g : C\to \lambda\) be
		the inverse of \(f\). Let \(D = g_*(U)\). Clearly \(U = f_*(D)\). We
		must show that \(D\) is Ketonen. Suppose \(W\sE D\). We claim
		\(f_*(W)\sE U\). Given this, it follows that \(\delta_{f_*(W)} <
		\gamma\) and hence \(\delta_W < \lambda\). It follows that \(D\) is a
		\(\sE\)-minimal element of the set of countably complete uniform
		ultrafilters on \(\lambda\), so \(D\) is Ketonen.
		
		We finally verify \(f_*(W)\sE U\). (The proof will show that if \(f :
		\lambda\to \gamma\) is an order preserving function, then the
		pushforward map \(f_*\) is Ketonen order preserving.) Fix \(I\in D\) and
		\(\langle W_\alpha : \alpha\in I\rangle\) such that \(W =
		D\text{-}\lim_{\alpha\in I} W_\alpha\) and \(\delta_\alpha \leq \alpha\)
		for all \(\alpha \in I\). Let \(J = f[I]\), so \(J\in U\) and moreover:
		\[f_*(W) = U\text{-}\lim_{\beta\in J}f_*(W_{g(\beta)})\] Moreover
		\(\delta_{f_*(W_{g(\beta)})} \leq \sup f[\delta_{W_{g(\beta)}}] \leq
		\sup f[g(\beta)] \leq \beta\). Thus the sequence \(\langle
		f_*(W_{g(\beta)}) : \beta\in J\rangle\) witnesses \(f_*(W)\sE U\), as
		desired.
	\end{proof}
\end{lma}

\subsection{Introducing \(\mathscr K_\lambda\)}
Under the Ultrapower Axiom, the Ketonen order is linear, so there is a canonical
Ketonen ultrafilter on each Fr\'echet cardinal \(\lambda\):
\begin{defn}[UA]\index{\(\mathscr K_\lambda\)}\index{\(\mathscr K_\lambda\)|seealso{Ketonen ultrafilter}}\index{Ketonen ultrafilter!\(\mathscr K_\lambda\)}
For any Fr\'echet cardinal \(\lambda\), the {\it least ultrafilter on
\(\lambda\)}, denoted by \(\mathscr K_\lambda\), is the unique Ketonen
ultrafilter on \(\lambda\).
\end{defn}
The analysis of supercompactness under UA proceeds by first completely analyzing
the ultrafilters \(\mathscr K_\lambda\) and then propagating the structure of
\(\mathscr K_\lambda\) to all ultrafilters.

Let us begin with some simple examples. Let \(\kappa_0\) be the least measurable
cardinal. Then without assuming UA, it is easy to prove that an ultrafilter on
\(\kappa_0\) is Ketonen if and only if it is normal. Assuming UA, \(\mathscr
K_{\kappa_0}\) is the unique normal ultrafilter on \(\kappa_0\). 

Moving up to the second measurable cardinal \(\kappa_1\), it is {\it not}
provable in ZFC that the Ketonen ultrafilters on \(\kappa_1\) are normal, or
even that there is a normal Ketonen ultrafilter on \(\kappa_1\). This is because
it is consistent that \(\kappa_0\) is \(\kappa_1\)-strongly compact. Under this
assumption, if \(U\) is a normal ultrafilter on \(\kappa_1\), \(\kappa_0\) is
\(j_U(\kappa_1)\)-strongly compact in \(M_U\), and hence \(U\) concentrates on
ordinals that carry \(\kappa_0\)-complete uniform ultrafilters. In fact, under
this hypothesis, if \(W\) is a Ketonen ultrafilter on \(\kappa_1\), then \(j_W\)
is \((\kappa_1,\delta)\)-tight for some \(\delta < j_W(\kappa)\), and hence
witnesses the \(\kappa_1\)-strong compactness of \(\kappa_0\).

Of course, under UA, \(\kappa_0\) is not \(\kappa_1\)-strongly compact, since by
\cref{KunenUA}, every countably complete ultrafilter in \(V_{\kappa_1}\) is
isomorphic to \(\mathscr K_{\kappa_0}^n\) for some \(n < \omega\). In fact, once
again \(\mathscr K_{\kappa_1}\) is the unique normal ultrafilter on
\(\kappa_1\). To see this, one can apply \cref{MuDichotomy} and the following
lemma:
\begin{lma}[UA]\label{KetonenIrreducible}
For any regular cardinal \(\lambda\), \(\mathscr K_\lambda\) is an irreducible
ultrafilter.\index{Ketonen ultrafilter!\(\mathscr K_\lambda\)!irreducibility}
\begin{proof}
Suppose \(D\sD \mathscr K_\lambda\). Then since \(D\sRK \mathscr K_\lambda\) and
\(\mathscr K_\lambda\) is weakly normal, \(\lambda_D < \lambda\). Therefore by
\cref{UFContinuity}, \[j_D(\lambda) = \sup j_D[\lambda]\] Assume towards a
contradiction that \(D\) is nonprincipal. Then by \cref{Pushdown}, \(\tr D
{\mathscr K_\lambda}\sE j_D(\mathscr K_\lambda)\), so \(\delta_{\tr D {\mathscr
K_\lambda}} < j_D(\lambda)\) by \cref{KetonenMinimality} applied in \(M_D\). But
\(\mathscr K_\lambda = j_D^{-1}[\tr D {\mathscr K_\lambda}]\), so
\[\delta_{\mathscr K_\lambda} = \min \{\delta : j_D(\delta) > \delta_{\tr D
{\mathscr K_\lambda}}\} < \lambda\] This contradicts that \(\mathscr K_\lambda\)
is a uniform ultrafilter on \(\lambda\).
\end{proof}
\end{lma}
We do not know whether this lemma is provable in ZFC, although it does follow
from \cref{ACC}. 

If \(\lambda\) is singular, then \(\mathscr K_\lambda\) is not necessarily
irreducible. (In fact, we will show under UA that for strong limit singular
cardinals  \(\lambda\), \(\mathscr K_\lambda\) is {\it never} irreducible.) For
example, suppose \(\lambda_0\) is the least singular cardinal that carries a
uniform countably complete ultrafilter. Of course, assuming just ZFC, one cannot
prove much about \(\lambda_0\): it is consistent that \(\lambda_0 =
\kappa_0^{+\kappa_0}\), or that \(\lambda_0\) is not a limit of regular
cardinals that carry uniform countably complete ultrafilters.

Assuming UA, it is not hard to give a complete analysis of \(\lambda_0\) and
\(\mathscr K_{\lambda_0}\). 
Let \(\langle\kappa_\alpha: \alpha < \kappa_0\rangle\) enumerate the first
\(\kappa_0\) measurable cardinals in increasing order. Then \(\lambda_0 =
\sup_{\alpha < \kappa_0} \kappa_\alpha\), and \[\mathscr K_{\lambda_0} =
\mathscr K_{\kappa_0}\text{-}\lim_{\alpha < \kappa_0} \mathscr
K_{\kappa_\alpha}\restriction \lambda_0\] The sets \(A_\alpha = \kappa_\alpha
\setminus \sup_{\beta < \alpha}\kappa_\beta\) witness that the sequence
\(\langle \mathscr K_{\kappa_\alpha}\restriction \lambda_0:\alpha <
\kappa_0\rangle\) is discrete, so \(\mathscr K_{\kappa_0}\sD \mathscr
K_{\lambda_0}\). In other words, \(\mathscr K_{\lambda_0}\) is produced by the
iterated ultrapower \(\langle \mathscr K_\kappa, \mathscr K^{M_{\mathscr
K_\kappa}}_{\lambda_0}\rangle\). 

Of course all this is closely related to \cref{FrechetSingularChar}. For
singular cardinals \(\lambda\), \(\mathscr K_\lambda\) is of greatest interest
if \(\lambda\) is not a limit of Fr\'echet cardinals, since in this case
\(\mathscr K_\lambda\) cannot be represented in terms of ultrafilters on smaller
cardinals.

\subsection{The universal property of \(\mathscr K_\lambda\)} \label{UniversalSection}
The main result of this section is a universal property of the least ultrafilter
\(\mathscr K_\lambda\) on a regular Fr\'echet cardinal:
\begin{thm}[UA]\label{UniversalProperty}\index{Ketonen ultrafilter!\(\mathscr K_\lambda\)!universal property}
Suppose \(\lambda\) is a regular Fr\'echet cardinal. Let \(j :V\to M\) be the
ultrapower of the universe by \(\mathscr K_\lambda\). Suppose \(i : V\to N\) is
an ultrapower embedding. Then the following are equivalent:
\begin{enumerate}[(1)]
\item There is an internal ultrapower embedding \(k : M\to N\) such that
\(k\circ j = i\).
\item \(\sup i[\lambda]\) carries no tail uniform ultrafilter in \(N\).
\item \(\textnormal{cf}^{N}(\sup i[\lambda])\) is not Fr\'echet in \(N\).
\end{enumerate}
\end{thm}

While the proof is quite simple, the result has profound consequences for the
structure of the ultrafilters \(\mathscr K_\lambda\). In fact, this universal
property is ultimately responsible for all of our results on supercompactness
under UA. 

Before proving \cref{UniversalProperty} (which is not that difficult), let us
show how it can be used to give a complete analysis of the internal ultrapower
embeddings of \(M_{\mathscr K_\lambda}\) when \(\lambda\) is regular.
\begin{thm}[UA]\label{EmbeddingChar}\index{Ketonen ultrafilter!\(\mathscr K_\lambda\)!internal ultrapowers}
Suppose \(\lambda\) is a regular Fr\'echet cardinal.  Let \(j :V \to M\) be the
ultrapower of the universe by \(\mathscr K_\lambda\). Suppose \(k : M  \to N\)
is an ultrapower embedding. Then the following are equivalent:
\begin{enumerate}[(1)]
\item \(k\) is an internal ultrapower embedding.
\item \(k\) is continuous at \(\sup j[\lambda]\).
\item \(k\) is continuous at \(\textnormal{cf}^M(\sup j[\lambda])\).
\end{enumerate}
\begin{proof}
{\it (1) implies (2):} Since \(\sup j[\lambda]\) carries no tail uniform
countably complete ultrafilter in \(M\), every elementary embedding of \(M\)
that is close to \(M\) is continuous at \(\sup j[\lambda]\). In particularly,
every internal ultrapower embedding of \(M\) is continuous at \(\sup
j[\lambda]\).

{\it (2) implies (1):} Let \(i = k\circ j\). Then \(\sup i[\lambda] = \sup
k[\sup j[\lambda]] = k(\sup j[\lambda])\) since \(k\) is continuous at \(\sup
j[\lambda]\). It follows that \(\sup i[\lambda]\) carries no tail uniform
countably complete ultrafilter in \(N\). Therefore by \cref{UniversalProperty},
there is an internal ultrapower embedding \(k' : M\to N\) such that \(k'\circ j
= i\). 

We claim \(k' = k\). First of all, \(k'\circ j = k\circ j\). In other words,
\(k'\restriction j[V] = k\restriction j[V]\). Moreover since \(k'\) is
\(M\)-internal \(k'(\sup j[\lambda]) = \sup i[\lambda] = k(\sup j[\lambda])\).
But \(M = H^M(j[V]\cup \{\id_{\mathscr K_\lambda}\}) = H^M(j[V]\cup \{\sup
j[\lambda]\})\) since \(\mathscr K_\lambda\) is weakly normal. Since we have
shown \(k'\restriction j[V]\cup \{\sup j[\lambda]\} =k\restriction j[V]\cup
\{\sup j[\lambda]\}\), it follows that \(k' = k\).

Since \(k'\) is an internal ultrapower embedding, so is \(k\), as desired.

The equivalence of (2) and (3) is trivial (and does not require UA).
\end{proof}
\end{thm}

The notion of {\it indecomposable ultrafilters} is an important part of infinite
combinatorics. We will need the following relativized version of this concept:
\begin{defn}\index{Indecomposable ultrafilter}\index{Decomposable ultrafilter|seealso{Indecomposable ultrafilter}}
	Suppose \(M\) is a transitive model of ZFC and \(U\) is an \(M\)-ultrafilter
	on \(X\). Suppose \(\delta\) is an \(M\)-cardinal. Then \(U\) is {\it
	\(\delta\)-indecomposable} if for any partition \(\langle X_\alpha : \alpha
	< \delta\rangle\in M\) of \(X\), there is some \(S\subseteq \delta\) in
	\(M\) with \(|S|^M < \delta\) and \(\bigcup_{\alpha\in S} X_\alpha\in U\).
\end{defn}

As a corollary of \cref{EmbeddingChar}, every \(\lambda\)-indecomposable
ultrafilter is internal to \(\mathscr K_\lambda\):

\begin{cor}[UA] Suppose \(\lambda\) is a regular Fr\'echet cardinal. Suppose
	\(D\) is a countably complete \(\lambda\)-indecomposable ultrafilter, then
	\(D\I \mathscr K_\lambda\). In particular, if \(D\) is a countably complete
	ultrafilter such that \(\lambda_D < \lambda\), then \(D\I \mathscr
	K_\lambda\).
	\begin{proof}
		Let \(j :  V\to M\) be the ultrapower of the universe by \(\mathscr
		K_\lambda\). To show \(D\I \mathscr K_\lambda\), we need to show that
		\(j_D\restriction M\) is an internal ultrapower embedding of \(M_U\). By
		\cref{PushUlt}, \(j_D\restriction M\) is an ultrapower embedding. Since
		\(D\) is \(\lambda\)-indecomposable, \(j_D\) is continuous at all
		ordinals of cofinality \(\lambda\), and in particular, \(j_D\) is
		continuous at \(\sup j[\lambda]\). Thus \(j_D\restriction M\) is an
		ultrapower embedding of \(M\) that is continuous at \(\sup j[\lambda]\),
		and it follows from \cref{EmbeddingChar} that \(j_D\restriction M\) is
		an internal ultrapower embedding of \(M\), as desired.
	\end{proof}
\end{cor}
 This fact is highly reminiscent of \cref{LipMO}, the theorem that analyzes
 which ultrafilters lie Mitchell below a Dodd solid ultrafilter. In fact, we
 will  show that \(\mathscr K_\lambda\) gives rise to a supercompact ultrapower
 precisely by leveraging the fact that so many ultrapower embeddings are
 internal to it. (See \cref{IndependentSection}, \cref{KSuperSection}, and
 especially \cref{GeneralStrong}.)

If \(M\) is a transitive model of ZFC and \(\delta\) is an \(M\)-regular
cardinal, then an \(M\)-ultrafilter \(U\) is \(\delta\)-indecomposable if and
only if \(j_U^M\) is continuous at \(\delta\), and we therefore obtain the
following combinatorial characterization of the countably complete
\(M\)-ultrafilters that belong to \(M\) when \(M\) is the ultrapower of the
universe by a Ketonen ultrafilter on a regular cardinal:

\begin{thm}[UA]\label{UFAmenableChar}
Suppose \(\lambda\) is a regular Fr\'echet cardinal. Let \(j : V\to M\) be the
ultrapower of the universe by \(\mathscr K_\lambda\). Let \(\delta =
\textnormal{cf}^M(\sup j[\lambda])\). Suppose \(U\) is a countably complete
\(M\)-ultrafilter. Then the following are equivalent:
\begin{enumerate}[(1)]
	\item \(U\) is \(\delta\)-indecomposable. 
	\item \(U\in M\).
\end{enumerate}
In particular, if \(U\) is a countably complete \(M\)-ultrafilter on a cardinal
\(\gamma < \delta\), then \(U\in M\).\qed
\end{thm}

In summary, the universal property of \(\mathscr K_\lambda\) is a powerful tool
for analyzing the model \(M_{\mathscr K_\lambda}\). Let us therefore prove it:
\begin{proof}[Proof of \cref{UniversalProperty}] {\it (1) implies (2):} First,
\(k(\sup j[\lambda])\) carries no tail uniform countably complete ultrafilter in
\(N\) by elementarity, since \(\sup j[\lambda]\) carries no tail uniform
countably complete ultrafilter in \(M\). Note also that \(k : M\to N\) is
continuous at \(\sup j[\lambda]\) since \(\sup j[\lambda]\) carries no tail
uniform countably complete ultrafilter in \(M\). Therefore \(k(\sup j[\lambda])
= \sup k\circ j[\lambda] = \sup i[\lambda]\). Hence \(\sup i[\lambda]\) carries
no tail uniform countably complete ultrafilter in \(N\).

{\it (2) implies (1):} Let \((e,h) : (M,N)\to P\) be an internal ultrapower
comparison of \((j,i)\). Then 
\[e(\sup j[\lambda]) = \sup e\circ j[\lambda] = \sup h\circ i[\lambda] = h(\sup
i[\lambda])\] The theorem is now an immediate consequence of \cref{RFSEquiv}:
\(M = H^M(j[V]\cup \{\sup j[\lambda]\})\) and \((e,h)\) witnesses \((j,\sup
j[\lambda]) =_S (i,\sup i[\lambda])\), so there is an internal ultrapower
embedding \(k : M\to N\) such that \(k\circ j = i\).

The equivalence of (2) and (3) is trivial (and does not require UA).
\end{proof}

\subsection{Independent families and the Hamkins properties}\label{IndependentSection}\index{Hamkins Properties}\index{Covering property}\index{Approximation property}
A key intuition in the theory of forcing is that forcing does not create new
large cardinals. The Levy-Solovay Theorem \cite{LevySolovay} establishes this
for small forcing, but various counterintuitive forcing constructions of the
next few decades demonstrate that in general, the intuition is just not correct.
The earliest example, due to Kunen, shows that it is consistent that there is a
forcing that makes a measurable cardinal out of a cardinal that is not even
weakly compact. Woodin's \(\Sigma_2\)-Resurrection Theorem (\cite{Larson},
Theorem 2.5.10) yields even more striking examples: for example, if there is a
proper class of Woodin cardinals and there is a huge cardinal, then arbitrarily
large cardinals can be forced to be huge cardinals.

 Hamkins isolated two closure properties of inner models: the {\it approximation
 and covering properties}, or collectively the {\it Hamkins properties}. If an
 inner model \(M\) has the Hamkins properties, then many of the large cardinal
 properties of the ambient universe of sets are downwards absolute to \(M\). For
 many forcing extensions \(V[G]\), the universe \(V\) satisfies the Hamkins
 properties inside \(V[G]\), and therefore the large cardinals of \(V[G]\)
 ``already exist" in \(V\).

Somewhat unexpectedly, the Hamkins properties have turned out to be relevant
outside forcing, in the domain of inner model theory. Woodin has shown that any
inner model that inherits a supercompact cardinal \(\kappa\) from the ambient
universe in a natural way necessarily satisfies the Hamkins properties at
\(\kappa\), and therefore inherits all large cardinals from the ambient
universe. Such models are called {\it weak extender models} for the
supercompactness of \(\kappa\). A canonical inner model with a  supercompact
cardinal is expected to be a weak extender model, and therefore Woodin
conjectures that if there is a canonical inner model with a supercompact
cardinal, in fact this is the {\it ultimate inner model}, a canonical inner
model that satisfies all true large cardinal axioms.

In our work on UA, the Hamkins properties rear their heads once again. Here they
arise in relation with (generalizations of) the Mitchell order, which can be
seen as yet another instantiation of the downwards absoluteness of large
cardinal properties to inner models. Recall that we are trying to show that the
ultrapower of the universe by \(\mathscr K_\lambda\) has closure properties. All
we know so far is that this ultrapower absorbs many countably complete
ultrafilters (\cref{UFAmenableChar}). To transform this into a model theoretic
closure property of the ultrapower, for example closure under
\(\lambda\)-sequences, we prove a converse to Hamkins and Woodin's absoluteness
theorems for models with the Hamkins properties. This converse says that any
inner model that inherits enough ultrafilters from the ambient universe must
satisfy the Hamkins properties. In our context, this will lead to a proof that
the ultrapower \(\mathscr K_\lambda\) is (roughly) closed under
\(\lambda\)-sequences.

The ultrapowers we consider do not satisfy the (relevant) Hamkins properties in
full, but rather satisfy local versions of these properties, introduced here for
the first time:
\begin{defn}
Suppose \(M\) is an inner model, \(\kappa\) is a cardinal, and \(\lambda\) is an
ordinal.
\begin{itemize}
\item \(M\) has the {\it \(\kappa\)-covering property at \(\lambda\)} if every
\(\sigma\in P_\kappa(\lambda)\) there is some \(\tau\in P_\kappa(\lambda)\cap
M\) with \(\sigma\subseteq \tau\).
\item \(M\) has the {\it \(\kappa\)-approximation property at \(\lambda\)} if
any \(A\subseteq \lambda\) with \(A\cap \sigma\in M\) for all \(\sigma\in
P_\kappa(\lambda)\cap M\) is an element of \(M\).
\end{itemize}
We say \(M\) has the {\it \(\kappa\)-covering property} if \(M\) has the
\(\kappa\)-covering property at all \(M\)-cardinals, and \(M\) has the {\it
\(\kappa\)-approximation property} if \(M\) has the \(\kappa\)-approximation
property at all \(M\) cardinals.
\end{defn}
In this section, we identify necessary and sufficient conditions for the
\(\kappa\)-covering and approximation properties that involve the absorption of
filters. We are working in slightly more generality than we will need, but we
think the results are quite interesting and hopefully lead to a clearer
exposition than would arise by working in a more specific case.

The condition equivalent to the covering property essentially comes from
Woodin's proof of the covering property for weak extender models:
\begin{prp}
Suppose \(M\) is an inner model. Then \(M\) has the \(\kappa\)-covering property
at \(\lambda\) if and only if there is a \(\kappa\)-complete fine filter on
\(P_\kappa(\lambda)\) that concentrates on \(M\).
\begin{proof}
First assume there is a \(\kappa\)-complete fine filter \(\mathcal F\) on
\(P_\kappa(\lambda)\) that concentrates on \(M\). Fix \(\sigma\in
P_\kappa(\lambda)\), and we will find \(\tau\in P_\kappa(\lambda)\cap M\) such
that \(\sigma\subseteq \tau\). For each \(\alpha < \lambda\), let \(A_\alpha =
\{\tau\in P_\kappa(\lambda) : \alpha\in \tau\}\), so that \(A_\alpha\in \mathcal
F\) by the definition of a fine filter. Then suppose \(\sigma\in
P_\kappa(\lambda)\). The set \[\{\tau\in P_\kappa(\lambda) : \sigma\subseteq
\tau\} = \bigcap_{\alpha\in \sigma} \{A_\alpha : \alpha\in \sigma\}\in \mathcal
F\] since \(\mathcal F\) is \(\kappa\)-complete. Since \(\mathcal F\)
concentrates on \(M\), \(\{\tau\in P_\kappa(\lambda) : \sigma\subseteq
\tau\}\cap M\in \mathcal F\), and in particular this set is nonempty. Any
\(\tau\) that belongs to this set satisfies \(\tau\in P_\kappa(\lambda)\cap M\)
and \(\sigma\subseteq \tau\), as desired.

Conversely, assume \(M\) has the \(\kappa\)-covering property at \(\lambda\).
Let \({\mathcal B} = \{A_\alpha \cap M: \alpha < \lambda\}\). Then \({\mathcal
B}\) is a \(\kappa\)-complete filter base: for any \(S\subseteq {\mathcal B}\)
with \(|S| < \kappa\), we have \(S = \{A_\alpha \cap M : \alpha\in \sigma\}\)
for some \(\sigma\in P_\kappa(\lambda)\), and so fixing \(\tau\in
P_\kappa(\lambda)\cap M\) such that \(\sigma\subseteq \tau\), we have \(\tau\in
\bigcap_{\alpha\in \sigma} A_\alpha\). Therefore \({\mathcal B}\) extends to a
\(\kappa\)-complete filter \(\mathcal G\). Let \[\mathcal F = \mathcal
G\restriction P_\kappa(\lambda) = \{A\subseteq P_\kappa(\lambda) : A\cap M\in
\mathcal G\}\] be the canonical extension of \(\mathcal G\) to an filter on
\(P_\kappa(\lambda)\). Then  \(\mathcal F\) is \(\kappa\)-complete and
concentrates on \(M\). Moreover, \(A_\alpha\in \mathcal F\) for all \(\alpha <
\lambda\), so \(\mathcal F\) is fine. Thus we have produced a
\(\kappa\)-complete fine filter on \(P_\kappa(\lambda)\) that concentrates on
\(M\), as desired.
\end{proof}
\end{prp}

One ultrafilter theoretic characterization of the approximation property is
given by the following theorem:

\begin{thm}
Suppose \(\kappa\) is strongly compact and \(M\) is an inner model with the
\(\kappa\)-covering property. Then \(M\) has the \(\kappa\)-approximation
property if and only if every \(\kappa\)-complete \(M\)-ultrafilter belongs to
\(M\).
\end{thm}

We will actually prove a local version of this theorem that requires no large
cardinal assumptions. The locality of this theorem is important in our analysis
of the ultrafilters \(\mathscr K_\lambda\). For the statement, we need use the
following definition:
\begin{defn}
Suppose \(X\) is a set and \(\Sigma\) is an algebra of subsets of \(X\). A set
\(U\subseteq \Sigma\) is said to be an {\it ultrafilter over \(\Sigma\)} if
\(U\) is closed under intersections and for any \(A\in \Sigma\), \(A\in U\) if
and only if \(X\setminus A\notin U\). An ultrafilter \(U\) over \(\Sigma\) is
said to be \(\kappa\)-complete if for any \(\sigma\in P_\kappa(U)\), \(\bigcap
\sigma\neq \emptyset\).
\end{defn}

What we call an ultrafilter over \(\Sigma\) is commonly referred to as an {\it
ultrafilter on the Boolean algebra \(\Sigma\),} but we are being a bit pedantic:
we do not want to confuse this with an ultrafilter with underlying set
\(\Sigma\), which in our terminology is a family of subsets of \(\Sigma\) rather
than a subset of \(\Sigma\). Notice that for us a \(\kappa\)-complete
ultrafilter over \(\Sigma\) is the same thing as an ultrafilter over \(\Sigma\)
that is a \(\kappa\)-complete filter base. (It is not the same thing as being
\(\kappa\)-complete ultrafilter on the Boolean algebra \(\Sigma\).)

\begin{thm}\label{ApproxThm}
Suppose \(M\) is an inner model, \(\kappa\) is a cardinal, \(\lambda\) is an
\(M\)-cardinal, and \(M\) has the \(\kappa\)-covering property at \(\lambda\).
Then the following are equivalent:
\begin{enumerate}[(1)]
\item \(M\) has the \(\kappa\)-approximation property at \(\lambda\).
\item Suppose \(\Sigma\in M\) is an algebra of sets of \(M\)-cardinality
\(\lambda\). Then every \(\kappa\)-complete ultrafilter over \(\Sigma\) belongs
to \(M\).
\end{enumerate}
\end{thm}

To simplify notation, we use the following lemma (analogous to \cref{XTight})
characterizing the approximation property at \(\lambda\):
\begin{lma}\label{XApprox}
Suppose \(M\) is an inner model, \(\kappa\) is a cardinal, and \(\lambda\) is an
\(M\)-cardinal. Then the following are equivalent:
\begin{enumerate}[(1)]
\item \(M\) has the \(\kappa\)-approximation property at \(\lambda\)
\item For all \(\Sigma\in M\) such that \(|\Sigma|^M \leq \lambda\), for all
\(B\subseteq \Sigma\) such that \(B\cap \sigma\in M\) for all \(\sigma\in
P_\kappa(\Sigma)\cap M\), \(B\in M\).
\item For some \(\Sigma\in M\) such that \(|\Sigma|^M = \lambda\), for all
\(B\subseteq \Sigma\) such that \(B\cap \sigma\in M\) for all \(\sigma\in
P_\kappa(\Sigma)\cap M\), \(B\in M\). \qed
\end{enumerate}
\end{lma}

The following notation will be convenient (although of course it is a bit
ambiguous):
\begin{defn}
Suppose \(X\) is a set and \(\sigma\) is a family of subsets of \(X\). Then the
{\it dual of \(\sigma\) in \(X\)} is the family \(\sigma^* = \{X\setminus A :
A\in \sigma\}\).
\end{defn}
We point out that the dualizing operation depends implicitly on the underlying
set \(X\).
\begin{defn}\index{Independent family}
Suppose \(\kappa\) is a cardinal and \(X\) is a set. A family \(\Gamma\) of
subsets of \(X\) is {\it \(\kappa\)-independent} if for any disjoint sets
\(\tau_0,\tau_1\in P_\kappa(\Gamma)\), \(\bigcap \tau_0\cap \bigcap \tau_1^*\neq
\emptyset\).
\end{defn}
Equivalently, \(\Gamma\) is \(\kappa\)-independent if for any disjoint sets
\(X,Y\subseteq \Gamma\), the collection \(X\cup Y^*\) is a \(\kappa\)-complete
filter base. Note that a \(\kappa\)-complete family of subsets of \(X\) is never
an algebra of sets, since if \(A\in \Gamma\), then \(X\setminus A\notin
\Gamma\).

\begin{thm}[Hausdorff]\label{Hausdorff0}
Suppose \(\kappa\) and \(\lambda\) are cardinals. Then there is a
\(\kappa\)-independent family of subsets of \(X = \{(\sigma,s):\sigma\in
P_\kappa(\lambda)\text{ and }s\in P_\kappa(P(\sigma))\}\) of cardinality
\(2^\lambda\).
\begin{proof}
Define \(f : P(\lambda)\to P(X)\) by \[f(A) = \{(\sigma,s) \in X: \sigma\cap
A\in s\}\] Let \(\Gamma = \text{ran}(f)\). Suppose \(\tau_0,\tau_1\in
P_\kappa(P(\lambda))\) are disjoint. We claim that the set
\[S = \bigcap f[\tau_0] \cap \bigcap f[\tau_1]^*\] is nonempty. This
simultaneously shows that \(f\) is injective and \(\Gamma\) is
\(\kappa\)-independent. Therefore \(\Gamma\) is a \(\kappa\)-independent family
of cardinality \(2^\lambda\).

Let \(\sigma\in P_\kappa(\lambda)\) be large enough that \(\sigma\cap A_0\neq
\sigma\cap A_1\) for any \(A_0\in \tau_0\) and \(A_1\in \tau_1\). Let \[s =
\{\sigma\cap A : A\in \tau_0\}\] We claim that \((\sigma, s) \in S\).

First we show that \((\sigma,s)\in \bigcap f[\tau_0]\). Suppose \(A\in \tau_0\).
We will show that \((\sigma,s)\in f(A)\). By the definition of \(s\), since
\(A\in \tau_0\), \(\sigma\cap A\in s\). Therefore by the definition of \(f\),
\((\sigma,s)\in f(A)\), as desired. This shows \((\sigma,s)\in \bigcap
f[\tau_0]\).

Next we show that \((\sigma,s)\in \bigcap f[\tau_1]^*\). Suppose \(B\in
\tau_1\), and we will show that \((\sigma,s)\in X\setminus B\). By the choice of
\(\sigma\), \(\sigma\cap B\neq \sigma \cap A\) for any \(A\in \tau_0\).
Therefore by the definition of \(s\), \(\sigma\cap B\notin s\). Finally, by the
definition of \(f\), it follows that \((\sigma,s)\notin f(B)\), or in other
words, \((\sigma,s)\in X\setminus B\). Hence \((\sigma,s)\in \bigcap
f[\tau_1]^*\).

Since \((\sigma,s)\in \bigcap f[\tau_0]\) and \((\sigma,s)\in \bigcap
f[\tau_1]^*\), it follows that \((\sigma,s) \in S\). Thus \(S\) is nonempty,
which completes the proof.
\end{proof}
\end{thm}

Computing cardinalities, Hausdorff's theorem implies the existence of
\(\kappa\)-independent sets that are as large as possible:

\begin{cor}[Hausdorff]\label{Hausdorff1}
Suppose \(\kappa\) and \(\lambda\) are cardinals such that \(\lambda^{<\kappa} =
\lambda\). Then there is a \(\kappa\)-independent family of subsets of
\(\lambda\) of cardinality \(2^\lambda\).
\begin{proof}
Let \(X= \{(\sigma,s) : \sigma\in P_\kappa(\lambda)\text{ and }s\in
P_\kappa(P(\sigma))\}\). In other words, \[X = \coprod_{\sigma\in
P_\kappa(\lambda)} P_\kappa(P(\sigma))\] Thus \[|X| = |P_\kappa(\lambda)|\cdot
\sup_{\sigma\in P_\kappa(\lambda)} |P_\kappa(P(\sigma))| = \lambda^{<\kappa}
\cdot (2^{<\kappa})^{<\kappa} =  \lambda^{<\kappa} = \lambda\] By
\cref{Hausdorff0}, there is a \(\kappa\)-independent family of subsets of \(X\)
of cardinality \(2^\lambda\), and therefore there is a \(\kappa\)-independent
family of subsets of \(\lambda\) of cardinality \(2^\lambda\).
\end{proof}
\end{cor}
We now establish our characterization of the approximation property.
\begin{proof}[Proof of \cref{ApproxThm}] {\it (1) implies (2):} Assume (1), and
we will prove (2). Suppose \(\Sigma\in M\) is an algebra of subsets of \(X\) of
\(M\)-cardinality \(\lambda\) and \(U\) is a \(\kappa\)-complete ultrafilter
over \(\Sigma\). Fix \(\sigma\in P_\kappa(\Sigma)\cap M\) and we will show that
\(\sigma\cap U\in M\). Since \(U\) is \(\kappa\)-complete, \[S = \bigcap \{A :
A\in \sigma\cap U\}\cap \bigcap \{X\setminus A:A\in \sigma\setminus U\}\] is
nonempty. Therefore fix \(a\in X\) with \(a\in S\). By the choice of \(a\),
\(\sigma \cap U = \{A\in \sigma : a\in A\}\). Thus \(\sigma\cap U\in M\).

By the \(\kappa\)-approximation property at \(\lambda\) (using \cref{XApprox}),
it follows that \(U\in M\).

{\it (2) implies (1):} Fix \(\Gamma\in M\) such that \(M\) satisfies that
\(\Gamma\) is a \(\kappa\)-independent family of subsets of some set \(X\) and
\(|\Gamma|^M = \lambda\). Suppose \(C\subseteq \Gamma\) is such that \(C\cap
\sigma\in M\) for all \(\sigma\in P_\kappa(\Gamma)\cap M\). We will show that
\(C\in M\). This verifies the condition of \cref{XApprox} (3), and so implies
that \(M\) satisfies the \(\kappa\)-approximation property at \(\lambda\).

Let \[\mathcal B = C\cup (\Gamma\setminus C)^*\] We claim that \(\mathcal B\) is
a \(\kappa\)-complete filter base on \(X\). Suppose \(\sigma\in
P_\kappa(\mathcal B)\). We must show that \(\bigcap\sigma\neq \emptyset\). Using
the \(\kappa\)-covering property at \(\lambda\), fix \(\tau\in
P_\kappa(\Gamma)\cap M\) such that \(\sigma\subseteq \tau\cup \tau^*\).

 By our assumption on \(C\), \(\tau\cap C\in M\). Let \(\tau_0 = \tau\cap C\)
and let \(\tau_1 = \tau \setminus C = \tau \setminus \tau_0 \in M\). Since
\(\sigma\subseteq \mathcal B = C\cup (\Sigma\setminus C)^*\), we have \(\sigma
\subseteq \tau_0\cup \tau_1^*\). Since \(\Gamma\) is \(\kappa\)-independent in
\(M\),
\[\textstyle\bigcap\tau_0\cap\bigcap \tau_1^*\neq \emptyset\] But
\(\textstyle\bigcap\tau_0\cap\bigcap \tau_1^* =\bigcap (\tau_0\cup
\tau_1)\subseteq \bigcap \sigma\), and hence \(\bigcap \sigma\neq\emptyset\), as
desired. This shows \(\mathcal B\) is a \(\kappa\)-complete filter base.

Let \(\Sigma\) be the algebra on \(X\) generated by \(\Gamma\) and let \(U\) be
the ultrafilter over \(\Sigma\) generated by \(\mathcal B\). Then \(U\) is
\(\kappa\)-complete because \(\mathcal B\) is \(\kappa\)-complete. Therefore
\(U\in M\) by our assumption on \(M\). But \(C = \Gamma\cap \mathcal B =
\Gamma\cap U\), so \(C\in M\), as desired. Thus \(M\) has the
\(\kappa\)-approximation property at \(\lambda\).
\end{proof}

The proof of \cref{ApproxThm} has the following corollary, which will be
important going forward:
\begin{prp}\label{ApproxPrp}
Suppose \(M\) is an inner model, \(\kappa\) is a cardinal, \(\lambda\) is an
\(M\)-cardinal, and \(M\) has the \(\kappa\)-covering property at \(\lambda\).
Then the following are equivalent:
\begin{enumerate}[(1)]
\item \(M\) has the \(\kappa\)-approximation property at \(\lambda\).
\item There is a \(\kappa\)-independent family \(\Gamma\) of \(M\) with
\(M\)-cardinality \(\lambda\) such that every \(\kappa\)-complete ultrafilter
over the algebra generated by \(\Gamma\) belongs to \(M\).\qed
\end{enumerate}
\end{prp}
\subsection{The strength and supercompactness of \(\mathscr K_\lambda\)}\label{KSuperSection}
\begin{defn}
For any Fr\'echet cardinal \(\lambda\), \(\kappa_\lambda\) denotes the
completeness of \(\mathscr K_\lambda\).\index{\(\kappa_\lambda\) (completeness
of \(\mathscr K_\lambda\))}
\end{defn}
In other words, \(\kappa_\lambda = \textsc{crt}(j_{\mathscr K_\lambda})\). In
\cref{FrechetSection}, we will prove the following theorem:
\begin{thm}[UA] Suppose \(\lambda\) is a Fr\'echet cardinal that is either a
successor cardinal or a strongly inaccessible cardinal. Then \(\kappa_\lambda\)
is \(\lambda\)-strongly compact.
\end{thm}
This is one of the harder theorems of this chapter, so we will just work under
this hypothesis for a while. The following theorem begins to show why it is a
useful assumption:
\begin{thm}[UA]\label{KetonenStrong}
Suppose \(\lambda\) is a regular Fr\'echet cardinal and \(\kappa_\lambda\) is
\(\lambda\)-strongly compact. Let \(j : V\to M\) be the ultrapower of the
universe by \(\mathscr K_\lambda\). Then \(P(\gamma)\subseteq M\) for all
\(\gamma < \lambda\).
\end{thm}
Because we will occasionally need to use this argument in a more general
context, let us instead prove the following:

\begin{prp}\label{GeneralStrong}\index{Ketonen ultrafilter!\(\mathscr K_\lambda\)!supercompactness}
Suppose \(\kappa \leq \gamma\) are cardinals, \(\kappa\) is \(\gamma\)-strongly
compact, and \(M\) is an inner model that is closed under
\({<}\kappa\)-sequences. Assume every \(\kappa\)-complete ultrafilter on
\(\gamma\) is amenable to \(M\). Then \(P(\gamma)\subseteq M\). Moreover if
\(\textnormal{cf}(\gamma)\geq \kappa\) then \(P(\eta)\subseteq M\) for all
\(\eta\leq 2^\gamma\) such that \(\kappa\) is \(\eta\)-strongly compact.
\end{prp}
\begin{proof}
We may assume by induction that \(P(\alpha)\subseteq M\) for all ordinals
\(\alpha < \gamma\). Let \(\nu = \text{cf}(\gamma)\). 

Assume first that \(\nu < \kappa\). Let \(\langle \gamma_\alpha : \alpha <
\nu\rangle\in M\) be an increasing sequence cofinal in \(\gamma\). Suppose
\(A\subseteq \gamma\). Let \(A_\alpha = A\cap \gamma_\alpha\), so \(A_\alpha\in
M\) for all \(\alpha < \nu\) by our inductive assumption. Then \(\langle
A_\alpha : \alpha < \nu\rangle\in M\) since \(M\) is closed under
\({<}\kappa\)-sequences. Therefore \(A = \bigcup_{\alpha < \nu} A_\alpha\in M\).
It follows that \(P(\gamma)\subseteq M\), which finishes the proof in this case.

Therefore we may assume that \(\nu\geq \kappa\). We claim that \(\kappa\) is
\(\gamma\)-strongly compact in \(M\). Fix an ordinal \(\alpha \in
[\kappa,\gamma]\) such that \(\text{cf}^M(\alpha)\geq \kappa\). Then
\(\text{cf}(\alpha)\geq \kappa\) since \(M\) is closed under
\(\kappa\)-sequences. Since \(\kappa\) is \(\gamma\)-strongly compact, there is
a \(\kappa\)-complete tail uniform ultrafilter \(U\) on \(\alpha\). But \(U\cap
M\in M\), so in \(M\) there is a tail uniform \(\kappa\)-complete ultrafilter on
\(\alpha\). In particular, every \(M\)-regular cardinal \(\iota\in
[\kappa,\lambda]\) carries a \(\kappa\)-complete ultrafilter in \(M\), so by
\cref{KetonenThm}, \(\kappa\) is \(\gamma\)-strongly compact in \(M\).

Therefore by \cref{SolovaySCH}, \((\gamma^{<\kappa})^M = \gamma\), so by
\cref{Hausdorff1}, \(M\) satisfies that there is a \(\kappa\)-independent family
of subsets of \(\gamma\) of cardinality \((2^\gamma)^M\). 

Let \(\Gamma\in M\) be such that \(M\) satisfies that \(\Gamma\) is a
\(\kappa\)-independent family of subsets of \(\gamma\) of cardinality
\(\gamma\).  Let \(\Sigma\) be the algebra of subsets of \(\gamma\) generated by
\(\Gamma\). If \(U_0\) is a \(\kappa\)-complete ultrafilter over \(\Sigma\),
then \(U_0\) extends to a \(\kappa\)-complete ultrafilter \(U\) on \(\gamma\) by
\cref{StrongCChar}, since \(\kappa\) is \(\gamma\)-strongly compact and \(U_0\)
is a \(\kappa\)-complete filter base of cardinality \(\gamma\). It follows from
\cref{ApproxPrp} that \(M\) has the \(\kappa\)-approximation property at
\(\gamma\). Since \(M\) is closed under \({<}\kappa\)-sequences, it follows from
this that \(P(\gamma)\subseteq M\).

We can now find larger independent families: since \(P(\gamma)\subseteq M\),
\((2^\gamma)^M\geq 2^\gamma\), and in particular,  \(M\) satisfies that there is
a \(\kappa\)-independent family of subsets of \(\gamma\) of cardinality
\((2^\gamma)^V\). 

Assume finally that \(\delta \leq 2^\gamma\) is a cardinal and \(\kappa\) is
\(\delta\)-strongly compact. Then let \(\Gamma\in M\) be a
\(\kappa\)-independent family of subsets of \(\gamma\) in \(M\) with cardinality
\(\delta\). As in the previous paragraph, any \(\kappa\)-complete ultrafilter
over the algebra generated by \(\Gamma\) belongs to \(M\), so \(M\) has the
\(\kappa\)-approximation property at \(\delta\) by \cref{ApproxPrp}. Since \(M\)
is closed under \({<}\kappa\)-sequences, it follows from this that
\(P(\delta)\subseteq M\).
\end{proof}

\begin{proof}[Proof of \cref{KetonenStrong}] By \cref{UFAmenableChar}, every
countably complete \(M\)-ultrafilter \(U\) on \(\gamma <\lambda\) belongs to
\(M\). Therefore if \(\gamma < \lambda\), our strong compactness assumption on
\(\kappa_\lambda\) implies the hypotheses of \cref{GeneralStrong} hold at
\(\gamma\), and so \(P(\gamma)\subseteq M\).
\end{proof}

Having proved that \(\mathscr K_\lambda\) has some strength, let us now turn to
the supercompactness properties of \(\mathscr K_\lambda\). 

\begin{thm}\label{KetonenTight}\index{Ketonen ultrafilter!\(\mathscr K_\lambda\)!tightness}
Suppose \(\lambda\) is a regular Fr\'echet cardinal and \(\kappa_\lambda\) is
\(\lambda\)-strongly compact. Let \(j : V\to M\) be the ultrapower of the
universe by \(\mathscr K_\lambda\). Then 
\begin{itemize}
	\item \(j\) is \(\lambda\)-tight.
	\item \(j\) is \(\gamma\)-supercompact for all \(\gamma < \lambda\).
\end{itemize}
In other words, \(M^\gamma\subseteq M\) for all \(\gamma < \lambda\) and \(M\)
has the \({\leq}\lambda\)-covering property.
\begin{proof}
		Suppose towards a contradiction that \(j\) is not \(\lambda\)-tight. By
		\cref{KetonenCov}, it follows that \(\delta = \textnormal{cf}^M(\sup
		j[\lambda])>\lambda\). By \cref{UFAmenableChar}, any countably complete
		\(M\)-ultrafilter \(U\) on \(\lambda\) belongs to \(M\). But then by
		\cref{GeneralStrong}, \(P(\lambda)\subseteq M\). But then \(\mathscr
		K_\lambda\) itself is a countably complete \(M\)-ultrafilter on
		\(\lambda\), so \(\mathscr K_\lambda\in M\). This contradicts the
		irreflexivity of the Mitchell order (\cref{MOStrict}).
		
		Now that we know \(j\) is \(\lambda\)-tight, let us show that \(j\) is
		\(\gamma\)-supercompact for all \(\gamma < \lambda\). We may assume by
		induction that \(j\) is \({<}\gamma\)-supercompact. Then if \(\gamma\)
		is singular, it is easy to see that \(j\) is \(\gamma\)-supercompact.
		Therefore assume \(\gamma\) is regular. Let \(\gamma' = \text{cf}^M(\sup
		j[\gamma])\). Then \(\gamma' \leq \lambda\) since \(j\) is
		\(\lambda\)-tight and hence \(j\) is \((\gamma,\lambda)\)-tight. Since
		\(\gamma <\lambda\), in fact \(\gamma' < \lambda\). Thus
		\(P(\gamma')\subseteq M\)  by \cref{KetonenStrong}. By
		\cref{KetonenCov}, \(j\) is \((\gamma,\gamma')\)-tight, so fix \(A\in
		M\) with \(|A|^M = \gamma'\) and \(j[\gamma]\subseteq A\). Note that
		since \(|A|^M = \gamma'\), \(P(A)\subseteq M\). Therefore
		\(j[\gamma]\subseteq M\). Therefore \(j\) is \(\gamma\)-supercompact, as
		desired.
		
		That \(M^\gamma\subseteq M\) for all \(\gamma < \lambda\) is an
		immediate consequence of  \cref{UltrapowerSC}. That \(M\) has the
		\({\leq}\lambda\)-covering property is an immediate consequence of
		\cref{UltrapowerStrongC}.
\end{proof}
\end{thm}

Finally, if \(\lambda\) is not a strongly inaccessible cardinal, we can show
that \(j_{\mathscr K_\lambda}\) is precisely as supercompact as it should be:
\begin{thm}[UA]\label{AccessibleSupercompact}
Suppose \(\lambda\) is a regular Fr\'echet cardinal and \(\kappa_\lambda\) is
\(\lambda\)-strongly compact. Let \(j : V\to M\) be the ultrapower of the
universe by \(\mathscr K_\lambda\). If \(\lambda\) is not strongly inaccessible
then \(j\) is \(\lambda\)-supercompact.
\begin{proof}
Let \(\kappa = \kappa_\lambda\) for ease of notation. We split into two cases:
\begin{case}\label{HighCofCase}
	For some \(\gamma < \lambda\) with \(\text{cf}(\gamma)\geq\kappa\),
	\(2^\gamma\geq\lambda\).
\end{case}
\begin{proof}[Proof in \cref{HighCofCase}] Since \(\gamma < \lambda\), by
	\cref{UFAmenableChar} every countably complete \(M\)-ultrafilter on
	\(\gamma\) belongs to \(M\). Since \(\text{cf}(\gamma)\geq\kappa\),
	\(\lambda\leq 2^\gamma\), and \(\kappa\) is \(\lambda\)-strongly compact, we
	can therefore apply the second part of \cref{GeneralStrong} to conclude that
	\(P(\lambda)\subseteq M\). 
	
	Given that \(j\) is \(\lambda\)-tight by \cref{KetonenTight}, it now follows
	easily that \(j\) is \(\lambda\)-supercompact: fix \(A\in M\) with \(|A|^M =
	\lambda\) and \(j[\lambda]\subseteq A\); then \(P(A)\subseteq M\) so
	\(j[\lambda]\in M\), as desired.
\end{proof}

\begin{case}\label{LowCofCase}
	For all \(\gamma < \lambda\) with \(\text{cf}(\gamma)\geq\kappa\),
	\(2^\gamma< \lambda\).
\end{case}
\begin{proof}[Proof in \cref{LowCofCase}] Since \(\lambda\) is not inaccessible,
	there is some \(\eta <\lambda\) such that \(2^\eta\geq\lambda\). Let
	\(\gamma = \eta^{<\kappa}\). Then \(\text{cf}(\gamma) \geq \kappa\) and
	\(2^{\gamma}\geq 2^\eta\geq \lambda\). Therefore by our case hypothesis,
	\(\lambda\leq\gamma\). By \cref{KetonenTight}, \(j\) is
	\(\eta\)-supercompact. By \cref{SmallCfCompact}, \(j\) is
	\(\eta^{<\kappa}\)-supercompact. Therefore \(j\) is \(\lambda\)-supercompact
	as desired.
	\end{proof}
Thus in either case \(j\) is \(\lambda\)-supercompact, which completes the
proof.
\end{proof}
\end{thm}
\section{Fr\'echet cardinals}\label{FrechetSection}
\subsection{The Fr\'echet successor}
Given the results of \cref{KSuperSection}, to analyze \(\mathscr K_\lambda\)
when \(\lambda\) is a regular Fr\'echet cardinal, it would be enough to show
that its completeness \(\kappa_\lambda\) is \(\lambda\)-strongly compact.
\begin{figure}
	\center
	\includegraphics[scale=.6]{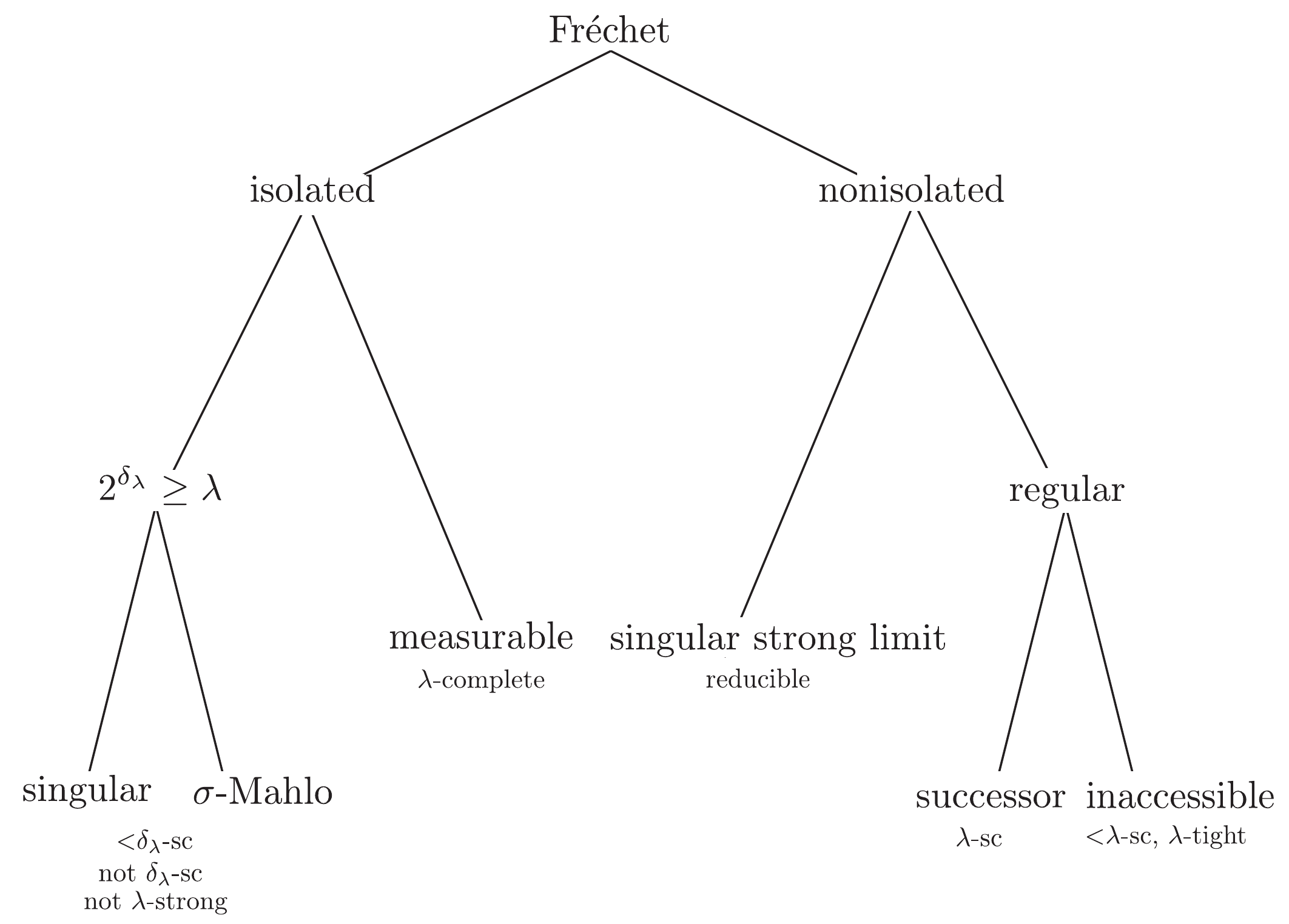}
	\caption{Types of Fr\'echet cardinals.}
\end{figure}
The following easy generalization of Ketonen's Theorem (\cref{KetonenThm})
reduces this to the analysis of Fr\'echet cardinals in the interval
\([\kappa_\lambda,\lambda]\):
\begin{prp}\label{KetonenFrechet}
	Suppose \(\lambda\) is a regular Fr\'echet cardinal. Suppose \(j : V\to M\)
	is the ultrapower of the universe by a Ketonen ultrafilter \(U\) on
	\(\lambda\).  Suppose \(\kappa \leq\lambda\) is a cardinal and every regular
	cardinal in the interval \([\kappa,\lambda]\) is Fr\'echet. Then \(j\) is
	\((\lambda, \delta)\)-tight for some \(\delta < j(\kappa)\). In particular,
	if \(\kappa = \textsc{crt}(j)\) then \(\kappa\) is \(\lambda\)-strongly
	compact.
	\begin{proof}
		Since \(U\) is Ketonen, the \(M\)-cardinal \(\delta = \text{cf}^M(\sup
		j[\lambda])\) is not Fr\'echet in \(M\). Therefore by elementarity
		\(\delta\notin j([\kappa,\lambda])\). Since \(\delta < j(\lambda)\), we
		must have \(\delta < j(\kappa)\). \cref{KetonenCov} implies that \(j\)
		is \((\lambda,\delta)\)-tight, proving the proposition.
	\end{proof}
\end{prp}

Suppose \(\lambda\) is a regular Fr\'echet cardinal. To obtain that every
regular cardinal in the interval \([\kappa_\lambda,\lambda)\) is Fr\'echet, it
actually suffices to show that every {\it successor} cardinal in the interval
\((\kappa_\lambda,\lambda]\) is Fr\'echet. (See \cref{SuccessorSuffices}.) Our
approach to this problem is as follows. Fix an ordinal \(\gamma\in
[\kappa_\lambda,\lambda)\). We consider the {\it Fr\'echet successor of
\(\gamma\):}
\begin{defn}\index{Fr\'echet cardinal!The Fr\'echet successor of \(\gamma\) (\(\gamma^\sigma\))}\index{\(\gamma^\sigma\)}\index{\(\gamma^\sigma\)|seealso{Fr\'echet cardinal}}
Suppose \(\gamma\) is an ordinal. Then the {\it Fr\'echet successor of
\(\gamma\)}, denoted \(\gamma^\sigma\), is the least Fr\'echet cardinal strictly
greater than \(\gamma\).
\end{defn}
We will attempt to use the fact that \(\gamma\) lies in the interval
\([\kappa_\lambda,\lambda)\) to show that \(\gamma^\sigma = \gamma^+\). Since
\(\gamma^\sigma\) is Fr\'echet by definition, this would show \(\gamma^+\) is
Fr\'echet. In this way, we we would establish that every successor cardinal in
the interval \((\kappa_\lambda,\lambda]\) is Fr\'echet, as desired.

The following classic result of Prikry \cite{Prikry} shows in particular that
there is nontrivial structure to the Fr\'echet cardinals even if we do not
assume UA:
\begin{thm}[Prikry]\label{PrikryThm}\index{Indecomposable ultrafilter!Prikry's Theorem}\index{Prikry's Theorem}
	Suppose \(\lambda\) is a cardinal and \(U\) is a \(\lambda^+\)-decomposable
	ultrafilter. Then \(U\) is \(\textnormal{cf}(\lambda)\)-decomposable.\qed
\end{thm}

A key part of our analysis of Fr\'echet cardinals is the following
generalization of \cref{PrikryThm}:
\begin{prp}\label{SuccessorPrp}
Suppose \(\eta\) is a cardinal such that \(\eta^+\) is Fr\'echet. Either
\(\eta\) is Fr\'echet or \(\eta\) is a singular cardinal and all sufficiently
large regular cardinals below \(\eta\) are Fr\'echet.
\begin{proof}
Suppose \(\gamma^\sigma = \eta^+\). We will show that either \(\eta\) is
Fr\'echet or \(\eta\) is a limit of Fr\'echet cardinals. Fix a countably
complete uniform ultrafilter \(U\) on \(\eta^+\), and let \(j : V\to M\) be the
ultrapower of the universe by \(U\). Let \[U_* = \{A\in j(P(\eta^+)) :
j^{-1}[A]\in U\}\] Thus \(U_*\) is an \(M\)-ultrafilter. Note that
\(\lambda_{U_*} < j(\eta^+)\) since for example \(\sup j[\eta^+]\in U_*\). Thus
\(\lambda_{U_*} \leq j(\eta)\). 

The proof now splits into two cases:
\begin{case}\label{UniformCase}
	\(\lambda_{U_*}\geq \sup j[\eta]\).
\begin{proof}[Proof in \cref{UniformCase}] Let \(\lambda = \lambda_{U_*}\). Then
		\(\sup j[\eta]\leq \lambda \leq j(\eta)\). Let \(W_*\) be an
		\(M\)-ultrafilter on \(j(\eta)\) that concentrates on \(\lambda\) and is
		isomorphic to \(U_*\). In other words, there is a set \(X\in U_*\) and a
		bijection \(f : \lambda\to X\) with \(f\in M\) such that \(W_* =
		\{f^{-1}[A] : A\in U_*\}\). All we need about \(W_*\) is that
		\(\lambda_{W_*} = \lambda\geq \sup j[\eta]\). Let \[W = j^{-1}[W_*]\]
		Then \(W\) is a countably complete ultrafilter on \(\eta\). 
		
		We claim that \(W\) is uniform. Suppose \(A\in W\). Then \(j(A) \in
		W_*\) so \(|j(A)|^M = \lambda\). In particular, since \(\lambda\geq \sup
		j[\eta]\), for any cardinal \(\kappa < \eta\), \(|j(A)|^M > j(\kappa)\),
		and therefore \(|A| > \kappa\). It follows that \(|A| \geq \eta\). Thus
		\(W\) is uniform.
\end{proof}
\end{case}

\begin{case}\label{NonuniformCase}
	\(\lambda_{U_*} < \sup j[\eta]\).
	\begin{proof}[Proof in \cref{NonuniformCase}] Fix \(\kappa < \eta\) and
	\(B\in U_*\) such that letting \(\delta = |B|^M\), we have \(\delta <
	j(\kappa)\). Let \(A = j^{-1}[B]\). Then \(A\in U\) so \(|A| = \eta^+\)
	since \(U\) is a uniform ultrafilter on \(\eta^+\). Since \(j[A]\subseteq
	B\), it follows that \(j\) is \((\eta^+,\delta)\)-tight. 
	
	We claim that \(j\) is discontinuous at every regular cardinal \(\iota\) in
	the interval \([\kappa,\eta^+]\). To see this, note that \(j(\iota) >
	\delta\) is a regular cardinal of \(M\). On the other hand, \(j[\iota]\) is
	contained in a set \(C\in M\) such that \(|C|^M \leq\delta\) since \(j\) is
	\((\iota,\delta)\)-tight. Therefore \(C\) is not cofinal in \(j(\iota)\),
	and hence neither is \(j[\iota]\). It follows that \(j\) is discontinuous at
	\(\iota\).
	
	Since \(j\) is discontinuous at every regular cardinal in the interval
	\([\kappa,\eta^+]\), which contains \(\eta\), it follows that either
	\(\eta\) is a regular Fr\'echet cardinal or \(\eta\) is a singular cardinal
	and all sufficiently large regular cardinals below \(\eta\) are Fr\'echet.
	\end{proof}
\end{case}
Thus in either case, the conclusion of the proposition holds.
\end{proof}
\end{prp}

An interesting feature of \cref{SuccessorPrp} is that it {\it does not} seem to
show that every \(\eta^+\)-decomposable ultrafilter \(U\) is either
\(\eta\)-decomposable or \(\iota\)-decomposable for all sufficiently large
\(\iota < \eta\). Instead the proof shows that this is true of \(U^2\). (Under
UA, we can in fact prove that every \(\eta^+\)-decomposable countably complete
ultrafilter \(U\) is either \(\eta\)-decomposable or \(\iota\)-decomposable for
all sufficiently large \(\iota < \eta\).)

\cref{SuccessorPrp} has two important consequences. The first is our claim above
that one need only show that all successor cardinals in
\([\kappa_\lambda,\lambda]\) are Fr\'echet to conclude that all regular
cardinals in \([\kappa_\lambda,\lambda]\) are. (This is really just a
consequence of \cref{PrikryThm}.)

\begin{cor}\label{SuccessorSuffices}
Suppose \(\kappa \leq \lambda\) are cardinals and every successor cardinal in
the interval \((\kappa,\lambda]\) is Fr\'echet. Then every regular cardinal in
the interval \([\kappa,\lambda)\) is Fr\'echet.
\begin{proof}
	Suppose \(\iota\) is a regular cardinal in the interval
	\([\kappa,\lambda)\). Then \(\iota^+\in (\kappa,\lambda]\), so \(\iota^+\)
	is a Fr\'echet cardinal. Therefore \(\iota\) is a Fr\'echet cardinal by
	\cref{SuccessorPrp}.
\end{proof}
\end{cor}

The consequence of \cref{SuccessorPrp} that is ultimately most important here is
a constraint on the Fr\'echet successor operation:
\begin{cor}\label{SuccessorDichotomy}
Suppose \(\gamma\) is an ordinal and \(\gamma^\sigma\) is a successor cardinal.
Then \(\gamma^\sigma = \gamma^+\).
\begin{proof}
	Suppose towards a contradiction that \(\gamma^\sigma = \eta^+\) for some
	cardinal \(\eta > \gamma\). Since \(\eta^+\) is Fr\'echet, by
	\cref{SuccessorPrp}, \(\eta\) is either Fr\'echet or a limit of Fr\'echet
	cardinals. Either way, there is a Fr\'echet cardinal in the interval
	\((\gamma,\eta^+)\). But the definition of \(\gamma^\sigma\) implies that
	there are no Fr\'echet cardinals in \((\gamma,\gamma^\sigma)\). This is a
	contradiction.
	
	Thus \(\gamma^\sigma = \eta^+\) for some cardinal \(\eta \leq \gamma\). In
	other words, \(\gamma^\sigma = \gamma^+\).
\end{proof}
\end{cor}

The problematic cases in the analysis of the Fr\'echet successor function
therefore occur when \(\gamma^\sigma\) is a limit cardinal:

\begin{defn}\index{Isolated cardinal}
A cardinal \(\lambda\) is {\it isolated} if the following hold:
\begin{itemize}
\item \(\lambda\) is Fr\'echet.
\item \(\lambda\) is a limit cardinal.
\item \(\lambda\) is not a limit of Fr\'echet cardinals.	
\end{itemize}	
\end{defn}
By \cref{SuccessorPrp}, \(\lambda\) is isolated if and only if \(\lambda =
\gamma^\sigma\) for some ordinal \(\gamma\) such that \(\gamma^+ < \lambda\).
Our analysis of Fr\'echet cardinals would be essentially complete if we could
prove the following conjecture:

\begin{conj}[UA]\label{IsolationConj}
A cardinal \(\lambda\) is isolated if and only if \(\lambda\) is a measurable
cardinal, \(\lambda\) is not a limit of measurable cardinals, and no cardinal
\(\kappa < \lambda\) is \(\lambda\)-supercompact.
\end{conj}

\cref{IsolatedStrongLimit} below shows that \cref{IsolationConj} is a
consequence of UA + GCH, so to some extent this problem is solved in the most
important case. But assuming UA alone, we do not know how to rule out, for
example, the existence of singular isolated cardinals. Enacting an analysis of
isolated cardinals under UA that is as complete as possible allows us to prove
our main results without cardinal arithmetic assumptions.

\subsection{The strong compactness of \(\kappa_\lambda\)}\label{kappaStrongCSection}
In this section we will prove the following theorem:
\begin{thm}[UA]\label{NonisolatedCompact}
Suppose \(\lambda\) is a nonisolated regular Fr\'echet cardinal. Then
\(\kappa_\lambda\) is \(\lambda\)-strongly compact.
\end{thm}
This yields the following corollary, which gives a complete analysis of
Fr\'echet successor cardinals:
\begin{cor}[UA]\label{SuccessorThm}\index{Ketonen ultrafilter!\(\mathscr K_\lambda\)!supercompactness}
Suppose \(\lambda\) is a Fr\'echet successor cardinal. Then \(\kappa_\lambda\)
is \(\lambda\)-supercompact and in fact the ultrapower embedding associated to
\(\mathscr K_\lambda\) is \(\lambda\)-supercompact.
\begin{proof}
This is an immediate consequence of \cref{NonisolatedCompact} and
\cref{AccessibleSupercompact}.
\end{proof}
\end{cor}
In general, we only obtain
\begin{prp}[UA]\label{GeneralThm}
Suppose \(\lambda\) is a nonisolated regular Fr\'echet cardinal. Then
\(\kappa_{\lambda}\) is \({<}\lambda\)-supercompact and \(\lambda\)-strongly
compact. In fact, the ultrapower embedding associated to \(\mathscr K_\lambda\)
is \({<}\lambda\)-supercompact and \(\lambda\)-tight.
\begin{proof}
	This is an immediate consequence of \cref{NonisolatedCompact} and
	\cref{KetonenTight}.
\end{proof}
\end{prp}

As we have sketched above, the proof of \cref{NonisolatedCompact} will follow
from an analysis of Fr\'echet cardinals in the interval
\([\kappa_\lambda,\lambda]\):

\begin{lma}\label{ImplicationLma}
	Suppose \(\kappa \leq \lambda\) are cardinals and there are no isolated
	cardinals in the interval \((\kappa,\lambda]\). Suppose that for all
	\(\gamma\in [\kappa,\lambda)\), there is a Fr\'echet cardinal in the
	interval \((\gamma,\lambda]\). Then every regular cardinal in the interval
	\([\kappa,\lambda)\) is Fr\'echet.
	\begin{proof}
			Since \(\lambda\) is Fr\'echet, we need only show that every regular
			cardinal in the interval \([\kappa,\lambda)\) is Fr\'echet. By
			\cref{SuccessorSuffices}, for this it is enough to show that every
			successor cardinal in the interval \((\kappa,\lambda]\) is
			Fr\'echet. In other words, it suffices to show that for any ordinal
			\(\gamma\in [\kappa,\lambda)\), \(\gamma^+\) is Fr\'echet. Therefore
			fix \(\gamma\in [\kappa,\lambda)\). By assumption, \(\gamma^\sigma
			\in (\gamma,\lambda]\), so in particular \(\gamma^\sigma\) is not
			isolated. Therefore \(\gamma^\sigma\) is not a limit cardinal. It
			follows that \(\gamma^\sigma\) is a successor cardinal, so by
			\cref{SuccessorPrp}, \(\gamma^\sigma = \gamma^+\), as desired.
	\end{proof}
\end{lma}

Our goal now it to prove the following lemma:
\begin{lma}[UA]\label{NonisolatedInterval}
	Suppose \(\lambda\) is a Fr\'echet cardinal that is either regular or
	isolated. Then there are no isolated cardinals in the interval
	\([\kappa_\lambda,\lambda)\). 
\end{lma}

Given this, we could complete the proof of \cref{NonisolatedCompact} as follows:
\begin{proof}[Proof of \cref{NonisolatedCompact} assuming
	\cref{NonisolatedInterval}] By \cref{NonisolatedInterval}, there are no
	isolated cardinals in the interval \([\kappa_\lambda,\lambda)\). Since
	\(\lambda\) is not isolated, there are no isolated cardinals in the interval
	\([\kappa_\lambda,\lambda]\). Therefore applying \cref{ImplicationLma},
	every regular cardinal in the interval \([\kappa_\lambda,\lambda]\) is
	Fr\'echet. By \cref{KetonenFrechet}, it follows that \(\kappa_\lambda\) is
	\(\lambda\)-strongly compact.
\end{proof}

We now proceed to the proof of \cref{NonisolatedInterval}. We will first need to
improve our understanding of isolated cardinals. The first step is to provide
some criteria that guarantee a cardinal's nonisolation:

\begin{lma}\label{Nonisolation1}
	Suppose \(\eta\) is a limit cardinal. Suppose \(U\) is a countably complete
	uniform ultrafilter on \(\eta\). Suppose \(W\) is a countably complete
	ultrafilter such that \(j_W\) is discontinuous at \(\eta\) and \(U\I W\).
	Then \(\eta\) is a limit of Fr\'echet cardinals.
	\begin{proof}
		Let \(i : V\to N\) be the ultrapower of the universe by \(W\). Let \[U_*
		= s_W(U) = \{B\in i(P(\lambda)) : i^{-1}[B]\in U\}\] By \cref{PushUlt},
		\(U_*\in N\).
		
		\begin{case}\label{U*Case}
			\(\lambda_{U_*} \geq\sup i[\eta]\)
			\begin{proof}[Proof in \cref{U*Case}] Working in \(N\),
				\(\lambda_{U_*}\) is a Fr\'echet cardinal \(\lambda\) with
				\(\sup i [\eta]\leq \lambda < i(\eta)\). It follows that for any
				\(\kappa < \eta\), \(N\) satisfies that there is a Fr\'echet
				cardinal strictly between \(i(\kappa)\) and \(i(\eta)\), and so
				by elementarity there is a Fr\'echet cardinal strictly between
				\(\kappa\) and \(\eta\). It follows that \(\eta\) is a limit of
				Fr\'echet cardinals.
			\end{proof}
		\end{case}
	
		\begin{case}\label{TightCase}
			\(\lambda_{U_*} < \sup i[\eta]\)
		\begin{proof}[Proof in \cref{TightCase}] Fix \(\kappa < \eta\) and
			\(B\in U_*\) such that letting \(\delta = |B|^N\), \(\delta <
			i(\kappa)\). Let \(A = i^{-1}[B]\). Then \(A\in U\), so \(|A| =
			\eta\) by the uniformity of \(U\). Since \(|A| = \eta\) and
			\(i[A]\subseteq B\), \(i\) is \((\eta,\delta)\)-tight by
			\cref{KetonenCov}. It follows that \(i\) is discontinuous at every
			regular cardinal in the interval \([\kappa,\eta]\). (See the proof
			of \cref{SuccessorPrp}.) In particular, \(\eta\) is a limit of
			Fr\'echet cardinals.
		\end{proof}
		\end{case}
	In either case, \(\eta\) is a limit of Fr\'echet cardinals, as desired.
	\end{proof}
\end{lma}

The second nonisolation lemma brings in a bit more of the theory of the internal
relation:
\begin{lma}[UA]\label{Nonisolation2}
		Suppose \(\eta\) is a Fr\'echet limit cardinal. Suppose there is a
		countably complete ultrafilter \(W\) such that \(\mathscr K_\eta\I W\)
		but \(W\not \I \mathscr K_\eta\). Then \(\eta\) is a limit of Fr\'echet
		cardinals.
\end{lma}
\begin{proof}
	By \cref{Nonisolation1}, if \(j_W\) is discontinuous at \(\eta\), then
	\(\eta\) is a limit of Fr\'echet cardinals. Therefore assume without loss of
	generality that \(j_W\) is continuous at \(\eta\).
	
	By the basic theory of the internal relation (\cref{IChar}), since
	\(\mathscr K_\eta\I W\), the translation \(\tr W {\mathscr K_\eta}\) is
	equal to the pushforward \(s_W(\mathscr K_\eta)\). 
	
	Since \(W\not \I \mathscr K_\eta\), the theory of the internal relation
	(\cref{IChar}) implies that in \(M_W\), \(\tr W {\mathscr K_\eta}\sE
	j_W(\mathscr K_\eta)\). Since \(M_W\) satisfies that \(j_W(\mathscr
	K_\eta)\) is the \(\sE\)-least uniform ultrafilter on \(j_W(\eta)\), it
	follows that \[\lambda_{\tr W {\mathscr K_\eta}} < j_W(\eta)\] But \(\tr W
	{\mathscr K_\eta} = s_W(\mathscr K_\eta)\) and \( j_W(\eta) = \sup
	j_W[\eta]\) by our assumption that \(j_W\) is continuous at \(\eta\). Thus
	\[\lambda_{s_ W (\mathscr K_\eta)} < \sup j_W[\eta]\]
	
	Fix \(\kappa < \eta\) and \(B\in s_W(\mathscr K_\eta)\) such that \(\delta =
	|B|^{M_W} < j_W(\kappa)\). Let \(A = j_W^{-1}[B]\). Then \(A\in \mathscr
	K_\eta\), so \(|A| = \eta\). Moreover \(j_W[A]\subseteq B\in M_W\), so
	\(j_W\) is \((\eta,\delta)\)-tight. In particular, \(j_W\) is discontinuous
	at every regular cardinal in the interval \([\kappa,\eta]\). (See the proof
	of \cref{SuccessorPrp}.) Therefore \(\eta\) is a limit of Fr\'echet
	cardinals.
\end{proof}

Finally, we need a version of \cref{EmbeddingChar} that applies at singular
cardinals.

We use a lemma that follows immediately from the ultrafilter sum construction:

\begin{lma}
	Suppose \(U\) is a countably complete ultrafilter on a cardinal \(\lambda\)
	and \(U'\) is a countably complete \(M_U\)-ultrafilter with \(\lambda_{U'}
	\leq j_U(\lambda)\). Then there is a countably complete ultrafilter \(W\) on
	\(\lambda\) such that \(j_W = j_{U'}^{M_U}\circ j_U\).\qed
\end{lma}

\begin{prp}[UA]\label{IsolatedInternal}
	Suppose \(\lambda\) is an isolated cardinal. Then \(\mathscr K_\lambda\) is
	\(\lambda\)-internal.
	\begin{proof}
		Suppose \(D\) is a countably complete ultrafilter on a cardinal \(\gamma
		< \lambda\). We will show \(D\I \mathscr K_\lambda\). Since \(\lambda\)
		is isolated, by increasing \(\gamma\), we may assume \(\lambda =
		\gamma^\sigma\).
		
		Assume towards a contradiction that in \(M_D\), \[\tr D {\mathscr
		K_\lambda} \sE j_D(\mathscr K_\lambda)\] Then \(\lambda_{\tr D {\mathscr
		K_\lambda}} < j_D(\lambda)\), and so since \(\lambda_{\tr D {\mathscr
		K_\lambda}}\) is a Fr\'echet cardinal of \(M_D\), \(\lambda_{\tr D
		{\mathscr K_\lambda}} \leq j_D(\gamma)\). Therefore, there is an
		ultrafilter \(W\) on \(\gamma\) such that \[j_W = j_{\tr D {\mathscr
		K_\lambda}}^{M_D}\circ j_D = j_{\tr {\mathscr K_\lambda} D}^{M_{\mathscr
		K_\lambda}}\circ j_{\mathscr K_\lambda}\] It follows from the basic
		theory of the Rudin-Keisler order (\cref{RKChar}) that \(\mathscr
		K_\lambda\RK W\), which contradicts that \(\lambda_{\mathscr K_\lambda}
		= \lambda > \gamma\geq\lambda_W\).
		
		Thus our assumption was false, and in fact, \(j_D(\mathscr K_\lambda)\E
		\tr D {\mathscr K_\lambda}\) in \(M_D\). By the theory of the internal
		relation (\cref{IChar}), this implies that \(D\I \mathscr K_\lambda\).
	\end{proof}
\end{prp}

In \cref{IsolatedUFSection}, we prove a much stronger version of this theorem
that constitutes a complete generalization of \cref{UniversalProperty} to
isolated cardinals.

\begin{lma}[UA]\label{Nonoverlapping}
		Suppose \(\eta < \lambda\) is are Fr\'echet cardinals that are regular
		or isolated. Then either \(\eta < \kappa_\lambda\) or \(\mathscr
		K_\lambda\not\I \mathscr K_\eta\).
		\begin{proof}
			By \cref{EmbeddingChar} or \cref{IsolatedInternal}, \(\mathscr
			K_\eta\) and \(\mathscr K_\lambda\) are \(\lambda\)-internal. 
			
			Assume \(\mathscr K_\lambda \I \mathscr K_\eta\). Note that we also
			have \(\mathscr K_\eta \I \mathscr K_\lambda\) since \(\mathscr
			K_\lambda\) is \(\lambda\)-uniform. By \cref{IUniformPrp}, \(\eta <
			\kappa_\lambda\).
		\end{proof}
\end{lma}

We can finally prove \cref{NonisolatedInterval}.

\begin{proof}[Proof of \cref{NonisolatedInterval}] Suppose towards a
	contradiction that \(\eta\in [\kappa_\lambda,\lambda)\) is isolated. Then by
	\cref{Nonoverlapping}, \(\mathscr K_\lambda\not\I \mathscr K_\eta\).
	Therefore by \cref{Nonisolation2}, \(\eta\) is a limit of Fr\'echet
	cardinals, contrary to the assumption that \(\eta\) is isolated.
\end{proof}

Since we will use it repeatedly, it is worth noting that \(\kappa_\lambda\) can
be characterized in terms of isolated cardinals:
\begin{lma}[UA]\label{KappaChar}
	Suppose \(\lambda\) is a nonisolated regular Fr\'echet cardinal. Then
	\(\kappa_\lambda\) is the supremum of the isolated cardinals less than
	\(\lambda\).
	\begin{proof}
		Let \(\kappa\) be the supremum of the isolated cardinals less than
		\(\lambda\). By \cref{NonisolatedInterval}, there are no isolated
		cardinals in the interval \([\kappa_\lambda,\lambda)\), so \(\kappa \leq
		\kappa_\lambda\). 
		
		Since there are no isolated cardinals in the interval
		\((\kappa,\lambda]\), \cref{ImplicationLma} implies that every regular
		cardinal in the interval \([\kappa,\lambda]\) is Fr\'echet. By
		\cref{KetonenFrechet}, it follows that \(\kappa_\lambda\leq \kappa\).
		Thus \(\kappa_\lambda = \kappa\), as desired.
	\end{proof}
\end{lma}

\subsection{The first supercompact cardinal}\label{FirstSuperSection}
In this subsection, we show how the theory of the internal relation can be used
to characterize the least supercompact cardinal (and its local instantiations).
\begin{thm}[UA]\label{omega1SC}
	Suppose \(\lambda\) is a successor cardinal and \(\kappa\) is the least
	\((\omega_1,\lambda)\)-strongly compact cardinal. Then \(\kappa\) is
	\(\lambda\)-supercompact. In fact, \(\kappa = \kappa_\lambda\).
	\begin{proof}
		Since \(\kappa\) is \((\omega_1,\lambda)\)-strongly compact, every
		regular cardinal in the interval \([\kappa,\lambda]\) is Fr\'echet. By
		\cref{KetonenFrechet}, \(\kappa_\lambda \leq \kappa\). By
		\cref{SuccessorThm}, \(\kappa_\lambda\) is \(\lambda\)-supercompact. In
		particular, \(\kappa_\lambda\) is \((\omega_1,\lambda)\)-strongly
		compact. Therefore \(\kappa \leq \kappa_\lambda\), and hence \(\kappa =
		\kappa_\lambda\). Thus \(\kappa\) is \(\lambda\)-supercompact, as
		desired.
	\end{proof}
\end{thm}
\begin{cor}[UA] Suppose \(\lambda\) is a successor cardinal and \(\kappa\) is
		the least \(\lambda\)-strongly compact cardinal. Then \(\kappa\) is
		\(\lambda\)-supercompact. In fact, \(\kappa = \kappa_\lambda\).\qed
\end{cor}

\begin{cor}[UA] The least \((\omega_1,\textnormal{Ord})\)-strongly compact
	cardinal \(\kappa\) is supercompact.
	\begin{proof}
		No cardinal \(\delta < \kappa\) is \((\omega_1,\kappa)\)-strongly
		compact. In particular, for any successor cardinal \(\lambda >\kappa\),
		\(\kappa\) is the least \((\omega_1,\lambda)\)-strongly compact
		cardinal. Therefore \(\kappa\) is \(\lambda\)-supercompact by
		\cref{omega1SC}.
	\end{proof}
\end{cor}

\begin{thm}[UA]\label{StrongSuper}\index{Strongly compact cardinal!equivalence with supercompactness}\index{Identity crisis}
	The least strongly compact cardinal is supercompact.\qed
\end{thm}

\cref{ArbHigh} identifies the following ordinals as key thresholds in the
structure theory of countably complete ultrafilters:
\begin{defn}
	The {\it ultrapower threshold} is the least cardinal \(\kappa\) such that
	for all \(\alpha\), there is an ultrapower embedding \(j : V\to M\) such
	that \(j(\kappa) > \alpha\).\index{Ultrapower threshold}
	
	Suppose \(\gamma\) is an ordinal. The {\it \(\gamma\)-threshold} is the
	least ordinal \(\kappa \leq \gamma\) such that for all \(\alpha < \gamma\)
	is an ultrapower embedding \(j : V\to M\) such that \(j(\kappa) > \alpha\). 
\end{defn}

The ultrapower threshold cannot be proved to exist without large cardinal
assumptions, but for any ordinal \(\gamma\), the \(\gamma\)-threshold exists and
is less than or equal to \(\gamma\). 

\begin{lma}\label{IdemThreshold}
	Suppose \(\kappa\) is a cardinal. If \(\kappa\) is the \(\gamma\)-threshold
	for some ordinal \(\gamma\) then \(\kappa\) is the \(\kappa\)-threshold.
	\begin{proof}
		We may assume without loss of generality that \(\kappa < \gamma\). Let
		\(\nu\leq \kappa\) be the \(\kappa\)-threshold.
		
		 We claim that for any \(\alpha < \gamma\), there is an ultrapower
		 embedding \(h : V\to N\) such that \(h(\nu) > \alpha\).  Fix \(\alpha <
		 \gamma\). Let \(j :V\to M\) be such that \(j(\kappa) > \alpha\). In
		 \(M\), \(j(\nu)\) is the \(j(\kappa)\)-threshold, so since \(\alpha <
		 j(\kappa)\), there is an internal ultrapower embedding \(i : M\to N\)
		 such that \(i(j(\nu)) > \alpha\). Let \(h = i\circ j\). Then \(h : V\to
		 N\) is an ultrapower embedding such that \(h(\nu) > \alpha\), as
		 desired.
		
		By the minimality of the \(\gamma\)-threshold, \(\kappa \leq \nu\).
		Hence \(\kappa = \nu\) as desired.
	\end{proof}
\end{lma}

\begin{thm}[UA]\label{ThreshThm}
	Suppose \(\lambda\) is a strong limit cardinal and \(\kappa < \lambda\) is
	the \(\lambda\)-threshold. Then \(\kappa\) is \(\gamma\)-supercompact for
	all \(\gamma < \lambda\).
\end{thm}

The proof uses the following lemma, an often-useful approximation to
\cref{IsolationConj}:
\begin{lma}[UA]\label{DoubleSigma}
	Suppose \(\lambda_0 \) is an isolated cardinal and \(\lambda_1 =
	(\lambda_0)^\sigma\). Then \(\lambda_1\) is measurable.
	\begin{proof}
		Note that \(\kappa_{\lambda_1} >\lambda_0\): otherwise \(\lambda_0 \in
		[\kappa_{\lambda_1},\lambda_1)\) contrary to the fact that there are no
		isolated cardinals in the interval \([\kappa_{\lambda_1},\lambda_1)\) by
		\cref{NonisolatedInterval}. Since \(\kappa_{\lambda_1}\) is measurable,
		\(\kappa_{\lambda_1}\) is Fr\'echet. Hence \(\lambda_1 =
		(\lambda_0)^\sigma\leq \kappa_{\lambda_1}\). Obviously
		\(\kappa_{\lambda_1}\leq \lambda_1\), so \(\kappa_{\lambda_1} =
		\lambda_1\). Therefore \(\lambda_1\) is measurable.
	\end{proof}
\end{lma}

	\begin{proof}[Proof of \cref{ThreshThm}] By induction, we may assume that
		the theorem holds for all strong limit cardinals \(\bar \lambda <
		\lambda\).
		
		Suppose \(\alpha < \lambda\). We claim that there is a countably
		complete ultrafilter \(D\) with \(\lambda_D < \lambda\) such that
		\(j_D(\kappa) > \alpha\). To see this, fix an ultrapower embedding \(j :
		V\to M\) such that \(j_D(\kappa) > \alpha\). Then by \cref{Exponential},
		one can find a countably complete ultrafilter \(D\) such that
		\(\lambda_D \leq 2^{|\alpha|} < \lambda\) and an elementary embedding
		\(k : M_D\to M\) such that \(k\circ j_D = j\) and \(\textsc{crt}(k) >
		\alpha\). Since \(k(j_D(\kappa)) = j(\kappa) >\alpha = k(\alpha)\), by
		the elementarity of \(k\), \(j_D(\kappa) >\alpha\).
		
		Next, we show that \(\lambda\) is a limit of Fr\'echet cardinals.
		Suppose \(\delta\) is a cardinal with \(\kappa\leq \delta < \lambda\).
		We will find a Fr\'echet cardinal in the interval \((\delta,\lambda)\).
		By the previous paragraph, there is a countably complete ultrafilter
		\(D\) such that \(j_D(\kappa) \geq (2^\delta)^+\) and \(\lambda_D <
		\lambda\). On the other hand \(\delta < \lambda_D\) since \(2^\delta <
		|j_D(\kappa)| \leq \kappa^{\lambda_D} = 2^{\lambda_D}\). Thus
		\(\lambda_D\) is a Fr\'echet cardinal in the interval
		\((\delta,\lambda)\), as desired.
		
		We claim that every regular cardinal in the interval
		\([\kappa,\lambda)\) is Fr\'echet. By \cref{ImplicationLma}, it suffices
		to show that there are no isolated cardinals in the interval
		\([\kappa,\lambda)\). Suppose \(\lambda_0 \in [\kappa,\lambda)\) is
		isolated. Let \(\lambda_1= (\lambda_0)^\sigma\). \cref{DoubleSigma}
		implies that \(\lambda_1\) is measurable. Since \(\lambda\) is a limit
		of Fr\'echet cardinals, \(\lambda_1 < \lambda\). Note that for all
		\(\alpha < \lambda_1\), there is an ultrapower embedding \(j : V\to M\)
		such that \(j(\kappa) >\alpha\), so the \(\lambda_1\)-threshold
		\(\kappa'\) is less than \(\lambda_1\). By our induction hypothesis,
		\(\kappa'\) is \(\gamma\)-supercompact for all \(\gamma < \lambda_1\).
		This contradicts that \(\lambda_1 = (\lambda_0)^\sigma\) is not a limit
		of Fr\'echet cardinals.
		
		We finally claim that \(\kappa\) is \(\delta\)-supercompact for any
		successor cardinal \(\delta\in (\kappa,\lambda)\), which proves the
		theorem. Suppose \(\delta \in (\kappa,\lambda)\) is a successor
		cardinal. Then \(\kappa_\delta\) is \(\delta\)-supercompact by
		\cref{SuccessorThm}.  Since \(\kappa_\delta\) is the limit of the
		isolated cardinals below \(\delta\) (\cref{KappaChar}), \(\kappa_\delta
		\leq \kappa\). On the other hand, by \cref{IdemThreshold}, \(\kappa\) is
		the \(\kappa\)-threshold, so in particular, no \(\nu < \kappa\) is
		\(\kappa\)-supercompact. Hence \(\kappa_\delta\not < \kappa\). It
		follows that \(\kappa =\kappa_\delta\), as desired.
\end{proof}

\subsection{The number of countably complete ultrafilters}
We close this section with a result that just barely uses the analysis of
\(\mathscr K_\lambda\) given by \cref{EmbeddingChar} and
\cref{IsolatedInternal}. Recall that \(\mathscr B(X)\) denotes the set of
countably complete ultrafilters on \(X\). The main result is a bound on the
cardinality of \(\mathscr B(X)\):
\begin{thm}[UA]\label{UFCounting}
	For any set \(X\), \(|\mathscr B(X)| \leq (2^{|X|})^+\).
\end{thm}

The theorem is proved by a generalizing Solovay's \cref{GeneralizedCounting}. To
do this, we need to generalize the notion of the Mitchell rank of an
ultrafilter:

\begin{defn}
	Suppose \(\delta\) is an ordinal and \(W\) is a countably complete
	ultrafilter on \(\delta\).
	\begin{itemize}
		\item \(\mathscr B_W(\delta)\) denotes the set of countably complete
		ultrafilters \(U\) on \(\delta\) such that \(U\sE W\).
		\item \(\sigma(W)\) denotes the rank of \((\mathscr B_W(\delta),\sE)\).
		\item \(\sigma(\delta)\) denotes the rank of \((\mathscr
		B(\delta),\sE)\).
		\item \(\sigma({<}\delta) = \sup_{\alpha < \delta} \sigma(\alpha)+1\).
	\end{itemize}
\end{defn}

Since the Ultrapower Axiom implies that the Ketonen order is linear, the rank of
an ultrafilter completely determines its position in the Ketonen order:

\begin{lma}[UA]\label{RankKet}
	Suppose \(U\) and \(W\) are countably complete ultrafilters on ordinals.
	Then \(U\E W\) if and only if \(\sigma(U) \leq \sigma(W)\).\qed
\end{lma}

The following lemma relates \(\sigma^V\) to \(\sigma^{M_U}\):

\begin{lma}[UA]\label{SigmaTrans}
	Suppose \(U\) is a countably complete ultrafilter and \(W\) is a countably
	complete ultrafilter on an ordinal \(\delta\). Then \(\sigma(W)\leq
	\sigma^{M_U}(\tr U W)\).
	\begin{proof}
		It suffices to show that there is a Ketonen order preserving embedding
		from \(\mathscr B_W(\delta)\) to \(\mathscr B^{M_U}_{\tr U
		W}(j_U(\delta))\). By \cref{OrderPreserving}, the translation function
		\(t_U\) restricts to such a function.
	\end{proof}
\end{lma}

We briefly mention that a version of \cref{SigmaTrans} is provable in ZFC.
Suppose \(Z\) is a countably complete ultrafilter and \(W\) is an ultrafilter on
an ordinal \(\delta\). If \(\langle W_i : i\in I\rangle\) is sequence of
countably complete ultrafilters on \(\delta\) such that \(W =
Z\text{-}\lim_{i\in I} W_i\), then \[\sigma(W)\leq [\langle \sigma(W_i) : i \in
I\rangle]_Z\] We omit the proof, which is an application of \cref{StrongTrans}.

\begin{cor}[UA]\label{FixI}
	Suppose \(U\) is a countably complete ultrafilter and \(W\) is a countably
	complete ultrafilter on an ordinal. If \(j_U(\sigma(W)) = \sigma(W)\) then
	\(U\I W\).
	\begin{proof}
		Assume \(j_U(\sigma(W)) = \sigma(W)\). Then \[\sigma^{M_U}(j_U(W)) =
		j_U(\sigma(W)) = \sigma(W)\leq \sigma^{M_U}(\tr U W)\] For the final
		inequality, we use \cref{SigmaTrans}. By \cref{RankKet}, it follows that
		\(j_U(W) \E \tr U W\) in \(M_U\). By the theory of the internal relation
		(\cref{IChar}), this implies \(U\I W\).
	\end{proof}
\end{cor}

\begin{lma}[UA]\label{Moving}
	Suppose \(\gamma\) is an ordinal. Then for any ordinal \(\xi \in
	[\sigma({<}\gamma),\sigma(\gamma))\), there is a countably complete tail
	uniform ultrafilter \(U\) on \(\gamma\) with \(j_U(\xi) > \xi\).
	\begin{proof}
		Let \(U\) be unique element of \(\mathscr B(\gamma)\) with \(\sigma(U) =
		\xi\). Since \(\xi \geq \eta\), \(U\) does not concentrate on \(\alpha\)
		for any \(\alpha < \eta\). Therefore \(U\) is a nonprincipal tail
		uniform ultrafilter on \(\gamma\). %(As a consequence of
		\cref{PrincipalOrder}, the principal ultrafilter \(D\) on \(\gamma\)
		concentrated at \(\alpha < \gamma\) satisfies \(\sigma(D) =
		\sigma(\alpha)< \sigma({<}\gamma)\).) Since \(U\) is nonprincipal,
		\(U\not \I U\). Therefore \(j_U(\sigma(U)) > \sigma(U)\) by \cref{FixI}.
		In other words, \(j_U(\xi) > \xi\).
	\end{proof}
\end{lma}

The following fact is ultimately equivalent to \cref{IPoint} below:

\begin{lma}[UA]\label{UnFix}
	Suppose \(\xi\) and \(\delta\) are ordinals and \(U\) is the \(\sE\)-minimum
	countably complete ultrafilter on \(\delta\) such that \(j_U(\xi) > \xi\).
	Then for any countably complete ultrafilter \(D\) such that \(j_D(\xi) =
	\xi\), \(D\I U\).
	\begin{proof}
		Since \(j_D\) is elementary and \(j_D(\xi) = \xi\), \(j_D(U)\) is the
		\(\sE^{M_D}\)-minimum countably complete ultrafilter \(Z\) of \(M_D\) on
		\(j_D(\delta)\) such that \(j_Z^{M_D}(\xi) > \xi\). On the other hand,
		\(\tr D U\) is a countably complete ultrafilter of \(M_D\) on
		\(j_D(\delta)\) such that \[j_{\tr D U}^{M_D}(\xi) = j_{\tr D
		U}^{M_D}(j_D(\xi)) = j_{\tr U D}^{M_U}(j_U(\xi)) \geq j_U(\xi) > \xi\]
		Hence by the linearity of the Ketonen order, \(j_D(U)\E \tr D U\) in
		\(M_D\). Now the basic theory of the internal relation (\cref{IChar})
		implies that \(D\I U\).
	\end{proof}
\end{lma}

The central combinatorial argument of \cref{UFCounting} appears in the following
proposition:

\begin{prp}[UA]\label{FrechetBound}
	Suppose \(\lambda\) is a Fr\'echet cardinal. Then for any ordinal \(\gamma <
	\lambda\), \(|\mathscr B(\gamma)| \leq 2^\lambda\).
	\begin{proof}
		Assume towards a contradiction that \(\lambda\) is the least Fr\'echet
		cardinal at which the theorem fails. In particular, \(\lambda\) is not a
		limit of Fr\'echet cardinals, so by \cref{EmbeddingChar} or
		\cref{IsolatedInternal}, \(\mathscr K_\lambda\) is \(\lambda\)-internal.
		Let \(\gamma < \lambda\) be the least ordinal such that \(|\mathscr
		B(\gamma)| > 2^\lambda\). Then in particular, \(\gamma\) is the least
		ordinal such that \(\sigma(\gamma) \geq (2^\lambda)^+\), so
		\(\sigma({<}\gamma) < (2^\lambda)^+\). 
		
		Let \(\xi\) be an ordinal with the following properties:
		\begin{itemize}
			\item \(\sigma({<}\gamma)\leq \xi < (2^\lambda)^+\).
			\item For all \(\alpha < \gamma\), for all \(D \in \mathscr
			B(\alpha)\), \(j_D(\xi) = \xi\).
			\item \(j_{\mathscr K_\lambda}(\xi) = \xi\).
		\end{itemize}
		To see that such an ordinal \(\xi\) exists, let \(S = \bigcup_{\alpha <
		\gamma} \mathscr B(\alpha)\cup \{\mathscr K_\lambda\}\). Note that \(|S|
		\leq 2^\lambda\) by the minimality of \(\gamma\). For each \(D\in S\),
		the collection of fixed points of \(j_D\) is \(\omega\)-closed unbounded
		in \((2^\lambda)^+\). Therefore the intersection of the fixed points of
		\(j_D\) for all \(D\in S\) is \(\omega\)-closed unbounded in
		\((2^\lambda)^+\).
		
		 Since \(\xi\in [\sigma({<}\gamma),\sigma(\gamma))\), \cref{Moving}
		 implies that there is a countably complete tail uniform ultrafilter
		 \(U\) on \(\gamma\) with \(j_U(\xi) > \xi\). Let \(U\) be the
		 \(\sE\)-least countably complete ultrafilter on \(\gamma\) such that
		 \(j_U(\xi) > \xi\). By \cref{UnFix}, \(U\) is \(\lambda\)-internal, and
		 moreover \(\mathscr K_\lambda \I U\). 
		
		Since \(\lambda_U < \lambda\), \(U\I \mathscr K_\lambda\). Thus \(U\I
		\mathscr K_\lambda\) and \(\mathscr K_\lambda\I U\), so by
		\cref{Commute}, \(U\) and \(\mathscr K_\lambda\) commute. Since \(U\) is
		\(\lambda_U\)-internal and \(\mathscr K_\lambda\) is
		\(\lambda\)-internal, we can apply the converse to Kunen's commuting
		ultrapowers lemma (\cref{IUniformPrp}) to obtain \(U\in
		V_{\kappa_\lambda}\). (Obviously \(\mathscr K_\lambda\) is not in
		\(V_\kappa\) where \(\kappa\) is the completeness of \(U\).) In
		particular \(\gamma < \kappa_\lambda\). But then \(|\mathscr
		B(\gamma)|^+ < \kappa_\lambda\) since \(\kappa_\lambda\) is
		inaccessible. This contradicts that \(\kappa_\lambda \leq\lambda <
		(2^\lambda)^+ \leq \sigma(\gamma) < |\mathscr B(\gamma)|^+\).
	\end{proof}
\end{prp}

The proof above is a bit mysterious, and the situation can be clarified by doing
a bit more work than the bare minimum required to prove the theorem. In fact one
can prove the following. Suppose \(\lambda\) is a Fr\'echet cardinal that is
either regular or isolated. Let \(\xi\) be the first fixed point of
\(j_{\mathscr K_\lambda}\) above \(\kappa_\lambda\). Then for any \(D\sE
\mathscr K_\lambda\), \(j_D(\xi) = \xi\). The \(\sE\)-minimum countably complete
ultrafilter \(U\) on \(\lambda\) such that \(j_U(\xi) > \xi\), if it exists, is
isomorphic the \(\mo\)-least normal fine ultrafilter \(\mathcal U\) on
\(P_{\kappa_\lambda}(\lambda)\) such that \(\mathscr K_\lambda\mo \mathcal U\).
This is related to \cref{TwoUFs}.

Incidentally, \cref{FrechetBound} yields an alternate proof of instances of GCH
from UA plus large cardinals. For example, assume \(|\mathscr B(\kappa)| =
2^{2^\kappa}\), \(|\mathscr B(\kappa^+)|> 2^{(\kappa^+)}\), and \(\kappa^{++}\)
is Fr\'echet. Then \[ 2^{2^\kappa} = |\mathscr B(\kappa)| \leq 2^{(\kappa^+)} <
|\mathscr B(\kappa^+)| \leq 2^{(\kappa^{++})}\] Thus \(2^{2^\kappa} <
2^{(\kappa^{++})}\), and in particular \(2^\kappa < \kappa^{++}\). In other
words, \(2^\kappa = \kappa^+\). (This result is not as strong as
\cref{SuccGCH}.)

\cref{FrechetBound} admits the following self-improvement:

\begin{thm}[UA]\label{FrechetBound2}
	For any Fr\'echet cardinal \(\lambda\), for any \(W\in \mathscr
	B(\lambda)\), \(|\mathscr B_W(\lambda)| \leq 2^\lambda\). Hence
	\(\sigma(\lambda) \leq (2^\lambda)^+\).
	\begin{proof}
		For \(\alpha \leq \lambda\), let \(\mathscr B(\lambda,\alpha) = \{Z\in
		\mathscr B(\lambda):\delta_Z\leq \alpha\}\). By the definition of the
		Ketonen order, every element of \(\mathscr B_W(\lambda)\) is of the form
		\(W\text{-}\lim_{\alpha\in I}U_\alpha\) for some \(I\in W\) and
		\(\langle U_\alpha : \alpha \in I\rangle\in \prod_{\alpha\in I}\mathscr
		B(\lambda,\alpha)\). Thus \(|\mathscr B_W(\lambda)| \leq
		|\coprod_{I\subseteq \lambda}\prod_{\alpha\in I}\mathscr
		B(\lambda,\alpha)|\). It therefore suffices show that the cardinality of
		\(\coprod_{I\subseteq \lambda}\prod_{\alpha\in I}\mathscr
		B(\lambda,\alpha)\) is at most \(2^\lambda\). Since
		\[\textstyle\left|\coprod_{I\subseteq \lambda}\prod_{\alpha\in
		I}\mathscr B(\lambda,\alpha)\right|= 2^\lambda\cdot\sup_{I\subseteq
		\lambda} \prod_{\alpha\in I}|\mathscr B(\lambda,\alpha)|\] it suffices
		to show that \(|\mathscr B(\lambda,\alpha)| \leq 2^\lambda\) for all
		\(\alpha < \lambda\). But there is a one-to-one correspondence between
		\(\mathscr B(\lambda,\alpha)\) and \(\mathscr B(\alpha)\), and by
		\cref{FrechetBound}, \(|\mathscr B(\alpha)| \leq \sigma(\alpha) <
		(2^\lambda)^+\). Thus \(|\mathscr B(\lambda,\alpha)| \leq 2^\lambda\),
		which completes the proof.
	\end{proof}
\end{thm}

We finally prove \(|\mathscr B(X)|\leq (2^{|X|})^+\):

\begin{proof}[Proof of \cref{UFCounting}] For \(A\subseteq X\), let \[\mathscr
	B(X,A) = \{U\in \mathscr B(X) : A\in U\}\] For any \(A\subseteq X\) of
	cardinality \(\lambda\), \(|\mathscr B(X,A)| = |\mathscr B(A)| = |\mathscr
	B(\lambda)|.\) Since every ultrafilter \(U\) concentrates on a set whose
	cardinality is a Fr\'echet cardinal, we have \[\mathscr B(X) = \bigcup
	\{\mathscr B(X,A) : A\subseteq X\text{ and }|A|\text{ is Fr\'echet}\}\]
	Hence \begin{equation}\label{FrechetUnion}|\mathscr B(X)| \leq 2^{|X|}\cdot
	\sup\{|\mathscr B(\lambda)| : \lambda \leq |X| \text{ is
	Fr\'echet}\}\end{equation} By \cref{FrechetBound2}, for any Fr\'echet
	cardinal \(\lambda\) such that \(\lambda\leq |X|\), \(|\mathscr B(\lambda)|
	\leq (2^\lambda)^+ \leq (2^{|X|})^+\). Hence by \cref{FrechetUnion},
	\(|\mathscr B(X)|  \leq 2^{|X|}\cdot (2^{|X|})^+ =  (2^{|X|})^+\), as
	desired.
\end{proof}

\section{Isolation}\label{IsolationSection}
In this section, we study isolated cardinals more deeply. The main objects of
study are nonmeasurable isolated cardinals, yet we have the feeling that a
clever argument might prove that these objects simply do not exist. (See
\cref{IsolationConj}.) So far, however, we have been unable to rule them out. In
this section, we show that there are significant constraints on their structure,
and this allows us to prove the linearity of the Mitchell order on normal fine
ultrafilters from UA without using any cardinal arithmetic assumptions
(\cref{ULinearity}).
\subsection{Isolated measurable cardinals}
Recall that isolated cardinals are Fr\'echet limit cardinals that are not limits
of Fr\'echet cardinals. We begin by giving a full analysis of cardinals that are
limits of Fr\'echet cardinals. This will lead us to a proof of
\cref{IsolationConj} for strong limit isolated cardinals.
\begin{figure}
	\center
	\includegraphics[scale=.65]{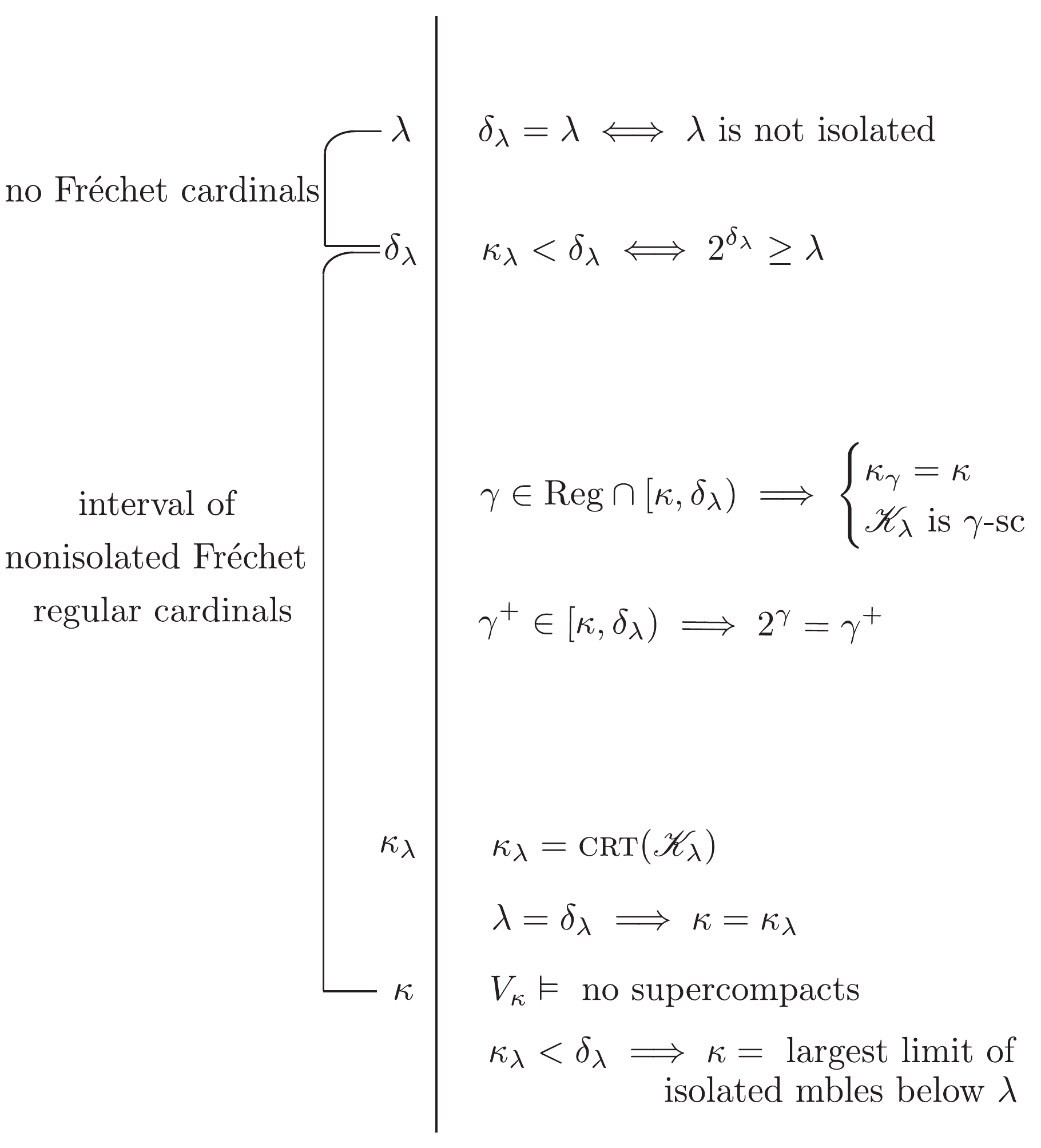}
	\caption{The interval of nonisolated cardinals below a Fr\'echet cardinal.}
\end{figure}
\begin{lma}[UA]\label{KappaChar2}\index{Fr\'echet cardinal!limit of Fr\'echet cardinals}
	Suppose \(\lambda\) is a limit of Fr\'echet cardinals. Let \(\kappa\) be the
	supremum of the isolated cardinals less than \(\lambda\), and assume
	\(\kappa < \lambda\). Then \(\kappa\) is \(\gamma\)-supercompact for all
	\(\gamma < \lambda\). In fact, \(\kappa = \kappa_\iota\) for all regular
	cardinals \(\iota\in [\kappa,\lambda)\).
	\begin{proof}
		Since there are no isolated cardinals in the interval
		\([\kappa,\lambda)\), \cref{ImplicationLma} implies that every regular
		cardinal in the interval \([\kappa,\lambda)\) is Fr\'echet. Assume
		\(\iota\in  [\kappa,\lambda)\) is a regular cardinal. Then \(\iota\) is
		a nonisolated Fr\'echet cardinal. Since \(\kappa\) is the supremum of
		the isolated cardinals below \(\iota\), \(\kappa= \kappa_\iota\) by
		\cref{KappaChar}. Now by \cref{GeneralThm}, \(\kappa\) is
		\(\gamma\)-supercompact for all \(\gamma < \iota\). Since \(\iota\) was
		an arbitrary regular cardinal in \([\kappa,\lambda)\) and \(\lambda\) is
		a limit cardinal, \(\kappa\) is \(\gamma\)-supercompact for all \(\gamma
		< \lambda\).
	\end{proof}
\end{lma}

\begin{cor}[UA]\label{LimitLma}\index{Fr\'echet cardinal!limit of Fr\'echet cardinals}
	Suppose \(\lambda\) is a cardinal. Then the following are equivalent:
	\begin{enumerate}[(1)]
		\item \(\lambda\) is a limit of Fr\'echet cardinals.
		\item Either \(\lambda\) is a limit of isolated measurable cardinals or
		some \(\kappa < \lambda\) is \(\gamma\)-supercompact for all \(\gamma <
		\lambda\).	
	\end{enumerate}
	\begin{proof}
		{\it (1) implies (2):} First assume \(\lambda\) is a limit of isolated
		cardinals. Then by \cref{DoubleSigma}, \(\lambda\) is a limit of
		isolated measurable cardinals.
		
		Assume instead that \(\lambda\) is not a limit of isolated cardinals and
		let \(\kappa < \lambda\) be the supremum of the isolated cardinals below
		\(\lambda\). By \cref{KappaChar2}, \(\kappa\) is \(\gamma\)-supercompact
		for all \(\gamma < \lambda\).
		
		{\it (2) implies (1):} Trivial.
	\end{proof}
\end{cor}

In particular, it follows that every limit of Fr\'echet cardinals is a strong
limit cardinal: if \(\lambda\) is a limit of measurable cardinals, this is
immediate; on the other hand, if some \(\kappa < \lambda\) is
\(\gamma\)-supercompact for all \(\gamma < \lambda\), then \cref{MainTheorem}
implies that for all \(\gamma\in [\kappa,\lambda)\), \(2^\gamma = \gamma^+\).

\begin{lma}[UA]\label{IsoThresh}
	Suppose \(\lambda\) is a strong limit cardinal such that no cardinal
	\(\kappa < \lambda\) is \(\gamma\)-supercompact for all \(\gamma <
	\lambda\). Then for all ultrapower embeddings \(j :V\to M\),
	\(j[\lambda]\subseteq \lambda\). In fact, no ordinal \(\kappa < \lambda\)
	can be mapped arbitrarily high below \(\lambda\) by ultrapower embeddings.
	\begin{proof}
		This is immediate from \cref{ThreshThm}.
	\end{proof}
\end{lma}

The following proposition shows that all the mystery of isolated cardinals falls
away if one assumes the Generalized Continuum Hypothesis.
\begin{prp}[UA]\label{IsolatedStrongLimit}\index{Isolated cardinal!strong limit}
	Suppose \(\lambda\) is cardinal. Then the following are equivalent:
	\begin{enumerate}[(1)]
		\item \(\lambda\) is a strong limit isolated cardinal.
		\item \(\lambda\) is a measurable cardinal, \(\lambda\) is not a limit
		of measurable cardinals, and no cardinal \(\kappa < \lambda\) is
		\(\lambda\)-supercompact.
	\end{enumerate}
	\begin{proof}
		{\it (1) implies (2):} Since \(\lambda\) is not a limit of Fr\'echet
		cardinals, clearly \(\lambda\) is not a limit of measurable cardinals
		and no \(\kappa < \lambda\) is \(\lambda\)-supercompact. It remains to
		show that \(\lambda\) is measurable. Let \(j : V\to M\) be the
		ultrapower of the universe by \(\mathscr K_\lambda\). Note that
		\(j[\lambda]\subseteq\lambda\) by \cref{IsoThresh}. By
		\cref{IsolatedInternal}, \(D\I \mathscr K_\lambda\) for all \(D\) with
		\(\lambda_D < \lambda\). Therefore by \cref{InternalComplete},
		\(\mathscr K_\lambda\) is \(\lambda\)-complete. Since there is a
		\(\lambda\)-complete uniform ultrafilter on \(\lambda\), \(\lambda\) is
		measurable.
		
		{\it (2) implies (1):} Since \(\lambda\) is measurable, \(\lambda\) is a
		strong limit cardinal. It remains to show that \(\lambda\) is isolated.
		Note that no cardinal \(\kappa < \lambda\) is \(\gamma\)-supercompact
		for all \(\gamma < \lambda\): since \(\lambda\) is measurable, this
		would imply \(\kappa\) is \(\lambda\)-supercompact, contrary to
		assumption. \cref{LimitLma} now implies that \(\lambda\) is not a limit
		of Fr\'echet cardinals. Therefore \(\lambda\) is isolated.
	\end{proof}
\end{prp}

The main application of isolated measurable cardinals is factoring ultrapower
embeddings:

\begin{thm}[UA]\label{Cutpoint}
	Suppose \(\kappa\) is a strong limit cardinal that is not a limit of
	Fr\'echet cardinals. Suppose \(U\) is a countably complete ultrafilter. Then
	there is a countably complete ultrafilter \(D\) such that \(\lambda_D <
	\kappa\) admitting an internal ultrapower embedding \(h : M_D\to M_U\) such
	that \(h\circ j_D = j_U\) and \(\textsc{crt}(h) \geq \kappa\) if \(h\) is
	nontrivial.
	\begin{proof}
		Fix \(\gamma < \kappa\) such that \(\kappa = \gamma^\sigma\). By
		 \cref{Exponential}, one can find a countably complete ultrafilter \(D\)
		 such that \(\lambda_D < \kappa\) and there is an elementary embedding
		 \(e : M_D\to M_U\) such that \(\textsc{crt}(e) > \beth_{10}(\gamma)\)
		 and \(e\circ j_D = j_U\).  Let \(\lambda = \lambda_D\). We may assume
		 without loss of generality that \(\lambda\) is the underlying set of
		 \(D\). Since \(\lambda < \kappa\) is a Fr\'echet cardinal, \(\lambda
		 \leq \gamma\). Let \(\lambda' = j_D(\lambda)\). Then \(\lambda'  <
		 (2^\lambda)^+\), so \(2^{2^{\lambda'}} < \beth_{10}(\gamma)\). Since
		 \(e : M_D\to M_U\) has critical point above \(2^{2^{\lambda'}}\),
		 \[P(P(\lambda')) \cap M_D = P(P(\lambda'))\cap M_U\] Thus the following
		 hold where \(\mathscr B(X)\) denotes the set of countably complete
		 ultrafilters on \(X\):
		 \begin{itemize}
		 	\item \(\lambda' = j_D(\lambda) = e(\lambda')  = j_U(\lambda)\).
			\item \(\mathscr B^{M_D} ({\lambda'})= \mathscr
			B^{M_U}({\lambda'})\).
			\item \({\E^{M_D}}\restriction\mathscr B^{M_D}({\lambda'}) =
			{\E^{M_U}} \restriction  \mathscr B^{M_U}({\lambda'}) \)
			\item \(j_D\restriction P(\lambda) = j_U\restriction P(\lambda)\).
		\end{itemize}
	 
	 By \cref{Reciprocity}, \(\tr D D\) is the \(\E^{M_D}\)-least element
	 \(D'\in\mathscr B^{M_D} ({\lambda'})\) such that \(j_D^{-1}[D'] = D\). By
	 \cref{Reciprocity}, \(\tr U D\) is the \(\E^{M_U}\)-least element \(D'\in
	 \mathscr B^{M_U}({\lambda'})\) such that \(j_U^{-1}[D'] = D\). By the
	 agreement set out in the bullet points above, it therefore follows that
	 \(\tr D D = \tr U D\). On the other hand, by \cref{TransRF0}, \(\tr D D\)
	 is principal in \(M_D\). Thus \(\tr U D\) is principal in \(M_U\).
	 Therefore by \cref{TransRF0}, \(D\D U\). 
	 
	 Let \(h : M_D \to M_U\) be the internal ultrapower embedding with \(h\circ
	 j_D = j_U\). By \cref{MinDefEmb}, \(h(\alpha) \leq e(\alpha)\) for all
	 \(\alpha\in \text{Ord}\), so \(\textsc{crt}(h)\geq \textsc{crt}(e) >
	 \beth_{10}(\gamma) > j_D(\gamma)\). Since \(h\) is an internal ultrapower
	 embedding of \(M_D\), if \(h\) is nontrivial then \(\textsc{crt}(h)\) is a
	 measurable cardinal of \(M_D\) above \(j_D(\gamma)\). Since there are no
	 measurable cardinals in the interval \((\gamma,\kappa)\), there are no
	 measurable cardinals of \(M_D\) in the interval
	 \((j_D(\gamma),j_D(\kappa))\). Therefore if \(h\) is nontrivial, then
	 \(\textsc{crt}(h) \geq \kappa\).
 \end{proof}
\end{thm}

\subsection{Ultrafilters on an isolated cardinal}\label{IsolatedUFSection}
In this subsection, which is perhaps the most technical of this dissertation, we
enact a very detailed analysis of the countably complete ultrafilters on an
isolated cardinal \(\lambda\). One of the goals is to prove the following
theorem:
\begin{lma}\label{IsoDecomp}
	Suppose \(\lambda\) is an isolated cardina and \(W\) is a countably complete
	ultrafilter. Then \(\mathscr K_\lambda\D W\) if and only if \(W\) is
	\(\lambda\)-decomposable and \(W\I \mathscr K_\lambda\) is and only if \(W\)
	is \(\lambda\)-indecomposable.
\end{lma}

This should be seen as a generalization of the universal property for \(\mathscr
K_\lambda\) when \(\lambda\) is regular to isolated cardinals \(\lambda\).
 
We begin with the following fact:
\begin{thm}[UA]\label{IsolatedWN}
	Suppose \(\lambda\) is an isolated cardinal. Then \(\mathscr K_\lambda\) is
	the unique  countably complete weakly normal ultrafilter on \(\lambda\).
\end{thm}

It turns out to be easier to prove something that is a priori slightly stronger.
Recall the notion of the Dodd parameter \(p(j)\) of an elementary embedding
\(j\), defined in \cref{DoddParamDef1} in the general context of elementary
embeddings, and once again in \cref{DoddParamDef2} in the more relevant special
case of ultrapower embeddings.
\begin{prp}[UA]\label{OneGenerator}
	Suppose \(\lambda\) is an isolated cardinal. Then \(\mathscr K_\lambda\) is
	the unique countably complete incompressible ultrafilter \(U\) on
	\(\lambda\) such that \(|p(j_U)| = 1\).
	\begin{proof}
		Suppose towards a contradiction that the proposition fails. Let \(U\) be
		the \(\sE\)-least countably complete incompressible ultrafilter on
		\(\lambda\) such that \(p(j_U) = 1\) and \(U \neq \mathscr K_\lambda\).
		Since \(\mathscr K_\lambda\) is the \(\sE\)-least uniform ultrafilter on
		\(\lambda\), \(\mathscr K_\lambda \sE U\).
		
		Let \(j : V\to M\) be the ultrapower of the universe by \(\mathscr
		K_\lambda\) and let \(\nu = a_{\mathscr K_\lambda}\). Let \(i : V\to N\)
		be the ultrapower of the universe by \(U\) and let \(\xi = \id_U\). By
		the incompressibility of \(U\), \(p(j_U) = \{\xi\}\).
		
		Let \((k,h) : (M,N)\to P\) be the pushout of \((j,i)\). Since \(\mathscr
		K_\lambda \sE U\), 
		\begin{equation}\label{SeedBelow}k(\nu) < h(\xi)\end{equation}
		
		We claim that \(h(\xi)\) is a generator of \(k: M\to P\), or in other words that \[h(\xi)\notin H^P(k[M]\cup h(\xi))\] Since \(\xi\) is a generator of \(i\), \(h(\xi)\) is a generator of \(h\circ i\) by \cref{CompositionGenerators}. Since \(k\circ j = h\circ i\), \(h(\xi)\) is a generator of \(k\circ j\).  Since \(M = H^M(j[V]\cup \{\nu\})\), \[H^P(k[M]\cup h(\xi)) = H^P(k\circ j[V]\cup \{k(\nu)\}\cup h(\xi)) = H^P(k\circ j[V]\cup h(\xi))\]
		The final equality follows from \cref{SeedBelow}. Therefore since \(h(\xi)\notin H^P(k\circ j[V]\cup h(\xi))\), \(h(\xi)\notin H^P(k[M]\cup h(\xi)) \), as desired.
		
		Let \(Z = \tr {\mathscr K_\lambda} U\), so \(Z\) is the \(M\)-ultrafilter on \(j(\lambda)\) derived from \(k\) using \(h(\xi)\). Then \(Z\) is a countably complete ultrafilter on \(j(\lambda)\) and \(\id_Z = h(\xi)\) is a generator of \(j_Z^M = k\). 
		
		We claim that \(Z\) is an incompressible ultrafilter on \(j(\lambda)\) in \(M\). Since \(\id_Z = h(\xi)\) is a generator of \(j_Z^M\), it suffices to show that \(Z\) is tail uniform, or in other words, \(\delta_Z = j(\lambda)\). Since \(\id_Z\) is a generator of \(j_Z^M\), \(\delta_Z =\lambda_Z\) is a Fr\'echet cardinal in \(M\). 
		By \cref{SeedBelow}, \(\delta_Z > \id_{\mathscr K_\lambda}\). Since \(U\) is on \(\lambda\), \(\xi < i(\lambda)\), so \(h(\xi) < h(i(\lambda)) = k(j(\lambda))\), which implies \(\delta_Z \leq j(\lambda)\). 
		Thus \(\delta_Z\in (\id_{\mathscr K_\lambda},j(\lambda)]\).
		Since \(\lambda\) is isolated, no Fr\'echet cardinal of \(M\) lies in the interval \([\sup j[\lambda],j(\lambda))\). Therefore \(\delta_Z = j(\lambda)\), as desired.
		
		It follows that in \(M\), \(Z\) is a countably complete incompressible ultrafilter on \(j(\lambda)\). Moreover \(p(j_Z^M) = \{h(\xi)\}\) by \cref{DoddGenerators}, so \(p(j_Z^M)\) has cardinality 1. 
		
		We claim that \(Z \neq j(\mathscr K_\lambda)\). The reason is that \(j^{-1}[Z] = U\) (since \(Z = \tr {\mathscr K_\lambda} U\)) while \(j^{-1}[j(\mathscr K_\lambda)] = \mathscr K_\lambda\).
		
		Thus we have shown that in \(M\), \(Z\) is a countably complete incompressible ultrafilter on \(j(\lambda)\) such that \(|p(j_Z^M)| = 1\) and \(Z \neq j(\mathscr K_\lambda)\). By the definition of \(U\) and the elementarity of \(j\), it follows that \(j(U)\E Z\) in \(M\). \cref{IChar} now implies that \(\mathscr K_\lambda \I U\). But \(j_U\) is discontinuous at \(\lambda\) since \(\lambda_U = \lambda\). Thus by \cref{Nonisolation1}, \(\lambda\) is not isolated. This is a contradiction.
		\end{proof}
\end{prp}

\begin{proof}[Proof of \cref{IsolatedWN}]
	If \(U\) is a countably complete weakly normal ultrafilter on \(\lambda\), then \(U\) is incompressible and \(p(j_U) = \{\id_U\}\) by \cref{WeaklyNormalGenerator}. Therefore we can apply \cref{OneGenerator}. 
\end{proof}

We now investigate the iterated ultrapowers of \(\mathscr K_\lambda\).

\begin{defn}
	If \(\lambda\) is an isolated cardinal, then the {\it iterated ultrapower of \(\mathscr K_\lambda\)} is the iterated ultrapower \[\mathcal I_\lambda = \langle M^\lambda_n,j^\lambda_{mn},U^\lambda_m : m \leq n < \omega\rangle\] formed by setting \(U^\lambda_m = j^\lambda_{0m}(\mathscr K_\lambda)\) for all \(m < \omega\).
	For \(n < \omega\), let \(p^{n\lambda} = p(j^\lambda_{0n})\), and let \(\mathscr K^n_\lambda\) be the ultrafilter on \([\lambda]^\ell\) derived from \(j_{0n}^\lambda\) using \(p^{n\lambda}\) where \(\ell = |p^{n\lambda}|\).
\end{defn}

Thus \(j_{0,n}^\lambda: V\to M_n^\lambda\) is the ultrapower of the universe by \(\mathscr K^n_\lambda\). We now analyze the parameters \(p^n\):

\begin{lma}\label{IsolatedParameter}
	Suppose \(\lambda\) is an isolated cardinal. Let  \(\langle M_n,j_{m,n},U_m : m \leq n < \omega\rangle\) be the iterated ultrapower of \(\mathscr K_\lambda\). For \(n < \omega\), let \(p^n = p^{n\lambda}\).
	Then for all \(n < \omega\), \(|p^n| = n\) and
	\begin{align}\label{PowerParameter}
		p^{n+1} \restriction n &= j_{01}(p(j_{0n})) \\
	p^{n+1}_n &= j_{01}(j_{0n})(\id_{\mathscr K_\lambda})\nonumber\end{align}
	
	\begin{proof}
		Note that the conclusion of the lemma holds when \(n = 0\). Assume \(m \geq 1\) and that the conclusion of the lemma holds when \(n = m - 1\). We will prove that the conclusion of the lemma holds when \(n = m\).
		\begin{figure}
			\center
			\includegraphics[scale=.6]{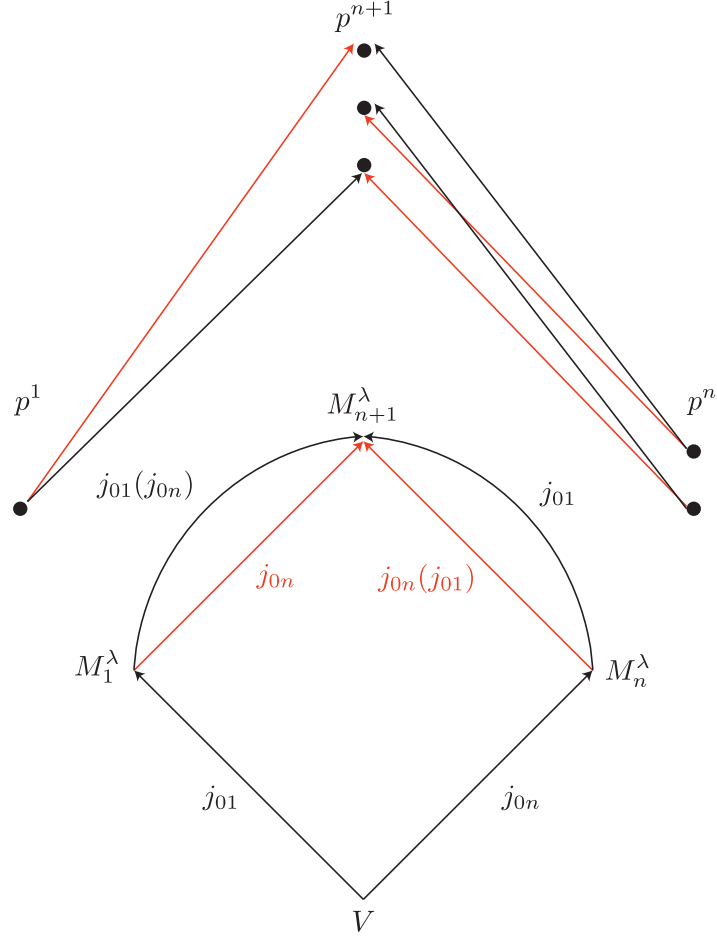}
			\caption{The iterated ultrapower of \(\mathscr K_\lambda\).}
		\end{figure}
		
		Note that \[j_{0m+1} = j_{1m+1}\circ j_{01} = j_{01}(j_{0m})\circ j_{01} = j_{01}\circ j_{0m}\]
		
		Since \(p^m\) is the Dodd parameter of \(j_{0m}\), \(p^{m+1}\) is the Dodd parameter of \(j_{0m+1}\), and \(j_{0m+1} = j_{01}\circ j_{0m}\), \(j_{01}(p^m) < p^{m+1}\) in the parameter order, and hence \(j_{01}(p^m) \leq p^{m+1}\restriction m\).
		
		Let \(W\) be the \(M_m\)-ultrafilter derived from \(j_{01}\) using \(j_{01}(j_{0m})(\id_{\mathscr K_\lambda})\). Then by the basic theory of the internal relation (\cref{PushUlt}),
		\[W = s_{\mathscr K_\lambda^m}(\mathscr K_\lambda)\] 
		\(j_W^{M_m} = j_{01}\restriction M_n\), \(\id_W = j_{01}(j_{0m})(\id_{\mathscr K_\lambda})\), and 
		\[M_{n+1} = H^{M_{n+1}}(j_{0n+1}[V]\cup  j_{01}(p^m)\cup \{j_{01}(j_{0m})(\id_{\mathscr K_\lambda}) \})\]
		It follows from the minimality of the Dodd parameter that \begin{equation}\label{ParamMinimal}p^{m+1}\leq j_{01}(p^m)\cup \{j_{01}(j_{0m})(\id_{\mathscr K_\lambda}) \}\end{equation}
		 By our induction hypothesis, \[\min p^m =
		j_{01}(j_{0m-1})(\id_{\mathscr K_\lambda}) = j_{1m}(\id_{\mathscr
		K_\lambda}) \geq \sup j_{1m}[\sup j_{01}[\lambda]] = \sup
		j_{0m}[\lambda]\] Therefore \(j_{0m}(\xi) < \min p^m\) for all \(\xi <
		\lambda\), so by \L o\'s's Theorem, \[j_{01}(j_{0m})(\id_{\mathscr
		K_\lambda}) < \min j_{01}(p^m)\] Combining this with
		\cref{ParamMinimal}, we can conclude that \(p^{m+1}\restriction m\leq
		j_{01}(p^m)\).
		
		Putting these two inequalities together, we have shown
		\(p^{m+1}\restriction m = j_{01}(p^m)\).
		
		By \cref{DoddGenerators}, to show \cref{PowerParameter}, it suffices to
		show that \(j_{01}(j_{0m})(\id_{\mathscr K_\lambda})\) is the largest
		\(j_{01}(p^m)\)-generator of \(j_{0m+1}\). By \cref{ParamMinimal}, it in
		fact suffices to show that \(j_{01}(j_{0m})(\id_{\mathscr K_\lambda})\)
		is a \(j_{01}(p^m)\)-generator of \(j_{0m+1}\). 
		
		We claim that an ordinal \(\xi\) is a \(j_{01}(p^m)\)-generator of
		\(j_{0m+1}\) if and only if \(\xi\) is a generator of
		\(j_{01}\restriction M_m\). This follows immediately from the fact that
		\[H^{M_{m+1}}(j_{0m+1}[V]\cup j_{01}(p^m) \cup \xi) =
		H^{M_{m+1}}(j_{01}[j_{0m}[V]\cup p^m]\cup \xi) =
		H^{M_{m+1}}(j_{1m}[M_m]\cup \xi)\] for any ordinal \(\xi\).
		
		Thus to finish, we must show that \(j_{01}(j_{0m})(\id_{\mathscr
		K_\lambda})\) is a generator of \(j_{01}\restriction M_m\).
		
		Let \(\lambda' = \sup j_{0m}[\lambda]\). We first show that \(\lambda_W
		= \lambda'\). By the definition of \(s_{\mathscr K^n_\lambda}(\mathscr
		K_\lambda)\), \(\lambda'\in W\): note that \(j_{0m}^{-1}[\lambda'] =
		\lambda\in \mathscr K_\lambda\). It follows that \(\lambda_W
		\leq\lambda'\). Thus we are left to show that \(\lambda'\leq
		\lambda_W\). Assume to the contrary that there is a set \(B\in W\) such
		that for some \(\kappa < \lambda\), letting \(\delta = |B|^{M_m}\),
		\(\delta < j_{0m}(\kappa)\). Then \(j_{0m}\) is
		\((\lambda,\delta)\)-tight, and it follows that \(j_{0m}\) is
		discontinuous at all regular cardinals in the interval
		\([\kappa,\lambda]\). (See the proof of \cref{SuccessorPrp}.) Therefore
		\(\lambda\) is a limit of Fr\'echet cardinals, which contradicts that
		\(\lambda\) is isolated.
		
		Since \(\lambda_W = \lambda'\), \(j_{01}\restriction M_m\) must have a
		generator in the interval \([\sup j_{01}[\lambda'],j_{01}(\lambda'))\).
		Let \(\xi\) be the least such generator. Then \[\xi \leq
		j_{01}(j_{0m})(\id_{\mathscr K_\lambda})\] by \cref{ParamMinimal}. Let
		\(U\) be the ultrafilter derived from \(j_{0m+1}\) using \(\xi\) and let
		\(k :M_U\to M_{m+1}\) be the factor embedding with \(k\circ j_U=
		j_{0m+1}\) and \(k(\id_U) = \xi\). Clearly \(\xi\) is a generator of
		\(j_{0m+1}\), so \(\id_U\) is a generator of \(j_U\). Therefore
		\(\lambda_U = \delta_U = \lambda\), since \(\sup j_{0m+1}[\lambda] =
		\sup j_{01}[\lambda']\leq \xi\). Thus \(U\) is a uniform countably
		complete ultrafilter on \(\lambda\) and \(p(j_U) = \{\id_U\}\), so by
		\cref{OneGenerator}, \(U = \mathscr K_\lambda\).
		
		 Thus \(k\) and \(j_{01}(j_{0m})\) are elementary embeddings from
		 \(M_1\) to \(M_{m+1}\). Since \(j_{01}(j_{0m})\) is an internal
		 ultrapower embedding of \(M_1\), it follows from \cref{MinDefEmb} that
		 \[j_{01}(j_{0m})(\id_{\mathscr K_\lambda}) \leq k(\id_{\mathscr
		 K_\lambda}) = k(\id_U) = \xi\] Putting these inequalities together,
		 \(\xi = j_{01}(j_{0m})(\id_{\mathscr K_\lambda})\), and therefore
		 \(j_{01}(j_{0m})(\id_{\mathscr K_\lambda})\) is a generator of
		 \(j_{01}\restriction M_m\), as desired.
	\end{proof}
\end{lma}

A key parameter in the theory of Fr\'echet cardinals is the strict cardinal
supremum of a cardinal's Fr\'echet predecessors:
\begin{defn}\label{DeltaLambda}
	For any cardinal \(\lambda\), \index{\(\delta_\lambda\) (Fr\'echet
	supremum)}\index{Fr\'echet cardinal!\(\delta_\lambda\)} \(\delta_\lambda =
	\sup \{\eta^+ : \eta < \lambda\text{ and }\eta\text{ is Fr\'echet}\}\).
\end{defn}
If \(\lambda\) is a Fr\'echet cardinal, then \(\lambda\) is isolated if and only
if \(\delta_\lambda < \lambda\).

We have the following immediate corollary:

\begin{lma}[UA]\label{IsolatedParameter2}
	Suppose \(\lambda\) is an isolated cardinal. Let  \(\langle M_n,j_{m,n},U_m
	: m \leq n < \omega\rangle\) be the iterated ultrapower of \(\mathscr
	K_\lambda\). Suppose \(i :V\to N\) is an ultrapower embedding of the form
	\(i = d\circ j_{0n}\) where \(d : M_n\to N\) is the ultrapower of \(M_n\) by
	a countably complete ultrafilter \(D\) of \(M_n\) with \(\lambda_D <
	j_{0n}(\delta_\lambda)\). Then \(p(i)\setminus i(\delta_\lambda) =
	d(p^{n,\lambda})\).
\end{lma}

The following theorem amounts to a complete analysis of the ultrafilters on an
isolated cardinal:
\begin{thm}[UA]\label{IsoFactor}
	Suppose \(\lambda\) is an isolated cardinal. Let  \(\langle M_n,j_{m,n},U_m
	: m \leq n < \omega\rangle\) be the iterated ultrapower of \(\mathscr
	K_\lambda\). Suppose \(i: V\to N\) is the ultrapower by a countably complete
	ultrafilter on \(\lambda\). Then for some \(n < \omega\), \(i = d\circ
	j_{0n}\) where \(d : M_n\to N\) is the ultrapower of \(M_n\) by a countably
	complete ultrafilter \(D\) of \(M_n\) with \(\lambda_D <
	j_{0n}(\delta_\lambda)\).
	\begin{proof}
		Suppose \(U\) is a countably complete ultrafilter on \(\lambda\). Assume
		by induction that the proposition holds when \(i  = j_W\) for an
		ultrafilter \(W\sE U\).
		
		Let \(i : V\to N\) be the ultrapower of the universe by \(U\), and we
		will show that the theorem is true for \(i\).
		
		If \(\lambda_U < \lambda\), then the theorem is vacuously true.
		Therefore we may assume \(\lambda_U = \lambda\). 
		
		Let \(j : V\to M\) be the ultrapower of the universe by \(\mathscr
		K_\lambda\). Let \(\nu = \id_{\mathscr K_\lambda}\).
		
		Let \((k,h) : (M,N)\to P\) be the pushout of \((j,i)\). Since \((k,h)\)
		is the pushout of \((j,i)\), \(k\) is the ultrapower embedding of \(M\)
		associated to \(\tr {\mathscr K_\lambda} U\). Since \(\lambda_U =
		\lambda\), \(j_U\) is discontinuous at \(\lambda\). Hence by
		\cref{Nonisolation1}, \(\mathscr K_\lambda \not \I U\). Therefore by
		\cref{IChar}, \(\tr {\mathscr K_\lambda} U \sE j(U)\) in \(M\). We can
		now apply our induction hypothesis, shifted by \(j\) to \(M\), to the
		ultrafilter \(\tr {\mathscr K_\lambda} U\) of \(M\). We conclude that
		for some \(\ell < \omega\), \(k = d\circ j(j_{0\ell})\) where \(d:
		j(M_\ell)\to P\) is the ultrapower of \(j(M_\ell)\) by a countably
		complete ultrafilter \(D\) of \(j(M_\ell)\) such that \(\lambda_D <
		j(j_{0\ell})(j(\delta_\lambda))\). Let \[n = \ell+1\] Then
		\(j(j_{0\ell}) = j_{1n}\), \(j(M_\ell)  = M_n\), and
		\(j(j_{0\ell})(j(\delta_\lambda)) = j_{1n}(j(\delta_\lambda)) =
		j_{0n}(\delta_\lambda)\). Thus \[k = d\circ j_{1n}\] where \(d : M_n\to
		P\) is the ultrapower of \(M_n\) by a countably complete ultrafilter
		\(D\) of \(M_n\) such that \(\lambda_D < j_{0n}(\delta_\lambda)\). Note
		that \[k\circ j = d\circ j_{0n}\] has the form we want to show that
		\(i\) has.
		
		Let \(p^\ell = p^{\ell,\lambda}\) and let \(p^n = p^{n,\lambda}\), so
		that
		\begin{equation}\label{IsoParamEq}p^n = j(p^\ell)\cup \{j_{1n}(\nu)\}\end{equation} by \cref{IsolatedParameter} and the fact that \(j_{1n} = j_{01}(j_{0\ell})\).
		
		Let \(q' = p(k)\setminus k(j(\delta_\lambda))\). By \cref{IsolatedParameter2} applied in \(M\), \(q' = d(j(p^{\ell}))\). Since \(k\circ j = d\circ j_{0n}\),
		 \begin{align}
		 	p(k\circ j ) \setminus k(j(\delta_\lambda)) &= p(d\circ j_{0n}) \setminus d (j_{0n}(\delta_\lambda)) \nonumber\\
		 	&= d(p^n) \label{IsoParam2App}\\
		 	&= d(j(p^{\ell})\cup \{j_{1n}(\nu)\})\label{IsoParam1App}\\ 
		 	&= q'\cup \{k(\nu)\}\label{LikeTermsApp}
		 \end{align}
	 	Here \cref{IsoParam2App} follows from \cref{IsolatedParameter2}; \cref{IsoParam1App} follows from \cref{IsoParamEq}; \cref{LikeTermsApp} follows from the fact that \(d(j(p^{\ell})) = q'\) and \(d\circ j_{1n} = k\).
		
		Let \(\xi\) be the least generator of \(i\) such that \(\sup i[\lambda]\leq \xi < i(\lambda)\). 
		
		\begin{clm} \label{DivisionClm}\(k(\nu) = h(\xi)\). \end{clm}
		\begin{proof}[Proof of \cref{DivisionClm}]
		By \cref{WeaklyNormalGenerator}, the ultrafilter derived from \(i\) using \(\xi\) is a countably complete weakly normal ultrafilter on \(\lambda\), and hence is equal to \(\mathscr K_\lambda\). Let \(e : M \to N\) be the factor embedding with \(e(\nu) = \xi\) and \(e\circ j = i\). The comparison \((e,\text{id})\) witnesses \((j,\nu) \E (i,\xi)\).
		Since \((j,\nu) \E (i,\xi)\), we must have \(k(\nu) \leq h(\xi)\). 
		
		Assume towards a contradiction that \(k(\nu) \neq h(\xi)\), so \(k(\nu) < h(\xi)\).
		
		Let \(q = p(i)\setminus \sup i[\lambda]\).
		 We claim that \(h(q) = p(k)\restriction |q|\). The proof is by induction. Assume \(m < |q|\) and \(h(q)\restriction m = p(k)\restriction m\). 
		 By \cref{DoddGenerators}, \(q_m\) is the largest \(q\restriction m\)-generator of \(i\). 
		 Hence \(h(q_m)\) is the largest \(h(q\restriction m)\)-generator of \(h\circ i\). 
		 Replacing like terms, \(h(q_m)\) is the largest \(p(k)\restriction m\)-generator of \(k\circ j\). Since \(q_m\) is a generator of \(i\) above \(\sup i[\lambda]\), \(q_m \geq \xi\). Hence \(h(q_m) \geq h(\xi) > k(\nu)\) by our assumption that \(h(\xi) > k(\nu)\). 
		 Therefore \(h(q_m)\) is not only a \(p(k)\restriction m\)-generator of \(k\circ j\) but also a \(p(k)\restriction m\cup \{k(\nu)\}\)-generator of \(k\circ j\). In other words, \(h(q_m)\) is a \(p(k)\restriction m\)-generator of \(k\), and it must therefore be the largest \(p(k)\restriction m\)-generator of \(k\). By \cref{DoddGenerators}, \(h(q_m) = p(k)_m\).
		 
		 Since \(q\) has no elements below \(\sup i[\lambda]\), in particular, \(q\) has no elements below \(i(\delta_\lambda)\). Therefore \(h(q)\) has no elements below \(h(i(\delta_\lambda)) = k(j(\delta_\lambda))\). Since \(h(q)\subseteq p(k)\) by the previous paragraph, it follows that \(h(q)  \subseteq p(k)\setminus k(j(\delta_\lambda))= q'\).
		
		We now claim that \(k(\nu)\) is a generator of \(h\). To show this, it suffices to show that \(k(\nu)\) is a \(h(p(i))\)-generator of \(h\circ i\). Let \(r = p(i)\cap \sup i[\lambda]\). Thus \(p(i) = q\cup r\). Note that \(h(r)\subseteq \sup h\circ i[\lambda] = \sup k\circ j[\lambda]\leq k(\nu)\), since \(\sup j[\lambda]\leq \nu\). Hence \(h(r)\subseteq k(\nu)\). Thus to show that \(k(\nu)\) is a \(k(p(i))\)-generator of \(h\circ i\), it suffices to show that \(k(\nu)\) is a \(h(q)\)-generator of \(h\circ i\). Since \(h(q) \subseteq q'\), it suffices to show that \(k(\nu)\) is a \(q'\)-generator of \(k\circ j\). This is an immediate consequence of \cref{LikeTermsApp}: by \cref{DoddGenerators}, \(k(\nu)\) is the largest \(q'\)-generator of \(k\circ j\).
		
		Thus \(k(\nu)\) is a generator of \(h\). Let \(W\) be the tail uniform \(N\)-ultrafilter derived from \(h\) using \(k(\nu)\). Then \(W\) is an incompressible ultrafilter. We have \(\sup i[\lambda]\leq \delta_W\) since \(\sup h[\sup i [\lambda]] = \sup k\circ j[\lambda] \leq k(\nu)\). Moreover \(\delta_W \leq \xi\) since \(k(\nu) < h(\xi)\). Since \(W\) is incompressible, \(\lambda_W = \delta_W\). But \(\lambda_W\) is a Fr\'echet cardinal of \(N\) and \(\sup i[\lambda]\leq \lambda_W \leq \xi <i(\lambda)\). This contradicts that \(\lambda\) is isolated.
		
		It follows that our assumption that \(k(\nu) \neq h(\xi)\) was false. This proves \cref{DivisionClm}.
		\end{proof}
	Since \(k(\nu) = h(\xi)\), it follows from \cref{SeedEq} that \(\mathscr K_\lambda\D U\). Let \(k' :M\to N\) be an internal ultrapower embedding. Then \((k',\text{id})\) is a pushout of \((j,i)\). By the uniqueness of pushouts, \(k = k'\) and \(h = \text{id}\). Hence \(k : M\to N\) is the unique internal ultrapower embedding such that \(k\circ j = i\). Thus \(i = k\circ j = d\circ j_{0n}\). Since \(d : M_n\to N\) is the ultrapower of \(M_n\) by a countably complete ultrafilter \(D\) of \(M_n\) with \(\lambda_D < j_{0n}(\lambda)\), this proves the proposition.
	\end{proof}
\end{thm}

\begin{cor}[UA]\label{IsoFactor2}
Suppose \(\lambda\) is an isolated cardinal.  Let  \(\langle M_n,j_{m,n},U_m : m \leq n < \omega\rangle\) be the iterated ultrapower of \(\mathscr K_\lambda\). Then for any ultrapower embedding \(k : V\to P\), there is some \(n < \omega\) such that \[ k = h\circ d\circ j_{0n}\] where \(M_n\stackrel d\longrightarrow N\stackrel h\longrightarrow P\) are ultrapower embeddings with the following properties:
\begin{itemize}
	\item \(d : M_n\to N\) is the ultrapower of \(M_n\) by a countably complete ultrafilter \(D\) of \(M_n\) with \(\lambda_D < j_{0n}(\delta_\lambda)\).
	\item \(h : N\to P\) is an internal ultrapower embedding of \(N\) with \(\textsc{crt}(h) > d(j_{0n}(\lambda))\) if \(h\) is nontrivial.
\end{itemize}
\begin{proof}
	We claim there is a strong limit cardinal \(\kappa > \lambda\) such that there are no Fr\'echet cardinals in the interval \((\lambda,\kappa)\). If there are no Fr\'echet cardinals above \(\lambda\), let \(\kappa = \beth_\omega(\lambda)\). Otherwise, let \(\kappa = \lambda^\sigma\). By \cref{DoubleSigma}, \(\kappa\) is measurable, and in particular, \(\kappa\) is a strong limit cardinal.
	
	By \cref{Cutpoint}, there is a countably complete ultrafilter \(U\) with \(\lambda_U < \kappa\) such that there is an internal ultrapower embedding \(h : M_U\to P\) with \(h\circ j_U = j\) and \(\textsc{crt}(h) \geq \kappa \). Since \(\lambda_U < \kappa\) is Fr\'echet and there are no Fr\'echet cardinals in the interval \([\lambda, \kappa]\), \(\lambda_U \leq \lambda\). Therefore we may assume that \(U\) is a countably complete ultrafilter on \(\lambda\). In particular \( \textsc{crt}(h)\geq \kappa > j_U(\lambda)\).
	
	 Let \(i = j_U\). By \cref{IsoFactor}, for some \(n < \omega\), \(i = d\circ j_{0n}\) where \(d : M_n\to N\) is the ultrapower of \(M_n\) be a countably complete ultrafilter \(D\) of \(M_n\) with \(\lambda_D < j_{0n}(\delta_\lambda)\). Putting everything together, \[j= h\circ d\circ j_{0n}\] and this proves the corollary.
\end{proof}
\end{cor}

It is not a priori obvious that \(p^n\) contains all the generators \(\xi\) of \(\mathscr K_\lambda^n\) with \(\xi \geq\sup j_{0n}[\lambda]\). In fact this is true:

\begin{prp}[UA]\label{IsoPowerChar}
	Suppose \(\lambda\) is an isolated cardinal and \(n < \omega\). Then \(\mathscr K_\lambda^n\) is the unique countably complete ultrafilter \(W\) on \([\lambda]^n\) such that \(\id_W\) is the set of generators  \(\xi\) of \(j_W\) with \(\xi \geq\sup j_W[\lambda]\).
	\begin{proof}
		Assume by induction that the corollary is true when \(n = m\), and we will prove it when \(n = m + 1\).
		
		Therefore assume \(W\) is a countably complete ultrafilter on \([\lambda]^{m+1}\) such that \(\id_W\) is the set of generators \(\xi\) of \(j_W\) with \(\xi\geq \sup j_W[\lambda]\). Let \(q\) be the first \(m\) generators of \(j_W\) above \(\sup j_W[\lambda]\). Let \(U\) be the ultrafilter derived from \(j_W\) using \(q\). Then by our induction hypothesis, \(U = \mathscr K_\lambda^{m}\). Let \(d : M_{m}\to M_W\) be the factor embedding with \(d\circ j_{0m} = j_W\) and \(d(p^{m}) = q\). By \cref{IsoFactor}, there is an internal ultrapower embedding \(d' : M_m\to M_W\). Note that \(d'(p^m)\) is a set of generators of \(d'\circ j_{0m} = d\circ j_{0m}\), so \(d'(p^m) \geq d(p^m)\). On the other hand, \(d'(p^m) \leq d(p^m)\) by \cref{MinDefEmb}. Thus \(d' (p^m) = d(p^m)\). Since \(d'\circ j_{0m} = d\circ j_{0m}\), we have \(d' = d\). Thus \(d\) is an internal ultrapower embedding. 
		
		Let \(\xi\) be the largest generator of \(j_W\). Thus \(d(q)\subseteq \xi\), so \(\xi\) is a \(d(q)\)-generator of \(j_W\) and hence \(\xi\) is a generator of \(d\). Let \(Z\) be the tail uniform \(M_m\)-ultrafilter derived from \(d\) using \(\xi\). Then \(Z\) is an incompressible ultrafilter of \(M_m\) and \(\delta_Z\in [\sup j_{0m}[\lambda],j_{0m}(\lambda)]\). Since \(\delta_Z = \lambda_Z\) is a Fr\'echet cardinal of \(M_m\), it follows that \(\delta_Z = j_{0m}(\lambda)\). Therefore by \cref{IsolatedWN}, \(Z = j_{0m}(\mathscr K_\lambda)\).
		
		Since \(M_W = H^{M_W}(j_W[V]\cup q\cup \{\xi\}) = H^{M_W}(d[M_m]\cup \{\xi\})\), we have \(d = j^{M_m}_Z = j_{mm+1}\). Thus \(d\circ j_{0m} = j_{0m+1}\). Thus \(j_W = j_{0m+1}\).
		
		Since \(p^{m+1}\) consists solely of generators of \(j_{0m+1}\) above \(\sup j_{0m+1}[\lambda]\), \(p^{m+1}\subseteq \id_W\). Since \(|\id_W| = |p^{m+1}|\), it follows that \(p^{m+1} = \id_W\). Therefore \(W = \mathscr K^{m+1}_\lambda\), as desired.
	\end{proof}
\end{prp}

\begin{prp}[UA]
	Suppose \(\lambda\) is an isolated cardinal. Then \(\mathscr K_\lambda^n\) is the unique countably complete ultrafilter \(W\) on \([\lambda]^n\) such that \(\id_W\) is a set of generators of \(j_W\) disjoint from \(\sup j_W[\lambda]\).
	\begin{proof}
		Suppose \(W\) is such an ultrafilter. Let \(p\) be the set of all generators of \(\xi\) of \(j_W\) with \(\xi \geq \sup j_W[\lambda]\). Let \(m = |p|\). By \cref{IsoPowerChar}, the ultrafilter derived from \(j_W\) using \(p\) is \(\mathscr K_\lambda^m\). It follows that \(j_W = j_{0m}\) and \(p = p^m\). Therefore \(p \leq \id_W\) by the minimality of the Dodd parameter. On the other hand, \(\id_W \subseteq p\) since \(p\) consists of all the generators of \(j_W\) above \(\sup j_W[\lambda]\). Therefore \(\id_W = p\). Hence \(m = n\) and \(W = \mathscr K_\lambda^n\), as desired.
	\end{proof}
\end{prp}

We also have an analog of \cref{UFAmenableChar} at isolated cardinals:

\begin{thm}[UA]\label{IsolatedUFAmenableChar}
	Suppose \(\lambda\) is an isolated cardinal. Let \(j : V\to M\) be the ultrapower of the universe by \(\mathscr K_\lambda\) and let \(\nu = \id_{\mathscr K_\lambda}\). Suppose \(Z\) is a countably complete \(M\)-ultrafilter that is \(\delta\)-indecomposable for all \(M\)-cardinals \(\delta\in [\sup j[\lambda],\nu]\). Then \(Z\in M\).
	\begin{proof}
		Let \(e : M\to P\) be the ultrapower of \(M\) by \(Z\). Then \(e(\nu)\) is a generator of \(e\circ j\) by \cref{CompositionGenerators}. Since \(Z\) is \(\delta\)-indecomposable for all \(\delta\in [\sup j[\lambda],\nu]\), \(e\) has no generators in the interval \([\sup e\circ j[\lambda], e(\nu)]\). In other words, \(e\circ j\) has no \(e(\nu)\)-generators in the interval \([\sup e\circ j[\lambda], e(\nu)]\). 
		
		Let \(k = e\circ j\). Applying \cref{IsoFactor2}, there is some \(n < \omega\) such that \[ k = h\circ d\circ j_{0n}\] where \(M_n\stackrel d\longrightarrow N\stackrel h\longrightarrow P\) are ultrapower embeddings with the following properties:
		\begin{itemize}
			\item \(d : M_n\to N\) is the ultrapower of \(M_n\) by a countably complete ultrafilter \(D\) of \(M_n\) with \(\lambda_D < j_{0n}(\delta_\lambda)\).
			\item \(h : N\to P\) is an internal ultrapower embedding of \(N\) with \(\textsc{crt}(h) > d(j_{0n}(\lambda))\) if \(h\) is nontrivial.
		\end{itemize}
	
		Let \(e' = h\circ d\circ j_{1n}\), so that \(e' : M\to P\) is an internal ultrapower embedding with \(e' \circ j = k = e\circ j\).
		We claim \(e'(\nu) = e(\nu)\). By \cref{MinDefEmb}, \(e'(\nu) \leq e(\nu)\). 
		
		Suppose towards a contradiction \(e'(\nu) < e(\nu)\). Then \(e'(\nu)\) is not an \(e(\nu)\)-generator of \(e\circ j = e'\circ j\). 
		Note that \(h(d(j_{1n}(\nu))) = e'(\nu)\) and \(h(e(\nu)) = e(\nu)\), so \(d(j_{1n}(\nu))\) is not an \(e(\nu)\)-generator of \(d\circ j_{0n}\). But consider the ultrafilter \(U\) on \([\lambda]^2\) derived from \(d\circ j_{0n}\) using \(\{d(j_{1n}(\nu)),e(\nu)\}\). Since \(d(j_{1n}(\nu))\) and \(e(\nu)\) are  generators \(\xi\) of \(d\circ j_{0n}\) with \(\xi \geq \sup d\circ j_{0n}[\lambda]\), \(\id_U\) consists of generators \(\xi\) of \(j_U\) with \(\xi \geq \sup j_U[\lambda]\). Thus \(U = \mathscr K_\lambda^2\). But then by \cref{IsolatedParameter}, \(\min(\id_U)\) is a \(j_U\)-generator. This contradicts that \(d(j_{1n}(\nu))\) is not an \(e(\nu)\)-generator.
	\end{proof}
\end{thm}

The main application of \cref{IsolatedUFAmenableChar} is the following fact:
\begin{lma}[UA]\label{IsolatedUFAmenableChar2}
	Assume \(\lambda\) is isolated and let \(j : V\to M\) be the ultrapower of the universe by \(\mathscr K_\lambda\). Either \(j[\lambda]\subseteq \lambda\) or \(\mathscr K_\lambda\cap M\in M\).
\begin{proof}
	Assume \(\sup j[\lambda] > \lambda\). Then \(\mathscr K_\lambda\cap M\) is not \(\gamma\)-decomposable for any \(M\)-cardinal \(\gamma\in [\sup j[\lambda],j(\lambda))\). Therefore \(\mathscr K_\lambda\cap M\in M\).
\end{proof}
\end{lma}
\cref{IsolatedUFAmenableChar} gives a coarse bound on the strength of \(\mathscr K_\lambda\) when \(\lambda\) is isolated.
\begin{prp}[UA]\label{IsoStrength}
	Suppose \(\lambda\) is isolated and let \(j : V\to M\) be the ultrapower of the universe by \(\mathscr K_\lambda\). Then \(P(\lambda)\subseteq M\) if and only if \(\mathscr K_\lambda\) is \(\lambda\)-complete.
	\begin{proof}
	Assume \(P(\lambda)\subseteq M\). Since \(\mathscr K_\lambda\notin M\), \(\mathscr K_\lambda\cap M\notin M\), so by \cref{IsolatedUFAmenableChar2}, \(\sup j[\lambda]\subseteq \lambda\). By the Kunen Inconsistency Theorem (\cref{KunenInconsistency}), this implies \(\textsc{crt}(j)\geq \lambda\). In other words, \(\mathscr K_\lambda\) is \(\lambda\)-complete.
\end{proof}
\end{prp}
Assume \(\lambda\) is a nonmeasurable isolated cardinal. One would expect to get a much better bound on the strength of \(\mathscr K_\lambda\) than \(\lambda\). When \(\delta_\lambda\) is a successor cardinal, one can in fact prove that \(P(\delta_\lambda)\not\subseteq M_{\mathscr K_\lambda}\), which determines the strength of \(j_{\mathscr K_\lambda}\) exactly (since \(j_{\mathscr K_\lambda}\) is \({<}\delta_\lambda\)-supercompact by \cref{DeltaLemma} below).  When \(\delta_\lambda\) is inaccessible, however, we do not know whether \(P(\delta_\lambda)\subseteq M_{\mathscr K_\lambda}\) is possible.
\subsection{Isolated cardinals and the GCH}
By \cref{IsolatedStrongLimit}, the existence of isolated cardinals that are not measurable is paired with failures of the Generalized Continuum Hypothesis. In this section, we study precisely how GCH fails below a nonmeasurable isolated cardinal. Here the cardinal \(\delta_\lambda\) (see \cref{DeltaLambda}) takes the center stage.

\begin{prp}[UA]\label{DeltaLemma}
	Suppose \(\lambda\) is an isolated cardinal that is not measurable.  Let \(j : V\to M\) be the ultrapower of the universe by \(\mathscr K_\lambda\). Let \(\kappa = \kappa_\lambda\) and \(\delta = \delta_\lambda\). Then the following hold:
	\begin{enumerate}[(1)]
		\item Every regular cardinal \(\iota\in [\kappa,\delta)\) is Fr\'echet.
		\item \(j\) is \({<}\delta\)-supercompact.
		\item If \(\delta\) is a limit cardinal then \(\delta\) is strongly inaccessible.
		\item Otherwise \(\delta\) is the successor of a cardinal \(\gamma\) of cofinality at least \( \kappa_\lambda\). In fact, no cardinal in the interval \((\textnormal{cf}(\gamma),\gamma)\) is \(\gamma\)-strongly compact.
	\end{enumerate}
\begin{proof}
	We first prove (1). Let \(\eta \in [\iota,\delta)\) be a Fr\'echet cardinal. Then for any \(\gamma\in [\kappa,\eta)\), then there is a Fr\'echet cardinal in \((\gamma,\eta]\). By \cref{NonisolatedInterval}, there are no isolated cardinals in \([\kappa,\lambda)\). \cref{ImplicationLma} implies that every regular cardinal in \([\kappa,\eta)\) is Fr\'echet. In particular, \(\iota\) is Fr\'echet.

	We now prove (2). Fix a regular cardinal \(\iota\in [\kappa,\delta)\), and we will show that \(j\) is \(\iota\)-supercompact. (This suffices since the  Recall that there are no isolated cardinals in \([\kappa,\lambda)\) (\cref{NonisolatedInterval}). Thus \(\kappa_\iota \leq \kappa\) as a consequence of \cref{KappaChar}. Moreover, by \cref{NonisolatedCompact}, \(\kappa_\iota\) is \(\iota\)-strongly compact. We can therefore apply our technique for converting amenability of ultrafilters into strength (\cref{GeneralStrong}) to conclude that \(P(\iota)\subseteq M\): \(\kappa_\iota\) is \(\iota\)-strongly compact, \(M\) is closed under \(\kappa_\iota\)-sequences, and every countably complete ultrafilter on \(\iota\) is amenable to \(M\) (\cref{IsolatedInternal}), so \(P(\iota)\subseteq M\).
	
	 By \cref{KetonenTight}, \(j_{\mathscr K_\iota}\) is \(\iota\)-tight. Moreover \(j_{\mathscr K_\iota}(\kappa) \geq j_{\mathscr K_\iota}(\kappa_\iota) > \iota\). By \cref{IsolatedInternal}, \(\mathscr K_\iota\I \mathscr K_\lambda\). We now use the following fact:
	 \begin{lma*}
	 	Suppose \(\kappa\leq \iota\) are cardinals, \(U\) and \(W\) are countably complete ultrafilters, \(U\) is \(\iota\)-tight, \(j_U(\kappa) > \iota\), \(W\) is \(\kappa\)-complete, and \(U\I W\). Then \(j_W\) is \(\iota\)-tight.
	 	\begin{proof}
	 		Since \(j_U(W)\) is \(\iota\)-complete in \(M_U\), \(\text{Ord}^\iota\cap M_U\subseteq M_{j_U(W)}^{M_U} = j_U(M_W)\subseteq M_W\). (The final containment uses \(U\I W\).) Therefore since \(M_U\) has the \({\leq}\iota\)-covering property, so does \(M_W\). Thus \(j_W\) is \(\iota\)-tight.
	 	\end{proof}
	 \end{lma*}
	  We can apply the fact to \(U = \mathscr K_\iota\) and \(W = \mathscr K_\lambda\). It follows that \(j\) is \(\iota\)-tight. Since \(j\) is \(\iota\)-tight and \(P(\iota)\subseteq M\), \(j\) is \(\iota\)-supercompact.
	 
	 We now prove (3). Suppose towards a contradiction that \(\delta\) is singular. Then by (2), \(j\) is \(\delta\)-supercompact. If \(\text{cf}(\delta) \geq \kappa_\lambda\), it follows that \(\delta\) is Fr\'echet, contrary to the definition of \(\delta_\lambda\). Therefore \(\text{cf}(\delta) < \kappa_\lambda\). But then by \cref{SmallCfSC}, \(j\) is \(\delta^+\)-supercompact. Then \(\delta^+\) is Fr\'echet. The definition of \(\delta\) implies that no cardinal in \([\delta,\lambda)\) is Fr\'echet, so it must be that \(\delta^+ = \lambda\). This contradicts that \(\lambda\) is isolated (and in particular is a limit cardinal).
	 
	 For (4), assume towards a contradiction that some cardinal \(\nu\) in the interval \((\textnormal{cf}(\gamma),\gamma)\) is \(\gamma\)-strongly compact. Then \(\nu\) it is \(\gamma^+\)-strongly compact by \cref{SmallCfSC}. But \(\gamma^+ = \delta\) is not Fr\'echet, and this is contradiction.
\end{proof}
\end{prp}

Suppose \(\lambda\) is an isolated cardinal, and let \(\delta = \delta_\lambda\). Must \(2^{<\delta} = \delta\)? By \cref{DeltaLemma} (3), this is true if \(\delta\) is a limit cardinal, but we are unable to answer the question when \(\delta\) is a successor. The following bound is sufficient for most applications:
\begin{thm}\label{LowerBound}
	Suppose \(\lambda\) is isolated and \(\delta = \delta_\lambda\). Then \(2^{<\delta} < \lambda\).
	\end{thm}
\begin{proof}
	 Assume by induction that the theorem holds for all isolated cardinals below \(\lambda\). Let \(j: V\to M\) be the ultrapower of the universe by \(\mathscr K_\lambda\). Then \(j\) is \({<}\delta\)-supercompact (\cref{DeltaLemma}). Thus \(2^{<\delta}\leq (2^{<\delta})^M\), so it suffices to show that \((2^{<\delta})^M < \lambda\).
	 
	 \begin{clm}\label{SmallSuccessorClm}
	 \((\delta^\sigma)^M \leq \lambda\).
	 \end{clm}
 	\begin{proof}[Proof of \cref{SmallSuccessorClm}]
 	There are two cases. 
 	
 	First assume \(\sup j[\lambda] = \lambda\). Since \(j\) is \({<}\delta\)-supercompact, Kunen's Inconsistency Theorem (\cref{Kunen}) implies that there is a measurable cardinal \(\iota <\delta\) such that \(j(\iota) > \delta\). Now \(j(\iota) < \lambda\) is a measurable cardinal of \(M\), so \((\delta^\sigma)^M \leq j(\iota) < \lambda\), as desired. 
 	
 	Assume instead that \(\lambda < \sup j[\lambda]\). Then \(\mathscr K_\lambda\cap M\in M\) by \cref{IsolatedUFAmenableChar}. Thus \(\lambda\) is Fr\'echet in \(M\), so \((\delta^\sigma)^M \leq \lambda\).
 	\end{proof}
 
 	If \(\delta^{+M}\) is Fr\'echet in \(M\), then \((2^{<\delta})^M = \delta\) by \cref{SuccGCH}. Assume therefore that \(\delta^{+M}\) is not Fr\'echet in \(M\). Let \(\eta = (\delta^\sigma)^M\). Then \(\eta\) is isolated in \(M\) by \cref{SuccessorPrp}. Moreover \(\eta \leq \lambda < j(\lambda)\), so our induction hypothesis shifted to \(M\) applies at \(\eta\). Notice that \(\delta \leq (\delta_\eta)^M\): indeed, by \cref{DeltaLemma}, \(M\) is correct about cardinals below \(\delta\), and by \cref{IsolatedInternal}, all sufficiently large cardinals below \(\delta\) are Fr\'echet in \(M\). Thus \[(2^{<\delta})^M \leq (2^{{<}\delta_\eta^M})^M < \eta \leq \lambda\] In particular \((2^{<\delta})^M < \lambda\), as desired.
\end{proof}

The following closely related fact can be seen as an ultrafilter-theoretic version of SCH:
\begin{prp}[UA]\label{IsoFix}
	Suppose \(\lambda\) is a regular isolated cardinal. Suppose \(D\) is a countably complete ultrafilter such that \(\lambda_D < \lambda\). Then \(j_D(\lambda) = \lambda\).
	\begin{proof}
		Suppose towards a contradiction that the theorem fails, and let \(\lambda\) be the least counterexample. Let \(j : V\to M\) be the ultrapower by \(\mathscr K_\lambda\). Let \(\delta = \delta_\lambda\) be the strict supremum of the Fr\'echet cardinals below \(\lambda\). By \cref{DeltaLemma}, \(M^{<\delta}\subseteq M\), and by \cref{IsolatedInternal}, \(M\) satisfies that there is a countably complete ultrafilter \(D\) is \(\lambda_D < \delta\) such that \(j_D(\lambda) \neq \lambda\).
		
		Suppose first that \(\lambda < \sup j[\lambda]\). Then \(\mathscr K_\lambda\cap M\in M\) by \cref{IsolatedUFAmenableChar}. Therefore \(\lambda\) is a regular Fr\'echet cardinal in \(M\). Clearly \(\lambda\) is a limit cardinal in \(M\). Since \(\lambda < j(\lambda)\), \(\lambda\) is not a counterexample to the proposition in \(M\). Therefore \(\lambda\) is not isolated in \(M\), so \(\lambda\) is strongly inaccessible in \(M\) by \cref{LimitLma}. But this contradicts that there is a countably complete ultrafilter \(D\) is \(\lambda_D < \delta\) such that \(j_D(\lambda) \neq \lambda\).
		
		Suppose instead that \(\lambda = \sup j[\lambda]\). Let \(\kappa = \kappa_\lambda\).  We claim that for any countably complete ultrafilter \(U\in V_\kappa\), \(j_U(\lambda) = \lambda\). Fix such an ultrafilter \(U\). Since \(2^{<\delta} < \lambda\), \(j_U(\delta) < \lambda\). By elementarity there are no Fr\'echet cardinals of \(M_U\) in the interval \([j_U(\delta),j_U(\lambda))\). But \(\mathscr K_\lambda \I U\) (by Kunen's commuting ultrapowers lemma, \cref{KunenCommute}), so \(\mathscr K_\lambda\cap M_U\in M_U\), and hence \(\lambda\) is Fr\'echet in \(M_U\). Thus \(\lambda\) is a Fr\'echet cardinal of \(M_U\) in the interval  \([j_U(\delta),j_U(\lambda)]\), so we must have \(j_U(\lambda) = \lambda\), as claimed.
		
		Let \(\eta\) be the least ordinal such that for some ultrafilter \(D\) with \(\lambda_D <\lambda\), \(j_Z(\eta) > \lambda\). (Note that \(\eta\) exists since \(\lambda\) is regular.) %We claim that for any \(\alpha < \lambda\), there is a countably complete ultrafilter \(Z\) on \(\delta\) such that \(j_Z(\eta) \in (\alpha, \lambda)\). Assume not, towards a contradiction. Let \(\xi\) be a fixed point of \(j\) such that there is no \(Z\) on \(\delta\) with \(j_Z(\eta)\in (\xi, \lambda)\). Let \(D\) be the \(\sE\)-least ultrafilter such that \(j_D(\eta) > \xi\). Then by our assumption on \(\xi\), \(j_D(\eta) > \lambda\).) We must have that \(j^M_{\tr {\mathscr K_\lambda} D}(j_D(\eta)) = j^{M_D}_{\tr D {\mathscr K_\lambda}}\circ j_D(\eta) \geq j_D(\eta) > \xi =  j(\xi)\), and hence \(j(D) \E \tr {\mathscr K_\lambda} D\) in \(M\) by the definition of \(D\). It follows that \(\mathscr K_\lambda \I D\) by \cref{IChar}. Since \(j_D(\eta) > \lambda\), in particular \(j_D(\lambda) > \lambda\), and this contradicts the previous paragraph. Thus our assumption was false, and indeed for any \(\alpha < \lambda\), there is a countably complete ultrafilter \(Z\) such that \(j_Z(\eta) \in (\alpha, \lambda)\). 
		We claim \(j(\eta) = \eta\). To see this, note that if \(D\) is an ultrafilter with \(\delta\), then \(j_D[\eta]\subseteq \eta\). (Otherwise we would contradict the minimality of \(\eta\) as in \cref{IdemThreshold}.) If \(j(\eta) > \eta\), however, then since \(j(\eta) < \lambda\), there is an ultrafilter \(D\) on \(\delta\) such that \(j_D(\eta) > j(\eta)\). Since \(Z\I \mathscr K_\lambda\), this contradicts that \(M\) thinks \(j(\eta)\) is closed under ultrapower embeddings associated to ultrafilters on \(j(\delta)\).
		
		Suppose \(\xi\) is a fixed point of \(j\). Let \(\gamma\) be the least cardinal that carries a countably complete ultrafilter \(U\) such that \(j_U(\eta) > \xi\). Then \(\gamma < \delta\) by assumption. We claim \(j(\gamma) = \gamma\). The reason is that \(M\) is closed under \(\gamma\) sequences and contains every ultrafilter on \(\gamma\), so \(M\) satisfies that there is an ultrafilter \(U\) on \(\gamma\) such that \(j_U(\eta) > \xi\). Since \(j(\eta) = \eta\) and \(j(\xi) = \xi\), it follows that \(j(\gamma)\) is the least \(M\) cardinal carrying such an ultrafilter \(U\), and hence \(j(\gamma) = \gamma\). Since \(\gamma\) is a fixed point of \(j\) below its supercompactness, \(\gamma < \kappa\) by the Kunen inconsistency theorem.
		
		It follows that \(\eta\) is mapped arbitrarily high below \(\lambda\) by ultrafilters in \(V_\kappa\). Since \(\lambda\) is regular, there must be a single ultrafilter \(U\in V_\kappa\) such that \(j_U(\eta) \geq \lambda\). This contradicts that for all \(U\in V_\kappa\), \(j_U(\lambda) = \lambda\).
	\end{proof} 
\end{prp}

Our next theorem shows that the problematic isolated cardinals \(\lambda\) suffer a massive failure of GCH precisely at \(\delta_\lambda\):

\begin{thm}\label{UpperBound}
	Suppose \(\lambda\) is a nonmeasurable isolated cardinal and \(\delta = \delta_\lambda\). Then \(2^\delta \geq \lambda\).
\end{thm}

It is not clear whether it is possible that \(2^\delta = \lambda\). This of course implies that \(\lambda\) is regular and hence weakly Mahlo by \cref{SigmaMahlo} below.

This theorem requires an analysis of indecomposable ultrafilters due to Silver. His analysis can be seen as an improvement of \cref{Exponential} in a the key special case of indecomposable ultrafilters.

\begin{thm}[Silver]\label{Silver}\index{Indecomposable ultrafilter!Silver's Theorem}
Suppose \(\delta\) is a regular cardinal and \(U\) is a countably complete ultrafilter that is \(\lambda\)-indecomposable for all \(\lambda\in [\delta,2^\delta]\). Then there is an ultrafilter \(D\) with \(\lambda_D < \delta\) such that there is an elementary embedding \(k : M_D\to M_U\) with \(j_U = k\circ j_D\) and \(\textsc{crt}(k) > j_D((2^{\delta})^+)\) if \(k\) is nontrivial.
\end{thm}

The proof does not really use that \(U\) is countably complete, and this was important in Silver's original work. Since we only need the theorem when \(U\) is countably complete, we make this assumption. (This is for notational convenience: the notion of the critical point of \(k\) does not really make sense if \(M_D\) is illfounded.) 

We begin by describing a correspondence between partitions of ultrafilters and points in the ultrapower embedding that is implicit in Silver's proof.

\begin{defn}
	Suppose \(P\) is a partition of a set \(X\) and \(A\) is a subset of \(X\). Then the {\it restriction of \(P\) to \(A\)} is the partition \(P\restriction A\) defined by \[P\restriction A = \{A\cap S : S\in P\text{ and }A\cap S\neq \emptyset\}\]
\end{defn}

\begin{defn}
	Suppose \(U\) is an ultrafilter on a set \(X\) and \(\lambda\) is a cardinal.
	\begin{itemize}
		\item \(\mathbb Q_U\) denotes the preorder on the collection of partitions of \(X\) defined by setting \(P\leq Q\) if there exists some \(A\in U\) such that \(Q\restriction A\) refines \(P\restriction A\).
		\item \(\mathbb Q_U^\lambda\subseteq \mathbb Q_U\) consists of those \(P\) such that \(|P\restriction A| < \lambda\) for some \(A\in U\).
		\item  \(\mathbb P_U\) denotes the preorder on \(M_U\) defined by setting \(x \leq y\) if \(x\) is definable in \(M\) from \(y\) and parameters in \(j_U[V]\).
		\item \(\mathbb P_U^\lambda\subseteq \mathbb P_U\) is the restriction of \(\mathbb P_U\) to \(H^{M_U}(j_U[V]\cup\sup j_U[\lambda])\). 
	\end{itemize}
\end{defn}

The following lemma, which is ultimately just an instance of the correspondence between partitions of \(X\) and surjective functions on \(X\), shows that the preorders \(\mathbb Q_U\) and \(\mathbb P_U\) are equivalent preorders:

\begin{lma}\label{Correspondence}
	Suppose \(U\) is an ultrafilter on a set \(X\). For \(P\in \mathbb Q_U\), let \(\Phi(P)\) be the unique \(S\in j_U(P)\) such that \(\id_U\in S\). Then the following hold:
	\begin{enumerate}[(1)]
		\item \(\Phi\) is order-preserving: for any \(P,Q\in \mathbb Q_U\), \(P\leq Q\) if and only if \(\Phi(P)\leq \Phi(Q)\).
		\item \(\Phi\) is surjective on equivalence classes: for any \(x\in \mathbb P_U\), there is some \(P\in \mathbb Q_U\) such that \(x\) and \(\Phi(P)\) are equivalent in \(\mathbb P_U\).
		\item For any cardinal \(\lambda\), \(\Phi[\mathbb Q^\lambda_U]\subseteq \mathbb P_U^\lambda\).
		\item Suppose \(P\in \mathbb Q_U\). Let \(D = \{A\subseteq P : \bigcup A\in U\}\). Then there is a unique elementary embedding \(k : M_D\to M_U\) such that \(k\circ j_D = j_U\) and \(k(\id_D) = \Phi(D)\).
	\end{enumerate}
	\begin{proof}
		{\it Proof of (1):} Suppose \(P,Q\in \mathbb Q_U\) and \(P\leq Q\). Fix \(A\in U\) such that \(Q\restriction A\) refines \(P\restriction A\). Then \(\Phi(P)\) is definable in \(M_U\) from the parameters \(\Phi(Q), j_U(P), j_U(A)\) as the unique \(S\in j_U(P)\) such that \(\Phi(Q)\cap j_U(A)\subseteq S\cap j_U(A)\). In other words, \(\Phi(P)\leq \Phi(Q)\).
		
		Conversely suppose \(\Phi(P)\leq \Phi(Q)\), so that \(\Phi(P) = j_U(f)(\Phi(Q))\) for some \(f: Q\to P\). Let \(A\subseteq X\) consist of those \(x\in X\) such that \(x\in f(S)\) where \(S\) is the unique element of \(Q\) with \(x\in S\). Then \(A\in U\) since \(\id_U\in j_U(f)(S)\) where \(S = \Phi(Q)\) is the unique \(S\in j_U(Q)\) such that \(\id_U\in S\). Moreover for any \(S\in Q\), \(S\cap A\subseteq f(S)\cap A\), so \(Q\restriction A\) refines \(P\restriction A\). In other words, \(P\leq Q\).
		
		{\it Proof of (2):} Fix \(x\in \mathbb P_U\). Fix \(f : X\to V\) such that \(x = j_U(f)(\id_U)\). Let \[P = \{f^{-1}[\{y\}] : y\in \text{ran}(f)\}\] Then \(\Phi(P)\) is interdefinable with \(x\) over \(M_U\) using parameters in \(j_U[V]\): \(\Phi(P)\) is the unique \(S\in j_U(P)\) such that \(x\in j_U(f)[S]\); and since \(j_U(f)[\Phi(P)] = \{x\}\), \(x = \bigcup j_U(f)[\Phi(P)]\).
		
		{\it Proof of (3):} Suppose \(P\in \mathbb Q^\lambda_U\). Fix \(\delta < \lambda\) and a surjection \(f : \delta\to P\). Then \(\Phi(P) = j_U(f)(\xi)\) for some \(\xi < j_U(\delta) \leq \sup j_U[\lambda]\). Hence \(\Phi(P) \in H^{M_U}(j_U[V]\cup \sup j_U[\lambda])\), as desired.
		
		{\it Proof of (4):} Define \(g : X\to P\) by setting \(g(a)\) equal to the unique \(S\in P\) such that \(a\in S\). Then \(g_*(U) = D\) and \(j_U(g)(\id_U) = \Phi(P)\). Therefore by the basic theory of the Rudin-Keisler order (\cref{RFChar}), there is a unique elementary embedding \(k : M_D\to M_U\) with \(k\circ j_D = j_U\) and \(k(\id_D) = \Phi(P)\).
	\end{proof}
\end{lma}

For Silver's theorem, it is useful to reformulate indecomposability in terms of \(\mathbb Q_U\):

\begin{lma}
	Suppose \(U\) is an ultrafilter on \(X\) and \(\lambda\) is a cardinal. Then \(U\) is {\it \(\lambda\)-indecomposable} if every partition of \(X\) into \(\lambda\) pieces is equivalent in \(\mathbb Q_U\) to a partition of \(X\) into fewer than \(\lambda\) pieces.\qed
\end{lma}

We now prove Silver's theorem.

\begin{proof}[Proof of \cref{Silver}]
		Let \((\mathbb Q,\leq) = \mathbb Q^{(2^\delta)^+}_U\) be the preorder of \(U\)-refinement on the set of partitions of \(X\) of size at most \(2^\delta\). Let \(\preceq\) be the preorder of refinement on \(\mathbb Q\), so \(P \preceq Q\) implies \(Q\) refines \(P\). Thus \((\mathbb Q,\leq)\) extends \((\mathbb Q,\preceq)\).
		
		Note that \(\preceq\) is \({\leq}\delta\)-directed. Indeed, suppose \(\mathcal S \subseteq \mathbb Q\) has cardinality \(\delta\). Then \[P = \left\{\bigcap \mathcal C : \mathcal C\text{ meets each element of }\mathcal S\text{ and }\bigcap \mathcal C\neq\emptyset\right\}\] refines every partition in \(\mathcal S\), and \(|P| \leq |\prod \mathcal S| \leq 2^\delta\). The partition \(P\) is called the {\it least common refinement} of \(\mathcal S\).
		
		We claim that \((\mathbb Q,\leq)\) has a maximum element (up to equivalence). Since \((\mathbb Q,\preceq)\) is directed, \((\mathbb Q,\leq)\) is directed, and thus it suffices to show that \((\mathbb Q,\leq)\) has a maximal element. Assume the contrary, towards a contradiction. Then since \((\mathbb Q,\preceq)\) is \({\leq}\delta\)-directed, we can produce a sequence \(\langle P_\alpha : \alpha \leq \delta\rangle\) of elements of \(\mathbb Q\) such that for all \(\alpha < \beta \leq \delta\), \(P_\alpha \preceq P_\beta\) and \(P_\beta \not \leq P_\alpha\).
		
		Since \(U\) is \(\lambda\)-indecomposable for all \(\lambda\in [\delta,2^\delta]\), there is some \(A\in U\) such that \(|P_\delta\restriction A| < \delta\).
		For each \(\alpha \leq \delta\), let \(Q_\alpha = P_\alpha\restriction A\). We use the following general fact:
		\begin{clm*} Suppose \(\delta\) is a regular cardinal, \(A\) is a set of size less than \(\delta\), and \(\langle Q_\alpha: \alpha < \delta\rangle\) is a sequence of partitions of \(A\) such that for all \(\alpha < \beta < \delta\), \(Q_\beta\) refines \(Q_\alpha\). Then for all sufficiently large \(\alpha < \beta < \delta\), \(Q_\alpha = Q_\beta\).
			\begin{proof}
				  Let \(Q\) be the least common refinement of \(\{Q_\alpha : \alpha < \delta\}\). Suppose \(S\in Q\). We claim that \(S\in Q_\alpha\) for some \(\alpha < \delta\). Consider the sequence \(\langle S_\alpha : \alpha < \delta\rangle\) where \(S_\alpha\in Q_\alpha\) is the unique element of \(Q_\alpha\) containing \(S\). Thus \(S = \bigcap_{\alpha < \delta} S_\alpha\). Note that \(\langle S_\alpha : \alpha < \delta\rangle\) is a decreasing sequence of sets, each of cardinality less than \(\delta\). Thus for all sufficiently large \(\alpha < \delta\), \(S_\alpha = S\), and in particular, \(S\in Q_\alpha\).
			  
			  	For each \(S\in Q\), fix \(\alpha_S < \delta\) such that \(S\in Q_{\alpha_S}\). Let \(\gamma  = \sup_{S\in Q}\alpha_S\). Then \(\gamma < \delta\) since \(|Q| < \delta\) and \(\delta\) is regular. By definition, \(Q\subseteq Q_\gamma\), so \(Q_\gamma = Q\), If \(\alpha \in [\gamma,\delta)\), then \(Q\) refines \(Q_\alpha\) which refines \(Q_\gamma = Q\), and hence \(Q = Q_\alpha\). This proves the claim.
		  \end{proof}
		\end{clm*}
		Thus for all sufficiently large \(\alpha < \beta < \delta\), \(Q_\alpha = Q_\beta\), or in other words, \(P_\alpha \restriction A = P_\beta \restriction A\). It follows that \(P_\beta \leq P_\alpha\), and this contradicts our choice of \(P_\beta\). Thus our assumption that \((\mathbb Q,\leq)\) has no maximum element was false. 
		
		Let \(P\) be a maximum element of \((\mathbb Q,\leq)\). By the indecomposability of \(U\), we may assume \(|P| < \delta\) by replacing \(P\) with an equivalent element of \((\mathbb Q,\leq)\). We now apply \cref{Correspondence}. Let \(D\) be the ultrafilter corresponding to \(P\) as in \cref{Correspondence} (4): \[D = \{\mathcal A\subseteq P : \textstyle\bigcup \mathcal A\in U\}\]
		Let \(k : M_D\to M_U\) be unique elementary embedding with \(k\circ j_D = j_U\) and \(k(\id_D) = \Phi(P)\). 
		We have \(\lambda_D < \delta\) since \(|P| < \delta\). 
		
		Let \(\eta = (2^\delta)^+\). We will show \(\textsc{crt}(k) > j_U(\eta)\) if \(k\) is nontrivial, or in other words, that \(j_U(\eta)\subseteq k[M_D]\). 
		Since \(P\) is a maximum element of \(\mathbb Q^{\eta}_U\), \cref{Correspondence} (1), (2), and (3) imply that \(\Phi(P)\) is a maximum element of \(\mathbb P^{\eta}_U\). In other words, if \(x\in H^{M_U}(j_U[V]\cup \sup j_U[\eta])\), then \(x\) is definable in \(M_U\) from \(\Phi(P)\) and parameters in \(j_U[V]\), or in other words \(x\in k[M_D]\). In particular, \(\sup j_U[\eta]\subseteq k[M_D]\).
		
		We finish by showing that \(\sup j_U[\eta] = j_U(\eta)\). Suppose not. Then since \(\eta\) is regular, \(U\) is \(\eta\)-decomposable. Since \(\eta = (2^\delta)^+\),  \cref{PrikryThm} implies that \(U\) is \(\lambda\)-decomposable where \(\lambda = \text{cf}(2^\delta)\). But by K\"onig's Theorem, \(\lambda\in [\delta,2^\delta]\), and this is a contradiction.
\end{proof}

We can finally prove \cref{UpperBound}:

\begin{proof}[Proof of \cref{UpperBound}]
	Let \(j : V\to M\) be the ultrapower of the universe by \(\mathscr K_\lambda\). Assume \(2^{\delta} < \lambda\). We will show that \(\textsc{crt}(j) \geq \lambda\), so \(\mathscr K_\lambda\) is a \(\lambda\)-complete uniform ultrafilter on \(\lambda\), and hence \(\lambda\) is measurable. 
	
	Since \(\lambda\) is isolated and \(2^\delta < \lambda\), \(\mathscr K_\lambda\) is \(\gamma\)-indecomposable for all cardinals in the interval \([\delta,2^{\delta}]\). By \cref{DeltaLemma}, \(\delta\) is regular. Therefore we can apply \cref{Silver}. Fix \(D\) with \(\lambda_D < \delta\) such that there is an elementary embedding embedding \(k : M_D\to M_{\mathscr K_\lambda}\) with \(k\circ j_D = j\) and \(\textsc{crt}(k) > j_D(\delta)\) if \(k\) is nontrivial. 
	
	By \cref{IsolatedInternal}, \(D\I \mathscr K_\lambda\). Therefore \(j_D\restriction \delta\in M\). But \(j_D\restriction \delta = j\restriction \delta\) since \(\textsc{crt}(k)  > j_D(\delta)\). It follows that \(j\) is \(\delta\)-supercompact. 
	Since \(\delta\) is regular and \(\mathscr K_\lambda\) is \(\delta\)-indecomposable, \(j(\delta) = \sup j[\delta]\). Since \(j\) is \(\delta\)-supercompact and \(j(\delta) = \sup j[\delta]\), the Kunen Inconsistency (\cref{KunenInconsistency0}) implies that \(\textsc{crt}(j)\geq \delta\). There are no measurable cardinals in the interval \([\delta,\lambda)\) since in fact there are no Fr\'echet cardinals in \([\delta,\lambda)\). The fact that \(\textsc{crt}(j)\geq \delta\) therefore implies \(\textsc{crt}(j) \geq \lambda\), as desired.
\end{proof}

\cref{Silver} can be combined with \cref{UFCounting} to prove a strengthening of \cref{Cutpoint}:

\begin{thm}[UA]\label{SilverCor}
	Suppose \(\delta\) is a regular cardinal and \(U\) is a countably complete ultrafilter that is \(\lambda\)-indecomposable for all \(\lambda\in [\delta,2^\delta]\). Then there is an ultrafilter \(D\) with \(\lambda_D < \delta\) such that there is an internal ultrapower embedding \(h : M_D\to M_U\) with \(h\circ j_D=j_U\) and \(\textsc{crt}(h) > j_D(\delta)\) if \(h\) is nontrivial.
	\begin{proof}
		Using Silver's theorem, fix a uniform countably complete ultrafilter \(D\) on a cardinal \(\eta < \delta\) such that there is an elementary embedding \(k : M_D\to M_U\) with \(k\circ j_D = j_U\) and \(\textsc{crt}(k) > j_D((2^\delta)^+)\) if \(k\) is nontrivial.
		
		Recall that \(\mathscr B(X)\) denotes the set of countably complete ultrafilters on \(X\). \cref{UFCounting} implies that \(|\mathscr B(\eta)| \leq (2^{\eta})^+\). Thus \(j_D(\mathscr B(\eta))\) has cardinality less than or equal to \(j_D((2^{\eta})^+)\) in \(M_U\). Since \(\textsc{crt}(k) > j_D((2^{\eta})^+)\), \(k\) restricts to an isomorphism from \(j_D(\mathscr B(\eta),\sE)\) to \(j_U(\mathscr B(\eta),\sE)\). Moreover, for any \(Z\in j_D(\mathscr B(\eta))\), \[j_D^{-1}[Z] = j_U^{-1}[k(Z)]\]
		
		We now use the fact that \(k\) is an isomorphism conjugating \(j_D^{-1}\) to \(j_U^{-1}\) to conclude that \(k(\tr D D) = \tr U D\). By \cref{Reciprocity}, \(\tr D D\) is the least element of \(j_D(\mathscr B(\eta),\sE)\) with \(j_D^{-1}[Z] = D\). Therefore since \(k\) is an order-isomorphism that conjugates \(j_D^{-1}\) to \(j_U^{-1}\), \(k(\tr D D)\) is the least element \(Z\) of \(j_U(\mathscr B(\eta),\sE)\) with \(j_U^{-1}[Z] = D\). But by \cref{Reciprocity}, the least such \(Z\) is equal to \(\tr U D\). Thus \(k(\tr D D) = \tr U D\).
		
		Recall the characterization of the Rudin-Frol\'ik order in terms of translation functions (\cref{TransRF0}): if \(W\) and \(Z\) are countably complete ultrafilters, then \(W\D Z\) if and only if \(\tr Z W\) is principal in \(M_Z\). Applying this characterization in one direction to \(D\D D\), \(\tr D D\) is principal in \(M_D\). Therefore \(\tr U D = k(\tr D D)\) is principal in \(M_U\), so applying the characterization in the other direction, it follows that \(D\D U\).
		
		Let \(h: M_D\to M_U\) be the unique internal ultrapower embedding such that \(h \circ j_D = j_U\). By \cref{TransRF0}, \(\tr D D\) is the principal ultrafilter concentrated at \(\id_D\) and \(\tr U D\) is the principal ultrafilter concentrated at \(h(\id_D)\). Since \(k(\tr D D) = \tr U D\), it follows that \(k(\id_D) = h(\id_D)\). Since \(k\circ j_D = j_U\), in fact \(k\restriction j_D[V]\cup \{\id_D\} = h\restriction j_D[V]\cup \{\id_D\}\). Thus \(k = h\), since \(M_D = H^{M_D}(j_D[V]\cup \{\id_D\})\). It follows that \(h : M_D\to M_U\) is an internal ultrapower embedding with \(h\circ j_D = j_U\) and \(\textsc{crt}(h) > j_D(\delta)\) if \(h\) is nontrivial.
	\end{proof}
\end{thm}

Our work on isolated cardinals leads to some relatively simple criteria for the completeness of an ultrafilter in terms of a local version of irreducibility that will become important when we analyze larger supercompact cardinals:
\begin{defn}\label{LambdaIrredDef}\index{Ultrafilter!\(\lambda\)-irreducible}
	Suppose \(\lambda\) is a cardinal and \(U\) is a countably complete ultrafilter.
	\begin{itemize}
		\item \(U\) is {\it \(\lambda\)-irreducible} if for all \(D\D U\) with \(\lambda_D < \lambda\), \(D\) is principal.
		\item \(U\) is {\it \({\leq}\lambda\)-irreducible} if \(U\) is \(\lambda^+\)-irreducible.
	\end{itemize}
\end{defn}
Note that \(U\) is \({\leq}\lambda\)-irreducible if and only if \(U\) is \(\lambda^\sigma\)-irreducible.

At isolated cardinals, we have the following fact which is often useful:

\begin{thm}[UA]\label{IsoComplete}
	Suppose \(\lambda\) is a cardinal and \(U\) is a countably complete ultrafilter. 
	\begin{enumerate}[(1)]
		\item If \(\lambda\) is a strong limit cardinal that is not a limit of Fr\'echet cardinals, then \(U\) is \(\lambda\)-irreducible if and only if \(U\) is \(\lambda\)-complete.
		\item If \(\lambda\) is a strong limit cardinal and no cardinal \(\kappa < \lambda\) is \(\gamma\)-supercompact for all \(\gamma < \lambda\), then \(U\) is \(\lambda\)-irreducible if and only if \(U\) is \(\lambda\)-complete.
		\item If \(\lambda\) is isolated, then \(U\) is \(\lambda^+\)-irreducible if and only if \(U\) is \(\lambda^+\)-complete.
	\end{enumerate}
	\begin{proof}
		(1) is immediate from \cref{Cutpoint}. 
		
		(2) follows from (1). By \cref{LimitLma}, either \(\lambda\) is not a limit of Fr\'echet cardinals or \(\lambda\) is a limit of isolated cardinals. The former case is precisely (1). In the latter case, we can apply (1) at each isolated cardinal below \( \lambda\). Thus we conclude that \(U\) is \(\bar \lambda\)-complete for all isolated cardinals \(\bar \lambda < \lambda\). It follows that \(U\) is \(\lambda\)-complete as desired.
		
		(3) also follows from (1). Since \(U\) is \(\lambda^+\)-irreducible, \(U\) is \(\lambda^\sigma\)-irreducible, and by \cref{DoubleSigma}, \(\lambda^\sigma\) is measurable. Thus \(U\) is \(\lambda^\sigma\)-complete by (1) and in particular, \(U\) is \(\lambda^+\)-complete.
	\end{proof}
\end{thm}

Working in a bit more generality but with a stronger irreducibility assumption, we have the following completeness result:

\begin{thm}[UA]\label{RegCompleteThm}
	Suppose \(\delta\) is a regular cardinal such that no cardinal \(\kappa \leq \delta\) is \(\delta\)-supercompact. Then a countably complete ultrafilter \(U\) is \(\delta^+\)-complete if and only if it is \({\leq}2^\delta\)-irreducible.
	\begin{proof}
		The forward direction is trivial, so let us prove the converse. 
		
		Suppose that \(U\) is \({\leq}2^\delta\)-irreducible.  We claim that \(U\) is \(\lambda\)-irreducible where \(\lambda > \delta\) is a strong limit cardinal that is not a limit of Fr\'echet cardinals. An immediate consequence of the factorization theorem for isolated measurable cardinal (\cref{Cutpoint}) is that any \(\lambda\)-irreducible ultrafilter is \(\lambda\)-complete, and this proves the theorem.
		
		If \(\delta^\sigma\) does not exist, then the \({\leq}\delta\)-irreducibility of \(U\) implies that \(U\) itself is principal, so \(U\) is \(\lambda\)-irreducible and \(\lambda\)-complete for any cardinal \(\lambda\). Thus assume \(\delta^\sigma\) exists. 
		
		There are two cases. Suppose first that \(\delta\) is a Fr\'echet cardinal. Let \(\lambda = \delta^\sigma\). Since \(U\) is \({\leq}\delta\)-irreducible, \(U\) is \(\lambda\)-irreducible. We claim that \(\lambda\) is an isolated measurable cardinal. First note that \(\lambda > \delta^+\) since otherwise \(\kappa_{\delta^+}\) is \(\delta\)-supercompact by \cref{SuccessorThm}. Thus by \cref{SuccessorPrp}, \(\lambda\) is isolated. Assume towards a contradiction that \(\lambda\) is not measurable. Then by \cref{DeltaLemma}, \(\kappa_{\lambda}\) is \({<}\delta_\lambda\)-supercompact. But \(\delta < \delta_\lambda\) since \(\delta < \lambda\) is Fr\'echet, and hence \(\kappa_\lambda\) is \(\delta\)-supercompact, a contradiction. Hence \(\lambda\) is measurable.
		
		Suppose instead that \(\delta\) is not a Fr\'echet cardinal. If \(\delta^\sigma\) is measurable, let \(\lambda = \delta^\sigma\). Suppose \(\delta^\sigma\) is not measurable. By \cref{UpperBound}, \(\delta^\sigma \leq 2^\delta\), so in particular \(U\) is \({\leq}\delta^\sigma\)-irreducible. Let \(\lambda = \delta^{\sigma\sigma}\). (If \(\delta^{\sigma\sigma}\) does not exist, then again since \(U\) \({\leq}\delta^\sigma\)-irreducible,  \(U\) is principal.) By \cref{DoubleSigma}, \(\lambda\) is measurable; here, one must check that \(\delta^\sigma\) is isolated. Since \(U\) is \({\leq}2^\delta\)-irreducible, \(U\) is \({\leq}\delta^\sigma\)-irreducible, so \(U\) is \({<}\lambda\)-irreducible.
	\end{proof}
\end{thm}

One might expect a strengthening of this theorem to be true: if \(U\) is just \({\leq}\delta\)-irreducible and no \(\kappa \leq \delta\) is \(\delta\)-supercompact, then \(U\) should be \(\delta^+\)-complete. The main issue is that if \(\lambda = \delta^\sigma\) is an isolated nonmeasurable cardinal, then \(U = \mathscr K_\lambda\) is a counterexample. If instead \(\lambda\) is measurable, then \({\leq}\delta\)-irreducibility indeed suffices. What about \(\delta\)-irreducibility? If \(\delta\) is the least cardinal such that \(\mathscr K_\delta\) exists and does not have a \(\delta\)-supercompact ultrapower, then \(U = \mathscr K_\delta\) is a counterexample.

A similar theorem is true for singular cardinals:

\begin{thm}[UA]
	Suppose \(U\) is a countably complete ultrafilter and \(\gamma\) is a singular cardinal such that no \(\kappa \leq \gamma\) is \(\gamma^+\)-supercompact. Then \(U\) is \(\gamma^+\)-complete if and only if \(U\) is \({\leq}2^\delta\)-irreducible for all \(\delta < \gamma\).
	\begin{proof}
		Let \(\delta = \sup \{\gamma^+ : \gamma < \lambda \text{ is a Fr\'echet cardinal}\}\). 
		
		Suppose first that \(\delta\) is regular. Since \(\delta\) is not Fr\'echet, no cardinal \(\kappa \leq \delta\) is \(\delta\)-supercompact. Since \(U\) is \({\leq}2^\delta\)-irreducible, we are in a position to apply \cref{RegCompleteThm}. We can conclude that \(U\) is \(\delta^+\)-complete. Since there are no measurable cardinals in the interval \((\delta,\gamma)\), it follows that \(U\) is \(\gamma^+\)-complete.
		
		Suppose instead that \(\delta\) is singular. If \(\delta^\sigma\) does not exist, then it is easy to see that \(U\) is principal, and thus we are done. Therefore assume \(\delta^\sigma\) exists, and let \(\lambda = \delta^\sigma\). Note that \(\lambda > \delta^+\): if \(\delta < \gamma\) this follows from the fact that \(\delta^+\) is not Fr\'echet, while if \(\delta = \gamma\), this follows from the fact that no cardinal is \(\gamma^+\)-supercompact. Thus \(\lambda\) is isolated. Note that \(\delta_\lambda = \delta\) is singular. Therefore by \cref{DeltaLemma}, \(\lambda\) is measurable. Since \(U\) is \({\leq}\delta\)-irreducible,  \(U\) is \({<}\lambda\)-irreducible, and therefore as an immediate consequence of the factorization theorem for isolated measurable cardinal (\cref{Cutpoint}), \(U\) is \(\lambda\)-complete.
	\end{proof}
\end{thm}

Let us close this subsection with a remark about the size of regular isolated cardinals.
\begin{defn}
	A regular cardinal \(\kappa\) is {\it \(\sigma\)-Mahlo} if there is a countably complete  weakly normal ultrafilter on \(\kappa\) that concentrates on regular cardinals.
\end{defn}

\begin{prp}
	If \(\kappa\) is \(\sigma\)-Mahlo then \(\kappa\) is weakly Mahlo.\qed
\end{prp}

In fact, \(\sigma\)-Mahlo cardinals are ``greatly weakly Mahlo." A theorem of Gitik shows that it is consistent that there is a \(\sigma\)-Mahlo cardinal that does not have the tree property.

\begin{thm}[UA]\label{SigmaMahlo}
	Suppose \(\kappa\) is a regular isolated cardinal. Then \(\kappa\) is \(\sigma\)-Mahlo. In fact, \(\mathscr K_\kappa\) concentrates on regular cardinals.
	\begin{proof}
		Let \(j : V\to M\) be the ultrapower of the universe by \(\mathscr K_\kappa\). Let \(\kappa_* = \sup j[\kappa] \). Let \(\delta = \text{cf}^{M}(\kappa_*)\). By \cref{KetonenTight}, \(j\) is \((\kappa, \delta)\)-tight, so \(j\) is discontinuous at any regular cardinal \(\iota\leq \kappa\) such that \(\delta < j(\iota)\). Since \(\kappa\) is isolated, \(j\) is continuous at all sufficiently large cardinals less than \(\kappa\). Putting these observations together, it follows that there are no regular cardinals \(\iota < \kappa\) such that \(j(\iota) > \delta\). In other words, \(\sup j[\kappa]\leq\delta\). Thus \(\kappa_* = \delta\), so \(\kappa_*\) is regular. Since \(\mathscr K_\kappa\) is weakly normal, \(\kappa_* =  \id_{\mathscr K_\kappa}\), so by \L o\'s's Theorem, \(\mathscr K_\kappa\) concentrates on regular cardinals.
	\end{proof} 
\end{thm}

This fact has a converse: assuming UA, any \(\sigma\)-Mahlo cardinal that is not measurable is isolated. It is not clear that singular Fr\'echet cardinals must be very large. For example, we do not know how to rule out that the least Fr\'echet cardinal \(\lambda\) that is neither measurable nor a limit of measurables is in fact equal to \(\kappa^{+\kappa}\) for some measurable \(\kappa < \lambda\).
\subsection{The linearity of the Mitchell order without GCH}
\cref{GCHLinear} states that assuming UA + GCH, the Mitchell order is linear on normal fine ultrafilters on \(P_\text{bd}(\lambda)\), the collection of bounded subsets of \(\lambda\). Here we prove essentially the same theorem using UA alone. Instead of using \(P_\text{bd}(\lambda)\) as our underlying set, we use the following variant:
\begin{defn}
	For any cardinal \(\lambda\), let \(P_*(\lambda)= \{\sigma\in P_\text{bd}(\lambda):|\sigma|^+ < \lambda\}\).\index{\(P_*(\lambda)\)}
\end{defn}
The following obvious characterization of \(P_*(\lambda)\) is often more useful than the definition above: 
\[P_*(\lambda) = 
\begin{cases}
P_\text{bd}(\lambda) &\text{if \(\lambda\) is a limit cardinal}\\ P_{\gamma}(\lambda) &\text{if \(\lambda\) is a successor cardinal and \(\gamma\) is its cardinal predecessor}\end{cases}\]

\begin{defn}
	For any cardinal \(\lambda\), let \(\mathscr U_\lambda\) denote the set of normal fine ultrafilters on \(P_*(\lambda)\). Let \(\mathscr U = \bigcup_{\lambda \in \text{Card}} \mathscr U_\lambda\).
\end{defn}

The main theorem of this subsection is the following:
\begin{thm}[UA]\label{ULinearity}\index{Generalized Mitchell order!linearity!on normal fine ultrafilters}
	The class \(\mathscr U\) is linearly ordered by the Mitchell order.
\end{thm}

Due to the following fact, \cref{ULinearity} can be seen as a precise formulation of the (literally false) statement that the Mitchell order is linear on normal fine ultrafilters:

\begin{prp}\label{UIso}
	Every normal fine ultrafilter is isomorphic to a unique element of \(\mathscr U\).
	\begin{proof}
		Recall that for any cardinal \(\lambda\), \(\mathscr N_\lambda\) denotes the set of normal fine ultrafilters on \(P_\text{bd}(\lambda)\) and \(\mathscr N = \bigcup_{\lambda\in \text{Card}} \mathscr N_\lambda\). Also recall \cref{NIso}, which states that every normal fine ultrafilter is isomorphic to a unique element of \(\mathscr N\). Therefore to prove the proposition, it suffices to show that there is a bijection \(\phi : \mathscr N\to \mathscr U\) such that \(\phi(\mathcal U)\cong \mathcal U\) for all \(\mathcal U\in \mathscr N\). 
		
		In fact, if \(\mathcal U\in \mathscr N_\lambda\), we will just set \(\phi(\mathcal U) = \mathcal U\restriction P_*(\lambda)\). 
		It is clear that \(\phi\) is as desired as long as \(P_*(\lambda)\in \mathcal U\). We now establish that this holds.
		Let \(j : V\to M\) be the ultrapower of the universe by \(\mathcal U\). Then \(\id_\mathcal U = j[\lambda]\) by \cref{NormalFineChar}. 
		Of course \(|j[\lambda]|^M = \lambda\), but note also that \(\lambda^{+M} < j(\lambda)\): by \cref{Kunen}, there is an inaccessible cardinal \(\kappa \leq\lambda\) such that \(\lambda < j(\kappa)\), so \(\lambda^{+M} < j(\kappa) \leq j(\lambda)\). 
		Thus \(|j[\lambda]|^{+M} < j(\lambda)\). 
		By \L o\'s's Theorem, it follows that \(\{\sigma\in P_\text{bd}(\lambda) : |\sigma|^+ < \lambda\}\in \mathcal U\).
		That is, \(P_*(\lambda)\in \mathcal U\).
	\end{proof}
\end{prp}

The reason we use \(P_*(\lambda)\) as an underlying set rather than sticking with \(P_\text{bd}(\lambda)\) is that without assuming GCH, we cannot prove \(|P_\text{bd}(\lambda)|= \lambda\). Therefore \(P_\text{bd}(\lambda)\) may be too large to use as an underlying set. On the other hand, we can prove \(|P_*(\lambda)| = \lambda\) in the relevant cases:
\begin{prp}[UA]\label{P*Lemma}
	Suppose \(\lambda\) is a cardinal such that \(\mathscr U_\lambda\) is nonempty. Then \(|P_*(\lambda)| = \lambda\).
	\begin{proof}
		Since \(\mathscr U_\lambda\) is nonempty, there is a normal fine ultrafilter on \(P_*(\lambda)\), and hence there is a cardinal \(\kappa \leq \lambda\) that is \(\lambda\)-supercompact.
		
		There are now two cases. 
		
		Suppose first that \(\lambda\) is a limit cardinal. Then \(P_*(\lambda) = P_\text{bd}(\lambda)\). Moreover by \cref{MainTheorem}, \(2^{<\lambda} = \lambda\). Thus \(|P_*(\lambda)| = |P_\text{bd}(\lambda)| = 2^{<\lambda} = \lambda\).
		
		Suppose instead that \(\lambda\) is a successor cardinal. Let \(\gamma\) be the cardinal predecessor of \(\lambda\). Then \(P_*(\lambda) = P_{\gamma}(\lambda)\), so \(|P_*(\lambda)|  = \lambda^{<\gamma}\). Since \(\lambda\) is regular, \(\lambda^{<\gamma} = \lambda \cdot \gamma^{<\gamma}\). To finish, it therefore suffices to show \(\gamma^{<\gamma} \leq \lambda\). By \cref{SuccGCH}, \(2^{<\gamma} = \gamma\). If \(\gamma\) is singular, then \(\gamma\) is a singular strong limit cardinal, so by Solovay's Theorem on SCH above a strongly compact cardinal (\cref{SolovaySCH}), \(\gamma^{<\gamma} \leq \gamma^\gamma = \gamma^+ = \lambda\). Otherwise, \(\gamma^{<\gamma} = 2^{<\gamma} = \gamma\).  
	\end{proof}
\end{prp}

Recall that an ultrafilter \(U\) on a set \(X\) is {\it hereditarily uniform} if \(|\text{tc}(X)| = \lambda_U\). We observed that the generalized Mitchell order is well-behaved on hereditarily uniform ultrafilters: for example it is isomorphism invariant (\cref{HeredLemma}) and transitive (\cref{StrongTransitivity}).
Under UA, it follows that the Mitchell order is well-behaved on \(\mathscr U\):

\begin{lma}[UA]\label{UHered}
	Every ultrafilter in \(\mathscr U\) is hereditarily uniform.
	\begin{proof}
		Suppose \(\mathcal U\in \mathscr U\). Fix a cardinal \(\lambda\) with \(\mathcal U\in \mathscr U_\lambda\). Since \(P_*(\lambda)\) is the underlying set of \(\mathcal U\), to show that \(\mathcal U\) is hereditarily transitive, we must show that \(|\text{tc}(P_*(\lambda))| \leq \lambda_\mathcal U\). Of course, \(\text{tc}(P_*(\lambda)) = P_*(\lambda)\cup \lambda\), which has cardinality \(\lambda\) by \cref{P*Lemma}. Since \(j_\mathcal U\) is \(\lambda\)-supercompact, \cref{UFSuperBound} implies that \(\lambda \leq \lambda_\mathcal U\).  Thus \(|\text{tc}(P_*(\lambda))| \leq \lambda_\mathcal U\), as desired.
	\end{proof}
\end{lma}

Recall that an ultrafilter \(U\) on a cardinal \(\lambda\) is {\it isonormal} if \(U\) is weakly normal and \(j_U\) is \(\lambda\)-supercompact. Recall \cref{IsoNormalThm}, which states that every normal fine ultrafilter is isomorphic to an isonormal ultrafilter. Combined with the isomorphism invariance of the Mitchell order on hereditarily uniform ultrafilters, the following theorem therefore easily implies \cref{ULinearity}:

\begin{thm}[UA]\label{IsoIPoint}
	Suppose \(U\) is an isonormal ultrafilter. Then for any \(D\sE U\), \(D\mo U\). In particular, the Mitchell order is linear on isonormal ultrafilters.
\end{thm}

Note that a strong version of this theorem (\cref{LipMO}) follows from GCH.
Let us explain in full detail how to prove the linearity of the Mitchell order on \(\mathscr U\) (\cref{ULinearity}) from \cref{IsoIPoint}:
\begin{proof}[Proof of \cref{ULinearity}]
	Suppose \(\mathcal U_0\) and \(\mathcal U_1\) are elements of \(\mathscr U\). We must show that either \(\mathcal U_0 \mo \mathcal U_1\), \(\mathcal U_0 = \mathcal U_1\), or \(\mathcal U_0\gmo \mathcal U_1\). Since every normal fine ultrafilter is isomorphic to an isonormal ultrafilter (\cref{IsoNormalThm}), there are isonormal ultrafilters \(U_0\) and \(U_1\) such that \(U_0 \cong \mathcal U_0\) and \(U_1\cong \mathcal U_1\). By \cref{IsoIPoint}, either \(U_0 = U_1\), \(U_0 \mo U_1\), or \(U_0\gmo U_1\). If \(U_0 = U_1\), then \(\mathcal U_0 \cong \mathcal U_1\). Therefore by the uniqueness clause of \cref{UIso}, \(\mathcal U_0 = \mathcal U_1\). If \(U_0\mo U_1\), then since the Mitchell order is isomorphism invariant on hereditarily uniform ultrafilters (\cref{HeredLemma}), \(\mathcal U_0\mo \mathcal U_1\). (All the ultrafilters we are considering are hereditarily uniform; the nontrivial part of this is \cref{UHered}.) Similarly, if \(U_0\gmo U_1\), then \(\mathcal U_0\gmo \mathcal U_1\).
\end{proof}

We therefore proceed to the proof of \cref{IsoIPoint}. This requires a general fact from the theory of the internal relation which is of independent interest. Here is the idea. Since no nonprincipal ultrafilter \(U\) satisfies \(U\mo U\), under UA there is a least \(W\) in the Ketonen order such that \(W\not \mo U\). What is the relationship between \(U\) and \(W\)? Perhaps \(W\D U\), but this is an open question. It turns out that one can make some headway if one considers instead the \(\sE\)-least \(W\) such that \(W\not \I U\). (\cref{MitchellPoint} shows that this actually defines the same ultrafilter.)
\begin{thm}[UA]\label{IPoint}
	Suppose \(U\) is a nonprincipal countably complete ultrafilter and \(W\) is the \(\sE\)-least countably complete uniform ultrafilter on an ordinal such that \(W\not \I U\). Then for any \(D\I U\), \(D\I W\).
\end{thm}

To prove \cref{IPoint}, we use the following closure property of the internal relation:

\begin{lma}\label{InternalClosed}
	Suppose \(D\I U\) is an ultrafilter on a set \(X\) and \(\langle W_i : i\in X\rangle\) is a sequence of ultrafilters on a set \(Y\) such that \(W_i\I U\) for all \(i\in X\). Then \(D\text{-}\sum_{i\in X}W_i\I U\) and \(D\text{-}\lim_{i\in X}W_i\I U\).
	\begin{proof}
		Since \(D\text{-}\lim_{i\in X}W_i\RK D\text{-}\sum_{i\in X}W_i\), if we show show that \(D\text{-}\sum_{i\in X}W_i\I U\), we obtain \(D\text{-}\lim_{i\in X}W_i\I U\) as a consequence of \cref{RKInternal}.
		
		Let \(j : V\to N\) be the ultrapower of the universe by \(D\). Let \(W = [\langle W_i : i\in X\rangle]_D\) and let \(k : N\to P\) be the ultrapower of \(M\) by \(W\). Thus \(k \circ j\) is the ultrapower embedding associated to \(D\text{-}\sum_{i\in X}W_i\), so to prove the lemma, we must show that \(k\circ j\restriction M_U\) is an internal ultrapower embedding of \(M_U\). 
		
		Since \(D\I U\), \(j\) is an internal ultrapower embedding of \(M_U\). Therefore to show \(k\circ j\restriction M_U\) is an internal ultrapower embedding of \(M_U\), it suffices to show that \(k\restriction j(M_U)\) is an internal ultrapower embedding of \(j(M_U)\). Note that by the elementarity of \(j : V\to N\), \(j(M_U) = (M_{j(U)})^N\). Since \(k = (j_W)^N\), to show that \(k \restriction j(M_U)\) is an internal ultrapower embedding of \(j(M_U)\), it suffices to show that \(W\I j(U)\) in \(N\). But \(W_i \I U\) for all \(i\in X\), so \(W\I j(U)\) in \(N\) by \L o\'s's Theorem.
	\end{proof}
\end{lma}

\begin{proof}[Proof of \cref{IPoint}]
		Suppose \(D\I U\). By \cref{IChar}, \(\tr D U = j_D(U)\). We claim \(j_D(W)\E^{M_D}\tr D W\). Suppose towards a contradiction that this fails, so \(\tr D W\sE^{M_D} j_D(W)\). Let \(X\) be the underlying set of \(D\), and fix \(\langle W_i : i\in I\rangle\) such that \(\tr D W = [\langle W_i : i\in X\rangle]_D\). Then since \(W_i \sE W\) for all \(i\in X\), in fact \(W_i\I U\) for all \(i\in X\). Thus \(D\text{-}\lim_{i\in X}W_i\I U\) by \cref{InternalClosed}. But \(W = D\text{-}\lim_{i\in I}W_i\), and this contradicts the definition of \(W\).
\end{proof}

\begin{proof}[Proof of \cref{IsoIPoint}]
	Let \(\lambda = \lambda_U\). If \(2^{<\lambda} = \lambda\), then \(U\) is Dodd sound (\cref{IsoSound}), so for all \(D\sE U\), we have \(D \mo U\) (\cref{LipMO}), and thus we are done. We therefore assume that \(2^{<\lambda} > \lambda\). (It is not clear whether this assumption is consistent. We will not try to reach a contradiction, however, but rather to prove that the theorem is true even if \(2^{<\lambda} > \lambda\).)
	
	Since \(j_U\) witnesses that some cardinal \(\kappa \leq \lambda\) is \(\lambda\)-supercompact, the local version of the theorem that GCH holds above a supercompact under UA (\cref{MainTheorem}) implies that \(2^{<\lambda} = \lambda\) if \(\lambda\) is a limit cardinal. Therefore by our assumption that \(2^{<\lambda} > \lambda\), \(\lambda\) is a successor cardinal. 
	
	Let \(\gamma\) be the cardinal predecessor of \(\lambda\). To simplify notation, we will from now on refer to \(\lambda\) only as \(\gamma^+\). We therefore reformulate our  assumption that \(2^{<\lambda} > \lambda\):
	\begin{equation}\label{BadAssumption}2^\gamma > \gamma^+\end{equation} Since \(\gamma^+\) is Fr\'echet, our local result on GCH (\cref{SuccGCH}) yields that \(2^{<\gamma} = \gamma\).  If \(\gamma\) is singular, then since \(2^{<\gamma} = \gamma\), \(\gamma\) is a singular strong limit cardinal, so the fact that \(2^\gamma > \gamma^+\) contradicts the local version of Solovay's Theorem that SCH holds above a strongly compact cardinal (\cref{SolovaySCH}). Therefore \(\gamma\) is regular.
	
	\begin{clm}\label{MCommandClm}
		\(M_U\) satisfies that \(2^{2^\gamma} = (2^\gamma)^+\)
	\begin{proof}
		Let \(\mathcal D\) be the normal fine ultrafilter on
		\(P_\text{bd}(\gamma)\) derived from \(j_U\) using \(j_U[\gamma]\).
		Since \(M_U\) is closed under \(\gamma^+\)-sequences, every ultrafilter
		on \(\gamma\) belongs to \(M_U\) (\cref{MitchellLemma}). Therefore since
		\(P_\text{bd}(\gamma)\) has hereditary cardinality \(2^{<\gamma} =
		\gamma\), we have \(\mathcal D\in M_U\). Therefore by a generalization
		of Solovay's argument that a \(2^\kappa\)-supercompact cardinal carries
		\(2^{2^\kappa}\) normal ultrafilters (\cref{MuGen}), \(M_U\) satisfies
		that every subset of \(P(\gamma)\) belongs to \(M_\mathcal W\) for some
		normal fine ultrafilter \(\mathcal W\) on \(P_\text{bd}(\gamma)\). By
		\cref{GeneralizedCounting} applied in \(M_U\), \(M_U\) satisfies that
		\(2^{2^\gamma} = (2^\gamma)^+\). (Alternately one can use
		\cref{UFCounting}.)
	\end{proof}
	\end{clm}
	Now let \(\eta = ((\gamma^+)^\sigma)^{M_U}\) be the least Fr\'echet cardinal
	above \(\gamma^+\) as computed in \(M_U\). The following claim is a
	consequence of our analysis of isolated cardinals:
	\begin{clm}\label{EtaMeasurableClm}
	\(\eta\) is a measurable cardinal of \(M_U\).
	\begin{proof}
		Since \(P(\gamma)\subseteq M_U\), \((2^{\gamma})^{M_U} \geq
		(2^{\gamma})^V > \gamma^+\), and therefore \(M_U\) satisfies that
		\(2^\gamma > \gamma^+\). 
		
		We now work in \(M_U\) to avoid a profusion of superscripts. We cannot
		have \(\eta = \gamma^{++}\): otherwise \(\gamma^{++}\) is Fr\'echet and
		hence \(2^{\gamma} = \gamma^+\) by \cref{SuccGCH}, contradicting that
		\(2^\gamma > \gamma^+\). Therefore \(\eta > \gamma^+\) and so by
		\cref{SuccessorPrp}, \(\eta\) is isolated. 
		
		Let \(\delta = \delta_\eta\). Then since \(\eta = (\gamma^+)^\sigma\),
		\(\delta \leq \gamma^{++}\leq 2^\gamma\). The final inequality uses the
		fact that \(2^\gamma > \gamma^+\). Thus \(2^{\delta} \leq 2^{2^\gamma} =
		(2^\gamma)^+\) by \cref{MCommandClm}. But \(2^\gamma < \eta\) by our
		results on the continuum function below an isolated cardinal
		(\cref{LowerBound}). Therefore \((2^\gamma)^+ < \eta\) since \(\eta\) is
		isolated (and therefore is a limit cardinal). It follows that \(2^\delta
		< \eta\). Therefore \(\eta\) is measurable by \cref{UpperBound}.
	\end{proof}
	\end{clm}

	Let \(W\) be the \(\sE\)-least countably complete ultrafilter on an ordinal
	such that \(W\not \I U\). Then \(W\E U\). To prove the theorem, we must show
	\(U = W\).
	
	Since every countably complete ultrafilter on \(\gamma\) belongs to \(M_U\)
	and hence is internal to \(U\), we have \(\lambda_W = \gamma^+\). Let
	\[(k,h) : (M_W,M_U)\to N\] be the pushout of \((j_W,j_U)\). 
	\begin{clm} If \(h\) is nontrivial, then \(\textsc{crt}(h)\geq \eta\).\end{clm}
	\begin{proof}
	Let \(W' = \tr U W\), so \(h : M_U\to N\) is the ultrapower of \(M_U\) by
	\(W'\). 
	
	Suppose \(D\) is a countably complete ultrafilter of \(M_U\) with
	\(\lambda_D < \eta\). We claim that \(D\I W'\) in \(M_U\). Since
	\(\lambda_D\) is a Fr\'echet cardinal of \(M_U\) below \(\eta =
	((\gamma^+)^{\sigma})^{M_U}\), \(\lambda_D \leq \gamma^+\). We may therefore
	assume that the underlying set of \(D\) is \(\gamma^+\). Since \(j_U\) is
	\(\gamma^+\)-supercompact, \(P(\gamma^+)\subseteq M_U\). Thus \(D\) is an
	ultrafilter on \(\gamma^+\) (in \(V\)). Since \(M_U\) is closed under
	\(\gamma^+\)-sequences, \(j_D\restriction M_U = j_D^{M_U}\), so in fact
	\(D\I U\). By \cref{IPoint}, \(D\I W\). Thus \(j_D\restriction N\) is
	amenable to both \(M_U\) and \(M_W\). By our characterization of the
	internal ultrapower embeddings of a pushout (\cref{PushoutInternal}),
	\(j_D\restriction N\) is an internal ultrapower embedding of \(N\).
	Equivalently \(j_D^{M_U}\restriction N\) is an internal ultrapower embedding
	of \(N\), or in other words, \(D\I W'\) in \(M_U\).
	
	By \cref{IsoThresh}, \(h[\eta]\subseteq \eta\). Working in \(M_U\), \(\eta\)
	is a strong limit cardinal, \(h[\eta]\subseteq \eta\), and for all \(D\)
	with \(\lambda_D < \eta\), \(D\I W'\). Applying in \(M_U\) our criterion for
	the completeness of an ultrafilter in terms of the internal relation
	(\cref{InternalComplete}), it follows that \(W'\) is \(\eta\)-complete.
	Since \(h\) is the ultrapower of \(M_U\) by \(W'\), if \(h\) is nontrivial
	then \(\textsc{crt}(h) \geq \eta\).
	\end{proof}
	Since \(U\) is a weakly normal ultrafilter on \(\gamma^+\), \(\id_U = \sup
	j_U[\gamma^+]\) (\cref{RegWeaklyNormal}). Since \(h\) is the identity or
	\(\textsc{crt}(h) \geq \eta > \gamma^+\), \(h\) is continuous at ordinals of
	\(M_U\)-cofinality \(\gamma^+\). Since \(M_U\) is closed under
	\(\gamma^+\)-sequences, \(\sup j_U[\gamma^+]\) is on ordinal of
	\(M_U\)-cofinality \(\gamma^+\). Therefore 
	\[h(\id_U) = h(\sup j_U[\gamma^+]) = \sup h\circ j_U[\gamma^+] \leq \sup
	k\circ j_W[\gamma^+] \leq k(\id_W)\] The final inequality follows from the
	fact that \(\lambda_W = \gamma^+\) and hence \( \sup j_W[\gamma^+]\leq
	\id_W\). Therefore \((k,h)\) witnesses that \(U\E W\). Since \(U\E W\) and
	\(W\E U\), \(U = W\), as desired.
\end{proof}

\chapter{Higher Supercompactness}\label{SCChapter2}
\section{Introduction}
\subsection{Obstructions to the supercompactness analysis}
The main result of \cref{SCChapter1} is that under UA, the first strongly
compact is supercompact. What about the second? What about all of the other
strongly compact cardinals? This chapter answers all these questions and more.
In this introductory section, we explain in broad strokes the obstructions to
generalizing the theory of \cref{SCChapter1} and the technique that ultimately
overcomes them.
\subsection{Menas's Theorem}\label{MenasSection}
The first obstruction to generalizing the results of \cref{SCChapter1} is that
not every strongly compact cardinal is supercompact. This is a consequence of
the following theorem of Menas:

\begin{thm}[Menas]\label{Menas}\index{Menas's Theorem}
	The least strongly compact limit of strongly compact cardinals is not
	supercompact.
\end{thm}

In order to explain the proof, we introduce an auxiliary notion:

\begin{defn}
	Suppose \(\kappa\) and \(\lambda\) are cardinals. A cardinal \(\kappa\) is
	{\it almost \(\lambda\)-strongly compact} if for any \(\alpha < \kappa\),
	there is an elementary embedding \(j : V\to M\) such that \(\textsc{crt}(j)
	> \alpha\) and \(M\) has the \(({\leq}\lambda,{<}j(\kappa))\)-covering
	property; \(\kappa\) is {\it almost strongly compact} if \(\kappa\) is
	almost \(\lambda\)-strongly compact for all cardinals \(\lambda\).
\end{defn}

As in \cref{StrongCChar}, there is a characterization of almost strong
compactness in terms of fine ultrafilters:

\begin{lma}
	A cardinal \(\kappa\) is almost \(\lambda\)-strongly compact if and only if
	for every \(\alpha < \kappa\), there is an \(\alpha^+\)-complete fine
	ultrafilter on \(P_\kappa(\lambda)\).\qed
\end{lma}

Unlike strongly compact cardinals, it is easy to see that almost strongly
compact cardinals form a closed class:

\begin{lma} \label{AlmostSCClosed}
	Any limit of almost \(\lambda\)-strongly compact cardinals is almost
	strongly compact. In particular, every limit of strongly compact cardinals
	is almost strongly compact.\qed 
\end{lma}
The following proposition shows that almost strongly compact cardinals really
almost are strongly compact:
\begin{prp}\label{MeasurableAlmostSC}
	A cardinal \(\kappa\) is \(\lambda\)-strongly compact if and only if
	\(\kappa\) is measurable and almost \(\lambda\)-strongly compact.
	\begin{proof}
		Since \(\kappa\) is measurable, there is a  \(\kappa\)-complete uniform
		ultrafilter \(U\) on \(\kappa\). Since \(\kappa\) is almost strongly
		compact, for each \(\alpha < \kappa\), there is an \(\alpha^+\)-complete
		fine ultrafilter \(\mathcal W_\alpha\) on \(P_\kappa(\lambda)\). Let
		\[\mathcal W = U\text{-}\lim_{\gamma < \kappa} \mathcal W_\gamma\] It is
		immediate that \(\mathcal W\) is a fine ultrafilter on
		\(P_\kappa(\lambda)\). 
		
		We claim that \(\mathcal W\) is \(\kappa\)-complete. Suppose \(\nu <
		\kappa\) and \(\{A_i : i < \nu\}\subseteq \mathcal W\). For each \(i <
		\nu\), let \(S_i  = \{\alpha < \kappa : A_i\in \mathcal W_\alpha\}\).
		Since \(A_i\in \mathcal W\), \(S_i\in U\) by the definition of an
		ultrafilter limit. Since \(U\) is \(\kappa\)-complete, \(\bigcap_{i <
		\nu} S_i\) belongs to \(U\). Since \(U\) is uniform, \(\bigcap_{i < \nu}
		S_i\setminus \nu\in U\). Suppose \(\alpha\in \bigcap_{i < \nu}
		S_i\setminus \nu\). Then \(\{A_i : i  < \nu\}\in \mathcal W_\alpha\).
		Therefore since \(\mathcal W_\alpha\) is \(\alpha^+\)-complete,
		\(\bigcap_{i< \nu} A_i\in \mathcal W_\alpha\). Thus \[\bigcap_{i < \nu}
		S_i\setminus \nu\subseteq \left\{\alpha < \kappa : \bigcap_{i < \nu}
		A_i\in \mathcal W_\alpha\right\}\] It follows that \(\{\alpha < \kappa :
		\bigcap_{i < \nu} A_i\in \mathcal W_\alpha\}\in U\). In other words,
		\(\bigcap_{i < \nu} A_i\in \mathcal W\).
	\end{proof}
\end{prp}

\begin{cor}[Menas]\label{MenasBig}
	Every measurable limit of strongly compact cardinals is strongly
	compact.\qed
\end{cor}

The least strongly compact limit of strongly compact cardinals is therefore in a
sense accessible from below:

\begin{lma}[Menas]\label{MenasSmall}
	Let \(\kappa\) be the least strongly compact limit of strongly compact
	cardinals. Then the set of measurable cardinals below \(\kappa\) is
	nonstationary in \(\kappa\). Therefore \(\kappa\) has Mitchell rank \(1\).
	In particular, \(\kappa\) is not \(\mu\)-measurable, let alone
	\(2^\kappa\)-strong, let alone \(2^\kappa\)-supercompact.
	\begin{proof}
		Let \(C\) be the set of limits of strongly compact cardinals less than
		\(\kappa\). Since \(\kappa\) is a regular limit of strongly compact
		cardinals, \(C\) is unbounded in \(\kappa\). Moreover, \(C\) is closed
		by definition. We claim that \(C\) contains no measurable cardinals.
		Suppose \(\delta\in C\) is measurable. Then by \cref{MenasBig},
		\(\delta\) is strongly compact. This contradicts that \(\kappa\) is the
		least strongly compact limit of strongly compact cardinals. It follows
		that the class of measurable cardinals is nonstationary in \(\kappa\).
	\end{proof}
\end{lma}

A strongly compact cardinal \(\kappa\) always carries \(2^{2^\kappa}\)-many
\(\kappa\)-complete ultrafilters. But Menas's theorem shows that the Mitchell
order may be trivial on \(\kappa\). Under UA, this has the following surprising
consequence:
\begin{thm}[UA] The least strongly compact limit of strongly compact cardinals
	carries a unique normal ultrafilter.
	\begin{proof}
		Let \(\kappa\) be the least strongly compact limit of strongly compact
		cardinals. By Menas's Theorem (\cref{MenasSmall}), the rank of the
		Mitchell order on normal ultrafilters on \(\kappa\) is 1. By
		\cref{UAMO}, the Mitchell order linearly orders these ultrafilters.
		Therefore \(\kappa\) carries exactly one normal ultrafilter.
	\end{proof}
\end{thm}

\subsection{Complete UA}
The second obstruction to generalizing the results of \cref{SCChapter1} is much
more subtle: UA alone does not seem to suffice to enact a direct generalization
of the structure of the least supercompact cardinal to the higher ones. In order
to shed light on the underlying issue, we introduce a principle called the
Complete Ultrapower Axiom (CUA), which does suffice.

\begin{defn}
	Suppose \(\kappa\) is an uncountable cardinal. Then
	\(\textnormal{UA}(\kappa)\) denotes the following statement. Suppose \(j_0 :
	V\to M_0\) and \(j_1: V\to M_1\) are ultrapower embeddings with
	\(\textsc{crt}(j_0)\geq \kappa\) and \(\textsc{crt}(j_1)\geq \kappa\). Then
	there is a comparison \((i_0,i_1) : (M_0,M_1)\to N\) of \((j_0,j_1)\) such
	that \(\textsc{crt}(i_0) \geq \kappa\) and \(\textsc{crt}(i_1)\geq \kappa\).
\end{defn}

Thus the usual Ultrapower Axiom is equivalent to \(\text{UA}(\omega_1)\). Notice
that \(\textnormal{UA}(\kappa)\) is equivalent to the assertion that the
Rudin-Frol\'ik order is directed on \(\kappa\)-complete ultrafilters.
\begin{cua}
	\(\text{UA}(\kappa)\) holds for all uncountable cardinals
	\(\kappa\).\index{Ultrapower Axiom!Complete Ultrapower Axiom}
\end{cua}

Assuming CUA, one can generalize all the proofs in the previous section to
obtain results about the higher supercompact cardinals. In fact, one does not
even need to dig into the details to see that this is possible:

\begin{prp}[CUA] Suppose \(\kappa\) is strongly compact. Either \(\kappa\) is
	supercompact or \(\kappa\) is a limit of supercompact cardinals.
	\begin{proof}[Sketch] Suppose first that \(\kappa\) is not a limit of
		strongly compact cardinals. We will show that \(\kappa\) is
		supercompact. Let \(\delta < \kappa\) be the supremum of the strongly
		compact cardinals below \(\kappa\). Let \(G\subseteq
		\text{Col}(\omega,\delta)\) be \(V\)-generic. Then in \(V[G]\),
		\(\kappa\) is the least strongly compact cardinal. Moreover, since
		\(\textnormal{UA}(\delta^+)\) holds in \(V\), UA holds in \(V[G]\).
		Therefore by the analysis of the least strongly compact cardinal under
		UA (\cref{StrongSuper}), \(\kappa\) is supercompact in \(V[G]\). It
		follows that \(\kappa\) is supercompact in \(V\), as desired.
		
		Suppose instead \(\kappa\) is a limit of strongly compact cardinals.
		Then by the result of the previous paragraph, every successor strongly
		compact cardinal below \(\kappa\) is supercompact, so \(\kappa\) is a
		limit of supercompact cardinals.
	\end{proof}
\end{prp}

 The issue now is that there is no inner model theoretic reason whatsoever to
 believe that CUA is consistent with very large cardinals, but it cannot be that
 easy to refute:
\begin{prp}[UA] Suppose \(j_0:V\to M_0\) and \(j_1:V\to M_1\) witness that
	\textnormal{CUA} is false and \(\lambda = \min
	\{\textsc{crt}(j_0),\textsc{crt}(j_1)\}\). Then some cardinal \(\kappa <
	\lambda\) is \(\lambda\)-supercompact.
	\begin{proof}[Sketch] Since \(\lambda\) is measurable, it suffices to show
		that some \(\kappa < \lambda\) is \(\gamma\)-supercompact for all
		\(\gamma < \lambda\). Assume towards a contradiction that no cardinal
		below \(\lambda\) has this property. Then by \cref{LimitLma}, for any
		ultrapower embedding \(i : V\to N\), \(i[\lambda]\subseteq \lambda\).
		Let \[(i_0,i_1) : (M_0,M_1)\to N\] be the pushout of \((j_0,j_1)\). Let
		\(W\) be a countably complete ultrafilter such that \(M_W = N\) and
		\(j_W = i_0\circ j_0 = i_1\circ j_1\). By the analysis of ultrafilters
		internal to a pushout (\cref{PushoutInternal}), \(W\) is
		\(\lambda\)-internal. Thus \(j_W[\lambda]\subseteq \lambda\) and \(W\)
		is  \(\lambda\)-internal, so the internal relation theoretic criterion
		for completeness (\cref{InternalComplete}) implies that \(W\) is
		\(\lambda\)-complete. Thus \(\textsc{crt}(i_0)\geq\textsc{crt}(i_0\circ
		j_0) = \textsc{crt}(j_W) = \kappa\), and similarly \(\textsc{crt}(i_1)
		\geq \kappa\). This contradicts that \(j_0\) and \(j_1\) witness the
		failure of CUA.
	\end{proof}
\end{prp}

One can do a bit better using the following fact, whose proof we omit:
\begin{prp}[UA] Suppose \textnormal{CUA} fails. Then there are irreducible
	ultrafilters \(U_0\) and \(U_1\) such that \(j_{U_0}\) and \(j_{U_1}\)
	witness the failure of \textnormal{CUA}.\qed
\end{prp}
Since UA implies the linearity of the Mitchell order on normal ultrafilters
(\cref{UAMO}), CUA cannot fail for a pair of normal ultrafilters, and hence the
analysis of normality and irreducible ultrafilters (\cref{MuMeasure}) implies
that \(\min\{\textsc{crt}(j_{U_0}),\textsc{crt}(j_{U_1})\}\) is a
\(\mu\)-measurable cardinal. One can push this quite a bit further, but not far
enough to answer the following question:
\begin{qst}
	Is CUA consistent with the existence of cardinals \(\kappa < \lambda\) that
	are both \(\lambda^+\)-supercompact?
\end{qst}
The most interesting possibility is that large cardinals refute CUA. In any
case, unless one can prove CUA from UA (or Weak Comparison), it is far from
well-justified. The analysis of the second strongly compact cardinal therefore
requires a different approach.
\subsection{Irreducible ultrafilters and supercompactness}
Given the techniques of the previous chapter, the obvious approach is to study
the \(\kappa\)-complete generalizations of Fr\'echet cardinals and the
ultrafilters \(\mathscr K_\lambda\).

\begin{defn}\index{Fr\'echet cardinal!\(\kappa\)-Fr\'echet cardinal}
	Suppose \(\kappa \leq \lambda\) are uncountable cardinals. Then \(\lambda\)
	is {\it \(\kappa\)-Fr\'echet} if there is a \(\kappa\)-complete uniform
	ultrafilter on \(\lambda\).
\end{defn}

\begin{defn}[UA]\index{Ketonen ultrafilter!\(mathscr K^\kappa_\lambda\)}
	Suppose \(\lambda\) is a \(\kappa\)-Fr\'echet cardinal. Then \(\mathscr
	K_\lambda^\kappa\) denotes the minimum \(\kappa\)-complete uniform
	ultrafilter on \(\lambda\) in the Ketonen order.
\end{defn}

Most of the key properties of \(\mathscr K_\lambda\) do not directly generalize
to \(\mathscr K_\lambda^\kappa\): the proofs seem to require
\(\text{UA}(\kappa)\). Essentially the only nontrivial UA result that lifts is
\cref{KetonenIrreducible}, the fact that \(\mathscr K_\lambda\) is irreducible
for regular \(\lambda\).
\begin{lma}[UA] Suppose \(\kappa\leq \lambda\) and \(\lambda\) is
	\(\kappa\)-Fr\'echet. Then \(\mathscr K_\lambda^\kappa\) is weakly normal.
	\begin{proof}
		Recall \cref{WeaklyNormalrRK}, which asserts that a uniform ultrafilter
		\(U\) on a cardinal \(\lambda\) is weakly normal if and only if for all
		\(W\rRK U\), \(\lambda_W < \lambda\). We will show that this holds for
		\(U = \mathscr K_\lambda^\kappa\). Suppose \(W\rRK \mathscr
		K_\lambda^\kappa\). Since \(W\RK \mathscr K_\lambda^\kappa\), \(W\) is
		\(\kappa\)-complete, and since \(W\rRK \mathscr K_\lambda^\kappa\),
		\(W\sE \mathscr K_\lambda^\kappa\). By the minimality of \(\mathscr
		K_\lambda^\kappa\), \(\lambda_W < \lambda\). 
	\end{proof}
\end{lma}

\begin{prp}[UA]\label{KetonenIrreducible2}
	Suppose \(\nu < \lambda\) and \(\lambda\) is a \(\nu^+\)-Fr\'echet regular
	cardinal.\footnote{It is necessary here to restrict to consideration of
	\(\mathscr K_\lambda^{\nu^+}\), rather than considering \(\mathscr
	K_\lambda^\kappa\) in general. In fact, \(\mathscr K_\lambda^\kappa\) is
	irreducible if and only if there is some \(\nu < \kappa\) such that
	\(\mathscr K_\lambda^\kappa = \mathscr K_\lambda^{\nu^+}\). This is closely
	related to Menas's Theorem (\cref{Menas}).} Then \(\mathscr
	K_\lambda^{\nu^+}\) is irreducible.
	\begin{proof}
		Let \(\mathscr K = \mathscr K_\lambda^{\nu^+}\). Suppose \(D\sD \mathscr
		K\). We must show that \(D\) is principal. Since \(\mathscr K\) is
		\(\nu^+\)-complete and \(D\RK \mathscr K\), \(D\) is \(\nu^+\)-complete,
		and in particular \(j_D(\nu) = \nu\). Since \(\mathscr K\) is weakly
		normal and \(D\sRK \mathscr K\), \(\lambda_D < \lambda\) by
		\cref{MinimalIncompressible}. Let \(j : V\to M\) be the ultrapower of
		the universe by \(\mathscr K\) and let \(h : M_D\to M\) be the unique
		internal ultrapower embedding such that \(h\circ j_D = j\). Then \(h\)
		is the ultrapower of \(M_D\) by \(\tr D {\mathscr K}\), and
		\(\textsc{crt}(h) \geq \textsc{crt}(j) > \nu = j_D(\nu)\). Thus \(\tr D
		{\mathscr K}\) is \(j_D(\nu^+)\)-complete in \(M_D\).
		
		Assume towards a contradiction that \(D\) is nonprincipal. By
		\cref{Pushdown}, \(\tr D {\mathscr K} \sE j_D(\mathscr K)\) in \(M_D\).
		Since  \(\tr D {\mathscr K}\) is \(j_D(\nu^+)\)-complete,  the
		\(\sE^{M_D}\)-minimality of \(j_D(\mathscr K)\) among
		\(j_D(\nu^+)\)-complete uniform ultrafilters on \(j_D(\lambda)\) implies
		that \(\lambda_{\tr D {\mathscr K}} < j_D(\lambda)\). Since
		\(j_D(\lambda)\) is \(M_D\)-regular, it follows that \(\delta_{\tr D
		{\mathscr K}} < j_D(\lambda)\). Since \(\lambda_D < \lambda\) and
		\(\lambda\) is regular, \(j_D(\lambda) = \sup j_D[\lambda]\) by
		\cref{UFContinuity}. Therefore there is some ordinal \(\alpha <
		\lambda\) such that \(\delta_{\tr D {\mathscr K}} < j_D(\alpha)\). But
		\(\alpha \in  j_D^{-1}[\tr D {\mathscr K}] = \mathscr K\), contradicting
		that \(\mathscr K\) is uniform. Thus our assumption was false, and in
		fact \(D\) is principal. This shows that \(\mathscr K\) is irreducible,
		as desired.
	\end{proof}
\end{prp}

Beyond \cref{KetonenIrreducible2}, the ultrafilters \(\mathscr
K^\kappa_\lambda\) turn out to be a bit of a red herring. The analysis of higher
supercompact cardinals does not proceed by generalizing the theory of \(\mathscr
K_\lambda\) to the ultrafilters \(\mathscr K^\kappa_\lambda\) but instead by
propagating the \(\lambda\)-supercompactness of \(\mathscr K_\lambda\) itself to
arbitrary irreducible ultrafilters. Recall that an ultrafilter \(U\) is
\(\lambda\)-irreducible if every ultrafilter \(D\D U\) such that \(\lambda_D <
\lambda\) is principal. The main theorems of this chapter, to which we refer
collectively as the {\it Irreducibility Theorem}, show that supercompactness and
irreducibility are equivalent:

\begin{repthm}{Irred1}[UA]\index{Irreducibility Theorem}
	Suppose \(\lambda\) is a successor cardinal or a strong limit singular
	cardinal and \(U\) is a countably complete uniform ultrafilter on
	\(\lambda\). Then the following are equivalent:
	\begin{enumerate}[(1)]
		\item \(j_U\) is \(\lambda\)-irreducible.
		\item \(j_U\) is \(\lambda\)-supercompact.
	\end{enumerate}
\end{repthm}

It does not seem to be possible to generalize this to the case that \(\lambda\)
is inaccessible, and instead we obtain the following theorem:

\begin{repthm}{Irred2}[UA] Suppose \(\lambda\) is an inaccessible cardinal and
	\(U\) is a countably complete ultrafilter on \(\lambda\). Then the following
	are equivalent:
	\begin{enumerate}[(1)]
		\item \(j_U\) is \(\lambda\)-irreducible.
		\item \(j_U\) is \({<}\lambda\)-supercompact and \(\lambda\)-tight.
	\end{enumerate}
\end{repthm}

We will use these two theorems to give a complete characterization of strongly
compact cardinals assuming UA:
\begin{repthm}{MenasUA}[UA]\index{Strongly compact cardinal!equivalence with supercompactness}
	Suppose \(\kappa\) is a strongly compact cardinal. Either \(\kappa\) is a
	supercompact cardinal or \(\kappa\) is a measurable limit of supercompact
	cardinals
\end{repthm}
\subsection{Outline of \cref{SCChapter2}}
We now outline the rest of this chapter.\\

\noindent {\sc\cref{IrredSection}.} We prove the main structural result of the section,
called the {\it Irreducibility Theorem}, from which all the other theorems flow.
The Irreducibility Theorem refers to a cluster of results (especially
\cref{SuccessorIrredThm} and \cref{LimitIrredThm}) that show an equivalence
between irreducibility and supercompactness.\\

\noindent {\sc\cref{IdentitySection}.} We use the Irreducibility Theorem to resolve the
Identity Crisis for strongly compact cardinals under UA. We also use it in
\cref{InternalSectionII} to completely characterize the internal relation in
terms of the Mitchell order.\\

\noindent {\sc\cref{VLCSection}.} We discuss the relationship between UA and very large
cardinals. We begin by (partially) analyzing the relationship between hugeness
and non-regular ultrafilters under UA (\cref{HugeThm}). We then turn to the
topic of cardinal preserving embeddings. We show that UA rules out such
embeddings (\cref{NoCardinalPreserve}), and more generally that local cardinal
preservation hypotheses are equivalent to rank-into-rank large cardinal large
cardinal axioms under UA (\cref{I3Thm}). Finally in \cref{PathologicalSection},
we discuss the structure of supercompactness at inaccessible cardinals, and in
particular the prospect that the local equivalence of strong compactness and
supercompactness breaks down there.\\
\section{The Irreducibility Theorem}\label{IrredSection}
In this section, we prove the central Irreducibility Theorem (\cref{Irred1} and
\cref{Irred2}). We begin in \cref{PseudocompactSection} by proving the forward
implication from supercompactness to irreducibility. This raises a central open
question (\cref{TightQ}) that will be discussed at greater length in
\cref{PathologicalSection}. The next two sections are devoted to proving the
preliminary lemmas necessary for the proof of the Irreducibility Theorem. In
\cref{KComparisonSection}, we prove two key lemmas regarding the comparison of
\(\mathscr K_\lambda\) with an arbitrary ultrafilter. In the very short
\cref{NormalComboSection}, we prove two theorems on the combinatorics of normal
ultrafilters that show up in the proof of the Irreducibility Theorem. Finally,
\cref{IrredProofSection} contains the proof of the Irreducibility Theorem as
well as some slightly more general theorems.
\subsection{Pseudocompactness and irreducibility}\label{PseudocompactSection}
In this short subsection, we prove the easy direction of the irreducibility
theorem: \(\lambda\)-supercompactness implies \(\lambda\)-irreducibility. In
fact, we will prove something slightly stronger. The following property is  a
priori somewhat weaker than \(\lambda\)-supercompactness, but already implies
\(\lambda\)-irreducibility.
\begin{defn}\index{Pseudocompact embedding}
	Suppose \(\lambda\) is a cardinal. An elementary embedding \(j : V\to M\) is
	said to be {\it \(\lambda\)-pseudocompact} if \(j\) is \(\gamma\)-tight for
	every cardinal \(\gamma \leq \lambda\).
\end{defn}

\begin{lma}
	An ultrapower embedding \(j :V \to M\) is \(\lambda\)-pseudocompact if and
	only if \(M\) has the \({\leq}\gamma\)-covering property for all
	\(\gamma\leq \lambda\).
	\begin{proof}
		This is an immediate consequence of the self-strengthening of tightness
		that holds for ultrapower embeddings (\cref{UltrapowerStrongC}).
	\end{proof}
\end{lma}

\begin{prp}\label{TightIrred}
	Suppose \(\lambda\) is a cardinal and \(U\) is a countably complete
	ultrafilter. If \(j_U\) is \(\lambda\)-pseudocompact, then \(U\) is
	\(\lambda\)-irreducible.
	\begin{proof}
	Suppose \(D\D U\) and \(\lambda_D < \lambda\). We must show that \(D\) is
	principal. We first show that \(j_D\) is \(\lambda\)-pseudocompact. Since
	\(j_U\) is \(\lambda\)-pseudocompact, \(M_U\) has the
	\({\leq}\gamma\)-covering property for all \(\gamma \leq \lambda\). Since
	\(D\D U\), \(M_U\subseteq M_D\). It follows that \(M_D\) has the
	\({\leq}\gamma\)-covering property for all \(\gamma \leq \lambda\): suppose
	\(\gamma\leq \lambda\) and \(A\) is a set of ordinals of cardinality
	\(\gamma\); then \(A\) is contained in a set \(B\in M_U\) such that
	\(|B|^{M_U} \leq \gamma\), and since \(M_U\subseteq M_D\), we have \(B\in
	M_D\) and \(|B|^{M_D} \leq \gamma\), as desired. Thus \(j_D\) is
	\(\lambda\)-pseudocompact. 
	
	In particular, since \(\lambda_D < \lambda\), \(D\) is
	\(\lambda_D^+\)-tight. Assume towards a contradiction that \(D\) is
	nonprincipal. By \cref{Ineqs}, \(j_D(\lambda_D) > \lambda_D^+\). Thus \(D\)
	is \((\lambda_D^+,\delta)\)-tight where \(\delta = \lambda_D^+ <
	j_D(\lambda_D)\). This contradicts \cref{UFStrongCBound}, which states that
	if \(\eta\) is a cardinal and \(Z\) is a nonprincipal countably complete
	ultrafilter such that \(\lambda_Z < \eta\), then \(Z\) is not
	\((\eta,\delta)\)-tight for any \(\delta < j_Z(\eta)\). Thus \(D\) is
	principal, as desired.
\end{proof}
\end{prp}

The only known instances of \(\lambda\)-pseudocompact elementary embeddings that
are not \(\lambda\)-supercompact come from large cardinal axioms at the level of
rank-into-rank cardinals. Specifically, assume the axiom \(I_2\). Thus there is
a cardinal \(\lambda\) and an elementary embedding \(j : V\to M\) such that
\(\textsc{crt}(j) < \lambda\), \(j(\lambda) = \lambda\), and
\(V_\lambda\subseteq M\). The embedding \(j\) is not \(\lambda\)-supercompact by
the Kunen Inconsistency Theorem, but \(j\) is trivially
\(\lambda\)-pseudocompact since \(j[\lambda]\subseteq \lambda\). In fact, \(j\)
is \(\lambda^{+\alpha}\)-pseudocompact for all \(\alpha < \textsc{crt}(j)\). On
the other hand, there are no known examples of {\it ultrapower} embeddings that
are \(\lambda\)-pseudocompact but not \(\lambda\)-supercompact. In fact, it is
not known whether it is consistent that such an example exists:

\begin{qst}[ZFC]\label{TightQ}
Suppose \(\lambda\) is a cardinal and \(j : V\to M\) is a
\(\lambda\)-pseudocompact ultrapower embedding. Must \(j\) be
\(\lambda\)-supercompact?
\end{qst}

The natural inclination is to conjecture that the answer is no: typically large
cardinal properties formulated in terms of covering do not imply
supercompactness in ZFC. But the problem turns out to be much more subtle than
one might expect. 

We highlight below the most basic instance of this problem (in simple English):
\begin{qst}
	Suppose \(j : V\to M\) is an elementary embedding with critical point
	\(\kappa\) such that \(\textnormal{cf}^M(\sup j[\kappa^+]) = \kappa^+\).
	Must \(j[\kappa^+]\) belong to \(M\)?
\end{qst}

We devote the final section of this dissertation (\cref{PathologicalSection}) to
the relationship between \cref{TightQ} and the Inner Model Problem. 

On this topic, let us mention an interesting way in which tightness can act as a
stand-in for strength:
\begin{lma}\label{TightContinuum}
	Suppose \(j : V\to M\) is an elementary embedding, \(\lambda\) is a
	cardinal, \(\delta\) is an \(M\)-cardinal, and \(j\) is
	\((\lambda,\delta)\)-tight. Then \(2^\lambda\leq (2^\delta)^M\).
	\begin{proof}
		Fix \(B\in M\) such that \(|B|^M = \delta\) and \(j[\lambda]\subseteq
		B\). Then the map \(f : P(\lambda)\to P(B)\cap M\) defined by \(f(S) =
		j(S)\cap B\) is an injection: if \(S\neq T\), then fix \(\alpha\in
		S\mathrel{\triangle} T\), and note that since \(j[\lambda]\subseteq B\),
		\(j(\alpha)\in (j(S)\mathrel{\triangle} j(T))\cap B =
		f(S)\mathrel{\triangle} f(T)\). Since \(|P(B)\cap M|^M = (2^\delta)^M\)
		it follows that \(2^\lambda \leq |(2^\delta)^M| \leq (2^\delta)^M\).
	\end{proof}
\end{lma}
As a sample application (and a brief diversion), suppose \(\kappa\) is a
cardinal such that for all cardinals \(\lambda \geq \kappa\), there is a
\(\lambda\)-tight embedding \(j: V\to M\) such that \(j(\kappa) > \lambda\).
Then the Generalized Continuum Hypothesis cannot fail first above \(\kappa\). To
see this, assume that for all cardinals \(\gamma < \kappa\), \(2^\gamma =
\gamma^+\). Fix \(\lambda\geq \kappa\). Let \(j : V\to M\) be a
\(\lambda\)-tight embedding with \(j(\kappa) >\lambda\). Then in \(M\),
\(2^\lambda = \lambda^+\). Therefore \(2^\lambda \leq (2^\lambda)^M \leq
(\lambda^+)^M \leq \lambda^+\), so \(2^\lambda = \lambda^+\).
\subsection{Translations of \(\mathscr K_\lambda\)}\label{KComparisonSection}
Suppose \(U\) is a \(\lambda\)-irreducible uniform ultrafilter on a successor
cardinal \(\lambda\). The Irreducibility Theorem asserts that \(j_U\) is
\(\lambda\)-supercompact. The proof proceeds by analyzing the pushout comparison
of \((j_{\mathscr K_\lambda},j_U)\) where \(\lambda\) is a Fr\'echet successor
cardinal. In this section, we will prove a number of lemmas regarding this
pushout that amount to pieces of this analysis.

The universal property of \(\mathscr K_\lambda\) (\cref{UniversalProperty})
identifies the pushout of \((j_{\mathscr K_\lambda},j_U)\) when
\(\text{cf}^{M_U}(\sup j_U[\lambda])\) is not Fr\'echet in \(M_U\): in fact,
\(\mathscr K_\lambda\D U\), so the pushout is given by the unique internal
ultrapower embedding \(h : M_{\mathscr K_\lambda}\to M_U\). It turns out that
the universal property is powerful enough to yield an analysis of this
comparison even when \(\text{cf}^{M_U}(\sup j_U[\lambda])\) is a Fr\'echet
cardinal of \(M_U\). The following lemma tells us which ultrafilter is hit on
the \(M_U\)-side of the comparison:

\begin{lma}[UA]\label{KetonenTranslation}
	Suppose \(\lambda\) is a regular Fr\'echet cardinal and \(U\) is a countably
	complete ultrafilter. Let \(\delta = \textnormal{cf}^{M_U}(\sup
	j_U[\lambda])\). 
	\begin{itemize}
		\item Suppose \(\delta\) is not Fr\'echet in \(M_U\). Then \(\tr U
		{\mathscr K_\lambda}\) is principal in \(M_U\).
		\item Suppose \(\delta\) is Fr\'echet in \(M_U\). Then \(\tr U {\mathscr
		K_\lambda} \cong (\mathscr K_\delta)^M\).
	\end{itemize}
	\begin{proof}
		The first bullet point is immediate from the universal property of
		\(\mathscr K_\lambda\)  (\cref{UniversalProperty}): we have \(\mathscr
		K_\lambda\D U\), so by \cref{TransRF0}, \(\tr U {\mathscr K_\lambda}\)
		is principal in \(M_U\). Therefore assume instead that \(\delta\) is
		Fr\'echet in \(M_U\). 
		
		Let \(Z = \tr U {\mathscr K_\lambda}\). We claim that in \(M_U\), \(Z\)
		is a \(\sE\)-minimal element of the set of countably complete
		ultrafilters \(W\) on \(j_U(\lambda)\) with \(\delta_W \geq \sup
		j_U[\lambda]\). Clearly \(\delta_Z\geq \sup j_U[\lambda]\), since
		otherwise \(\delta_{j_U^{-1}[Z]} < \lambda\) contradicting that
		\(j_U^{-1}[Z] = \mathscr K_\lambda\). Suppose \(W\in j_U(\mathscr
		B(\lambda))\) and \(W\sE Z\) in \(M_U\), and we will show \(\delta_W <
		\sup j_U[\lambda]\). Let \(\bar W = j_U^{-1}[W]\). Then \(\tr U {\bar W}
		\E W\sE Z = \tr U {\mathscr K_\lambda}\). By \cref{OrderPreserving}, it
		follows that \(\bar W \sE \mathscr K_\lambda\). Since \(\mathscr
		K_\lambda\) is the \(\sE\)-least uniform ultrafilter on the regular
		cardinal \(\lambda\), \(\delta_{\bar W} < \lambda\). But
		\(j_U(\delta_{\bar W})\in W\), so \(\delta_W \leq j_U(\delta_{\bar W}) <
		\sup j_U[\lambda]\).
		
		Applying the analysis of \(\sE\)-minimal tail uniform ultrafilters
		(\cref{KetonenOrdinals}) in \(M_U\), it follows that in \(M_U\), there
		is a Ketonen ultrafilter \(D\) on \(\text{cf}^{M_U}(\sup j_U[\lambda]) =
		\delta\) that is isomorphic to \(Z\). Applying UA in \(M_U\), \(D =
		\mathscr K_\delta\), the unique Ketonen ultrafilter on \(\delta\).
	\end{proof}
\end{lma}

The analysis of the \(M_{\mathscr K_\lambda}\)-side of the comparison is much
more subtle, and uses the following fact:
\begin{lma}[UA]\label{KFactor}
	Suppose \(\lambda\) is a nonisolated regular Fr\'echet cardinal. Let \(M =
	M_{\mathscr K_\lambda}\). Suppose \(i : M \to N\) is an internal ultrapower
	embedding. Then there is a countably complete ultrafilter \(D\) of \(M\)
	with \(\lambda_D < \lambda\) such that there is an internal ultrapower
	embedding \(h : (M_D)^M\to N\) with \(h\circ j_D = i\) and \(\textsc{crt}(h)
	> j_D(\lambda)\).
\end{lma}
The proof uses an analysis of \(\lambda^\sigma\) in \(M_{\mathscr K_\lambda}\)
which is similar to \cref{EtaMeasurableClm} of \cref{IsoIPoint}:
\begin{lma}[UA]\label{KLambdaSigma}
	Suppose \(\lambda\) is a nonisolated regular Fr\'echet cardinal. Let \(j :
	V\to M\) be the ultrapower of the universe by \(\mathscr K_\lambda\). Then
	\((\lambda^\sigma)^M\) is a measurable cardinal of \(M\).
	\begin{proof}
		 By \cref{KetonenTight} and \cref{NonisolatedCompact}, \(j\) is
		 \(\lambda\)-tight and therefore \(\text{cf}^M(\sup j[\lambda]) =
		 \lambda\). Therefore by the definition of \(\mathscr K_\lambda\) (or
		 more precisely, \cref{KetonenChar}), \(\lambda\) is not Fr\'echet in
		 \(M\). 
				
		Let \(\eta = (\lambda^\sigma)^M\). Assume towards a contradiction that
		\(\eta\) is not measurable. Let \(i : M\to N\) be the ultrapower of
		\(M\) by \((\mathscr K_\eta)^M\) and let \(a = \id_{(\mathscr
		K_\eta)^M}\). 
		
		We claim that every countably complete \(N\)-ultrafilter \(D\) on
		\(\lambda\) belongs to \(M\). For any such \(D\), \(j_D^N\circ i\) is
		continuous at \(\lambda\): \(i\) is continuous at \(\lambda\) because
		\(i\) is internal to \(M\) and \(\lambda\) is not Fr\'echet in \(M\),
		while \(j_D^N\) is continuous at \(i(\lambda)\) since \(i(\lambda)\) is
		an \(N\)-regular cardinal with \(i(\lambda) > \lambda \geq \lambda_D\),
		and combining these observations: \[j_D^N(i(\lambda)) = \sup
		j_D^N[i(\lambda)] = \sup j_D^N[\sup i[\lambda] ] = \sup j_D^N\circ
		i[\lambda]\] Thus by the characterization of internal ultrapower
		embeddings of \(M_{\mathscr K_\lambda}\) (\cref{EmbeddingChar}),
		\(j_D^N\circ i\) an internal ultrapower embedding of \(M\). Since
		\(j_D^N\) can be defined at a typical element of \(N\) by setting
		\[j_D^N([f]_{(\mathscr K_\eta)^M}) = j_D^N\circ i(f)(j_D^N(a))\] it
		follows that \(j_D^N\) is definable over \(M\). Thus \(D\in M\).
		Applying inside \(M\) the characterization of countably complete
		ultrafilters amenable to an isolated ultrapower
		(\cref{IsolatedUFAmenableChar}), we have that \(D\in N\).
		
		\cref{GeneralStrong} states that if \(\kappa\) is \(\lambda\)-strongly
		compact and \(Q\) is a \({<}\kappa\)-closed inner model such that every
		\(\kappa\)-complete ultrafilter \(U\) on \(\lambda\) is amenable to
		\(Q\), then \(P(\lambda)\subseteq Q\). By \cref{NonisolatedCompact},
		\(\kappa_\lambda\) is \(\lambda\)-strongly compact. Moreover \(N\) is a
		\({<}\kappa_\lambda\)-closed (indeed \({<}\lambda\)-closed by
		\cref{DeltaLemma}) inner model such that every countably complete
		\(N\)-ultrafilter on \(\lambda\) belongs to \(\lambda\). It therefore
		follows that \(P(\lambda)\subseteq N\). But then \(\mathscr K_\lambda\)
		itself is an \(N\)-ultrafilter, so \(\mathscr K_\lambda\in N\). Since
		\(N\subseteq M\), this implies \(\mathscr K_\lambda\in M = M_{\mathscr
		K_\lambda}\), so \(\mathscr K_\lambda\mo \mathscr K_\lambda\),
		contradicting the irreflexivity of the Mitchell order (\cref{MOStrict}).
		
	\end{proof}
\end{lma}

\begin{proof}[Proof of \cref{KFactor}] By \cref{KLambdaSigma}, \(\eta =
	(\lambda^\sigma)^M\) is a measurable cardinal that is not a limit of
	Fr\'echet cardinals. The theorem follows by applying in \(M\) the fact that
	ultrapower embeddings can be factored across strong limit cardinals that are
	not limits of Fr\'echet cardinals (\cref{Cutpoint}).
\end{proof}

\cref{KFactor} has the following curious and sometimes useful corollary:

\begin{lma}[UA]\label{AlwaysTight}\index{Tightness of an elementary embedding!at inaccessible cardinals}
	Suppose \(\lambda\) is a strongly inaccessible cardinal such that one of the
	following holds: 
	\begin{itemize}
		\item \(\lambda\) is Fr\'echet.
		\item \(\lambda^\sigma\) is measurable. 
	\end{itemize}
	Then every ultrapower embedding is \(\lambda\)-tight.
	\begin{proof}
		Suppose \(U\) is a countably complete ultrafilter. We will show that
		\(j_U\) is \(\lambda\)-tight. 
		
		Assume first that \(\lambda\) is not Fr\'echet. Then by assumption
		\(\eta = \lambda^\sigma\) is measurable. By \cref{Cutpoint}, there is a
		countably complete ultrafilter \(D\) with \(\lambda_D < \eta\) such that
		there is an elementary embedding \(k : M_D\to M_U\) with \(k\circ j_D =
		j_U\) and \(\textsc{crt}(k) \geq \eta\). Since \(\lambda_D < \eta\), in
		fact \(\lambda_D < \lambda\), so \(j_D(\lambda) = \lambda\) since
		\(\lambda\) is inaccessible. But since \(\textsc{crt}(k) >
		j_D(\lambda)\),  \(j_U(\lambda) = j_D(\lambda) = \lambda\). Therefore
		\(j_D\) is vacuously \(\lambda\)-tight.
		
		Assume instead that \(\lambda\) is Fr\'echet. Let \((h,i) :
		(M_U,M_{\mathscr K_\lambda})\to N\) be the pushout of \((j_U,j_{\mathscr
		K_\lambda})\). Applying \cref{KFactor}, \(i\) factors in such a way that
		we can conclude that \(i(\lambda) = \lambda\) by the argument of the
		previous paragraph. Since \(\mathscr K_\lambda\) is \(\lambda\)-tight by
		\cref{GeneralThm} and \(i\) is vacuously \(\lambda\)-tight, \(i\circ
		j_{\mathscr K_\lambda}\) is \(\lambda\)-tight. In other words, \(N\) has
		the \({\leq}\lambda\)-covering property. Since \(N\subseteq M_U\) and
		\(N\) has the \({\leq}\lambda\)-covering property, \(M_U\) has the
		\({\leq}\lambda\)-covering property. Therefore \(j_U\) is
		\(\lambda\)-tight, as desired.
	\end{proof}
\end{lma}

\subsection{Elementary embeddings and normal filters}\label{NormalComboSection}
In this short subsection, we prove some combinatorial constraints on comparisons
involving normal filters. Suppose \(U\) and \(W\) are countably complete
ultrafilters on a cardinal \(\kappa\). A question that often arises in the
context of UA is what sort of \(M_W\)-ultrafilters \(Z\) on \(j_W(\kappa)\) pull
back to \(U\) in the sense that \(U = j_W^{-1}[Z]\). Such \(M_W\)-ultrafilters
arise from any comparison of \((j_U,j_W)\). Focusing on a more specific
question, assume \(U\) is normal, and suppose \(Z\) is a tail uniform
\(M_W\)-ultrafilter on \(j_W(\kappa)\) with \(j_W^{-1}[Z] = U\). Must \(Z =
j_W(U)\)? The following lemma, which has almost certainly been discovered
before, tells us that the answer is yes:
\begin{lma}\label{NormalGeneration}
	Suppose \(\mathcal F\) is a normal fine filter on a set \(Y\), and \(W\) is
	an ultrafilter on \(X = \bigcup Y\). Then \(j_W(\mathcal F)\) is the unique
	\(M\)-filter on \(j_W(Y)\) that extends \(j_W[F]\) and concentrates on
	\(\{\sigma\in j_W(Y) : \id_W\in \sigma\}\). In particular, \(j_W(\mathcal
	F)\) is the unique fine \(M\)-filter on \(j_W(Y)\) extending \(j_W[F]\).
	\begin{proof}
		Suppose \(A\in j_W(\mathcal F)\). We will find \(B\in \mathcal F\) such
		that \[j_W(B)\cap \{\sigma\in j_W(Y) : \id_W\in \sigma\}\subseteq A\]
		Fix a function \(G : X\to \mathcal F\) such that \(A = j_W(G)(\id_W)\).
		Let \[B = \triangle_{x\in X}G(x)\] Suppose \(\tau\in j_W(B)\cap
		\{\sigma\in j_W(Y) : \id_W\in \sigma\}\). We will show that \(\tau\in
		A\). Since \(\tau\) belongs to \(j_W(B) = \triangle_{x\in
		j_W(X)}j_W(G)(x)\), the definition of the diagonal intersection
		operation implies that \(\tau\in j_W(G)(x)\) for all \(x\in \tau\). But
		\(\id_W\in \tau\), and hence \(\tau\in j_W(G)(\id_W) = A\).
	\end{proof}
\end{lma}

In general, one must adjoin the set \(\{\sigma\in j(Y) : \id_W\in \sigma\}\) in
order to generate all of \(\mathcal F\). Suppose \(\mathcal F\) is a normal fine
ultrafilter on \(Y\) and \(W\) is an ultrafilter on \(X = \bigcup Y\). Then
\(j_W[\mathcal F]\) generates \(j_W(\mathcal F)\) if and only if there is some
\(\tau\in Y\) such that \(W\) concentrates on \(\tau\) and \(\mathcal F\)
concentrates on \(\{\sigma\in Y : \tau\subseteq \sigma\}\).

To better explain how this lemma is related to UA, we offer a sample corollary:
\begin{cor}[UA]\label{LeastExtension}
	Suppose \(F\) is a normal filter on a cardinal \(\kappa\). Let \(U\) be the
	\(\sE\)-least countably complete ultrafilter on \(\kappa\) that extends
	\(F\). Then for all \(D\sE U\), \(D\I U\).
	\begin{proof}
		Suppose \(D\sE U\). We claim \(j_D(U)\E \tr D U\) in \(M_D\), which
		implies \(D\I U\) by the theory of the internal relation (\cref{IChar}).
		Since \(j_D(U)\) is the \(\sE^{M_D}\)-least countably complete
		ultrafilter of \(M_D\) that extends \(j_D(F)\), it suffices to show that
		\(j_D(F)\subseteq \tr D U\). Of course \(j_D[F]\subseteq \tr D U\) since
		\(j_D^{-1}[\tr D U] = U\). Moreover since \(D\sE U\), we must have that
		\(\tr D U\) concentrates on ordinals greater than \(\id_D\) (since
		otherwise \(\tr D U\) witnesses \(U\E D\)). In other words, \(\{\alpha <
		j_D(\kappa): \id_D\in \alpha\}\in \tr D U\). Therefore by
		\cref{NormalGeneration}, \(j_D(F)\subseteq \tr D U\), as desired.
	\end{proof}
\end{cor}
Here is an intriguing consequence of \cref{LeastExtension}. Suppose \(\kappa\)
is a regular cardinal and \(F\) is the \(\omega\)-club filter on \(\kappa\).
Suppose \(F\) extends to a countably complete ultrafilter. Mitchell
\cite{MitchellSkies} showed that this hypothesis is equiconsistent with a
measurable cardinal of Mitchell order \(\omega\), but assuming UA, it implies
that there is a \(\mu\)-measurable cardinal and quite a bit more. The reason is
that \cref{LeastExtension} shows that the \(\sE\)-least extension of \(F\) is
irreducible; clearly it is not normal, so we can apply the dichotomy between
normal ultrafilters and \(\mu\)-measurability (\cref{MuDichotomy}).
\begin{qst} Let \(T\) be the theory  ZFC + UA +
	there is a regular cardinal \(\delta\) that carries a countably complete ultrafilter
	extending the closed unbounded filter.
	What is the consistency strength of \(T\)?
	Does \(T\) imply that there is an inner model with a superstrong cardinal?
	Does \(T\) imply that there is a superstrong cardinal?
\end{qst}
The only models of \(T\) that we know of are Jensen's canonical inner models
in the vicinity of a subcompact cardinal.

As a corollary of \cref{NormalGeneration}, we have a similar unique extension
theorem for isonormal ultrafilters on regular cardinals. We begin with a
corollary of Solovay's Lemma (\cref{SolovayLemma}) that explains the statement
of \cref{IsoExtension}:
\begin{lma}\label{WNStationary}
	Suppose \(\lambda\) is a regular cardinal and \(W\) is a countably complete
	weakly normal ultrafilter on \(\lambda\). Suppose \(\langle S_\xi : \xi <
	\lambda\rangle\) is a partition of \(S^\lambda_\omega\) into stationary
	sets. Then for any \(\xi < \lambda\), \(W\) concentrates on the set of
	\(\alpha < \lambda\) such that \(S_\xi\) is stationary in \(\alpha\).
	\begin{proof}
		Let \(j : V\to M\) be the ultrapower of the universe by \(W\). Then
		since \(W\) is weakly normal, \(\id_W = \sup j[\lambda]\). Let \(\langle
		T_\xi : \xi < j(\lambda)\rangle = j(\langle S_\alpha : \alpha <
		\lambda\rangle)\). By Solovay's Lemma (\cref{StationaryPartition}),
		\[j[\lambda] = \{\xi < j(\lambda) : T_\xi\text{ is stationary in }\sup
		j[\lambda]\}\] In particular, if \(\xi < \lambda\), then \(M\) satisfies
		that \(T_{j(\xi)}\) is stationary in \(\id_W\), and so by \L o\'s's
		Theorem, \(W\) concentrates on the set of \(\alpha < \lambda\) such that
		\(S_\xi\) is stationary in \(\alpha\).
	\end{proof}
\end{lma}

\begin{lma}\label{IsoExtension}
	Suppose \(\lambda\) is a regular cardinal, \(W\) is an isonormal ultrafilter
	on \(\lambda\), and \(D\) is a countably complete ultrafilter on
	\(\lambda\). Let \(\langle S_\xi : \xi < \lambda\rangle\) be a partition of
	\(S^\lambda_\omega\) into stationary sets, and let \(\langle T_\xi : \xi <
	j_D(\lambda)\rangle = j_D(\langle S_\xi : \xi < \lambda\rangle)\). Let \[A =
	\{\alpha < j_D(\lambda) : M_D\vDash T_{\id_D}\text{ is stationary in
	}\alpha\}\] Then \(j_D(W)\) is the unique \(M_D\)-filter on \(j_D(\lambda)\)
	that extends \(j_D[W]\) and concentrates on \(A\).
	\begin{proof}
		Let \(\mathcal U\) be the normal fine ultrafilter on \(P(\lambda)\)
		isomorphic to \(W\). Let \(g : P(\lambda)\to \lambda+1\) be the sup
		function \[g(\sigma) = \sup \sigma\] By Solovay's Lemma
		(\cref{SolovayCor}), \(g_*(\mathcal U) = W\). Let \(f : \lambda\to
		P(\lambda)\) be the function defined by \[f(\alpha) = \{\xi < \lambda :
		S_\xi\text{ is stationary in }\alpha\}\] By the proof of Solovay's
		Lemma, for any \(A\subseteq \lambda\), \(f[A]\) and \(g^{-1}[A]\) are
		equal modulo \(\mathcal U\). Thus since \(W = g_*(\mathcal U)\), \[W =
		\{A\subseteq \lambda : f[A]\in \mathcal U\}\]
		
		By \cref{NormalGeneration}, \(j_D(\mathcal U)\) is the unique
		\(M_D\)-filter on \(j_D(Y)\) that extends \(j_D[\mathcal U]\) and
		concentrates on \(\{\sigma\in j_D(P(\lambda)) : \id_D \in \sigma\}\).
		Since \(j_D(W) = \{A\subseteq \lambda : j_D(f)[A]\in j_D(\mathcal
		U)\}\), it follows that \(j_D(W)\) is the unique \(M_D\)-filter on
		\(j_D(\lambda)\) that extends \(\{A : j_D(f)[A]\in j_D[\mathcal U]\}\)
		and concentrates on \[\{\alpha < \lambda: \id_D\in j_D(f)(\alpha)\} =
		\{\alpha < \lambda : M_D\vDash T_{\id_D}\text{ is stationary in
		}\alpha\}\] In other words, \(j_D(W)\) is the unique \(M_D\)-filter on
		\(j_D(\lambda)\) that extends \(j_D[W]\) and concentrates on \(A\), as
		desired.
	\end{proof}
\end{lma}

Let us include one more useful combinatorial fact, this time about pullbacks of
weakly normal ultrafilters. To state the lemma in the generality we will need,
we introduce a relativized version of the notion of a weakly normal ultrafilter.
\begin{defn}
	Suppose \(M\) is a transitive model of ZFC, \(\lambda\) is an \(M\)-regular
	cardinal,  and \(F\) is an \(M\)-filter on \(\lambda\). Then \(F\) is {\it
	weakly normal} if for all sequences \(\langle A_\alpha : \alpha <
	\lambda\rangle\in M\) of subsets of \(\lambda\) such that \(A_\alpha\in F\)
	for all \(\alpha < \lambda\) and \(A_\alpha\supseteq A_\beta\) for all
	\(\alpha \leq \beta < \lambda\), the diagonal intersection
	\(\triangle_{\alpha < \lambda} A_\alpha\) belongs to \(F\).
\end{defn}
We will really only need this notion for \(M\)-ultrafilters, in which case it
has the following familiar formulation:
\begin{lma}
	If \(M\) is a transitive model of \textnormal{ZFC}, \(\lambda\) is an
	\(M\)-regular cardinal, and \(U\) is an \(M\)-ultrafilter on \(\lambda\),
	then \(U\) is weakly normal if and only if \(\id_U = \sup
	j^M_U[\lambda]\).\qed
\end{lma}

\begin{lma}\label{WNPullback}
	Suppose \(\lambda\) is a regular cardinal and \(j :  V\to M\) is an
	elementary embedding that is continuous at \(\lambda\). Suppose \(F\) is a
	weakly normal \(M\)-filter on \(j(\lambda)\). Then \(j^{-1}[F]\) is a weakly
	normal filter on \(\lambda\).
	\begin{proof}
		Suppose \(\langle A_\alpha : \alpha < \lambda\rangle\) is a decreasing
		sequence of subsets of \(\lambda\) such that \(A_\alpha\in j^{-1}[F]\)
		for all \(\alpha < \lambda\). We must show that \(\triangle_{\alpha <
		\lambda} A_\alpha\in j^{-1}[F]\). Let \(\langle B_\beta : \beta <
		j(\lambda)\rangle = j(\langle A_\alpha : \alpha < \lambda\rangle)\).
		Since \(j(\triangle_{\alpha < \lambda} A_\alpha) = \triangle_{\beta <
		j(\lambda)} B_\beta\), it suffices to show that \(\triangle_{\beta <
		j(\lambda)} B_\beta\in F\).
		
		By the elementarity of \(j\), \(\langle B_\beta: \beta <
		j(\lambda)\rangle\) is a decreasing sequence of subsets of
		\(j(\lambda)\). We claim that for all \(\beta < j(\lambda)\),
		\(B_\beta\in F\). To see this, fix \(\beta < j(\lambda)\). Since \(j\)
		is continuous at \(\lambda\), there is some \(\alpha < \lambda\) such
		that \(\beta \leq j(\alpha)\). Now \(B_{j(\alpha)} = j(A_\alpha) \in F\)
		since \(A_\alpha\in j^{-1}[F]\). But \(B_{j(\alpha)}\subseteq B_\beta\)
		since \(\beta \leq j(\alpha)\) and \(\langle B_\beta: \beta <
		j(\lambda)\rangle\) is a decreasing sequence. Therefore \(B_\beta\in
		F\), as claimed. Since \(F\) is weakly normal, it follows that
		\(\triangle_{\beta < j(\lambda)} B_\beta\in F\). 
	\end{proof}
\end{lma}

\subsection{Proof of the Irreducibility Theorem}\label{IrredProofSection}
We will obtain the Irreducibility Theorem as an immediate consequence of the
following slightly more general fact:
\begin{thm}[UA]\label{SuccessorIrredThm}\index{Irreducibility Theorem}
	Suppose \(U\) is a countably complete ultrafilter and \(\lambda\) is a
	Fr\'echet successor cardinal. Then there is a countably complete ultrafilter
	\(D\) with \(\lambda_D < \lambda\) and an internal ultrapower embedding \(e
	: M_D\to M_U\) that is \(j_D(\lambda)\)-supercompact in \(M_D\).
	\begin{proof}
		Let \(j : V\to M\) be the ultrapower of the universe by \(\mathscr
		K_\lambda\) and let \(i : V\to N\) be the ultrapower of the universe of
		by \(U\). Let \[(i_*,j_*) : (M,N)\to P\] be the pushout of \((j,i)\).
		{\it Note that \(i_*\) denotes the embedding on the \(M\)-side of the
		comparison and \(j_*\) denotes the embedding on the \(N\)-side of the
		comparison.} The proof amounts to an analysis of \((i_*,j_*)\).
		\begin{figure}
			\center
			\includegraphics[scale=.6]{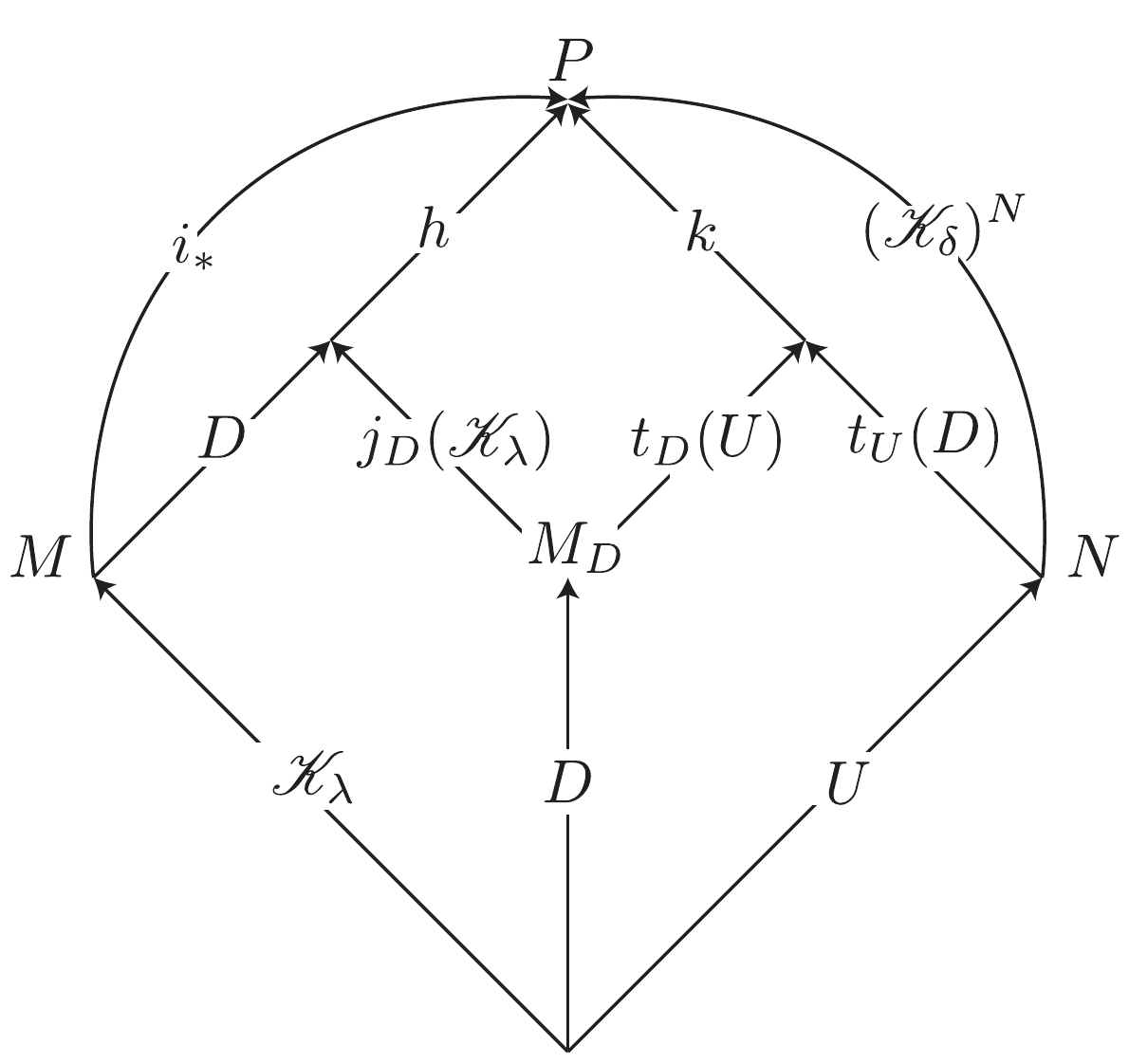}
			\caption{Diagram of the Irreducibility Theorem.}
		\end{figure}
		
		We first characterize \(j_*\). By definition (\cref{TransChar}), \(j_*\)
		is the ultrapower of \(N\) by \(\tr U {\mathscr K_\lambda}\). Let
		\[\delta = \text{cf}^N(\sup i[\lambda])\] By the analysis of
		translations of \(\mathscr K_\lambda\) (\cref{KetonenTranslation}), one
		of the following holds in \(N\): 
		\begin{itemize}
			\item \(\delta\) is not Fr\'echet and \(\tr U {\mathscr K_\lambda}\)
			is principal.
			\item \(\delta\) is Fr\'echet and  \(\tr U {\mathscr K_\lambda}\) is
			isomorphic to \((\mathscr K_\delta)^N\).
		\end{itemize}
		
		The hard part of the proof is the analysis of \(i_*\), the embedding on
		the \(M\)-side of the comparison of \((j,i)\). Let \(\eta\) be the least
		measurable cardinal of \(M\) above \(\lambda\). Applying \cref{KFactor},
		let \(D\) be a countably complete ultrafilter of \(M\) with \(\lambda_D
		< \lambda\) such that there is an internal ultrapower embedding \(h :
		(M_D)^M\to P\) with \(\textsc{crt}(h)\geq \eta\) and \(i_* = h\circ
		j_D^M\). We may assume without loss of generality that the underlying
		set of \(D\) is the cardinal \(\lambda_D\). Recall \cref{SuccessorThm},
		which states that \(M^\lambda\subseteq M\). In particular,
		\(P(\gamma)\subseteq M\), so \(D\) truly is an ultrafilter.
		
		The following are the two key claims:
		\begin{clm}\label{AgreementClm}
			\(\delta = j_D(\lambda)\) and \( \textnormal{Ord}^{j_D(\lambda)}\cap
			N = \textnormal{Ord}^{j_D(\lambda)} \cap P =
			\textnormal{Ord}^{j_D(\lambda)} \cap M_D\).
		\end{clm}
		\begin{clm}\label{DUClm}\(D\D U\).\end{clm}
		
		Assuming these claims, the conclusion of the theorem is immediate: by
		\cref{DUClm}, let \(e : M_D\to N\) be the unique internal ultrapower
		embedding such that \(e\circ j_D = i\); then \(e\) is
		\({j_D(\lambda)}\)-supercompact in \(M_D\) since \(
		\textnormal{Ord}^{j_D(\lambda)}\cap N =
		\textnormal{Ord}^{j_D(\lambda)}\cap M_D\) by \cref{AgreementClm}.
		
		We therefore focus on proving these two claims.
		\begin{proof}[Proof of \cref{AgreementClm}] We begin by showing
			\(\text{Ord}^{j_D(\lambda)}\cap M_D = \text{Ord}^{j_D(\lambda)}\cap
			P\). Since \(j: V\to M\) is a \(\lambda\)-supercompact ultrapower
			embedding, \(\text{Ord}^\lambda = \text{Ord}^\lambda\cap M\). By the
			elementarity of \(j_D\), \[\text{Ord}^{j_D(\lambda)}\cap M_D =
			\text{Ord}^{j_D(\lambda)}\cap j_D(M) = \text{Ord}^{j_D(\lambda)}\cap
			(M_D)^M\] The final equality follows from the fact that \(M\) is
			closed under \(\lambda\)-sequences and hence correctly computes the
			ultrapower of \(M\) by \(D\). But \(h : (M_D)^M\to P\) is an
			internal ultrapower embedding such that \(\textsc{crt}(h)\geq \eta >
			j_D(\lambda)\). Hence \(\text{Ord}^{j_D(\lambda)}\cap (M_D)^M =
			\text{Ord}^{j_D(\lambda)}\cap P\). Putting all this together, we
			have shown \[\text{Ord}^{j_D(\lambda)}\cap M_D =
			\text{Ord}^{j_D(\lambda)}\cap P\]
			
			One consequence of the agreement between \(M_D\) and \(P\), which we
			set down now for future use, is that \(j_D(\lambda)\) is a successor
			cardinal of \(P\): \(\lambda\) is a successor cardinal, so by
			elementarity, \(j_D(\lambda)\) is a successor cardinal of \(M_D\),
			and therefore since \(\text{Ord}^{j_D(\lambda)}\cap M_D =
			\text{Ord}^{j_D(\lambda)}\cap P\), \(j_D(\lambda)\) is a successor
			cardinal of \(P\).
			
			Next, we show that \(\delta = j_D(\lambda)\). To do this, we
			calculate the \(P\)-cofinality of the ordinal \(\sup j_*\circ
			i[\lambda]\) in two different ways.
			
			On the one hand, we claim
			\begin{equation}\label{CofCalc1}\text{cf}^P(\sup j_*\circ i[\lambda]
			) = j_D(\lambda)\end{equation} We have that \( j_*\circ i = h \circ
			j_D \circ j = h\circ j_D(j)\circ  j_D\). Since \(\lambda_D <
			\lambda\), and \(\lambda\) is regular, \(j_D(\lambda) = \sup
			j_D[\lambda]\) (\cref{UFContinuity}). Now we calculate:
			\begin{align*}
				\textnormal{cf}^P(\sup h\circ j_D(j)\circ  j_D[\lambda]) 
					&= \textnormal{cf}^P(\sup h\circ j_D(j)[\sup j_D[\lambda]]) \\
					&= \textnormal{cf}^P(\sup h\circ j_D(j)[j_D(\lambda)]) \\
					&= \textnormal{cf}^{M_D}(\sup h\circ j_D(j)[j_D(\lambda)]) \\
					&=  j_D(\lambda)
			\end{align*} 
			The second-to-last equality uses the fact that
			\(\text{Ord}^{j_D(\lambda)}\cap P = \text{Ord}^{j_D(\lambda)}\cap
			M_D\). The final equality uses the fact that \(h\circ j_D(j)\) is
			increasing and definable over \(M_D\) and \(j_D(\lambda)\) is
			regular in \(M_D\). Putting everything together yields that
			\(\text{cf}^P(\sup j_*\circ i[\lambda] ) = j_D(\lambda)\), as
			claimed.
			
			On the other hand, we claim
			\begin{equation}\label{CofCalc2}\text{cf}^P(\sup j_*\circ i
			[\lambda]) = \text{cf}^P(\sup j_*[ \delta])\end{equation} Since
			\(\delta = \text{cf}^N(\sup i[\lambda])\), there is an increasing
			cofinal function \(f : \delta\to \sup i[\lambda]\) with \(f\in N\).
			Now \(\sup j_*\circ i [\lambda] =\sup j_*[\sup f[\delta]] = \sup
			j_*(f)[\sup j_*[\delta]]\). Thus \(j_*(f)\in P\) restricts to an
			increasing cofinal function from \(\sup j_*[\delta]\) to \(\sup
			j_*\circ i [\lambda]\). It follows that \(\text{cf}^P(\sup j_*\circ
			i [\lambda]) = \text{cf}^P(\sup j_*[ \delta])\), as desired.
			
			Combining \cref{CofCalc1} and \cref{CofCalc2}, we have shown
			\(\text{cf}^P(\sup j_*[ \delta]) = j_D(\lambda)\). To show \(\delta
			= j_D(\lambda)\), we must show \(\text{cf}^P(\sup j_*[ \delta]) =
			\delta\). In other words (applying the easy direction of
			\cref{KetonenTight}), we must show \(j_*\) is \(\delta\)-tight.
			
			Recall that \(j_*\) is the ultrapower of \(N\) by \(\tr U {\mathscr
			K_\lambda}\). If \(\tr U {\mathscr K_\lambda}\) is principal, then
			trivially \(j_*\) is \(\delta\)-tight. Therefore assume \(\tr U
			{\mathscr K_\lambda}\) is nonprincipal. By the second paragraph of
			this proof, \(N\) satisfies that \(\delta\) is Fr\'echet and \(\tr U
			{\mathscr K_\lambda}\) is isomorphic to \((\mathscr K_\delta)^N\). 
			
			It suffices to show that \(\delta\) is not isolated in \(N\). Then
			applying in \(N\) the analysis of \(\mathscr K_\delta\) at
			nonisolated cardinals \(\delta\) (\cref{GeneralThm}), \(j_*\) is
			\(\delta\)-tight. 
			
			Thus assume towards a contradiction that \(\delta\) is isolated in
			\(N\). In particular, \(\delta\) is a regular limit cardinal in
			\(N\). Moreover, by \cref{SigmaMahlo}, \((\mathscr K_\delta)^N\)
			concentrates on \(N\)-regular cardinals, so by \L o\'s's Theorem,
			\(\id_{(\mathscr K_\delta)^N} = \sup j_*[\delta]\) is regular in
			\(P\). Thus by \cref{CofCalc2}, \(\text{cf}^P(\sup j_*\circ
			i[\lambda]) = \sup j_*[\delta]\), and so by \cref{CofCalc1}, \(\sup
			j_*[\delta] = j_D(\lambda)\). Since \(\delta\) is a limit cardinal
			of \(N\), \(\sup j_*[\delta]\) is a limit cardinal of \(P\). This
			contradicts the fact (set down earlier) that \(j_D(\lambda)\) is a
			successor cardinal of \(P\). Thus our assumption that \(\delta\) is
			isolated in \(N\) was false. It follows that \(\delta\) is not
			isolated and hence \(j_*\) is \(\delta\)-tight, and hence
			\(\text{cf}^P(\sup j_*[\delta]) = \delta\), and hence by
			\cref{CofCalc1} and \cref{CofCalc2}, \(j_D(\lambda) = \delta\).
			
			We finally show that \(\text{Ord}^\delta \cap N = \text{Ord}^\delta
			\cap P\). If \(\tr U {\mathscr K_\lambda}\) is principal then \(P =
			N\), so this is obvious. If not, then \(j_* : N \to P\) is the
			ultrapower embedding associated to \((\mathscr K_\delta)^N\). Note
			that \(\delta = j_D(\lambda)\) is a successor cardinal of \(P\), and
			so since \(P\subseteq N\), \(\delta\) is a successor cardinal of
			\(N\). Thus by the analysis of Ketonen ultrafilters on successor
			cardinals (\cref{SuccessorThm}) applied in \(N\), \(j_*\) is
			\(\delta\)-supercompact. In particular, \(\text{Ord}^\delta \cap N =
			\text{Ord}^\delta \cap P\).
		\end{proof}
		We now turn to the proof that \(D\D U\).
		\begin{proof}[Proof of \cref{DUClm}] To show \(D\D U\), it suffices (by
		the definition of translation functions, or \cref{TransRF0}) to show
		that \(\tr U D\) is principal in \(N\).
			
		Let us first show that \[\tr U D\D \tr U {\mathscr K_\lambda}\] in
		\(N\). Note that \[(h\circ j_D(j),j_*) : (M_D,N)\to P\] is an internal
		ultrapower comparison of \((j_D,i)\). Since \[(j_{\tr D U}^{M_D}, j_{\tr
		U D}^{M_U}) : (M_D,N)\to (M_{\tr U D})^N\] is the pushout of \((j_D,i)\)
		(by \cref{TransChar}), it follows that there is an internal ultrapower
		embedding \(k : (M_{\tr U D})^N\to P\) such that \(k\circ j_{\tr D
		U}^{M_D} = h\circ j_D(j)\) and \(k\circ j_{\tr U D}^{M_U} = j_* = j_{\tr
		U {\mathscr K_\lambda}}^N\). The latter equation is equivalent to the
		statement that \(\tr U D\D \tr U {\mathscr K_\lambda}\) in \(N\).
		
		 Since \(\tr U {\mathscr K_\lambda}\) is either principal or isomorphic
		 to the ultrafilter \(({\mathscr K_{j_D(\lambda)}})^N\), which is
		 irreducible by \cref{KetonenIrreducible}, one of the following must
		 hold:
		\begin{enumerate}[(1)]
			\item \({j_D(\lambda)}\) is Fr\'echet in \(N\) and \(N\vDash \tr U D
			\cong ({\mathscr K_{j_D(\lambda)}})^N\).
			\item \(\tr U D\) is principal in \(N\).
		\end{enumerate}
		Our goal is to show that (2) holds, so to finish the proof of the claim,
		it suffices  to show that (1) fails. Towards this, we will prove the
		following subclaim:
		\begin{sclm}\label{FinalAgreementClm}
			Assume \({j_D(\lambda)}\) is Fr\'echet in \(N\). Then \((\mathscr
			K_{j_D(\lambda)})^N = j_D(\mathscr K_\lambda)\).
		\end{sclm}
		\begin{proof}[Proof of \cref{FinalAgreementClm}] We plan to prove the
			claim by applying our unique extension lemma for isonormal
			ultrafilters. By \cref{SuccessorThm}, \(\mathscr K_\lambda\) is an
			isonormal ultrafilter on \(\lambda\). By \cref{AgreementClm},
			\((\mathscr K_{j_D(\lambda)})^N\) is an \(M_D\)-filter on
			\(j_D(\lambda)\). Let \(\langle S_\xi : \xi < \lambda\rangle\) be a
			partition of \(S^\lambda_\omega\) into stationary sets. Let
			\(\langle T_\xi : \xi < j_D(\lambda)\rangle = j_D(\langle S_\xi :
			\xi < \lambda\rangle)\). By \cref{IsoExtension}, to show that
			\(j_D(\mathscr K_\lambda) = (\mathscr K_{j_D(\lambda)})^N\), it
			suffices to show that the following hold:
			\begin{enumerate}[(i)]
				\item \(\{\alpha < j_D(\lambda) : M_D \vDash T_{\id_D}\text{ is
				stationary in }\alpha\}\in (\mathscr K_{j_D(\lambda)})^N\).
				\item \(j_D[\mathscr K_\lambda]\subseteq (\mathscr
				K_{j_D(\lambda)})^N\).
			\end{enumerate}
			
			(i) will be proved by applying \cref{WNStationary}. Note that
			\(\langle T_\xi : \xi < j_D(\lambda)\rangle\) belongs to \(N\) and
			\(N\) satisfies that \(\langle T_\xi : \xi < j_D(\lambda)\rangle\)
			is a stationary partition of \(S^{j_D(\lambda)}_\omega\): this
			follows from the fact that \(P(j_D(\lambda))\cap N =
			P(j_D(\lambda))\cap M_D\) by \cref{AgreementClm} and \(\langle T_\xi
			: \xi < j_D(\lambda)\rangle\) is a stationary partition of
			\(S^{j_D(\lambda)}_\omega\) in \(M_D\). Since \((\mathscr
			K_{j_D(\lambda)})^N\) is a countably complete weakly normal
			ultrafilter of \(N\), \cref{WNStationary} implies that \((\mathscr
			K_{j_D(\lambda)})^N\) concentrates on \(\{\alpha < j_D(\lambda) :
			M_D \vDash T_{\xi}\text{ is stationary in }\alpha\}\) for any \(\xi
			< j_D(\lambda)\), and in particular \(\{\alpha < j_D(\lambda) : M_D
			\vDash T_{\id_D}\text{ is stationary in }\alpha\}\in (\mathscr
			K_{j_D(\lambda)})^N\), as desired.
			
			Towards (ii), let \(W = j_D^{-1}[(\mathscr K_{j_D(\lambda)})^N]\).
			It suffices to show that \(W = \mathscr K_\lambda\). It is clear
			that \(W\) is a countably complete uniform ultrafilter on
			\(\lambda\). Recall that \(\mathscr K_\lambda\) is the unique
			Ketonen ultrafilter on \(\lambda\). Let \(A\) be the set of ordinals
			below \(\lambda\) that carry no countably complete tail uniform
			ultrafilter. By the definition of a Ketonen ultrafilter on a regular
			cardinal (\cref{KetonenRegDef}), to show \(W = \mathscr K_\lambda\),
			it suffices to show that the following hold:
			\begin{itemize}
				\item \(A\in W\).
				\item \(W\) is weakly normal.
			\end{itemize}
			
			Let us show that \(A\in W\). In other words, we must show that
			\(j_D(A)\in (\mathscr K_\delta)^N\). Note that \(j_D(A)\) is the set
			of ordinals less than \(j_D(\lambda) = \delta\) that carry no
			countably complete tail uniform ultrafilter in \(M_D\). By the
			definition of a Ketonen ultrafilter on a regular cardinal
			(\cref{KetonenRegDef}) applied in \(N\), \((\mathscr K_\delta)^N\)
			concentrates on the set of ordinals less than \(\delta\) that carry
			no countably complete tail uniform ultrafilter in \(N\). Thus to
			show that \(j_D(A)\in (\mathscr K_\delta)^N\), it suffices to show
			that if an ordinal less than \(\delta\) carries no countably
			complete tail uniform ultrafilter in \(N\), then it carries no
			countably complete tail uniform ultrafilter in \(M_D\). In fact we
			will show that for any ordinal \(\alpha < \delta\), \[\mathscr
			B^{M_D}(\alpha) = \mathscr B^N(\alpha)\] where \(\mathscr B(X)\)
			denotes the set of countably complete ultrafilters on \(X\).
			
			This is an application  \cref{MitchellLemma}, which asserts that if
			\(\gamma\) is a cardinal and \(Q\) is an ultrapower of the universe
			that is closed under \(\gamma\)-sequences, then for any ordinal
			\(\alpha < \gamma\), \(\mathscr B(\alpha) = \mathscr B^Q(\alpha)\).
			Fix an ordinal \(\alpha < \delta\). Applying \cref{MitchellLemma}
			in \(M_D\) to the ultrapower \(P\) of \(M_D\), which satisfies
			\(\textnormal{Ord}^\delta\cap P = \textnormal{Ord}^\delta\cap M_D\)
			by \cref{AgreementClm}, \[\mathscr B^{M_D}(\alpha) = \mathscr
			B^P(\alpha)\] Similarly, applying \cref{MitchellLemma} and
			\cref{AgreementClm} in \(N\) to \(P\), \[\mathscr B^N(\alpha) =
			\mathscr B^P(\alpha)\] Hence \(\mathscr B^{M_D}(\alpha) = \mathscr
			B^{N}(\alpha)\), as desired. This shows \(A\in W\).
			
			We now show that \(W\) is weakly normal. We do this by applying
			\cref{WNPullback}. Note that \((\mathscr K_{j_D(\lambda)})^N\) is a
			weakly normal \(M_D\)-ultrafilter since it is a weakly normal
			ultrafilter of \(N\) and \(P(j_D(\lambda))\cap M_D =
			P(j_D(\lambda))\cap N\). Therefore since \(j_D : V\to M_D\) is
			continuous at \(\lambda\), \cref{WNPullback} implies that
			\(j_D^{-1}[(\mathscr K_{j_D(\lambda)})]\) is weakly normal. In other
			words, \(W\) is weakly normal.
			
			Thus we have shown that \(W\) is a Ketonen ultrafilter on
			\(\lambda\), so \(W = \mathscr K_\lambda\). This implies (ii).
			
			As we explained above, (i), (ii), and \cref{IsoExtension} together
			imply \((\mathscr K_{j_D(\lambda)})^N = j_D(\mathscr K_\lambda)\),
			which proves the subclaim.
		\end{proof}
		Using \cref{FinalAgreementClm}, we show that (1) above does not hold. If
		\({j_D(\lambda)}\) is not Fr\'echet in \(N\), then obviously (1) does
		not hold, so assume instead that \({j_D(\lambda)}\) is Fr\'echet in
		\(N\). Let \(\mathscr K = (\mathscr K_{j_D(\lambda)})^N = j_D(\mathscr
		K_\lambda)\). Thus \(\mathscr K\in M_D\cap N\). 
		
		Recall that \(M_{\tr U D}^N\) is the target model of the pushout of
		\((j_D,i)\). Thus by the analysis of ultrafilters amenable to a pushout
		(\cref{PushoutInternal}), \(\mathscr K\cap M_{\tr U D}^N\in M_{\tr U
		D}^N\). On the other hand, we will show that \(\mathscr K\cap P\notin
		P\). By the strictness of the Mitchell order on nonprincipal
		ultrafilters (\cref{MOStrict}), \[\mathscr K\notin M_\mathscr K^{M_D} =
		j_D(M_{\mathscr K_\lambda}) = j_D(M)\] Recall that \(h : j_D(M)\to P\)
		is an internal ultrapower embedding, so in particular \(P\subseteq
		j_D(M)\), and hence \(\mathscr K\notin P\) since \(\mathscr K\notin
		j_D(M)\). Since \(P({j_D(\lambda)})\cap N = P({j_D(\lambda)})\cap P\) by
		\cref{AgreementClm}, it follows that \(\mathscr K = \mathscr K\cap P\),
		and so \(\mathscr K\cap P\notin P\). 
		
		We have \(\mathscr K\cap M_{\tr U D}^N\in M_{\tr U D}^N\) and \(\mathscr
		K\cap P\notin P\), so \(M_{\tr U D}^N\neq P\). Since \(P = M_{\mathscr
		K}^N\), it follows that \(\tr U D\) and \(\mathscr K\) are not
		isomorphic in \(N\): they have different ultrapowers. In other words,
		(1) above fails.
		
		Thus (2) holds, which proves \(D\D U\), establishing the claim.
		\end{proof}
		Having proved \cref{AgreementClm} and \cref{DUClm}, the theorem follows,
		as we explained after the statement of \cref{DUClm}.
	\end{proof}
\end{thm}

An immediate corollary of \cref{SuccessorIrredThm} is the following fact, which
will imply the Irreducibility Theorem:

\begin{cor}[UA]\label{SuccessorIrredCor}
	Suppose \(\lambda\) is a Fr\'echet successor cardinal and \(U\) is a
	\(\lambda\)-irreducible ultrafilter. Then \(j_U\) is
	\(\lambda\)-supercompact.
	\begin{proof}
		We begin with the case that \(\lambda\) is a successor cardinal. By
		\cref{SuccessorIrredThm}, there is an ultrafilter \(D\) with \(\lambda_D
		< \lambda\) such that there is an internal ultrapower embedding \(e :
		M_D \to M_U\) with \(e \circ j_D = j_U\) that is
		\(j_D(\lambda)\)-supercompact in \(M_D\). Since \(U\) is
		\(\lambda\)-irreducible, \(D\) is principal, and hence \(j_U= e\circ j_D
		= e\) is \(\lambda\)-supercompact as desired.
	\end{proof}
\end{cor}

\begin{cor}[UA]\label{LimitIrredThm}\index{Irreducibility Theorem}
	Suppose \(\lambda\) is a strong limit cardinal and \(U\) is a
	\(\lambda\)-irreducible ultrafilter. Then \(j_U\) is
	\({<}\lambda\)-supercompact. If \(\lambda\) is singular, then \(j_U\) is
	\(\lambda\)-supercompact. If \(\lambda\) is regular and Fr\'echet, then
	\(j_U\) is \(\lambda\)-tight. 
	\begin{proof}
		We start by showing that \(j_U\) is \({<}\lambda\)-supercompact. Fix a
		successor cardinal \(\delta < \lambda\). If \(\delta\) is Fr\'echet,
		then \(j_U\) is \(\delta\)-supercompact by \cref{SuccessorIrredCor}. If
		\(\delta\) is not Fr\'echet, then we can apply \cref{RegCompleteThm}: no
		cardinal \(\kappa\leq \delta\) is \(\delta\)-supercompact and \(U\) is
		\({\leq}2^\delta\)-irreducible, so \(U\) is \(\delta^+\)-complete and
		\(j_U\) is vacuously \(\delta\)-supercompact. Thus \(j_U\) is
		\({<}\lambda\)-supercompact.
		
		Since \(j_U\) is a \({<}\lambda\)-supercompact ultrapower embedding,
		\((M_U)^{<\lambda}\subseteq M_U\). If \(\lambda\) is singular, this
		immediately implies \((M_U)^\lambda\subseteq M_U\). Therefore \(j_U\) is
		\(\lambda\)-supercompact.
		
		If \(\lambda\) is regular and Fr\'echet, we can apply \cref{AlwaysTight}
		to conclude that \(j_U\) is \(\lambda\)-tight.
	\end{proof}
\end{cor}

As a corollary, we can finally prove the Irreducibility Theorem.
\begin{thm}[UA]\label{Irred1}\index{Irreducibility Theorem}
	Suppose \(\lambda\) is a successor cardinal or a strong limit singular
	cardinal and \(U\) is a countably complete uniform ultrafilter on
	\(\lambda\). Then the following are equivalent:
	\begin{enumerate}[(1)]
		\item \(j_U\) is \(\lambda\)-irreducible.
		\item \(j_U\) is \(\lambda\)-supercompact.
	\end{enumerate}
\end{thm}
\begin{proof}
	{\it (1) implies (2):} Follows from \cref{SuccessorIrredThm} and
	\cref{LimitIrredThm}.
	
	{\it (2) implies (1):} Follows from \cref{TightIrred}.
\end{proof}

\begin{thm}[UA]\label{Irred2}
	Suppose \(\lambda\) is an inaccessible cardinal and \(U\) is a countably
	complete ultrafilter on \(\lambda\). Then the following are equivalent:
	\begin{enumerate}[(1)]
		\item \(j_U\) is \(\lambda\)-irreducible.
		\item \(j_U\) is \({<}\lambda\)-supercompact and \(\lambda\)-tight.
	\end{enumerate}
\end{thm}
\begin{proof}
	{\it (1) implies (2):} Follows from \cref{LimitIrredThm} and
	\cref{AlwaysTight}.
	
	{\it (2) implies (1):} Follows from \cref{TightIrred}.
\end{proof}

It is sometimes easier to use a version of the Irreducibility Theorem in the
form of \cref{SuccessorIrredThm}. This follows from \cref{LimitIrredThm} using
the structure of the Rudin-Frol\'ik order (\cref{ACC}).
\begin{lma}[UA] Suppose \(U\) is a countably complete ultrafilter and
	\(\lambda\) is a cardinal. Then there is a countably complete ultrafilter
	\(D\D U\) with \(\lambda_D < \lambda\) such that \(\tr D U\) is
	\(\lambda_*\)-irreducible in \(M_D\) where \(\lambda_* = \sup
	j_D[\lambda]\).
	\begin{proof}
		 By the local ascending chain condition for the Rudin-Frol\'ik order
		 (\cref{ACC}), there is an \(\D\)-maximal \(D\D U\) such that
		 \(\lambda_D < \lambda\). Let \(i :  M_D\to M_U\) be the unique internal
		 ultrapower embedding such that \(i\circ j_D = j_U\). Then \(i\) is the
		 ultrapower of \(M_D\) by \(\tr D U\). 
		
		Suppose towards a contradiction that \(\tr D U\) is not
		\(\lambda_*\)-irreducible in \(M_D\). Fix a cardinal \(\gamma <
		\lambda\) and a countably complete ultrafilter \(Z\) of \(M_D\) on
		\(j_D(\gamma)\) such that \(Z\D \tr D U\). Then the iteration \(\langle
		D, W\rangle\) is given by an ultrafilter \(D'\) on \(\lambda_D\cdot
		\gamma\). Now \(\lambda_{D'} \leq \lambda_D\cdot \gamma < \lambda\) but
		\(D\sD D'\D U\). This contradicts the maximality of \(D\).
	\end{proof}
\end{lma}

Combining this with the Irreducibility Theorem immediately yields the following
fact:

\begin{cor}[UA]\label{SupercompactFactor}\index{Irreducibility Theorem}
	Suppose \(U\) is a countably complete ultrafilter.
	\begin{itemize}
		\item If \(\lambda\) is a Fr\'echet successor cardinal, then there is an
		ultrafilter \(D\D U\) with \(\lambda_D < \lambda\) such that the unique
		internal ultrapower embedding \(h: M_D\to M_U\) with \(h\circ j_D =
		j_U\) is \(j_D(\lambda)\)-supercompact in \(M_D\).
		\item If \(\lambda\) is a Fr\'echet inaccessible cardinal, then there is
		an ultrafilter \(D\D U\) with \(\lambda_D < \lambda\) such that the
		unique internal ultrapower embedding \(h: M_D\to M_U\) with \(h\circ j_D
		= j_U\) is \({<}\lambda\)-supercompact and \(\lambda\)-tight in \(M_D\).
		\item If \(\lambda\) is a strong limit singular cardinal, then there is
		an ultrafilter \(D\D U\) with \(\lambda_D < \lambda\) such that the
		unique internal ultrapower embedding \(h: M_D\to M_U\) with \(h\circ j_D
		= j_U\) is \(\lambda_*\)-supercompact in \(M_D\) where \(\lambda_* =
		\sup j_D[\lambda]\).\qed
	\end{itemize}
\end{cor}
\section{Resolving the identity crisis}\label{IdentitySection}\index{Identity crisis}
In this section, we characterize all strongly compact cardinals assuming UA.
This begins with an analysis of the \(\kappa\)-complete analog of \(\mathscr
K_\lambda\), denoted \(\mathscr K_\lambda^\kappa\). 
\subsection{The equivalence of strong compactness and supercompactness}
Recall that if \(\lambda\) is \(\kappa\)-Fr\'echet, then \(\mathscr
K_\lambda^\kappa\) is the \(\sE\)-least \(\kappa\)-complete uniform ultrafilter
on \(\lambda\).  Applying the Irreducibility Theorem, \cref{KetonenIrreducible2}
yields a generalization of our analysis of \(\mathscr K_\lambda\) for successor
\(\lambda\) (\cref{SuccessorThm}) to these more complete ultrafilters:

\begin{cor}[UA]\label{SuccessorThm2}
	Suppose \(\kappa < \lambda\) and \(\lambda\) is a \(\kappa^+\)-Fr\'echet
	successor cardinal. Let \(j : V\to M\) be the ultrapower of the universe by
	\(\mathscr K_\lambda^{\kappa^+}\). Then \(M^\lambda\subseteq M\). 
	\begin{proof}
		By \cref{KetonenIrreducible2}, \(\mathscr K = \mathscr
		K_\lambda^{\kappa^+}\) is irreducible. Since \(\lambda_{\mathscr K} =
		\lambda\), \(\mathscr K\) is \(\lambda\)-irreducible. Therefore by the
		Irreducibility Theorem (\cref{SuccessorIrredCor}), \(M^\lambda\subseteq
		M\).
	\end{proof}
\end{cor}

\begin{cor}[UA] Suppose \(\kappa < \lambda\) and \(\lambda\) is a
	\(\kappa^+\)-Fr\'echet successor cardinal. Then there is a
	\(\lambda\)-supercompact cardinal \(\delta\) such that \(\kappa < \delta <
	\lambda\).\qed
\end{cor}

As in the case of the first supercompact cardinal, if \(\lambda\) is strongly
inaccessible, it is not clear whether \(\mathscr K_\lambda^{\kappa^+}\)
witnesses full \(\lambda\)-supercompactness:
\begin{cor}[UA]\label{InaccessibleThm2}
	Suppose \(\kappa < \lambda\) and \(\lambda\) is a \(\kappa^+\)-Fr\'echet
	inaccessible cardinal. Let \(j : V\to M\) be the ultrapower of the universe
	by \(\mathscr K_\lambda^{\kappa^+}\). Then \(M^{<\lambda}\subseteq M\) and
	\(M\) has the \({\leq}\lambda\)-covering property. 
	\begin{proof}
		By \cref{KetonenIrreducible2}, \(\mathscr K = \mathscr
		K_\lambda^{\kappa^+}\) is irreducible. Since \(\lambda_{\mathscr K} =
		\lambda\), \(\mathscr K\) is \(\lambda\)-irreducible. Therefore by the
		Irreducibility Theorem (\cref{LimitIrredThm}), \(M^{<\lambda}\subseteq
		M\) and \(M\) has the \({\leq}\lambda\)-covering property. 
	\end{proof}
\end{cor}

Let us now analyze \(\mathscr K_\lambda^\kappa\) for general \(\kappa\):

\begin{thm}[UA]\label{HigherKetonen}
	Suppose \(\kappa \leq \lambda\) and \(\lambda\) is a \(\kappa\)-Fr\'echet
	regular cardinal. Let \(\mathscr K = \mathscr K^\kappa_\lambda\).
	\begin{enumerate}[(1)]
		\item Suppose  \(\kappa\) is not a measurable limit of
		\(\lambda\)-strongly compact cardinals. Then \(\mathscr K\) is
		irreducible.
		\item Suppose \(\kappa\) is a measurable limit of \(\lambda\)-strongly
		compact cardinals. Let \(D\) be the \(\mo\)-least normal ultrafilter on
		\(\kappa\). Then \(D\D \mathscr K\).
	\end{enumerate}
	\begin{proof}
		{\it Proof of (1):} Assume first that \(\kappa\) is not measurable.
		Since \(\mathscr K\) is \(\kappa\)-complete, it is
		\(\kappa^+\)-complete. Hence \(\mathscr K = \mathscr
		K^{\kappa^+}_\lambda\) is irreducible by \cref{KetonenIrreducible2}. 
		
		Assume instead that \(\kappa\) is not a limit of \(\lambda\)-strongly
		compact cardinals. Let \(\nu < \kappa\) be the supremum of the
		\(\lambda\)-strongly compact cardinals below \(\kappa\). Note that
		\(\lambda\) is \(\nu^+\)-Fr\'echet since \(\lambda\) is
		\(\kappa\)-Fr\'echet and \(\nu \leq \kappa\). Moreover, \(\mathscr K\)
		is \(\nu^+\)-complete since \(\nu^+\leq \kappa\). Since \( \mathscr
		K^{\nu^+}_\lambda\) is the \(\sE\)-least \(\nu^+\)-complete uniform
		ultrafilter on \(\lambda\),  \( \mathscr K^{\nu^+}_\lambda\E \mathscr
		K^\kappa_\lambda\). On the other hand, \(\mathscr K^{\nu^+}_\lambda\) is
		\(\kappa\)-complete: by \cref{SuccessorThm2} and
		\cref{InaccessibleThm2}, the completeness of \(\mathscr
		K^{\nu^+}_\lambda\) is a \(\lambda\)-strongly compact cardinal in the
		interval \((\nu,\lambda)\), and by choice of \(\nu\), the completeness
		is at least \(\kappa\). Since \(\mathscr K^{\nu^+}_\lambda\) is a
		\(\kappa\)-complete uniform ultrafilter on \(\lambda\) and \(\mathscr K
		=\mathscr K^\kappa_\lambda\) is the \(\sE\)-least such ultrafilter,
		\(\mathscr K \E \mathscr K^{\nu^+}_\lambda\). By the antisymmetry of the
		Ketonen order, \(\mathscr K = \mathscr K^{\nu^+}_\lambda\), and in
		particular \(\mathscr K\) is irreducible by \cref{KetonenIrreducible2}.
		
		{\it Proof of (2):} Let \(j : V\to M\) be the ultrapower of the universe
		by \(\mathscr K\). 
		
		We first claim that \(\kappa\) is not measurable in \(M\). Since
		\(\mathscr K\) is \(\kappa\)-complete, \(\textsc{crt}(j)\geq \kappa\).
		Therefore if \(\delta < \kappa\) is \(\lambda\)-strongly compact, then
		\(\delta\) is \(j(\lambda)\)-strongly compact in \(M\). Suppose towards
		a contradiction that \(\kappa\) is measurable in \(M\). Then \(\kappa\)
		is a measurable limit of \(j(\lambda)\)-strongly compact cardinals in
		\(M\), so \(\kappa\) is \(j(\lambda)\)-strongly compact in \(M\) by
		Menas's Theorem (\cref{MenasBig}). But by the minimality of \(\mathscr
		K_\lambda^\kappa\) (see \cref{KetonenExistence}), \(\text{cf}^M(\sup
		j[\lambda])\) is not \(\kappa\)-Fr\'echet in \(M\), contradicting that
		\(\kappa\) is \(\text{cf}^M(\sup j[\lambda])\)-strongly compact in
		\(M\). Thus our assumption was false and so \(\kappa\) is not measurable
		in \(M\).
		
		Since \(\kappa\) is measurable in \(V\) but not in \(M\), it follows
		that \(\textsc{crt}(j) \leq \kappa\), so \(\textsc{crt}(j) = \kappa\).
		Let \(D\) be the ultrafilter on \(\kappa\) derived from \(j\) using
		\(\kappa\). Since \(D\) is a normal ultrafilter and \(\kappa\) is not
		measurable in \(M_D\), \(D\) is the \(\mo\)-least ultrafilter on
		\(\kappa\) (by the linearity of the Mitchell order, \cref{UAMO}). Recall
		that our analysis of derived normal ultrafilters (\cref{MuMeasure})
		implies that either \(D\mo \mathscr K\) or \(D\D \mathscr K\). Since
		\(\kappa\) is not measurable in \(M = M_\mathscr K\), it cannot be that
		\(D\mo \mathscr K\), and therefore we can conclude that \(D\D \mathscr
		K\).
	\end{proof}
\end{thm}

It is not hard to show that in the situation of \cref{HigherKetonen} (2), in
fact \(\mathscr K^\lambda_\kappa\) is one of the ultrafilters defined in the
proof of Menas's Theorem (\cref{MenasBig}):
\[\mathscr K^\lambda_\kappa = D\text{-}\lim_{\alpha < \kappa} \mathscr
K_\lambda^{\alpha^+}\] Moreover, there is a set \(I\in D\)  such that the
sequence \(\langle \mathscr K_\lambda^{\alpha^+} : \alpha \in I\rangle\) is
discrete, which explains why \(D\D \mathscr K^\lambda_\kappa\).

We now characterize the critical point of \(\mathscr K_\lambda^\nu\).
\begin{defn}
	Suppose \(\nu\leq \lambda\) are uncountable cardinals and \(\lambda\) is
	\(\nu\)-Fr\'echet. Then \(\kappa^\nu_\lambda\) denotes the completeness of
	\(\mathscr K_\lambda^\nu\).
\end{defn}

To analyze \(\kappa^\nu_\lambda\), we use the following generalization of
\cref{KetonenFrechet}:

\begin{lma}\label{KetonenFrechet2}
	Suppose \(\nu\leq \lambda\) are cardinals and \(\lambda\) is regular.
	Suppose \(\kappa\leq\lambda\) is the least \((\nu,\lambda)\)-strongly
	compact cardinal. Suppose \(j : V\to M\) is an elementary embedding such
	that \(\textnormal{cf}^M(\sup j[\lambda])\) is not \(j(\nu)\)-Fr\'echet in
	\(M\). Then \(j\) is \((\lambda,\delta)\)-tight for some \(M\)-cardinal
	\(\delta < j(\kappa)\).
	\begin{proof}
		Since \(\kappa\) is \((\nu,\lambda)\)-strongly compact, every cardinal
		in the interval \([\kappa,\lambda]\) is \(\nu\)-Fr\'echet. Thus in
		\(M\), every cardinal in the interval \(j([\kappa,\lambda])\) is
		\(j(\nu)\)-Fr\'echet. Let \(\delta = \textnormal{cf}^M(\sup
		j[\lambda])\). By \cref{KetonenCov}, \(j\) is
		\((\lambda,\delta)\)-tight. Moreover \(\delta \leq \sup j[\lambda] \leq
		j(\lambda)\) and \(\delta\notin j([\kappa,\lambda])\) since \(\delta\)
		is not \(j(\nu)\)-Fr\'echet. Thus \(\delta < j(\kappa)\). This proves
		the lemma.
	\end{proof}
\end{lma}

The following proposition shows that under UA, all the ultrafilter-theoretic
generalizations of strong compactness collapse to a single concept:

\begin{prp}[UA]\label{KappaNuLambda}
	Suppose \(\nu \leq \kappa\leq \lambda\) are cardinals, \(\lambda\) is a
	regular cardinal, and \(\kappa\) is the least \((\nu,\lambda)\)-strongly
	compact cardinal. Then \(\kappa = \kappa^\nu_\lambda\) and \(\kappa\) is
	\(\lambda\)-strongly compact.
	\begin{proof}
		Since there is a \((\nu,\lambda)\)-strongly compact cardinal \(\kappa
		\leq \lambda\), there is some cardinal below \(\lambda\) that is
		\(\lambda\)-supercompact. Thus if \(\lambda\) is a limit cardinal then
		\(\lambda\) is strongly inaccessible by our results on GCH
		(\cref{MainTheorem}). In particular, we are in a position to apply the
		Irreducibility Theorem.
		
		By \cref{HigherKetonen}, either \(\kappa^\nu_\lambda\) is a measurable
		limit of \(\lambda\)-strongly compact cardinals or \(\mathscr
		K^\nu_\lambda\) is irreducible. In the former case
		\(\kappa^\nu_\lambda\) is \(\lambda\)-strongly compact by \cref{Menas}.
		In the latter case, \(\mathscr K^\nu_\lambda\) witnesses that
		\(\kappa^\nu_\lambda\) is \({<}\lambda\)-supercompact and
		\(\lambda\)-strongly compact by the Irreducibility Theorem
		(\cref{SuccessorIrredCor} and \cref{LimitIrredThm}).
		
		In particular, \(\kappa^\nu_\lambda\) is \((\nu,\lambda)\)-strongly
		compact, so \(\kappa \leq \kappa^\nu_\lambda\). 
		
		On the other hand, by \cref{KetonenFrechet2}, \(\kappa \leq
		\kappa^\nu_\lambda\).
		
		Thus \(\kappa = \kappa^\nu_\lambda\), and in particular \(\kappa\) is
		\(\lambda\)-strongly compact.
	\end{proof}
\end{prp}

\begin{cor}[UA]\label{LocalMenasUA}
		Suppose \(\kappa\leq \lambda\) are cardinals, \(\lambda\) is a successor
		cardinal, and \(\kappa\) is \(\lambda\)-strongly compact. Then either
		\(\kappa\) is \(\lambda\)-supercompact or \(\kappa\) is a measurable
		limit of \(\lambda\)-supercompact cardinals.
		\begin{proof}
			Assume by induction that the theorem is true for \(\bar \kappa <
			\kappa\). By \cref{KappaNuLambda}, \(\kappa=
			\kappa^\kappa_\lambda\). By \cref{HigherKetonen}, either \(\kappa\)
			is a measurable limit of \(\lambda\)-strongly compact cardinals or
			\(\mathscr K^\kappa_\lambda\) is irreducible. If \(\kappa\) is a
			limit of \(\lambda\)-strongly compact cardinals, then by our
			induction hypothesis, \(\kappa\) is a measurable limit of
			\(\lambda\)-supercompact cardinals. If instead \(\mathscr
			K_\lambda^\kappa\) is irreducible, then by \cref{SuccessorIrredThm},
			\(\mathscr K_\lambda^\kappa\) witnesses that \(\kappa\) is
			\(\lambda\)-supercompact.
		\end{proof}
\end{cor}

This implies our converse to Menas's Theorem, stating that under UA, a strongly
compact cardinal is either a supercompact cardinal or a measurable limit of
supercompact cardinals:

\begin{thm}[UA]\label{MenasUA}\index{Strongly compact cardinal!equivalence with supercompactness}
	Suppose \(\kappa\) is a strongly compact cardinal. Either \(\kappa\) is a
	supercompact cardinal or \(\kappa\) is a measurable limit of supercompact
	cardinals
\end{thm}
\begin{proof}
	Suppose \(\kappa\) is strongly compact. By the Pigeonhole Principle, there
	is a cardinal \(\gamma\geq \kappa\) such that a cardinal \(\bar \kappa\leq
	\kappa\) is supercompact if and only if \(\bar \kappa\) is
	\(\gamma\)-supercompact. Since \(\kappa\) is \(\gamma^+\)-strongly compact,
	\cref{LocalMenasUA} implies that either \(\kappa\) is
	\(\gamma^+\)-supercompact or \(\kappa\) is a limit of
	\(\gamma^+\)-supercompact cardinals. By our choice of \(\gamma\), it follows
	that either \(\kappa\) is supercompact or \(\kappa\) is a limit of
	supercompact cardinals, as desired.
\end{proof}
The use of the Pigeonhole Principle is unnecessary here, since the cardinal
\(\gamma\) turns out to equal \(\kappa\); a more careful argument appears in the
proof of \cref{ThreshThm2}.

Before generalizing our results on ultrapower thresholds (\cref{ThreshThm}), it
is worth noting that our large cardinal assumptions now put us in a local GCH
context. For example, we have the following lemma:
\begin{lma}[UA]\label{TwoSCGCH}
	Suppose \(\lambda\) is a regular Fr\'echet cardinal. Suppose \(\lambda\) is
	also \(\kappa^+_\lambda\)-Fr\'echet.\footnote{By the proof of the lemma,
	this hypothesis can be reformulated as the statement that there are distinct
	\(\lambda\)-strongly compact cardinals.} Then for all cardinals \(\gamma\in
	[\kappa_\lambda,\lambda]\), \(2^\gamma = \gamma^+\).
	\begin{proof}
		Let \(\kappa = \kappa_\lambda\). Let \(\mathscr K_0 = \mathscr
		K_\lambda\) and \(\mathscr K_1 = \mathscr K^{\kappa^+}_\lambda\). Then
		\(\mathscr K_1\) is \(\lambda\)-decomposable yet since \(\mathscr K_1\)
		is \(\kappa^+\)-complete, \(\mathscr K_0\not \D \mathscr K_1\).
		Therefore \cref{IsoFactor} implies that \(\lambda\) is not isolated. It
		follows that \(\kappa\) is \({<}\lambda\)-supercompact. In particular,
		applying our results on GCH (namely \cref{MainTheorem}), either
		\(\lambda\) is a successor cardinal or \(\lambda\) is a strongly
		inaccessible cardinal. Thus we are in a position to apply
		\cref{SuccessorThm2} and \cref{InaccessibleThm2}. 
		
		A weak consequence of the conjunction of these two theorems is that
		there is an elementary embedding \(j : V\to M\) such that
		\(\textsc{crt}(j) > \kappa\), \(j(\lambda) > \lambda^{++M}\), and \(j\)
		is \(\lambda\)-pseudocompact (or in other words, \(j\) is
		\(\gamma\)-tight for all \(\gamma \leq \lambda\)). Since \(j(\kappa) =
		\kappa\) and \(j(\lambda) > \lambda^{++M}\), \(\kappa\) is
		\(\lambda^{++M}\)-supercompact in \(M\). Thus by our results on GCH
		(\cref{MainTheorem}) applied in \(M\), \(M\) satisfies that for all
		\(\gamma\in [\kappa,\lambda]\), \(2^\gamma = \gamma^+\). But for all
		\(\gamma \leq \lambda\), the \(\gamma\)-tightness of \(j\) implies that
		\(2^\gamma \leq (2^\gamma)^M\) (by \cref{TightContinuum}), and hence
		\[2^\gamma \leq (2^\gamma)^M\leq \gamma^{+M} \leq \gamma^+\] as desired.
	\end{proof}
\end{lma}

\begin{defn}\index{Ultrapower threshold!\((\nu,\lambda\)-threshold}
	Suppose \(\nu \leq \lambda\) are uncountable cardinals. The {\it
	\((\nu,\lambda)\)-threshold} is the least ordinal \(\kappa\) such that for
	all \(\alpha < \lambda\), there is an ultrapower embedding \(j :V \to M\)
	such that \(\textsc{crt}(j) \geq \nu\) and \(j(\kappa) >\alpha\).
\end{defn}

The following theorem is proved in ZFC and has nothing to do with UA.

\begin{thm}\label{StronglyTallLocal}\index{Strongly tall cardinal}
	Suppose \(\kappa\leq\lambda\) are cardinals, \(\lambda\) is regular, and
	\(\kappa\) is the \((\nu,{\leq}\lambda^+)\)-threshold. Assume \(2^\gamma =
	\gamma^+\) for all cardinals \(\gamma\in [\kappa,\lambda]\). Then \(\kappa\)
	is \((\nu,\lambda)\)-strongly compact.
	\begin{proof}
		Let \(U\) be a \(\nu\)-complete ultrafilter such that \(j_U(\kappa) \geq
		\lambda^+\). Suppose \(\gamma\) is a regular cardinal in the interval
		\([\kappa,\lambda]\). Suppose towards a contradiction that \(U\) is
		\(\gamma\)-indecomposable and \(\gamma^+\)-indecomposable. Since
		\(2^\gamma = \gamma^+\), we can apply Silver's Theorem (\cref{Silver}).
		This yields an ultrafilter \(D\) with \(\lambda_D < \gamma\) such that
		there is an elementary embedding \(k: M_D\to M_U\) with \(k\circ j_D =
		j_U\) and \(\textsc{crt}(k) > j_D(\gamma^+)\). Since \(j_D(\kappa) \leq
		j_D(\gamma^+)\), \[ j_D(\kappa) = k(j_D(\kappa)) = j_U(\kappa) \geq
		\lambda^+\] But \(j_D(\kappa) < (\kappa^{\lambda_D})^+ \leq
		(\gamma^{<\gamma})^+ = \gamma^+\leq \lambda^+\), which is a
		contradiction.
		
		Therefore \(U\) is either \(\gamma\)-decomposable or
		\(\gamma^+\)-decomposable. But if \(U\) is \(\gamma^+\)-decomposable,
		then since \(\gamma\) is regular, in fact, \(U\) is
		\(\gamma\)-decomposable (by Prikry's Theorem \cite{Prikry}, or the proof
		of \cref{SuccessorPrp}). In particular, every regular cardinal in the
		interval \([\kappa,\lambda]\) carries a \(\nu\)-complete uniform
		ultrafilter, which implies that \(\kappa\) is \((\nu,\lambda)\)-strongly
		compact.
	\end{proof}
\end{thm}

Let us point out that this answers a question of Hamkins \cite{Hamkins} assuming
GCH. Hamkins defines a cardinal \(\kappa\) to be {\it strongly tall} if
\(\kappa\) is the \((\kappa,\text{Ord})\)-threshold, and asks about the
relationship between strongly tall and strongly compact cardinals:

\begin{thm}[GCH] If \(\kappa\) is strongly tall, then \(\kappa\) is strongly
	compact.\qed
\end{thm}

\begin{thm}[UA]\label{ThreshThmLocal}
	Suppose \(\lambda\) is a regular Fr\'echet cardinal. Suppose \(\kappa \leq
	\lambda\) is the \((\nu,{\leq}\lambda^+)\)-threshold for some \(\nu >
	\kappa_\lambda\). Then \(\kappa\) is \(\lambda\)-strongly compact.
	\begin{proof} The following is the main claim:
		\begin{clm*} \(\lambda\) is \(\nu\)-Fr\'echet.\end{clm*}
		\begin{proof}[Sketch] We first claim that there is some
		\(\nu\)-Fr\'echet cardinal in the interval \([\lambda,2^\lambda]\).
		Assume towards a contradiction that this fails.  Fix \(U\) such that
		\(j_U(\kappa) \geq \lambda^+\). By Silver's Theorem (\cref{Silver}),
		there is an ultrafilter \(D\) with \(\lambda_D < \lambda\) such that
		there is an elementary embedding \(k : M_D\to M_U\) with
		\(\textsc{crt}(k) > j_D((2^\lambda)^+)\). In particular, \(j_D(\kappa)
		\geq \lambda^+\). In particular, it follows that \(\lambda\) is not
		isolated by \cref{IsoFix}. Let \(\gamma = \lambda_D\). We claim that
		\(2^\gamma = \gamma^+\). If \(\gamma\) is singular, this follows from
		\cref{MainThm}:  note that \(\gamma\in [\kappa_\lambda,\lambda]\) so
		some cardinal is \(\gamma\)-supercompact by \cref{NonisolatedCompact},
		and hence \(2^\gamma = \gamma^+\) by \cref{MainThm}. If \(\gamma\) is
		regular, then this follows from \cref{TwoSCGCH} since by
		\cref{KappaChar}, \(\kappa_\gamma \leq \kappa_\lambda\leq \nu\). Thus
		\(2^\gamma = \gamma^+\) in either case. From this (and \cref{MainThm})
		it follows that \(\lambda^\gamma = \lambda\). This contradicts that
		\(j_D(\lambda) \geq \lambda^+\). Thus our assumption was false, so there
		is a \(\nu\)-Fr\'echet cardinal in the interval \([\lambda,2^\lambda]\).
		
		Now let \(\lambda'\) be the least \(\nu\)-Fr\'echet cardinal greater
		than or equal to \(\lambda\). Suppose towards a contradiction that
		\(\lambda' > \lambda\). 
		
		We claim \(\lambda'\) is an isolated cardinal. Clearly \(\lambda'\) is
		Fr\'echet. By the proof of \cref{SuccessorPrp}, \(\lambda'\) is a limit
		cardinal. Finally, \(\lambda'\) is not a limit of Fr\'echet cardinals:
		otherwise by \cref{LimitLma}, \(\lambda'\) is a strong limit cardinal,
		contradicting that \(\lambda < \lambda' \leq 2^{\lambda}\). Thus
		\(\lambda'\) is isolated, as claimed.
		
		 \cref{IsoFactor} implies \(\mathscr K_{\lambda'}\D \mathscr
		 K_{\lambda'}^\nu\), which implies that \(\mathscr K_{\lambda'}\) is
		 \(\nu\)-complete, or in other words \(\kappa_{\lambda'}\geq \nu\).
		 Since \(\lambda \geq \kappa_{\lambda'}\), \cref{Nonoverlapping} implies
		 \(\mathscr K_{\lambda'}\not \I \mathscr K_\lambda\). By the
		 characterization of internal ultrapower embeddings of \(M_{\mathscr
		 K_\lambda}\) (\cref{EmbeddingChar}), \(\mathscr K_{\lambda'}\) must be
		 discontinuous at \(\lambda\). But this implies \(\lambda\) is
		 \(\kappa_{\lambda'}\)-Fr\'echet, and hence \(\lambda\) is
		 \(\nu\)-Fr\'echet. This contradicts our assumption that \(\lambda'
		 >\lambda\) is the least  \(\nu\)-Fr\'echet cardinal greater than or
		 equal to \(\lambda\).
		\end{proof}
	Since \(\lambda\) is \(\nu\)-Fr\'echet and \(\nu > \kappa_\lambda\), we are
	in the situation of \cref{TwoSCGCH}. Therefore for all cardinals \(\gamma\in
	[\kappa_\lambda,\lambda]\), \(2^\lambda = \lambda^+\). This yields the
	cardinal arithmetic hypothesis of \cref{StronglyTallLocal}, so we can
	conclude that \(\kappa\) is the least \((\nu,\lambda)\)-strongly compact
	cardinal. By \cref{KappaNuLambda}, it follows that \(\kappa\) is
	\(\lambda\)-strongly compact.
	\end{proof}
\end{thm}

Of course, if one works below a strong limit cardinal, one obtains the complete
generalization of \cref{ThreshThm}:

\begin{cor}[UA]\label{ThreshThm2}
	If \(\lambda\) is a strong limit cardinal and \(\kappa < \lambda\) is the
	\((\nu,\lambda)\)-threshold, then \(\kappa\) is \(\gamma\)-strongly compact
	for all \(\gamma < \lambda\). Therefore one of the following holds:
	\begin{itemize}
		\item \(\kappa\) is \(\gamma\)-supercompact for all \(\gamma <
		\lambda\).
		\item \(\kappa\) is a measurable limit of cardinals that are
		\(\gamma\)-supercompact for all \(\gamma < \lambda\).
	\end{itemize}
	\begin{proof}
		Let \(\kappa_0\) be the \(\lambda\)-threshold. By \cref{ThreshThm},
		\(\kappa_0\) is \({<}\lambda\)-supercompact. If \(\nu \leq \kappa_0\),
		then \(\kappa_0\) is the \((\nu,\lambda)\)-threshold, so \(\kappa =
		\kappa_0\), which proves the corollary. 
		
		Therefore assume \(\nu > \kappa\). Suppose \(\delta \in
		[\kappa,\lambda]\) is  a regular cardinal. By the proof of
		\cref{ThreshThm}, \(\kappa_0 = \kappa_\delta\). Moreover \(\kappa\) is
		the \((\nu,\delta^+)\)-threshold by \cref{IdemThreshold}. Therefore we
		can apply \cref{ThreshThmLocal} to obtain that \(\kappa\) is
		\(\delta\)-strongly compact. 
		
		The final two bullet points are immediate from \cref{LocalMenasUA}.
		Suppose \(\kappa\) is not \(\delta\)-supercompact for some \(\delta <
		\lambda\). By \cref{LocalMenasUA}, \(\kappa\) is a measurable limit of
		\(\gamma\)-supercompact cardinals for all \(\gamma\in
		[\delta,\lambda)\). Now suppose \(\kappa_0 < \kappa\) is
		\(\kappa\)-supercompact. We claim \(\kappa_0\) is
		\(\gamma\)-supercompact for all \(\gamma < \lambda\). Fix \(\gamma <
		\lambda\). There is some \(\kappa_1 \in (\kappa_0,\kappa]\) that is
		\(\gamma\)-supercompact. But \(\kappa_0\) is \(\kappa_1\)-supercompact,
		so in fact, \(\kappa_0\) is \(\gamma\)-supercompact, as desired.
	\end{proof}
\end{cor}

\subsection{Level-by-level equivalence at singular cardinals}\index{Level-by-level equivalence}
A well-known theorem of Apter-Shelah \cite{Apter} shows the consistency of {\it
level-by-level equivalence} of strong compactness and supercompactness: it is
consistent with very large cardinals that for all regular \(\lambda\), a
cardinal \(\kappa\) is \(\lambda\)-strongly compact if and only if it is
\(\lambda\)-supercompact or a measurable limit of \(\lambda\)-supercompact
cardinals. (By \cref{MenasBig}, this is best possible.) We showed this is a
consequence of UA assuming \(\lambda\) is a successor cardinal; when \(\lambda\)
is inaccessible, we ran into the usual problems.

When \(\lambda\) is singular, level-by-level equivalence is in general false.
This is a consequence of the following observation:
\begin{lma}\label{SingularNonEQ}Suppose \(\kappa \leq \lambda\) are cardinals.
	\begin{itemize}
	\item If \(\textnormal{cf}(\lambda) < \kappa\), then \(\kappa\) is
	\(\lambda\)-strongly compact if and only if \(\kappa\) is
	\(\lambda^+\)-strongly compact.
	\item If \(\kappa \leq \textnormal{cf}(\lambda) < \lambda\), then \(\kappa\)
	is \(\lambda\)-strongly compact if and only if \(\kappa\) is
	\({<}\lambda\)-strongly compact.
	\end{itemize}
\end{lma}
The first bullet point shows that if level-by-level equivalence holds at
successor cardinals, it also holds at singular cardinals of small cofinality.
But by the second bullet point, it need not hold at singular cardinals of larger
cofinality:
\begin{prp}\label{DumbExample}
	Suppose \(\kappa\) is the least cardinal \(\delta\) that
	\(\beth_\delta(\delta)\)-strongly compact. Then \(\kappa\) is not
	\(\beth_\kappa(\kappa)\)-supercompact.
	\begin{proof}
		In fact, if \(\delta\) is \(\beth_\delta(\delta)\)-supercompact, then
		\(\delta\) is a limit of cardinals \(\bar\delta < \delta\) that are
		\(\beth_{\bar \delta}(\bar \delta)\)-strongly compact. To see this, let
		\(j : V\to M\) be an elementary embedding such that \(\textsc{crt}(j) =
		\delta\), \(j(\delta) > \beth_\delta(\delta)\), and
		\(M^{\beth_\delta(\delta)}\subseteq M\). Then \(\delta\) is
		\({<}\beth_\delta(\delta)\)-supercompact in \(M\). It follows from
		\cref{SingularNonEQ} that \(\delta\) is
		\(\beth_\delta(\delta)\)-strongly compact in \(M\). Therefore by the
		usual reflection argument, \(\delta\) is a limit of cardinals \(\bar
		\delta < \delta\) that are \(\beth_{\bar \delta}(\bar \delta)\)-strongly
		compact.
	\end{proof}
\end{prp}
Upon further thought, however, \cref{DumbExample} {\it does not} rule out that a
version of level-by-level equivalence that holds at singular cardinals, but
rather shows that the conventional localization of strong compactness
degenerates at singular cardinals of large cofinality. We therefore introduce an
alternate localization of strong compactness:

\begin{defn}
	A cardinal \(\kappa\) is {\it \(\lambda\)-club compact} if there is a
	\(\kappa\)-complete ultrafilter on \(P_\kappa(\lambda)\) that extends the
	closed unbounded filter.
\end{defn}

If \(\kappa\) is \(\lambda\)-supercompact, then \(\kappa\) is \(\lambda\)-club
compact: a normal fine ultrafilter always extends the closed unbounded filter.
On the other hand, if every \(\kappa\)-complete filter on \(P_\kappa(\lambda)\)
extends to a \(\kappa\)-complete ultrafilter, then in particular, the closed
unbounded filter on \(P_\kappa(\lambda)\) extends to a \(\kappa\)-complete
ultrafilter, so \(\kappa\) is \(\lambda\)-club compact.

\begin{qst}[ZFC] Suppose \(\lambda\) is a regular cardinal and \(\kappa\) is
	\(\lambda\)-strongly compact. Must \(\kappa\) be \(\lambda\)-club compact?
\end{qst}

Gitik \cite[Theorem 7]{GitikClubCompact} answers this question positively under the assumption that 
\(2^\lambda = \lambda^+\). (The statement of Gitik's theorem omits the hypothesis 
that \(\lambda\) is regular, which is presumably a typo.)

To state stronger results, we introduce the Bagaria-Magidor versions of club
compactness as well:

\begin{defn}
	A cardinal \(\kappa\) is {\it \((\nu,\lambda)\)-club compact} if there is a
	\(\nu\)-complete ultrafilter on \(P_\kappa(\lambda)\) that extends the
	closed unbounded filter, and \(\kappa\) is {\it almost \(\lambda\)-club compact}
	if \(\kappa\) is \((\nu,\lambda)\)-club compact for all \(\nu < \kappa\).
\end{defn}

As is typical in the Bagaria-Magidor notation, if \(\kappa\) is
\((\nu,\lambda)\)-club compact, then every cardinal greater than \(\kappa\) is
\((\nu,\lambda)\)-club compact.

Menas's Theorem (\cref{MenasBig}) carries over to club compactness:

\begin{lma}\label{MenasClub}
	Suppose \(\lambda\) is a cardinal. Any limit of \(\lambda\)-club compact
	cardinals is almost \(\lambda\)-club compact. An almost \(\lambda\)-club
	compact cardinal is \(\lambda\)-club compact if and only if it is
	measurable. Thus every measurable limit of \(\lambda\)-club compact
	cardinals is \(\lambda\)-club compact.\qed
\end{lma}

The main theorem of this section is that under UA, level-by-level equivalence
holds for club compactness at singular cardinals.

\begin{thm}[UA]\label{ClubEQ}
	Suppose \(\kappa \leq \lambda\) are cardinals and \(\lambda\) is singular.
	Then the following are equivalent:
	\begin{enumerate}[(1)]
		\item \(\kappa\) is \(\lambda\)-club compact.
		\item \(\kappa\) is the least \((\nu,\lambda)\)-club compact cardinal
		for some \(\nu \leq \kappa\).
		\item \(\kappa\) is \(\lambda\)-supercompact or a measurable limit of
		\(\lambda\)-supercompact cardinals.
	\end{enumerate}
\end{thm}

For the proof, we use the following much more general lemma:
\begin{defn}
	The {\it Kat\v etov order} is defined on filters \(F\) and \(G\) by setting
	\(F\leq_{\text{Kat}} G\) if there is a function \(f\) on a set in \(G\) such
	that \(F\subseteq f_*(G)\).
\end{defn}
Thus \(F\leq_{\text{Kat}} G\) if and only if there is an extension \(F'\) of
\(F\) below \(G\) in the Rudin-Keisler order.

\begin{lma}[UA] Suppose \(\nu < \lambda\) are cardinals. Suppose \(\mathcal F\)
	is a normal fine filter on a set \(Y\) such that \(\lambda\subseteq
	Y\subseteq P(\lambda)\). Suppose \(A\) is a set of ordinals and \(U\) is the
	\(\sE\)-least \(\nu^+\)-complete ultrafilter on \(A\) such that \(\mathcal
	F\leq_{\textnormal{Kat}} U\). Then \(U\) is \(\lambda\)-irreducible.
	\begin{proof}
		Suppose \(D\D U\) and \(\lambda_D < \lambda\). We must show that \(D\)
		is principal. To do this, we will show that \(j_D(U) \E \tr D U\) in
		\(M_D\). By \cref{Pushdown}, it then follows that \(D\) is principal. As
		usual, to show \(j_D(U) \E \tr D U\) in \(M_D\), we verify that the
		properties for which \(U\) is minimal hold for \(\tr D U\) with the
		parameters shifted by \(j_D\). In other words, we show that \(M_D\)
		satisfies the following:
		\begin{itemize}
			\item \(\tr D U\) is a \(j_D(\nu^+)\)-complete ultrafilter on
			\(j_D(A)\).
			\item \(j_D(\mathcal F) \leq_{\textnormal{Kat}}\tr D U\).
		\end{itemize}
		The first bullet point is rather easy. By definition, \(\tr D U\) is an
		ultrafilter on \(j_D(A)\). Moreover, \(\tr D U\) is
		\(j_D(\nu^+)\)-complete in \(M_D\) since \[\textsc{crt}(k) \geq
		\textsc{crt}(j) > \nu = j(\nu) \geq j_D(\nu)\]
		
		The second bullet point is a bit more subtle. Since \(\mathcal
		F\leq_{\text{Kat}} U\), there is some \(B\in M_U\) such that \(\mathcal
		F\) is contained in the ultrafilter derived from \(j_U\) using \(B\). In
		other words, for all \(S\in \mathcal F\), \(B\in j_U(S)\). Note that for
		any \(f: \lambda\to \lambda\), \(B\) is closed under \(j_U(f)\): by
		normality, \(\{\sigma\in Y: \sigma \text{ is closed under }f\}\in
		\mathcal F\), and hence \(B\in j_U(\{\sigma\in Y: \sigma \text{ is
		closed under }f\})\), or in other words, \(B\) is closed under
		\(j_U(f)\). We will use this fact in an application of \cref{ClubLemma}.
		
		Let \(k : M_D\to M\) be the unique internal ultrapower embedding with
		\(k\circ j_D = j_U\). Thus \(k\) is the ultrapower of \(M_D\) by \(\tr D
		U\). Let \[\mathcal W= \{S\in j_D(P(Y)) : B\in k(S)\}\] Thus \(\mathcal
		W\) is the \(M_D\)-ultrafilter on \(j_D(Y)\) derived from \(k\) using
		\(B\). In particular, \(\mathcal W\RK \tr D U\) by the characterization
		of the Rudin-Keisler order in terms of derived embeddings
		(\cref{RKChar}). We claim that \(j_D(\mathcal F)\subseteq \mathcal W\).
		Clearly \(j_D[\mathcal F]\subseteq \mathcal W\). The key point is that
		by \cref{ClubLemma}, \(k(\id_D)\in B\). In other words, \[\{\sigma\in
		j_D(Y): \id_D\in \sigma\}\in \mathcal W\] Therefore by our unique
		extension lemma for normal filters (\cref{NormalGeneration}),
		\(j_D(\mathcal F)\subseteq\mathcal W\), as desired.
		
		Now \(j_D(\mathcal F)\subseteq \mathcal W\RK \tr D U\), or in other
		words \(j_D(\mathcal F)\leq_{\text{Kat}}\tr D U\).
	\end{proof}
\end{lma}

\begin{proof}[Proof of \cref{ClubEQ}] {\it (1) implies (2):} Trivial.
	
	{\it (2) implies (3):} Clearly \(\lambda\) is a limit of Fr\'echet
	cardinals, so by \cref{LimitLma}, \(\lambda\) is a strong limit cardinal.
	
	We first handle the case in which there is some \(\nu < \kappa\) such that
	\(\kappa\) is the least \((\nu,\lambda)\)-club compact cardinal. Note that
	\(\nu\) is either not measurable or not almost \(\lambda\)-club compact,
	since otherwise \(\nu\) would be the least \((\nu,\lambda)\)-club compact
	cardinal. If \(\nu\) is not almost \(\lambda\)-club compact, then there is
	some \(\bar \nu < \nu\) such that \(\kappa\) is the least \((\bar
	\nu^+,\lambda)\)-club compact cardinal. If \(\nu\) is not measurable, then
	\(\kappa\) is the least \((\nu^+,\lambda)\)-club compact cardinal. In either
	case, we can fix \(\eta < \kappa\) such that \(\kappa\) is the least
	\((\eta^+,\lambda)\)-club compact cardinal. 
	
	Let \(\mathcal F\) be the closed unbounded filter on \(P_\kappa(\lambda)\).
	Let \(U\) be the least \(\eta^+\)-complete ultrafilter on an ordinal such
	that \(\mathcal F\leq_\text{Kat} U\). Then \(U\) is \(\lambda\)-irreducible.
	Since \(\lambda\) is a singular strong limit cardinal, by
	\cref{LimitIrredThm}, \((M_U)^\lambda\subseteq M_U\). Thus
	\(\textsc{crt}(j_U)\) is \(\lambda\)-supercompact. Note that
	\(\textsc{crt}(j_U)\leq \kappa\) since \(\mathcal F\leq_{\text{Kat}} U\) and
	\(\mathcal F\) is not \(\kappa^+\)-complete. On the other hand
	\(\textsc{crt}(j_U) > \eta\), so \(\textsc{crt}(j_U)\) is an
	\((\eta^+,\lambda)\)-club compact cardinal, and hence \(\textsc{crt}(j_U)
	\leq \kappa\). Thus \(\kappa = \textsc{crt}(j_U)\) is
	\(\lambda\)-supercompact.
	
	We now handle the case in which \(\kappa\) is \((\kappa,\lambda)\)-club
	compact but there is no \(\nu < \kappa\) such that \(\kappa\) is the least
	\((\nu,\lambda)\)-club compact cardinal. Since \(\kappa\) is
	\((\nu,\lambda)\)-club compact for all \(\nu < \kappa\), it follows that for
	each \(\nu < \lambda\), the least \((\nu,\lambda)\)-club compact cardinal
	lies strictly below \(\kappa\). Thus by the previous case, \(\kappa\) is a
	limit of \(\lambda\)-supercompact cardinals. Moreover, \(\kappa\) is
	measurable since \(\kappa\) is \((\kappa,\lambda)\)-club compact. Thus
	\(\kappa\) is a measurable limit of \(\lambda\)-club compact cardinals, as
	desired.
	
	{\it (3) implies (1):} This follows from \cref{MenasClub}.
\end{proof}
\subsection{The Mitchell order, the internal relation, and coherence}\label{InternalSectionII}
Assume UA and suppose \(U\) is a normal ultrafilter on \(\kappa\). Can
\(P(\kappa^+)\subseteq M_U\)? The question remains open in general, but the
following theorem shows that if \(\kappa^+\) is Fr\'echet, this cannot occur: 
\begin{thm}[UA]\label{StrongToSuper}
	Suppose \(\lambda\) is a Fr\'echet cardinal. Suppose \(U\) is a countably
	complete ultrafilter such that \(P(\lambda)\subseteq M_U\). Then
	\((M_U)^\lambda\subseteq M_U\). 
	\begin{proof}
		Assume by induction that the theorem holds for cardinals below
		\(\lambda\).  If \(\lambda\) is a limit of Fr\'echet cardinals, we then
		have \((M_U)^{<\lambda}\subseteq M_U\). In particular, if \(\lambda\) is
		a singular limit of Fr\'echet cardinals, then \((M_U)^\lambda\subseteq
		M_U\). Thus we may assume that \(\lambda\) is either regular or
		isolated. This puts the analysis of \(\mathscr K_\lambda\) (especially
		\cref{EmbeddingChar} and \cref{IsolatedInternal}) at our disposal.
		
		We first show that \(U\) is \(\lambda\)-irreducible. Suppose towards a
		contradiction that there is a uniform ultrafilter \(D\D U\) on an
		infinite cardinal \(\gamma < \lambda\). Since \(M_U\subseteq M_D\), so
		in particular \(P(\lambda)\subseteq M_D\). A general bound on the
		strength of ultrapowers (\cref{AmenabilityLemma}) implies that \[\lambda
		< j_D(\gamma)\] 
		
		Assume first that \(\lambda\) is isolated. By \cref{IsolatedInternal},
		\(D\I \mathscr K_\lambda\), and by \cref{DeltaLemma},
		\(P(\gamma)\subseteq M_{\mathscr K_\lambda}\). Thus
		\[P(\lambda)\subseteq j_D(P(\lambda_D))\subseteq M_{\mathscr
		K_\lambda}\] Therefore by our bound on the strength of \(j_{\mathscr
		K_\lambda}\) for nonmeasurable isolated cardinals \(\lambda\)
		(\cref{IsoStrength}), \(\lambda\) is measurable. Since \(\lambda\) is a
		strong limit, \(D\in H(\lambda)\subseteq M_D\), and this is a
		contradiction.
		
		Assume instead that \(\lambda\) is a nonisolated regular cardinal. We
		use an argument similar to the one from the local proof of GCH
		(\cref{MainThm}). Let \(M = M_{\mathscr K_\lambda}\) and let \(N =
		(M_D)^{M}\). Consider the embedding \(j_{\mathscr K_\lambda}^{N}\circ
		j^M_D\). (Note: \(j_{\mathscr K_\lambda}^N\) denotes the ultrapower
		formed by using functions in \(N\) modulo the \(N\)-ultrafilter
		\(\mathscr K_\lambda\), {\it not} the ultrafilter \((\mathscr
		K_\lambda)^N\), which we have not proved to exist.) This is an
		ultrapower embedding from \(M\), and we claim that it is internal to
		\(M\). By our analysis of internal ultrapower embeddings of \(M\)
		(\cref{EmbeddingChar}), it suffices to show that \(j_{\mathscr
		K_\lambda}^{N}\circ j^M_D\) is continuous at \(\text{cf}^M(\sup
		j[\lambda]) = \lambda\). (To compute the cofinality of \(\sup
		j[\lambda]\) in \(M\), we use \cref{GeneralThm}.) Clearly
		\(j_D^M(\lambda) = \sup j_D^M[\lambda]\) since \(\lambda\) is regular
		and \(D\) lies on \(\gamma < \lambda\). Moreover \(j_D^M(\lambda)\) is
		regular in \(N\) and is larger than \(\lambda\) since \(j_D^M(\gamma) =
		j_D(\gamma) > \lambda\). Thus \(j_{\mathscr
		K_\lambda}^{N}(j^M_D(\lambda)) = \sup j_{\mathscr
		K_\lambda}^N[j^M_D(\lambda)]\). Putting it all together, \[j_{\mathscr
		K_\lambda}^N\circ j_D^M(\lambda) = \sup j_{\mathscr
		K_\lambda}^N[j_D^M(\lambda)] = \sup  j_{\mathscr K_\lambda}^N[\sup
		j_D^M[\lambda]] = \sup j_{\mathscr K_\lambda}^N\circ j_D^M[\lambda]\]
		Thus \(j_{\mathscr K_\lambda}^N\circ j_D^M\) is an internal ultrapower
		embedding of \(M\).
		
		In fact, \(j_{\mathscr K_\lambda}^N\) itself is definable over \(M\):
		for any \(f\in M^\gamma\),
		\[j_{\mathscr K_\lambda}^N([f]_D^M) = j_{\mathscr K_\lambda}^N\circ
		j_D^M(f)(\id_{\mathscr K_\lambda}^N)\] Thus \(j_{\mathscr K_\lambda}^N\)
		is definable over \(M\). Since \(P(\lambda)\subseteq N\), we have
		\(\mathscr K_\lambda = \{A\subseteq \lambda : \id_{\mathscr
		K_\lambda}^N\in j_{\mathscr K_\lambda}^N(A)\}\). Thus \(\mathscr
		K_\lambda\) is definable over \(M\), and it follows that \(\mathscr
		K_\lambda\in M\), or in other words, \(\mathscr K_\lambda\mo \mathscr
		K_\lambda\). This is a contradiction. 
		
		Thus our assumption was false, and in fact \(U\) is
		\(\lambda\)-irreducible. 
		
		To finish the proof, we break once again into cases.
		
		Suppose first that \(\lambda\) is a nonmeasurable isolated cardinal. We
		will show that \(U\) is \(\lambda^+\)-complete. We claim that \(\mathscr
		K_\lambda \not \D U\): otherwise, \(P(\lambda)\subseteq M_U\subseteq
		M_{\mathscr K_\lambda}\), and hence \(\mathscr K_\lambda\) is
		\(\lambda\)-complete by \cref{IsoStrength}, contradicting that
		\(\lambda\) is not measurable. Since \(\mathscr K_\lambda\not \D U\),
		our factorization theorem for isolated cardinals (\cref{IsoFactor})
		implies that \(U\) is \(\lambda^+\)-irreducible. Therefore by
		\cref{IsoComplete}, \(U\) is \(\lambda^+\)-complete, as claimed.
		
		If \(\lambda\) is not a nonmeasurable isolated cardinal, then
		\(\lambda\) is either a Fr\'echet successor cardinal or a Fr\'echet
		inaccessible cardinal. Since \(U\) is \(\lambda\)-irreducible, the
		Irreducibility Theorem (\cref{SuccessorIrredCor} and
		\cref{LimitIrredThm}) implies that \(j_U[\lambda]\) is contained in a
		set \(A\in M_U\) such that \(|A|^{M_U} = \lambda\). Since
		\(P(\lambda)\subseteq M_U\) and \(|A|^{M_U} = \lambda\), in fact
		\(P(A)\subseteq M_U\). In particular, the subset \(j_U[\lambda]\subseteq
		A\) belongs to \(M_U\), so \(j_U\) is \(\lambda\)-supercompact, and
		hence \((M_U)^\lambda\subseteq M_U\).
	\end{proof}
\end{thm}

A consequence of the coincidence of strength and supercompactness at Fr\'echet
cardinals is that under UA, the generalized Mitchell order is very well-behaved.

\begin{thm}[UA]\label{UltrapowerCorrectness}
	Suppose \(U\) and \(W\) are countably complete ultrafilters such that \(U\mo
	W\). Then \((j_U)^{M_W} = j_U\restriction M_W\). In fact,
	\((M_W)^{\lambda_U}\subseteq M_W\).
	\begin{proof}
		Let \(\lambda= \lambda_U\). Fix \(A\in U\) with \(|A| = \lambda\). Since
		\(U\in M_W\), \(P(A)\subseteq M_W\), and hence \(P(\lambda)\subseteq
		M_W\). Since \(\lambda = \lambda_U\), \(\lambda\) is Fr\'echet. Hence
		\((M_W)^{\lambda}\subseteq M_W\) by \cref{StrongToSuper}. By
		\cref{USupercompact}, this implies \((j_U)^{M_W} = j_U\restriction
		M_W\). 
	\end{proof}
\end{thm}

As a consequence, UA implies that the internal relation and the seed order
extend the Mitchell order:

\begin{cor}[UA]\label{MOExtensions}\index{Generalized Mitchell order!vs. the internal relation}\index{Ketonen order!vs. the generalized Mitchell order}\index{Internal relation!vs. the generalized Mitchell order}
	Suppose \(U\) and \(W\) are countably complete ultrafilters such that \(U\mo
	W\). Then \(U\I W\). Assume moreover that \(\lambda_U\) is the underlying
	set of \(U\) and \(W\) concentrates on ordinals. Then \(U\swo W\).
	\begin{proof}
		By \cref{UltrapowerCorrectness}, \(U\I W\). Moreover, \(j_{W}\) is
		\(\lambda_U\)-supercompact, so by \cref{UFSuperBound}, \(\lambda_U \leq
		\lambda_W\). Thus if \(\lambda_U\) is the underlying set of \(U\) and
		\(W\) concentrates on ordinals, then \[\delta_U = \lambda_U \leq
		\lambda_W\leq \delta_W\] Therefore by \cref{InternalKet}, we have
		\(U\swo W\).
	\end{proof}
\end{cor}

Using the Irreducibility Theorem, we prove some converses of \cref{MOExtensions}
that demystify the internal relation. This requires an argument we have seen
before but which we now make explicit:
\begin{lma}\label{InternalToSuper}
	Suppose \(W\) is a countably complete ultrafilter such that \(j_W\) is
	\({<}\lambda\)-strong and \(\lambda\)-tight.\footnote{Equivalently, \(j_W\)
	is \({<}\lambda\)-supercompact and \(\lambda\)-tight.} Suppose that there is
	a countably complete ultrafilter \(U\) on \(\lambda\) such that \(U\I W\)
	and \(\sup j_U[\lambda] < j_U(\lambda)\). Then \(j_W\) is
	\(\lambda\)-supercompact.
	\begin{proof}
		We first show that \(P(\lambda)\subseteq M_W\). Since \(W\) is
		\({<}\lambda\)-strong, \(P(\alpha)\subseteq M_W\) for all \(\alpha <
		\lambda\). Therefore by the elementarity of \(j_U\), \(M_U\) satisfies
		that \(P(\sup j_U[\lambda])\subseteq j_U(M_W)\). In other words,
		\(P^{M_U}(\sup j_U[\lambda])\subseteq j_W(M_U)\). Since \(U\I W\),
		\(j_U(M_W)\subseteq M_W\), and therefore \(P^{M_U}(\sup
		j_U[\lambda])\subseteq M_W\). Now fix \(A\subseteq \lambda\). We have
		\(j_U(A)\cap \sup j_U[\lambda] \in P^{M_U}(\sup j_U[\lambda])\subseteq
		M_W\). Moreover \(j_U\restriction \lambda\in M_W\) since \(U\I W\).
		Hence \[A = j_U^{-1}[j_U(A)\cap \sup j_U[\lambda]]\in M_W\] This shows
		that \(P(\lambda)\subseteq M_W\), as claimed.
		
		Now suppose \(B\) is a subset of \(M_W\) of cardinality at most
		\(\lambda\). Since \(j_W\) is \(\lambda\)-tight, there is a set \(C\in
		M_W\) of \(M_W\)-cardinality at most \(\lambda\) such that \(B\subseteq
		C\). Since \(P(\lambda)\subseteq M_W\) and \(|C|^{M_W}\leq \lambda\),
		\(P(C)\subseteq M_W\). Thus \(B\in M_W\). It follows that \(j_W\) is
		\(\lambda\)-supercompact.
	\end{proof} 
\end{lma}

\begin{thm}[UA]\label{XIrredInternal}
	Suppose \(W\) is a countably complete ultrafilter and \(U\) is a countably
	complete uniform ultrafilter on a set \(X\subseteq M_W\). Then the following
	are equivalent:
	\begin{enumerate}[(1)]
		\item \(U\mo W\).
		\item \(U\I W\) and \(W\) is \(|X|\)-irreducible.
	\end{enumerate}
	\begin{proof}
		Let \(\lambda = \lambda_U = |X|\).
		
		{\it (1) implies (2):} Suppose \(U\mo W\). Then \(j_W\) is
		\(\lambda\)-supercompact by \cref{StrongToSuper}, so \(W\) is
		\(\lambda\)-irreducible by \cref{TightIrred}. Moreover by
		\cref{MOExtensions}, \(U\I W\). This shows that (2) holds.
		
		{\it (2) implies (1):} Suppose \(U\I W\) and \(W\) is
		\(\lambda\)-irreducible.
		
		Suppose first that \(\lambda\) is an isolated cardinal. We claim that
		\(W\) is \(\lambda^+\)-complete. Note that \(j_W\) must be continuous at
		\(\lambda\) by \cref{Nonisolation1}. It follows that \(W\) is
		\(\lambda^+\)-irreducible. Hence \(W\) is
		\(\lambda^\sigma\)-irreducible. But \(\lambda^\sigma\) is measurable (by
		\cref{DoubleSigma}), so by \cref{RegCompleteThm} it follows that \(W\)
		is \(\lambda^+\)-complete. As an immediate consequence, \(U\mo W\).
		
		Suppose instead that \(\lambda\) is not isolated. We can then apply the
		Irreducibility Theorem (\cref{SuccessorIrredCor} and
		\cref{LimitIrredThm}) to conclude that \(W\) is
		\({<}\lambda\)-supercompact and \(\lambda\)-tight. Since \(U\I W\),
		\cref{InternalToSuper} yields that \(j_W\) is \(\lambda\)-supercompact.
		In particular, \(P(\lambda)\subseteq M_W\), so \(U\mo W\), as desired.
	\end{proof}
\end{thm}

We can reformulate \cref{XIrredInternal} slightly to characterize the internal
relation in terms of the Mitchell order:

\begin{thm}[UA]\label{InternalChar}
	Suppose \(U\) and \(W\) are hereditarily uniform irreducible ultrafilters.
	Then the following are equivalent:
	\begin{enumerate}[(1)]
		\item \(U\I W\).
		\item Either \(U\mo W\) or \(W \in V_\kappa\) where \(\kappa =
		\textsc{crt}(j_U)\).
	\end{enumerate}
\end{thm}
For this, we need the following theorem, which shows that the notions of
\(\lambda\)-irreducible, \(\lambda\)-Mitchell, and \(\lambda\)-internal
ultrafilters (\cref{LambdaIrredDef}, \cref{LambdaMitchellDef},
\cref{LambdaInternalDef} respectively) coincide under UA:
\begin{thm}[UA]\label{LambdaIrredChar}\index{Ultrafilter!\(\lambda\)-internal}\index{Ultrafilter!\(\lambda\)-irreducible}\index{Ultrafilter!\(\lambda\)-Mitchell}
	Suppose \(U\) is an ultrafilter and \(\lambda\) is a cardinal. Then the
	following are equivalent:
	\begin{enumerate}[(1)]
		\item \(U\) is \(\lambda\)-irreducible.
		\item \(U\) is \(\lambda\)-Mitchell.
		\item \(U\) is \(\lambda\)-internal.
	\end{enumerate}

	\begin{proof}
		{\it (1) implies (2):} Assume \(U\) is \(\lambda\)-irreducible. We may
		assume by induction that for all \(U'\E U\) and \(\lambda'\leq \lambda\)
		with \(U'\sE U\) or \(\lambda' < \lambda\),  if \(U'\) is \(
		\lambda'\)-irreducible then \(U'\) is \( \lambda'\)-Mitchell. Thus \(U\)
		is \(\lambda'\)-Mitchell for all \(\lambda' < \lambda\). In particular,
		\(U\) is automatically \(\lambda\)-Mitchell unless \(\lambda\) is a
		successor cardinal and and the cardinal predecessor \(\gamma\) of
		\(\lambda\) is Fr\'echet. Therefore we can assume \(\lambda = \gamma^+\)
		where \(\gamma\) is a Fr\'echet cardinal.
		
		We may also assume that \(\gamma^\sigma\) exists, since otherwise the
		\(\lambda\)-irreducibility of \(U\) implies \(U\) is principal, so (2)
		holds automatically. Let \(\eta = \gamma^\sigma\). 
		
		Assume first that \(\eta = \gamma^+\). Then \(\gamma^+\) is Fr\'echet,
		so by the Irreducibility Theorem (\cref{SuccessorIrredCor}), \(U\) is
		\(\gamma^+\)-supercompact. Therefore every countably complete
		ultrafilter on \(\gamma\) belongs to \(M_U\) by \cref{MitchellLemma}. In
		other words, \(U\) is \(\gamma^+\)-Mitchell.
		
		This leaves us with the case that \(\eta > \gamma^+\). In other words,
		by \cref{SuccessorPrp}, \(\eta\) is isolated.
		
		Assume first that \(\mathscr K_\eta\not \D U\). Then by
		\cref{IsoFactor}, \(U\) is \(\eta\)-indecomposable, and so in particular
		\(U\) is \(\eta^+\)-irreducible. By \cref{IsoComplete} (3), \(U\) is
		\(\eta^+\)-complete, which easily implies that \(U\) is
		\(\gamma^+\)-Mitchell.
		
		Assume finally that \(\mathscr K_\eta\D U\). Let \(j : V\to M\) be the
		ultrapower of the universe by \(\mathscr K_\eta\). Let \(h : M\to M_U\)
		be the unique internal ultrapower embedding with \(h\circ j = j_U\). 
		
		Recall that \(\tr {\mathscr K_\eta} U\) is the canonical ultrafilter
		\(Z\) of \(M\) such that \(j_Z^M = h\). We claim that \(\tr {\mathscr
		K_\eta} U\) is \(\gamma^+\)-irreducible in \(M\).  Suppose \(M\)
		satisfies that \(D\) is an ultrafilter on \(\gamma\) with \(D\D \tr
		{\mathscr K_\eta} U\). Let \(i : (M_D)^M \to M_U\) be the unique
		internal ultrapower embedding such that \[i\circ j_D^M = h\] We will
		show \(D\) is principal by showing that \(D\D U\). By \cref{DeltaLemma},
		\(M\) is closed under \(\gamma\)-sequences. In particular,
		\(P(\gamma)\subseteq M\), so \(D\) really is an ultrafilter on
		\(\gamma\), and hence the question of whether \(D\D U\) makes sense.
		Moreover \(j_D \restriction M = j_D^M\), and so \(j_D^M\circ j =
		j_D(j)\circ j_D\). Now \[i\circ j_D(j)\circ j_D = i \circ j_D^M\circ j =
		h\circ j = j_U\] Thus \(i\circ j_D(j) : M_D\to M_U\) is an internal
		ultrapower embedding witnessing \(D\D U\). It follows that \(D\) is
		principal since \(U\) is \(\gamma^+\)-irreducible.
		
		Thus \(\tr {\mathscr K_\eta} U\) is \(\gamma^+\)-irreducible in \(M\).
		Moreover by \cref{Pushdown}, \(\tr {\mathscr K_\eta} U\sE j(U)\) in
		\(M\). Our induction hypothesis yields that for all \(U'\sE U\) and all
		\( \lambda '\leq \gamma^+\), if \(U'\) is \( \lambda '\)-irreducible
		then \(U'\) is \( \lambda '\)-Mitchell. Shifting this hypothesis by the
		elementary embedding \(j: V\to M\), we have that for all \(U'\sE j(U)\)
		and all \( \lambda '\leq j(\gamma^+)\), if \(U'\) is \( \lambda
		'\)-irreducible in \(M\) then \(U'\) is \( \lambda '\)-Mitchell in
		\(M\). Applying this with \(U' = \tr {\mathscr K_\lambda} U\) and \(
		\lambda ' = \gamma^+\), it follows that \(\tr {\mathscr K_\lambda} U\)
		is \(\gamma^+\)-Mitchell in \(M\). Thus every countably complete
		ultrafilter of \(M\) on \(\gamma\) belongs to \((M_{\tr {\mathscr
		K_\eta} U})^M = M_U\). But by \cref{IsolatedInternal} and
		\cref{DeltaLemma}, every countably complete ultrafilter on \(\gamma\)
		belongs to \(M\). Hence every countably complete ultrafilter on
		\(\gamma\) belongs to \(M_U\). In other words, \(U\) is
		\(\gamma^+\)-Mitchell as desired.
		
		{\it (2) implies (3):} Immediate from \cref{MOExtensions}.
		
		{\it (3) implies (4):} Assume \(U\) is \(\lambda\)-internal. Suppose
		\(D\D U\) and \(\lambda_D < \lambda\). We will show \(D\) is principal.
		Since \(\lambda_D < \lambda\), \(D\I U\). Thus \(D\D U \gI D\), so \(D\I
		D\) by \cref{RFInternal}. Since the internal relation is irreflexive on
		nonprincipal ultrafilters, \(D\) is principal.
\end{proof}
\end{thm}

\begin{proof}[Proof of \cref{InternalChar}] {\it (1) implies (2):} Suppose \(U\I
	W\). 
	
	Assume first that \(\lambda_U\leq \lambda_W\). Then since \(W\) is
	irreducible, \(W\) is \(\lambda_U\)-irreducible. By \cref{XIrredInternal},
	\(U\mo W\). 
	
	Assume instead that \(\lambda_W < \lambda_U\). Then by
	\cref{LambdaIrredChar}, \(W\I U\). Since \(U\I W\) and \(W \I U\),
	\cref{Commute} implies that \(U\) and \(W\) are commuting ultrafilters in
	the sense of Kunen's commuting ultrapowers lemma (\cref{KunenCommute}).
	Moreover, again by \cref{LambdaIrredChar}, \(U\) is \(\lambda_U\)-internal
	and \(W\) is \(\lambda_W\)-internal. We can therefore apply our converse to
	Kunen's commuting ultrapowers lemma, from which it follows that \(W \in
	V_\kappa\) where \(\kappa = \textsc{crt}(j_U)\).
	
	{\it (2) implies (1):} If \(U\mo W\), then \(U\I W\) by \cref{MOExtensions}.
	If \(W\in V_\kappa\) where \(\kappa = \textsc{crt}(j_U)\), then \(U\I W\) by
	Kunen's commuting ultrapowers lemma (\cref{KunenCommute}).
\end{proof}

We now reformulate UA in terms of a form of coherence:

\begin{defn}
	Suppose \(C\) is a class of countably complete ultrafilters. 
	\begin{itemize}
		\item Suppose \(\mathcal I = \langle M_n,j_{nm}, U_n : n< m
		\leq\ell\rangle\) is a finite iterated ultrapower.
		\begin{itemize}
			\item A countably complete ultrafilter \(U\) is {\it given by
			\(\mathcal I\)} if \(j_U = j_{0\ell}\).
			\item \(\mathcal I\) is a {\it \(C\)-iteration} if \(U_n\in
			j_{0n}(C)\) for all \(n < \ell\).
		\end{itemize}
		\item \(C\) is {\it cofinal} if the class of ultrafilters given by
		\(C\)-iterations is Rudin-Frol\'ik cofinal. 
		\item \(C\) is {\it coherent} if for any distinct ultrafilters \(U\) and
		\(W\) of \(C\),  either \(U\in j_W(C)\) and \((M_W)^{\lambda_U}\subseteq
		M_W\), or \(W\in j_U(C)\) and \((M_U)^{\lambda_W}\subseteq M_U\).

	\end{itemize}
\end{defn}

\begin{thm}\label{SOtoUA}\index{Ultrapower Axiom!and coherence}
The following are equivalent:
\begin{enumerate}[(1)]
	\item There is a coherent cofinal class of countably complete ultrafilters.
	\item The Ultrapower Axiom holds.
\end{enumerate}
\end{thm}

For one direction of the theorem, we show that under UA, there is a canonical
coherent cofinal class of ultrafilters:
\begin{defn}\index{Mitchell point}
	An ultrafilter \(D\) is a {\it Mitchell point} if for all uniform countably
	complete ultrafilters \(U\), if \(U\sE D\), then \(U\mo D\).
\end{defn}
Dodd sound ultrafilters are Mitchell points by \cref{LipMO}. Under UA, isonormal
ultrafilters are Mitchell points by \cref{IsoIPoint}. The following fact is
trivial:
\begin{lma}[UA]\label{MPSO}
	The Mitchell points form a coherent class of ultrafilters.
	\begin{proof}
		Let \(C\) be the class of Mitchell points. Since the Ketonen order is
		linear, \(C\) is linearly ordered by \(\sE\), and hence by the
		definition of a Mitchell point, \(C\) is linearly ordered by the
		Mitchell order. The property of being a Mitchell point is absolute, so
		if \(U\mo W\) are Mitchell points, then \(U\in j_W(C)\). Moreover
		\cref{UltrapowerCorrectness}, \((M_W)^{\lambda_U}\subseteq M_W\). Thus
		\(C\) is coherent.
	\end{proof}
\end{lma}

We next show that under UA, the Mitchell points form a cofinal class. The first
step is to give an alternate characterization in terms of the internal relation:

\begin{prp}[UA]\label{MitchellPoint}
	Suppose \(U\) is a nonprincipal countably complete tail uniform ultrafilter
	on an ordinal \(\delta\). The following are equivalent:
	\begin{enumerate}[(1)]
		\item For all countably complete uniform ultrafilters \(U\), if \(U\sE
		D\), then \(U\I D\).
		\item \(D\) is a Mitchell point
		\item For all Mitchell points \(D'\), if \(D' \sE D\), then \(D'\mo D\).
	\end{enumerate}
	\begin{proof}
		{\it (1) implies (2):} Note that (1) implies in particular that \(U\) is
		\(\delta\)-internal. Thus \(U\) is a uniform ultrafilter on \(\delta\).
		There are two cases. Suppose first that \(D = \mathscr K_\delta\). Then
		\cref{LambdaIrredChar}, \(D\) is \(\delta\)-Mitchell, which is what (2)
		asserts. Assume instead that \(D\neq \mathscr K_\delta\), so \(\mathscr
		K_\delta \sE D\) since \(\mathscr K_\delta\) is the least uniform
		ultrafilter on \(\delta\). By (1), \(\mathscr K_\delta\I U\), and in
		particular by \cref{Nonisolation1}, \(\delta\) is not isolated. By
		\cref{LambdaIrredChar}, \(U\) is \(\delta\)-irreducible, and therefore
		by the Irreducibility Theorem, \(U\) is \({<}\delta\)-supercompact and
		\(\delta\)-tight. Since \(\mathscr K_\delta\I D\),
		\cref{InternalToSuper} yields that \(j_D\) is \(\delta\)-supercompact.
		In particular, \(P(\delta)\subseteq M_D\), and so for any countably
		complete ultrafilter \(U\) on \(\delta\) with \(U\I D\), \(U\mo D\).
		Given (1), this implies (2).
		
		{\it (2) implies (3):} Immediate.
		
		{\it (3) implies (1):} Let \(D'\) be the \(\sE\)-least tail uniform
		ultrafilter that is not internal to \(D\). To show that (1) holds, we
		must show \(D' = D\).  Clearly \(D'\E D\) (since a nonprincipal
		ultrafilter is never internal to itself). By \cref{MOExtensions}, the
		internal relation extends the Mitchell order, so \(D'\not \mo D\).
		\cref{IPoint} asserts that \(D'\) has the following property: for any
		\(U\I D\), in fact \(U\I D'\). In particular, for any \(U\sE D'\), by
		the minimality of \(D'\), we have \(U\I D\), and so we can conclude that
		\(U\I D'\). Since we have shown that (1) implies (3), we can conclude
		that \(D'\) is a Mitchell point. Since \(D'\) is a Mitchell point and
		\(D'\not \mo D\), (3) implies that \(D'\not \sE D\). Since \(D'\E D\),
		it follows that \(D = D'\), as desired.
	\end{proof}
\end{prp}

\begin{defn}
	For any countably complete ultrafilter \(W\), the {\it Mitchell point of
	\(W\)}, denoted \(D(W)\), is the \(\sE\)-least tail uniform ultrafilter
	\(D\) such that \(D\not \mo W\).
\end{defn}
The proof of \cref{MitchellPoint} yields the following fact:
\begin{thm}[UA] Suppose \(W\) is a nonprincipal countably complete ultrafilter
	and \(D = D(W)\). Then the following hold:
	\begin{itemize}
		\item \(D\) is a Mitchell point.
		\item \(\{U : U \mo W\} = \{U : U\mo D\}\).
		\item If \(U\) is a countably complete ultrafilter such that \(U\I W\),
		then \(U\I D\).
		\item \(D\not \I W\).		\qed
	\end{itemize}
\end{thm}

\begin{thm}[UA]\label{Cofinal1}
The Mitchell points form a cofinal class of ultrafilters.
	\begin{proof}[Sketch] Suppose \(U\) is a countably complete ultrafilter. We
		will show that there is an ultrafilter \(U'\) given by a Mitchell point
		iteration such that \(U\D U'\). By induction, we may assume that this
		statement is true for all \(\bar U \sE U\). Let \(D = D(U)\). Since
		\(D\not \I U\), \(\tr D U\sE j_D(U)\) in \(M_D\). Therefore by our
		induction hypothesis, \(M_D\) satisfies that there is an ultrafilter
		\(W'\) given by a Mitchell point iteration of such that \(\tr D U\D
		W'\). Let \(U'\) be such that \(j_{U'} = j_{W'}^{M_D}\circ j_D\). It is
		easy to see that \(U'\) is given by a Mitchell point iteration and \(U\D
		U'\).
	\end{proof}
\end{thm}

We now turn to the other direction of \cref{SOtoUA}. It would be enough to prove
the following fact:
\begin{prp}\label{SODirected}
Suppose \(C\) is a coherent class of countably complete ultrafilters. Then the
restriction of the Rudin-Frol\'ik order to the  class of ultrafilters given by
\(C\)-iterations is directed.
\begin{proof}
	The idea of the proof is that the ultrafilters in \(C\) can be compared by
	the comparisons given by the internal relation \cref{InternalComparison},
	and then this can be propagated to compare arbitrary \(C\)-iterations by
	recursion. This is quite easy to see (given the right definition of a
	coherent class), but we nevertheless include a very detailed
	proof.\footnote{We caution, however, that as usual it may be easier for the
	reader work out the details than to read them.}
	
	We use the following convention: if \(\mathcal I\) is an iterated ultrapower
	of length \(\ell\), then \(j^{\mathcal I} = j_{0\ell}^{\mathcal I}\).
	
	We begin with a one-step claim:
	\begin{clm}\label{1StepClm}Suppose \(D\in C\). For any \(C\)-iteration
	\(\mathcal I\),  there is a \(C\)-iteration \(\mathcal J\) such that
	\(U^{\mathcal J}_0 = D\)  and a \(C\)-iteration \(\mathcal I'\) extending
	\(\mathcal I\) such that \(j^{\mathcal I'} = j^{\mathcal J}\).
	\end{clm}
	\begin{proof}[Proof of \cref{1StepClm}]The proof is by induction on the
	length of \(\mathcal I\). 
	
	If \(U^{\mathcal I}_0 = D\), then we can take \(\mathcal I = \mathcal J\). 
	
	Therefore assume \(U^\mathcal I_0\neq D\). Since \(C\) is coherent, either
	\(D\mo  U^\mathcal I_0\) or \(U^\mathcal I_0\mo D\). Define
	\[D_* =\begin{cases} j^\mathcal I_{01}(D)&\text{if }U^\mathcal I_0\mo D\\
	D&\text{if }D\mo U^\mathcal I_0\end{cases}\] and
	\[U^*_1 =\begin{cases} U^\mathcal I_0&\text{if }U^\mathcal I_0\mo D\\
	j_D(U^\mathcal I_0)&\text{if }D\mo U^\mathcal I_0\end{cases}\] The key point
	is  that by the definition of a coherent class of ultrafilters, \(D_*\in
	j^\mathcal I_{01}(C)\), \(U_1^*\in j_D(C)\), and \[j^{M_D}_{U_1^*}\circ j_D
	= j_{D_*}^{M^\mathcal I_1}\circ j^\mathcal I_{01}\] 
	
	Let  \(\mathcal I_* = \mathcal I\restriction [1,\infty)\), which is a
	\(j^{\mathcal I}_{01}(C)\)-iteration of \(M^\mathcal I_1\). By our induction
	hypothesis applied in \(M^\mathcal I_1\) to the  and the ultrafilter
	\(D_*\in j^\mathcal I_{01}(C)\), there is a \(j^\mathcal
	I_{01}(C)\)-iteration \(\mathcal J_*\) with \(U_0^{\mathcal J_*} = D_*\) and
	a \(j^\mathcal I_{01}(C)\)-iteration \(\mathcal I_*'\) extending \(\mathcal
	I_*\) such that \(j^{\mathcal I_*'} = j^{\mathcal J_*}\). 
	
	Let \(\mathcal I'\) be the iterated ultrapower of \(V\) given by
	\(U_0^\mathcal I\) followed by \(\mathcal I_*'\). Clearly \(\mathcal I'\) is
	a \(C\)-iteration extending \(\mathcal I\). Let \(\ell = \text{lth}(\mathcal
	J_*)\), and define a \(C\)-iteration \(\mathcal J\) of length \(\ell+1\) in
	terms of the ultrafilters \(U^\mathcal J_n\):
	\begin{align*}
		U_0^{\mathcal J} &= D\\
		U_1^{\mathcal J} &= U_1^*\\
		U_n^{\mathcal J} &=U_{n-1}^{\mathcal J_*}
	\end{align*}
	Then \[j^\mathcal J = j^{\mathcal J_*}_{1\ell-1}\circ j^{M_D}_{U^*_1}\circ
	j_D = j^{\mathcal J_*}_{1\ell-1}\circ j^{M^\mathcal I_1}_{D_*}\circ
	j^\mathcal I_{01} = j^{\mathcal J_*}\circ j^\mathcal I_{01} =  j^{\mathcal
	I'_*}\circ j^\mathcal I_{01} = j^{\mathcal I'}\] This verifies the induction
	step, and proves the claim.\end{proof}
	
	We now turn to the multi-step claim:
	\begin{clm}\label{MultiClm}For any \(C\)-iteration \(\mathcal H\), for any \(C\)-iteration \(\mathcal I\), there are \(C\)-iterations \(\mathcal H^*\) and \(\mathcal I^*\) extending \(\mathcal H\) and \(\mathcal I\) respectively such that \(j^{\mathcal H^*} = j^{\mathcal I^*}\).\end{clm}
	\begin{proof}[Proof of \cref{MultiClm}] The proof is by induction on the
	length \(\ell\) of \(\mathcal H\): thus our induction hypothesis is that for
	any \(C\)-iteration \(\bar {\mathcal H}\) of length less \(\ell\), for any
	\(C\)-iteration \(\mathcal I\), there are \(C\)-iterations \(\bar{\mathcal
	H}^*\) and \(\mathcal I'\) extending \(\bar{\mathcal H}\) and \(\mathcal I\)
	respectively such that \(j^{\bar {\mathcal H}^*} = j^{\mathcal I'}\). 
	
	Let \(D = U_0^\mathcal H\). By our first claim, there is a \(C\)-iteration
	\(\mathcal J\) such that \(U^{\mathcal J}_0 = D\) and a \(C\)-iteration
	\(\mathcal I'\) extending \(\mathcal I\) such that \(j^{\mathcal I'} =
	j^{\mathcal J}\). Now we work in \(M_D\). Let \(\bar {\mathcal H} = \mathcal
	H\restriction [1,\infty)\). Thus \(\bar {\mathcal H}\) is a
	\(j_D(C)\)-iteration of \(M_D\) of length less than \(\ell\). Let \(\bar
	{\mathcal J} = \mathcal J\restriction [1,\infty)\), so that \(\bar J\) is
	also a \(j_D(C)\)-iteration of \(M_D\). 
	
	By our induction hypothesis applied in \(M_D\), there are
	\(j_D(C)\)-iterations \(\bar {\mathcal H}^*\) and \(\bar {\mathcal J}^*\) of
	\(M_D\) extending \(\bar {\mathcal H}\) and \(\bar {\mathcal J}\)
	respectively such that \(j^{\bar{\mathcal H}^*} = j^{\bar{\mathcal J}^*}\).
	Define 
	\begin{align*}\mathcal H^* &= D^\frown \bar {\mathcal H}^*\\ 
		\mathcal I^* &= \mathcal I'^\frown \mathcal K\end{align*}
	 where \(\mathcal K\) is the iteration such that \(\bar {\mathcal J}^* =
	 \bar {\mathcal J}^\frown \mathcal K\).
	
	Obviously \(\mathcal H^*\) and \(\mathcal I^*\) are \(C\)-iterations
	extending \(\mathcal H\) and \(\mathcal I\) respectively. Moreover
	\[j^{\mathcal H'} = j^{\bar {\mathcal H}'} \circ j_D =  j^{\bar {\mathcal
	J}'} \circ j_D =  j^{\mathcal K}\circ j^{\bar {\mathcal J}} \circ j_D =
	j^{\mathcal K}\circ j^{\mathcal J} = j^{\mathcal K}\circ j^{\mathcal I'} =
	j^{\mathcal I'}\] This proves the claim
	\end{proof}
	It follows easily from \cref{MultiClm} that the restriction of the
	Rudin-Frol\'ik order to the class of ultrafilters given by \(C\)-iterations
	is directed.
\end{proof}
\end{prp}

We finally prove our characterization of UA in terms of coherent cofinal
sequences.

\begin{proof}[Proof of \cref{SOtoUA}] {\it (1) implies (2):} This is immediate
	from \cref{MPSO} and \cref{Cofinal1}.
	
	{\it (2) implies (1):} Let \(C\) be a coherent cofinal class of
	ultrafilters. Since \(C\) is coherent, \cref{SODirected} implies that the
	restriction of the Rudin-Frol\'ik order to the class of ultrafilters \(C'\)
	given by \(C\)-iterations is directed. Since \(C\) is cofinal, \(C'\) is
	cofinal in the Rudin-Frol\'ik order. Since the Rudin-Frol\'ik order has a
	cofinal directed subset, the Rudin-Frol\'ik order is itself directed. This
	implies that the Ultrapower Axiom holds (by \cref{RFDirected}).
\end{proof}
\section{Very large cardinals}\label{VLCSection}
\subsection{Huge cardinals}
The notion of {\it \((\kappa,\lambda)\)-regularity} is a two cardinal
generalization of \(\kappa^+\)-incompleteness that has already shown up
implicitly in this dissertation:
\begin{defn}\index{Regular ultrafilter}
	Suppose \(\kappa\leq \lambda\) are cardinals. An ultrafilter \(U\) is {\it
	\((\kappa,\lambda)\)-regular} if there is a set \(F\subseteq U\) of
	cardinality \(\lambda\) such that \(\bigcap \sigma\notin U\) for any
	\(\sigma\subseteq F\) of cardinality at least \(\kappa\).
\end{defn}
The combinatorial definition of \((\kappa,\lambda)\)-regularity defined above
obscures its true significance:
\begin{lma}\label{RegularityChar}
	Suppose \(\kappa\leq \lambda\) are cardinals and \(U\) is an ultrafilter.
	Then the following are equivalent:
	\begin{enumerate}[(1)]
		\item \(U\) is \((\kappa,\lambda)\)-regular.
		\item For some fine ultrafilter \(\mathcal U\) on \(P_\kappa(\lambda)\),
		\(\mathcal U\RK U\).
		\item \(j_U\) is \((\lambda,\delta)\)-tight for some \(M_U\)-cardinal
		\(\delta < j_U(\kappa)\).
	\end{enumerate}
	\begin{proof}
		{\it (1) implies (2):} Fix a set \(F\subseteq U\) of cardinality
		\(\lambda\) such that \(\bigcap \sigma\notin U\) for any
		\(\sigma\subseteq F\) of cardinality at least \(\kappa\). Let \(X\) be
		the underlying set of \(U\). Define \(f : X\to P_\kappa(F)\) by setting
		\(f(x) = \{A\in F: x\in A\}\). Let \(\mathcal U = f_*(U)\). We claim
		\(\mathcal U\) is a fine ultrafilter on \(P_\kappa(F)\). Suppose \(A\in
		F\). We must show \(\{\sigma\in P_\kappa(F) : A\in \sigma\}\in \mathcal
		U\). But by the definition of \(f\), \(A\in f(x)\) if and only if \(x\in
		A\). Thus \[f^{-1}[\{\sigma\in P_\kappa(F) : A\in \sigma\}] = A\in U\]
		and so \(\{\sigma\in P_\kappa(F) : A\in \sigma\}\in \mathcal U\).
		
		{\it (2) implies (3):} Fix a fine ultrafilter \(\mathcal U\) on
		\(P_\kappa(\lambda)\) such that \(\mathcal U \RK U\). Let \(A =
		\id_\mathcal U\). Then \(j_\mathcal U[\lambda]\subseteq A\) by
		\cref{NormalFineChar}, and \(|A|^{M_{\mathcal U}} < j_\mathcal
		U(\kappa)\) by \L o\'s's Theorem. Let \(k : M_\mathcal U\to M_U\) be an
		elementary embedding such that \(k\circ j_\mathcal U = j_U\). Then
		\(j_U[\lambda] = k[j_\mathcal U[\lambda]]\subseteq k(A)\) and
		\(|k(A)|^{M_U} < k(j_\mathcal U(\kappa)) = j_U(\kappa)\). Let \(\delta =
		|k(A)|^{M_U}\). Then \(k(A)\) witnesses that \(j_U\) is
		\((\lambda,\delta)\)-tight, as desired.
		
		{\it (3) implies (1):} Fix \(A\in M_U\) such that \(|A|^{M_U} <
		j_U(\kappa)\) and \(j_U[\lambda]\subseteq A\). Let \(f\) be a function
		such that \(A = [f]_U\). By \L o\'s's Theorem, there is a set \(X\in U\)
		such that \(f[X]\subseteq P_\kappa(\lambda)\). Let \(S_\alpha = \{x \in
		X : \alpha\in f(x)\}\). Let \(F = \{S_\alpha : \alpha < \lambda\}\). We
		claim that \(\bigcap_{\alpha \in \sigma} S_\alpha = \emptyset\) for any
		\(\sigma\subseteq \lambda\) of cardinality at least \(\kappa\). Suppose
		towards a contradiction that \(x\in \bigcap_{\alpha\in \sigma}
		S_\alpha\). Then \(\sigma\subseteq f(x)\), so \(|f(x)|\geq \kappa\),
		contradicting that \(f(x)\in P_\kappa(\lambda)\). Thus \(F\) witnesses
		that \(U\) is \((\kappa,\lambda)\)-regular.
	\end{proof}
\end{lma}
Another way of stating (2) above is to say that the minimum fine filter on
\(P_\kappa(\lambda)\) lies below \(U\) in the Kat\v etov order.

\begin{defn}
	If \(\kappa \leq \lambda\) are cardinals, then \(P^\kappa(\lambda)\) denotes
	the collection of subsets of \(\lambda\) of cardinality exactly \(\kappa\).
\end{defn}

Thus \(P^\kappa(\lambda) = P_{\kappa^+}(\lambda)\setminus P_\kappa(\lambda)\).
\begin{defn}
	A cardinal \(\kappa\) is {\it huge} if there is an elementary embedding \(j
	: V\to M\) with critical point \(\kappa\) such that \(M^{j(\kappa)}\subseteq
	M\).
\end{defn}

A question raised in \cite{Ketonen} is the relationship between nonregular
ultrafilters and huge cardinals. Assuming UA, we can almost show an equivalence:
\begin{thm}[UA]\label{HugeThm}\index{Huge cardinal}
	Suppose \(\kappa < \lambda\) are cardinals and \(\lambda\) is regular. The
	following are equivalent:
	\begin{enumerate}[(1)]
		\item There is a countably complete fine ultrafilter on
		\(P^\kappa(\lambda)\) that cannot be pushed forward to a fine
		ultrafilter on \(P_\kappa(\lambda)\).
		\item There is a countably complete ultrafilter that is
		\((\kappa^+,\lambda)\)-regular but not \((\kappa,\lambda)\)-regular. 
		\item There is an elementary embedding \(j : V\to M\) such that
		\(j(\kappa) = \lambda\), \(M^{<\lambda}\subseteq M\), and \(M\) has the
		\({\leq}\lambda\)-covering property.
	\end{enumerate}
	 If \(\lambda\) is a successor cardinal, then we can add to the list: 
	 \begin{enumerate}[(1)]
	 \item[(4)] There is an elementary embedding \(j : V\to M\) such that
	 \(j(\kappa) = \lambda\) and \(M^{\lambda}\subseteq M\).
	 \item[(5)] There is a normal fine ultrafilter on \(P^\kappa(\lambda)\).
	 \end{enumerate}
	\begin{proof}
		
		The equivalence of (1) and (2) is immediate from \cref{RegularityChar}.
		We now turn to the equivalence of (2) and (3). Before we begin,  we
		point out that the property of being \((\kappa^+,\lambda)\)-regular but
		not \((\kappa,\lambda)\)-regular can be reformulated in terms of
		ultrapowers: 
		\begin{quote}\(U\) is \((\kappa^+,\lambda)\)-regular but not
		\((\kappa,\lambda)\)-regular if and only if \(\text{cf}^{M_U}(\sup
		j_U[\lambda]) = j_U(\kappa)\).\end{quote} This is an immediate
		consequence of \cref{RegularityChar} (3) and Ketonen's analysis of tight
		embeddings in terms of cofinality (\cref{KetonenCov}).
		
		{\it (2) implies (3):} Let \(U\) be the \(\sE\)-least countably complete
		ultrafilter concentrating on ordinals that is
		\((\kappa^+,\lambda)\)-regular but not \((\kappa,\lambda)\)-regular.
		
		We claim that \(U\) is \(\lambda\)-irreducible. (In fact, \(U\) is an
		irreducible weakly normal ultrafilter on \(\lambda\), but this is not
		relevant to the proof.) Suppose \(D\D U\) and \(\lambda_D < \lambda\).
		We must show that \(D\) is principal. We claim \(\tr D U\) is
		\((j_D(\kappa^+),j_D(\lambda))\)-regular but not
		\((j_D(\kappa),j_D(\lambda))\)-regular. Let \(i : M_D\to M_U\) be the
		unique internal ultrapower embedding with \(i\circ j_D = j_U\). Thus \(i
		: M_D\to M_U\) is the ultrapower of \(M_D\) by \(\tr D U\). Therefore to
		show that \(\tr D U\) is \((j_D(\kappa),j_D(\lambda))\)-regular it
		suffices (by our remark at the beginning of the proof) to show that
		\(\text{cf}^{M_U}(\sup i[j_D(\lambda)]) = i(j_D(\kappa))\). Since
		\(\lambda_D < \lambda\), by \cref{UFContinuity}, \[\sup i[j_D(\lambda)]
		= \sup i\circ j_D[\lambda] = \sup j_U[\lambda]\] Furthermore, since
		\(U\) is \((j_D(\kappa^+),j_D(\lambda))\)-regular but not
		\((j_D(\kappa),j_D(\lambda))\)-regular, applying our remark at the
		beginning of the proof again, \[\text{cf}^{M_U}(\sup j_U[\lambda]) =
		j_U(\kappa) = i(j_D(\kappa))\] Thus \(\text{cf}^{M_U}(\sup
		i[j_D(\lambda)]) = i(j_D(\kappa))\), as desired.
		
		By elementarity \(j_D(U)\) is the \(\sE\)-least ultrafilter that is
		\((j_D(\kappa^+),j_D(\lambda))\)-regular but not
		\((j_D(\kappa),j_D(\lambda))\)-regular. Hence \(j_D(U) \E \tr D U\).
		Recall \cref{Pushdown}, which states that if \(D\) is nonprincipal and
		\(D\D U\), then \(\tr D U\sE j_D(U)\). It follows that \(D\) is
		principal.
		
		Since \(U\) is \(\lambda\)-irreducible, and now we would like to apply
		the Irreducibility Theorem. For this, we need that \(\lambda\) is either
		a successor cardinal or an inaccessible cardinal. Assume \(\lambda\) is
		a limit cardinal, and we will show that \(\lambda\) is a strong limit
		cardinal. Since \(\kappa < \lambda\), we have \(\kappa^+ < \lambda\).
		Since \(U\) is \((\kappa^+,\lambda)\)-regular, \(U\) is
		\(\delta\)-decomposable for all regular cardinals \(\delta \in
		[\kappa^+,\lambda]\). Therefore \(\lambda\) is a limit of Fr\'echet
		cardinals, and hence by \cref{LimitLma}, \(\lambda\) is a strong limit
		cardinal, as desired.
		
		To summarize, \(j_U : V\to M_U\) is an elementary embedding such that
		\(j_U(\kappa) = \lambda\), \(M_U^{<\lambda}\subseteq M_U\) and \(M_U\)
		has the \({\leq}\lambda\)-covering property. This shows that (3) holds.
		
		{\it (3) implies (2):} Let \(U\) be the ultrafilter on \(\lambda\)
		derived from \(j\) using \(\sup j[\lambda]\), and let \(k : M_U\to M\)
		be the factor embedding with \(k\circ j_U = j\) and \(k(\id_U) = \sup
		j[\lambda]\). Then \(\id_U = \sup j_U[\lambda]\), and
		\(k(\text{cf}^{M_U}(\id)) = \text{cf}^M(\sup j[\lambda]) = \lambda =
		j(\kappa) = k(j_U(\kappa))\). By the elementarity of \(k\),
		\[\text{cf}^{M_U}(\sup j_U[\lambda]) = \text{cf}^{M_U}(\id_U) =
		j_U(\kappa)\] Thus by our remark at the beginning of the proof, \(U\) is
		\((\kappa^+,\lambda)\)-regular but not \((\kappa,\lambda)\)-regular.
		This shows that (1) holds.
		
		Assuming \(\lambda\) is a successor cardinal, the argument that (2)
		implies (3) shows that in fact (2) implies (4), since the Irreducibility
		Theorem leads to full \(\lambda\)-supercompactness in the case that
		\(\lambda\) is a successor cardinal.
		
		Finally, (4) and (5) are equivalent (in general) by an easy argument
		using derived ultrafilters and ultrapowers (\cref{DerivedNF}).  
	\end{proof}
\end{thm}
We cannot show that \(M^\lambda\subseteq M\) in the key case that \(\lambda\) is
inaccessible, which blocks proving the equivalence between huge cardinals and
nonregular countably complete ultrafilters.
\subsection{Cardinal preserving elementary embeddings}
In this section, we turn to even larger large cardinal axioms.
\begin{defn}\index{Cardinal preserving embedding}
	An elementary embedding \(j : V\to M\) is {\it weakly cardinal preserving}
	if whenever \(\kappa\) is a cardinal, \(j(\kappa)\) is also a cardinal.
\end{defn}
The following question, due to Caicedo, essentially asks whether the Kunen
Inconsistency Theorem can be strengthened to rule out cardinal preserving
elementary embeddings:
\begin{qst}
	Is it consistent that there is a nontrivial weakly cardinal preserving
	elementary embedding?
\end{qst}
Under UA, we will show that there are no nontrivial weakly cardinal preserving
embeddings. 
\begin{lma}[UA]\label{BadSuccessor}
	Suppose \(U\) is a countably complete uniform ultrafilter on \(\kappa^+\)
	such that \(j_U[\kappa]\subseteq \kappa\). Either \(\kappa\) is
	\(\kappa^+\)-supercompact or \(\kappa\) is a limit of
	\(\kappa^+\)-supercompact cardinals.
	\begin{proof}
	By \cref{SupercompactFactor}, there is some \(D\D U\) with \(\lambda_D <
	\kappa^+\) such that there is an internal ultrapower embedding \(i : M_D\to
	M_U\) with \(i\circ j_D = j_U\) that is \(j_D(\kappa^+)\)-supercompact in
	\(M_D\). Note that \(\sup j_D[\kappa]\subseteq \kappa\) and \(\sup
	i[\kappa]\subseteq \kappa\), since both \(i\) and \(j_D\) are bounded on the
	ordinals by \(j_U\).
	
	We claim that \( \textsc{crt}(i)\in [\kappa, j_D(\kappa)]\). To see this,
	note that \(\sup i[\kappa]\subseteq \kappa\) and \(i\) is
	\(\kappa\)-supercompact, so by the Kunen Inconsistency Theorem
	(\cref{KunenInconsistency}), \(\textsc{crt}(i)\geq \kappa\). On the other
	hand, \(i\) is given by an ultrafilter on \(j_D(\kappa^+)\), so
	\(\textsc{crt}(i) \leq j_D(\kappa)\). 
	
	Now \(i\) witnesses that \(\textsc{crt}(i)\) is
	\(j_D(\kappa^+)\)-supercompact in \(M_D\). If \(\textsc{crt}(i) =
	j_D(\kappa)\), then \(\kappa\) is \(\kappa^+\)-supercompact by elementarity.
	Otherwise \(\sup j_D[\kappa] = \kappa \leq \textsc{crt}(i) < j_D(\kappa)\),
	so \(\kappa\) is a limit of \(\kappa^+\)-supercompact cardinals by a
	standard reflection argument.  
	\end{proof}
\end{lma}

The following observation is due to Caicedo:

\begin{lma}\label{CaicedoLemma}
	Suppose \(j : V\to M\) and \(\gamma\) is a cardinal. If \(j(\gamma^+) \neq
	\gamma^+\) and \(j\) is continuous at \(\gamma^+\), then \(j(\gamma^+)\) is
	not a cardinal.
	\begin{proof}
		Note that \(j(\gamma^+)\) is a singular ordinal since \(j[\gamma^+]\) is
		cofinal in \(j(\gamma^+)\). Moreover \(j(\gamma) < j(\gamma^+) =
		j(\gamma)^{+M}\leq j(\gamma)^+\). There are no singular cardinals
		between \(j(\gamma)\) and \(j(\gamma)^+\), so \(j(\gamma^+)\) is not a
		cardinal.
	\end{proof}
\end{lma}

\begin{lma}[UA]\label{NoCardinalPreserve}
	Suppose \(j : V\to M\) is a nontrivial elementary embedding with critical
	point \(\kappa\). Let \(\gamma\) be a cardinal above \(\kappa\) with
	\(j(\gamma) = \gamma\). Then \(j\) is continuous at \(\gamma^{+\kappa+1}\)
	and therefore \(j(\gamma^{+\kappa+1})\) is not a cardinal.
	\begin{proof}
		We begin the proof by making some general observations about the action
		of \(j\) on cardinals in the vicinity of \(\gamma\). First, for all
		\(\alpha < \kappa\), \(j(\gamma^{+\alpha}) = (\gamma^{+\alpha})^M \leq
		\gamma^{+\alpha}\). It follows that \(j(\gamma^{+\alpha}) =
		\gamma^{+\alpha}\). Hence \(\sup j[\gamma^{+\kappa}] =
		\gamma^{+\kappa}\).
		
		Next, we claim that \((\gamma^{+\kappa+1})^M = \gamma^{+\kappa+1}\).
		This is proved by following the argument of \cref{Ineqs}: fix \(\alpha <
		\gamma^{+\kappa+1}\), and we will show that \(\alpha <
		(\gamma^{+\kappa+1})^M\). Let \((\gamma^{+\kappa},\prec)\) be a
		wellorder of order type \(\alpha\). Then \((\gamma^{+\kappa},j(\prec))\)
		is a wellorder of \(\gamma^{+\kappa}\) that belongs to \(M\). Since
		\(j[\gamma^{+\kappa}]\subseteq \gamma^{+\kappa}\), \(j\) embeds
		\((\gamma^{+\kappa},\prec)\) into \((\gamma^{+\kappa},j(\prec))\), so
		\[\alpha \leq
		\text{ot}(\gamma^{+\kappa},\prec)\leq\text{ot}(\gamma^{+\kappa},j(\prec))<
		(\gamma^{+\kappa+1})^M\] as desired.
		
		It follows that \[j(\gamma^{+\kappa+1}) > j(\gamma^{+\kappa}) =
		(\gamma^{+j(\kappa)})^M > (\gamma^{+\kappa+1})^M = \gamma^{+\kappa+1}\]
		Thus to prove \(j(\gamma^{+\kappa+1})\) is not a cardinal, by
		\cref{CaicedoLemma} it suffices to show \(j\) is continuous at
		\(\gamma^{+\kappa+1}\).
		
		Suppose towards a contradiction that \(j\) is discontinuous at
		\(\gamma^{+\kappa+1}\). Let \(U\) be the ultrafilter on
		\(\gamma^{+\kappa+1}\) derived from \(j\) using \(\sup
		j[\gamma^{+\kappa+1}]\). Then \(U\) is a countably complete uniform
		ultrafilter on \(\gamma^{+\kappa+1}\). Moreover, \[\sup
		j_U[\gamma^{+\kappa}]\leq \sup j[\gamma^{+\kappa}] = \gamma^{+\kappa}\]
		Therefore by \cref{BadSuccessor}, \(\gamma^{+\kappa}\) is either
		\(\gamma^{+\kappa+1}\)-supercompact or else a limit of
		\(\gamma^{+\kappa+1}\)-supercompact cardinals. This is impossible since
		there are no inaccessible cardinals in the interval
		\((\gamma,\gamma^{+\kappa}]\). Thus our assumption was false, and in
		fact \(j\) is continuous at \(\gamma^{+\kappa+1}\). 
		
		Now \(j\) is continuous at \(\gamma^{+\kappa+1}\) and
		\(j(\gamma^{+\kappa+1}) > \gamma^{+\kappa+1}\). Therefore by
		\cref{CaicedoLemma}, \(j(\gamma^{+\kappa+1})\) is not a cardinal.
	\end{proof}
\end{lma}
\begin{cor}[UA] Any weakly cardinal preserving elementary embedding of the
	universe is the identity.
\end{cor} 

We now investigate the relationship between cardinal preservation and
rank-into-rank axioms.
\begin{thm}[UA]\label{I3Thm}\index{Rank-into-rank cardinal}
	Assume \(\lambda\) is an ordinal, \(M\subseteq V_\lambda\) is a transitive
	set, and \(j : V_\lambda\to M\) is an elementary embedding with critical
	point \(\kappa\) that has no fixed points above \(\kappa\). Suppose that
	\(\textnormal{Card}^M\cap \lambda = \textnormal{Card}\cap \lambda\). Then
	\(M = V_\lambda\).
\end{thm}
If the assumption that \(\textnormal{Card}^M\cap \lambda = \textnormal{Card}\cap
\lambda\) is weakened to the assumption that \(j\) is weakly cardinal preserving
below \(\lambda\) (or in other words that \(j[\text{Card}\cap \lambda]\subseteq
\text{Card}\cap \lambda\)), then the resulting statement is false. Let us
provide a counterexample. Suppose \(j : V\to M\) is an elementary embedding with
critical point \(\kappa\). Let \(\lambda\) be the first cardinal fixed point of
\(j\) above \(\kappa\). Assume \(V_\lambda\subseteq M\), so \(j\) witnesses the
axiom \(I_2\). Suppose \(U\) is a \(\kappa\)-complete ultrafilter on \(\kappa\).
Then by \cref{RanFix2}, \(j^M_U\circ j : V\to (M_U)^M\) has the property that
\(j^M_U\circ j\restriction \text{Ord} = j\restriction \text{Ord}\), so in
particular \(j^M_U\circ j[\text{Card}\cap \lambda] = j[\text{Card}\cap
\lambda]\subseteq \text{Card}\cap \lambda\). But of course \((M_U)^M\) does not
contain \(V_\lambda\).

One of the key lemmas is the following curiosity, a close cousin of
\cref{AlwaysTight}:
\begin{lma}[UA]\label{SuccPreservation}
	Suppose \(U\) is a countably complete ultrafilter and \(\delta\) is a
	successor cardinal. Then \(\textnormal{cf}^{M_U}(\sup j_U[\delta])\) is a
	successor cardinal of \(M_U\).
	\begin{proof}
	If \(\sup j_U[\delta] = j_U(\delta)\), then \(\sup j_U[\delta]\) is itself a
	successor cardinal of \(M_U\), so of course its \(M_U\)-cofinality (which is
	again \(\sup j_U[\delta]\)) is a successor cardinal of \(M_U\). We may
	therefore assume that \(\sup j_U[\delta] < j_U(\delta)\). 
		
	 Hence \(\delta\) is Fr\'echet, and so we are in a position to apply
	 \cref{SuccessorIrredThm}. By \cref{SuccessorIrredThm}, there is an
	 ultrafilter \(D\) with \(\lambda_D < \delta\) such that there is an
	 internal ultrapower embedding \(h : M_D\to M_U\) such that \(h\) is
	 \(j_D(\delta)\)-supercompact in \(M_D\). Since \(\lambda_D < \delta\),
	 \(j_D(\delta) = \sup j_D[\delta]\) by \cref{UFContinuity}. Thus
	 \[\text{cf}^{M_U}(\sup j_U[\delta]) = \text{cf}^{M_U}(\sup h[j_D(\delta)])
	 = \text{cf}^{M_D}(j_D(\delta)) = j_D(\delta)\] Since \(j_D(\delta)\) is a
	 successor cardinal of \(M_D\), and \(\text{Ord}^{j_D(\delta)}\cap
	 M_D\subseteq M_U\), \(j_D(\delta)\) is a successor cardinal of \(M_U\). 
	\end{proof}
\end{lma}

We now turn to the proof of \cref{I3Thm}. 
\begin{proof}[Proof of \cref{I3Thm}] For \(n < \omega\), let \(\kappa_n =
	\kappa_n^j\) be the \(n\)th element of the critical sequence of \(j\)
	(\cref{CriticalSequence}), and note that \(\lambda = \sup_{n < \omega}
	\kappa_n\) since \(j\) has no fixed points above \(\kappa\).
	
	Let us make some preliminary remarks about the interaction between
	ultrapowers and the structure \(V_\lambda\). Suppose that \(U\) is a
	countably complete ultrafilter on a set \(X\in V_\lambda\). Then any
	function \(f : X\to V_\lambda\) is bounded on a set in \(U\). In particular,
	\[j_U(V_\lambda) = \{[f]_U : f\in V_\lambda\text{ and dom}(f) = X\}\] In
	other words, \(V_\lambda\) correctly computes the ultrapower by \(U\). We
	will go to great lengths, however, not to work inside \(V_\lambda\), which
	we have not yet proved to be a model of ZFC.
	
	Suppose \(X\in V_\lambda\), \(a\in j(X)\), and \(U\) is the ultrafilter on
	\(X\) derived from \(j\) using \(a\). Then \(U\) is countably complete, so
	the remark of the previous paragraph applies. Thus we can define a factor
	embedding \(k : j_U(V_\lambda)\to M\) by setting \(k([f]_U) = j(f)(a)\)
	whenever \(f\in V_\lambda\) is a function on \(X\). The usual argument shows
	that \(k\) is well-defined and elementary. Moreover, \(k\circ
	(j_U\restriction V_\lambda) = j\) and \(k(\id_U) = a\).
	
	Suppose \(\delta < \lambda\) is a successor cardinal. Let \(U\) be the
	uniform ultrafilter derived from \(j\) using \(\sup j[\delta]\), and let \(k
	: j_U(V_\lambda)\to M\) be the factor embedding. We claim: 
	\begin{itemize}
		\item \(\text{cf}^{M}(\sup j[\delta]) = \delta\).
		\item \(j_U\) is \(\delta\)-tight.
		\item \(k(\delta) = \delta\).
	\end{itemize}
		
	 By \cref{SuccPreservation}, \(\sup j_U[\delta]\) is a successor cardinal of
	 \(M_U\). Thus \(\sup j_U[\delta]\) is a successor cardinal of
	 \(j_U(V_\lambda)\), so \(k(\sup j_U[\delta]) = \text{cf}^M(\sup
	 j[\delta])\) is a successor cardinal of \(M\). Since \(M\) is correct about
	 cardinals below \(\lambda\), \(\text{cf}^{M}(\sup j[\delta])\) is a
	 successor cardinal (in \(V\)). In particular, \(\text{cf}^{M}(\sup
	 j[\delta])\) is regular. Thus \(\text{cf}^{M}(\sup j[\delta]) =
	 \text{cf}(\text{cf}^{M}(\sup j[\delta])) = \text{cf}(\sup j[\delta]) =
	 \delta\), as desired.
	
	It follows that \(j_U\) is \(\delta\)-tight: \[\text{cf}^{M_U}(\sup
	j_U[\delta]) = \text{cf}^{j_U(V_\lambda)}(\sup j_U[\delta]) \leq
	k(\text{cf}^{j_U(V_\lambda)}(\sup j_U[\delta])) = \text{cf}^M(\sup
	j[\delta]) = \delta\] so \(\text{cf}^{M_U}(\sup j_U[\delta]) = \delta\), and
	hence \(j_U\) is \(\delta\)-tight by \cref{KetonenCov}.
	
	Repeating the same argument, it now follows that \(k(\delta) = \delta\):
	\[k(\delta) = k(\text{cf}^{M_U}(\sup j_U[\delta])) =
	k(\text{cf}^{j_U(V_\lambda)}(\sup j_U[\delta])) = \text{cf}^M(\sup
	j[\delta]) = \delta\]
	
	We recall an argument due to Caicedo-Woodin (\cite{Caicedo}) that shows that
	\(\kappa_n\) is strongly inaccessible for all \(n < \omega\). Suppose by
	induction that \(\kappa_n\) is strongly inaccessible, and we will show that
	\(\kappa_{n+1}\) is strongly inaccessible. Since \(\kappa_n\) is strongly
	inaccessible and \(j(\kappa_n) = \kappa_{n+1}\), \(M\) satisfies that
	\(\kappa_{n+1}\) is strongly inaccessible. Since \(M\) is cardinal correct
	below \(\lambda\) and satisfies that \(\kappa_{n+1}\) is a limit cardinal,
	\(\kappa_{n+1}\) is a limit cardinal (in \(V\)). It therefore suffices to
	show that for all successor cardinals \(\delta < \kappa_{n+1}\), \(2^\delta
	< \kappa_{n+1}\). Let \(U\) be the ultrafilter on derived from \(j\) using
	\(\sup j[\delta]\), and let \(k : j_U(V_\lambda)\to M\) be the factor
	embedding. Since \(j_U\) is \(\delta\)-tight, \( 2^\delta \leq
	(2^\delta)^{M_U}\) by \cref{TightContinuum}. Since \(\delta < \kappa_{n+1}\)
	and \(\kappa_{n+1}\) is strongly inaccessible in \(M\), \(M\) satisfies that
	\(2^\delta\) exists, and \((2^\delta)^M < \kappa_{n+1}\). Since \(k(\delta)
	= \delta\), by elementarity \(j_U(V_\lambda)\) satisfies that \(2^\delta\)
	exists, and hence \((2^\delta)^{j_U(V_\lambda)} = (2^\delta)^{M_U}\). Thus
	\[ 2^\delta \leq (2^\delta)^{M_U}  = (2^\delta)^{j_U(V_\lambda)} \leq k(
	(2^\delta)^{j_U(V_\lambda)}) = (2^\delta)^M < \kappa_{n+1}\] Thus \(\lambda
	= \sup_{n <\omega} \kappa_n\) is a limit of strongly inaccessible cardinals.
	
	Suppose \(\eta\in [\kappa,\lambda]\) is a strongly inaccessible cardinal. We
	will show \(V_\eta\subseteq M\). 
	
	Let \(U\) be the ultrafilter on \(\eta\) derived from \(j\) using \(\sup
	j[\eta]\). Let \(k : j_U(V_\lambda)\to M\) be the factor embedding.

	 By \cref{SupercompactFactor}, there is an ultrafilter \(D\) with
	 \(\lambda_D < \eta\) such that there is an internal ultrapower embedding
	 \(h : M_D\to M_U\) with \(h\circ j_D = j_U\) that is
	 \({<}j_D(\eta)\)-supercompact in \(M_D\). Since \(\lambda_D < \eta\) and
	 \(\eta\) is strongly inaccessible, \(j_D(\eta) = \eta\).
	 
	 In particular, we have that \(V_\eta \cap M_D = V_\eta \cap M_U\). We can
	 therefore define an elementary embedding \(i : V_\eta\to V_\eta\cap M\):
	 for \(x\in V_\eta\), set \(i(x) = k(j_D(x))\). Note that \(i\) is a weakly
	 cardinal preserving embedding of \(V_\eta\): \[\text{Card}^{V_\eta\cap M} =
	 \text{Card}^M\cap \eta = \text{Card}\cap \eta\] The second-order structure
	 \((V_\eta,V_{\eta+1})\) is a model of \(\text{NBG + UA}\), so we can apply
	 \cref{NoCardinalPreserve} in \((V_\eta,V_{\eta+1})\) to conclude that \(i\)
	 is the identity. In particular, \(V_\eta\cap M = V_\eta\), and therefore
	 \(V_\eta\subseteq M\), as desired.
	 
	 Since \(\eta < \lambda\) was an arbitrary inaccessible cardinal and
	 \(\lambda\) is a limit of inaccessible cardinals, \(V_\lambda\subseteq M\).
	 Hence \(M= V_\lambda\), as desired.
\end{proof}

The following question remains open:
\begin{qst}
	Suppose there is a weakly cardinal preserving elementary embedding from
	\(V_\lambda\) into a transitive set \(M\subseteq V_\lambda\). Must there be
	an elementary embedding \(j : V_\lambda\to V_\lambda\)?
\end{qst}
This cannot be entirely trivial: an application of \cref{RanFix2} shows that a
weakly cardinal preserving embedding itself need not have target model
\(V_\lambda\). Suppose \(\kappa < \lambda\) are cardinals, \(j : V\to M\) is an
elementary embedding with critical point \(\kappa\), \(j(\lambda) = \lambda\),
and \(V_\lambda\subseteq M\). Let \(U\) a \(\kappa\)-complete ultrafilter on
\(\kappa\). Then \(j_U\circ j\restriction \text{Ord} = j\) by \cref{RanFix2},
and it follows that \(j_U\circ j\restriction V_\lambda\) is weakly cardinal
preserving, even though its target model is \(M_U\cap V_\lambda\) and not
\(V_\lambda\).
\subsection{Supercompactness at inaccessible cardinals}\label{PathologicalSection}
The following are probably the most interesting questions left open by our work:
\begin{qst}[UA]\index{Supercompactness!at inaccessible cardinals}
	Suppose \(\lambda\) is an inaccessible cardinal and \(\kappa\) is the least
	\(\lambda\)-strongly compact cardinal. Must \(\kappa\) be
	\(\lambda\)-supercompact? More generally, if \(\kappa\) is
	\(\lambda\)-strongly compact, must \(\kappa\) be \(\lambda\)-supercompact or
	a measurable limit of \(\lambda\)-supercompact cardinals?
\end{qst}
This final chapter consists of some inconclusive observations regarding this
problem.

The whole question, it turns out, reduces to the analysis of \(\mathscr
K_\lambda\):
\begin{lma}[UA]\label{InaccessibleContainment}
	Assume \(\lambda\) is an inaccessible Fr\'echet cardinal. Let \(j : V\to M\)
	be the ultrapower of the universe by \({\mathscr K_\lambda}\), and let
	\(\kappa\) be the least measurable cardinal of \(M\) above \(\lambda\). Then
	for any \(\lambda\)-irreducible ultrafilter \(U\),
	\(\textnormal{Ord}^\kappa\cap M\subseteq M_U\). 
	\begin{proof}
		Let \((k,h) : (M,M_U)\to P\) be the pushout of \((j,j_U)\), and let
		\(W\) be such that \(P = M_W\). By the analysis of ultrafilters internal
		to a pushout, for any \(D\) with \(\lambda_D < \lambda\), since \(D\I
		U\) and \(D\I \mathscr K_\lambda\), in fact, \(D \I W\). In particular,
		\(W\) is \(\lambda\)-irreducible, so \(V_\lambda\subseteq M_W=P\) by
		\cref{LimitIrredThm}. By our factorization lemma for  embeddings of
		\(M\) (\cref{KFactor}), it follows that \(\textsc{crt}(k) \geq \kappa\).
		(Otherwise \(k\) would factor through an ultrapower by an ultrafilter in
		\(V_\lambda\), contrary to the fact that \(V_\lambda\subseteq P\).)
		Therefore \(\text{Ord}^\kappa\cap M\subseteq P\subseteq M_U\), as
		desired.
		\end{proof}
\end{lma}

\begin{cor}[UA]\label{KlambdaReduction}
	Suppose \(\lambda\) is a Fr\'echet inaccessible cardinal. Let \(M\) be the
	ultrapower of the universe by \(\mathscr K_\lambda\), and assume \(M\) is
	closed under \(\lambda\)-sequences. Then for any \(\lambda\)-irreducible
	ultrafilter \(U\), \(M_U\) is closed under \(\lambda\)-sequences.
	\begin{proof}
		By \cref{InaccessibleContainment}, \(\text{Ord}^\lambda =
		\text{Ord}^\lambda\cap M \subseteq M_U\), so \(M_U\) is closed under
		\(\lambda\)-sequences.
	\end{proof}
\end{cor}

We now show that the \(\sE\)-second irreducible ultrafilter on an inaccessible
cardinal \(\lambda\) always witnesses \(\lambda\)-supercompactness. This is a
bit surprising given that we cannot prove the supercompactness of \(\mathscr
K_\lambda\).

We use the following lemma, extracted from Ketonen's proof that the Ketonen
order is wellfounded on weakly normal ultrafilters.
\begin{lma}\label{ClubExtend}
	Suppose \(\lambda\) is a regular cardinal. Suppose \(W\) is a countably
	complete ultrafilter on \(\lambda\) that extends the closed unbounded
	filter. Suppose \(U\sE W\). Then \(\delta_{\tr U W} = j_U(\lambda)\). In
	fact, \(\tr U W\) extends the closed unbounded filter on \(j_U(\delta)\).
	\begin{proof}
		Let \(F\) be the closed unbounded filter on \(\lambda\). Clearly
		\(j_U[F]\subseteq \tr U W\). Moreover \(\{\alpha < j_U(\delta) :
		\id_U\in \alpha\}\in \tr U W\) since \[j_{\tr U W}^{M_U}(\id_U) < j_{\tr
		W U}^{M_W}(\id_W) = \id_{\tr U W}\] Thus by \cref{NormalGeneration},
		\(j_U(F)\subseteq \tr U W\), as claimed.
	\end{proof}
\end{lma}
We choose not to cite the Irreducibility Theorem in the proof of the following
proposition since it predates the Irreducibility Theorem and is really much
easier:
\begin{prp}[UA]\label{TwoUFs}
	Suppose \(\lambda\) is a regular cardinal. The following are equivalent: 
	\begin{enumerate}[(1)]
		\item \(\lambda\) carries distinct uniform irreducible ultrafilters.
		\item There is a countably complete uniform ultrafilter \(U\) such that
		\(\mathscr K_\lambda\not \D U\) and \(U\not \I \mathscr K_\lambda\).
		\item \(\lambda\) carries a countably complete weakly normal ultrafilter
		that concentrates on ordinals that carry countably complete tail uniform
		ultrafilters.
		\item \(\lambda\) carries distinct countably complete weakly normal
		ultrafilters.
		\item \(\lambda\) carries distinct countably complete ultrafilters
		extending the closed unbounded filter.
		\item There is a a normal fine \(\kappa_\lambda\)-complete ultrafilter
		\(\mathcal U\) on \(P_{\kappa_\lambda}(\lambda)\) such that \(\mathscr
		K_\lambda \mo \mathcal U\).
	\end{enumerate}
	\begin{proof}
		{\it (1) implies (2):} Suppose \(U\neq \mathscr K_\lambda\) is an
		irreducible ultrafilter on \(\lambda\). By irreducibility, \(\mathscr
		K_\lambda\not \D U\). Since \(\sup j_{\mathscr K_\lambda}[\lambda]\)
		carries no countably complete tail uniform ultrafilter in \(M_{\mathscr
		K_\lambda}\), \(j_U\restriction M_{\mathscr K_\lambda}\) is not internal
		to \(M_{\mathscr K_\lambda}\), since it is discontinuous at \(\sup
		j_{\mathscr K_\lambda}[\lambda]\). In other words \(U\not \I \mathscr
		K_\lambda\).
		
		{\it (2) implies (3):} Suppose \(U\) is a countably complete ultrafilter
		such that \(\mathscr K_\lambda\not \D U\) and \(U\not \I \mathscr
		K_\lambda\). Since \(U\not \I \mathscr K_\lambda\), by the
		characterization of internal ultrapower embeddings of \(M_{\mathscr
		K_\lambda}\) (\cref{EmbeddingChar}), \(j_U\) must be discontinuous at
		\(\lambda\). Since \(\mathscr K_\lambda \not\D U\), by the universal
		property of \(\mathscr K_\lambda\), \(\sup j_U[\lambda]\) carries a
		countably complete tail uniform ultrafilter in \(M_U\). Let \(W\) be the
		ultrafilter on \(\lambda\) derived from \(j_U\) using \(\sup
		j_U[\lambda]\). Then \(W\) is weakly normal (by \cref{DerivedWN}) and
		\(W\) concentrates on ordinals carrying countably complete tail uniform
		ultrafilters by the definition of a derived ultrafilter.
		
		{\it (3) implies (4):} If \(\lambda\) carries a countably complete
		uniform ultrafilter, then \(\lambda\) carries a countably complete
		weakly normal ultrafilter that {\it does not} concentrate on ordinals
		carrying countably complete tail uniform ultrafilters (by
		\cref{KetonenExistence}); in the context of UA, this is \(\mathscr
		K_\lambda\). Thus if (3) holds, \(\lambda\) carries distinct countably
		complete weakly normal ultrafilters.
		
		{\it (4) implies (5):} Immediate given the fact that weakly normal
		ultrafilters extend the closed unbounded filter.
		
		{\it (5) implies (6):} Assume (5) holds. Let \(U\) be the \(\sE\)-least
		countably complete ultrafilter that extends the closed unbounded filter
		on \(\lambda\) and is not equal to \(\mathscr K_\lambda\). We claim that
		for all \(D\sE U\), \(D\I U\). We will verify the criterion for showing
		\(D\I U\) given by \cref{IChar} by showing that \(j_D(U) \E \tr D U\) in
		\(M_D\).
		
		Let \(U' = \tr D U\). By \cref{ClubExtend}, \(U'\) extends the closed
		unbounded filter on \(j_D(\lambda)\). Moreover we claim that
		\(j_D(\mathscr K_\lambda) \neq U'\). To see this, note that
		\[j_D^{-1}[j_D(\mathscr K_\lambda)]  = \mathscr K_\lambda \neq U =
		j_D^{-1}[U']\] Thus \(j_D(\mathscr K_\lambda)\neq U'\), as claimed.
		
		By elementarity, in \(M_{D}\), \(j_{D}(U)\) is the \(\E\)-least
		countably complete ultrafilter that extends the closed unbounded filter
		on \(j_{D}(\lambda)\) and is not equal to \(j_{D}(\mathscr K_\lambda)\).
		It follows that \(j_{D}(U)\E U'\) in \(M_D\). \cref{IChar} now implies
		that \(D \I U\), as claimed.
		
		Let \(\kappa = \kappa_\lambda\). Since \(\lambda\) is not isolated, by
		\cref{KappaChar}, \(\kappa\) is a limit of isolated cardinals. By
		\cref{IsoThresh}, for all isolated cardinals \(\gamma < \kappa\),
		\(j_U[\gamma]\subseteq \gamma\), and hence \(j_U[\kappa]\subseteq
		\kappa\). \cref{InternalComplete} states that if \(\kappa\) is a strong
		limit cardinal such that \(j_U[\kappa]\subseteq \kappa\) and \(D\I U\)
		for all countably complete ultrafilters \(D\) with \(\lambda_D <
		\kappa\), then \(U\) is \(\kappa\)-complete. Thus \(U\) is
		\(\kappa\)-complete. In particular, \(\text{Ord}^\kappa\subseteq M_U\).
		Since \(\mathscr K_\lambda\I U\), \(j_{\mathscr
		K_\lambda}(\text{Ord}^\kappa) = \text{Ord}^{j_{\mathscr
		K_\lambda}(\kappa)}\cap M_{\mathscr K_\lambda}\subseteq M_U\). As
		\(j_{\mathscr K_\lambda}(\kappa) > \lambda\) by \cref{KetonenFrechet},
		it follows that \(\text{Ord}^\lambda \cap M_{\mathscr
		K_\lambda}\subseteq M_U\). 
		
		Now suppose \(A\in \text{Ord}^\lambda\). Then \(j_{\mathscr
		K_\lambda}[A]\) is contained in a set \(B\in [\text{Ord}]^\lambda \cap
		M_{\mathscr K_\lambda}\). Hence \(B\in M_U\). We may assume \(B\subseteq
		j_{\mathscr K_\lambda}(A)\), so that \(j_{\mathscr K_\lambda}^{-1}[B] =
		A\). Since \(\mathscr K_\lambda\I U\), \(j_{\mathscr
		K_\lambda}\restriction \alpha\in M_U\) for all ordinals \(\alpha\).
		Hence \(A = j_{\mathscr K_\lambda}^{-1}[B]\in M_U\). Thus
		\(\text{Ord}^\lambda\subseteq M_U\).
		
		If \(Z\) is a countably complete ultrafilter extending the closed
		unbounded filter on \(\lambda\) such that \(Z\mo \mathcal U_0\), then
		\(Z\I U\) so \(Z\sE U\) by \cref{InternalSeed} and consequently by the
		minimality of \(U\), \(Z = \mathscr K_\lambda\). In particular, no
		cardinal less than or equal to \(\lambda\) can be
		\(2^\lambda\)-supercompact in \(M_U\). It follows that \(j_U(\kappa) >
		\lambda\): otherwise \(j_U(\kappa) < \lambda\) is
		\(j_U(\lambda)\)-supercompact and since \(2^\lambda < j_U(\lambda)\), we
		contradict the previous sentence.
		
		Thus \(U\) is \(\kappa\) complete and \(j_U(\kappa) > \lambda\). Let
		\(\mathcal U\) be the normal fine \(\kappa\)-complete ultrafilter on
		\(P_\kappa(\lambda)\) derived from \(j_U\) using \(j_U[\lambda]\). It is
		easy to see that \(\mathscr K_\lambda\mo \mathcal U\) (and in fact \(U
		\cong \mathcal U\)). This completes the proof.
	\end{proof}
\end{prp}

We now turn to the question of pseudocompact cardinals first raised in
\cref{PseudocompactSection}. Recall that an elementary embedding is
\(\lambda\)-pseudocompact if it is \(\gamma\)-tight for all \(\gamma
\leq\lambda\). Our main question asked whether \(\lambda\)-pseudocompactness and
\(\lambda\)-supercompactness coincide below rank-into-rank cardinals. If
\(\lambda\) is the least cardinal where this fails, then it has the following
property:
\begin{defn}\index{Pathological cardinal}
	A cardinal \(\lambda\) is said to be {\it pathological} if there is an
	elementary embedding \(j : V\to M\) that is \({<}\lambda\)-supercompact and
	\(\lambda\)-tight but not \(\lambda\)-supercompact. The embedding \(j\) is
	said to {\it witness the pathology of \(\lambda\).} 
\end{defn}
Equivalently, \(j: V\to M\) witnesses the pathology of \(\lambda\) if
\(H(\lambda)\subseteq M\) and \(j[\lambda]\) can be covered by a set of size
\(\lambda\) in \(M\), yet \(j[\lambda]\notin M\). The axiom \(I_2(\lambda)\)
\index{\(I_2(\lambda)\)} asserts that there is an elementary embedding \(j :
V\to M\) with critical point less than \(\lambda\) such that \(j(\lambda) =
\lambda\) and \(V_\lambda\subseteq M\). By the Kunen Inconsistency Theorem
(\cref{KunenInconsistency}), \(j[\lambda]\notin M\). Thus if \(I_2(\lambda)\)
holds, then \(\lambda\) is pathological. 
\begin{qst}
	Suppose \(\lambda\) is pathological. Must \(\text{cf}(\lambda) = \omega\)?
	Must \(I_2(\lambda)\)?
\end{qst}
Our guess is that the answer is no.

We begin by establishing a dichotomy: pathological cardinals are either regular
or of countable cofinality. For the proof we use the following fact, a
generalization of the Kunen inconsistency theorem that is a slight improvement
on an observation due to Foreman \cite{VickersWelch}. (The proof of this theorem
included in the original version of this thesis was nonsense.)

\begin{thm}\label{ForemanInconsistency}\index{Kunen Inconsistency Theorem!Foreman's Theorem}
	Suppose \(\lambda\) is a cardinal. Suppose \(Q\) is a transitive set 
	that is closed under countable sequences and satisfies \(\textnormal{Ord}\cap Q = \lambda\). 
	Suppose \(k : Q\to H(\lambda)\) is a nontrivial elementary embedding. Let \(\gamma\) be
	the supremum of the critical sequence of \(k\). Then \(\lambda = \gamma^+\).
	\begin{proof}
		Since \(\gamma\) has countable cofinality and \(Q\) is closed under
		countable sequences, \(\gamma \in Q\), and in particular \(\gamma <
		\lambda\). The closure of \(Q\) under countable sequences also easily
		implies that \(k(\gamma) = \gamma\). 

		Assume towards a contradiction that \(\gamma^+ < \lambda\). We claim
		that \(k[\gamma]\) is definable over \(H(\lambda)\) from the ordinal
		\(\sup k[\gamma^{+Q}]\) and parameters in \(k[Q]\). This follows from
		the stationary splitting argument, which actually implies that if
		\(\langle T_\alpha: \alpha < \gamma\rangle\) is any stationary splitting
		of \(\{\alpha < \gamma^+ : \text{cf}(\alpha) = \omega\}\) that lies in
		the range of \(k\), then \(k[\gamma] = \{\alpha < \gamma : T_\alpha\cap
		\sup k[\gamma^{+Q}]\text{ is stationary}\}\). We omit the proof.

		We now split into two cases, each of which leads to a contradiction.

		\begin{case}
		\(\gamma^{+Q} = \gamma^+\)
		\end{case}
		In this case, \(\sup k[\gamma^{+Q}] = \gamma^+\), so \(k[\gamma] =
		\{\alpha < \gamma : T_\alpha\text{ reflects to }\gamma^+\} = \gamma\).
		In other words, \(k\restriction \gamma\) is the identity, contrary to
		the fact that \(\gamma > \kappa_0 = \textsc{crt}(k)\). (This is just
		Woodin's proof of the Kunen inconsistency.) Given this contradiction, we
		turn to our second case.

		\begin{case}
		\(\gamma^{+Q} < \gamma^+\).
		\end{case}
		In this case, we will use Solovay's argument that SCH holds above a
		strongly compact cardinal to show that \(\gamma^\omega = \gamma^+\).
		This immediately leads to a contradiction: by elementarity,
		\(\gamma^{+Q} = (\gamma^\omega)^Q\); by the closure of \(Q\) under
		countable sequences \((\gamma^\omega)^Q \geq \gamma^\omega\); and hence
		\(\gamma^{+Q} \geq \gamma^\omega = \gamma^+\) (again using the closure
		of \(Q\) under countable sequences), contrary to our case hypothesis.

		We finish by showing that \(\gamma^\omega = \gamma^+\). Suppose not,
		towards our final contradiction. 

		Let \(U\) be the \(Q\)-ultrafilter on \(\gamma^{+Q}\) derived from \(k\)
		using \(\sup k[\gamma^{+Q}]\). Let \(j : Q\to M\) be the ultrapower
		embedding and \(i : M\to H(\lambda)\) the factor embedding. Since
		\(k[\gamma]\) is definable from elements of \(\{\sup k[\gamma^{+Q}]\}
		\cup k[Q]\subseteq \text{ran}(i)\), we have that \(k[\gamma] = i(S)\)
		for some \(S\in M\). But \(S = i^{-1}[k[\gamma]] = j[\gamma]\). This
		shows that \(j[\gamma]\in M\). 

		Since \(\gamma^+ < k(\gamma^+) < \lambda\) and \(k(\gamma^+)\) is a
		cardinal, \(\gamma^{++} < \lambda\). Therefore every subset of
		\(P_{\omega_1}(\gamma)\) of cardinality \(\gamma^{++}\) belongs to
		\(H(\lambda)\). By elementarity, we can fix a set \(A\subseteq
		P_{\omega_1}(\gamma)\) with \(A\in Q\) and \(|A|^Q = \gamma^{++Q}\). Now
		\(j[A]\in M\): indeed, \[j[A] = \{\sigma\in j(A) : \sigma\subseteq
		j[\gamma]\}\] The forwards inclusion is immediate, and the reverse
		inclusion follows from the fact that \(Q\) is closed under countable
		seqences.

		Now let \(f : A\to \gamma^{++Q}\) be a surjection that lies in \(Q\).
		Then \(j(f)[j[A]] = j[\gamma^{++Q}]\), so \(j[\gamma^{++Q}]\in M\).
		Since \(j(\gamma^{++Q}) > j(\gamma^{+Q}) \geq \gamma^{++Q}\) and
		\(j(\gamma^{++Q}) = \gamma^{++M}\) is an \(M\)-regular cardinal,
		\(j[\gamma^{++Q}]\) cannot be cofinal in \(j(\gamma^{++Q})\). It follows
		that \(j\) is discontinuous at \(\gamma^{++Q}\). This contradicts that
		\(j\) is the ultrapower embedding associated to a \(Q\)-ultrafilter on
		\(\gamma^{+Q}\): in general, the ultrapower embedding of a model \(N\)
		associated to an \(N\)-ultrafilter on an ordinal \(\delta\) is
		continuous at every \(N\)-regular cardinal above \(\gamma\).
	\end{proof}
\end{thm}

\begin{lma}\label{PathoUF}
	Suppose \(\lambda\) is a pathological cardinal of uncountable cofinality and
	\(j : V\to M\) witnesses the pathology of \(\lambda\). Let \(A\in M\) be a
	cover of \(j[\lambda]\) of \(M\)-cardinality \(\lambda\), and let \(\mathcal
	U\) be the fine ultrafilter on \(P(\lambda)\) derived from \(j\) using
	\(A\). Let \(k : M_\mathcal U\to M\) be the factor embedding. Then
	\(\textsc{crt}(k) > \lambda\) and therefore \(j_\mathcal U\) witnesses the
	pathology of \(\lambda\).
	\begin{proof}
		 Let \(k : M_\mathcal U\to M\) be the factor embedding. We must show
		 that \(\textsc{crt}(k) > \lambda\). Let \(\bar A = \id_{\mathcal U}\),
		 so \(k(\bar A) = A\). Clearly \(j_{\mathcal U}[\lambda]\subseteq \bar
		 A\), so \(|\bar A|^{M_{\mathcal U}}\geq |\bar A|\geq \lambda\). On the
		 other hand, \(|\bar A|^{M_{\mathcal U}} \leq k(|\bar A|^{M_{\mathcal
		 U}}) = |A|^M = \lambda\). Thus \(|\bar A|^{M_{\mathcal U}} = \lambda\),
		 so \[k(\lambda) =k(|\bar A|^{M_{\mathcal U}}) = |A|^M = \lambda\]
		
		Assume towards a contradiction that \(\textsc{crt}(k) < \lambda\). Since
		\(j\) is \({<}\lambda\)-supercompact, \(j\) is \({<}\lambda\)-strong,
		and therefore \(H(\lambda)\cap M = H(\lambda)\). Thus \(k\) restricts to
		a nontrivial elementary embedding \(k : H(\lambda)\cap M_{\mathcal U}\to
		H(\lambda)\). Since \(M_{\mathcal U}\) is closed under countable
		sequences, we can apply Foreman's inconsistency theorem.  Since
		\(\lambda\) has uncountable cofinality and \(k(\lambda) = \lambda\),
		\(k\) has a fixed point in the interval \((\textsc{crt}(k),\lambda)\).
		Therefore by Foreman's theorem (\cref{ForemanInconsistency}), \(\lambda
		= \gamma^+\) where \(\gamma\) is the supremum of the critical sequence
		of \(k\). But \(j\) is \(\gamma\)-supercompact and \(j\) is continuous
		at \(\gamma\), so by \cref{SmallCfCompact}, \(j\) is
		\(\gamma^+\)-supercompact. Since \(j\) witnesses the pathology of
		\(\lambda\), \(j\) is not \(\lambda\)-supercompact. This contradicts
		that \(\lambda  = \gamma^+\).
		
		Thus our assumption was false, and in fact \(\textsc{crt}(k) \geq
		\lambda\). Since \(k(\lambda) = \lambda\), it follows that
		\(\textsc{crt}(k) > \lambda\). We finally show that this implies
		\(j_{\mathcal U}\) witnesses the pathology of \(\lambda\). 
		
		The set \(\bar A\) witnesses that \(j_{\mathcal U}\) is
		\(\lambda\)-tight. 
		
		Assume towards a contradiction that \(j_{\mathcal U}\) is
		\(\lambda\)-supercompact. Since \(\textsc{crt}(k) > \lambda\),
		\(k(j_{\mathcal U}[\lambda]) = k\circ j_{\mathcal U}[\lambda] =
		j[\lambda]\), so \(j\) is \(\lambda\)-supercompact, which is a
		contradiction.
		
		We finally show that \(j_{\mathcal U}\) is \({<}\lambda\)-supercompact.
		Since \(j_{\mathcal U}\) is an ultrapower embedding, it suffices to show
		that \(j_{\mathcal U}\) is \(\delta\)-supercompact for all regular
		cardinals \(\delta < \lambda\). To do this, it is enough to show that
		\(j[\delta]\in k[M_{\mathcal U}]\), since then \(k^{-1}(j[\delta]) =
		j_{\mathcal U}[\delta]\) belongs to \(M_{\mathcal U}\). By Solovay's
		Lemma (\cref{StationaryPartition}), \(j[\delta]\) is definable in \(M\)
		from \(\sup j[\delta]\) and parameters in \(j[V]\). Since
		\(j[V]\subseteq k[M_{\mathcal U}]\) and \(k[M_{\mathcal U}]\) is closed
		under definability in \(M\), to show \(j[\delta]\in k[M_{\mathcal U}]\),
		it suffies to show that \(\sup j[\delta]\in k[M_{\mathcal U}]\). To
		finish, we show that \(k(\sup j_{\mathcal U}[\delta]) = \sup
		j[\delta]\), or in other words that \(k\) is continuous at \(\sup
		j_{\mathcal U}[\delta]\). Since \(\textsc{crt}(k) > \lambda\), it is
		enough to show that \(\text{cf}^{M_{\mathcal U}}(\sup j_{\mathcal
		U}[\delta]) \leq \lambda\). Since \(j_{\mathcal U}\) is
		\(\lambda\)-tight, \(j_{\mathcal U}\) is \((\delta,\lambda)\)-tight, so
		by the easy direction of \cref{KetonenCov}, \(\text{cf}^{M_{\mathcal
		U}}(\sup j_{\mathcal U}[\delta]) \leq \lambda\), as desired.
	\end{proof}
\end{lma}

As a corollary, we eliminate many pathologies which a priori might have seemed
plausible:

\begin{cor}\label{PathologicalDichotomy}
	Suppose \(\lambda\) is a pathological cardinal. Either \(\lambda\) is
	regular or \(\lambda\) has countable cofinality.
	\begin{proof}
		Assume \(\lambda\) has uncountable cofinality, and we will show that
		\(\lambda\) is regular. By \cref{PathoUF}, the pathology of \(\lambda\)
		is witnessed by an ultrapower embedding \(i : V\to N\). Since \(i\) is a
		\({<}\lambda\)-supercompact ultrapower embedding, \(N\) is closed under
		\({<}\lambda\)-sequences. If \(\lambda\) is singular, it follows that
		\(N\) is closed under \(\lambda\)-sequences, contradicting that \(i\) is
		not \(\lambda\)-supercompact. Therefore \(\lambda\) is regular.
	\end{proof}
\end{cor}

\begin{cor}\label{PathoReg}
	Suppose \(\lambda\) is a regular pathological cardinal. Suppose \(j : V\to
	M\) witnesses the pathology of \(\lambda\). Let \(U\) be the ultrafilter on
	\(\lambda\) derived from \(j\) using \(\sup j[\lambda]\), and let \(k :
	M_U\to M\) be the factor embedding. Then \(\textsc{crt}(k) >\lambda\) and
	\(j_U\) witnesses the pathology of \(\lambda\).
	\begin{proof}
		Since \(\lambda\) is regular and \(j\) is \(\lambda\)-tight,
		\(\text{cf}^M(\sup j[\lambda]) = \lambda\). Note that \(\id_U = \sup
		j_U[\lambda]\), so \(k(\sup j_U[\lambda]) = \sup j[\lambda]\). We have
		\(\text{cf}^{M_U}(\sup j_U[\lambda])\geq \text{cf}(\sup j_U[\lambda]) =
		\lambda\) on the one hand, and \(\text{cf}^{M_U}(\sup j_U[\lambda])\leq
		k(\text{cf}^{M_U}(\sup j_U[\lambda])) = \text{cf}^M(\sup j[\lambda]) =
		\lambda\) on the other. Thus \(\text{cf}^{M_U}(\sup j_U[\lambda]) =
		\lambda\). It follows that \(k(\lambda) = k(\text{cf}^{M_U}(\sup
		j_U[\lambda])) = \text{cf}^M(\sup j[\lambda]) = \lambda\). 
		
		Given that \(k(\lambda)= \lambda\), one can finish the proof as in
		\cref{PathoUF}. Instead of redoing this proof, however, we note that the
		corollary follows from an {\it application} of \cref{PathoUF}. Using
		\cref{KetonenCov}, fix a cover \(\bar A\subseteq \sup j_U[\lambda]\) of
		\(j_U[\lambda]\) of \(M_U\)-cardinality \(\lambda\). Let \(A = k(\bar
		A)\). Thus \(|A|^M = k(\lambda) = \lambda\). Moreover, it is easy to see
		that 
		\[H^M(j[V]\cup \{\sup j[\lambda]\}) = H^M(j[V]\cup \{A\})\] The
		left-to-right inclusion follows from the fact that \(\sup j[\lambda] =
		\sup A\) is definable from \(A\) in \(M\), while the right-to-left
		inclusion follows from the fact that \(A = k(\bar A)\) and \(k[M_U] =
		H^M(j[V]\cup \{\sup j[\lambda]\}) \). Therefore \(j_U = j_\mathcal U\)
		and the factor embeddings from \(M_\mathcal U\) into \(M\) is equal to
		\(k\). Therefore by \cref{PathoUF}, \(\textsc{crt}(k) >\lambda\) and
		\(j_U\) witnesses the pathology of \(\lambda\).
	\end{proof}
\end{cor}

Pathological cardinals of countable cofinality, on the other hand, have a
property that is a lot like \(I_2(\lambda)\):
\begin{prp}\label{PathoCountable}
	Suppose \(\lambda\) is a pathological cardinal of countable cofinality. Then
	there is a countably complete fine ultrafilter \(\mathcal U\) on
	\(P(\lambda)\) such that there is a nontrivial elementary embedding \(k :
	M_{\mathcal U}\to M\) such that \(k\circ j_\mathcal U = j\) and \(\lambda\)
	is the supremum of the critical sequence of \(k\).
\begin{proof}
	Immediate from the proof of \cref{PathoUF}.
\end{proof}
\end{prp}
If the ultrafilter \(\mathcal U\) of the previous lemma is principal, then
\(I_2(\lambda)\) holds. Under UA, there is a way to make this conclusion:
\begin{thm}[UA] Suppose \(\lambda\) is a pathological cardinal of countable
	cofinality. Then \(I_2(\lambda)\).
	\begin{proof}
		Let \(j : V\to M\) witness the pathology of \(\lambda\). Then \(j\)
		witnesses that some cardinal \(\kappa < \lambda\) is
		\(\gamma\)-supercompact for all \(\gamma < \lambda\). In particular, by
		our results on GCH (\cref{MainThm}), \(\lambda\) is a strong limit
		cardinal.
		
		Applying \cref{PathoCountable}, fix a countably complete fine
		ultrafilter \(\mathcal U\) on \(P(\lambda)\) and a nontrivial elementary
		embedding \(k : M_{\mathcal U}\to M\) such that \(k\circ j_\mathcal U =
		j\) and \(\lambda\) is the supremum of the critical sequence of \(k\).
		
		By \cref{SupercompactFactor}, fix a countably complete ultrafilter \(D\)
		with \(\lambda_D < \lambda\) and an elementary embedding \(k : M_D\to
		M_\mathcal U\) such that \(k\) is \(j_D(\lambda)\)-supercompact in
		\(M_D\). Since \(\lambda\) is a strong limit cardinal of countable
		cofinality, \(j_D(\lambda) = \lambda\). In particular, \(V_\lambda \cap
		M_D = V_\lambda\cap {M_{\mathcal U}}\). Since \((j_D\restriction
		V_\lambda) : V_\lambda\to V_\lambda\cap M_D\) and \((k\restriction
		V_\lambda\cap M_D) : V_\lambda\cap M_D\to V_\lambda\) are elementary
		embeddings, 
		\[i = (k\restriction V_\lambda\cap M_D)\circ (j_D\restriction
		V_\lambda)\] is an elementary embedding from \(V_\lambda\) to
		\(V_\lambda\). Moreover, suppose \(A\subseteq V_\lambda\) is a
		wellfounded relation. Then \(i(A) = \bigcup_{\alpha < \lambda} i(A\cap
		V_\alpha)\) is also wellfounded since \(i(A) = k(j_D(A))\), and \(k\)
		and \(j_D\) preserve wellfoundedness. Thus \(i\) extends to an
		elementary embedding \(i^* : V\to N\) where \(N\) is wellfounded, and it
		follows that \(I_2(\lambda)\) holds.
	\end{proof}
\end{thm}

Under UA, regular pathological cardinals are inaccessible:
\begin{prp}[UA]\label{SuccPathology}
	Suppose \(\lambda\) is a regular pathological cardinal. Then \(\lambda\) is
	strongly inaccessible and \(\mathscr K_\lambda\) witnesses the pathology of
	\(\lambda\).
	\begin{proof}
		By \cref{PathoUF}, there is a countably complete ultrafilter \(U\) such
		that \(j_U\) witnesses the pathology of \(\lambda\). In particular,
		\(j_U\) is \({<}\lambda\)-supercompact and \(\lambda\)-tight. It follows
		that \(j_U\) is \(\lambda\)-pseudocompact, since this just means \(j_U\)
		is \(\gamma\)-tight for all cardinals \(\gamma\leq\lambda\). In
		particular, \(U\) is \(\lambda\)-irreducible by \cref{TightIrred}. 
		
		Note that \(j_U\) witnesses that \(\lambda\) is Fr\'echet. Suppose
		towards a contradiction that \(\lambda\) is a successor cardinal. Then
		by the Irreducibility Theorem (\cref{SuccessorIrredCor}), \(j_U\) is
		\(\lambda\)-supercompact, contradicting that \(U\) witnesses the
		pathology of \(\lambda\).
			
		Thus \(\lambda\) is a limit cardinal. But \(j_U\) is
		\({<}\lambda\)-supercompact, so by our results on GCH (\cref{MainThm}),
		\(\lambda\) is a strong limit cardinal. Therefore \(\lambda\) is
		strongly inaccessible.
		
		Finally we show that \(\mathscr K_\lambda\) witnesses the pathology of
		\(\lambda\). Let \(j : V\to M\) be the ultrapower of the universe by
		\(\mathscr K_\lambda\). It suffices to show that \(j\) is not
		\(\lambda\)-supercompact, since by \cref{KetonenTight}, \(j\) is
		\({<}\lambda\)-supercompact and \(\lambda\)-tight. Suppose towards a
		contradiction that \(j\) is \(\lambda\)-supercompact. Then by
		\cref{KlambdaReduction}, {\it every} ultrapower by a
		\(\lambda\)-irreducible ultrafilter is \(\lambda\)-supercompact,
		contradicting that \(j_U\) is not \(\lambda\)-supercompact. Thus
		\(\mathscr K_\lambda\) witnesses the pathology of \(\lambda\).
	\end{proof}
\end{prp}

To summarize, under UA, if a cardinal is pathological, it is pathological for
good reason:
\begin{thm}[UA] If \(\lambda\) is a pathological cardinal, then one of the
	following holds:
	\begin{itemize}
		\item \(\lambda\) is a strong limit singular cardinal of countable
		cofinality and \(I_2(\lambda)\) holds.
		\item \(\lambda\) is a strongly inaccessible cardinal and \(\mathscr
		K_\lambda\) witnesses the pathology of \(\lambda\).\qed
	\end{itemize}
\end{thm}

We now turn to the question of whether regular pathological  cardinals can exist
at all (without assuming UA). This is equivalent to the existence of
pseudocompact embeddings that are not supercompact:
\begin{prp}\label{PathoPseudo}
	Suppose \(j : V\to M\) is an elementary embedding such that \(M\) has the
	\({\leq}\gamma\)-covering property for all \(\gamma \leq \lambda\). Then one
	of the following holds:
	\begin{itemize}
		\item \(M^\lambda\subseteq M\).
		\item \(j\) witnesses the pathology of a regular cardinal \(\delta \leq
		\lambda\).
	\end{itemize}
	\begin{proof}
		We claim that if \(\gamma \leq \lambda\) and \(j\) is
		\(\gamma\)-supercompact, then \(M^\gamma\subseteq M\). (In fact, it
		suffices that \(P(\gamma)\subseteq M\) which follows from
		\(\gamma\)-supercompactness by \cref{SupercompactClosure}.)  To see
		this, suppose \(A\subseteq M\) and \(|A|\leq \gamma\). Using the
		\({\leq}\gamma\)-covering property, fix \(B\in M\) with \(A\subseteq B\)
		and \(|B|^M \leq \gamma\). Then since \(P(\gamma)\subseteq M\),
		\(P(B)\subseteq M\), and hence \(A\in M\), as desired.
		
		Therefore let \(\delta\) be the least cardinal such that \(j\) is not
		\(\delta\)-supercompact. Note that \(\delta\) is the least cardinal such
		that \(M^\delta\not \subseteq M\), and therefore \(\delta\) is regular.
		If \(\delta \leq \lambda\), then \(j\) witnesses that \(\delta\) is
		pathological. Otherwise \(\delta > \lambda\), and hence \(j\) is
		\(\lambda\)-supercompact, so \(M^\lambda\subseteq M\) by the previous
		paragraph.
	\end{proof}
\end{prp}

Recall Woodin's Ultimate \(L\) Conjecture, which in a weak form states that if
\(\delta\) is extendible then there is an inner model with the
\(\delta\)-approximation and \(\delta\)-covering properties that satisfies the
axiom \(V = \text{Ultimate }L\). The motivation is that the canonical inner
model with a supercompact cardinal should give rise to such an inner model. The
same intuition motivates the {\it UA Conjecture}, which we now define.

\begin{defn}\index{Ultrapower Axiom!UA Hypothesis}
	We say that the {\it UA Hypothesis} holds at a cardinal \(\delta\) if there
	is a inner model of UA with the \(\delta\)-cover and
	\(\delta\)-approximation properties.
\end{defn}

\begin{conj}[UA Conjecture]\index{Ultrapower Axiom!UA Conjecture}
	\textnormal{ZFC} proves that if \(\delta\) is an extendible cardinal, then
	the \textnormal{UA Hypothesis} holds at \(\delta\).
\end{conj}
It is a plausible conjecture that the axiom \(V = \text{Ultimate }L\) implies
UA. If this is the case, then the Ultimate \(L\) Conjecture implies the UA
Conjecture. On the other hand, the UA Conjecture implies the HOD Conjecture.

Our next theorem, due to Woodin in the case that \(\lambda\) is strongly
inaccessible, shows that the pathologies we are studying are in a sense
absolute.
\begin{defn}
	Suppose \(\kappa < \lambda\) are cardinals. Then \(\lambda\) is {\it
	\(\kappa\)-pathological} if there is an elementary embedding \(j : V\to M\)
	with critical point \(\kappa\) that witnesses the pathology of \(\lambda\).
\end{defn}

\begin{thm}\label{WoodinAbsoluteness}
	Suppose \(\delta < \kappa < \lambda\) are cardinals and \(\lambda\) is
	regular. Suppose \(N\) is an inner model with the \(\delta\)-cover and
	\(\delta\)-approximation properties. Suppose \(\lambda\) is
	\(\kappa\)-pathological. Then \(\lambda\) is \(\kappa\)-pathological in
	\(N\).
\end{thm}
The proof uses several facts from the remarkable theory of models with the
approximation and covering properties. First, we will need Hamkins's uniqueness
theorem for models with the approximation and covering properties:
\begin{thm}[Hamkins]\label{UniquenessThm}
	Suppose \(\delta\) is a cardinal and \(N_0\) and \(N_1\) are inner models of
	\textnormal{ZFC} with the \(\delta\)-approximation and \(\delta\)-covering
	properties at an ordinal \(\alpha\). If \(N_0\cap H(\delta^+) = N_1\cap
	H(\delta^+)\) then \(N_0\cap P(\alpha) = N_1\cap P(\alpha)\).\qed
\end{thm}
Second we need the Hamkins-Reitz theorem on the propagation of the covering
property:
\begin{thm}[Hamkins-Reitz]\label{PropagationThm}
	Suppose \(\delta\) is a cardinal and \(N\) is an inner model of
	\textnormal{ZFC} with the \(\delta\)-approximation and \(\delta\)-covering
	properties. Then \(N\) has the \(\lambda\)-covering property for every
	cardinal \(\lambda\geq \delta\).\qed
\end{thm}

With these tools in hand, we can prove \cref{WoodinAbsoluteness}.

\begin{proof}[Proof of \cref{WoodinAbsoluteness}] Applying \cref{PathoReg}, let
		\(U\) be a \(\kappa\)-complete weakly normal ultrafilter on \(\lambda\)
		such that \(j_U\) witnesses the pathology of \(\lambda\). Let \(W =
		U\cap N\). By \cref{ApproxThm}, \(W\) belongs to \(N\) (and in fact
		every \(\delta\)-complete \(N\)-ultrafilter belongs to \(N\)).
		
		Let \(j: V\to M\) be the ultrapower of the universe by \(U\). Let \(i:
		N\to P\) be the ultrapower of \(N\) by \(W\). Let \(k : P\to j(N)\) be
		the factor embedding, defined by \[k(i(f)(\sup i[\lambda])) = j(f)(\sup
		j[\lambda])\] Thus \(k\circ i  = j\restriction N\) and \(k(\sup
		i[\lambda]) = \sup j[\lambda]\).
		
		We now show that \[\text{cf}^{j(N)}(\sup j[\lambda])= \lambda\] Since
		\(j\) is \(\lambda\)-tight, \(\text{cf}^M(\sup j[\lambda]) = \lambda\).
		Since \(j(N)\) has the \(\delta\)-approximation and \(\delta\)-covering
		properties in \(M\), in fact \(j(N)\) has the \(\lambda\)-covering
		property in \(M\) (by \cref{PropagationThm}).  Therefore \(j(N)\)
		correctly computes the cofinality of \(\sup j[\lambda]\) in \(M\), and
		it follows that \(\text{cf}^{j(N)}(\sup j[\lambda])= \lambda\).
		
		We now claim that \(k(\lambda) = \lambda\) and \(\text{cf}^N(\sup
		i[\lambda]) = \lambda\). The argument is by now familiar. Since \(k(\sup
		i[\lambda]) = \sup j[\lambda]\), 
		\[k(\text{cf}^N(\sup i[\lambda])) = \text{cf}^{j(N)}(\sup j[\lambda]) =
		\lambda\] Since \(\lambda \leq \text{cf}^N(\sup i[\lambda])\), \(\lambda
		\leq k(\lambda) \leq k(\text{cf}^N(\sup i[\lambda])) = \lambda\). Thus
		\(k(\lambda) = \lambda\). Similarly \(\lambda \leq \text{cf}^N(\sup
		i[\lambda]) \leq k(\text{cf}^N(\sup i[\lambda]))\leq \lambda\), so
		\(\text{cf}^N(\sup i[\lambda])  = \lambda\).
		
		We claim \(j(N)\cap P(\alpha) = N\cap P(\alpha)\) for all \(\alpha <
		\lambda\). The argument is due to Hamkins \cite{Hamkins}. Fix \(\alpha <
		\lambda\). Since \(j\) is \({<}\lambda\)-supercompact, \(M\cap P(\alpha)
		= P(\alpha)\). By elementarity, \(j(N)\) has the
		\(\delta\)-approximation and \(\delta\)-covering properties in \(M\),
		and in particular \(j(N)\) has the \(\delta\)-approximation and
		\(\delta\)-covering properties at \(\alpha\). Similarly, \(N\) has the
		\(\delta\)-approximation and \(\delta\)-covering properties at
		\(\alpha\). But \(N\cap H(\delta^+) = j(N)\cap H(\delta^+)\) since
		\(\textsc{crt}(j) > \delta\). Thus \(j(N)\cap P(\alpha) = N\cap
		P(\alpha)\) by the uniqueness theorem (\cref{UniquenessThm}).
		
		We claim \(j[\alpha]\in j(N)\) for all \(\alpha < \lambda\). This
		follows from the \(\delta\)-approximation property for \(j(N)\) in \(M\)
		and the \(\delta\)-covering property for \(N\) in \(V\). Let \(\alpha_*
		= \sup j[\alpha]\). Suppose \(\sigma\in P_\delta(\alpha_*)\cap j(N)\).
		Fix \(\tau\in P_\delta(\alpha)\) such that \(j^{-1}[\sigma]\subseteq
		\tau\). Then \(\sigma\cap j[\alpha] = \sigma\cap j[\tau] = \sigma\cap
		j(\tau)\in j(N)\).
		
		Since \(P(\alpha)\cap N = P(\alpha)\cap j(N)\) for all \(\alpha <
		\lambda\), we have that \(H(\lambda)\cap N = H(\lambda) \cap j(N)\). Let
		\(H = H(\lambda)\cap N\) and let \(Q = H(\lambda) \cap P\). We claim
		that \(k\restriction Q\in N\). The proof is a generalization of the
		proof of Woodin's Universality Theorem for models with the approximation
		property. Since \(Q\) is transitive, \(k\restriction Q\) is the inverse
		of the transitive collapse of \(k[Q]\), and therefore it suffices to
		show that \(k[Q]\in N\). Since \(N\) satisfies the
		\(\delta\)-approximation property, it suffices to show that \(k[Q]\cap
		\sigma\in N\) for any \(\sigma \in P_\delta(H)\cap N\). Fix such a
		\(\sigma\). Since \(N\) has the \(\delta\)-cover property, there is some
		\(\tau\in P_\delta(Q)\cap N\) with \(k^{-1}[\sigma]\subseteq \tau\).
		Since \(P^\delta\cap N\subseteq P\), \(\tau\in P\), and hence \(\tau \in
		Q\). Since \(\textsc{crt}(k) \geq \delta\), \(k(\tau) = k[\tau]\). Thus
		\[k[Q]\cap \sigma = k[\tau]\cap \sigma = k(\tau) \cap \sigma \in N\]
		Thus \(k[Q]\cap \sigma\in N\). By the \(\delta\)-approximation property,
		\(k[Q]\in N\), and hence \(k\restriction Q\in N\).
		
		We now apply Foreman's inconsistency result in \(N\). Assume towards a
		contradiction that \(\textsc{crt}(k) < \lambda\). Note that \(k\)
		restricts to an elementary embedding from \(Q\) to \(H(\lambda)\cap N\)
		that belongs to \(N\). Moreover \(Q\) is closed under
		\(\omega\)-sequences in \(N\). Then it follows from
		\cref{ForemanInconsistency} applied in \(N\) that \(\lambda = \gamma^+\)
		where \(\gamma\) is the supremum of the critical sequence of \(k\). It
		follows that \(\lambda\) is the successor of a singular cardinal
		\(\gamma\) of countable cofinality in \(N\).
		
		Recall that \(j[\gamma]\in j(N)\) since \(j[\alpha]\in j(N)\) for all
		\(\alpha < \lambda\). Note that \[j[P^N_\kappa(\gamma)] = \{ j[\sigma] :
		\sigma\in P^N_\kappa(\gamma)\} = P^{ j(N)}_\kappa(j[\gamma])\] In
		particular, \(j[P^N_\kappa(\gamma)]\in j(N)\) since it is definable over
		\(j(N)\) from \(j[\gamma]\in j(N)\). Recall that \(\lambda  =
		\gamma^{+N}\). Therefore by K\"onig's theorem, there is a surjection \(f
		: P^N_\kappa(\gamma)\to \lambda\) in \(N\). Then
		\(j(f)[j[P^N_\kappa(\gamma)]] = j[\lambda]\). Thus \(j[\lambda]\in
		j(N)\). In particular, \(j[\lambda]\in M\). This contradicts the fact
		that \(j\) is not \(\lambda\)-supercompact.
		
		Thus our assumption that \(\textsc{crt}(k) < \lambda\) was false.
		Moreover since \(k(\lambda) = \lambda\), it follows that
		\(\textsc{crt}(k) >\lambda\).
		
		We established the following:
		\begin{itemize}
			\item \(k : P \to j(N)\) is an elementary embedding with critical
			point above \(\lambda\)
			\item \(j[\alpha]\in j(N)\) for all \(\alpha < \lambda\).
			\item \(j[\lambda]\notin j(N)\).
		\end{itemize}
		Since \(\textsc{crt}(k)>\lambda\), \(\textsc{crt}(i) = \kappa\). The
		proof that \(i[\alpha]\in P\) for all \(\alpha < \lambda\), and
		\(i[\lambda]\notin P\) now proceeds exactly as in \cref{PathoUF}. Thus
		\(i\) witnesses that \(\lambda\) is \(\kappa\)-pathological in \(N\).
		\end{proof}

The UA Conjecture rules out certain kinds of pathological cardinals that are not
obviously ruled out in ZFC alone:

\begin{thm}\label{UAHNonPathology}
	Suppose \(\delta\) is a cardinal and the \textnormal{UA Hypothesis} holds at
	\(\delta\). If \(\gamma\) is a singular cardinal, then \(\gamma^+\) is not
	\(\kappa\)-pathological for any \(\kappa > \delta\).
	\begin{proof}
		Let \(N\) be an inner model of UA with the \(\delta\)-covering and
		\(\delta\)-approximation properties. Assume towards a contradiction that
		\(\gamma^+\) is \(\kappa\)-pathological for some \(\kappa> \delta\). By
		\cref{WoodinAbsoluteness}, \(\gamma^+\) is pathological in \(N\). But by
		\cref{PropagationThm}, \(\gamma^+\) is a successor cardinal in \(N\):
		this follows from the fact that every ordinal in the interval
		\([\gamma,\gamma^+]\) has cofinality less than \(\gamma\) in \(N\) by
		the \(\gamma\)-covering property. Since \(N\) satisfies UA, by
		\cref{SuccPathology}, no successor cardinal is pathological in \(N\).
		This is a contradiction.
	\end{proof}
\end{thm}

\begin{qst}[ZFC] Can the successor of a singular cardinal be pathological?
\end{qst}

Finally let us tie this all back up with the question of whether UA implies that
\(\lambda\)-irreducible ultrafilters are \(\lambda\)-supercompact.

\begin{thm}
	Suppose \(\delta\) is a cardinal and assume there is an inner model \(N\)
	with the \(\delta\)-approximation and \(\delta\)-covering properties that
	has no regular pathological cardinals above \(\delta\). Suppose \(\lambda>
	\delta\) is a cardinal. Suppose \(j :V\to M\) is an elementary embedding
	with critical point above \(\delta\) such that \(M\) has the
	\({\leq}\gamma\)-covering property for all \(\gamma \leq \lambda\). Then
	\(M^\lambda \subseteq M\).
	\begin{proof}
		Suppose not. Then by \cref{PathoPseudo}, \(j\) witnesses that some
		regular cardinal \(\gamma \leq\lambda\) is
		\(\textsc{crt}(j)\)-pathological. By \cref{WoodinAbsoluteness},
		\(\gamma\) is  \(\textsc{crt}(j)\)-pathological in \(N\), contrary to
		our assumption that \(N\) has no regular pathological cardinals above
		\(\delta\).
	\end{proof}
\end{thm}

Thus granting the UA Conjecture, either pseudocompactness principles are
(eventually) equivalent to supercompactness or else UA is consistent with
regular pathological cardinals. it seems more reasonable to make the following
conjecture:
\begin{conj}
	\textnormal{UA} is consistent with the existence of a regular pathological
	cardinal.
\end{conj}

\printindex
\bibliography{Bibliography}{}

\begin{thebibliography}{10}

\bibitem{Kanamori}
Robert~M. Solovay, William~N. Reinhardt, and Akihiro Kanamori.
\newblock Strong axioms of infinity and elementary embeddings.
\newblock {\em Ann. Math. Logic}, 13(1):73--116, 1978.

\bibitem{Schimmerling}
Ernest Schimmerling.
\newblock Combinatorial principles in the core model for one {W}oodin cardinal.
\newblock {\em Ann. Pure Appl. Logic}, 74(2):153--201, 1995.

\bibitem{KunenLU}
Kenneth Kunen.
\newblock Some applications of iterated ultrapowers in set theory.
\newblock {\em Ann. Math. Logic}, 1:179--227, 1970.

\bibitem{Menas}
Telis~K. Menas.
\newblock On strong compactness and supercompactness.
\newblock {\em Ann. Math. Logic}, 7:327--359, 1974/75.

\bibitem{Mitchell}
William~J. Mitchell.
\newblock Sets constructible from sequences of ultrafilters.
\newblock {\em J. Symbolic Logic}, 39:57--66, 1974.

\bibitem{Neeman}
Itay Neeman.
\newblock Inner models in the region of a {W}oodin limit of {W}oodin cardinals.
\newblock {\em Ann. Pure Appl. Logic}, 116(1-3):67--155, 2002.

\bibitem{Schlutzenberg}
Farmer~Salamander Schlutzenberg.
\newblock {\em Measures in mice}.
\newblock ProQuest LLC, Ann Arbor, MI, 2007.
\newblock Thesis (Ph.D.)--University of California, Berkeley.

\bibitem{Larson}
Paul~B. Larson.
\newblock {\em The stationary tower}, volume~32 of {\em University Lecture
  Series}.
\newblock American Mathematical Society, Providence, RI, 2004.
\newblock Notes on a course by W. Hugh Woodin.

\bibitem{MitchellSteel}
William~J. Mitchell and John~R. Steel.
\newblock {\em Fine structure and iteration trees}, volume~3 of {\em Lecture
  Notes in Logic}.
\newblock Springer-Verlag, Berlin, 1994.

\bibitem{Woodin}
W.~Hugh Woodin.
\newblock In search of {U}ltimate-{$L$}: the 19th {M}idrasha {M}athematicae
  {L}ectures.
\newblock {\em Bull. Symb. Log.}, 23(1):1--109, 2017.

\bibitem{Ketonen}
Jussi Ketonen.
\newblock Strong compactness and other cardinal sins.
\newblock {\em Ann. Math. Logic}, 5:47--76, 1972/73.

\bibitem{JechStationary}
Thomas Jech.
\newblock Stationary subsets of inaccessible cardinals.
\newblock In {\em Axiomatic set theory ({B}oulder, {C}olo., 1983)}, volume~31
  of {\em Contemp. Math.}, pages 115--142. Amer. Math. Soc., Providence, RI,
  1984.

\bibitem{Ketonen2}
Jussi Ketonen.
\newblock Open problems in the theory of ultrafilters.
\newblock In {\em Essays on mathematical and philosophical logic ({P}roc.
  {F}ourth {S}candinavian {L}ogic {S}ympos. and {F}irst {S}oviet-{F}innish
  {L}ogic {C}onf., {J}yv\"askyl\"a, 1976)}, volume 122 of {\em Synthese
  Library}, pages 227--247. Reidel, Dordrecht-Boston, Mass., 1979.

\bibitem{MitchellHandbook}
William~J. Mitchell.
\newblock Beginning inner model theory.
\newblock In {\em Handbook of {S}et {T}heory. {V}ols. 1, 2, 3}, pages
  1449--1495. Springer, Dordrecht, 2010.

\bibitem{Steel}
John~R. Steel.
\newblock An outline of inner model theory.
\newblock In {\em Handbook of {S}et {T}heory. {V}ols. 1, 2, 3}, pages
  1595--1684. Springer, Dordrecht, 2010.

\bibitem{CummingsHandbook}
James Cummings.
\newblock Iterated forcing and elementary embeddings.
\newblock In {\em Handbook of {S}et {T}heory. {V}ols. 1, 2, 3}, pages 775--883.
  Springer, Dordrecht, 2010.

\bibitem{Hamkins}
Joel~D. Hamkins.
\newblock Tall cardinals.
\newblock {\em MLQ Math. Log. Q.}, 55(1):68--86, 2009.

\bibitem{MagidorBagaria}
Joan Bagaria and Menachem Magidor.
\newblock On {$\omega_1$}-strongly compact cardinals.
\newblock {\em J. Symb. Log.}, 79(1):266--278, 2014.

\bibitem{Gitik}
Moti Gitik.
\newblock Some constructions of ultrafilters over a measurable cardinal.
\newblock To appear.

\bibitem{VanWesep}
Robert Van~Wesep.
\newblock Wadge degrees and descriptive set theory.
\newblock In {\em Cabal {S}eminar 76--77 ({P}roc. {C}altech-{UCLA} {L}ogic
  {S}em., 1976--77)}, volume 689 of {\em Lecture Notes in Math.}, pages
  151--170. Springer, Berlin, 1978.

\bibitem{Solovay}
Robert~M. Solovay.
\newblock Strongly compact cardinals and the {GCH}.
\newblock In {\em Proceedings of the {T}arski {S}ymposium ({P}roc. {S}ympos.
  {P}ure {M}ath., {V}ol. {XXV}, {U}niv. {C}alifornia, {B}erkeley, {C}alif.,
  1971)}, pages 365--372. Amer. Math. Soc., Providence, R.I., 1974.

\bibitem{Kunen}
Kenneth Kunen.
\newblock Elementary embeddings and infinitary combinatorics.
\newblock {\em J. Symbolic Logic}, 36:407--413, 1971.

\bibitem{Kanamori2}
Akihiro Kanamori.
\newblock {\em The {H}igher {I}nfinite}.
\newblock Springer Monographs in Mathematics. Springer-Verlag, Berlin, second
  edition, 2009.
\newblock Large cardinals in set theory from their beginnings, Paperback
  reprint of the 2003 edition.

\bibitem{SteelWellfounded}
J.~R. Steel.
\newblock The well-foundedness of the {M}itchell order.
\newblock {\em J. Symbolic Logic}, 58(3):931--940, 1993.

\bibitem{Negrepontis}
W.~W. Comfort and S.~Negrepontis.
\newblock {\em The theory of ultrafilters}.
\newblock Springer-Verlag, New York-Heidelberg, 1974.
\newblock Die Grundlehren der mathematischen Wissenschaften, Band 211.

\bibitem{Donder}
Hans-Dieter Donder.
\newblock Regularity of ultrafilters and the core model.
\newblock {\em Israel J. Math.}, 63(3):289--322, 1988.

\bibitem{AbrahamMagidor}
Uri Abraham and Menachem Magidor.
\newblock Cardinal arithmetic.
\newblock In {\em Handbook of {S}et {T}heory. {V}ols. 1, 2, 3}, pages
  1149--1227. Springer, Dordrecht, 2010.

\bibitem{Frolik}
Zden\v{e}k Frol\'ik.
\newblock Sums of ultrafilters.
\newblock {\em Bull. Amer. Math. Soc.}, 73:87--91, 1967.

\bibitem{KunenRK}
Kenneth Kunen.
\newblock Ultrafilters and independent sets.
\newblock {\em Trans. Amer. Math. Soc.}, 172:299--306, 1972.

\bibitem{LevySolovay}
A.~L\'{e}vy and R.~M. Solovay.
\newblock Measurable cardinals and the continuum hypothesis.
\newblock {\em Israel J. Math.}, 5:234--248, 1967.

\bibitem{KunenText}
Kenneth Kunen.
\newblock {\em Set theory}, volume~34 of {\em Studies in Logic (London)}.
\newblock College Publications, London, 2011.

\bibitem{Magidor}
Menachem Magidor.
\newblock How large is the first strongly compact cardinal? or {A} study on
  identity crises.
\newblock {\em Ann. Math. Logic}, 10(1):33--57, 1976.

\bibitem{Prikry}
Kenneth Kunen and Karel Prikry.
\newblock On descendingly incomplete ultrafilters.
\newblock {\em J. Symbolic Logic}, 36:650--652, 1971.

\bibitem{MitchellSkies}
William Mitchell.
\newblock Indiscernibles, skies, and ideals.
\newblock In {\em Axiomatic set theory ({B}oulder, {C}olo., 1983)}, volume~31
  of {\em Contemp. Math.}, pages 161--182. Amer. Math. Soc., Providence, RI,
  1984.

\bibitem{Apter}
Arthur~W. Apter and Saharon Shelah.
\newblock On the strong equality between supercompactness and strong
  compactness.
\newblock {\em Trans. Amer. Math. Soc.}, 349(1):103--128, 1997.

\bibitem{GitikClubCompact}
Moti Gitik.
\newblock A note on sequences witnessing singularity, following {M}agidor and
  {S}inapova.
\newblock {\em MLQ Math. Log. Q.}, 64(3):249--253, 2018.

\bibitem{Caicedo}
Andr\'es~Eduardo Caicedo.
\newblock Cardinal preserving elementary embeddings.
\newblock In {\em Logic {C}olloquium 2007}, volume~35 of {\em Lect. Notes
  Log.}, pages 14--31. Assoc. Symbol. Logic, La Jolla, CA, 2010.

\bibitem{VickersWelch}
J.~Vickers and P.~D. Welch.
\newblock On elementary embeddings from an inner model to the universe.
\newblock {\em J. Symbolic Logic}, 66(3):1090--1116, 2001.

\end{thebibliography}
\bibliographystyle{unsrt}
\end{document}